\documentclass[10pt]{amsart}

\usepackage{xr-hyper}

\usepackage[utf8]{inputenc}         %
\usepackage[T1]{fontenc}            %

\usepackage{a4wide}                 %
\usepackage{amsmath, amssymb}       %
\usepackage{amsthm}                 %
\usepackage{appendix}               %

\usepackage{tikz}                   %
\usepackage{tikz-cd}                %
\tikzcdset{ %
arrow style=tikz
}
\usetikzlibrary{arrows}             %

\usepackage{graphicx}               %
\usepackage{stackrel}               %

\usepackage{color}                  %
\usepackage[normalem]{ulem}         %
\usepackage{enumitem}               %
\setlist[itemize]{leftmargin=*}
\setlist[enumerate]{leftmargin=*,label=\roman*),ref=\roman*)}
\newlist{subenumerate}{enumerate}{2}
\setlist[subenumerate]{leftmargin=*,label=\alph*),ref=\alph*)}

\usepackage{blkarray}               %

\definecolor{darkblue}{rgb}{0,0,0.6} %
\usepackage[ocgcolorlinks,colorlinks=true, citecolor=darkblue, filecolor=darkblue, linkcolor=darkblue, urlcolor=darkblue]{hyperref}         %

\usepackage[alphabetic]{amsrefs}

\usepackage{stix}                   %
\usepackage{bbm}                    %
\usepackage{eucal}                  %

\DeclareFontFamily{T1}{cbgreek}{}
\DeclareFontShape{T1}{cbgreek}{m}{n}{<-6>  grmn0500 <6-7> grmn0600 <7-8> grmn0700 <8-9> grmn0800 <9-10> grmn0900 <10-12> grmn1000 <12-17> grmn1200 <17-> grmn1728}{}
\DeclareSymbolFont{quadratics}{T1}{cbgreek}{m}{n}
\DeclareMathSymbol{\qoppa}{\mathord}{quadratics}{19}
\DeclareMathSymbol{\Qoppa}{\mathord}{quadratics}{21}

\makeatletter
\let\old@font@info\@font@info
\def\@font@info#1{%
\expandafter\ifx\csname\detokenize{#1}\endcsname\relax
  \old@font@info{#1}%
\fi
\expandafter\xdef\csname\detokenize{#1}\endcsname{}%
}
\makeatother

\usepackage{comment}                %
\usepackage{xspace}                 %

\newtheoremstyle{thms}
	{}{}{\itshape}{}{\bfseries }{.}{ }
	{\thmname{#1} \thmnumber{#2}. \thmnote{\bfseries{[#3]}}}
\newtheoremstyle{thms2}
	{}{}{\itshape}{}{\bfseries }{.}{ }
	{}
\newtheoremstyle{ithreethms}
	{}{}{\itshape}{}{\bfseries }{}{ }
	{\thmname{#1} \thmnumber{#2}. \thmnote{\bfseries{[#3]}}}

\newtheoremstyle{name}
	{}{}{\itshape}{}{\bfseries }{.}{ }
	{\thmname{#1}\thmnumber{#2}\thmnote{\bfseries{[#3]}}}

\newtheoremstyle{defs}
	{}{}{\normalfont}{}{\bfseries }{.}{ }
	{\thmname{#1} \thmnumber{#2}. \thmnote{\bfseries{(#3)}}}
\newtheoremstyle{defs2}
	{}{12pt}{\normalfont}{}{\bfseries }{.}{ }
	{\thmname{#1}\thmnumber{#2}. \thmnote{\bfseries{(#3)}}}
\newtheoremstyle{ithreedefs}
	{}{}{\normalfont}{}{\bfseries }{}{ }
	{\thmname{#1} \thmnumber{#2}. \thmnote{\bfseries{(#3)}}}

\newtheoremstyle{rmk}
	{}{}{\normalfont}{}{\itshape }{.}{ }
        {}

\newtheoremstyle{claim}
	{}{}{\normalfont}{}{\itshape}{.}{ }
        {\thmname{#1} \thmnumber{#2}. \thmnote{#3}}

\theoremstyle{thms2}

\theoremstyle{thms2}
\newtheorem*{ucorollary}{Corollary}

\newtheorem*{namedthm}{\namedthmname}
\newcounter{namedthm}

\makeatletter
\newenvironment{named}[1]
{\def\namedthmname{#1}%
\refstepcounter{namedthm}%
\namedthm\def\@currentlabel{#1}}%
{\endnamedthm}
\makeatother

\newtheorem{theoremintrotwo}{Theorem}

\newtheorem{corollaryintrotwo}[theoremintrotwo]{Corollary}

\theoremstyle{rmk}

\theoremstyle{ithreethms}

\theoremstyle{ithreedefs}

\theoremstyle{thms2}

\swapnumbers         %

\newcounter{rthree}
\theoremstyle{defs}
\newtheorem{definition-r-three}{Definition}[rthree]
\newtheorem{notation-r-three}[definition-r-three]{Notation}
\newtheorem{remark-r-three}[definition-r-three]{Remark}
\newtheorem{example-r-three}[definition-r-three]{Example}

\theoremstyle{thms}
\newtheorem{proposition-r-three}[definition-r-three]{Proposition}
\newtheorem{corollary-r-three}[definition-r-three]{Corollary}
\newtheorem{theorem-r-three}[definition-r-three]{Theorem}

\newcounter{rfour}
\theoremstyle{defs}
\newtheorem{definition-r-four}{Definition}[rthree]
\newtheorem{notation-r-four}[definition-r-four]{Notation}
\newtheorem{remark-r-four}[definition-r-four]{Remark}
\newtheorem{example-r-four}[definition-r-four]{Example}

\theoremstyle{thms}
\newtheorem{proposition-r-four}[definition-r-four]{Proposition}
\newtheorem{corollary-r-four}[definition-r-four]{Corollary}
\newtheorem{theorem-r-four}[definition-r-four]{Theorem}

\theoremstyle{thms}
\newtheorem{proposition}{Proposition}[subsection]
\newtheorem{theorem}[proposition]{Theorem}
\newtheorem{lemma}[proposition]{Lemma}
\newtheorem{corollary}[proposition]{Corollary}

\newtheorem{observation}[proposition]{Observation}

\theoremstyle{defs}
\newtheorem{definition}[proposition]{Definition}%
\newtheorem{notation}[proposition]{Notation}

\newtheorem{construction}[proposition]{Construction}
\newtheorem{example}[proposition]{Example}
\newtheorem{examples}[proposition]{Examples}
\newtheorem{remark}[proposition]{Remark}
\newtheorem{remarks}[proposition]{Remarks}

\theoremstyle{defs2}

\theoremstyle{rmk}
\newtheorem*{Rmk}{Remark}
\newtheorem*{warn}{Warning}

\theoremstyle{claim}

\newcommand{\defi}[1]{\emph{#1}}                 %

\newcommand{\spacefont}{\mathcal}                %

\newcommand{\alp}{\alpha}
\newcommand{\bet}{\beta}

\newcommand{\Del}{\Delta}

\newcommand{\Lam}{\Lambda}
\newcommand{\Om}{\Omega}

\newcommand{\Sig}{\Sigma}
\newcommand{\sig}{\sigma}
\newcommand{\vphi}{\varphi}

\newcommand{\lrar}{\longrightarrow}
\newcommand{\lto}{\longrightarrow}

\newcommand{\hrar}{\hookrightarrow}

\newcommand{\wtl}{\widetilde}

\newcommand{\st}{\stackrel}
\newcommand{\ovl}{\overline}

\newcommand{\adj}{\mathbin{%
\begin{tikzpicture}[baseline,thick] 
\coordinate (source) at (0ex,.5ex);
\coordinate (target) at (3ex,.5ex);
\draw[->] ([yshift=1ex]source) -- ([yshift=1ex]target); 
\draw[->] ([yshift=-.5ex]target) -- ([yshift=-.5ex]source);
\node at (1.5ex,.8ex) {$\scriptscriptstyle \perp$};
\end{tikzpicture}%
}}

\newcommand{\bigcircle}{\bigcirc}  %

\newcommand{\oneover}[1]{{\scriptstyle \frac{1}{#1}}}   %
\newcommand{\onehalf}{\oneover{2}}   %
\newcommand{\oneoverl}{\oneover{l}}   %
\newcommand{\comma}{,}                     %

\newcommand{\myperp}{\scalebox{0.5}{\rotatebox{270}{$\dashv$}}} %

\newcommand{\cocolon}{\nobreak \mskip6mu plus1mu \mathpunct{}\nonscript\mkern-\thinmuskip {:}\mskip2mu \relax}

\makeatletter
\newcommand{\dovl}[1]{\overline{\dbl@overline{#1}}}
\newcommand{\dbl@overline}[1]{\mathpalette\dbl@@overline{#1}}
\newcommand{\dbl@@overline}[2]{%
  \begingroup
  \sbox\z@{$\m@th#1\overline{#2}$}%
  \ht\z@=\dimexpr\ht\z@-2\dbl@adjust{#1}\relax
  \box\z@
  \ifx#1\scriptstyle\kern-\scriptspace\else
  \ifx#1\scriptscriptstyle\kern-\scriptspace\fi\fi
  \endgroup
}
\newcommand{\dbl@adjust}[1]{%
  \fontdimen8
  \ifx#1\displaystyle\textfont\else
  \ifx#1\textstyle\textfont\else
  \ifx#1\scriptstyle\scriptfont\else
  \scriptscriptfont\fi\fi\fi 3
}
\makeatother

\newcommand{\NN}{\mathbb{N}}               %
\newcommand{\RR}{\mathbb{R}}               %
\newcommand{\CC}{\mathbb{C}}               %
\newcommand{\ZZ}{\mathbb{Z}}               %
\newcommand{\QQ}{\mathbb{Q}}               %

\newcommand{\im}{\mathrm{im}}              %
\newcommand{\coker}{\mathrm{coker}}        %
\newcommand{\Tor}{\mathrm{Tor}}            %

\newcommand{\spec}{\mathrm{spec}}          %

\newcommand{\Hom}{\operatorname{Hom}}            %

\newcommand{\fib}{\operatorname{fib}}            %
\newcommand{\cof}{\operatorname{cof}}            %

\newcommand{\colim}{\mathop{\mathrm{colim}}}     %

\newcommand{\Map}{\operatorname{Map}}            %

\newcommand{\grpcr}{\iota}                       %
\newcommand{\core}{\mathrm{Cr}}                  %
\newcommand{\Core}{\core}                        %

\newcommand{\cor}{\mathrm{Cor}}                  %

\newcommand{\dec}{\mathrm{dec}}                  %

\newcommand{\lag}{\operatorname{lag}}            %

\newcommand{\idem}{\mathrm{idem}}                %
\newcommand{\mn}{\mathrm{min}}                  %

\newcommand{\nerv}{\operatorname N}              %

\newcommand{\Pro}{\operatorname{Pro}}            %
\newcommand{\Ind}{\operatorname{Ind}}            %
\newcommand{\Tate}{\operatorname{Tate}}          %
\newcommand{\Latt}{\operatorname{Latt}}          %

\newcommand{\op}{^\mathrm{op}}                   %

\newcommand{\Ar}{\operatorname{Ar}}              %
\newcommand{\Twar}{\operatorname{TwAr}^{\mathrm{r}}}          %
\newcommand{\target}{\mathrm{t}}                 %
\newcommand{\source}{\mathrm{s}}                 %

\newcommand{\cp}{\omega}                         %
\newcommand{\perf}{\mathrm{p}}                   %

\newcommand{\Mon}{\operatorname{Mon}}            %
\newcommand{\Grp}{\operatorname{Grp}}            %
\newcommand{\grp}{\mathrm{grp}}                  %

\newcommand{\Bs}{\mathrm{B}}                     %

\newcommand{\sd}{\operatorname{sd}}              %

\newcommand{\asscat}{\operatorname{asscat}}      %
\newcommand{\snerv}{\operatorname{N}}            %
\newcommand{\comp}{\operatorname{comp}}	         %

\newcommand{\pr}{\mathrm{pr}}                %
\newcommand{\incl}{\mathrm{inc}}             %
\newcommand{\id}{\mathrm{id}}                %
\newcommand{\const}{\mathrm{const}}          %

\newcommand{\Spa}{{\mathcal Sp}}             %
\newcommand{\Spaf}{\Spa^{\free}}             %
\newcommand{\PSp}{{\mathcal {PS}p}}          %
\newcommand{\Sps}{\mathcal S}                %
\newcommand{\ssmash}{\otimes}                  %
\renewcommand{\hom}{\operatorname{hom}}        %
\newcommand{\map}{\operatorname{hom}}          %

\renewcommand{\SS}{\mathbb{S}}                 %

\newcommand{\sph}{\mathrm{S}}                  %

\newcommand{\solid}{localisation\xspace}       %
\newcommand{\compatible}{compatible\xspace}    %

\newcommand{\h}{\mathrm{h}}                    %
\newcommand{\Ct}{\mathrm{C_2}}                 %
\newcommand{\BC}{\Bs \Ct}                      %
\newcommand{\hC}{{\h\Ct}}                      %
\newcommand{\tC}{{\mathrm{t}\Ct}}              %
\newcommand{\undl}{\mathrm{u}}                 %
\newcommand{\geofix}{{\varphi\Ct}}             %
\newcommand{\gC}{{\mathrm{g}\Ct}}              %
\newcommand{\Spagc}{\Spa^\gC}                  %
\newcommand{\gCt}{\gC}                         %

\newcommand{\Ab}{\mathcal{A}b}                   %

\newcommand{\Ring}{\mathrm{Ring}}                %

\newcommand{\Mod}{\operatorname{Mod}}            %

\newcommand{\fg}{{\mathrm{f}}}                   %

\newcommand{\Proj}{\mathrm{Proj}}                %

\newcommand{\Modcp}{\Mod^{\cp}}                  %

\newcommand{\Modp}[1]{\Modcp({#1})}              %
\newcommand{\Modf}[1]{\mathrm{Mod}^\fg({#1})}    %
\newcommand{\Modc}[2]{\mathrm{Mod}^{#1}({#2})}   %

\newcommand{\Unimod}{\mathrm{Unimod}}            %

\newcommand{\Ering}{\mathrm{E}}			 %
\newcommand{\Ek}[1]{{\Ering_{#1}}}               %
\newcommand{\Einf}{\Ek{\infty}}          	 %
\newcommand{\Eone}{\Ek{1}}               	 %

\newcommand{\Alg}{\mathrm{Alg}}                  %

\newcommand{\Ch}{\operatorname{Ch^b}}            %

\newcommand{\qiso}{\mathrm{qIso}}                %

\newcommand{\GEM}{\mathrm{H}}                    %

\newcommand{\Der}{\mathcal D}                    %
\newcommand{\Dfree}{{\mathcal D}^\fg}            %
\newcommand{\Dperf}{{\mathcal D}^\perf}          %

\newcommand{\free}{\mathrm{f}}                   %

\newcommand{\Set}{\mathrm{Set}}                  %

\newcommand{\Cat}{\mathrm{Cat}_\infty}                  %

\newcommand{\Catx}{\Cat^\ex}     %
\newcommand{\Catxi}{\mathrm{Cat}^{\ex,\natural}_{\infty}}                 %

\newcommand{\sCatx}{\mathrm{s}\Catx}             %

\newcommand{\sSps}{\mathrm{s}\Sps}	 %
\newcommand{\cSS}{\mathrm{cS}\Sps}	         %

\newcommand{\Cath}{\mathrm{Cat}^\mathrm{h}_\infty}      %
\newcommand{\sCath}{\mathrm{s}\Cath} 	                %
\newcommand{\Catp}{\mathrm{Cat}^{\mathrm p}_\infty}     %
\newcommand{\Catpi}{\mathrm{Cat}^{\mathrm p,\natural}_{\infty}}                 %

\newcommand{\sCatp}{\mathrm{s}\Catp}             %

\newcommand{\Span}{\operatorname{Span}}          %
\newcommand{\Hlgy}{\operatorname{Hlgy}}          %
\newcommand{\Cob}{\mathrm{Cob}}                  %

\newcommand{\Fun}{\operatorname{Fun}}            %
\newcommand{\Nat}{\operatorname{Nat}}            %
\newcommand{\nat}{\operatorname{nat}}            %

\newcommand{\ex}{\mathrm{ex}}                    %
\newcommand{\rT}{\mathrm{T}}                     %
\newcommand{\RFib}{\operatorname{RFib}}          %

\newcommand{\add}{\mathrm{add}}                  %
\newcommand{\vloc}{\mathrm{vloc}}                %
\newcommand{\kloc}{\mathrm{kloc}}                %

\newcommand{\Funadd}{\Fun^\add}                  %

\newcommand{\Funvloc}{\Fun^\vloc}                  %
\newcommand{\Funkloc}{\Fun^\kloc}                  %

\newcommand{\Funx}{\operatorname{Fun^{ex}}}      %
\newcommand{\Funq}{\operatorname{Fun^q}}         %
\newcommand{\Funs}{\operatorname{Fun^s}}         %
\newcommand{\Funh}{\operatorname{Fun^h}}         %

\newcommand{\Null}{\operatorname{Null}}          %
\newcommand{\Pairings}{\operatorname{Pair}}      %
\newcommand{\pair}{\mathrm{pair}}                %

\newcommand{\metstab}{\mathrm{stab}}             %
\newcommand{\ads}{\mathrm{ad}}                   %

\newcommand{\Motp}[1][]{{\if\relax\detokenize{#1}\mathrm{Mot^p}\relax\else\mathrm{Mot}_{#1}^{\mathrm{p}}\fi}}  %
\newcommand{\Motpun}[1][]{\mathrm{Mot^p_{un\if\relax\detokenize{#1}\relax\else{,}#1\fi}}}           %
\newcommand{\Motkar}[1][]{{\if\relax\detokenize{#1}\mathrm{Mot^{k}}\relax\else\mathrm{Mot^k_{#1}}\fi}}  %
\newcommand{\Motpbord}[1][]{{\if\relax\detokenize{#1}\mathrm{Mot^{pw}}\relax\else\mathrm{Mot}_{#1}^{\mathrm{pw}}\fi}} %
\newcommand{\Mot}[1][]{{\if\relax\detokenize{#1}\mathrm{Mot}\relax\else\mathrm{Mot_{#1}}\fi}}            %

\newcommand{\MotpE}[2]{\mathrm{Mot^p_{{#2}\if\relax\detokenize{#1}\relax\else{,}#1\fi}}}           %

\newcommand{\MotkwE}[2]{\mathrm{Mot^{kw}_{{#2}\if\relax\detokenize{#1}\relax\else{,}#1\fi}}}           %

\newcommand{\Motwun}[1][]{\mathrm{Mot^w_{un\if\relax\detokenize{#1}\relax\else{,}#1\fi}}}           %
\newcommand{\MotwE}[2]{\mathrm{Mot^w_{{#2}\if\relax\detokenize{#1}\relax\else{,}#1\fi}}}           %

\newcommand{\MotpwE}[2]{\mathrm{Mot^{pw}_{{#2}\if\relax\detokenize{#1}\relax\else{,}#1\fi}}}           %

\newcommand{\Motw}[1][]{{\if\relax\detokenize{#1}\mathrm{Mot^w}\relax\else\mathrm{Mot^w_{#1}}\fi}}  %

\newcommand{\Motpw}[1][]{{\if\relax\detokenize{#1}\mathrm{Mot^{pw}}\relax\else\mathrm{Mot^{pw}_{#1}}\fi}}  %

\newcommand{\Motkw}[1][]{{\if\relax\detokenize{#1}\mathrm{Mot^{kw}}\relax\else\mathrm{Mot^{kw}_{#1}}\fi}}  %

\newcommand{\MotkarE}[2]{\mathrm{Mot^k_{{#2}\if\relax\detokenize{#1}\relax\else{,}#1\fi}}}           %

\newcommand{\BO}{\mathrm{BO}}                    %
\newcommand{\SO}{\mathrm{SO}}                    %
\newcommand{\BSO}{\mathrm{BSO}}                  %
\newcommand{\Sp}{\mathrm{Sp}}                    %
\newcommand{\MTO}{\mathrm{MTO}}                  %
\newcommand{\MTSO}{\mathrm{MTSO}}                %
\newcommand{\MSO}{\mathrm{MSO}}                  %

\newcommand{\LA}{\mathrm{LA}}                    %

\newcommand{\Pic}{\operatorname{Pic}}            %

\newcommand{\Kspace}{\operatorname{\spacefont{K}}} %
\newcommand{\kk}{\Kspace}                        %
\newcommand{\K}{\operatorname K}                 %
\newcommand{\KK}{\mathbb K}                      %

\newcommand{\KR}{\operatorname{KR}}              %

\newcommand{\Lspace}{\operatorname{\spacefont{L}}}%
\renewcommand{\L}{\operatorname L}               %
\newcommand{\W}{\operatorname W}                 %
\newcommand{\Witt}{\W}                           %

\newcommand{\GWspace}{\operatorname{\spacefont{GW}}} %
\newcommand{\gw}{\GWspace}                       %
\newcommand{\GW}{\operatorname{GW}}              %
\newcommand{\KGW}{\mathbb{GW}}                   %
\newcommand{\GGW}{\mathbb{GW}}                   %
\newcommand{\GWgroup}{\operatorname{GW_0}}       %

\newcommand{\U}{\operatorname{U}}                %
\newcommand{\V}{\operatorname{V}}                %

\newcommand{\bord}{\mathrm{bord}}                %
\newcommand{\rbord}{{\mathrm{c}\bord}}		 	 %
\newcommand{\fgt}{\mathrm{fgt}}                  %

\newcommand{\CCob}{\mathbb{C}\mathrm{ob}}        %

\newcommand{\Q}{\operatorname{Q}}                %
\newcommand{\Qh}{\Q}                             %
\newcommand{\dualQ}{\mathrm d\Q}                 %
\newcommand{\cyl}{\mathrm{cyl}}                  %
\newcommand{\dMet}{\mathrm{d}\hMet}              %
\newcommand{\bcyl}{\mathrm{bcyl}}                %

\newcommand{\arr}{\mathrm{ar}}                   %
\newcommand{\hMet}{\mathrm{Null}}                %

\newcommand{\rS}{\mathrm{S}}                     %

\newcommand{\fpm}{\lambda}			 %
\newcommand{\dfpm}{\check{\fpm}}			 %
\newcommand{\fpmg}{Q}				 %
\newcommand{\gfpm}{{\g\fpm}}
\newcommand{\gdfpm}{{\g\dfpm}}

\newcommand{\Sym}{\mathrm{Sym}}                  %
\newcommand{\Quad}{\mathrm{Quad}}                %

\newcommand{\B}{\mathrm{B}}                      %

\newcommand{\Bil}{\mathrm{B}}               %
\newcommand{\Lin}{\Lambda}               %

\newcommand{\Dual}{\mathrm{D}}              %

\newcommand{\ev}{\mathrm{ev}}                    %

\newcommand{\qshift}[1]{^{[#1]}}                 %

\newcommand{\Hyp}{\operatorname{Hyp}}            %
\newcommand{\HypG}{{\mathcal{H}\mathrm{yp}}}     %
\newcommand{\gHyp}{\operatorname{gHyp}}          %
\newcommand{\ghyp}{{\mathrm{g}\hyp}}             %
\newcommand{\Met}{\operatorname{Met}}            %
\newcommand{\ilag}{\mathrm{can}}                 %
\newcommand{\met}{\mathrm{met}}                  %

\newcommand{\rN}{\mathrm{N}}                    %

\newcommand{\hyp}{\mathrm{hyp}}                  %

\newcommand{\cl}{\mathrm{cl}}                    %

\newcommand{\sym}{\mathrm{s}}                    %
\newcommand{\s}{\mathrm{s}}                      %
\newcommand{\vis}{\mathrm{v}}                    %
\newcommand{\qdr}{\mathrm{q}}                     %
\newcommand{\g}{\mathrm{g}}                      %
\newcommand{\gs}{\mathrm{gs}}                    %
\newcommand{\gq}{\mathrm{gq}}                    %
\newcommand{\gev}{\mathrm{ge}}                   %

\newcommand{\uni}{\mathrm{u}}                    %

\newcommand{\proj}{\mathrm{proj}}                %

\newcommand{\Even}{\mathrm{Ev}}                  %

\newcommand{\Poinc}{\mathrm{Pn}}                 %

\newcommand{\Poi}{\operatorname{Poi}}            %
\newcommand{\catforms}{\mathrm{He}}              %

\newcommand{\spsforms}{\mathrm{Fm}}              %

\newcommand{\Surg}{\mathrm{Surg}}                %

\newcommand{\Cart}{\operatorname{Cart}}          %
\newcommand{\Cocart}{\operatorname{Cocart}}      %

\newcommand{\Ver}{\mathrm{Ver}}
\newcommand{\PVer}{\mathrm{PVer}}

\newcommand{\spl}{\mathrm{spl}} 		%
\newcommand{\qspl}{\mathrm{qspl}} 		%

\newcommand{\gd}{\mathrm{gd}} 		%
\newcommand{\Catgd}{\Cat^{\gd}}
\newcommand{\Dualb}{\Dual_{\mathrm{b}}}
\newcommand{\Waldgd}{\mathrm{Wald}^{\gd}_\infty}
\newcommand{\fl}{\mathrm{fl}}
\newcommand{\conj}{\mathrm{conj}}

\newcommand{\T}{\mathcal{T}}               %

\newcommand{\I}{\mathcal{I}}               %
\newcommand{\J}{\mathcal{J}}               %

\newcommand{\A}{\mathcal{A}}               %
\newcommand{\cB}{\mathcal{B}}              %
\newcommand{\E}{\mathcal{E}}               %
\newcommand{\Etwo}{\mathcal{E}'}           %

\newcommand{\cO}{\mathcal{O}}              %

\newcommand{\C}{\mathcal C}                %
\newcommand{\Ctwo}{{{\mathcal C}'}}        %
\newcommand{\D}{\mathcal{D}}               %
\newcommand{\Dtwo}{\mathcal{D}'}           %

\newcommand{\F}{\mathcal{F}}               %
\newcommand{\G}{\mathcal{G}}               %

\newcommand{\x}{x}                         %
\newcommand{\y}{y} 			   			   %
\newcommand{\z}{z}                         %

\newcommand{\QF}{\Qoppa}                   %
\newcommand{\QFtwo}{{\QF'}}                %
\newcommand{\QFD}{\Phi}                    %
\newcommand{\QFDtwo}{{\QFD'}}              %
\newcommand{\QFE}{\Psi}                    %
\newcommand{\QFEtwo}{\Psi'}                %

\newcommand{\Qgen}[2]{\QF^{\geq #1}_{#2}}  %

\newcommand{\qone}{q}                      %
\newcommand{\qtwo}{q'}                     %

\newcommand{\Lag}{\mathcal{L}}             %

\makeatletter

\renewcommand{\tocsection}[3]{%
\indentlabel{\@ifnotempty{#2}{\parbox[b]{3ex}{\bfseries\ignorespaces#1 #2}}}\bfseries#3} 

\renewcommand{\tocsubsection}[3]{%
\indentlabel{\@ifnotempty{#2}{\hspace{1.6em}\parbox[b]{5ex}{\ignorespaces#1 #2}}}#3}

\renewcommand{\tocsubsubsection}[3]{%
\indentlabel{\@ifnotempty{#2}{\hspace{3.9em}\parbox[b]{5ex}{\ignorespaces#1 #2}}}#3}

\makeatother

\DeclareRobustCommand{\SkipTocEntry}[5]{} 

\newcommand{\introsubsection}[1]{\addtocontents{toc}{\SkipTocEntry}\subsection*{#1}}
\newcommand{\bibsubsubsection}[1]{\addtocontents{toc}{\SkipTocEntry}\subsubsection*{#1}}

\BibSpec{collection.article}{%
    +{}  {\PrintAuthors}                {author}
    +{,} { \textit}                     {title}
    +{.} { }                            {part}
    +{:} { \textit}                     {subtitle}
    +{,} { \PrintContributions}         {contribution}
    +{,} { \PrintConference}            {conference}
    +{}  {\PrintBook}                   {book}
    +{,} { }                            {booktitle}
    +{,} { }                            {series}
    +{,} { \voltext}                    {volume}
    +{,} { }                            {publisher}
    +{,} { }                            {organization}
    +{,} { }                            {address}
    +{,} { \PrintDateB}                 {date}
    +{,} { \eprintpages}                {pages}
    +{,} { }                            {status}
    +{,} { available at \eprint}        {eprint}
    +{}  { \parenthesize}               {language}
    +{}  { \PrintTranslation}           {translation}
    +{;} { \PrintReprint}               {reprint}
    +{.} { }                            {note}
    +{.} {}                             {transition}
    +{}  {\SentenceSpace \PrintReviews} {review}
}

\BibSpec{article}{%
    +{}  {\PrintAuthors}                {author}
    +{.} { \textit}                     {title}
    +{.} { }                            {part}
    +{:} { \textit}                     {subtitle}
    +{,} { \PrintContributions}         {contribution}
    +{.} { \PrintPartials}              {partial}
    +{,} { }                            {journal}
    +{}  { \textbf}                     {volume}
    +{} { \PrintDate}                  {date}
    +{,} { \issuetext}                  {number}
    +{,} { \eprintpages}                {pages}
    +{,} { }                            {status}
    +{,} { available at \eprint}        {eprint}
    +{}  { \parenthesize}               {language}
    +{}  { \PrintTranslation}           {translation}
    +{;} { \PrintReprint}               {reprint}
    +{.} { }                            {note}
    +{.} {}                             {transition}
    +{}  {\SentenceSpace \PrintReviews} {review}
}

\newcommand{\reftwo}[1]{\ref{#1}}
\newcommand{\eqreftwo}[1]{\eqref{#1}}
\newcommand{\reftwoitem}[1]{\ref{#1}}

\newcommand{\refone}[1]{\cite{Part-one}.\ref{I-#1}}

\newcommand{\refoneitem}[1]{\ref{I-#1}}

\newcommand{\refthree}[1]{\cite{Part-three}.\ref{III-#1}}

\newcommand{\reffour}[1]{\cite{Part-four}.\ref{IV-#1}}

\newcommand{\paperone}{Paper~\cite{Part-one}\xspace}

\newcommand{\paperthree}{Paper~\cite{Part-three}\xspace}
\newcommand{\paperfour}{Paper~\cite{Part-four}\xspace}

\newcommand{\refundefined}[1]{\textcolor{red}{Undefined ref}}

\title[Hermitian K-theory for stable $\infty$-categories II: Cobordism categories and additivity]{Hermitian K-theory for stable $\infty$-categories II:\\
Cobordism categories and additivity}

\author[Calmès]{Baptiste Calmès}
\address{Université d'Artois, Laboratoire de Mathématiques de Lens, Lens, France}
\email{baptiste.calmes@univ-artois.fr}

\author[Dotto]{Emanuele Dotto}
\address{University of Warwick, Mathematics Institute, Coventry, United Kingdom}
\email{emanuele.dotto@warwick.ac.uk}

\author[Harpaz]{Yonatan Harpaz}
\address{Université Paris Cité, Sorbonne Université, Paris, France}
\email{harpaz@imj-prg.fr}

\author[Hebestreit]{Fabian Hebestreit}
\address{Universität Bielefeld, Fakultät für Mathematik, Bielefeld, Germany}
\email{hebestreit@math.uni-bielefeld.de}

\author[Land]{Markus Land}
\address{LMU München, Mathematisches Institut, Germany}
\email{markus.land@math.lmu.de}

\author[Moi]{Kristian Moi}
\address{KTH, Institutionen för matematik, Stockholm, Sweden}
\email{kristian.moi@gmail.com}

\author[Nardin]{Denis Nardin}
\address{Universität Regensburg, Mathematisches Institut, Regensburg, Germany}

\author[Nikolaus]{Thomas Nikolaus}
\address{WWU Münster, Mathematisches Institut, Münster, Germany}
\email{nikolaus@uni-muenster.de}

\author[Steimle]{Wolfgang Steimle}
\address{Universität Augsburg, Institut für Mathematik, Augsburg, Germany}
\email{wolfgang.steimle@math.uni-augsburg.de}

\dedicatory{To Andrew Ranicki.}

\date{\today}

\begin{document}

\begin{abstract}
We define Grothendieck-Witt spectra in the setting of Poincaré $\infty$-categories and show that they fit into an extension with a $\K$- and an $\L$-theoretic part. As consequences, we deduce localisation sequences for Verdier quotients and generalisations of Karoubi's fundamental and periodicity theorems for rings in which $2$ need not be invertible. Our set-up allows for the uniform treatment of such algebraic examples alongside homotopy-theoretic generalisations: For example, the periodicity theorem holds for complex oriented $\Eone$-rings, and we show that the Grothendieck-Witt theory of parametrised spectra recovers Weiss and Williams' $\LA$-theory.

Our Grothendieck-Witt spectra are defined via a version of the hermitian $\Q$-construction, and a novel feature of our approach is to interpret the latter as a cobordism category. This perspective also allows us to give a hermitian version -- along with a concise proof -- of the theorem of Blumberg, Gepner and Tabuada, and provides a cobordism theoretic description of the aforementioned $\LA$-spectra. 

\end{abstract}

\maketitle
\tableofcontents

\section*{Introduction}
\introsubsection{Overview}

Unimodular quadratic and symmetric forms are ubiquitous objects in mathematics appearing in contexts ranging from norm constructions in number theory to surgery obstructions in geometric topology. Their classification, however, even over simple rings such as the integers, remains out of reach. A simplification, following ideas of Grothendieck for the study of projective modules, suggests to consider for a commutative ring \(R\) (for ease of exposition) the abelian groups \(\GW^{\qdr}_0(R)\) and $\GW^\sym_0(R)$ given as the group completion of the monoid of isomorphism classes of finitely generated projective \(R\)-modules \(P\), equipped with a unimodular quadratic or symmetric form \(q\), with addition the orthogonal sum
\[
[P,q] + [P',q'] = [P \oplus P',q \perp q'].
\]
This group, commonly known as the (quadratic or symmetric) Grothendieck-Witt group of $R$, was given a homotopy-theoretical refinement at the hands of Karoubi and Villamayor in \cite{karoubi-villamayor}, by adapting Quillen's approach to higher algebraic \(\K\)-theory. In modern terms, one organises the collection of unimodular quadratic or symmetric forms \((P,q)\) into a groupoid $\Unimod^{r}(R)$ for $r \in \{\qdr,\sym\}$, which may be viewed as an
\(\Einf\)-monoid
using the symmetric monoidal structure on $\Unimod^{r}(R)$ arising from orthogonal sum. One can then take the group completion to obtain an \(\Einf\)-group 
\[
\GWspace_\cl^{r}(R) = \Unimod^{r}(R)^\grp,
\]
the classical Grothendieck-Witt space, whose group of components is the Grothendieck-Witt group described above. By definition the higher Grothendieck-Witt groups of $R$ are then the higher homotopy groups of $\GWspace^{r}_\cl(R)$.

There are variants of this construction also for forms over non-commutative rings $R$ taking values in what we call an invertible module $M$ with involution over $R$, which in particular covers forms over rings with involution. In this generality, one can also simply change the sign of the involution on $M$ to pass from symmetric and quadratic to skew-symmetric and skew-quadratic forms; we will denote the resulting module with involution by $-M$. Finally, we shall also consider even forms in the sequel, that is forms admitting a quadratic refinement but not equipped with one. For $R$ commutative they for example encode symplectic, rather than just skew-symmetric, forms via $M=-R$ also in characteristic $2$. The various kinds of forms are then connected by polarisation maps
\[
\GWspace_\cl^\qdr(R,M) \lrar \GWspace_\cl^\ev(R,M) \lrar \GWspace_\cl^\sym(R,M)
\]
which are equivalences if $2$ is a unit in $R$, but not in general. \\

In the present paper we establish a decomposition of any of the above Grothendieck-Witt spaces into a $\K$-theoretic and an $\L$-theoretic part, and establish a tighter relation between them than was previously expected. To motivate the decomposition into $\K$- and $\L$-theoretic pieces, we recall for $r \in \{\qdr,\ev,\sym\}$ the following classical exact sequence 
\[
\K_0(R,M)_{\Ct} \xrightarrow{\hyp} \GW^r_0(R,M) \longrightarrow \Witt^r_0(R,M) \longrightarrow 0;
\]
here, the $\Ct$-coinvariants on the group $\K_0(R,M) = \K_0(R)$ on the left are formed with respect to the dualisation action $P \mapsto P^* \coloneq \Hom_R(P,M)$, and $\Witt^r_0(R,M)$ is the Witt group, defined as the monoid of isomorphism classes of unimodular $M$-valued $r$-forms modulo those that admit a Lagrangian. The middle map is then simply the projection and the left hand one, the hyperbolisation, takes $P$ to $P \oplus P^*$ equipped with its evaluation form. 

The principal goal of the present paper is to extend this to a long exact sequence with $\L$-groups playing the role of higher Witt groups. Such results are well-known principally from the work of Karoubi and Schlichting if $2$ is a unit in $R$, and have lead to a good understanding of Grothendieck-Witt theory relative to $\K$-theory since the relevant $\L$-groups are often rather accessible. For instance, by work of Ranicki \cite{Ranickiblue} when $2$ is invertible in $R$, one has $\Witt_0(R,M) = \L_0(R,M)$ and in addition $\L_{i+2}(R,M) \cong \L_i(R,-M)$ for all $i \in \mathbb Z$, so that the $\L$-groups are in particular $4$-periodic. There is furthermore a large body of work computing Witt-groups, culminating for example in Voevodsky’s solution to Milnor's conjecture that Witt groups of fields not of characteristic $2$ admit complete filtrations with filtration quotients given by Galois cohomology \cite{VoevodskyMilnor,Voe-Witt}. An older result of Kato achieves a similar description in characteristic $2$, see \cite{Kato-Witt}. The resulting  calculations for Grothendieck-Witt groups of fields can then often be leveraged for calculations over other rings by means of various localisation sequences  if $2$ is invertible in $R$~\cite{hornbostel-schlichting, schlichting-mv}. Another overarching goal of our paper series is to provide such localisation sequences for Grothendieck-Witt groups if $2$ is not assumed invertible in $R$ and to describe the extent to which periodicity statements for $\L$-groups still hold. The latter turns out to be surprisingly closely related to the behaviour of the polarisation maps connecting the various Grothendieck-Witt spaces.\\

\introsubsection{Historical background}

The study of Grothendieck-Witt spaces begins by comparing them to Quillen's algebraic $\K$-space $\Kspace(R)$ defined as the group completion of the groupoid of finitely generated projective modules over $R$. To this end, one has maps
\[
\fgt\colon \GWspace_\cl^\sym(R,M) \lrar \Kspace(R) \quad \text{and} \quad \hyp \colon \Kspace(R) \lrar \GWspace_\cl^\qdr(R,M),
\]
the former extracting the underlying module of a unimodular form, the latter induced by the hyperbolisation construction explained above. 

In his fundamental papers \cite{Karoubi-Le-theoreme-fondamental, karoubi-quillen} Karoubi analysed the case in which $2$ is a unit in $R$ (so no distinction between the three flavours of Grothendieck-Witt groups is necessary). He considered the spaces
\[
\mathcal U_\cl(R,M) = \fib(\Kspace(R) \xrightarrow{\hyp} \GWspace_\cl(R,M)) \quad \text{and} \quad \mathcal V_\cl(R,M) = \fib(\GWspace_\cl(R,M) \xrightarrow{\fgt} \Kspace(R)),
\]
and the cokernels $\mathrm W_i(R,M)$ of the maps $\K_i(R) \xrightarrow{\hyp} \GW_i(R,M)$. He then showed the following two results nowadays known as his fundamental and periodicity theorems, namely that
\[
\Omega \mathcal U_\cl(R,-M) \simeq \mathcal V_\cl(R,M), \quad \text{and} \quad \mathrm W_i(R,M)[1/2] \cong \mathrm W_{i+2}(R,-M)[1/2].
\]
In fact, Karoubi showed the right hand identification without the assumption that $2$ be invertible in $R$. These results permit one to inductively deduce results on higher Grothendieck-Witt groups from information about algebraic \(\K\)-theory on the one hand and about \(\mathrm W_{i}(R,\pm M)\) for \(i=0,1\) on the other. 

To control the behaviour of the $2$-torsion in the cokernel of the hyperbolisation map Kobal in \cite{Kobal} introduced refinements of the hyperbolic and forgetful maps: By the invertibility assumption on $M$, the functor taking $M$-valued duals induces an action of the group $\Ct$ on the algebraic $\K$-spectrum and we denote the arising $\Ct$-spectrum by $\K(R,M)$ and similarly for the $\K$-space $\Kspace(R,M)$. The maps above then refine to maps
\[
\Kspace(R,M)_\hC \xrightarrow{\hyp} \GWspace_\cl^\qdr(R,M) \longrightarrow \GWspace_\cl^\sym(R,M) \xrightarrow{\fgt} \Kspace(R,M)^\hC,
\]
where the orbits on the left are formed in $\Einf$-groups (and not on underlying spaces), and the composite is the norm of $\Kspace(R,M)$. Kobal used these refinements to show that, if $2$ is invertible in $R$, the cofibre of $\hyp \colon \Kspace(R,M)_\hC \rightarrow \GWspace_\cl(R,M)$ is $4$-periodic on the nose. 

The next major steps forward were then taken by Schlichting in \cite{schlichting-derived}, who introduced Grothendieck-Witt spectra for differential graded categories with duality in which $2$ is invertible. Applied to categories of chain complexes over $R$ he obtained (usually non-connective) spectra $\GW_\cl(R,M)$ with $\Omega^\infty \GW_\cl(R,M) \simeq \GWspace_\cl(R,M)$ and used this machinery to establish a fibre sequence
\[
\GW_\cl(R,M\qshift{-1}) \xrightarrow{\fgt} \K(R,M) \xrightarrow{\hyp} \GW_\cl(R,M),
\]
which he termed the Bott sequence, where $(-)\qshift{i}$ refers to the $i$-fold shift of the duality afforded by his set-up. He also proved that $\GW_\cl(R,M\qshift{i}) \simeq \GW_\cl(R,-M\qshift{i+2})$, yielding in total a new conceptual proof of Karoubi's fundamental theorem. Still assuming $2$ invertible in $R$ he furthermore showed that the ($4$-periodic) homotopy groups of the cofibre of the refined hyperbolic map $\hyp \colon \K(R,M)_\hC \rightarrow \GW_\cl(R,M)$ are indeed given by the Witt groups $\Witt_0(R,\pm M)$ in even degrees and by Witt groups of formations in odd degrees.

This led to the folk theorem that if $2$ is a unit in $R$ the cofibre of $\hyp \colon  \K(R,M)_\hC \rightarrow \GW_\cl(R,M)$ is given by Ranicki's $\L$-spectum $\L(R,M)$ from \cite{Ranickiblue}, whose homotopy groups are well-known to match Schlichting's results, though as far as we are aware no account at the level of spectra has appeared in the literature. 
Without the assumption on the invertibility of $2$, however, many of the methods employed break down. In particular, the relation between Grothendieck-Witt-groups and $\L$-groups remained mysterious: For example, it is well-known that the cofibre of the hyperbolisation map in the quadratic case cannot directly relate to quadratic $\L$-spectra. If $2$ is not invertible in $R$, there are indeed many flavours of $\L$-groups and as far as we are aware not even a precise conjecture had been put forward. In contrast to this situation, Karoubi predicted in \cite{karoubi-periodicity} a precise form in which his fundamental theorems should extend to general rings: It is not only the sign on $M$ that changes when passing from $\U$ to $\V$ but also the form parameter. A similar suggestion was made by Giffen, see \cite{williams-quadratic}. In what is hopefully evident notation they conjectured
\[
\Omega \mathcal U_\cl^\qdr(R,-M) \simeq \mathcal V_\cl^\ev(R,M) \quad \text{and} \quad \Omega \mathcal U_\cl^\ev(R,-M) \simeq \mathcal V_\cl^\sym(R,M).
\]

In the present paper, we obtain extensions of Karoubi's periodicity and fundamental theorem to arbitrary rings, affirming in particular the conjecture of Karoubi and Giffen, and also determine the cofibre of the hyperbolisation map as an $\L$-spectrum (which, however, need not be $4$-periodic in general). \\

Karoubi's fundamental and periodicity theorem have also been instrumental in concrete computations, for instance they have lead to an 
almost complete computation of the Grothendieck-Witt groups of $\ZZ[\onehalf]$ in \cite{berrick-karoubi} and to great structural insight by controlling the $2$-adic behaviour of the forgetful map $\GW_\cl \rightarrow \K^\hC$ in \cite{BKSO-fixed-point} (still under the assumption that $2$ is a unit). In \paperthree, we will employ the results of this paper to prove extensions of these computational and structural results to general rings of integers in number fields, and in particular to $\ZZ$.

\introsubsection{Main results}

Our approach is based on placing Grothendieck-Witt- and $\L$-theory into a
common general framework, namely that of Poincar\'e
$\infty$-categories, which were introduced by Lurie in his lectures on $\L$-theory  \cite{Lurie-L-theory},
and developed in detail in \paperone. A Poincaré
$\infty$-category is a small stable $\infty$-category $\C$ together with a
certain functor $\QF\colon \C\op\to \Spa$ which encodes
the type of forms (such as, quadratic, even or symmetric) under consideration. 
The requirements on $\QF$ are such that there is an associated duality equivalence $\Dual_\QF \colon \C\op \rightarrow \C$ and as in Schlichting's setup, Poincar\'e structures can be shifted via $\QF\qshift{i} = \SS^i \otimes \QF$ with associated duality $\SS^i \otimes \Dual_\QF$.

As mentioned above, Lurie defined $\L$-spectra $\L(\C,\QF)$ for Poincaré $\infty$-categories, and it is by now standard to view $\K$-spectra as a functor on stable $\infty$-categories. The duality $\Dual_\QF$ induces a $\Ct$-action on the $\K$-spectrum of a Poincaré $\infty$-category and we will denote the resulting $\Ct$-spectrum by $\K(\C,\QF)$. In the present paper we define a Grothendieck-Witt spectrum $\GW(\C,\QF)$ for every Poincar\'e $\infty$-category and as the main result we provide extensions of Karoubi's periodicity theorem and Schlichting's Bott sequence for arbitrary Poincar\'e $\infty$-categories:

\begin{named}{Main Theorem}
\label{theorem:main-intro-two}%
For every Poincaré $\infty$-category $(\C,\QF)$, there is a fibre sequence
\[
\K(\C,\QF)_\hC \xrightarrow{\hyp} \GW(\C,\QF) \xrightarrow{\bord} \L(\C,\QF),
\]
which canonically splits after inverting $2$, and a fibre sequence
\[
\GW(\C,\QF\qshift{-1}) \xrightarrow{\fgt} \K(\C) \xrightarrow{\hyp} \GW(\C,\QF).
\]
\end{named}

To explain how this generalises the results of Karoubi and Schlichting on Grothendieck-Witt spectra of discrete rings mentioned above, take $\C = \Dperf(R)$, the stable subcategory of the derived $\infty$-category $\D(R)$ spanned by the perfect complexes over $R$. As part of \paperone, we constructed Poincaré structures
\[
\QF^\qdr_M \Longrightarrow \QF^{\gq}_M {\Longrightarrow} \QF^{\gev}_M \Longrightarrow \QF^{\gs}_M \Longrightarrow \QF^\sym_M
\]
connected by maps as indicated: Roughly, the outer two assign to a chain complex its spectrum of homotopy coherent quadratic or symmetric $M$-valued forms, whereas the middle three are the more subtle animations, or in more classical terminology non-abelian derivations, of the functors 
\[
\Quad_M, \Even_M, \Sym_M \colon \Proj(R)\op \rightarrow \Ab
\]
parametrising ordinary $M$-valued quadratic, even and symmetric forms, respectively. We shall refer to the outer two as the homotopical and inner ones as the genuine Poincar\'e structures on $\Dperf(R)$. The comparison maps between all five Poincar\'e structures are equivalences if $2$ is a unit in $R$, but in general they give rise to five distinct Grothendieck-Witt spectra. Moreover,
from the main result of \cite{comparison}, we find that it is the middle three, which relate to the classical Grothendieck-Witt spaces, i.e.\ we have
\[
\Omega^\infty\GW(\Dperf(R),\QF^{\mathrm{g}r}_M) \simeq \GWspace_\cl^r(R,M)
\]
for $r \in \{\qdr,\ev,\sym\}$, whereas it is the outer two Poincar\'e structures $\QF_M^\qdr$ and $\QF_M^\sym$ that give rise to Ranicki's $4$-periodic quadratic and symmetric $\L$-spectra.
This mismatch (which is also the reason for carrying the subscript $\cl$ through the introduction) explains much of the subtlety that arose in previous attempts to connect Grothendieck-Witt- and $\L$-theory. 

For consistency of notation we shall allow ourselves to write $\GW_\cl^r(R,M)$ for $\GW(\Dperf(R),\QF^{\mathrm{g}r}_M)$ also if $2$ is not necessarily invertible in $R$; if it is, these spectra really do agree with Schlichting's constructions, as we verify in the second appendix. The identifications at the heart of his proofs of Karoubi's fundamental theorem then extend to identifications
\[
(\Dperf(R),(\QF^\qdr_M)\qshift{2}) \simeq (\Dperf(R),\QF^\qdr_{-M}) \quad \text{and} \quad (\Dperf(R),(\QF^\sym_M)\qshift{2}) \simeq (\Dperf(R),\QF^\sym_{-M}),
\]
so if $2$ is a unit in $R$, we, in particular, recover his results and extend the identification of the cofibre of the hyperbolisation map to the spectrum level. More importantly, however, even if $2$ is not necessarily invertible, we also find
\[
(\Dperf(R),(\QF^{\gs}_M)\qshift{2}) \simeq (\Dperf(R),\QF^{\gev}_{-M}) \quad \text{and} \quad (\Dperf(R),(\QF^{\gev}_M)\qshift{2}) \simeq (\Dperf(R),\QF^{\gq}_{-M}),
\]
whence the second part of the \reftwo{theorem:main-intro-two} settles the conjecture of Giffen and Karoubi. %
Explicitly, in hopefully evident notation, we obtain:

\begin{ucorollary}
For a discrete ring $R$ and an invertible $R$-module $M$ with involution there are canonical equivalences
\[
\U_\cl^\qdr(R,-M) \simeq \SS^1 \otimes \V_\cl^\ev(R,M) \quad \text{and} \quad \U_\cl^\ev(R,-M) \simeq \SS^1 \otimes \V_\cl^\sym(R,M).
\]
\end{ucorollary}

As a consequence of the first part of our \reftwo{theorem:main-intro-two} we also obtain a direct relation between the Grothendieck-Witt spectra for different form parameters. As an implementation of Ranicki's $\L$-theoretic periodicity results, Lurie produced canonical equivalences
$
\L(\C,\QF\qshift{1}) \simeq \SS^1 \otimes \L(\C,\QF).
$
Applying this four times we obtain a stabilisation map
\[
\mathrm{stab} \colon \SS^4 \otimes \L(\Dperf(R),\QF_M^\gs) \simeq \L(\Dperf(R),\QF_M^\gq)  \xrightarrow{\mathrm{pol}} \L(\Dperf(R),\QF_M^\gs)
\]
which is an equivalence if $2$ is invertible in $R$. As another articulation of periodicity we have:

\begin{ucorollary}
The natural map $\mathrm{pol} \colon \GW_\cl^\qdr(R,M) \rightarrow \GW_\cl^\sym(R,M)$ fits into a commutative diagram
\[
\begin{tikzcd}
\K(R,M)_\hC \ar[d,"\id"] \ar[r,"\hyp"] & \GW_\cl^\qdr(R,M) \ar[r] \ar[d,"\mathrm{pol}"] & \SS^4  \otimes \L(\Dperf(R),\QF^\gs_M) \ar[d,"\mathrm{stab}"] \\
\K(R,M)_\hC \ar[r,"\hyp"] & \GW_\cl^\sym(R,M) \ar[r,"\bord"] & \L(\Dperf(R),\QF^\gs_M)
\end{tikzcd}
\]
of fibre sequences.
\end{ucorollary}

In particular, this result ties the behaviour of the polarisation map $\GW_\cl^\qdr(R,M) \rightarrow \GW_\cl^\sym(R,M)$ directly to the periodicity properties $\L(\Dperf(R),\QF^\gs_M)$. In the body of the text, we shall derive more general versions of both corollaries concerning Grothendieck-Witt groups associated to suitable pairs of Bak's form parameters. In particular, we settle Conjectures~1 and~2 of \cite{karoubi-periodicity} in full.%

In the third part of this series \paperthree, we investigate the spectra $\L(\Dperf(R),\QF^{\g r}_M)$ and $\GW(\Dperf(R),\QF^{\g r}_M)$ for general form parameters in detail. In particular, we refine work of Ranicki to show that for Dedekind rings, and thus in particular number rings, the comparison map $\L(\Dperf(R),\QF^\gs_M) \rightarrow \L(\Dperf(R),\QF^\s_M)$ is an equivalence in non-negative degrees. In particular, in this case the non-negative genuine symmetric $\L$-groups are still $4$-periodic, and surprisingly it is symmetric and not quadratic $\L$-groups which govern the cofibre of the hyperbolisation map $\K(R,M)_\hC \rightarrow \GW^\qdr(R,M)$ outside small degree. We use this analysis in conjuction with our \reftwo{theorem:main-intro-two} to obtain dévissage results for Grothendieck-Witt spectra of Dedekind rings, establish finite generation results and confirm Thomason's homotopy limit problem for number rings, and ultimately provide an almost complete computation of the Grothendieck-Witt groups of the integers. \\

\introsubsection{Hermitian $\K$-theory of Poincaré $\infty$-categories}

Let us now sketch the road to our main result. Besides the set-up of Poincaré $\infty$-categories the main novelty of our approach is its direct connection to the theory of cobordism categories of manifolds. To facilitate the discussion, recall that $\Cob_d$ has as objects $(d-1)$-dimensional closed, oriented, smooth manifolds, and spaces of cobordisms thereof as morphisms. The celebrated equivalence 
\[
|\Cob_d| \simeq \Omega^{\infty-1}\MTSO(d),
\]
established by Galatius, Madsen, Tillmann and Weiss in \cite{GMTW} then lies at the heart of much modern work on the homotopy types of diffeomorphism groups \cite{Stablemoduli}; here $\MTSO(d)$ denotes the Thom spectrum of $-\gamma_d \rightarrow \BSO(d)$, where $\gamma_d$ denotes the universal vector bundle over $\mathrm{BSO}(d)$.

Now, a Poincaré $\infty$-category $(\C,\QF)$ determines a space of Poincaré objects $\Poinc(\C,\QF)$ to be thought of as the higher categorical generalisation of the groupoid $\Unimod^r(R,M)$ of unimodular forms considered in the case of discrete rings above. By adapting Quillen's $\Q$-construction to the hermitian setup, we produce for every Poincaré $\infty$-category $(\C,\QF)$ an analogous cobordism category $\Cob(\C,\QF) \in \Cat$ with objects given by $\Poinc(\C,\QF\qshift{1})$ and morphisms given by spaces of Poincaré cobordisms, Ranicki style.

As the technical heart of the present paper, we show the following additivity theorem:

\begin{theoremintrotwo}
\label{theorem:additivity-intro-two}%
If $p \colon (\D,\QFD) \to (\E,\QFE)$ is a split Poincaré-Verdier projection with kernel $(\C,\QF)$, then it induces a bicartesian fibration of $\infty$-categories $\Cob(\D,\QFD) \to \Cob(\E,\QFE),$ whose fibre over $0 \in \Cob(\E,\QFE)$ is $\Cob(\C,\QF)$. In particular, one obtains a fibre sequence
\[
|\Cob(\C,\QF)| \longrightarrow |\Cob(\D,\QFD)| \longrightarrow |\Cob(\E,\QFE)|
\]
 of spaces.
\end{theoremintrotwo}

Here, $p$ is called a Poincaré-Verdier projection, if the sequence $(\C,\QF) \rightarrow (\D,\QFD) \rightarrow (\E,\QFE)$ is both a fibre and a cofibre sequence in $\Catp$, the $\infty$-category of Poincaré $\infty$-categories. We call such sequences Poincar\'e-Verdier sequences and they are split if the underlying functors both admit both adjoints. This is the Poincar\'e analogue of the notion of (split) Verdier sequences of stable $\infty$-categories and indeed the result is a hermitian analogue of Waldhausen's additivity theorem, as we shall momentarily explain. To any Poincar\'e category $(\C,\QF)$ we can associate a new Poincar\'e category $\Met(\C,\QF)$ whose underlying category is $\Ar(\C) = \Fun(\Delta^1,\C)$ and whose Poincar\'e objects encode Poincar\'e objects in $(\C,\QF)$ together with a Lagrangian. As such Lagrangians are the same thing as Ranicki's null-cobordisms, the cobordism category $\Cob^\partial(\C,\QF) \coloneq \Cob(\Met(\C,\QF)\qshift{1})$ plays the role of a cobordism category with boundary akin to the geometric situation. Extracting boundaries and viewing objects as having trivial boundary then induces a split Poincar\'e-Verdier sequence, the \emph{metabolic sequence}
\[
(\C,\QF\qshift{-1}) \lto \Met(\C,\QF) \stackrel{\met}{\lto} (\C,\QF),
\]
whose right hand map is underlain by the target projection $\mathrm{t} \colon \Ar(\C) \rightarrow \C$. Applied to this special case, Theorem~\reftwo{theorem:additivity-intro-two} gives a fibre sequence 
\[
|\Cob(\C,\QF\qshift{-1})| \longrightarrow |\Cob^\partial(\C,\QF\qshift{-1})| \stackrel{\partial}{\longrightarrow} |\Cob(\C,\QF)|
\]
which on the one hand, we view as an algebraic version of Genauer's fibre sequence
\[
|\Cob_{d+1}| \longrightarrow |\Cob^\partial_{d+1}| \stackrel{\partial}{\longrightarrow} |\Cob_{d}|
\]
from geometric topology \cite{Genauer}. On the other hand, it is a hermitian analogue of the (split) fibre sequence 
\[
|\Span(\C)| \longrightarrow |\Span(\Ar(\C))| \stackrel{\mathrm{t}}{\longrightarrow} |\Span(\C)|,
\]
which in view of the equivalence $\Omega|\Span(\C)| \simeq \mathcal K(\C)$ is equivalent to Waldhausen's classical additivity theorem $\Kspace(\Ar(\C)) \simeq \Kspace(\C)^2$. In contrast to this case, however, neither the geometric nor our hermitian fibre sequences split. 

Our proof of Theorem~\reftwo{theorem:additivity-intro-two} is in fact modelled on the recent proof of Genauer's fibre sequence at the hands of the ninth author \cite{Steimleadd} and is new even in the context of algebraic $\K$-theory. Similar additivity results are known in other set-ups and in varying degrees of generality, e.g.\ in \cite{schlichting-derived, Spitzweckreal}. The actual statement we prove in the body of the paper is, however, more general and also applies to an arbitrary additive functor $\F \colon \Catp \rightarrow \Sps$ replacing $\F=\Poinc$ in the definition of cobordism categories; here, we call such $\F$ additive if it takes split Poincar\'e-Verdier squares, that is, cartesian squares of Poincar\'e $\infty$-categories whose vertical maps are split Poincar\'e-Verdier projections, to cartesian squares of spaces. %

We then proceed on to define the Grothendieck-Witt spectrum $\GW(\C, \QF)$ by iterating the hermitian $\Q$-construction and thus also the Grothendieck-Witt space $\GWspace(\C,\QF) = \Omega^\infty \GW(\C,\QF)$. The refined version of Theorem~\reftwo{theorem:additivity-intro-two}, which in particular applies to such iterated constructions, gives us a basis for establishing:

\begin{theoremintrotwo}
\label{theorem:gwuniversal-intro-two}%
\begin{enumerate}
\item There is a equivalences 
\[
|\Cob(\C,\QF)| \simeq \Omega^{\infty-1}\GW(\C,\QF),
\]
natural in $(\C,\QF) \in \Catp$, so in particular, $\Omega|\Cob(\C,\QF)| \simeq \GWspace(\C,\QF)$.
\item 
\label{item:universal-properties-GW}%
The functors $\GW \colon \Catp \rightarrow \Spa$ and $\GWspace \colon \Catp \rightarrow \Grp_\Einf(\Sps)$ are the initial additive functors equipped with a transformation $\Poinc \rightarrow \GWspace \simeq \Omega^\infty\GW$.
\item The functor $\L \colon \Catp \rightarrow \Spa$ is the initial additive, bordism invariant functor equipped with a transformation $\Poinc \rightarrow \Omega^\infty\L$.  
\end{enumerate}
\end{theoremintrotwo}

Here, an additive functor $\Catp \to \Spa$ is called bordism invariant if it vanishes on $\Met(\C,\QF)$ for every Poincar\'e $\infty$-category $(\C,\QF)$, a notion that on account of Waldhausen additivity has no sensible counterpart in the context of ordinary algebraic $\K$-theory. Theorem~\reftwo{theorem:gwuniversal-intro-two} simultaneously gives an algebraic analogue of the theorem of Galatius, Madsen, Tillmann and Weiss above concerning the homotopy type of the cobordism category and the hermitian analogue of the theorem of Blumberg, Gepner and Tabuada from \cite{BGT} that $\K\colon \Catx \rightarrow \Spa$ is the initial additive functor with a transformation $\core \rightarrow \Omega^\infty \K = \Kspace$. \\

To establish the fibre sequence $\K(\C,\QF)_\hC \rightarrow \GW(\C,\QF) \rightarrow \L(\C,\QF)$ starting from Theorem \reftwo{theorem:gwuniversal-intro-two} we finally also exhibit the cofibre of the left hand map as the initial bordism invariant functor under $\GW$, so that Theorem \reftwo{theorem:gwuniversal-intro-two} applies to identify it with $\L(\C,\QF)$. To obtain the extension $\GW(\C,\QF\qshift{-1}) \rightarrow \K(\C) \rightarrow \GW(\C,\QF)$ of Schlichting's fibre sequence we use a version of Ranicki's algebraic Thom construction to derive a natural equivalence
\[
|\Cob^\partial(\C,\QF)| \simeq |\Span(\C)| \simeq \Omega^{\infty-1}\K(\C),
\]
analogous to Genauer's identification $|\Cob^\partial_d| \simeq \Omega^{\infty-1}\SS[\mathrm{BSO}(d)]$. This identification upgrades to an equivalence $\GW(\Met(\C,\QF)) \simeq \K(\C)$ and the desired fibre sequence is thus obtained by taking Grothendieck-Witt spectra of the metabolic sequence. In particular, it is an algebraic analogue of Genauer's results, and we therefore term it the Bott-Genauer sequence.

With the \reftwo{theorem:main-intro-two} established, we observe that since both $\L$- and $\K$-theory take arbitrary Poincar\'e-Verdier sequences in $\Catp$ to fibre sequences, we obtain the following localisation theorem for Grothendieck-Witt theory:

\begin{corollaryintrotwo}
\label{corollary:gwadditive-intro-two}%
The functor $\GW \colon \Catp \rightarrow \Spa$ is Verdier localising, i.e.\ it takes arbitrary Poincaré-Verdier sequences 
$
(\C,\QF) \to (\D,\QFD) \to (\E,\QFE)
$
to fibre sequences
\[
\GW(\C,\QF) \longrightarrow \GW(\D,\QFD) \longrightarrow \GW(\E,\QFE)
\]
of spectra.
\end{corollaryintrotwo}

An important example of Poincaré-Verdier sequences arises from localisations of rings. Applied in these cases, Corollary~\reftwo{corollary:gwadditive-intro-two} can for example be used to reproduce and generalise the classical localisation and Mayer-Vietoris sequences among (symmetric) Grothendieck-Witt spaces and spectra~\cite{hornbostel-schlichting,schlichting-mv}, see Corollary~\reftwo{corollary:Mayerviet-intro-two} below and the subsequent discussion. 

As the final piece of general structure, we recall that in \cite{hesselholt-madsen-real} Hesselholt and Madsen combined the symmetric Grothendieck-Witt and algebraic $\K$-spectra of a ring $R$ and an invertible module with involution $M$ into what they termed the real algebraic $\K$-spectrum, a genuine $\Ct$-spectrum with underlying spectrum $\K(R)$ and genuine fixed points given by the connective cover of $\GW_\cl^\sym(R,M)$. We likewise construct a functor $\KR \colon \Catp \rightarrow \Spa^\gCt$, where $\Spa^\gCt$ denotes the $\infty$-category of genuine $\Ct$-spectra. The underlying $\Ct$-spectrum is again given by $\K(\C,\QF)$ and in our case the genuine fixed points are $\GW(\C,\QF)$ (not just its connective cover). Our \reftwo{theorem:main-intro-two} can then be expressed as identifying the geometric fixed points of $\KR$ with $\L$-spectra, which combined with the comparison results of \cite{comparison} affirms the conjecture of Hesselholt and Madsen, that the geometric fixed points of the real algebraic $\K$-spectrum of a discrete ring are a version of Ranicki's $\L$-theory. In total one can identify the isotropy separation square of $\KR(\C,\QF)$ %
\[
\begin{tikzcd}
\KR(\C,\QF)^\gCt \ar[r] \ar[d] & \KR(\C,\QF)^\geofix \ar[d] \\
\KR(\C,\QF)^\hC \ar[r] & \KR(\C,\QF)^\tC
\end{tikzcd} \quad \text{as} \quad\begin{tikzcd}
\GW(\C,\QF) \ar[r,"\bord"] \ar[d,"\fgt"] & \L(\C,\QF) \ar[d,"\Xi"] \\
\K(\C,\QF)^\hC \ar[r,"\mathrm{can}"] & \K(\C,\QF)^\tC,
\end{tikzcd}
\]
where $\Xi \colon \L(\C,\QF) \rightarrow \K(\C,\QF)^\tC$ is an extension of a transformation into the Tate construction on algebraic $\K$-theory first constructed by Weiss and Williams in \cite{WWII}.

As the ultimate expression of periodicity, we then have the following statement in the language of genuine equivariant homotopy theory, which combines the periodicity results about $\L$-spectra and Karoubi's fundamental theorem into a single result:

\begin{theoremintrotwo}
\label{theorem:kr-intro-two}%
The boundary map of the metabolic Poincaré-Verdier sequence provides a canonical equivalence
\[
\KR(\C,\QF\qshift{1}) \simeq \SS^{1-\sigma} \otimes \KR(\C,\QF),
\]
where $\sigma$ denotes the sign representation of $\Ct$.
\end{theoremintrotwo}

Indeed, the periodicity properties of $\L$-spectra are now an instance of the identification $(\mathbb S^{1-\sigma} \otimes X)^\geofix \simeq \mathbb S ^1 \otimes X^\geofix$ and the Genauer-Bott sequence is recovered via the fibre sequence $X^{\gCt}\rightarrow X \rightarrow (\mathbb S^{1-\sigma} \otimes X)^{\gCt}$, both valid for general $X \in \Spa^\gCt$.

\introsubsection{Further applications to rings and parametrised spectra}
Let us now specialise the abstract results of the previous section to the classical Grothendieck-Witt spectrum of a (discrete) ring $R$ with coefficients in an invertible module with involution $M$.
We constructed in \paperone a sequence of Poincaré structures 
\[
\QF^\qdr_M = \Qgen{\infty}{M} \Longrightarrow\dots \Longrightarrow  \Qgen{m}{M} \Longrightarrow \Qgen{m-1}{M} \Longrightarrow  \dots \Longrightarrow \Qgen{-\infty}{M} = \QF_M^\sym
\]
on $\Dperf(R)$, all of whose transition maps are equivalences if $2$ is invertible in $R$.
The genuine Poincaré structures introduced earlier appear as $\QF^{\gs} = \Qgen{0}{M}$, $\QF^{\gev}_M = \Qgen{1}{M}$ and $\QF^{\gq}_M = \Qgen{2}{M}$. 
Furthermore, in \paperone we showed that the shift functors on $\Dperf(R)$ induce
equivalences of the form
\[
(\Dperf(R),(\Qgen{m}{M})\qshift{2}) \simeq (\Dperf(R),\Qgen{m+1}{-M}).
\]
Applying Theorem~\reftwo{theorem:kr-intro-two} we find the following result, which (for $2$ invertible in $R$) proves another unpublished conjecture of Hesselholt-Madsen:

\begin{corollaryintrotwo}[Genuine Karoubi periodicity]
\label{corollary:karper}%
For a (discrete) ring $R$, an invertible $R$-module $M$ with involution, 
and $m \in \ZZ \cup \{\pm \infty\}$ there are canonical equivalences
\[
\KR(\Dperf(R),\Qgen{m}{M}) \simeq \SS^{2-2\sigma} \otimes \KR(\Dperf(R),\Qgen{m-1}{-M}).
\]
In particular, the genuine $\Ct$-spectra $\KR(\Dperf(R),\QF^\sym_M)$ and $\KR(\Dperf(R),\QF^\qdr_M)$ are $(4-4\sigma)$-periodic and even $(2-2\sigma)$-periodic if $R$ has characteristic $2$. 
\end{corollaryintrotwo}

In particular, we find
\[
\KR(\Dperf(R),\QF^{\gq}_M) \simeq \SS^{2-2\sigma} \otimes \KR(\Dperf(R),\QF^{\gev}_{-M})\simeq \SS^{4-4\sigma} \otimes \KR(\Dperf(R),\QF^{\gs}_M).
\]
In fact, we show that the $(4-4\sigma)$- or $(2-2\sigma)$-fold periodicity of
$\KR(\Dperf(R),\QF^\sym_M)$ and $\KR(\Dperf(R),\QF^\qdr_M)$
in fact holds for any complex oriented or real oriented $\Eone$-ring $R$, respectively; we will deduce it in this generality in the body of the paper along with periodicity results for other ring spectra such as $\mathrm{ko}$ and $\mathrm{tmf}$ often with periods larger than $4$.

Let us now turn to a sample application regarding the behaviour of Grothendieck-Witt spectra under localisations of rings. As part of our general analysis we show that
\[
\begin{tikzcd} 
\biggl(\Dperf(R),\Qgen{m}{M}\biggr) \ar[r] \ar[d] & \biggl(\Der^{\im(f)}(R\bigl[\oneover{f}\bigr]),\Qgen{m}{M\bigl[\oneover{f}\bigr]}\biggr) \ar[d] \\
\biggl(\Der^{\im(g)}(R\bigl[\oneover{g}\bigr]),\Qgen{m}{M\bigl[\oneover{g}\bigr]}\biggr) \ar[r] & \biggl(\Der^{\im(fg)}(R\bigl[\oneover{fg}\bigr]),\Qgen{m}{M\bigl[\oneover{fg}\bigr]}\biggr)
\end{tikzcd}
\]
is a Poincar\'e-Verdier square for any commutative ring, invertible $R$-module $M$ with involution and $f,g \in R$ that span the unit ideal; here, for any $a \in R$ the superscript $\im(a)$ indicates the Poincar\'e subcategory of $\Dperf(R)$ generated by the image of the localisation functor $\Dperf(R) \to \Dperf(R[1/a])$. As a consequence of Corollary~\reftwo{corollary:gwadditive-intro-two} we thus find:

\begin{corollaryintrotwo}
\label{corollary:Mayerviet-intro-two}%
Ler $R$ be a (discrete) commutative ring, $M$ an invertible $R$-module with involution, and $f,g \in R$ elements spanning the unit ideal. Then applying any Verdier-localising functor $\F \colon \Catp \rightarrow \Spa$ such as $\K,\GW$ or $\L$ to the square above yields a cartesian diagram of spectra.
\end{corollaryintrotwo}

The categories $\D^{\im(a)}(R) \subseteq \Dperf(R)$ are dense (and in general proper) stable subcategories. By work of Thomason such subcategories are classified by subgroups $c \subseteq \K_0(R)$, and in the body of the text we treat the case of such control groups in general. In \paperfour we then extend Schlichting's hermitian cofinality theorem to general Poincar\'e $\infty$-categories, which in particular implies that passage to a dense Poincar\'e subcategory does not affect the positive Grothendieck-Witt groups. Combined with Corollary \reftwo{corollary:Mayerviet-intro-two} this extends Schlichting's Mayer-Vietoris sequences \cite{schlichting-mv} from higher symmetric to quadratic and even Grothendieck-Witt groups. %
A more systematic approach, also pursued in \paperfour, considers what we term Karoubi-Grothendieck-Witt spectra $\GGW(\C,\QF)$, a variant of Grothendieck-Witt spectra that is not only Verdier localising, but also invariant under passage to dense subcategories mimicking the passage from connective to non-connective $\K$-theory. Corollary~\reftwo{corollary:Mayerviet-intro-two} then
shows that $\GGW$ satisfies a Mayer-Vietoris principle at the spectrum level. 
In \cite{motives}, this is exploited to show that Karoubi-Grothendieck-Witt spectra of schemes satisfy Zariski (even Nisnevich) descent, recovering and extending Schlichting's results.  %
\medskip

Lastly, we turn to another class of examples of Poincaré \(\infty\)-categories, namely those formed by compact parametrised spectra over a space $B$. The relevance of these examples is already visible in the equivalences
\[
\mathrm A(B) \simeq \K((\Spa/B)^\omega)
\]
describing Waldhausen's $\K$-theory of spaces in the present framework. Given a stable spherical fibration $\xi$ over $B$, there are three important Poincaré structures on $(\Spa/B)^\omega$, the quadratic, symmetric and visible one, all of whose underlying duality is the Costenoble-Waner functor 
\[
E \mapsto \Hom_B(E \boxtimes E, \Delta_! \xi);
\]
here $\boxtimes$ is the exterior tensor product, $\Delta \colon B \rightarrow B \times B$ is the diagonal map, and the subscript $!$ denotes the left adjoint functor to its associated pullback. Then from our identification of the isotropy separation square of $\KR((\Spa/B)^\omega,\QF_\xi^r)$ with $r \in \{\qdr,\sym,\vis\}$ we find:

\begin{corollaryintrotwo}
There are canonical equivalences
\[
\GW((\Spa/B)^\omega,\QF^r_{\xi}) \simeq \LA^r(B,\xi)
\]
and in particular
\[
\Omega^{\infty-1}\LA^r(B,\xi) \simeq |\Cob((\Spa/B)^\omega,\QF^r_{\xi})|
\]
for $r \in \{\qdr,\sym,\vis\}$.
\end{corollaryintrotwo}

Here, $\LA^r(B,\xi)$ denotes the spectra constructed by Weiss and Williams (under the names $\mathrm{LA}_\bullet$, $\mathrm{LA}^\bullet$, and  $\mathrm{VLA}$) in their pursuit of a combination of surgery theory and pseudo-isotopy theory into a direct description of the spaces $\mathcal G(M)/\mathrm{Top}(M)$ for closed manifolds $M$, see \cite{WWIII}. This result unites their work with the recent approaches to the study of diffeomorphism groups at the hands of Galatius and Randal-Williams \cite{Stablemoduli}. In particular, the second part provides a cycle model for the previously rather mysterious spectra $\LA^r(B,\xi)$ that can be used to give a new construction of the celebrated map
\[
\widetilde{\mathrm{Top}}(M)/\mathrm{Top}(M) \longrightarrow \Omega^{\infty+1}(\mathrm{Wh}(M)_\hC)
\]
from \cite{WWIII} along with a new proof of the index theorems of Weiss and Williams set to appear in \cite{WWIV}. These results will appear in future work. 

\subsubsection*{Related recent work} 
Firstly, in \cite{Spitzweckreal} Heine, Spitzweck and Verdugo also construct a real algebraic $\K$-spectrum in greater generality than in the present paper (in particular, for not necessarily stable $\infty$-categories), with a version of Theorem \reftwo{theorem:gwuniversal-intro-two}, part \reftwoitem{item:universal-properties-GW} as their main result, albeit using a weaker notion of additivity than the one we use here (resulting in a logically incomparable result). We compare our construction of Grothendieck-Witt spectra to theirs in Appendix \reftwo{appendix:AppIIB}, but do not pursue a full comparison of the real algebraic $\K$-spectra here. 

Secondly, Schlichting announced results similar to the corollaries of our main theorem, and some of the applications we pursue in the third instalment of this series in \cite{Schlichting-integers}, but see also \cite{Schlichting-wrong}, during the completion of this work. He uses the set-up from \cite{SchlichtinghigherI} and a comparison to our construction follows from the main results of \cite{comparison}.

Neither of these papers systematically relates Grothendieck-Witt theory to $\L$-theory, a main thread of our work.

\introsubsection{Organisation of the paper}

In the next section we briefly summarise the necessary results of \paperone, providing in particular a guide to the requisite parts. In \S\reftwo{section:verdier}, we study (co)fibre sequences in $\Catp$ in detail and introduce additive and localising functors. The analogous results in the setting of stable $\infty$-categories, on which our results are based, are well-known but seem difficult to locate coherently in the literature. We therefore give a systematic account in Appendix \reftwo{appendix:AppIIA}, without any claim of originality.

The real work of the present paper then starts in \S\reftwo{section:cobcats}. It contains the definition of the hermitian $\Q$-construction and the algebraic cobordism category and proves Theorem \reftwo{theorem:additivity-intro-two} as \reftwo{theorem:additivity} and \reftwo{theorem:fulladditivity}. In \S\reftwo{section:structure-theory} we then analyse the behaviour of arbitrary additive functors $\Catp \rightarrow \Spa$ and $\Catp \rightarrow \Sps$. This leads to very general versions of Theorem \reftwo{theorem:gwuniversal-intro-two} in \reftwo{corollary:universal}, \reftwo{proposition:positive-om}, our \reftwo{theorem:main-intro-two} in \reftwo{corollary:bord-seq-unique} and Theorem~\reftwo{theorem:kr-intro-two} in \reftwo{theorem:periodicity}. We then obtain all other results of this introduction as simple consequences in \S\reftwo{section:GW}, where we specialise the discussion to the Grothendieck-Witt functor.

Finally, there is a second appendix which establishes two comparison results to other Grothendieck-Witt spectra, not immediate from the results of \cite{comparison}. They are not used elsewhere in the paper.

\introsubsection{Acknowledgements}
For useful discussions about our project,
we heartily thank
Tobias Barthel,
Clark Barwick,
Lukas Brantner,
Mauricio Bustamante,
Denis-Charles Cisinski,
Dustin Clausen,
Diarmuid Crowley,
Uriya First,
Rune Haugseng,
André Henriques,
Lars Hesselholt,
Gijs Heuts,
Geoffroy Horel,
Marc Hoyois,
Max Karoubi,
Daniel Kasprowski,
Ben Knudsen,
Manuel Krannich,
Achim Krause,
Henning Krause,
Sander Kupers,
Wolfgang Lück,
Ib Madsen,
Cary Malkiewich,
Mike Mandell,
Akhil Mathew,
Lennart Meier,
Irakli Patchkoria,
Nathan Perlmutter,
Maxime Ramzi,
Oscar Randal-Williams,
Andrew Ranicki,
George Raptis,
Marco Schlichting,
Peter Scholze,
Stefan Schwede,
Graeme Segal,
Markus Spitzweck,
Jan Steinebrunner,
Georg Tamme,
Ulrike Tillmann,
Michael Weiss,
Christoph Winges,
and
Maria Yakerson.

Furthermore, we owe a tremendous intellectual debt to Jacob Lurie for creating the framework we exploit here, and to Søren Galatius for originally spotting the overlap between various separate projects of ours; his insight ultimately led to the present collaboration.

\medskip

The authors would also like to thank the Hausdorff Center for Mathematics at the University of Bonn, the Newton Institute at the University of Cambridge, the University of Copenhagen and the Mathematical Research Institute Oberwolfach for hospitality and support while parts of this project were undertaken.
\medskip

BC was supported by the French National Centre for Scientific Research (CNRS) through a ``délégation'' at LAGA, University Paris 13, and by the french French National Research Agency (ANR) project ``Motivic homotopy, quadratic invariants and diagonal classes'' (ANR grant no.\ 21-CE40-0015) at the University of Burgundy.
ED was supported by the German Research Foundation (DFG) through the priority program ``Homotopy theory and Algebraic Geometry'' (DFG grant no.\ SPP 1786) at the University of Bonn and WS by the priority program ``Geometry at Infinity'' (DFG grant no.\ SPP 2026) at the University of Augsburg. ED was further supported by the Engineering and Physical Sciences Research Council (EPSRC) through the grant ``Characteristic polynomials for symmetric forms'' (EPSRC grant no.\ EP/W019620/1) at the University of Warwick.
YH was supported by the European Research Council (ERC) through the project ``Foundations of Motivic Real K-Theory'' (ERC grant no.\ 949583) at Paris Cité University.
YH and DN were supported by the ANR through the grant ``Chromatic Homotopy and K-theory'' (ANR grant no.\ 16-CE40-0003) at LAGA, University of Paris 13.
FH was a member of the Hausdorff Center for Mathematics (DFG grant no.\ EXC 2047 390685813) at the University of Bonn and TN of the cluster ``Mathematics Münster: Dynamics-Geometry-Structure'' (DFG grant no.\ EXC 2044 390685587) at the University of Münster. 
FH was further supported by the DFG through the collaborative research center ``Integral structures in Geometry and Representation Theory'' (DFG grant no.\ TRR 358 4913924) at the University of Bielefeld, TN through the centre ``Geometry: Deformations and Rigidity" (DFG grant no.\ SFB 1442 427320536) at the University of Münster and ML and DN through the centre ``Higher Invariants'' (DFG grant no.\ SFB 1085 224262486) at the University of Regensburg. 
FH was also supported by the ERC through the project ``Moduli spaces, Manifolds and Arithmetic'' (ERC grant no.\ 682922) of Søren Galatius and KM by the project ``$\K$-theory, $\L^2$-invariants, manifolds, groups and their interactions'' (ERC grant no.\ 662400) of Wolfgang Lück. 
FH, TN and WS were further supported by the EPSRC through the program ``Homotopy harnessing higher structures'' at the Isaac Newton Institute for Mathematical Sciences (EPSRC grants no.\ EP/K032208/1 and EP/R014604/1). 
ML was supported by the research fellowship ``New methods in algebraic K-theory'' (DFG grant no.\ 424239956) and by the Danish National Research Foundation (DNRF) through the Center for Symmetry and Deformation (DNRF grant no.\ 92) and the Copenhagen Centre for Geometry and Topology (DNRF grant no.\ 151) at the University of Copenhagen. 
KM was also supported by the K\&A Wallenberg Foundation at the University of Stockholm.

\section*{Recollection}

In the present section, we briefly recall the parts of \paperone that are most relevant for the considerations of the present paper. We first summarise the abstract features of the theory, and then spell out some examples.

\subsubsection*{Poincaré $\infty$-categories and Poincaré objects}
Recall from \S\refone{subsection:hermitian-and-poincare-cats} that an hermitian structure on a small stable $\infty$-category $\C$ is a reduced, quadratic functor $\QF \colon \C\op \rightarrow \Spa$; we will discuss such functor in detail momentarily. We call a pair $(\C,\QF)$ consisting of these data an hermitian \(\infty\)-category. These organise into an $\infty$-category $\Cath$ whose morphisms $(\C,\QF) \rightarrow (\D,\QFD)$ consist of what we term hermitian functors, that is pairs $(f,\eta)$ where $f \colon \C \rightarrow \D$ is an exact functor and $\eta \colon \QF \Rightarrow \QFD \circ f\op$ is a natural transformation.

To such an hermitian \(\infty\)-category is associated its $\infty$-category of hermitian forms $\catforms(\C,\QF)$, whose objects consist of pairs $(X,q)$ where $X \in \C$ and $q$ is a $\QF$-hermitian form on $X$, i.e.\ a point in $\Omega^\infty \QF(X)$, see \S\refone{subsection:hermitian-poincare-objects}. Morphisms are maps in $\C$ preserving the hermitian forms. The core of the $\infty$-category $\catforms(\C,\QF)$ is denoted $\spsforms(\C,\QF)$, and these assemble into functors
\[
\catforms \colon \Cath \lrar  \Cat \quad \text{and} \quad \spsforms \colon \Cath \lrar \Sps.
\]
In order to impose a non-degeneracy condition on the forms in $\spsforms(\C,\QF)$, one needs a non-degeneracy condition on the hermitian \(\infty\)-category $(\C,\QF)$ itself. To this end recall the classification of quadratic functors from Goodwillie calculus: Any reduced quadratic functor uniquely extends to a natural cartesian diagram
\begin{equation*}
\begin{tikzcd}
\QF(X) \ar[r] \ar[d] & \Lin_\QF(X) \ar[d]\\
\Bil_\QF(X,X)^\hC \ar[r] & \Bil_\QF(X,X)^\tC
\end{tikzcd}
\end{equation*}
where $\Lin_\QF \colon \C\op \rightarrow \Spa$ is linear (i.e.\ exact), and $\Bil_\QF \colon \C\op \times \C\op \rightarrow \Spa$ is bilinear (i.e.\ exact in each variable) and symmetric (i.e.\ comes equipped with a refinement to an element of $\Fun(\C\op \times \C\op,\Spa)^\hC$, with $\Ct$ acting by flipping the input variables); see \S\refone{subsection:classification}.

A hermitian structure $\QF$ is called Poincaré if there exists an equivalence $\Dual \colon \C\op \rightarrow \C$ such that
\[
\Bil_\QF(X,Y) \simeq \Hom_\C(X,\Dual Y)
\]
naturally in $X,Y \in \C\op$. By Yoneda's lemma, such a functor $\Dual$ is uniquely determined if it exists, so we refer to it as $\Dual_\QF$. By the symmetry of $\B_\QF$, the functor $\Dual_\QF$ then automatically satisfies $\Dual_\QF \circ \Dual_\QF\op \simeq \id_\C$. Any hermitian functor $(f,\eta) \colon (\C,\QF) \rightarrow (\D,\QFD)$ between Poincaré $\infty$-categories (i.e.\ hermitian $\infty$-categories whose hermitian structure is Poincaré) induces a canonical map
$
f \circ \Dual_\QF \Rightarrow \Dual_\QFD \circ f\op;
$
see \S\refone{subsection:hermitian-and-poincare-cats}.
We say that $(F,\eta)$ is a Poincaré functor if this transformation is an equivalence, and Poincaré $\infty$-categories together with Poincaré functors form a (non-full) subcategory $\Catp$ of $\Cath$. \\

Now, if $(\C,\QF)$ is Poincaré, then to any hermitian form $(X,q) \in \spsforms(\C,\QF)$ there is canonically associated a map
$
q_\sharp \colon X \to \Dual_\QF X
$
as the image of $q$ under
\[
\Omega^\infty \QF(X) \longrightarrow \Omega^\infty\Bil_\QF(X,X) \simeq \Hom_\C(X,\Dual_\QF X)
\]
and we say that $(X,q)$ is Poincaré if $q_\sharp$ is an equivalence. The full subspace of $\spsforms(\C,\QF)$ spanned by the Poincaré forms is denoted by $\Poinc(\C,\QF)$ and provides a functor
$
\Poinc \colon \Catp \rightarrow \Sps,
$
which we suggest to view in analogy with the functor $\core \colon \Catx \rightarrow \Sps$ taking a stable $\infty$-category to its groupoid core. Details about this functor are spelled out in \S\refone{subsection:hermitian-poincare-objects}.
\medskip

Hermitian structures can be shifted by forming $\QF\qshift{i} = \SS^i \otimes \QF$, and the simplest example of a Poincaré $\infty$-category to keep in mind is $\C = \Dperf(R)$, where $R$ is a discrete commutative ring and $\Dperf(R)$ is the $\infty$-category of perfect complexes over $R$ (i.e.\ finite chain complexes of finitely generated projective $R$-modules), together with the symmetric and quadratic Poincaré structures given by 
\[
\QF^\qdr_R(X) \simeq \mathbb \hom_{R}(X \otimes^{\mathbb L}_R X,R)_\hC \quad \text{and} \quad \QF^\sym_R(X) \simeq \mathbb \hom_{R}(X \otimes^{\mathbb L}_R X,R)^\hC,
\]
where $\hom_R$ denotes the mapping spectrum of the stable $\infty$-category $\Dperf(R)$ (in other words the spectrum underlying the derived mapping complex $\mathbb R\Hom_R$).
In either case, the bilinear part and duality are given by
\[
\Bil(X,Y) \simeq \hom_R(X \otimes_R^\mathbb L Y, R) \quad \text{and} \quad \Dual(X) \simeq \mathbb R\Hom_R(X,R),
\]
which makes both $\QF^\sym_R$ and $\QF^\qdr_R$ into Poincaré structures on $\Dperf(R)$.

\subsubsection*{Constructions of Poincaré $\infty$-categories}

We next collect a few important structural properties of the $\infty$-categories $\Cath$ and $\Catp$, which are the source of key constructions of Poincaré $\infty$-categories.
First of all, by the results of \S\refone{subsection:limits} they are both complete and cocomplete, and the inclusion $\Catp \rightarrow \Cath$ is detects equivalences.
Furthermore, the forgetful functors
\[
\Catp \longrightarrow \Cath \longrightarrow \Catx
\]
both possess both adjoints, so preserve both limits and colimits; the adjoints are constructed in \S\refone{subsection:bilinear-and-pairings} and \S\refone{subsection:thom}. For the right hand functor, both of them simply equip a stable $\infty$-category $\C$ with the trivial hermitian structure $0$. 
The right adjoint to the left hand functor is given by the pairings construction $(\D,\QFD) \mapsto \Pairings(\D,\QFD)$ studied in \S\refone{subsection:thom}.  Its underlying $\infty$-categoy is given by the bivariant unstraightening of $\Omega^\infty \Bil_\QFD \colon \D\op \times \D\op \rightarrow \Sps$, and its objects are thus given by triples $(X,Y,\beta)$ with $X \in \D,Y \in \D\op$ and $\beta \in \Omega^\infty\Bil_\QFD(X,Y)$. Its hermitian structure $\QFD_\pair$ takes $(X,Y,\beta)$ to the pullback of $\QFD(X) \rightarrow \Bil_\QFD(X,X) \leftarrow \hom_\D(X,Y)$, with right hand map given by pullback of $\beta$ and the adjunction counit simply projects $(X,Y,\beta)$ to $X$ with the evident hermitian refinement. The left adjoint is given by $(\D,\QFD) \mapsto \Pairings(\D,\QFD\qshift{-1})$, with unit underlain by $X \mapsto (\Omega X,0,0)$.

The following two special cases
are of great importance. By the above discussion, the left and right adjoint of the composite $\Catp \rightarrow \Catx$ agree. They are given by the hyperbolic construction $\C \rightarrow \Hyp(\C)$ with underlying $\infty$-category $\C \times \C\op$ and Poincaré structure $\hom_\C \colon (\C \times \C\op)\op \rightarrow \Spa$, see \S\refone{subsection:hyp-and-sym-poincare-objects}. The associated duality is given by $(X,Y) \mapsto (Y,X)$, and there is a natural equivalence
$
\core \C \simeq \Poinc \Hyp(\C)
$
implemented by $X \mapsto (X,X)$. 

The other important case is the composite of the inclusion $\Catp \rightarrow \Cath$ with either of its adjoints, see \S\refone{subsection:metabolic-and-L}: For $(\C,\QF)$ a Poincar\'e $\infty$-category the duality identifies the underlying category of $\Pairings(\C,\QF)$ with $\Ar(\C)$, and the Poincaré structure thereon becomes $\QF_\arr(X \rightarrow Y) \simeq \QF(X) \times_{\Bil_\QF(X,X)} \Bil_\QF(X,Y)$, resulting in the Poincar\'e $\infty$-category $\Ar(\C,\QF) = (\Ar(\C),\QF_\arr)$. The equivalence $\cof \colon \Ar(\C) \rightarrow \Ar(\C)$ then translates the shifted hermitian structure $(\QF_\arr)\qshift{-1}$ into
\[
\QF_\met(X \rightarrow Y) \simeq \fib(\QF(Y) \rightarrow \QF(X)).
\]
The left adjoint therefore takes $(\C,\QF)$ to $\Met(\C,\QF) = (\Ar(\C),\QF_\met)$, the metabolic Poincar\'e $\infty$-category of $(\C,\QF)$, which comes equipped with an equivalence $\cof \colon \Ar(\C,\QF)\qshift{-1} \simeq \Met(\C,\QF)$. The associated duality is
\[
\Dual_{\QF_\met}(X \rightarrow Y)  \simeq \fib(\Dual_\QF(Y) \rightarrow \Dual_\QF(X)) \longrightarrow \Dual_\QF(Y)
\]
and Poincaré objects in $\Met(\C,\QF)$ are best thought of as Poincaré objects with boundary in the Poincaré $\infty$-category $(\C,\QF\qshift{-1})$ which embeds into $\Met(\C,\QF)$ via $X \mapsto (X \rightarrow 0)$, i.e.\ as the objects with trivial boundary.

From the various adjunction units and counits there then arises a diagram
\[
\begin{tikzcd}[cramped]
 & \Hyp(\C) \ar[ld,equal] \ar[rd,"\hyp"] \ar[dd,"\ilag"] &  \\
\Hyp(\C) & & (\C,\QF) \\
 & \ar[lu,"\lag"] \Met(\C,\QF) \ar[ru,"\met"'] & 
\end{tikzcd}
\]
in $\Catp$ for every Poincaré $\infty$-category; the underlying functors pointing to the right are given by 
\[
\met(X \rightarrow Y) = Y
\quad \text{and} \quad 
\hyp(X,Y) = X \oplus \Dual_\QF Y,
\]
whereas the other two are given by extending the source and identity functors
\[
s \colon \Ar(\C) \longrightarrow \C
\quad \text{and} \quad
\id \colon \C \longrightarrow \Ar(\C)
\]
using the adjunction properties of $\Hyp$. Regarding the induced maps after applying $\Poinc$, one finds that an element $(X,q) \in \pi_0 \Poinc(\C,\QF)$ is in the image of $\met$ if it admits a Lagrangian, that is a map $f \colon L \rightarrow X$ such that there is an equivalence $f^*q \simeq 0$, whose associated nullhomotopy of the composite
\[
L \xrightarrow{f} X \simeq \Dual_\QF X \xrightarrow{\Dual_\QF f} \Dual_\QF L
\]
makes this sequence into a fibre sequence in $\C$. In particular, this is the case for all Poincar\'e objects in a metabolic Poincar\'e $\infty$-category. Similarly, $(X,q)$ lies in the image of $\hyp$ if there is an equivalence $X \simeq L \oplus \Dual_\QF L$ which translates the form $q$ into the canonical evaluation form on the target. 

Thus the Poincaré categories $\Hyp(\C)$ and $\Met(\C,\QF)$ encode the theory of metabolic and hyperbolic forms in $(\C,\QF)$, and the remainder of the diagram witnesses that any hyperbolic form has a canonical Lagrangian, from which it can be reconstructed.

One further property of these constructions that we shall need is that the duality $\Dual_\QF$ equips the underlying $\infty$-category of $(\C,\QF)$ with the structure of a homotopy fixed point in $\Catx$ under the $\Ct$-action given by taking $\C$ to $\C\op$, or in other words, the forgetful functor $\Catp \rightarrow \Catx$ is $\Ct$-equivariant for the trivial action on the source and the opponing action on the target; see \S\refone{subsection:bilinear-and-pairings}. As a formal consequence, its adjoint $\Hyp$ is equivariant as well, and thus the composite
\[
\Catp \xrightarrow{\fgt} \Catx \xrightarrow{\Hyp} \Catp
\]
lifts to a functor $\HypG \colon \Catp \rightarrow (\Catp)^\hC = \Fun(\mathrm{B}\Ct,\Catp)$, the $\infty$-category of $\Ct$-objects in $\Catp$. The action map on $\HypG(\C,\QF)$ is given by the composite
\[
\Hyp(\C) \xrightarrow{\mathrm{flip}} \Hyp(\C\op) \xrightarrow{\Hyp(\Dual_\QF)} \Hyp(\C)
\]
and the functor $\hyp \colon \Hyp(\C) \rightarrow (\C,\QF)$ is invariant under the action on the source. \\

We recall from \S\refone{subsection:monoidal-structure} that the $\infty$-category $\Cath$ admits a symmetric monoidal structure making the functor $\fgt \colon \Cath \rightarrow \Catx$ symmetric monoidal for Lurie's tensor product of stable $\infty$-categories on the target. While we do not use the monoidal structure much in the present paper, we heavily exploit the following: The monoidal structure on $\Cath$ is cartesian closed, i.e.\ $\Cath$ admits internal function objects, and also both tensors and cotensors over $\Cat$, see \S\refone{subsection:internal}, \S\refone{subsection:tensoring} and \S\refone{subsection:cotensoring}, respectively. More explicitly, to hermitian \(\infty\)-categories $(\C,\QF)$ and $(\D,\QFD)$ and an ordinary category $\I$, there are associated hermitian \(\infty\)-categories
\[
\Funx((\C,\QF),(\D,\QFD)), \quad (\C,\QF)_\I \quad \text{and} \quad (\C,\QF)^\I
\]
connected by natural equivalences
\[
\Funx((\C,\QF)_\I,(\D,\QFD)) \simeq \Funx((\C,\QF),(\D,\QFD))^\I \simeq \Funx((\C,\QF),(\D,\QFD)^\I)
\]
and whose underlying $\infty$-categories in the outer cases are given by 
\[
\Funx(\C,\D) \quad \text{and} \quad \Fun(\I,\C).
\]
This results in particular in equivalences
\[
\spsforms\Funx((\C,\QF),(\D,\QFD)) \simeq \Hom_{\Cath}((\C,\QF),(\D,\QFD)), \quad \Poinc\Funx((\C,\QF),(\D,\QFD)) \simeq \Hom_{\Catp}((\C,\QF),(\D,\QFD))
\]
\[
\text{and} \quad \catforms((\C,\QF)^\I) \simeq \Fun(\I,\catforms(\C,\QF)).
\]
The tensoring construction is unfortunately far less explicit, but 
for $\I$ a finite poset it is described in Proposition~\refone{proposition:tensor-strongly-finite}. Finally, we note that neither the tensor nor cotensor construction generally preserve Poincaré \(\infty\)-categories, though Lurie established sufficient criteria which we recorded in \S\refone{subsection:finite-complexes}.

\subsubsection*{Examples of Poincaré $\infty$-categories}

We now recall the most relevant examples in detail: first, Poincaré structures on module $\infty$-categories and then parametrised spectra. 
\medskip

Fix an $\Eone$-algebra $A$ over a base $\Einf$-ring spectrum $k$ and a subgroup $c \subseteq \K_0(A)$. Consider then the $\infty$-category of compact $A$-module spectra $\Modp{A}$ or more generally its full subcategory $\Mod^c(A)$ spanned by all those $X \in \Modp{A}$ with $[X] \in c \subseteq \K_0(A)$. For example, $\Modf{A} = \Mod^{\langle A \rangle}(A)$ is the stable subcategory of $\Mod(A)$ generated by $A$ itself. 
Any discrete ring $R$ gives rise to such data, via the functor $\GEM \colon \Ring \rightarrow \Alg_{\Eone}(\Mod({\GEM\ZZ}))$ induced by the lax symmetric monoidal Eilenberg-Mac Lane functor $\GEM \colon \Ab \rightarrow \Spa$.
There are, furthermore, equivalences
\[
\Modp{\GEM R} \simeq \Dperf(R) \quad \text{and} \quad \Modf{\GEM R} \simeq \Dfree(R),
\]
where $\Dperf(R)$ and $\Dfree(R)$ denote the full subcategory of the derived $\infty$-category $\D(R)$ of $R$ spanned by the finite chain complexes of finitely generated projective $R$-modules and of free $R$-modules, respectively.
In this regime, terms such as $\otimes_{\GEM \ZZ}$ or $\Hom_{\GEM R}$ correspond to the functors $\otimes^\mathbb L_\ZZ$ and $\mathbb R\Hom_R$.
\medskip

According to \S\refone{subsection:genuine-modules}, hermitian structures on the $\infty$-categories $\Mod^c(A)$ are classified by $A$-modules with genuine involution $(M,N,\alpha)$; the first entry $M$ is an $A \otimes_k A$-module, equipped with the structure of a homotopy fixed point in the $\infty$-category $\Mod({A \otimes_k A})$ under the $\Ct$-action flipping the two factors; see \S\refone{subsection:modules-with-involution}; we call it a module with (naive) involution over $A$.

The simplest example of such a structure over a discrete ring $R$ is given by a discrete $R \otimes_\ZZ R$-module $M$, and a selfmap $M \rightarrow M$ that squares to the identity on $M$, and is semilinear for the flip map of $R \otimes_\ZZ R$. If $R$ is equipped with an anti-involution $\sigma$, then $M = R$ with involution $\sigma$ or $-\sigma$ (or $\epsilon \sigma$ for any other central unit $\epsilon$ with $\sigma(\epsilon) = \epsilon^{-1}$) is a valid choice, using $\sigma$ to turn the usual $R \otimes_\ZZ R\op$-module structure on $R$ into an $R \otimes_\ZZ R$-module structure.

The second entry $N$ is an $A$-module spectrum and the third entry $\alpha$ is an $A$-linear map $N \rightarrow M^\tC$, where the $A$-module structure on $M^\tC$ is obtained via the Tate diagaonal $A \rightarrow (A \otimes_k A)^\tC$, which is a map of $\Eone$-ring spectra; see \S\refone{subsection:genuine-modules} for details. 

Even if only interested in discrete rings, one therefore has to leave the realm of discrete $R$-modules to form the Tate construction, and even the realm of derived categories 
in order to consider the Tate diagonal.
\medskip

The hermitian structure associated to a module with genuine involution $(M,N,\alpha)$
is given by the pullback
\[
\begin{tikzcd}
[column sep=3ex]
\QF_M^\alpha(X) \ar[rr] \ar[d] & & \hom_A(X,N) \ar[d,"{\alpha_*}"] \\
  \hom_{A \otimes_k A}(X \otimes_k X,M)^\hC \ar[r] & \hom_{A \otimes_k A}(X \otimes_k X,M)^\tC \ar[r,"\simeq"] & \hom_A(X,M^\tC)
\end{tikzcd}
\]
where the $\Ct$-action on $\hom_{A \otimes_k A}(X \otimes_k X,M)$ flips the factors in the source and applies the involution on $M$. It is a Poincaré structure on $\Modp A$ if $M$ restricts to an object of $\Modp{A}$ under either inclusion $A \rightarrow A \otimes_k A$, and furthermore $M$ is invertible, i.e.\ the natural map
$A \rightarrow \hom_A(M,M)$
is an equivalence. In this case, the associated duality is given by $X \mapsto \hom_A(X,M)$ regarded as an $A$-module via the extraneous $A$-module structure on $M$; see again \S\refone{subsection:modules-with-involution}. Given a subgroup $c \in \K_0(A)$, one obtains a Poincaré structure on $\Mod^c(A)$ if in addition $c$ is closed under the duality on $\K_0(A)$ induced by $M$. In the example of $\Modf{A}$ this translates to $M \in \Modf{A}$. \\

Let us give some concrete examples of such genuine structures in the case of a discrete ring $R$ and of a discrete invertible module with involution $M$ over $R$.
From this data are most easily defined the quadratic and symmetric Poincaré structures $\QF^\qdr_M$ and $\QF^\sym_M$ given by
\[
\QF^\qdr(X) = \hom_{R \otimes_\ZZ^\mathbb L R}(X \otimes_\ZZ^\mathbb L X,M)_\hC 
\quad \text{and} \quad 
\QF^\sym(X) = \hom_{R \otimes_\ZZ^\mathbb L R}(X \otimes_\ZZ^\mathbb L X,M)^\hC
\]
which correspond to the  modules with genuine involution
$(M,0,0)$ and $(M,M^\tC, \id)$,
respectively.
Interpolating between these, the genuine family of Poincaré structures $\Qgen{i}{M}$ correspond to the modules with genuine involution $(M,\tau_{\geq i}M^\tC, \tau_{\geq i}M^\tC \rightarrow M^\tC)$ for $i \in \ZZ$. 
These intermediaries are mostly important because they contain the following examples: The functors 
\[
\Quad_M, \quad \Even_M, \quad \text{and} \quad \Sym_M \colon \Proj(R)\op \longrightarrow \Ab
\]
assigning to a finitely generated projective module its abelian group of $M$-valued quadratic, even or symmetric forms, respectively, admit animations (or non-abelian derived functors in more classical terminology) which we term 
\[
\QF^\gq_M, \quad \QF^\gev_M \quad \text{and} \quad \QF^\gs_M \colon \Dperf(R)\op \longrightarrow \Spa,
\]
respectively. One of the main results of \paperone are equivalences
\[
\QF^\gq_M \simeq \Qgen{2}{M}, \quad \QF^\gev_M \simeq \Qgen{1}{M} \quad \text{and} \quad \QF^\gs_M \simeq \Qgen{0}{M}.
\]
No further members of the genuine family arise as animations of functors $\Proj(R) \rightarrow \Ab$; 
see \S\refone{subsection:discrete-rings}.
\medskip

Turning to a different kind of example, consider the $\infty$-categories $\Spa_B = \Fun(B,\Spa)$ for some $B \in \Sps$. Entirely parallel to the discussion above, one can derive hermitian structures on the compact objects of $\Spa_B$ from triples $(M,N, \alpha)$ with $M \in (\Spa_{B \times B})^\hC$ and $\alpha \colon N \rightarrow (\Delta^*M)^\tC$ a map in $\Spa_B$, where $\Delta \colon B \rightarrow B \times B$ is the diagonal, \S\refone{subsection:parametrised-spectra}. The most important examples of such functors are the visible Poincaré structures $\QF_\xi^\vis$ given by the triples 
\[
(\Delta_! \xi, \xi, u \colon \xi \rightarrow (\Delta^*\Delta_!\xi)^\tC),
\]
where $\xi \colon B \rightarrow \Pic(\SS)$ is some stable spherical fibration over $B$, where $\Delta_! \colon \Spa_B \rightarrow \Spa_{B \times B}$ is the left adjoint to $\Delta^*$ and where $u$ is the unit of this adjunction (which factors through $\xi \rightarrow (\Delta^*\Delta_!\xi)^\hC$ since $\Delta$ is invariant under the $\Ct$-action on $B \times B$). These hermitian structures are Poincaré with associated duality given by 
\[
X \longmapsto \hom_B(X,\Delta_!\xi),
\]
the Costenoble-Waner duality functor twisted by $\xi$. The reason these are so important is that any closed manifold $M$, in fact any Poincaré complex, defines a very interesting element, its visible symmetric signature, in $\Poinc(\Spa_M^\cp,\QF_{\nu}^\vis)$, where $\nu \colon M \rightarrow \Pic(\SS)$ is the stable normal bundle of a manifold or more generally the Spivak fibration of a Poincaré complex.

\section{Poincaré-Verdier sequences and additive functors}
\label{section:verdier}%

In this section, we study the analogue of (split) Verdier sequences in the context of Poincaré \(\infty\)-categories, as well as their analogue for idempotent complete Poincaré \(\infty\)-categories, which, following a suggestion of Clausen and Scholze, we call \defi{Karoubi sequences}. 
In particular, our terminology differs from that of Blumberg-Gepner-Tabuada~\cite{BGT}; see Appendix~\reftwo{appendix:AppIIA} for a thorough discussion. 

After developing the example of module \(\infty\)-categories in some detail, we proceed to introduce the notions of \defi{additive}, \defi{Verdier-localising}, and \defi{Karoubi-localising} functors \(\Catp \rightarrow \E\), encoding the preservation of an increasing number of such sequences, or rather, for not necessarily stable 
\(\E\), of a mild generalisation thereof in the form of certain cartesian and cocartesian squares in \(\Catp\). These three notions we introduce correspond loosely to satisfying Waldhausen's additivity theorem, Quillen's localisation theorem and Bass' strengthening thereof. 

The notion of an additive functor from \(\Catp\) to \(\Spa\) is central in our work, since it essentially abstracts the minimal property enjoyed by our main subject of interest, the functor \(\GW\colon \Catp \to \Spa\) (only to be defined in Definition~\reftwo{definition:gw-spectrum}), that permits us to develop a general theory decomposing it into simpler pieces; furthermore \(\GW\) is characterised as the universal such additive functor with a transformation from the functor \(\Poinc\), extracting the space of Poincaré forms. Analogously to \(\K\)-theory, the functor \(\GW\) turns out to be furthermore Verdier-localising, justifying as well the study of that notion. Finally, just as non-connective \(\K\)-theory relates to \(\K\)-theory, the search for a Karoubi-localising approximation of \(\GW\) will yield in \paperfour the Karoubi-Grothendieck-Witt spectrum functor \(\KGW\).

In the present section, we only give the very basic properties of such functors, as the only immediately interesting examples are the space valued functors \(\Core\) and \(\Poinc\). 
After \S\reftwo{section:cobcats}
introduces more interesting examples, we return to a detailed study of additive  functors in \S\reftwo{section:structure-theory}.
The study of Karoubi-localising functors will be taken up in \paperfour.

\subsection{Poincaré-Verdier sequences}
\label{subsection:poincare-verdier-sequences}%

As the basis for our study we require a rather detailed analysis of Verdier sequences in the set-up of stable \(\infty\)-categories. Essentially all of the results we need seem well-known, but have not been coherently organised. 
We have largely collected such statements and their proofs into Appendix~\reftwo{appendix:AppIIA}, the focus of the present section being on incorporating Poincaré structures.	

A sequence 
\[
\C \xrightarrow{f} \D \xrightarrow{p} \E
\]
in \(\Catx\) with vanishing composite is a \emph{Verdier sequence} (Definition~\reftwo{definition:verdier}) if it is both a fibre and a cofibre sequence in \(\Catx\), in which case we refer to \(f\) as a \emph{Verdier inclusion} and to \(p\) as a \emph{Verdier projection}. 
We also say that the sequence is \emph{split} if \(p\), or equivalently, \(f\), admits \emph{both} adjoints (see Definition~\reftwo{definition:split-verdier} and Corollary~\reftwo{corollary:equivalence-split}). 

\begin{definition}
\label{definition:poincare-verdier}%
A sequence 
\[
(\C,\QF) \xrightarrow{(f,\eta)} (\D, \QFD) \xrightarrow{(p,\vartheta)} (\E, \QFE)
\]
of Poincaré functors with vanishing composite is called a \defi{Poincaré-Verdier sequence} if it is both a fibre sequence and a cofibre sequence in \(\Catp\), in which case we call \((f,\eta)\) a \defi{Poincaré-Verdier inclusion} and \((p,\vartheta)\) a \defi{Poincaré-Verdier projection}. We shall say that the sequence is \defi{split} if the underlying Verdier sequence splits.
\end{definition}

We have collected examples of interest of (split) Poincaré-Verdier sequences into Section \reftwo{subsection:examples-poincare-verdier} below, and encourage the reader yearning for them to jump ahead to that section. In the present section we develop the theory governing such sequences.

\begin{remark}
\label{remark:poincare-verdier-fiber-square}%
By a sequence with vanishing composite we implicitly mean a square in \(\Catp\) whose (say) bottom left corner is the zero Poincaré \(\infty\)-category. We note that if a composable pair of Poincaré functors extends to such a square then it does so essentially uniquely, since the zero exact functor between any two Poincaré \(\infty\)-categories admits an essentially unique Poincaré structure.
\end{remark}

\begin{remark}
Since the forgetful and hyperbolic functors are both-sided adjoints to one another, we immediately find that the underlying sequence \(\C \rightarrow \D \rightarrow \E\) of a (split) Poincaré-Verdier sequence is a (split) Verdier sequence, and that the hyperbolisation of any (split) Verdier sequence is a (split) Poincaré-Verdier sequence (where we note that the action of \(\C \mapsto \Hyp(\C) = \C \times \C\op\) on underlying exact functors preserves the property of having \emph{both} adjoints). 
\end{remark}

We now proceed to consider Poincaré-Verdier sequences more closely. Recall that the inclusion \(\Catp \to \Cath\) preserves both limits and colimits (Proposition~\refone{proposition:Catp-cocomplete}), and since it is also conservative, it moreover detects limits and colimits. We may hence test if a given sequence of Poincaré \(\infty\)-categories is a (co)fibre sequence at the level of \(\Cath\).
Now, limits in \(\Cath\) are computed by first taking the limit \(\D\) of underlying stable \(\infty\)-categories, then pulling back all the quadratic functors to \(\D\op\), and finally calculating the limit of the resulting diagram in the \(\infty\)-category of quadratic functors on \(\D\op\); the dual version holds for colimits. 
Finally, we also note that limits and colimits in \(\Funq(\D\op,\Spa)\), i.e.\ of quadratic functors, can be computed in \(\Fun(\D\op,\Spa)\), see Remark~\refone{remark:closed}, and that the operation of left Kan extending along an exact functor preserves quadratic functors, see Lemma~\refone{lemma:kan-extension-exact-quadratic}(iii).

\begin{proposition}
\label{proposition:criterion-poincare}%
Let 
\[
(\C,\QF) \xrightarrow{(f,\eta)} (\D,\QFD) \xrightarrow{(p,\vartheta)} (\E,\QFE)
\]
be a sequence of functors in \(\Catp\) with vanishing composite. Then the following holds:
\begin{enumerate}
\item 
\label{item:fiberCatp}%
It is a fibre sequence in \(\Catp\) if and only if its image in \(\Catx\) is a fibre sequence and
\(\eta \colon \QF \rightarrow f^*\QFD\) is an equivalence. 
\item 
\label{item:cofiberCatp}%
It is a cofibre sequence in \(\Catp\) if and only if its image in \(\Catx\) is a cofibre sequence and
\(\vartheta\colon \QFD \to p^*\QFE\) exhibits \(\QFE \colon \E\op \rightarrow \Spa\) as the left Kan extension of \(\QFD\) along \(p\op\).
\end{enumerate}
In particular, it is a Poincaré-Verdier sequence if and only if its image in \(\Catx\) is a Verdier sequence, and the Poincaré structures on \(\C\) and \(\E\) are obtained from that of \(\D\) by restriction and left Kan extension, respectively.
\end{proposition}
\begin{proof}
Specialising the preceding discussion to the case of squares with one corner the zero Poincaré \(\infty\)-category the sequence from the statement is a fibre sequence in \(\Catp\) if and only if its image in \(\Catx\) is a fibre sequence and \(\QF \to f^*\QFD \to f^*p^*\QFE\) is a fibre sequence in \(\Fun(\C\op,\Spa)\), which, since \(f^*p^*\QFE \simeq 0\), just means that the map \(\QF \to f^*\QFD\) is an equivalence, showing \reftwoitem{item:fiberCatp}. 

Similarly, it is a cofibre sequence in \(\Catp\) if and only if its image in \(\Catx\) is a cofibre sequence and \(p_!f_!\QF \to p_!\QFD \to \QFE\) is a cofibre sequence of quadratic functors, which, since \(p_!f_!\QF \simeq 0\) just means that the map \(p_!\QFD \to \QFE\) is an equivalence, showing \reftwoitem{item:cofiberCatp}.
\end{proof}

Combining Proposition~\reftwo{proposition:criterion-poincare} with Proposition~\reftwo{proposition:inclusion-criterion}, which states that an exact functor \(\C \rightarrow \D\) between stable \(\infty\)-categories is a Verdier inclusion if and only if it is fully faithful and its essential image is closed under retracts in \(\D\), we get:
\begin{corollary}
\label{corollary:poincare-inclusion}%
A Poincaré functor \((f,\eta)\colon (\C,\QF) \to (\D,\QFD)\) is a Poincaré-Verdier inclusion if and only if \(f\) is fully faithful, its essential image is closed under retracts, and the map \(\eta \colon \QF \rightarrow f^*\QFD\) is an equivalence. 
\end{corollary}

To state the analogous corollary concerning Poincaré-Verdier projections, 
recall that the localisation \(\D[W^{-1}]\) of an \(\infty\)-category \(\D\) at a set \(W\) of morphisms 
is the initial \(\infty\)-category under \(\D\) in which the morphisms from \(W\) become invertible. We refer to localisations which are also themselves left/right adjoints as left/right Bousfield localisations (note the difference with the terminology of~\cite{HTT}*{Definition 5.2.7.2}, which calls localisation what we call left Bousfield localisations, see also Corollary~\reftwo{corollary:bousfieldgenII}).

Given an exact functor \(f \colon \C \to \D\), the \emph{Verdier quotient} \(\D/\C\) of \(\D\) by \(\C\) is the localisation of \(\D\) with respect to the collection of maps whose fibre lies in the smallest stable subcategory containing the essential image of \(f\) (see Definition~\reftwo{definition:quotient}). 
By \cite{NS}*{Theorem I.3.3(i)}, \(\D/\C\) is again a stable \(\infty\)-category and the canonical functor \(\D \rightarrow \D/\C\) is exact. For a further discussion of Verdier quotients, we refer the reader to \S\reftwo{section:appendix-verdier}. The main output of the discussion there is Corollary~\reftwo{corollary:criterion-projection}, which shows that an exact functor is a Verdier projection if and only if it is a localisation. Combining this with Proposition~\reftwo{proposition:criterion-poincare}, we get:
\begin{corollary}
\label{corollary:poincare-projection}%
A Poincaré functor \((p,\vartheta)\colon (\D,\QFD) \to (\E,\QFE)\) is a Poincaré-Verdier projection if and only if \(p\colon \D \to \E\) is a localisation and \(\QFD \to p^*\QFE\) exhibits \(\QFE\) as the left Kan extension of \(\QFD\) along \(p\op\).
\end{corollary}

\begin{remark}
\label{remark:compute}%
If $p\colon \D \to \E$ is a Verdier projection and $\QFD\colon \D\op \to \Sp$ a quadratic functor then the left Kan extension $p_!\QFD$ of $\QFD$ along $p\op$ is given by the explicit formula
\[
(p_!\QFD)(p(Y))\simeq  \colim_{[\beta\colon Y \to Z] \in (\ker(p)_{Y/})\op} \QFD(\fib(\bet)),
\]
see Remark~\reftwo{remark:kan-along-verdier-projection}.
\end{remark}

\begin{example}
\label{example:kan-poincare}%
If \(p\colon \D \to \E\) is a Verdier projection and \(\QFD\) is a Poincaré structure on \(\D\), then the hermitian structure \(p_!\QFD\) on \(\E\) and the canonical hermitian refinement of \(p\) are Poincaré if and only if \(\ker(p)\) is invariant under the duality, and in this case
\[
(\ker(p),\QFD) \longrightarrow (\D,\QFD) \longrightarrow (\E,p_!\QFD).
\]
is a Poincaré-Verdier sequence.

Indeed, if \(p_!\QFD\) and \(p\) are Poincaré, then it is immediate that \(\ker(p)\) is closed under the duality. Conversely if \(\ker(p)\) is closed under the duality, then, since the forgetful functor \(\Catp \to \Catx\) preserves colimits, the cofibre of the inclusion \((\ker(p),\QFD) \to (\D,\QFD)\) in \(\Catp\) must be equivalent to a Poincaré \(\infty\)-category of the form \((\E,\QFE)\) for some Poincaré structure on \(\E\) equipped with a Poincaré functor \((p,\vartheta)\colon (\D,\QFD) \to (\E,\QFE)\). The latter is then a Poincaré-Verdier projection by construction, and by Proposition~\reftwo{proposition:criterion-poincare} \reftwoitem{item:cofiberCatp}
the natural transformation \(p_!\QFD \to \QFE\) determined by \(\vartheta\) must be an equivalence, and so the desired properties of \(p_!\QFD\) follow.
\end{example}

Poincaré-Verdier sequences can be assembled into an \(\infty\)-category in various ways by letting the morphisms be some form of diagrams
\[
\begin{tikzcd}
(\C,\QF) \ar[r]\ar[d] & (\D,\QFD) \ar[r]\ar[d] & (\E,\QFE) \ar[d] \\
(\C',\QF')\ar[r] & (\D',\QFD')\ar[r] & (\E',\QFE')
\end{tikzcd}
\]
with horizontal Verdier sequences. While one obvious choice is requiring the right hand square in the diagram cartesian in $\Catp$, one obtains a well-behaved theory already by requiring only adjointability: Recall from Definition~\reftwo{definition:adjointable} that we call a square of stable \(\infty\)-categories \defi{adjointable} if it becomes both vertically and horizontally right adjointable after inductive completion, or equivalently (by Lemma~\reftwo{lemma:adjointable-equivalent}) horizontally and vertically, respectively, left adjointable after projective completion. We will say that a square in $\Catp$ is adjointable if the underlying square of stable \(\infty\)-categories is so; on account of the dualities  in $\Catp$ this underyling square is equivalent to the square of its opposites so any of the four adjointability conditions in fact implies the others in this case.
We then write 
\[
\PVer \subseteq \Fun(\Del^2,\Catp)
\]
for the (non-full) subcategory spanned by the Poincaré-Verdier sequences and maps as in the displayed diagram above in which the right hand square (or by Lemma~\reftwo{lemma:adj-verdier} equivalently the left hand one) is adjointable.

We then have the following relation between cartesian and adjointable squares,
 a Poincaré analogue of Proposition~\reftwo{proposition:verdier-square-char}.

\begin{proposition}
\label{proposition:poincare-verdier-square-char}%
For a commutative square
\[
\begin{tikzcd}
(\D,\QFD) \ar[r, "f"]\ar[d, "p"']& (\D',\QFD') \ar[d, "p'"] \\
(\E,\QFE) \ar[r, "f'"] & (\E',\QFE')
\end{tikzcd}
\]
in \(\Catp\) the following are equivalent:

\begin{enumerate}
\item
\label{item:square-cartesian}%
The square is cartesian and \(p'\) is a Poincaré-Verdier projection.
\item
\label{item:square-adjointable}%
The square is adjointable, both \(p\) and \(p'\) are Poincaré-Verdier projections and the induced functor \(\ker(p) \to \ker(p')\) is an equivalence in \(\Catp\).
\end{enumerate}

Similarly, the following are equivalent:
\begin{enumerate}
\item
\label{item:square-cocartesian}%
The square is cocartesian and \(f\) is a Poincaré-Verdier inclusion.
\item
\label{item:square-adjointable2}%
The square is adjointable, both \(f\) and \(f'\) are Poincaré-Verdier inclusions and the induced functor \(\coker(f) \to \coker(f')\) is an equivalence in \(\Catp\).
\end{enumerate}

\end{proposition}

\begin{corollary}
\label{corollary:PV-projpullbackstable}%
\label{corollary:PV-proj-pullback}%
The collection of (split) Poincaré-Verdier projections is closed under pullbacks and the collection of (split) Poincar\'e-Verdier inclusions is closed under pushouts.
\end{corollary}
\begin{proof}
The non-split case follows from Proposition~\reftwo{proposition:poincare-verdier-square-char}. To obtain the split case combine this with the fact that on the level of underlying stable \(\infty\)-categories split Verdier projections are closed under pullbacks, see Corollary~\reftwo{corollary:split-verdier-proj-pullback}.
\end{proof}

\begin{proof}[Proof of Proposition~\reftwo{proposition:poincare-verdier-square-char}]
We prove the first equivalence, the second one is obtained by dual arguments.
Let us denote the composite $p' \circ f \simeq f' \circ p$ by $h$ throughout. Assume \reftwoitem{item:square-cartesian} holds. Then clearly the induced functor \(\ker(p) \to \ker(p')\) is an equivalence in \(\Catp\). In addition, 
applying Proposition~\reftwo{proposition:verdier-square-char} on the level of underlying stable \(\infty\)-categories 
we also get that the square is adjointable and that the underlying exact functor of \(p\) is a Verdier projection. Now since we assume that the square is cartesian in \(\Catp\) the square of quadratic functors
\[
\begin{tikzcd}
\QFD \ar[r]\ar[d] & f^*\QFD'\ar[d] \\
p^*\QFE \ar[r] & h^*\QFE'
\end{tikzcd}
\]
is cartesian in \(\Funq(\D)\). By adjointability and Remark~\reftwo{remark:adj-beck-chevalley}, the maps \(p_!f^*\QFD'\Rightarrow (f')^*p'_!\QFD'\) and \(p_!h^*\QFE'\Rightarrow (f')^*p'_!(p')^*\QFE'\) are equivalences and since \(p'\) is assumed a Poincaré-Verdier projection we get from Corollary~\reftwo{corollary:poincare-projection} that the right vertical map is sent to an equivalence by \(p_!\). It follows that the left vertical map is also sent to an equivalence by \(p_!\), and so \(p\) is a Poincaré-Verdier projection by another application of Corollary~\reftwo{corollary:poincare-projection}.

Now assume that \reftwoitem{item:square-adjointable} holds. By Proposition~\reftwo{proposition:verdier-square-char} the square of underlying stable \(\infty\)-categories of the square of the proposition is cartesian, and so to show that it is also cartesian in \(\Catp\) it suffice to check that the square of hermitian structures is cartesian. Let \(\Lam\) be the total fibre of this square. Since the underlying square is cartesian in \(\Catx\) and all involved functors are duality preserving the square is also cartesian on the level of stable \(\infty\)-categories with duality. It follows that the square of the proposition is cartesian on the level of bilinear parts, and so the total fibre \(\Lam\) is exact. At the same time, by Remark~\reftwo{remark:adj-beck-chevalley} and the fact that the vertical maps are Poincaré-Verdier projections the vertical maps in this square are sent to equivalences by \(p_!\), and so \(p_!(\Lam) = 0\). By Corollary~\reftwo{corollary:criterion-projection} we get that \(\Lam\) is left Kan extended from its restriction to \(\ker(p)\). But the map $p^*\Psi \rightarrow h^*\Psi'$ restricts to \(0\) on \(\ker(p)\) and the condition that the Poincaré functor \(\ker(p) \to \ker(p')\) is an equivalence implies that $\Phi \rightarrow f^*\Phi'$ restricts to an equivalence. We conclude that \(\Lam|_{\ker(p)} = 0\) and hence \(\Lam =0\), as desired.
\end{proof}

\subsection{Split Poincaré-Verdier sequences and Poincaré recollements}
\label{subsection:poincare-recollements}%

We now turn to split Poincaré-Verdier sequences, which are by definition Poincaré-Verdier sequences in which the underlying Verdier sequence is split. Let us therefore mention from Lemma~\reftwo{lemma:quasi-split} that a sequence
\[
\C \xrightarrow{f} \D \xrightarrow{p} \E
\]
in \(\Catx\) with vanishing composite is a split Verdier sequence if and only if it is a fibre sequence and \(p\) admits fully faithful left and right adjoints, if and only if it is a cofibre sequence and \(f\) is fully faithful and admits left and right adjoints. Furthermore, this notion is equivalent to that of a stable recollement (see Proposition~\reftwo{proposition:criterion-split}).

\begin{observation}
\label{observation:adjoints-equivalent}%
The underlying functor \(p\) of a Poincaré functor admits a left adjoint if and only if it admits a right adjoint. Indeed, a left or right adjoint to \(p\) gives a right or left adjoint to \(p\op\), respectively, but \(p\) and \(p\op\) are naturally equivalent by means of the dualities in source and target.
\end{observation}

With this at hand, we derive the following criterion to recognise split Poincaré-Verdier sequences.
 
\begin{proposition}
\label{proposition:criterion-poincare-split}%
Let 
\[
(\C,\QF) \xrightarrow{(f,\eta)} (\D,\QFD) \xrightarrow{(p,\vartheta)} (\E,\QFE)
\]
be a sequence of functors in \(\Catp\) with vanishing composite. Then the following holds:
\begin{enumerate}
\item 
\label{item:fiber-left-adjoint}%
If it is a fibre sequence in \(\Catp\), then it is a split Poincaré-Verdier sequence if and only if \(p\) admits a fully faithful left adjoint \(g\) and the transformation 
\[
g^*\QFD  \stackrel{g^*\vartheta}{\Longrightarrow} g^*p^*\QFE \stackrel{u^*}{\Longrightarrow} \QFE
\]
is an equivalence, where \(u \colon \id_{\C} \Rightarrow pg\) denotes an adjunction unit.
\item
\label{item:cofiber-left-adjoint}%
If it is a cofibre sequence in \(\Catp\), then it is a split Poincaré-Verdier sequence if and only if \(f\) is fully faithful, admits a right adjoint, and \(\eta \colon \QF \rightarrow f^*\QFD\) is an equivalence.
\end{enumerate}
\end{proposition}

\begin{proof}
Assume that the sequence from the statement is a fibre sequence in \(\Catp\).  Then its image in \(\Catx\) is a fibre sequence as well. 
By the previous observation, the existence of a left adjoint to \(p\) implies that of a right adjoint, 
so the underlying sequence of stable \(\infty\)-categories is a split Verdier-sequence if and only if \(p\) admits a fully faithful left adjoint \(g\colon \E \to \D\). In this case, it follows from Remark~\reftwo{remark:compute} that \(g^*\QFD\) is a left Kan extension of \(\QFD\), and the transformation from the statement is the extension of \(\vartheta\). Thus, \(\QFE\) is a left Kan extension of \(\QFD\) if and only if this extension is an equivalence, which gives the claim by Proposition~\reftwo{proposition:criterion-poincare}.

The second item is immediate from Observation~\reftwo{observation:adjoints-equivalent} and Proposition~\reftwo{proposition:criterion-poincare} \reftwoitem{item:fiberCatp}.
\end{proof}

\begin{corollary}
\label{corollary:split-poincare-projection-inclusion}%
\ %
\begin{enumerate}
\item
\label{item:split-poinc-verdier-incl-concrete}%
A Poincaré functor \((f,\eta)\colon (\C,\QF) \to (\D,\QFD)\) is a split Poincaré-Verdier inclusion if and only if \(f\) is fully faithful, admits a right adjoint, and the map \(\eta \colon \QF \rightarrow f^*\QFD\) is an equivalence. 
\item
\label{item:split-poinc-verdier-proj-concrete}%
A Poincaré functor \((p,\vartheta)\colon (\D,\QFD)\to (\E,\QFE)\) is a split Poincaré-Verdier projection if and only if 
\(p\) admits a fully faithful left adjoint \(g\) and the composite transformation 
\(
\begin{tikzcd}
g^*\QFD  \ar[r,Rightarrow,"g^*\vartheta"] & g^*p^*\QFE \ar[r,Rightarrow,"\sim"] & \QFE
\end{tikzcd}
\)
is an equivalence.
\end{enumerate}
\end{corollary}

\begin{remark}
\label{remark:hermitian-adjoints}%
By means of the equivalence \(g^*\QFD \simeq \QFE\), the left adjoint \(g\) to a Poincaré-Verdier projection \(p \colon (\D,\QFD) \rightarrow (\E,\QFE)\) automatically becomes an hermitian functor \((\E,\QFE) \rightarrow (\D,\QFD)\) (which is usually not Poincaré). One readily checks that the unit gives an equivalence of hermitian functors \(\id_{(\E,\QFE)} \Rightarrow pg\), making \(g\) a section of \(p\) in \(\Cath\). 

In fact, granting that the \(\infty\)-categories \(\catforms(\Funx((\D,\QFD),(\E,\QFE)))\) provide a \(\Cat\)-enrichment to \(\Cath\) (a fact we will neither prove nor even make precise here), the adjunction between \(g\) and \(p\) is an enriched one, i.e.\ its unit \(\id_\E \Rightarrow gp\) and counit \(pg \Rightarrow \id_\D\) canonically promote to objects in \(\catforms(\Funx((\E,\QFE),(\E,\QFE)))\) and \(\catforms(\Funx((\D,\QFD),(\D,\QFD)))\), such that the triangle identities hold in  these \(\infty\)-categories. 

Conversely, the existence of such an enriched left adjoint to \(p\), whose unit is an equivalence, is readily checked to amount precisely to the conditions of Corollary~\reftwo{corollary:split-poincare-projection-inclusion} \reftwoitem{item:split-poinc-verdier-proj-concrete}.

Similarly, the existence of an enriched right adjoint with counit an equivalence boils down to precisely the conditions in \reftwoitem{item:split-poinc-verdier-incl-concrete} above, and therefore detects split Poincaré-Verdier inclusions; in particular, the counit always provides the right adjoint to a Poincaré-Verdier inclusion with an hermitian structure (which is again usually not Poincaré).

We warn the reader that the analogous statements involving the right adjoint to a Poincaré-Verdier projection and the left adjoint to a Poincaré-Verdier inclusion fail entirely; for instance in the metabolic Poincaré-Verdier sequence of Example~\reftwo{example:metabolicfseq} below, the only hermitian refinement of the right adjoint to the projection is null, and so certainly does not give rise to a splitting of \(p\). 
\end{remark}

\begin{example}
\label{example:metabolicfseq}%
For any Poincar\'e $\infty$-category $(\C,\QF)$ there are two natural Poincar\'e structures
\[
\QF_{\arr}(X \to Y) = \QF(X) \times_{\Bil_{\QF}(X,X)} \Bil_{\QF}(X,Y)\quad \text{and} \quad \QF_{\met}(X \rightarrow Y) = \fib(\QF(Y) \rightarrow \QF(X))
\]
on the category $\Ar(\C)$ resulting in Poincar\'e categories $\Ar(\C,\QF)$ and $\Met(\C,\QF)$, together with an equivalence $\cof \colon \Ar(\C,\QF)\qshift{-1} \rightarrow \Met(\C,\QF)$. We then have tautological Poincar\'e functors %
\[
\const_{(-)} \colon (\C,\QF) \lrar \Ar(\C,\QF) \quad\text{and}\quad \met\colon \Met(\C,\QF) \lrar (\C,\QF)
\]
whose underlying exact functors are given by \(X \mapsto \id_X\) and \([W \to X] \mapsto X\); indeed, in terms of the identifications of $\Met(\C,\QF)$ and $\Ar(\C,\QF)$ in terms of the pairing categories in \S\refone{subsection:thom}, they are the adjunction unit and counit of the latter, respectively. The left hand map is an example of a split Poincaré-Verdier inclusion and the right hand one an example of a split Poincaré-Verdier projection: The left and right adjoints of \(\const_{(-)}\) send an arrow \(X \to Y\) to its target and source, respectively, and the left and right adjoints of \(\met\) send \(X \in \C\) to \(0 \to X\) and \(X \to X\), respectively. Inspecting the relevant definitions we see that the structure transformations  
\[
\QF(X) \to \QF_{\arr}(X \xrightarrow{\id} X) = \QF(X) \times_{\Bil_{\QF}(X,X)} \Bil_{\QF}(X,X) \quad\text{and}\quad \QF(X) \to \QF_{\met}(0 \to X) = \fib(\QF(X) \to \QF(0))
\]
are indeed equivalences so Proposition~\reftwo{proposition:criterion-poincare-split} applies. The resulting split Poincaré-Verdier sequences are
\[
(\C,\QF) \lrar \Ar(\C,\QF) \xrightarrow{\cof} (\C,\QF\qshift{1})
\quad \text{or equivalently} \quad 
(\C,\QF\qshift{-1}) \lrar \Met(\C,\QF) \xrightarrow{\met} (\C,\QF).
\]
up to the indicated shift. We will use the term \emph{metabolic sequence} to refer to either one of them, and they will play a crucial role throughout the paper. %
\end{example}

Let $\PVer^{\perp} \subset \PVer$ denote the full subcategory spanned by the split Verdier sequences. Let us also remind the reader that we already recorded in Corollary \reftwo{corollary:PV-projpullbackstable} that split Poincar\'e-Verdier sequences are closed under pullbacks. 

Recall then from \S\reftwo{section:appendix-split-verdier} that for a Verdier inclusion \(f \colon\C \to \D\) we denote by \(\D^{\spl} \subseteq \D\) the full subcategory of \(\D\) spanned by those objects whose images under the right and left adjoints of \(f\) on the inductive or projective completions, respectively, lie in \(\C\). We call $\D^{\spl}$ the split core of $f$. By Lemma~\reftwo{lemma:split-core} the induced functor \(\C \to \D^{\spl}\) is a split Verdier inclusion and is universally characterised as being the terminal split Verdier inclusion under \(f\colon\C \to \D\) (where morphisms of Verdier inclusions are given by adjointable squares). We now adapt the split core construction to the Poincaré setting: given a Verdier inclusion \((\C,\QF) \to (\D,\QFD)\), we write \((\D,\QFD)^{\spl} := (\D^{\spl},\QFD|_{\D^{\spl}})\). The inclusion \(\D^{\spl} \subseteq \D\) is stable under the duality of \(\D\) by Observation~\reftwo{observation:adjoints-equivalent} and so \((\D,\QFD)^{\spl}\) is Poincaré, and the induced Poincaré functor \((\C,\QF) \to (\D,\QFD)^{\spl}\) is a split Poincaré-Verdier inclusion by Corollary~\reftwo{corollary:split-poincare-projection-inclusion}. We then also write \((\E,\QFE)^{\qspl} \coloneq (\D,\QFD)^{\spl}/(\C,\QF)\).

To state the universal property of the split core in full, define finally for two Poincar\'e-Verdier sequences $(\C,\QF) \rightarrow (\D,\QFD) \rightarrow (\E,\QFE)$ and $(\C',\QF') \rightarrow (\D',\QFD') \rightarrow (\E',\QFE')$ a new Poincar\'e $\infty$-category $\Fun^{\Ver}((\D,\QFD)\rightarrow (\E,\QFE),(\D',\QFD') \rightarrow (\E',\QFE))$ as the full subcategory of the pullback of
\[
\Funx((\D,\QFD),(\D',\QFD')) \longrightarrow \Funx((\C,\QF),(\D',\QFD'))\longleftarrow\Funx((\C,\QF),(\C',\QF'))
\]
spanned by those functors that give rise to adjointable squares; it is closed under the duality on the pullback since conjugation by the dualities in source and target swaps the (inductive/projective) left and right adjoints. By Corollary \refone{corollary:forms-in-functor-cats} we then find
\[
\Poinc(\Fun^{\Ver}((\D,\QFD)\rightarrow (\E,\QFE),(\D',\QFD') \rightarrow (\E',\QFE))) \simeq \Hom_{\PVer}((\D,\QFD)\rightarrow (\E,\QFE),(\D',\QFD') \rightarrow (\E',\QFE)).
\]

\begin{lemma}
\label{lemma:split-poincare-core}%
Let \((\C,\QF) \to (\D,\QFD) \to (\E,\QFE)\) 
be a Poincaré-Verdier sequences. Then post-composition with the inclusion of the split core
induces an equivalence 
\[
\Fun^{\Ver}((\D',\QFD')\to(\E',\QFE'),(\D,\QFD)^{\spl}\to(\E,\QFE)^{\qspl}) \longrightarrow \Fun^{\Ver}((\D',\QFD')\to(\E',\QFE'),(\D,\QFD)\to(\E,\QFE))
\]
for any split Poincaré-Verdier sequence \((\C',\QF')\to(\D',\QFD')\to(\E',\QFE')\). In particular, the formation of split cores assembles 
into a right adjoint 
to the inclusion \(\PVer^{\perp} \subseteq \PVer\).
\end{lemma}
\begin{proof}
By Lemma~\reftwo{lemma:adj-verdier} it suffices to note that in any adjointable square
\[
\begin{tikzcd}
\C' \ar[r]\ar[d] & \C\ar[d] \\
\D'\ar[r, "f"] & \D \ ,
\end{tikzcd}
\]
with \((\C',\QF') \to (\D',\QFD')\) a split Poincaré-Verdier inclusion the underlying exact functor of \(f\) must factor through \(\D^{\spl}\) by adjointability and this factorisation refines to an hermitian functor in an essentially unique manner since the Poincaré structure on \(\D^{\spl}\) is restricted from that of \(\D\). Finally, the resulting factorisation of \(f\) through \((\D,\QFD)^{\spl}\) is automatically duality preserving since the inclusion \(\D^{\spl} \subseteq \D\) is stable under the duality. 
The claim about the adjunction follows by taking Poincar\'e objects and applying Corollary \refone{corollary:forms-in-functor-cats}.
\end{proof}

\subsection{Poincaré-Karoubi sequences}
\label{subsection:poincare-karoubi-sequences}%

In this section, we study Poincaré-Karoubi sequences, the analogues of Poincaré-Verdier sequences in the setting of idempotent complete Poincaré \(\infty\)-categories. 
On the one hand, these are important in their own right when considering the hermitian analogue of non-connective \(\K\)-theory (we will do this in \paperfour), on the other, it is often easier to establish Poincaré-Verdier sequences in a two-step process: First one constructs a Poincaré-Karoubi sequence using the Thomason-Neeman localisation theorem \reftwo{theorem:indkaroubi}, or in modern guise, the equivalence between small stable \(\infty\)-categories, and compactly generated stable \(\infty\)-categories, and then in a second step isolates subcategories forming Poincaré-Verdier sequences; see Proposition~\reftwo{proposition:module-Verdier-example} for an example.

Let us establish some terminology: Recall that we denote by \(\C^\natural\) the idempotent completion of a small \(\infty\)-category \(\C\) and refer the reader to \cite{HTT}*{\S 5.1.4} for its construction. The \(\infty\)-category \(\C^\natural\) is stable if \(\C\) is, and the natural functor \(i \colon \C \rightarrow \C^\natural\) is fully faithful, exact and has dense essential image, 
where a full subcategory \(\D \subseteq \C\) is called \defi{dense} if every object of \(\C\) is a retract of one in \(\D\). A stable $\infty$-category $\C$ always contains a minimal stable dense subcategory $\C^\mathrm{min}$ spanned by all $X \in \C$ with $[X] = 0 \in \K_0(\C)$, see \reftwo{proposition:minimal-cats}, we say that $\C$ is minimal if $\C^\mathrm{min} = \C$. Recall also that we call an exact functor a \defi{Karoubi equivalence} if it is fully faithful with dense essential image, in other words if it induces an equivalence on the idempotent completions or minimalisations (cf.\ Definition~\reftwo{definition:karoubi-equivalence}).

\begin{definition}
\label{definition:poincare-karoubi-equivalence}%
A Poincaré \(\infty\)-category is \defi{idempotent complete} or \defi{minimal} if its underlying stable \(\infty\)-category is. We denote by \(\Catpi \subseteq \Catp\) the full subcategory spanned by the idempotent complete Poincaré 
\(\infty\)-categories. 
A Poincaré functor \((f, \eta)\colon (\C,\QF) \rightarrow (\D,\QFD)\) is a \defi{Karoubi equivalence} if \(f\) is a Karoubi equivalence and \(\eta\colon \QF\to f^*\QFD\) is an equivalence.
\end{definition}

\begin{proposition}
\label{proposition:idempotent-completion}%
Let \((\C,\QF)\) be a Poincaré \(\infty\)-category, \(i\colon\C\to \C^\natural\) its idempotent completion, and $j \colon \C^\mathrm{min} \rightarrow \C$ its minimalisation. Then \(i_!\QF\) and $j^*\QF$ are Poincaré structures on \(\C^\natural\) and $\C^\mathrm{min}$, respectively, and the canonical Poincar\'e functors \((\C,\QF)\to (\C^\natural,i_!\QF)\) and $(\C^\mathrm{min},j^*\QF) \to (\C,\QF)$ are Karoubi equivalences.

Moreover, for an idempotent complete Poincaré \(\infty\)-category \((\D,\QFD)\) and a minimal Poincar\'e $\infty$-category $(\E,\QFE)$ the functors
\[
\Funx((\C^\natural,i_!\QF),(\D,\QFD))\lrar \Funx((\C,\QF),(\D,\QFD)) \quad \text{and} \quad \Funx((\E,\QFE),(\C^\mathrm{min},j^\ast\QF)) \lrar \Funx((\E,\QFE),(\C,\QF))
\]
are equivalence of Poincaré \(\infty\)-categories. 
\end{proposition}

\begin{proof}
By Lemma~\refone{lemma:kan-extension-exact-quadratic} and Proposition~\refone{proposition:left-kan-bilinear-linear},
the functor \(i_!\QF\) is quadratic with bilinear part \((i\times i)_!\Bil_\QF\). To see that this is perfect, note first that it restricts back to \(\Bil_\QF\) since \(i\) is fully faithful. Now, the idempotent completion of the equivalence \(\Dual_\QF \colon \C\op \rightarrow \C\) is another equivalence \(\Dual \colon (\C^\natural)\op \simeq (\C\op)^\natural \rightarrow \C^\natural\), and by the previous observation, the functors $\Hom_{\C^\natural}(-,\Dual-) \quad \text{and} \quad \Bil_{i_!\QF}$ agree on \(\C\op \times \C\op\), and therefore on all of \((\C^\natural)\op \times (\C^\natural)\op\) by \cite{HTT}*{Proposition 5.1.4.9}. This shows that both \(i_!\QF\) and \(i\) are Poincaré.

To see the second claim, let us fix \((\D,\QFD)\) an idempotent-complete Poincaré \(\infty\)-category and consider the Poincaré functor
\[
i^* \colon \Funx((\C^\natural,i_!\QF),(\D,\QFD))\lrar \Funx((\C,\QF),(\D,\QFD))\,.
\]
By Proposition~\reftwo{proposition:minimal-cats}
this is an equivalence of the underlying stable \(\infty\)-categories, so it suffices to show that it induces also an equivalence on the corresponding quadratic functors. But for an exact functor \(f:\C^\natural\to \D\), this map is precisely the canonical equivalence
\(\nat(i_!\QF,f^*\QFD)\to \nat(\QF,i^*f^*\QFD)\).

The claims about minimalisation are evident since the duality of $\QF$ preserves $\C^\mathrm{min}$.
\end{proof}

\begin{remark}
\label{remark:idem-poincare}%
For $i\colon \C \to \C'$ a Karoubi equivalence, the adjunction \(i_!\colon \Funq(\C) \adj \Funq(\C')\cocolon i^*\) between hermitian structures on \(\C\) and hermitian structures on \(\C'\) is an equivalence, since \(i_!\) is fully faithful and \(i^*\) is conservative by the density of \(i\). By Proposition~\reftwo{proposition:idempotent-completion}, this equivalence restricts to an equivalence between Poincaré structures on \(\C\) and Poincaré structures on \(\C'\) whose associated duality preserves \(\C\).
\end{remark}

\begin{corollary}
\label{proposition:minimal-poincare}%
The localisation of \(\Catp\) at the Karoubi equivalences admits both a left and a right adjoint, the right adjoint is induced by \((\C,\QF) \mapsto (\C^\natural,i_!\QF)\), and the left adjoint by \((\C,\QF) \mapsto (\C^\mn,j^*\QF)\).
Consequently, regarding the idempotent completion as a functor \((-)^{\natural}\colon\Catp \to \Catpi\) it has both adjoints and thus preserves both limits and colimits.
\end{corollary}

We will write \((\C,\QF)^{\natural}\) and $(\C,\QF)^\mn$ for the idempotent completion and minimalisation of a Poincar\'e $\infty$-category $(\C,\QF)$. The analogous statement for the underlying stable \(\infty\)-categories is Proposition~\reftwo{proposition:minimal-cats}.

\begin{proof}
According to Lemma~\reftwo{lemma:bousfieldgenI}, this follows from Proposition \reftwo{proposition:idempotent-completion} upon investing that Poincar\'e objects in functor categories are Poincar\'e functors by Corollary \refone{corollary:forms-in-functor-cats}.
\end{proof}

\begin{corollary}
\label{corollary:poinc-karoubi-equiv-pullback}%
The collection of Karoubi equivalences is closed under pullbacks and pushouts in $\Catp$. %
\end{corollary}

Recall that a sequence \(\C \xrightarrow{f} \D \xrightarrow{p} \E\) of exact functors with vanishing composite is a \emph{Karoubi sequence} (Definition~\reftwo{definition:kar-seq}) if the sequence 
$
\C^\natural \to \D^\natural \rightarrow \E^\natural
$
is both a fibre and a cofibre sequence in \(\Catxi\). In this case, we refer to \(f\) as a \emph{Karoubi inclusion} and to \(p\) as a \emph{Karoubi projection}. 
In the same spirit, we put:
\begin{definition}
\label{definition:poinc-kar-seq}%
A sequence 
\[
(\C,\QF) \xrightarrow{(f,\eta)} (\D, \QFD) \xrightarrow{(p,\theta)} (\E, \QFE)
\]
of Poincaré functors with vanishing composite is a \defi{Poincaré-Karoubi sequence} if 
\[
(\C,\QF)^{\natural} \xrightarrow{(f,\eta)^{\natural}} (\D, \QFD)^{\natural} \xrightarrow{(p,\vartheta)^{\natural}} (\E, \QFE)^{\natural}
\]
is both a fibre sequence and a cofibre sequence in \(\Catpi\). We then call \((f,\eta)\) a \defi{Poincaré-Karoubi inclusion} and \((p,\vartheta)\) a \defi{Poincaré-Karoubi projection}.
\end{definition}

We warn the reader that, contrary to the situation for (Poincaré-)Verdier sequences, a (Poincaré-)Karoubi sequence is determined by its inclusion or its projection only up to idempotent completion of the third term.

We record a few simple consequences of the definition.

\begin{observation}
Since the forgetful and hyperbolic functors both preserve limits, colimits and Karoubi equivalences, it follows from Propositions~\reftwo{proposition:minimal-poincare} that they induce limit and colimit preserving functors between \(\Catpi\) and \(\Catxi\). 
It then follows that the sequence of stable \(\infty\)-categories underlying a Poincaré-Karoubi sequence is a Karoubi sequence and the hyperbolisation of a Karoubi sequence is a Poincaré-Karoubi sequence.
\end{observation}

By Proposition~\reftwo{proposition:criterion-karoubi}, any Verdier sequence is a Karoubi sequence. Analogously: 
\begin{proposition}
\label{proposition:poincare-verdier-is-karoubi}%
Every Poincaré-Verdier sequence is a Poincaré-Karoubi sequence.
\end{proposition}
\begin{proof}
A bifibre sequence in \(\Catp\) remains so in \(\Catpi\) after applying $(-)^\natural$ by Proposition~\reftwo{proposition:minimal-poincare}.
\end{proof}

\begin{proposition}
\label{proposition:poincare-karoubi-underlying}%
Let 
\[
(\C,\QF) \xrightarrow{(f,\eta)} (\D,\QFD) \xrightarrow{(p,\vartheta)} (\E,\QFE)
\]
be a sequence of Poincaré functors with vanishing composite. Then:
\begin{enumerate}
\item 
\label{item:poincare-karoubi-fiber}%
Its idempotent completion is a fibre sequence in \(\Catpi\) if and only if the idempotent completion of its underlying sequence is a fibre sequence in \(\Catxi\)  
and \(\eta\) induces an equivalence \(\QF \Rightarrow f^*\QFD\).
\item 
\label{item:poincare-karoubi-cofiber}%
Its idempotent completion is a cofibre sequence in \(\Catpi\) if and only if 
the idempotent completion of its underlying sequence is a cofibre sequence in \(\Catxi\) 
and \(\vartheta\) exhibits \(\QFE\) as the left Kan extension of \(\QFD\) along \(p\op\).
\end{enumerate}
In particular, it is a Poincaré-Karoubi sequence if and only if
its underlying sequence is a Karoubi sequence and
both \(\eta\) induces an equivalence \(\QF \Rightarrow f^*\QFD\)
and \(\vartheta\) exhibits \(\QFE\) as the left Kan extension of \(\QFD\) along \(p\).
\end{proposition}
\begin{proof}
By Proposition~\reftwo{proposition:minimal-poincare}, fibres in \(\Catpi\) are computed in \(\Catp\) while cofibres are computed as idempotent completions of cofibres in \(\Catp\). Similarly, by Proposition~\reftwo{proposition:minimal-cats} the analogous statement holds on the level of underlying stable \(\infty\)-categories. In addition, idempotent completion of Poincaré \(\infty\)-categories is compatible with idempotent completion on the underlying stable \(\infty\)-categories by Proposition~\reftwo{proposition:idempotent-completion}.
Thus, \reftwoitem{item:poincare-karoubi-fiber} and \reftwoitem{item:poincare-karoubi-cofiber} follow from Proposition~\reftwo{proposition:criterion-poincare} \reftwoitem{item:fiberCatp} and \reftwoitem{item:cofiberCatp}, respectively, using the equivalence between quadratic functors on \(\C\) and on \(\C^\natural\) explained in Remark~\reftwo{remark:idem-poincare} (as well as the ones for \(\D\) and \(\E\)). 
\end{proof}

Combining Proposition~\reftwo{proposition:poincare-karoubi-underlying}
and Corollary~\reftwo{corollary:characterisation-of-Karoubi-inclusions-and-projections} we obtain: 
\begin{corollary}
\label{corollary:criterion-poincare-karoubi-projection}%
\ %
\begin{enumerate}
\item
A Poincaré functor \((f,\eta)\colon (\C,\QF) \to (\D,\QFD)\) is a Poincaré-Karoubi inclusion if and only if \(f\) is fully faithful and the map \(\eta \colon \QF \Rightarrow f^*\QFD\) is an equivalence. 
\item
A Poincaré functor \((p,\vartheta)\colon (\D,\QFD) \to (\E,\QFE)\) is a Poincaré-Karoubi projection if and only if \(p\) has dense essential image, the induced functor \(\D \to p(\D)\) is a localisation and \(\vartheta\colon \QFD \Rightarrow p^*\QFE\) exhibits \(\QFE\) as the left Kan extension of \(\QFD\) along \(p\op\).
\end{enumerate}
\end{corollary}

Combining Corollary~\reftwo{corollary:criterion-poincare-karoubi-projection} with
Corollary~\reftwo{corollary:poincare-inclusion} and Corollary~\reftwo{corollary:poincare-projection} we obtain:

\begin{corollary}
\label{corollary:karoubi-composite}%
A Poincaré functor is a Poincaré-Karoubi inclusion if and only if it is a composite of a Karoubi equivalence followed by a Poincaré-Verdier inclusion. A Poincaré functor is a Poincaré-Karoubi projection if and only if it is a composite of a Poincaré-Verdier projection followed by a Karoubi equivalence.
\end{corollary}

\begin{corollary}
\label{corollary:criterion-poincare-karoubi-is-verdier}%
A Poincaré-Karoubi inclusion is Poincaré-Verdier if and only if the essential image of its underlying exact functor is closed under retracts. A Poincaré-Karoubi projection is Poincaré-Verdier if and only if its underlying exact functor is essentially surjective.
\end{corollary}

\begin{corollary}
\label{corollary:PK-projpullbackstable}%
The collection of Poincaré-Karoubi projections is closed under pullbacks and the collection of Poincar\'e-Karoubi inclusions is closed under pushouts.
\end{corollary}
\begin{proof}
This follows from Corollary~\reftwo{corollary:karoubi-composite}, Corollary~\reftwo{corollary:PV-proj-pullback}, and Corollary~\reftwo{corollary:poinc-karoubi-equiv-pullback}.
\end{proof}

\subsection{Poincaré-Verdier sequences among $\infty$-categories of modules} 
\label{subsection:examples-poincare-verdier}%

In this section, we will work out the important example of Poincar\'e-Verdier and Poincar\'e-Karoubi sequences of module \(\infty\)-categories in detail, where the hermitian structure is defined by means of a module with genuine involution, as introduced in \S\refone{subsection:genuine-modules} (see also the recollection section). We do so first in the generality of a map of \(\Eone\)-algebras \(\phi \colon A \rightarrow B\) over some base \(\Einf\)-ring \(k\) together with a map \(\eta\) of modules with genuine involution \(\phi_!(M,N,\alpha) \rightarrow (M',N',\beta)\) over \(B\), and eventually specialise to Ore localisations of discrete rings (and thus in particular localisations of commutative rings) with skew-involution in Corollary~\reftwo{corollary:ore-localisation}. 
The reader only interested in this case is invited to take \(k\) to be the (Eilenberg-Mac Lane spectrum of the) integers and \(A\) and \(B\) (and even \(M\) and \(M'\)) discrete from the start, though this does not simplify the discussion. Furthermore, it is important to allow \(N\) and \(N'\) to be non-discrete, so as to capture the genuine Poincaré structures. 

\begin{notation}
Throughout, unmarked tensor products are always over \(k\), and in case \(k = \GEM R\) translate to the derived tensor product \(\otimes_R^\mathbb L\). 
\end{notation}

We want to establish general conditions on \(\phi\) under which the hermitian functor \((\phi_!,\eta)\) becomes a Poincaré-Verdier or Poincaré-Karoubi projection. To obtain a 
Karoubi sequence on the underlying stable \(\infty\)-categories, the following conditions are necessary and sufficient: A map \(\phi\colon A \to B\) of \(\Eone\)-ring spectra is said to be a \defi{\solid} if the map
\[
B \otimes_A B \lrar B,
\]
induced by the multiplication of \(B\) is an equivalence of spectra. 
In the case of discrete rings, we call a map \(A \rightarrow B\) a \defi{derived localisation} if \(\GEM A \rightarrow \GEM B\) is a \solid in the sense above. 
For such a \solid of ring spectra denote by \(I \in \Mod(A)\) its fibre. Straight from the definition one finds that \(I\) belongs to \(\Mod(A)_B\), the kernel of \(\phi_! \colon \Mod(A) \rightarrow \Mod(B)\). We say that \(\phi\) has \defi{perfectly generated fibre} if \(I\) belongs to the smallest full subcategory of \(\Mod(A)_B\) which contains \(\Modp{A} \cap \Mod(A)_B\) and is closed under colimits. 
By Proposition~\reftwo{proposition:app-verdier-free} we then have that if \(\phi\colon A \to B\) is \solid of \(\Eone\)-rings with perfectly generated fibre, then for any subgroup \(\mathrm c \subseteq \K_0(A)\) the extension of scalars functors
\[
\phi_!^\cp \colon \Modp{A} \lrar \Modp{B}
\quad
\text{and}
\quad
\phi_!^\mathrm c\colon\Modc{\mathrm c}{A} \lrar \Modc{\phi(\mathrm c)}{B}
\]
are Karoubi and Verdier projections, respectively; here \(\Modc{\mathrm c}{A}\) denotes the full subcategory of \(\Modp{A}\) spanned by all those \(A\)-modules with \([A] \in \mathrm c \subseteq \K_0(A)\).

\begin{remark}
We warn the reader  
that for discrete rings, the notion of derived localisation as above 
differs from that of localisation as used in ordinary algebra: If \(A \rightarrow B\) is a localisation of discrete rings, then it is a derived localisation 
if and only if additionally \(\Tor_i^A(B,B) = 0\) for all \(i > 0\). This is automatic for commutative rings, or more generally if the localisation satisfies an Ore condition, but not true in general. Moreover, there are quotient maps \(A \rightarrow A/I\) which are derived localisations. 
\end{remark}
 
The following example essentially covers all of our immediate applications: 
\begin{example}
If \(A\) is an \(\Eone\)-ring spectrum and \(S \subseteq \pi_*A\) is a multiplicatively closed subset of homogeneous elements which satisfies the left or right Ore condition, 
then the map \(\phi \colon A \rightarrow A[S^{-1}]\) 
is a \solid with perfectly generated fibre by Example~\reftwo{example:ore}. Thus 
\(\phi_!^\cp \colon \Modp{A} \to \Modp{A[S^{-1}]} \)
is a Karoubi projection and 
\(\phi_!^\mathrm c \colon \Modc{\mathrm c}{A} \to \Modc{\im(\mathrm c)}{A[S^{-1}]}\)
is a Verdier projection for any \(\mathrm c \subseteq \K_0(A)\), see Corollary~\reftwo{corollary:ore-2}.
\end{example}

Let us now introduce hermitian structures into the picture. As discussed in section \S\refone{subsection:genuine-modules} and the last part of the recollection section, an invertible module with genuine involution \((M,N,\alpha \colon N \rightarrow M^\tC)\) over \(A\) gives rise to a Poincaré structure \(\QF^{\alpha}_M\) on \(\Modp A\); it restricts to a Poincaré structure on \(\Modc{\mathrm c}{A}\) provided that \(M\) belongs to \(\Modc{\mathrm c}{A}\) and provided \(\mathrm c\) is closed under the involution on \(\K_0(A)\) induced by \(M\). For example, if \(\mathrm c\) is the image of the canonical map \(\ZZ \rightarrow \K_0(A), 1 \mapsto A\), then \(\Modc{\mathrm c}{A} = \Modf{A}\) and this assumption is satisfied if also \(M \in \Modf{A}\). We computed the left Kan extension of this Poincaré structure along the functor \(\phi_!^\cp\colon\Modp{A} \to \Modp{B}\) in Corollary~\refone{corollary:induction-genuine-modules}: It is the hermitian structure associated to the module with genuine involution
\[
\phi_!(M,N,\alpha) = ((B \otimes B) \otimes_{A \otimes A} M, B \otimes_A N, \beta),
\]
over \(B\); here, \(\beta\) is the composition 
\[
B \otimes_A N \xrightarrow{\Delta \otimes \alpha} (B \otimes B)^\tC \otimes_A M^\tC \lrar ((B \otimes B) \otimes_{A \otimes A} M)^\tC
\]
where \(\Delta\) is the Tate diagonal. For example, by Remark~\reftwo{remark:compute}, the same formula then applies for the Kan extension along \(\phi_!^\mathrm c\colon\Modc{\mathrm{c}}{A} \to \Modc{\phi(\mathrm{c})}{B}\). 

In order to obtain Poincaré-Karoubi projections, we need a compatibility condition between \(\phi \colon A \rightarrow B\) and the module with involution \(M\) over \(A\):

\begin{definition}
\label{definition:module-compatible-with-induction}%
An invertible module with involution \(M\) over \(A\) is called \defi{\compatible} with a \solid of \(\Eone\)-rings \(A \rightarrow B\) if the composite
\[
B \otimes_A M \simeq (B \otimes A) \otimes_{A \otimes A} M \longrightarrow (B \otimes B) \otimes_{A \otimes A} M
\]
is an equivalence.
\end{definition}

\begin{examples}
\label{example:involution-induction-compatible}%
\ 
\begin{enumerate}
\item
\label{item:symmetric-modules-compatible}%
If \(A\) is an \(\Einf\)-ring and \(M\) an invertible \(A\)-module with \(A\)-linear involution (regarded as an \(A \otimes A\)-module via the multiplication map \(A \otimes A \rightarrow A\)), then compatibility is automatic, since in this case the map in question identifies with the evident one \(B \otimes_A M \rightarrow B \otimes_A B \otimes_A M\) which is an equivalence by the assumption that \(A \rightarrow B\) is a localisation.
\item If \(M\) is the module with involution over \(A\) associated to a Wall anti-structure \((\epsilon,\sigma)\) on a discrete ring \(A\) as in Example~\refone{example:wall-antistructure} (i.e., \(M=A\) regarded as an \(A \otimes A\)-module using the involution \(\sigma\), and then equipped with the involution \(\epsilon \sigma\), where \(\epsilon \in A^*\)) and 
\[
\phi \colon (A,\epsilon, \sigma) \longrightarrow (B,\delta,\tau)
\]
is a map of rings with anti-structure, then \(M\) is also automatically \compatible with \(\phi\) if the latter is a derived \solid: For in this case, it is readily checked that the maps 
\[
b \otimes b' \otimes a \longmapsto ba \otimes \tau(b') \quad \text{and} \quad b \otimes b' \longmapsto b \otimes b' \otimes 1
\]
give inverse equivalences
\[
B \otimes B \otimes_{A \otimes A} A \simeq B \otimes_A B,
\]
which translates the map in Definition~\reftwo{definition:module-compatible-with-induction} to the unit map \(B \rightarrow B \otimes_A B\) which is an equivalence since \(\phi\) is a derived \solid.
\item 
\label{item:ore-compatible}%
If \(\phi\) is an Ore localisation at the set \(S \subseteq \pi_*(A)\), and \(M\) is an invertible module with involution over \(A\), then \(M\) is \compatible with \(\phi\) if after inverting the action of \(S\) on \(M\) using the first \(A\)-module structure, \(S\) operates invertibly through the second one.
\item Combining the two previous examples, if \(M\) is the \(A\)-module associated to a Wall anti-structure \((\epsilon, \sigma)\) on \(A\), and \(S \subseteq A\) satisfies the Ore condition and is closed under the involution \(\sigma\), then \(M\) is \compatible with the localisation map \(A \rightarrow A[S^{-1}]\).
\end{enumerate}  
\end{examples}

\begin{proposition}
\label{proposition:module-Verdier-example}%
Let \(\phi\colon A \rightarrow B\) be a \solid of \(\Eone\)-ring spectra, with perfectly generated fibre and let
\((M,N,\alp)\) be an invertible module with genuine involution over \(A\), such that \(M\) is \compatible with \(\phi\).
Then \(\phi_!(M,N,\alpha)\) is invertible and 
\[
\phi_!^\cp\colon (\Modp{A},\QF^{\alp}_{M}) \lrar (\Modp{B},\QF^{\phi_!\alpha}_{\phi_!M}) \quad \text{and} \quad \phi_!^\mathrm c\colon (\Modc{\mathrm c}{A},\QF^{\alp}_{M}) \lrar (\Modc{\phi(\mathrm c)}{B},\QF^{\phi_!\alpha}_{\phi_!M}).
\]
are a Poincaré-Karoubi and Poincar\'e-Verdier projection, respectively, so long as \(\mathrm c \subseteq \K_0(A)\) is closed under the involution induced by \(M\).
\end{proposition}

\begin{proof}
The natural map 
\[
B \otimes_A \hom_A(X,M) \longrightarrow \hom_B(B \otimes_A X,B \otimes_A M)
\]
is an equivalence for \(X=A\) and thus for every compact \(A\)-module \(X\), in particular for \(X=M\), which shows that \(B \otimes_A M\) has \(B\) as its \(B\)-linear endomorphisms, and therefore by assumption \((B \otimes B) \otimes_{A \otimes A} M\) is invertible (or alternatively, one can apply Example~\reftwo{example:kan-poincare} together with Proposition~\refone{proposition:invertible-poincare}). Both functors \(\phi_!^\cp\) and \(\phi_!^\mathrm c\) are then Poincaré by Lemma~\refone{lemma:Poincresscalars}. 

By Corollary~\reftwo{corollary:poincare-projection}, the functor \(\phi_!^\mathrm c\) is a Poincaré-Verdier projection since the underlying map on module categories is a Verdier projection and by definition the Poincaré structure on the target is the left Kan extension of that on the source.
Similarly, the functor \(\phi_!^\cp\) is a Poincaré-Karoubi projection by Corollary~\reftwo{corollary:criterion-poincare-karoubi-projection}.
\end{proof}

\begin{corollary}
\label{corollary:quadratic-Verdier-sequence}%
Let \(\phi\colon A \rightarrow B\) be a \solid of \(\Eone\)-ring spectra, with perfectly generated fibre and let \(M\) be an invertible module with involution over \(A\) that is \compatible with \(\phi\).
Then 
\[
\phi_!^\cp\colon (\Modp{A},\QF^{\qdr}_{M}) \lrar (\Modp{B},\QF^{\qdr}_{\phi_!M}) 
\quad \text{and} \quad 
\phi_!^{\mathrm c}\colon (\Modc{\mathrm{c}}{A},\QF^{\qdr}_{M}) \lrar (\Modc{\phi(\mathrm{c})}{B},\QF^{\qdr}_{\phi_!M})
\]
are a Poincaré-Karoubi and Poincaré-Verdier projection (for \(\mathrm c \subseteq \K_0(A)\) closed under the duality), respectively. 
\end{corollary}

Symmetric Poincaré structures, however, are not generally preserved by left Kan extension: 

\begin{example}
\label{example:2-inverted-sphere}%
The map $p \colon \SS \rightarrow \SS[\tfrac 1 2]$ does \emph{not} induce an equivalence $p_!\QF_\SS^\sym \simeq \QF_{\SS[1/2]}^\sym$. Indeed, the linear part of the former is represented by \(\SS[\tfrac 1 2] \otimes \SS^\tC\) which by Lin's theorem is equivalent to \(\GEM\QQ_2\), whereas the linear part of the latter is represented by \(\SS[\tfrac 1 2]^\tC \simeq 0\). Consequently, the functor 
$
p_! \colon (\Modp{\SS},\QF^\sym) \to (\Modp{\SS[\frac 1 2]},\QF^\sym)
$
is not a Poincaré-Karoubi projection.
\end{example}

In the discrete case, an additional flatness assumption excludes such examples:

\begin{proposition}
\label{proposition:left-Kan-extension-specific-structures}%
Let \(\phi\colon A \rightarrow B\) be a derived localisation between discrete rings with perfectly generated fibre that furthermore makes \(B\) into a flat right module over \(A\) and let \(M\) be a discrete invertible module with involution over \(A\) that is \compatible with \(\phi\). Then for arbitrary \(m \in \ZZ \cup \{\pm \infty\}\), the maps
\[
\phi_!^\cp\colon (\Dperf(A),\QF^{\geq m}_M) \lrar (\Dperf(B),\QF^{\geq m}_{\phi_!M}) 
\quad \text{and} \quad 
\phi_!^\mathrm{c}\colon (\D^\mathrm{c}(A),\QF^{\geq m}_M) \lrar (\D^{\phi(\mathrm c)}(B),\QF^{\geq m}_{\phi_!M}).
\]
are a Poincaré-Karoubi and a Poincaré-Verdier projection, respectively, for every \(\mathrm c \subseteq \K_0(A)\) closed under the duality. 
\end{proposition}
Recall that $\QF^{\geq m}_M$ specialises to the quadratic and symmetric Poincar\'e structures when $m=\infty$ and $m=-\infty$, respectively, and that for $m=0,1,2$, it specialises to the animations of the functors on $\mathrm{Proj}(A)$ which extract genuine quadratic, even, and symmetric forms, respectively.
\begin{proof}
We use the following two inputs: firstly, \(B\) being a flat \(A\op\)-module implies that it can be written as filtered colimit of finitely generated free \(A\op\)-modules \(B_i\) and secondly, Tate cohomology commutes with filtered colimits of discrete modules in the coefficients. The former statement is a classical theorem of Lazard, see e.g.\ \cite{Lazard-platitude}*{Théorème 1.2} or \cite{stacks}*{Tag 058G}, and the second statement (for group cohomology) was discovered by Brown in \cite{Brown-finite}*{Theorem 3} for groups admitting a classifying space of finite type; given the \(2\)-periodicity of Tate cohomology for \(\Ct\) the case at hand also follows immediately from the same statement for group homology, which is obvious from the definitions. 

We start by considering the case \(m=\infty\). According to the general description of Kan extensions of modules with genuine involution above we need to show that \(B \otimes_A M^\tC \rightarrow ((B \otimes B) \otimes_{A \otimes A} M)^\tC\) is an equivalence. 
We can regard this as a natural transformation between spectrum valued functors
\[
X \longmapsto X \otimes_A M^\tC 
\quad \text{and} \quad 
X \longmapsto \big((X \otimes X) \otimes_{A \otimes A} M\big)^\tC
\]
on both the category of discrete \(A\op\)-modules and \(\D(A\op)\). From the latter case, we obtain that it is an equivalence for every perfect \(X\), as both sides are exact functors and the claim is evidently true for \(X = A\). In particular, the claim is true for all finitely generated projective \(A\op\)-modules since these are perfect when regarded in \(\D(A\op)\). Since filtered colimits are in particular sifted, regarding the two assignments as functors on the category of discrete \(A\op\)-modules the second fact makes them commute with filtered colimits of finitely generated free \(A\op\)-modules (for \(X\) a finitely generated free \(A\op\)-module, \(X \otimes X\) is a finitely generated free \((A \otimes A)\op\)-module, so \((X \otimes X) \otimes_{A \otimes A} M\) remains discrete, despite \(X \otimes X\) and \(A \otimes A\) potentially having higher homotopy). Taken together, the transformation is an equivalence for all flat \(A\op\)-modules, so in particular for \(X=B\) as desired. 

To obtain the case of the genuine Poincaré structures, just observe that the flatness of \(B\) also guarantees that the functor \(B \otimes_A - \colon \D(A) \rightarrow \D(B)\) commutes with the connective cover functors \(\tau_{\geq m}\) for all \(m \in \ZZ\).

The case of the quadratic Poincaré structure is Corollary~\reftwo{corollary:quadratic-Verdier-sequence} above.
\end{proof}

\begin{corollary}
\label{corollary:ore-localisation}%
Let \((A,\epsilon, \sigma)\) a ring with Wall anti-structure, and \(S \subseteq A\) a multiplicative subset satisfying the left Ore condition and closed under the involution \(\sigma\). Let \(M\) denote the module with involution over \(A\) given by endowing \(A\) with the \(A \otimes A\)-module structure arising from \(\sigma\) and the involution \(\epsilon \sigma\). Then for arbitrary \(m \in \ZZ \cup \{\pm\infty\}\), the maps 
\[
 \big(\Dperf(A),\QF^{\geq m}_M\big) \lrar \big(\Dperf(A[S^{-1}]),\QF^{\geq m}_{M[S^{-1}]}\big) \quad \text{and} \quad \big(\D^\mathrm c(A),\QF^{\geq m}_M\big) \lrar \big(\D^{\im(\mathrm c)}(A[S^{-1}]),\QF^{\geq m}_{M[S^{-1}]}\big)
\]
are a Poincar\'e-Karoubi and a Poincar\'e-Verdier projection, respectively, for every $c \subseteq \K_0(A)$ closed under the duality.
\end{corollary}
\begin{proof}
Note only that \(A[S^{-1}]\) is a derived localisation with perfectly generated fibre by Example~\reftwo{example:ore} and is furthermore 
flat as the construction of \(A[S^{-1}]\) as one-sided fractions displays it as a filtered colimit of free \(A\op\)-modules of rank \(1\).
Proposition~\reftwo{proposition:left-Kan-extension-specific-structures} therefore applies.
\end{proof}

\begin{remark}
The Ore condition is in fact often necessary to achieve flatness of the localisation: In \cite{teichner-ore}*{Main Theorem}, Teichner shows that if \(S\) is the set of elements that become invertible modulo a two-sided ideal \(I\), then flatness of \(A[S^{-1}]\) as a right \(A\)-module is equivalent to \(S\) being left Ore. \\
\end{remark}

\subsection{Construction and preservation of Poincaré-Verdier sequences} 
\label{subsection:preservation-poincare-verdier}%

In this section, we record several general constructions of (split) Poincar\'e-Verdier and Poincar\'e-Karoubi sequences arising from the tensor and cotensor construction as well as the closed symmetric monoidal structure of $\Catp$, that we will use throughout the sequel.

Before we dive in, let us give a simple recognition criterion for Poincaré-Verdier sequences involving hyperbolic Poincaré \(\infty\)-categories. Recall from Corollary~\refone{corollary:hyp-is-adjoint} and Remark~\refone{remark:units-and-counits} that \(\Hyp\) is both left and right adjoint to the underlying stable \(\infty\)-category functor \(\Catp \to \Catx\).

\begin{lemma}
\label{lemma:PV-hyperbolic}%
Let \(g \colon \C \rightarrow \D\) be an exact functor. Then, for a Poincaré structure \(\QF\) on \(\C\), the functor 
\[
g^\hyp \colon (\C,\QF) \longrightarrow \Hyp(\D)
\]
obtained by right adjointness of \(\Hyp\)
is a split Poincaré-Verdier projection if and only if \(g\) is a split Verdier projection and the restrictions of \(\QF\) to both the essential images of \(l\op\) and \(\Dual_\QF\op \circ r\) vanish, where \(l\) and \(r\) denote the adjoints of \(g\).

Similarly, for a Poincaré structure \(\QFD\) on \(\D\), the functor
\[
g_\hyp \colon \Hyp(\C) \longrightarrow (\D,\QFD)
\]
obtained by left adjointness of \(\Hyp\)
is a split Poincaré-Verdier inclusion if and only if \(g\) is a split Verdier inclusion and the restrictions of \(\QFD\) to both the essential images of \(g\op\) and \(\Dual_\QFD\op \circ g\) vanish.
\end{lemma}

\begin{proof}
Let us prove the first statement, the proof of the second is similar (and easier). 
First, it is easy to check that the functor \(g^\hyp\), which is given by
\[
(g, g\op \circ \Dual_\QF\op) \colon \C \longrightarrow \D \oplus \D\op
\]
admits a left adjoint \(l'\)
if and only if \(g\) admits both adjoints \(l,r\),
in which case the left adjoint \(l'\) 
is given by \((d,d') \mapsto ld \oplus \Dual_\QF(rd')\).
It is then clear from this description that \(l'\) is fully faithful if and only if \(l\) and \(r\) are fully faithful and the mapping spectra \(\map(ld,\Dual_\QF(rd')) = \map(d,g\Dual_{\QF}(rd'))\) and \(\map(\Dual_\QF(rd'), ld) = \map(g\Dual_{\QF}(ld),d')\) vanish for every \(d,d'\). By adjunction, this vanishing is equivalent to the vanishing of the mapping spectra \(\map(\Dual_{\QF}(rd'),rd)\) and \(\map(ld',\Dual_\QF(ld))\) for every \(d\), and hence to the vanishing of \(\Bil_{\QF}\) on both the image of \(l\) and the image of \(\Dual_{\QF}r\). By Proposition~\reftwo{proposition:criterion-poincare-split} \reftwoitem{item:fiber-left-adjoint}, to finish the proof it will hence suffice to show that in the case where \(l'\) is fully faithful, the map \(\QF(l'(d,d')) \to \QF_{\hyp}(d,d')\) is an equivalence if and only if \(\QF\) vanishes on the images of both \(l\) and \(\Dual_{\QF}r\). Indeed, unwinding the definitions, the fibre of this map is exactly \(\QF(ld) \oplus \QF(\Dual_{\QF}(rd'))\).
\end{proof}

Let us now consider examples involving the tensor and cotensor constructions as considered in \S\refone{subsection:cotensoring} and \refone{subsection:tensoring}; several of the results will be required for our analysis of the hermitian \(\Q\)-construction in the next section. For this, let us recall that for a Poincaré \(\infty\)-category \((\C,\QF)\) and an \(\infty\)-category \(\I\), the \(\I\)-tensor \((\C,\QF)_\I = (\C_\I,\QF_\I) =(\I \otimes \C,\QF_\I)\) and \(\I\)-cotensor \((\C,\QF)^{\I} = (\C^{\I},\QF^{\I}) = (\Fun(\I,\C),\QF^{\I})\) of \((\C,\QF)\) are generally only hermitian \(\infty\)-categories, and that for a given functor \(\alp\colon \I \to \J\), the induced functors \((\alp_*,\eta_{\alp})\colon (\C,\QF)_{\I} \to (\C,\QF)_{\J}\) and \((\alp^*,\eta^{\alp})\colon (\C,\QF)^{\J} \to (\C,\QF)^{\I}\) are generally only hermitian functors. Here we note that the underlying exact functor \(\alp^*\) is given by restriction along \(\alp\) and the underlying exact functor \(\alp_*\) is given by right Kan extension along \(\alp\op\) under an identification of the tensor \(\I \otimes \C\) with a certain full subcategory of \(\Fun(\I\op,\Pro(\C))\). When \(\I\) is a finite poset or more generally a strongly finite $\infty$-category (i.e. $|\I|$ and $\Hom_\I(i,j)$ are finite for all $i,j \in \I$) this full subcategory is just \(\Fun(\I\op,\C)\). In addition, if \(\I\) and \(\J\) are strongly finite and \(\bet\colon \I \to \J\) is a cofinal map we also have the exceptional hermitian functors \((\bet^*,\vartheta^{\beta})\colon (\C,\QF)_\J \to (\C,\QF)_\I\) and \((\bet_*,\vartheta_{\beta})\colon (\C,\QF)^{\I} \to (\C,\QF)^{\J}\) going in the opposite direction to the usual functoriality. Here the underlying exact functors \(\bet^*\) and \(\bet_*\) are again given by restriction and right Kan extension, respectively.
\begin{proposition}
\label{proposition:localizationforp4}%
\label{proposition:kan-extension-cofinal}%
Let \(\alp\colon \I \to \J\) be a functor of \(\infty\)-categories which exhibits \(\J\) as a localisation of \(\I\) with respect to a collection of maps \(W\), and let \(\C\) be a stable \(\infty\)-category. 
Let \(\ovl{W}\) be the collection of images in \(\C_\I\) of the arrows of the form \(\id_x \times \alp\) with \(x \in \C\) and \(\alp\) in \(W\). 
Then the following hold:
\begin{enumerate}
\item
\label{item:induced-functor-verdier}%
The induced functor \(\alp_*\colon \C_\I \to \C_\J\) is a Verdier projection whose kernel is the smallest stable subcategory of \(\C_\I\) closed under retracts and containing the fibres of the maps in \(\ovl W\). 
In addition, if \(\QF\) is a Poincaré structure on \(\C\) such that \((\C,\QF)_{\I}\) is Poincaré and \(\ker(\alp_*)\) 
is closed under its duality then 
\[
(\alp_*,\eta_{\alp})\colon (\C,\QF)_{\I} \lrar (\C,\QF)_{\J}
\]
is a Poincaré-Verdier projection.
\item 
\label{item:induced-functor-verdier-inclusion}%
The induced functor \(\alp^*\colon \C^\J \rightarrow \C^\I\) is a Verdier inclusion. In addition, if \(\QF\) is a Poincaré structure on \(\C\) such that \((\C,\QF)^\I\) is Poincaré and \(\alp^*(\C^\J)\) is closed under the duality then
\[
(\alp^*,\eta^{\alp})\colon (\C,\QF)^\J \lrar (\C,\QF)^{\I}
\]
is a Poincaré-Verdier inclusion. 
\item
\label{item:kan-extension-cofinal}%
If \(\I\) and \(\J\) 
are strongly finite and the exceptional pushforward \((\alp_*,\vartheta_{\alp})\colon (\C,\QF)^{\I} \to (\C,\QF)^{\J}\) is duality preserving then it is a split Poincaré-Verdier projection.
\end{enumerate}
\end{proposition}

\begin{proof}
The universal properties of \(\C_\I\) and \(\C_\J\) as tensors imply that for every stable \(\infty\)-category \(\D\), we have a diagram with invertible horizontal arrows 
\[
\begin{tikzcd}[column sep=2ex]
\Funx(\C_{\J},\D) \ar[rr]\ar[d] && \Fun^l(\C \times \J,\D) \ar[d]\ar[rr] && \Fun(\J,\Funx(\C,\D))\ar[d] \\
\Funx(\C_{\I},\D) \ar[rr] && \Fun^l(\C \times \I,\D) \ar[rr] && \Fun(\I,\Funx(\C,\D)) , \
\end{tikzcd}
\]
where \(\Fun^l(\C \times \I,\D)\) denotes the full subcategory of \(\Fun(\C \times \I,\D)\) spanned by those functors which are exact in the \(\C\) entry, and similarly for \(\Fun^l(\C \times \J,\D)\). The assumption that \(\I \to \J\) exhibits \(\J\) as the localisation of \(\I\) at \(W\) now implies that the right most vertical arrow is fully faithful, with essential image spanned by those diagrams \(\I \to \Funx(\C,\D)\) which send the arrows in \(W\) to equivalences. It then follows that the left most vertical arrow is fully faithful with essential image spanned by those exact functors \(\C_\I \to \D\) which send the arrows in \(\ovl{W}\) to equivalences. By Corollary~\reftwo{corollary:criterion-projection} this establishes the first part of \reftwoitem{item:induced-functor-verdier}.

Now suppose that \(\QF\) is a Poincaré structure on \(\C\) such that 
\((\C,\QF)_\I\)
is Poincaré and \(\E := \ker(\alp_*)\) 
is closed under the associated duality. 
Then \((\E,(\QF_{\I})|_{\E})\) is Poincaré as well. 
In light of Proposition \reftwo{proposition:criterion-poincare}, to show that 
\((\alp_*,\eta_\alp)\)
is a Poincaré-Verdier projection it thus suffices to show that extended by its kernel it is a cofibre sequence in \(\Cath\), since the forgetful functor \(\Catp \rightarrow \Cath\) preserves and detects colimits. For this, by the universal properties of \((\C,\QF)_\I\)
 and \((\C,\QF)_\J\) as tensors in \(\Cath\), we obtain for every hermitian \(\infty\)-category \((\D,\QFD)\) a diagram
\[
\begin{tikzcd}[column sep=2ex]
\Fun(\J,\Funh((\C,\QF),(\D,\QFD)))\ar[d] \ar[r,phantom,"\simeq"] & \Funh((\C,\QF)_\J,(\D,\QFD)) \ar[d]\ar[rr] && \Funx(\C_\J,\D) \ar[d]\ar[r,phantom,"\simeq"] & \Fun(\J,\Funx(\C,\D))\ar[d] \\
\Fun(\I,\Funh((\C,\QF),(\D,\QFD))) \ar[r,phantom,"\simeq"] & \Funh((\C,\QF)_\I,(\D,\QFD)) \ar[rr] && \Funx(\C_\I,\D) \ar[r,phantom,"\simeq"] & \Fun(\I,\Funx(\C,\D)) . \
\end{tikzcd}
\]
Since \(\I \to \J\) is a localisation at \(W\), the left most vertical arrow is fully faithful with essential image consisting of those diagrams \(\I \to \Funh((\C,\QF),(\D,\QFD))\) which send the arrows in \(W\) to equivalences, and similarly for the right most arrow. The conservativity of \(\Funh((\C,\QF),(\D,\QFD)) \to \Funx(\C,\D)\) implies that the outer square (and hence the middle one) are cartesian and since \(\C_{\I} \to \C_{\J}\) is a Verdier projection,
\[
\Funh((\C,\QF)_\J,(\D,\QFD)) \lrar \Funh((\C,\QF)_\I,(\D,\QFD))
\]
is fully faithful with essential image those hermitian functors \((\C,\QF)_\I \to (\D,\QFD)\) which invert \(\ovl W\), i.e.\ which vanish on the kernel of \(\C_\I \rightarrow \C_\J\). This shows the second part of \reftwoitem{item:induced-functor-verdier}.

For \reftwoitem{item:induced-functor-verdier-inclusion},
we recall that \(\C^\I \simeq \Fun(\I,\C)\), so the fact that \(\I \rightarrow \J\) is a localisation implies that \(\C^\J \rightarrow \C^\I\) is fully faithful, and the characterisation of the image as those functors that invert \(\W\) shows that the image is closed under retracts. The first part of \reftwoitem{item:induced-functor-verdier-inclusion} now follows from Proposition \reftwo{proposition:inclusion-criterion}.
For the second part of \reftwoitem{item:induced-functor-verdier-inclusion}, we have to check that \(\QF^\J \simeq p^*\QF^\J\), but this follows from the definition of the hermitian structure on cotensors as a limit, together with the fact that localisations are (co)final.

Finally, let us prove \reftwoitem{item:kan-extension-cofinal}.
First up, restriction along \(\alpha\) is left adjoint to \(\alpha_*\) and fully faithful since \(\alpha\) is a localisation. Thus, by \reftwo{corollary:poincare-projection}, we are left to prove that the natural map
\(
\QF^\J(\varphi) \to \QF^\I(\varphi \circ \alpha)
\)
is an equivalence. Indeed,
\[
\QF^\J(\varphi) \simeq \lim_{j \in \J\op} \QF(\varphi(j)) \simeq \lim_{i \in \I\op} \QF(\varphi \alpha(i)) \simeq \QF^\I(\varphi \circ \alpha)
\]
again since \(\alpha\op\) is final.
\end{proof}

We now consider a somewhat opposite case, where \(\alp\colon \I \to \J\) is very far from being a localisation. 
\begin{proposition}
\label{proposition:example-cosieve}%
\label{proposition:example-cosieve-2}%
Let \(\alp\colon \I \hrar \J\) be an upwards closed inclusion between strongly finite categories, and let \((\C, \QF)\) be a Poincaré \(\infty\)-category. 
\begin{enumerate}
\item
\label{item:cosieve-cotensor}%
If the hermitian \(\infty\)-categories \((\C,\QF)^{\I}\) and \((\C,\QF)^{\J}\) are Poincaré then 
\((\alp^*,\eta^{\alp})\colon (\C,\QF)^{\J} \to (\C,\QF)^\I\) is a split Poincaré-Verdier projection. 
\item
\label{item:cosieve-tensor}%
If \((\C,\QF)_\I\) and \((\C,\QF)_\J\) are Poincaré then \((\alp_*,\eta_{\alp})\colon (\C,\QF)_{\I} \to (\C,\QF)_\J\) is a split Poincaré-Verdier inclusion. 
\end{enumerate}
\end{proposition}
Here a fully faithful functor $\alpha \colon \I \rightarrow \J$ is called upwards closed, if $\Hom_\J(\alpha(i),j) \neq \emptyset$ for some $i \in \I$ and $j \in \J$ implies that $j$ lies in the essential image of $\alpha$.
\begin{proof}
We first prove~\reftwoitem{item:cosieve-cotensor}. For this, first note that \((\alp^*,\eta^{\alp})\) is a Poincaré functor by Proposition~\refone{proposition:induced-maps-duality-pres}. In addition, a fully faithful left adjoint to \(\alpha^*\) is given by the exact functor
\(\alpha_!\colon \C^{\I} \to \C^{\J}\)
performing left Kan extension. In fact, identifying \(\I\) with a subposet of \(\J\) via \(\alp\), the fact that \(\I\) is upwards closed means that this left Kan extension admits a very explicit formula: for a diagram \(\vphi\colon \I \to \C\) the value of \(\alpha_!\vphi\) is given by
\[
\alpha_!\vphi(j) = \left\{\begin{matrix} \vphi(j) & j \in \I \\ 0 &  j \notin \I \ . \\ \end{matrix}\right.
\]
Since \(\QF(0) \simeq 0\), the spectrum valued diagram \(j \mapsto \QF(r_!\vphi(j))\) is, for a similar reason, a right Kan extension of its restriction to \(\I\op\), and so the natural map
\[
\QF^{\J}(\alpha_!\vphi) \simeq \lim_{j \in \J\op} \QF(\alpha_!\vphi(j)) \to \lim_{i \in \I\op} \QF(\vphi(i)) = \QF^{\I}(\vphi)
\]
is an equivalence. The Poincaré functor
\((\alp^*,\eta^{\alp})\) is then Poincaré-Verdier projection by Corollary~\reftwo{corollary:split-poincare-projection-inclusion} \reftwoitem{item:split-poinc-verdier-proj-concrete}.

We now prove~\reftwoitem{item:cosieve-tensor}. Here \((\alp_*,\eta_{\alp})\) is Poincaré by Corollary~\refone{corollary:tensor-duality-preserving}. Now on the level of underlying stable \(\infty\)-categories the right Kan extension functor \(\alpha_*\) is fully faithful (since \(\alpha\) is) and admits a left adjoint given by restriction. To finish the proof, it will suffice by Corollary~\reftwo{corollary:split-poincare-projection-inclusion} \reftwoitem{item:split-poinc-verdier-incl-concrete} to show that for every diagram \(\vphi\colon \I\op \to \C\), the composite map
\[
\displaystyle\mathop{\colim}_{i \in \I}\QF(\vphi(i)) \lrar \displaystyle\mathop{\colim}_{i \in \I}\QF(\alpha^*\alpha_*\vphi(i)) \lrar \displaystyle\mathop{\colim}_{j \in \J}\QF(\alpha_*\vphi(j))
\]
is an equivalence. Here, the first map is an equivalence since \(\alpha_*\) is fully faithful. To see that the second map is an equivalence, we argue as in the proof of the first part
and observe that
\[
\alpha_*\vphi(j) = \left\{\begin{matrix} \vphi(j) & j\in \I \\ 0 &  j \notin \I \\ \end{matrix}\right.
\]
The spectrum valued diagram \(j \mapsto \QF(\alpha_*\vphi(j))\) is then, for a similar reason, the left Kan extension of its restriction to \(\I\), and so the second map above is an equivalence as well.
\end{proof}

\begin{proposition}
\label{proposition:(co)tensor-Verdier}%
Let us consider a (split) Poincaré-Verdier sequence \((\C,\QF) \rightarrow (\D,\QFD) \rightarrow (\E, \QFE)\) and a strongly finite $\infty$-category \(\I\), e.g.\ a finite poset, such that \((-)_\I \colon \Catp \rightarrow \Cath\) preserves Poincaré \(\infty\)-categories. Then the induced sequences
\[
(\C,\QF)^{\I}\longrightarrow (\D,\QFD)^{\I} \longrightarrow (\E,\QFE)^{\I}
\quad
\text{and}
\quad
(\C,\QF)_{\I}\longrightarrow (\D,\QFD)_{\I} \longrightarrow (\E,\QFE)_{\I}
\]

are (split) Poincaré-Verdier sequences.
\end{proposition}

Note that by \refone{corollary:I-preserve-poincare}, the functor \((-)^\I\colon \Cath \rightarrow \Cath\) preserves \(\Catp\) provided \((-)_\I\) does, which in turn is equivalent to \((\Spa^\omega,\QF^\uni)_\I\) being Poincaré. At the level of underlying stable $\infty$-categories, the assumption that $(-)_\I$ preserve Poincar\'e $\infty$-categories is of course unnecessary.

\begin{proof}[Proof of Proposition~\reftwo{proposition:(co)tensor-Verdier}]
Let us treat the tensoring, the argument for the cotensoring being entirely dual. As a left adjoint, the tensoring construction generally preserves colimits, and by \refone{corollary:limit-preservation} tensoring with a finite poset also preserves limits. This gives the part of the statement disregarding splittings. But for example from \refone{proposition:tensor-strongly-finite} we find that the operation \((-)_\I\colon \Catx \rightarrow \Catx\) preserves adjoints, which implies the split case.
\end{proof}

Recall from Theorem \refone{theorem:tensorpoincare} that $\Catp$ has a closed monoidal structure lifting Lurie's tensor product on $\Catx$. Its interaction with Poincar\'e-Verdier sequences is as follows:

\begin{proposition}
\label{lemma:Morita-exact-criterion}%
For a fixed Poincaré $\infty$-category $(\cB,\Pi)$ the functor $(\cB,\Pi) \otimes (-)\colon \Catp \to \Catp$ has the following properties:
\begin{enumerate}
\item
\label{item:preserves-K-equivalences}%
It preserves Karoubi equivalences.
\item
\label{item:preserves-PV-proj}%
It preserves Poincaré-Verdier projections.
\item 
\label{item:preserves-split-PV-sequences}%
It preserves split Poincar\'e-Verdier sequences.
\item
\label{item:sends-PV-incl-to-PK-inc}%
It preserves Poincaré-Karoubi inclusions.
\item
\label{item:sends-PV-squares-to-PK-squares}%
It preserves Poincaré-Karoubi sequences.
\end{enumerate}
\end{proposition}

\begin{proof}[Proof of Lemma~\reftwo{lemma:Morita-exact-criterion}]
\reftwoitem{item:preserves-K-equivalences} follows from the fact that the functor $\Funx((\cB,\Pi),-)\colon \Catp \to \Catp$, right adjoint to $(\cB,\Pi) \otimes (-)$, preserves idempotent complete Poincaré $\infty$-categories, and \reftwoitem{item:preserves-PV-proj} follows from 
the fact that $(\cB,\Pi)\otimes (-) $ is a left adjoint.

For \reftwoitem{item:preserves-split-PV-sequences}, we first consider the underlying stable $\infty$-categories.
Since the symmetric monoidal structure on $\Catx$ is closed we may view $\Catx$ as a $\Catx$-enriched category, in which case the functor $\cB \otimes (-)\colon \Catx \to \Catx$ canonically extends to a $\Catx$-enriched functor. It hence sends adjunctions to adjunctions, and adjunctions whose (co)units are equivalences to adjunctions whose (co)units are equivalences. We conclude that $\cB \otimes (-)$ preserves split Verdier inclusions and split Verdier projections. Since this functor also preserves colimits, we conclude that it preserves split Verdier sequences. It remains to show that $(\cB,\Pi) \otimes (-)$ preserves split Poincar\'e-Verdier inclusions. But using the linear-bilinear decomposition of Poincaré structures on tensor products given by Proposition~\refone{proposition:tensor}, we see that if $(i,\eta)\colon (\C,\QF) \to (\D,\QFD)$ is a split Poincaré-Verdier inclusion then the restriction of the Poincaré structure $\Pi \otimes \QFD$ of $(\cB,\Pi)\otimes(\D,\QFD)$ to $\cB \otimes \C$ coincides with $\Pi \otimes \QF$ as needed.

For \reftwoitem{item:sends-PV-incl-to-PK-inc}, we first show that if $(\C,\QF) \to (\D,\QFD)$ is a Karoubi inclusion then $\cB \otimes \C \to \cB \otimes \D$ is fully faithful. Indeed, to see this we may apply $\Ind(-)$, use that it is symmetric monoidal, and the fact that $\Ind(\C) \to \Ind(\D)$ is a (right) split inclusion. The same argument as above then implies that $\Ind(\cB) \otimes (-)$ preserves (right) split inclusions, and that the restriction of $\Pi \otimes \QFD$ to $\cB\otimes \C$ is given by $\Pi \otimes \QF$ as needed.
For \reftwoitem{item:sends-PV-squares-to-PK-squares}, since $(\cB,\Pi) \otimes (-)$ preserves colimits, it will suffice to show that it preserves Poincar\'e-Karoubi inclusions, which we have just shown.
\end{proof}

\begin{remark}
\label{remark:no-verdier}%
The functor $(\cB,\Pi) \otimes (-)\colon \Catp \to \Catp$ does \emph{not} preserve Poincaré-Verdier sequences in general. In fact, already on the level of stable $\infty$-categories the operation $\C \otimes (-)\colon \Catx \to \Catx$ does not preserve Verdier inclusions. To see this, let $R$ be a commutative ring with $\KK_{-1}(R) \neq 0$. For instance, consider a field $k$ and $R \subseteq k[t]$ the subring consisting of polynomials with common value on $0$ and $1$, i.e.\ the pullback of $k[t]$ over $k\times k$ along the diagonal $k \to k \times k$. Then the excision sequence \cite{Bass}*{Theorem XII.8.3} gives $\K_{-1}(R) \cong \ZZ$. Let us then denote by $\Dperf(R \,\mathrm{on}\, t)$ the kernel of the Karoubi projection $\Dperf(R[t]) \to \Dperf(R[t^{\pm 1}])$ and by $i \colon \Dperf(R \,\mathrm{on}\, t) \subseteq \Dperf(R[t])$ the inclusion. We then claim that the map 
\[
\Dperf(R\,\mathrm{on}\, t)_{\B\ZZ} \lrar \Dperf(R[t])_{\B\ZZ}
\]
obtained by tensoring $i$ with $\Spaf_{\B\ZZ}$ is \emph{not} a Verdier inclusion. As it is a Karoubi inclusion by Lemma~\reftwo{lemma:Morita-exact-criterion}, we have to show that the image is not closed under retracts. By Thomason's classification of dense subcategories, this is equivalent to the map 
\[
\KK_0(\Dperf(R\,\mathrm{on}\, t)_{\B\ZZ})/\K_0(\Dperf(R\,\mathrm{on}\, t)_{\B\ZZ}) \lrar \KK_0(\Dperf(R[t])_{\B\ZZ})/\K_0(\Dperf(R[t])_{\B\ZZ})
\]
not being injective. By the fundamental theorem in algebraic K-theory, see e.g.\ \cite{saunier}*{Theorem 1.3} for the exact version we use here, this map contains the map
\[
\KK_{-1}(\Dperf(R\,\mathrm{on}\, t)) \lrar \KK_{-1}(\Dperf(R[t]))
\]
as a direct summand. This map in turn participates in a long exact sequence
\[
\KK_0(R[t]) \lrar \KK_0(R[t^{\pm 1}]) \lrar \KK_{-1}(\Dperf(R\,\mathrm{on}\, t)) \lrar \KK_{-1}(\Dperf(R[t]))
\]
so it suffices to recall that the cokernel of the first map contains $\KK_{-1}(R)$ as a direct summand (by another application of the fundamental theorem).

We can also promote this example to the level of Poincaré $\infty$-categories by endowing $\Dperf(R[t])$ and $\Dperf(R[t^{\pm1}])$ with the symmetric Poincaré structures.
\end{remark}

Finally, we consider the analogous questions of $\Funx((\B,\Pi),-)$ instead of $(\cB,\Pi) \otimes (-)$:

\begin{lemma}
\label{lemma:and-smashing}%
For a fixed $(\cB,\Pi) \in \Catp$ the functor 
$
\Funx((\cB,\Pi),-)\colon \Catp \to \Catp
$
preserves split Poincaré-Verdier sequences. %
\end{lemma}
\begin{proof}
Since $\Funx((\cB,\Pi),-)$ preserves limits it will suffice to show that $\Funx((\cB,\Pi),-)$ sends split Poincaré-Verdier projections to split Poincaré-Verdier projections. For this, let $(p,\eta)\colon(\D,\QFD) \to (\E,\QFE)$ be a split Poincaré-Verdier projection and let $g,h \colon \E \to \D$ be the left and right adjoints of $p$. Then $g_*,h_*\colon \Funx(\cB,\E) \to \Funx(\cB,\D)$ are left and right adjoints of $p_*$, showing that $p_*$ is a split Verdier projection. To see that it is a Poincar\'e-Verdier projection, by Corollary~\reftwo{corollary:split-poincare-projection-inclusion} it then suffices to show
that for every exact functor $f\colon \cB \to \E$ the composite map of spectra
\[
\nat(\Pi,f^*g^*\QFD)) \lrar \nat(\Pi,f^*g^*p^*\QFE) \lrar \nat(\Pi,f^*\QFE)
\]
is an equivalence. This map is induced by the composite
$
g^*\QFD \Rightarrow g^*p^*\QFE \Rightarrow \QFE
$
which is an equivalence by the assumption that $(p,\eta)$ is a split Poincaré-Verdier projection (and Corollary~\reftwo{corollary:split-poincare-projection-inclusion} yet again).
\end{proof}

\begin{remark}
\label{remark:not-smashing}%
Unlike the case of $(\cB,\Pi) \otimes (-)$, the functor $\Funx((\cB,\Pi),-)$ does not preserve Poincaré-Karoubi sequences in general.
To see this, let us take $(\cB,\Pi) = (\Spaf,\QF^{\uni})_{\B\ZZ}$, so that 
\[
\Funx((\Spaf,\QF^{\uni})_{\B\ZZ},-) = (-)^{\B\ZZ}.
\]
We claim that the operation of cotensoring with $\B\ZZ$ preserves neither Verdier nor Karoubi projections, already on the level of underlying stable $\infty$-categories. For instance, $\Dperf(\ZZ) \to \Dperf(\QQ)$ is a (Poincar\'e)-Verdier projection, but the essential image of $p\colon \Fun(\B\ZZ,\Dperf(\ZZ)) \to \Fun(\B\ZZ,\Dperf(\QQ))$ is not dense. Indeed, let us consider the object $(\QQ,\cdot 2) \in \Fun(\B\ZZ,\Dperf(\QQ))$. We claim that it is not a retract of an object in the image of $p$. Indeed, assume there exists $(M,\phi) \in \Fun(\B\ZZ,\Dperf(\ZZ))$ such that $p(M,\phi)$ contains $(\QQ,\cdot 2)$ as a retract. Then the automorphism $\phi_0$ of $H_0(M)/\mathrm{tors}$ induced by $\phi$ is such that its rationalisation contains $(\QQ,\cdot 2)$ as a retract. It follows that $\phi_0$ is represented by an integral invertible matrix which has $2$ as eigenvalue. This is impossible: For any eigenvalue $\lambda \in \mathbb C$ of an integral matrix $A$ the number $\det(A)/\lambda \in \CC$ is a product of eigenvalues of $A$ (e.g.\ via the Jordan normal form) and thus an algebraic integer.
\end{remark}

\subsection{Additive and localising functors}
\label{subsection:additive-localizing}%

In this section, we establish the basic notions of additive, Verdier-localising and Karoubi-localising functors. They are based on a mild generalisation of Poincaré-Verdier and Poincaré-Karoubi sequences in the form of certain bicartesian squares. Sending these particular bicartesian squares to cartesian squares isolates the localisation properties enjoyed by Grothendieck-Witt theory axiomatically.

In the present paper, we focus almost exclusively on Verdier-localising (or even additive) functors. Together with the principal example of the Karoubi-Grothendieck-Witt functor, Karoubi-localising functors are studied thoroughly in \paperfour and we only briefly mention them here for the sake of completeness.

\begin{definition}
\label{definition:poincare-square}%

A \defi{(split) Poincaré-Verdier square} is a diagram
\[
\begin{tikzcd}
(\D,\QFD) \ar[r]\ar[d] & (\Dtwo,\QFDtwo) \ar[d] \\
(\E,\QFE) \ar[r] & (\Etwo,\QFEtwo) 
\end{tikzcd}
\]
in \(\Catp\) which is cartesian and whose vertical maps 
are (split) Poincaré-Verdier projections. 
We say that such a square is a \defi{Poincaré-Karoubi square} if it becomes cartesian after applying the idempotent completion functor of Proposition~\reftwo{proposition:idempotent-completion} and its vertical maps are Poincaré-Karoubi projections.

We will also speak of (split) Verdier squares and Karoubi squares for the evident analogs of the above definition for squares in $\Catx$.
\end{definition}

\begin{remarks}
\label{remark:cocartesian}%
\label{remark:pullbackofVerdier}%
\label{remark:verdier-is-karoubi}%
\ %
\begin{enumerate}
\item 
A (split) Poincaré-Verdier square with lower left corner \(0 \in \Catp\) is exactly a (split) Poincaré-Verdier sequence. The same holds for Poincaré-Karoubi sequences.
\item 
Any Poincaré-Verdier square is also cocartesian in \(\Catp\): Indeed, extend the defining square to a commutative rectangle
\[
\begin{tikzcd}
(\C,\QF) \ar[r]\ar[d] & (\D,\QFD) \ar[r]\ar[d] & (\Dtwo,\QFDtwo) \ar[d] \\
{0} \ar[r] & (\E,\QFE) \ar[r] & (\Etwo,\QFEtwo) 
\end{tikzcd}
\]
in which both squares are cartesian and the vertical maps are Poincaré-Verdier projections. Then the external rectangle is cartesian by the pasting lemma, see \cite{HTT}*{Lemma 4.4.2.1}, hence cocartesian since the right vertical map is a Verdier projection. For the same reason, the left square is cocartesian and so the right square is cocartesian by the pasting lemma. Similarly, every Poincaré-Karoubi square becomes cocartesian in \(\Catpi\) after applying idempotent completion.
\item
\label{item:PV-stable-by-pullback}%
The classes of split (Poincaré-)Verdier, (Poincaré-)Verdier, and (Poincaré-)Karoubi squares are all closed under pullbacks, see Corollaries \reftwo{corollary:split-verdier-proj-pullback}, \reftwo{corollary:pullback-non-split-verd-proj}, \reftwo{corollary:Karoubi-proj-pullback}, \reftwo{corollary:PV-proj-pullback} and \reftwo{corollary:PK-projpullbackstable}. In the definition of any of these six types of squares, one may therefore just require the right vertical map to be a projection of the given kind.
\item
\label{lemma:criterion-Verdier-squares}%
By Proposition~\reftwo{proposition:verdier-square-char}, a square as in Definition~\reftwo{definition:poincare-square} is a Verdier square if and only if it is adjointable in the sense of Definition~\reftwo{definition:adjointable}, has vertical arrows Verdier projections, and the induced map on vertical kernels is an equivalence in \(\Catx\). Similarly, by Proposition~\reftwo{proposition:poincare-verdier-square-char}, a square is a Poincaré-Verdier square if and only if it is adjointable, has vertical arrows Poincaré-Verdier projections, and the induced map on vertical kernels is an equivalence in \(\Catp\). The analogous characterisation in the split case follows immediately.
\item 
Proposition~\reftwo{proposition:minimal-poincare} implies that every Poincaré-Verdier square is a Poincaré-Karoubi square. Conversely, a Poincaré-Karoubi square involving idempotent complete Poincaré \(\infty\)-categories is a Poincaré-Verdier square if and only if its vertical maps are essentially surjective. Moreover, the idempotent completion of a Poincar\'e-Karoubi square is again Poincar\'e-Karoubi.
\end{enumerate}
\end{remarks}

\begin{definition}
\label{definition:additive}%
Let \(\E\) be an \(\infty\)-category which admits finite limits and \(\F\colon \Catp \to \E\) a functor. Recall that \(\F\) is said to be \defi{reduced} if \(\F({0})\) is a terminal object in \(\E\). We say that a reduced functor \(\F\) is \defi{additive}, \defi{Verdier-localising} or \defi{Karoubi-localising}, if it  takes split Poincaré-Verdier squares, arbitrary Poincaré-Verdier squares or  Poincaré-Karoubi squares to cartesian squares, respectively. 

We shall denote the \(\infty\)-categories of these functors by 
\[
\Funadd(\Catp,\E), \quad \Funvloc(\Catp,\E), \quad \text{and}\quad \Funkloc(\Catp,\E),
\]
respectively. 
\end{definition}

It follows from Remark~\reftwo{remark:verdier-is-karoubi} that there are inclusions
\[
\Funkloc(\Catp,\E) \subseteq \Funvloc(\Catp,\E) \subseteq \Funadd(\Catp,\E)
\]
as full subcategories. We note that additive, Verdier-localising and Karoubi-localising invariants are closed in \(\Fun(\Catp,\E)\) under limits (which are computed pointwise), such as taking loops. Colimits on the other hand are generally not computed pointwise (unless \(\E\) is stable), and in fact, we will discuss in detail the suspension operation in the next section as part of our analysis of the universal property of Grothendieck--Witt theory.

\begin{warn}
Here we follow the convention of the fifth author and Tamme to divorce the preservation of filtered colimits from the preservation of certain fibre sequences and squares. As a result, the \(\infty\)-categories appearing in the end of Definition~\reftwo{definition:additive} are (presumably) \emph{not} locally small. The reader who is adverse to non-locally small \(\infty\)-categories is invited to restrict attention only to accessible additive/Verdier-localising/Karoubi-localising functors; this will not affect any of the statements in this paper, nor their proofs. 

We also note that if one fixes a regular cardinal \(\kappa\) and restricts attention only to those additive/Verdier-localising/Karoubi-localising functors that preserve \(\kappa\)-filtered colimits, then the corresponding variants of \(\Funadd(\Catp,\E)\), \(\Funvloc(\Catp,\E)\) and \(\Funkloc(\Catp,\E)\) become presentable, and a reader who so prefers may fix at this moment once and for all a sufficiently large such \(\kappa\). At any rate, the most interesting examples of such functors that appear in this paper, such as the Grothendieck-Witt, \(\K\)-and \(\L\)-theory spectra, even preserve (\(\omega\)-)filtered colimits.
\end{warn}

\begin{proposition}
\label{proposition:squares-and-sequences}%
A reduced functor \(\F \colon \Catp \rightarrow \E\) with \(\E\) stable is additive, Verdier-localising or Karoubi-localising if and only if it takes split Poincaré-Verdier, Poincaré-Verdier or all Poincaré-Karoubi sequences to fibre sequences in \(\E\).
\end{proposition}

\begin{proof}
Apply \(\F\) to the rectangle in Remark~\reftwo{remark:cocartesian} and use the pasting lemma.
\end{proof}

For non-stable \(\E\) we expect, however, that the condition of being additive or Verdier-localising is strictly stronger than sending split Poincaré-Verdier or Poincaré-Verdier sequences to fibre sequences, and similarly for the condition of being Karoubi-localising. We will need the stronger variant in \S\reftwo{subsection:fibcob} with target \(\Sps\), when we discuss the additivity theorem for cobordism categories.

\begin{proposition}
\label{proposition:localizing-criterion}%
A functor \(\F\colon \Catp \to \E\) is Karoubi-localising if and only if it is Verdier-localising and invariant under Karoubi equivalences.
\end{proposition}

\begin{proof}
The ``only if'' part follows from Remark~\reftwo{remark:verdier-is-karoubi} and the fact that  
\[
\begin{tikzcd}
(\C,\QF) \ar[r] \ar[d] & {0}\ar[d] \\
 (\C,\QF)^{\natural} \ar[r] & {0}
\end{tikzcd}
\]
forms a Poincaré-Karoubi square. For the converse, assume $\F$ is Verdier-localising and invariant under Karoubi equivalences. It then suffices to show that $\F$ sends a Poincar\'e-Karoubi square all whose terms are idempotent complete to cartesian squares. Such squares are in turn the idempotent completion of cartesian squares where the vertical maps are Poincar\'e-Verdier projections: Indeed, simply replace the bottom terms with the essential image of the Poincar\'e-Karoubi projections and use Corollary~\reftwo{corollary:criterion-poincare-karoubi-projection}. We arrive at the desired conclusion as $\F$ is Verdier-localising.
\end{proof}

\begin{lemma}
\label{lemma:Funadd-semiadd}%
The categories $\Funadd(\Catp,\E)$, $\Funvloc(\Catp,\E)$, and$\Funkloc(\Catp,\E)$
are semi-additive and the forgetful functor 
$
\Funadd(\Catp,\Mon_{\Einf}(\E)) \to \Funadd(\Catp,\E)
$
and its localising analogues are equivalences.
\end{lemma}

Since the \(\infty\)-category \(\Catp\) is semi-additive, see Proposition~\refone{proposition:catp-pre-add}, but not additive, the analogous statement for \(\Grp_{\Einf}(\E)\) in place of \(\Mon_{\Einf}(\E)\) is not true. Noting the unfortunate clash in the use of the word \emph{additive} arising from this, we set:

\begin{definition}
\label{definition:grouplike}%
An additive functor \(\F \colon \Catp \rightarrow \E\) is called \defi{group-like} if the canonical lift of \(\F\) to \(\Mon_{\Einf}(\E)\) actually takes values in the full subcategory \(\Grp_{\Einf}(\E)\).
\end{definition}

Equivalently, this is the same as saying that for every Poincaré \(\infty\)-category \((\C,\QF)\) it takes the shear map \((\C,\QF) \times (\C,\QF) \to (\C,\QF) \times (\C,\QF)\) (given at the level of objects by \((\x,\y) \mapsto (\x,\x \oplus \y)\)) to an equivalence in \(\E\).

\begin{proof}[Proof of Lemma \reftwo{lemma:Funadd-semiadd}]
Given the semi-additivity of \(\Catp\), the \(\infty\)-category of product preserving functors \(\Catp \rightarrow \E\) is semi-additive by \cite{GepGroNik}*{Corollary 2.4}. But products of additive, Verdier or Karoubi localising functors are again such, which implies the first statement. The second follows from \cite{GepGroNik}*{Corollary 2.5 iii)} together with the observation that $\F\colon \Catp \to \Mon_{\Einf}(\E)$ is additive if and only if it is additive when viewed as taking values in $\E$ itself. 
\end{proof}

\begin{remark}
\label{remark:E-additive}%
If \(\E\) is additive  then any additive functor \(\F\colon \Catp \to \E\) is group-like, because the inclusion $\Grp_{\Einf}(\E) \subseteq \Mon_{\Einf}(\E)$ is an equivalence. 
\end{remark}

\begin{examples}
\label{examples:additive}%
\
\begin{enumerate}
\item The functors \(\Core\) and \(\Poinc \colon \Catp \rightarrow \Sps\) are Verdier-localising since, by virtue of being representable, they preserve all limits. They are neither group-like nor Karoubi-localising. 
\item Lurie's \(\L\)-theory functor \(\L \colon \Catp \rightarrow \Spa\) \cite{Lurie-L-theory} is Verdier localising, see \S\reftwo{subsection:L+tate} below for a discussion.
\item
\label{item:from-catx}%
Many examples can be obtained from functors \(\Catx \rightarrow \E\) satisfying the corresponding conditions for (non-Poincaré) stable \(\infty\)-categories: For example, the functors \(\Kspace\colon \Catp \to \Sps\) and \(\K \colon \Catp \to \Spa\), which associate to a Poincaré \(\infty\)-category the algebraic \(\K\)-theory space or spectrum of its underlying stable \(\infty\)-category are Verdier-localising and group-like; this essentially follows from Waldhausen's additivity and fibration theorems, as implemented in the setting of stable \(\infty\)-categories by Blumberg-Gepner-Tabuada~\cite{BGT}, we will review the situation in \S\reftwo{subsection:Kadd}. Similarly, the functor \(\KK\colon \Catp \to \Spa\) which associates to a Poincaré \(\infty\)-category the nonconnective \(\K\)-theory spectrum of its underlying stable \(\infty\)-category is Karoubi-localising by~\cite{BGT}.
\item The functor \(\K \circ (-)^\natural \colon \Catp \to \Spa\) (where \((-)^{\natural}\) is the idempotent completion functor of Proposition~\reftwo{proposition:idempotent-completion}) is an example of an additive, but non-Verdier-localising functor. By contrast, the cofinality theorem implies that \(\Kspace\circ (-)^\natural \colon \Catp \rightarrow \Sps\) is Karoubi-localising. It is one of the main purposes of this paper series to show that these \(\K\)-theoretic examples have hermitian analogues.
\item The functors \(\mathrm H\K_i(-) \colon \Catp \to \Spa\) are further examples of additive functors that are not Verdier-localising and so is \(\mathrm{H}(\K_0(-^\natural)/\K_0(-))\), where \(\mathrm H \colon \Ab \rightarrow \Spa\) denotes the Eilenberg-Mac Lane functor. The analogous claims for Grothendieck-Witt and \(\L\)-groups are, however, not correct: For example, it follows from the results of Section~\reftwo{section:GW} that the functor \(\mathrm H\GW_0 \colon \Catp \rightarrow \Spa\) can only take all metabolic Poincaré-Verdier sequences to fibre sequences if the map \(\fgt \colon \GW_0(\C,\QF) \rightarrow \K_0(\C)\) is injective for all $\C$. This fails in many cases, e.g.\ for \((\Dperf(\mathbb Z),\QF^\s)\), where it is a map \(\ZZ \oplus \ZZ \rightarrow \ZZ\).
\end{enumerate}
\end{examples}

Finally, we will make use of the following: 

\begin{proposition}
\label{proposition:retract}%
Let 
\[
(\C,\QF) \st{f}{\lrar} (\D,\QFD) \st{p}{\lrar} (\E,\QFE)
\]
be a (split) Poincaré-Verdier sequence and let \(\F\colon \Catp \to \E\) be a group-like (additive or) Verdier-localising functor. Assume that Poincar\'e functor $p$ admits a section \(s\colon  (\E,\QFE) \to (\D,\QFD)\) in \(\Catp\). Then \(f\) and \(s\) together induce an equivalence
\[
\F(\C,\QF) \oplus \F(\E,\QFE) {\lrar} \F(\D,\QFD).
\]
Similarly, if the Poincaré functor \(f\) admits a retraction \(r\colon (\D,\QFD)\to (\C,\QF)\) in \(\Catp\) then \(p\) and \(r\) together induce an equivalence
\[
\F(\D,\QFD) {\lrar} \F(\C,\QF) \oplus \F(\E,\QFE).
\]
In case both $s$ and $r$ exist, the above equivalences are inverse to each other if \(r \circ s\) is the zero functor.
\end{proposition}

Note that since \(\Catp\) is only semi-additive, but not additive, the middle term in a Poincaré-Verdier sequence admitting a Poincaré splitting as above, need not split as a direct sum before applying \(\F\).

Proposition~\reftwo{proposition:retract} will be deduced from the following version of the splitting lemma:

\begin{lemma}
\label{lemma:retract}%
Let \(\A\) be an additive \(\infty\)-category which admits fibres and cofibres and let 
\[
\x \st{i}{\lrar} \y \st{r}{\lrar} \x
\]
be a retract diagram. Then the following statements hold:
\begin{enumerate}
\item
\label{item:retract-y}%
The maps \(i\colon\x \to \y\) and \(\fib(r) \to \y\) induce an equivalence \(\x \oplus \fib(r) \to \y\). 
\item
\label{item:retract-x}%
The maps \(r\colon\y \to \x\) and \(\y \to \cof(i)\) induce an equivalence \(\y {\to} \x \oplus \cof(i)\). 
\item
\label{item:retract-fiber-cofiber}%
The fibre sequence \(\fib(r) \to \y \to \x\) is also a cofibre sequence.
\item
\label{item:retract-cofiber-fiber}%
The cofibre sequence \(\x \to \y \to \cof(i)\) is also a fibre sequence.
\item
\label{item:retract-composite-eq}%
The composite map \(\fib(r) \to \y \to \cof(i)\) is an equivalence.
\end{enumerate}
\end{lemma}
\begin{proof}
We first note that \reftwoitem{item:retract-x} and \reftwoitem{item:retract-cofiber-fiber} follow from \reftwoitem{item:retract-y} and \reftwoitem{item:retract-fiber-cofiber}, respectively, applied to the additive \(\infty\)-category \(\A\op\). For the remaining cases \reftwoitem{item:retract-y}, \reftwoitem{item:retract-fiber-cofiber}, and \reftwoitem{item:retract-composite-eq} using the Yoneda lemma, we may assume that $\A = \Grp_{\Einf}(\Sps)$, in which case the long exact sequence of homotopy groups of a fibre sequence reduces the desired statements to the analoguous and well-known statements about split short exact sequences of abelian groups.
\end{proof}

\begin{proof}[Proof of Proposition \reftwo{proposition:retract}]
The first part follows from \reftwo{lemma:retract}~\reftwoitem{item:retract-y}, the second from \reftwo{lemma:retract}~\reftwoitem{item:retract-fiber-cofiber} and ~\reftwoitem{item:retract-x}.
To see the final statement, note that two equivalences are inverse to each other if and only if they are one-sided inverses. Composing in one direction, we get the functor
\[
\F(\C,\QF) \oplus \F(\E,\QFE) \to \F(\C,\QF) \oplus \F(\E,\QFE)
\]
whose ``matrix components'' are \(\left(\begin{matrix} \id & r_*s_* \\ 0 & \id \end{matrix}\right)\), and which is hence homotopic to the identity as soon as \(r\circ s\) is the zero functor.  
\end{proof}

\subsection{The classification of Poincar\'e-Verdier sequences}
The goal of this final section is to show that Poincar\'e-Verdier sequences, split Poincar\'e-Verdier sequences and Poincar\'e-Karoubi sequences with given kernel can each be classified by means of a universal example, analogous to the structure theory for Verdier sequences developed in~\cite{Nikolaus-Verdier} and described in ~\S\reftwo{section:classverdier}.

As in \S\reftwo{section:appendix-verdier}, it will be convenient (and partially necessary) to extend the scope from small Poincaré (and hermitian) $\infty$-categories to locally small ones. Here, by a locally small hermitian $\infty$-category we simply mean a locally small $\infty$-category $\C$ equipped with an hermitian structure $\QF\colon \C\op \to \Sp$, and such a locally small hermitian $\infty$-category is Poincaré if \(\Bil_{\QF}(-,-) \simeq \map_{\C}(-,\Dual-)\) for some equivalence \(\Dual\colon \C\op \rightarrow \C\). To get started, note that the (Verdier) quotient of a locally small Poincaré $\infty$-category $(\D,\QFD)$ by a small full Poincaré subcategory $(\C,\QFD|_{\C})$ is again a locally small Poincaré $\infty$-category, with the Poincaré structure on the quotient given by the formula from Remark~\reftwo{remark:compute}, which is well-defined since the $\C$ is small.
By a locally small Poincaré-Verdier sequence we will then mean a null-composite sequence 
\[
(\C,\QF) \lrar (\D,\QFD) \lrar (\E,\QFE)
\]
of locally small Poincar\'e $\infty$-categories, with $\C$ small, satisfying the characterisation from Proposition \reftwo{proposition:criterion-poincare} (or equivalently being both a fibre and cofibre sequence in the evident $\infty$-category of locally small Poincar\'e categories). We then use the term locally small Poincaré-Verdier projections/inclusions for the projections and inclusions appearing in such sequences (so that by definition a locally small Poincaré-Verdier inclusion has a small domain and a locally small Poincaré-Verdier projection a small kernel). Recall also from \S\reftwo{section:appendix-verdier} that the notion of adjointability extends to squares of locally small stable $\infty$-categories with a pair of parallel legs being either locally small Verdier inclusions or projections. We denote by $\textsc{PVer}$ the ensuing enlargement of $\mathrm{Pver}$ consisting of locally small Poincaré-Verdier sequences, with morphisms defined as commutative diagrams made up of adjointable squares. With these definitions Proposition~\reftwo{proposition:poincare-verdier-square-char} characterising Poincar\'e-Verdier squares applies verbatim to locally small Poincaré $\infty$-categories so long as the vertical fibres are assumed small. 

Now recall from Theorem~\reftwo{theorem:universal-verdier} that to an \(\infty\)-category \(\cB\) we may associate a Verdier sequence
\[
\cB^\natural \lrar \Latt(\cB) \lrar \Tate(\cB),
\]
which is the universal Vedier sequence with fibre \(\cB^\natural\) in the sense that morphisms of Verdier sequences from any Verdier sequence \(\C \to \D \to \E\) to it correspond to exact functors \(\C \to \cB^\natural\). 
The \(\infty\)-category \(\Latt(\cB)\) is by definition the full subcategory of \(\Pro\Ind(\cB)\) spanned by those arrows whose source is in \(\Ind(\cB)\) and target in \(\Pro(\cB)\). The Verdier inclusion \(\cB^\natural \to \Latt(\cB)\) sends \(b\) to the identity arrow \(\id\colon b \to b\) (where \(b\) is considered as an object of both \(\Ind(\cB)\) and \(\Pro(\cB)\) via the Yoneda embeddings) and \(\Tate(\cB)\) can be identified 
with the smallest stable subcategory of \(\Pro\Ind(\C)\) containing both \(\Pro(\cB)\) and \(\Ind(\cB)\) with the Verdier projection $\Latt(\cB) \rightarrow \Tate(\cB)$ taking cofibres. We refer the reader to the end of \S\reftwo{section:appendix-verdier} for more details.

To equip $\Latt(\cB)$ with a Poincar\'e structure we note that there is a canonical equivalence 
\[
\Latt(\cB) \simeq \Pairings(\Ind(\cB),\Pro(\cB),\wtl{\map}_{\cB}),
\]
where 
$\wtl{\map}_\cB \colon \Ind(\cB)\op \times \Pro(\cB) \to \Sp$
is the essentially unique extension of $\map_{\cB}$ which preserves limits in each variable, and for a biexact functor $\Bil\colon \C\op \times \D \to \Sp$ we write $\Pairings(\C,\D,\Bil)$ for the bivariant unstraightening of the space valued functor $\Om^{\infty}\Bil$, see Remark~\reftwo{remark:concrete}. Now observe that for an hermitian \(\infty\)-category \((\cB, \Pi)\), \(\Ind(\cB)\) canonically inherits an hermitian structure 
\[
\widetilde \Pi\colon \Ind(\cB)\op\simeq \Pro(\cB\op)\xrightarrow{\Pro(\Pi)} \Pro(\Spa)\xrightarrow \lim \Spa,
\]
so that \((\Ind(\cB),\widetilde \Pi)\) is a locally small hermitian $\infty$-category.
Note that \((\Ind(\cB), \widetilde \Pi)\) is usually not Poincaré, even if \((\cB, \Pi)\) is: Instead, the duality functor \(\Dual_\Pi\colon \cB\op\to \cB\) extends to an equivalence 
\[
\wtl{\Dual}_{\Pi}\colon \Ind(\cB)\op\simeq \Pro(\cB\op) \lrar \Pro(\cB).
\]
giving rise to an identification 
\[
\Latt(\B)\simeq \Pairings(\Ind(\cB),\Pro(\cB),\widetilde{\hom}_\cB) \simeq \Pairings(\Ind(\cB),\Ind(\cB)\op,\Omega^\infty\Bil_{\widetilde \Pi}) 
\]
for Poincar\'e $(\cB,\Pi)$. 
Now the Poincar\'e refinement of the pairing construction 
\[
\Pairings(\D,\QFD) = (\Pairings(\D,\D\op,\Omega^\infty\Bil_\QFD),\QFD^\pair)
\]
 from \S\refone{subsection:thom} works equally well for locally small hermitian $\infty$-categories $(\D,\QFD)$, as does Proposition~\refone{proposition:first-adj-hermitian}, showing that it provides a right adjoint to the forgetful functor from locally small Poincaré to locally small hermitian $\infty$-categories. We thus put: %

\begin{definition}
For a Poincaré \(\infty\)-category \((\cB, \Pi)\) we define $
\Latt(\cB, \Pi) \coloneq \Pairings(\Ind(\cB), \widetilde \Pi)$.%
\end{definition}

In particular, \(\Latt(\cB,\Pi)\) is a locally small Poincaré \(\infty\)-category and its underlying stable \(\infty\)-category is indeed equivalent to \(\Latt(\cB)\).%
We shall unpack this definition in Proposition \reftwo{proposition:TateandLattarewhattheyaresupposedtobe} below.
For now the universal property of the pairing construction implies that the hermitian functor \((\cB, \Pi) \to (\Ind(\cB), \widetilde \Pi)\) lifts uniquely to a Poincaré functor (marked as dashed)
\[
\begin{tikzcd}
(\cB,\Pi)\ar[r,equal]\ar[d,"i",dashed] & (\cB,\Pi) \ar[d, "j"] \\
\Latt(\cB,\Pi) \ar[r] & (\Ind(\cB),\widetilde \Pi).
\end{tikzcd}
\]
In addition, the Poincaré functor \(i\)  
is a Poincaré-Verdier inclusion since it can be written as a composite of Poincaré-Verdier inclusions
\[
(\cB,\Pi) \lrar \Ar(\cB,\Pi) = \Pairings(\cB,\Pi) \lrar \Pairings(\Ind(\cB),\widetilde\Pi) = \Latt(\cB,\Pi),
\]
see Examples~\refone{examples:typical-pairings}\;\refoneitem{item:correspondence-pairing} for the equivalence between the arrow and pairing constructions for \((\cB,\Pi)\) Poincaré.

\begin{definition}
\label{definition:tatepoinc}%
For a Poincaré \(\infty\)-category \((\cB,\Pi)\) we put $\Tate(\cB,\Pi) = \Latt(\cB,\Pi\qshift{-1})/(\cB,\Pi\qshift{-1})$.
\end{definition}

In particular, $\Tate(\cB,\Pi)$ is again a locally small Poincar\'e $\infty$-category for every $(\cB,\Pi) \in \Catp$, and in total, we thus obtain a locally small Poincaré-Verdier sequence
\[
(\cB, \Pi)^\natural \lrar  \Latt(\cB, \Pi) \lrar  \Tate(\cB, \Pi\qshift{1})
\]
refining the universal Verdier sequence with kernel \(\cB^\natural\) from Theorem~\reftwo{theorem:universal-verdier}. 

\begin{remark}
\label{remark:explicit-q-i}%
Under the identification of \(\Latt(\cB)\) with the full subcategory of \(\Ar(\Ind\Pro(\cB))\) spanned by arrows from inducitve to projective systems, the counit \(q \colon \Latt(\cB) \rightarrow \Ind(\cB)\) sends an arrow \(x \to y\) to \(x \in \Ind(\cB)\), while \(i \colon \cB \rightarrow \Latt(\cB)\) sends \(b\) to the arrow \(\id\colon b \to b\). In particular, \(q\) induces an equivalence \(\hom_{\Latt(\cB)}(i(b),z) \to \hom_{\Ind(\cB)}(qi(b)),q(z))\) for all $b \in \cB$ and $z \in \Latt(\cB)$. By Remark~\reftwo{remark:adj-concrete} the square above is vertically inductively adjointable.
\end{remark}

Let us now unpack the construction of the Poincar\'e structure in terms of the original definition of $\Latt(\cB) \subseteq \Ar(\Ind\Pro(\cB))$:

\begin{proposition}
\label{proposition:TateandLattarewhattheyaresupposedtobe}%
For $(\cB,\Pi)$ a Poincar\'e $\infty$-category, the Poincar\'e structure on $\Tate(\cB,\Pi)$ is given by
\[
\Pi_{\mathrm{tate}}\colon \Tate(\cB)\op \subseteq \Ind(\Pro(\cB))\op = \Pro\Ind(\cB\op) \xrightarrow{\Pro\Ind\Pi} \Pro\Ind(\Sp) \xrightarrow{\colim} \Pro(\Sp) \xrightarrow{\lim} \Sp,
\]
and that on $\Latt(\cB,\Pi)$ by
\[
\Pi_{\mathrm{latt}} \colon \Latt(\cB)\op \subseteq \Ar(\Tate(\cB))\op \xrightarrow{(\Pi_\mathrm{tate})_\arr} \Spa.
\]
\end{proposition}

\begin{proof}
Under the equivalence 
\[
\Pairings(\Ind(\cB),\Pro(\cB),\widetilde{\hom}) \simeq \Ar(\Ind\Pro(\cB)) \times_{\Ind\Pro(\cB) \times \Ind\Pro(\cB)}[\Ind(\cB) \times \Pro(\cB)]
\]
of Remark~\reftwo{remark:concrete}, the Poincaré structure $\widetilde{\Pi}_{\mathrm{latt}}$ is per construction given by the formula
\[
\widetilde{\Pi}_{\mathrm{latt}}(I \to P) = \wtl{\Pi}(I) \times_{\Bil_{\wtl{\Pi}}(I,I)}\map_{\Pro(\cB)}(P,\wtl{\Dual}_{\Pi}(I)).
\]
Denoting by $\widehat \Pi \colon \Tate(\cB)\op \rightarrow \Sp$ the hermitian structure from the statement which we want to show is $\Pi_\mathrm{tate}$, this formula agrees, by direct inspection, with the restriction of $\widehat\Pi_\arr \colon \Ar(\Tate(\cB))\op \rightarrow \Spa$ to $\Latt(\cB)\op$. To deduce that $\Pi_\mathrm{tate} \simeq \widehat \Pi$ consider the diagram
\[
\begin{tikzcd} 
(\cB,\Pi)^\natural \ar[r]\ar[d] & \Latt(\cB,\Pi) \ar[r,"\cof"]\ar[d] & \Tate(\cB,\Pi\qshift{1}) \ar[d]\\
(\Tate(\cB),\widehat \Pi) \ar[r] & \Ar(\Tate(\cB), \widehat\Pi) \ar[r,"\cof"] & (\Tate(\cB),\widehat\Pi\qshift{1})
\end{tikzcd}
\]
of horizontal Verdier sequences issuing from the discussion so far. As the left Kan extensions of the hermitian structures along the right hand horizontal maps can be computed by restriction along the $\Pro$-left adjoints, it suffices to show that the right hand square is adjointable. By Lemma \reftwo{lemma:adj-verdier} we may as well check that the left hand square is adjointable, which is true by inspection, as the horizontal (projective/inductive) adjoints simply take an arrow to its target or source.
\end{proof}

The following theorem then establishes the Poincar\'e-Verdier sequence formed by the lattice and Tate objects as a universal one:

\begin{theorem}
\label{theorem:universal-poincare-verdier}%
For every locally small Poincaré-Verdier sequence \((\C,\QF) \to (\D,\QFD) \to (\E,\QFE)\) and every (small) Poincar\'e $\infty$-category \((\cB,\Pi)\) the map
\[
\Fun^{\textsc{Ver}}((\D,\QFD) \to (\E,\QFE),\Latt(\cB, \Pi) \to \Tate(\cB, \Pi\qshift{1})) \longrightarrow \Funx((\C,\QF),(\cB,\Pi)^\natural)
\]
extracting kernels is an equivalence of (small) Poincar\'e $\infty$-categories. Consequently, there 
exists an essentially unique adjointable square
\[
\begin{tikzcd}
(\D,\QFD) \ar[r,"\vphi"]\ar[d, "p"'] & \Latt(\C,\QF) \ar[d] \\
(\E,\QFE) \ar[r,"\psi"] & \Tate(\C,\QF\qshift{1}) \ ,
\end{tikzcd}
\]
inducing the inclusion $(\C,\QF) \rightarrow (\C,\QF)^\natural$ on vertical fibres. This square is cartesian if and only if $\C$ is idempotent complete, so that in this case the Poincaré-Verdier projection \(p\) is pulled back from the universal Poincaré-Verdier projection with fibre $(\C,\QF)$ on the right. 
\end{theorem}

The Poincaré functor \(\psi\colon (\E,\QFE) \to \Tate(\C,\QF\qshift{1})\) is called the \defi{classifying functor} of the given Poincaré-Verdier sequence, and the square its \defi{classifying square}. Theorem~\reftwo{theorem:universal-poincare-verdier} then implies the following classification result for Poincaré-Verdier projections (among small Poincar\'e $\infty$-categories), for the statement of which we note that $\K(\Tate(\C,\QF\qshift{1})) \simeq \mathbb S^1 \otimes \K((\C,\QF)^\natural)$ in $\Spa^\hC$ since $\Latt(\C)$ admits an Eilenberg swindle, compare \cite{Hennion-Tate}*{Corollary 4.3}:
	
\begin{corollary}
Given Poincaré $\infty$-categories $(\C,\QF)$ and $(\E,\QFE)$, extracting classifying functors and the boundary map of $\K$-spectra induces an equivalence between the space of Poincaré-Verdier sequences $(\C,\QF) \rightarrow (\D,\QFD) \rightarrow (\E,\QFE)$ and the space of pairs $(\psi, \eta)$, where $\psi \colon (\E,\QFE) \rightarrow \Tate(\C,\QF\qshift{1})$ and $\eta \colon \K(\E,\QFE) \rightarrow \mathbb S^1 \otimes \K(\C,\QF)$ is a lift of
\[
\K(\E,\QFE) \xrightarrow{\psi_\ast} \K(\Tate(\C,\QF\qshift{1})) \xrightarrow{\partial} \mathbb S^1 \otimes \K((\C,\QF)^\natural)
\]
as maps in $\Spa^\hC$. Equivalently, it is the space of pairs $(\psi,s)$, where $\psi$ is as above, such that in the exact sequence%
\[
\K_0((\C,\QF)^\natural)/\K_0(\C,\QF) \longrightarrow \K_0\big(\Latt(\C,\QF) \times_{\Tate(\C,\QF\qshift{1})} (\E,\QFE)\big)/\K_0(\C,\QF) \longrightarrow \K_0(\E,\QFE) \longrightarrow 0
\]
the first map is injective and $s$ is a splitting of it as $\mathbb Z[\Ct]$-modules.
\end{corollary}

The deduction of this Corollary from Theorem~\reftwo{theorem:universal-poincare-verdier} is essentially the same as that of Corollary~\reftwo{corollary:classverdierfull}, so we refrain from repeating it here.  The proof of Theorem~\reftwo{theorem:universal-poincare-verdier} will use the following lemma. 

\begin{lemma}
\label{lemma:kan-extended}%
Let \((\D,\QFD)\) be a locally small hermitian \(\infty\)-category, \(\C \subseteq \D\) a small full subcategory and \((\E,\QFE)\) a locally small hermitian \(\infty\)-category with $\E$ admitting filtered colimits and \(\QFE\) preserving cofiltered limits. Denote then by \(\A \subseteq \Funx(\D,\E)\) be the full subcategory spanned by those exact funtors \(\D \to \E\) which are left Kan extended from their restriction to \(\C\). Then the composite hermitian functor
\[
(\A,{\Nat^{\QFD}_{\QFE}}|_{\A}) \longrightarrow \Funx((\D,\QFD),(\E,\QFE)) \longrightarrow \Funx((\C,\QFD|_{\C}),(\E,\QFE))
\]
is an equivalence of locally small hermitian \(\infty\)-categories.
\end{lemma}

\begin{proof}
Note that left Kan extension along \(\C \subseteq \D\) is well defined since $\C$ is small, and sends exact functors \(\C \to \E\) to exact functors \(\D \to \E\), see the first part of Appendix~\reftwo{appendix:AppIIA}.  
The composite 
\[
\A \longrightarrow \Funx(\D,\E) \longrightarrow \Funx(\C,\E)
\]
is thus an equivalence by~\cite{HTT}*{Proposition 4.3.2.15}. It is left to show that if \(f\colon \D \to \E\) is left Kan extended from \(\C\) then 
\[
\Nat(\QFD,f^*\QFE) \longrightarrow \Nat(\QFD|_{\C\op},(f^*\QFE)|_{\C\op})
\]
is an equivalence of spectra. Indeed, noting that the comma \(\infty\)-categories of the inclusion \(\C \subseteq \D\) are filtered, the assumptions that \(f\) is left Kan extended from \(f|_{\C}\) and that \(\QFE\) sends filtered colimits to limits together imply that \(f^*\QFE\) is right Kan extended from its restriction to \(\C\op\), and so the desired result follows.
\end{proof}

\begin{proof}[Proof of Theorem~\reftwo{theorem:universal-poincare-verdier}]
We will start by showing that the extraction of kernels gives an equivalence
\[
\Hom_{\textsc{PVer}}((\D,\QFD) \to (\E,\QFE),\Latt(\cB, \Pi) \to \Tate(\cB, \Pi\qshift{1})) \longrightarrow \Hom_{\Catp}((\C,\QF),(\cB,\Pi)^\natural).
\]
Equivalently, given Poincaré functor \((f,\eta)\colon (\C,\QF) \to (\cB,\Pi)^\natural\) we have to show that the fibre over \((f,\eta)\) is contractible.
By Lemma~\reftwo{lemma:adj-verdier}, this fibre can be described as the space of extensions of \((f,\eta)\) to commutative squares
\[
\begin{tikzcd}
(\C,\QF) \ar[r, "{(f,\eta)}"]\ar[d] & (\cB,\Pi)^\natural \ar[d, "i"] \\
(\D,\QFD) \ar[r, dashed] & \Latt(\cB,\Pi)
\end{tikzcd}
\]
whose underlying square is adjointable (where all arrows are fixed except the dashed one). %
Investing the universal property of the pairing construction (Proposition~\refone{proposition:first-adj-hermitian}) such squares correspond to squares
\[
\begin{tikzcd}
(\C,\QF) \ar[r, "{(f,\eta)}"]\ar[d] & (\cB,\Pi)^\natural \ar[d, "j"] \\
(\D,\QFD) \ar[r, "{(g,\theta)}", dashed] & (\Ind(\cB),\wtl{\Pi})
\end{tikzcd}
\]
where \((g,\theta)\) is now an hermitian functor. By Remark~\reftwo{remark:explicit-q-i} the vertical inductive adjointability of this latter square is equivalent to the (vertical inductive) adjointability of the former.
In this last square the underlying right vertical exact functor is just the canonical embedding, and hence by Remark~\reftwo{remark:adj-concrete} \reftwoitem{item:equivalent-def} (in the universal case of $\vphi$ being the embedding $j\colon \cB \to \Ind(\cB)$)
its vertical inductive adjointability is equivalent to the condition that \(g^*j_!(j) = g\colon \D \to \Ind(\cB)\) is left Kan extended from its restriction to \(\C\). We now obtain the claim by invoking Lemma~\reftwo{lemma:kan-extended}. The addendum already follows from this by Proposition \reftwo{proposition:poincare-verdier-square-char}.

The full statement about hermitian functor categories can be deduced as follows: On underlying stable $\infty$-categories it is an instance of Theorem \reftwo{theorem:universal-verdier}, and to check the Poincar\'e structures it then suffices to show that the functor in question induces an equivalence on (a priori not necessarily small) spaces of shifted hermitian forms on account of the general cartesian squares
\[
\begin{tikzcd}\Omega^{\infty-i} \mathrm P(X) \ar[r] \ar[d] & \spsforms(\mathcal A,\mathrm P\qshift{i}) \ar[d] \\
\{X\} \ar[r] & \grpcr(\mathcal A)\end{tikzcd}
\]
for any hermitian $\infty$-category $(\mathcal A,\mathrm P)$. But we have 
\[
\spsforms(-\qshift{i}) = \Hom_{\textsc{Cat}^\mathrm h_\infty}((\Sp^\omega,\QF^\uni)\qshift{-i},-) = \Hom_{\textsc{Cat}^\mathrm p_\infty}(\Met(\Sp^\omega,\QF^\uni)\qshift{-i},-)
\]
by Propositions \refone{proposition:corepresentability-of-poinc}, \refone{proposition:second-adj-hermitian} and Remark \refone{remark:second-counit}. 
Specialising to the situation at hand, we can further use the adjunction equivalences
\[
\Hom_{{\textsc{Cat}_\infty^\mathrm p}}((\Met(\Sp^\omega,\QF^\uni)\qshift{-i},\Funx(-,-)) \simeq \Hom_{{\textsc{Cat}_\infty^\mathrm p}}(\Met(\Sp^\omega,\QF^\uni)\qshift{-i} \otimes -, -)
\]
from Corollary \refone{proposition:internal} to rewrite the requisite statement as an instance of the claim already proven applied to the Poincar\'e-Verdier sequence
\[
[\Met(\Sp^\omega,\QF^\uni)\qshift{-i} \otimes (\C,\QF)]^r \longrightarrow \Met(\Sp^\omega,\QF^\uni)\qshift{-i} \otimes (\D,\QFD)\longrightarrow \Met(\Sp^\omega,\QF^\uni)\qshift{-i} \otimes (\E,\QFE),
\]
whose source is the retract closure of $\Met(\Sp^\omega,\QF^\uni)\qshift{-i} \otimes (\C,\QF)$, see Lemma \reftwo{lemma:Morita-exact-criterion}.
\end{proof}

\begin{remark}
\label{remark:explicit-classifying-hermitian}%
To give a somewhat more concrete description of the classifying functor \(\psi\), let us write \(r\colon \C \to \Ind(\D)\) and \(q \colon \D \to \Pro(\C)\) for the inductive right and projective left adjoints of \(f \colon \C \rightarrow \D\). 
Unwinding the proof of Theorem~\reftwo{theorem:universal-poincare-verdier} we see that 
the Poincaré functor \(\vphi\colon (\D,\QFD) \to \Latt(\C,\QF) = \Pairings(\Ind(\C),\widetilde \QF)\) is the essentially unique Poincaré functor determined by the canonical hermitian refinement \((r,\theta)\colon (\D,\QFD) \to (\Ind(\C),\widetilde \QF)\) of \(r\) 
given by 
\[
\eta\colon \QFD(X) \longrightarrow \lim_{Y \in \C_{/X}}\QFD(f(Y)) = \lim_{Y \in \C_{/X}}\QF(Y) =  \wtl{\QF}(r(X)).
\]
Under the identification of \(\Latt(\C)\) with the full subcategory of \(\Pro\Ind(\C)\) spanned by the arrows between inductive and projective diagrams, the underlying exact functor of \(\vphi\) then sends \(d\) to the essentially unique arrow \(r(X) \to q(X)\) mapping under \(\Pro\Ind(f)\) to the unit-counit composite map \(fr(X) \to X \to fq(X)\). 
Similarly, we may identify \(\Tate(\C)\) with the smallest stable subcategory of \(\Pro\Ind(\C)\) containing \(\Ind(\C)\) and \(\Pro(\C)\), in which case the underlying exact functor of \(\psi\) sends \(p(X) \in \E\) to \(\cof[r(X) \to q(X)] \in \Tate(\C)\).
\end{remark}

Next, we carry the classification over to split Poincar\'e-Verdier sequences. This is now easy by forming split cores, and we arrive at an hermitian analogue of the structure theory of split Verdier sequences described in \S\reftwo{section:appendix-split-verdier}. The role of the universal split Verdier sequence is played by the metabolic Poincaré-Verdier sequence above.

To this end note that for $\cB \in \Catx$ idempotent complete the split core of 
\(\cB \to \Latt(\cB) \to \Tate(\cB)\) is given by the full subcategories \(\Ar(\cB) \subseteq \Latt(\cB)\) and $\cB \subseteq \Tate(\cB)$, see the discussion towards the end of \S\reftwo{section:appendix-split-verdier}, and hence on the level of Poincaré \(\infty\)-categories the split core is given by the Poincaré full subcategory 
\[
\Ar(\cB,\Pi) = \Pairings(\cB,\Pi) \subseteq \Pairings(\Ind(\cB),\wtl{\Pi}) = \Latt(\cB,\Pi)
\]
mapping to $(\cB,\Pi\qshift{1}) \subseteq \Tate(\cB,\Pi\qshift{1})$. In total, the split core of the universal Verdier sequence then coincides with
the metabolic sequence
\[
(\cB,\Pi) \longrightarrow \Ar(\cB,\Pi) \longrightarrow (\cB,\Pi\qshift{1})
\]
of Example~\reftwo{example:metabolicfseq}.
This construction gives a universal split Poincaré-Verdier sequence, even when \(\cB\) is not necessarily idempotent complete:

\begin{corollary}
\label{corollary:universal-split-poincare-verdier}%
\label{corollary:structure-square}%
Let \((\cB,\Pi)\) be a Poincaré \(\infty\)-category. Then,
for every split Poincaré-Verdier sequence
\((\C,\QF) \to (\D,\QFD) \to (\E,\QFE)\)
the map
\[
\Fun^{\Ver}((\D,\QFD) \to (\E,\QFE),\Ar(\cB,\Pi) \to (\cB,\Pi\qshift{1})) \longrightarrow \Funx((\C,\QF),(\cB,\Pi))
\]
extracting the induced functor on kernels is an equivalence. Consequently, there exist essentially unique (equivalent) adjointable squares
\[
\begin{tikzcd}
(\D,\QFD) \ar[r,"\vphi"]\ar[d, "p"'] & \Ar(\C,\QF) \ar[d, "\cof"] && (\D,\QFD) \ar[r,"\vphi'"]\ar[d, "p"'] & \Met(\C,\QF\qshift{1}) \ar[d,"\met"] \\
(\E,\QFE) \ar[r,"\psi"] & (\C,\QF\qshift{1}) && (\E,\QFE) \ar[r,"\psi"] & (\C,\QF\qshift{1})
\end{tikzcd}
\]
inducing the identity on vertical fibres. In addition, both are cartesian, so that the split Poincaré-Verdier projection \((\D,\QFD) \to (\E,\QFE)\) is pulled back from the universal split Poincaré-Verdier projections \(\cof \colon \Ar(\C,\QF) \to (\C,\QF\qshift{1})\) or \(\met \colon \Met(\C,\QF\qshift{1}) \rightarrow (\C,\QF\qshift{1})\).
\end{corollary}

\begin{proof}
In the case where \(\cB\) is idempotent complete this follows from Theorem~\reftwo{theorem:universal-poincare-verdier} via Lemma~\reftwo{lemma:split-poincare-core}. For arbitrary \(\cB\), it suffices by Lemma~\reftwo{lemma:adj-verdier} to note that in any adjointable square
\[
\begin{tikzcd}
\C \ar[r, "\rho"]\ar[d] & \cB^{\natural} \ar[d,"\id_{(-)}"] \\
\D \ar[r, "\vphi"] & \Ar(\cB)^{\natural} \ ,
\end{tikzcd}
\]
the condition that \(\rho\) factors through \(\cB \subseteq \cB^{\natural}\) is equivalent to the condition that \(\vphi\) factors through \(\Ar(\cB)\). The final claim follows from Proposition~\reftwo{proposition:poincare-verdier-square-char}.
\end{proof}

The metabolic fibre sequence 
\[
(\C,\QF) \longrightarrow \Ar(\C,\QF) \longrightarrow (\C,\QF\qshift{1}) \quad \text{or equivalently} \quad (\C,\QF) \longrightarrow \Met(\C,\QF\qshift{1}) \longrightarrow (\C,\QF\qshift{1})
\]
is thus the universal split Poincaré-Verdier sequence with kernel $(\C,\QF)$. Note that, in contrast to the non-split case it consists entirely of small Poincar\'e $\infty$-categories.%
 
The Poincaré functor \(\psi\colon (\E,\QFE) \to (\C,\QF\qshift{1})\) is then again called the \defi{classifying functor} of the given split Poincaré-Verdier sequence, and the square the associated classifying square.
Just as in the non-split case Corollary~\reftwo{corollary:universal-split-poincare-verdier} in total implies the following classification result for split Poincaré-Verdier projections:

\begin{corollary}
Given two Poincar\'e $\infty$-categories $(\C,\QF)$ and $(\E,\QFE)$ extracting classifying functors induces an equivalence between the space of Poincar\'e-Verdier sequences $(\C,\QF) \rightarrow (\D,\QFD) \rightarrow (\E,\QFE)$ and $\Hom_{\Catp}((\E,\QFE), (\C,\QF\qshift{1}))$.
\end{corollary}%

\begin{remark}
\label{remark:explicit-classifying-split}%
\label{remark:explicit-met}%
The descriptions of \(\vphi\) and \(\psi\) given in Remark~\reftwo{remark:explicit-classifying-hermitian}
simplify in the split case. In particular, writing \(r,q\colon \D \to \C\) for the right and left adjoints of \(f\), the Poincaré functor \(\vphi\colon (\D,\QFD) \to \Ar(\C,\QF) = \Pairings(\C,\QF)\) is the one determined by the canonical hermitian refinement \((r,\tau)\colon (\D,\QFD) \to (\C,\QF)\) of \(r\) (see Remark~\reftwo{remark:hermitian-adjoints}).
The underlying exact functor of \(\vphi\) sends \(d \in \D\) to the unit-counit composite arrow \(r(d) \to q(d)\). The Poincaré structure is given by 
\[
\QFD(d) \to \QFD(fr(d)) \displaystyle\mathop{\times}_{\Bil_{\QFD}(fr(d),fr(d))} \Bil_{\QFD}(fr(d),fq(d))  = \QF(r(d)) \displaystyle\mathop{\times}_{\Bil_{\QF}(r(d),r(d))} \Bil_{\QF}(r(d),q(d)) = \QF_{\arr}([r(d) \to q(d)])
\]
where the first map is given by restriction along the counit $fr(d) \to d$ and by the composite 
\[
\QFD(d) \to \Bil_{\QFD}(d,d) = \map_\D(d,\Dual_\Phi d) \to \map_\D(fr(d),fr(\Dual_\Phi d)) = \map_\D(fr(d),\Dual_\Phi fq(d)) = \Bil_{\QFD}(fr(d),fq(d))
\]
using the fact that conjugating with the duality switches between the two adjoints of $f$.
Furthermore, if we write \(g,h\colon \E \to \D\) for the left and right adjoints of \(p\), respectively, then \(rh = qg = 0\) and so the underlying exact functor of \(\psi\) is given by \(\psi(e) = \cof\vphi(g(e)) = \cof\vphi(h(e)) = qh(e) = \Sig rg(e)\) with the Poincaré structure induced by that of $\vphi$.

Upon replacing $\Ar(\C,\QF)$ by $\Met(\C,\QF\qshift{1})$ %
by means of the equivalence $\Ar(\C,\QF) \to\Met(\C,\QF\qshift{1})$ taking $f\colon X \to Y$ to $Y \to \cof(f)$ the structure square of Corollary~\reftwo{corollary:structure-square} then has the %
underlying exact functor of $\vphi'$ given by sending $X \in \D$ to 
\[
[q(X) \to \cof[r(X)\to q(X)]] = [q(X) \to q(\cof[r(X) \to X])] = [q(X) \to qhp(X)],
\]
where $r,q$ are the right and left adjoint of $f$ and $h,g$ are the right and left adjoint of $p$, respectively. 
\end{remark}

Finally, we record the classification of Poincaré-Karoubi sequences, as another direct consequence of Theorem \reftwo{theorem:universal-poincare-verdier} combined with the universal property of idempotent completions. To state the result let $\Fun^\mathrm{Kar}(-,-)$ denote the evident generalisation of $\Fun^\Ver(-,-)$ allowing Poincar\'e-Karoubi sequences as input. We then have:

\begin{corollary}
\label{corollary:universal-poincare-karoubi}%
For $(\cB,\Pi)$ a Poincaré $\infty$-category, the restriction map
\[
\Fun^{\textsc{Kar}}((\C,\QF) \to (\D,\QFD) \to (\E,\QFE), (\cB,\Pi) \to \Latt(\cB,\Pi) \to \Tate(\cB, \Pi\qshift{1})^\natural) \longrightarrow \Funx((\C,\QF),(\cB,\Pi))
\]
is an equivalence for every locally small Poincaré Karoubi sequence \((\C,\QF) \to (\D,\QFD) \to (\E,\QFE)\), and there is consequently an essentially unique adjointable square
\[
\begin{tikzcd}
(\D,\QFD) \ar[r,"\vphi"]\ar[d, "p"'] & \Latt(\C,\QF) \ar[d] \\
(\E,\QFE) \ar[r,"\psi"] & \Tate(\C,\QF\qshift{1})^\natural \ ,
\end{tikzcd}
\]
inducing the inclusion $(\C,\QF) \rightarrow (\C,\QF)^\natural$ on vertical fibres. This square is cartesian if and only if $\C$ is idempotent complete, so that in this case the Poincaré-Karoubi projection \(p\) is pulled back from the universal Poincaré-Karoubi projection with fibre $(\C,\QF)^\natural$ on the right. 
\end{corollary}

\begin{corollary}
Given Poincar\'e $\infty$-categories $(\C,\QF)$ and $(\E,\Psi)$ extracting classifying functors and the boundary map of $\K$-spectra induces an equivalence between the space of Poincar\'e-Karoubi sequences $(\C,\QF) \rightarrow (\D,\QFD) \rightarrow (\E,\Psi)$ and the space of quadruples $(\psi, \E', \C', \eta)$ where $\psi \colon (\E,\QFE) \rightarrow \Tate(\C,\QF\qshift{1})^\natural$ is a Poincar\'e functor, $\E^\mathrm{min} \subseteq \E' \subseteq \E$ and $\C \subseteq \C' \subseteq \C^\natural$ are closed under the dualities and $\eta \colon \K(\E',\QFE) \rightarrow \mathbb S^1 \otimes \K(\C',\QF)$ is a lift of $\psi_\ast \colon \K(\E',\QFE) \rightarrow \K(\E,\QFE) \rightarrow \K(\Tate(\C,\QF\qshift{1})^\natural)$ along
\[
\mathbb S^1 \otimes \K(\C',\QF) \longrightarrow \mathbb S^1 \otimes \K((\C,\QF)^\natural) \xrightarrow{\partial} \K(\Tate(\C,\QF\qshift{1})) \longrightarrow \K(\Tate(\C,\QF\qshift{1})^\natural)
\]
as maps in $\Spa^\hC$. Equivalently, it is the space of quadruples $(\psi,e,c,s)$ where, 
\[
e \subseteq \mathrm{ker}(\psi_\ast \colon \K_0(\E) \longrightarrow \K_{0}(\Tate(\C)^\natural)) \quad \text{and} \quad \K_0(\C) \subseteq c \subseteq \K_0(\C^\natural)
\]
are subgroups closed under the duality, such that 
in the exact sequence
\[
\K_0((\C,\QF)^\natural)/c \longrightarrow \K_0\big(\Latt(\C,\QF) \times_{\Tate(\C,\QF\qshift{1})} (\E^e,\QFE)\big)/c \longrightarrow e \longrightarrow 0
\]
the left hand map is injective and $s$ is a splitting of it as $\mathbb Z[\Ct]$-modules, where $\E^e = \{X \in \E \mid [X] \in e\}$.
\end{corollary}

In particular, the space of Poincar\'e-Karoubi sequences $(\C,\QF) \rightarrow (\D,\QFD) \rightarrow (\E,\QFE)$ among idempotent complete (or minimal) Poincar\'e $\infty$-categories agrees with $\Hom_{\Catp}((\E,\QFE),\Tate(\C,\QF\qshift{1})^\natural)$.

\begin{proof}[Proof of Corollary~\reftwo{corollary:universal-poincare-karoubi}]
By Corollary~\reftwo{corollary:criterion-poincare-karoubi-projection} both $(\C,\QF) \xrightarrow{f} (\D,\QFD) \xrightarrow{p} (\E,\QFE)$ and $\ker(p) \rightarrow (\D,\QFD) \rightarrow \mathrm{im}(p)$ admit canonical pointwise dense inclusions into $\ker(p) \rightarrow (\D,\QFD) \rightarrow (\E,\QFE)$, where we write $\im(p)$ for the Poincaré-Verdier quotient of $(\D,\QFD)$ by $\ker(p)$.
These inclusions then induce equivalences on mapping spaces into $(\cB,\Pi)^{\natural} \rightarrow \Latt(\cB,\Pi) \rightarrow \Tate(\cB, \Pi\qshift{1})^\natural$, since the latter consists of idempotent complete Poincaré $\infty$-categories. Similarly, the inclusion $(\C,\QF) \subseteq \ker(p)$ induces an equivalence on mapping spaces into $(\cB,\Pi)^{\natural}$.
But $\ker(p) \rightarrow (\D,\QFD) \rightarrow \mathrm{im}(p)$ is a Poincaré-Verdier sequence by Corollary~\reftwo{corollary:poincare-inclusion},
and so combining the above  
with Theorem \reftwo{theorem:universal-poincare-verdier} we conclude that the map
\[
\Fun^{\textsc{Kar}}((\C,\QF) \to (\D,\QFD) \to (\E,\QFE), (\cB,\Pi)^{\natural} \to \Latt(\cB,\Pi) \to \Tate(\cB, \Pi\qshift{1})^\natural) \to \Funx((\C,\QF),(\cB,\Pi)^{\natural}) ,
\]
induced by extracting kernels, is an equivalence.
The categories from the statement then correspond precisely to $\Funx((\C,\QF),(\cB,\Pi)) \subseteq \Funx((\C,\QF),(\cB,\Pi)^{\natural})$ 
by inspection.
\end{proof}

\section{The hermitian $\Q$-construction and algebraic cobordism categories}
\label{section:cobcats}%

In this section, we introduce the main objects of study, namely the cobordism $\infty$-category constructed from a Poincaré $\infty$-category. To motivate our perspective, let $(\C,\QF)$ be a Poincaré \(\infty\)-category and $(X,\qone),(X',\qtwo)$ be two Poincaré objects in $\C$. A \defi{cobordism} from $(X,\qone)$ to $(X',\qtwo)$ is a span of the form
\[
X \xleftarrow{\alp} W \xrightarrow{\bet} X'
\]
together with a path $\eta\colon \alp^*\qone \rightarrow \bet^*\qtwo$ in the space $\Om^{\infty}\QF(W)$ of hermitian structures on $W$, such that $\eta$ satisfies the Poincaré-Lefschetz condition with respect to $X$ and $X'$, i.e.\ that a canonical map 
\[
\fib(W \rightarrow X) \simeq \fib(X' \rightarrow X \cup_W X') \longrightarrow \fib(X' \rightarrow \Dual_\QF W) \simeq \Omega\Dual_\QF(\fib(W \rightarrow X')),
\]
induced by $\eta$ is an equivalence. 
For example, if $W$ is an oriented cobordism between two closed oriented $d$-manifolds $M$ and $N$, we obtain a span of the form
\[
C^*(M) \longleftarrow C^*(W) \lrar C^*(N)
\]
and the fundamental class $[W]$ determines a path relating the pullbacks of the two symmetric Poincaré structures $\qone_M$ and $\qone_N$ on $C^*(M)$ and $C^*(N)$, respectively. Lefschetz duality for manifolds with boundary precisely implies that this path exhibits the span as a cobordism between the Poincaré objects $(C^*(M),q_{M})$ and $(C^*(N),q_{N})$ of $(\Dperf(\ZZ), {\QF_{\ZZ}^{s}}\qshift{-d})$ in the sense above.

Now, cobordisms can be composed in a natural way, by first forming the corresponding composition at the level of spans and then at the level of the paths between hermitian structures. This will allow us to define an $\infty$-category $\Cob(\C,\QF)$ whose objects are the Poincaré objects of $(\C,\QF\qshift{1})$ and whose morphisms are given by cobordisms; the choice in shifts here adheres to the usual convention from manifold theory that the category $\Cob_d$ have $(d-1)$-dimensional, closed, oriented, smooth manifolds as objects and $d$-dimensional cobordisms as morphisms. 

To make this idea precise, we interpret a cobordism in $(\C, \QF)$ as a Poincaré object in the diagram $\infty$-category $(\Fun(P, \C), \QF^P)$, where $P$ is the category $\bullet \leftarrow \bullet \rightarrow \bullet$, and $\QF^P$ is the Poincaré structure on the diagram $\infty$-category given by the limit of the values of $\QF$ on the diagram. This construction turns out to be the degree $1$ part of a simplicial Poincaré $\infty$-category $\Q(\C, \QF)$, whose Poincaré objects in degree $n$ may be interpreted as the datum of $n$ composable tuples of cobordisms. Varying $(\C,\QF)$, this construction gives rise to a functor
\[
\Q \colon \Catp \lrar \sCatp,
\]
our implementation of the hermitian $\Q$-construction; see \S\reftwo{subsection:hermitian-Q}. Considering the spaces of Poincaré objects,
we will obtain a complete Segal space and then extract $\Cob(\C,\QF) \in \Cat$ as the associated $\infty$-category in \S\reftwo{subsection:cobcat}. 

Then, we develop the two main tools that will allow us to analyse this cobordism $\infty$-category and its homotopy type. First, we show how to describe the cobordism $\infty$-category using Ranicki's algebraic surgery techniques from \cite{RanickiATS1}, adapted to the setting of Poincaré \(\infty\)-categories by Lurie in \cite{Lurie-L-theory}. Beside its uses in the present paper, this serves as a fundamental tool in \cite{comparison} to compare our definition of Grothendieck-Witt theory with the classical $\L$-, Witt and Grothendieck-Witt groups, and it is also used extensively in \paperthree. The second topic, in \S\reftwo{subsection:additvity} and \S\reftwo{subsection:fibcob}, is the additivity theorem, which says that the functor
\[
\vert \Cob-\vert=\vert \Poinc\Q(-)\vert \colon \Catp \lrar \Sps;
\]
is additive. This will be the basis for most of the structural results we prove about Grothendieck-Witt theory. 

As far as the additivity theorem is concerned, the only property of the functor $\Poinc$ that enters the proof, is that it is itself additive. In fact, we will show that the functor
\[
\vert \F\Q(-)\vert\colon \Catp\lrar \Sps
\]
is additive whenever $\F\colon \Catp\to \Sps$ is additive. This added layer of generality will be used to establish the additivity of the spectral version of the Grothendieck-Witt functor, defined via iteration of the hermitian $\Q$-construction, and also enters into the proof of its universal property.

Finally, in \S\reftwo{subsection:Kadd} we explain how our methods give rise to a new proof of the more classical additivity theorem for the algebraic $\K$-theory of stable $\infty$-categories.

\subsection{Recollection on Segal spaces and their associated categories}
\label{subsection:rezk}%
 
Before developing the hermitian $\Q$-construction and our cobordism categories, we shall collect the necessary results regarding the relationship between $\infty$-categories and Rezk's complete Segal spaces, originally established in \cite{RezkSegal}. 

There is an adjoint pair of functors
\[
\asscat \colon \sSps \adj \Cat \cocolon \snerv
\]
with the Rezk nerve as right adjoint, given by
\[
\snerv(\C)_n =\Hom_{\Cat}(\Delta^n, \C),
\]
and left adjoint given by left Kan extending the cosimplicial $\infty$-category
\[
\Delta^- \colon \Delta \longrightarrow \Cat
\]
along the Yoneda embedding $\Delta \rightarrow \sSps$. The nerve is fully faithful with essential image the complete Segal spaces $\cSS \subseteq \sSps$, in particular making $\Cat$ a left Bousfield localisation of $\sSps$ (at what we shall refer to as the categorical equivalences). Consequently, there is also a left adjoint $\comp\colon \sSps \rightarrow \cSS$ to the inclusion, often referred to as completion, and composing adjoints we find $\asscat \circ \comp = \asscat$; the point-set version was originally shown in \cite{JoyalTierney} and the exact formulation above was established in  \cite{LurieGoodwillie}*{Corollary 4.3.16}, but see also \cite{Hebestreit-Steinebrunner} for a short proof. Recall that a simplicial space $X$ is called a Segal space if the Segal maps 
\[
\mathrm{seg}_i \colon [1] \longrightarrow [n], \quad 0,1 \longmapsto i-1,i
\]
for $1 \leq i \leq n$ induce equivalences
\[
(\mathrm{seg}_1, \dots, \mathrm{seg}_n) \colon X_n \longrightarrow X_1 \times_{X_0} X_1 \times_{X_0} \dots \times_{X_0} X_1
\]
for all $n \geq 1$. Completeness can be characterised in many ways; for us, the most convenient criterion is that a Segal space is complete (i.e.\ lies in the essential image of $\snerv$) if and only if 
\[
\begin{tikzcd}
X_0 \ar[r,"s"] \ar[d,"\Delta"] & X_3 \ar[d,"{(d_{02},d_{13})}"] \\ 
X_0^2 \ar[r,"{(s,s)}"] & X_1^2
\end{tikzcd}
\]
is cartesian, where $d_{02} \colon [1] \rightarrow [3]$ is the unique injective map that misses $0$ and $2$ and similarly for $d_{13}$, see \cite{LurieGoodwillie}*{Proposition 1.1.13} or \cite{Rezk-ncat}*{Proposition 10.1}.

By \cite{HTT}*{Proposition 5.5.4.15}, the categorical equivalences are the saturation of the spine inclusions, which encode the Segal condition, and the map 
\[
\Delta^3/\Delta^{0,2},\Delta^{1,3} = \Delta^0 \cup_\Delta^{\{0,2\}} \Delta^3 \cup_\Delta^{\{1,3\}} \Delta^0 \longrightarrow \Delta^0
\]
that encodes the completeness condition. Thus, a colimit preserving functor $\sSps \rightarrow \E$ (with $\E$ cocomplete) factors (uniquely) through $\asscat \colon \sSps \rightarrow \Cat$ if and only if it inverts these maps, i.e.\ if its restriction $\Delta\op \rightarrow \E\op$ along the Yoneda embedding is a complete Segal object in $\E\op$.

\medskip

The restriction of the nerve functor to $\Sps\subset \Cat$ is given by the inclusion of constant diagrams $\Sps \rightarrow \sSps$, and passing to adjoints shows that $|\asscat(X)| \simeq |X|$ for every simplicial space $X$. In particular, $\pi_0|\asscat (X)|$ is always the coequaliser of the two boundary maps $\pi_0X_1 \rightarrow \pi_0X_0$. 
The category $\asscat(X)$ itself is rather difficult to access for a general simplicial space $X$, except for the trivial fact that the canonical map
\[
X_0\lrar  \grpcr \asscat(X),
\]
induced by the inclusion $X_0\to X$ on associated $\infty$-categories, is $\pi_0$-surjective \cite{Hebestreit-Steinebrunner}*{Corollary 3.12}; here, $\grpcr$ denotes the groupoid core of an arbitrary $\infty$-category (we reserve the use of the symbol $\core$ for those occasions where we regard the groupoid core as an additive functor on stable or Poincaré $\infty$-categories). This changes for a Segal space $X$: Then Joyal, Rezk and Tierney showed that the unit provides a canonical equivalence between $\Hom_X(x,y)$, the fibre over $(x,y)$ of $(d_1,d_0) \colon X_1 \rightarrow X_0 \times X_0$, and $\Hom_{\comp(X)}(x,y) \simeq \Hom_{\asscat(X)}(x,y)$, see \cite{Hebestreit-Steinebrunner}*{Lemma 4.1}

Any categorical construction can thus be translated to the setting of Segal spaces. We require the following instances. Firstly, the \defi{groupoid core} of a Segal space $X$ is the Segal space $X^\times$ where $X^\times_n$ is the collection of components of $X_n$ consisting of all those simplices whose edges become equivalences in $\asscat(X)$. Then $\nerv(\C)^\times \simeq \nerv(\grpcr(\C))$ and the mapping space formula implies that the canonical maps
\[
\asscat(X^\times)\longrightarrow |X^\times| \longrightarrow \grpcr\asscat(X)
\]
are equivalences. It follows that for $X$ complete, the map $X_0\to \grpcr\asscat(X)$ is an equivalence.  

Secondly, a map of Segal spaces $p\colon X\to Y$ is called an \defi{isofibration} if the map $X_1^\times \to Y_1^\times\times_{Y_0} X_0$ is $\pi_0$-surjective, where the pullback is formed using either $d_0$ or $d_1\colon Y_1^\times\to Y_0$; the two conditions are equivalent as a consequence of the equivalence $\Hom_{X^\times}(x,y)\simeq  \Hom_{X^\times}(y,x)$ e.g.\ obtained from taking inverses in $\grpcr(\asscat(X))$. We also note that this condition is automatically satisfied if $Y$ is complete (e.g.\ constant). While $\asscat$ does not commute with pullbacks of Segal spaces in general, we still have:

\begin{lemma}
\label{lemma:asscat_pullbacks}%
For any cospan $X\to Y\leftarrow Z$ of Segal spaces, the natural comparison map 
\[
\asscat(X\times_Y Z) \lrar \asscat(X)\times_{\asscat(Y)} \asscat(Z)
\]
is fully faithful. It is an equivalence if the map $X\to Y$ is an isofibration, so $\asscat$ in particular commutes with finite products of Segal spaces.
\end{lemma}

\begin{proof}
The map is fully faithful by the formula for mapping spaces in the associated category of a Segal space. For essential surjectivity, we recall that an object in the target is represented by some $x\in X_0$, some $z\in Z_0$, and an equivalence between their images in $\asscat(Y)$. By the same formula, this equivalence is represented by some $f\in Y_1$ with source and target the images of $x$ and $z$, respectively. Lifting this to a morphism $f'\in X_1$ starting at $x$, the target of $f'$ together with $z$ defines a preimage.
\end{proof}

Thirdly, a morphism $f \in X_1$ is called \defi{$p$-cartesian} if the following square is cartesian
\[
 \begin{tikzcd}
  \{f\}\times_{X_1} X_2 \ar[r, "{(d_0, d_1)}"] \ar[d, "p"]  & \{z\}\times_{X_0} X_1 \ar[d, "p"]\\
  \{p(f)\}\times_{Y_1} Y_2  \ar[r, "{(d_0, d_1)}"] & \{p(z)\} \times_{Y_0} Y_1,
 \end{tikzcd}
\]
where $z = d_0(f)$ and the fibre products on the left are taken along the boundary $d_0$. Passing to fibres over some $x\in X_0$ and $p(x)\in Y_0$, this condition translates to the condition that $f$ is $\asscat(p)$-cartesian in the usual sense. The map $p$ is called \defi{cartesian fibration} if the map 
\[
(p, d_0)\colon X_1^{p-\mathrm{cart}} \lrar Y_1\times_{Y_0} X_0
\]
is surjective on $\pi_0$, where $X_1^{p-\mathrm{cart}} \subseteq X_1$ is the subspace spanned by the $p$-cartesian edges. It follows that cartesian fibrations are isofibrations. Moreover, a functor $\varphi$ of $\infty$-categories is a cartesian fibration (of $\infty$-categories) if and only if $\nerv(\varphi)$ is a cartesian fibration (of Segal spaces), and if $p$ is a cartesian fibration of Segal spaces, then $\asscat(p)$ is a cartesian fibration of $\infty$-categories, but generally the condition on $p$ is stronger. The notion of $p$-cocartesian morphisms and cocartesian fibrations are defined by passing to opposite simplicial spaces.

Lastly, we require the Segal version of the slice construction. Let us recall generally that for any simplicial object $X$ in some $\infty$-category $\mathcal A$, its décalage is the simplicial object 
\[
\dec(X) \colon \Delta\op \rightarrow \mathcal A, \quad [n]\longmapsto X_{[0] \ast [n]}= X_{1+n}
\]
and that the d\'ecalage is a Segal object whenever $X$ is Segal. Considering the maps $d_0 \colon [n] \rightarrow [1+n]$ and $0 \colon [0] \rightarrow [1+n]$ as a natural transformations $\id_{\Delta} \Rightarrow [0] \ast \id_{\Delta}$ and $\const_{[0]} \Rightarrow [0] \ast \id_{\Delta}$ yields maps of simplicial objects $\dec(X) \to X$ and $\dec(X) \rightarrow \const X_0$, respectively, the former of which, for $X$ a Segal object, participates in a pullback square
\[
\begin{tikzcd}
	\const(X_1) \ar[r] \ar[d,"d_0"] & \dec(X) \ar[d,"\pi"] \\
	\const(X_0) \ar[r] & X
\end{tikzcd}
\]
We then call the Segal space 
\[
X_{x/}\coloneq\dec(X)\times_{\const(X_0)} \const(\{x\})
\]
the \defi{Segal slice} of $X$ under $x \in X_0$. This extends the categorical slice construction in the sense that there is a canonical natural equivalence $\nerv(\C)_{c/}\simeq \nerv(\C_{c/})$ induced by the canonical equivalence $(\Delta^n\times \Delta^1)\cup_{\Delta^n\times \{0\}} \{*\}\simeq \Delta^{1+n}$ of $\infty$-categories. We also have:

\begin{lemma}
\label{lemma:Segal_slice_asscat}%
The canonical map
\[
\asscat(X_{x/})\longrightarrow \asscat(X)_{x/}
\]
is an equivalence for any Segal space $X$ and $x \in X_0$.
\end{lemma}

\begin{proof}
To see that this is an equivalence, note that it suffices to show that $X_{x/} \rightarrow (\comp X)_{x/}$ is a categorical equivalence. The map $\asscat(X_{x/})\to \asscat((\comp X)_{x/})$ is essentially surjective because any object in the target is represented by some $f\in \Hom_{\comp X}(x,y)\simeq \Hom_X(x,y)$ for suitable $x,y\in X_0$. To see that it is fully faithful, we need to show that for further $g\in \Hom_X(x,z)$, the expression $\Hom_{X_{x/}}(f,g)$ is invariant under passing to completion. Now, $\Hom_{X_{x/}}(f,g)$ is the fibre over $(f,g)$ of the left-hand map in the chain
\[
X_2 \xrightarrow{(d_2, d_1)} X_1\times_{X_0} X_1 \xrightarrow{(\id, d_0)} X_1\times X_0.
\]
The right-hand fibre over $(f,z)$ is $\Hom_X(x,z)$ and the fibre of the composite map is (by the Segal condition) $\Hom_X(y,z)$. Since both individual fibres are invariant under passing to completion, the claim follows. 
\end{proof}

\subsection{The hermitian $\Q$-construction}
\label{subsection:hermitian-Q}%

Let $K$ be an $\infty$-category and $(\C,\QF)$ a hermitian \(\infty\)-category. We denote by $\Twar(K)$ the twisted arrow category, the conventions being such that $(s,t) \colon \Twar(K) \rightarrow K \times K\op$ give the right fibration classifying $\Hom_K \colon K\op \times K \rightarrow \Sps$.

\begin{definition}
\label{definition:hermitianQ}%
Let $\Q_K(\C,\QF)$ denote the following hermitian \(\infty\)-category: The underlying stable \(\infty\)-category is given as the full subcategory $\Q_K(\C)$ of $\Fun(\Twar(K),\C)$ spanned by those functors $F$ such that for every functor $[3] \rightarrow K$, say $i \rightarrow j \rightarrow k \rightarrow l \in K$, the square
\[
\begin{tikzcd}
F(i\rightarrow l) \ar[r]\ar[d] & F(j\rightarrow l) \ar[d] \\
F(i\rightarrow k) \ar[r]       & F(j\rightarrow k)
\end{tikzcd}
\]
is cartesian. The hermitian structure is given by restricting the quadratic functor 
\[
F \longmapsto \QF^{\Twar(K)}(F) = \lim_{\Twar(K)\op} \QF\circ  F\op
\]
from Proposition~\refone{proposition:basic-properties-herm-diagrams}.
\end{definition}

When $K=\Delta^n$, we will shorten notation, and denote $\Q_{K}(\C,\QF)$ by $\Q_n(\C,\QF)$ and $\QF^{\Twar(\Delta^n)}$ by $\QF_n$. Also, by definition, the hermitian \(\infty\)-category $(\Fun(\Twar(K), \C), \QF^{\Twar(K)})$ is the cotensor $(\C, \QF)^{\Twar(K)}$, in the sense of \S\refone{subsection:cotensoring}. It is usually not Poincaré, while $\Q_K(\C, \QF)$ is, as we will see below.

\begin{examples}
\label{examples:DescriptionsofQ}%
\ %
\begin{enumerate}
\item
\label{item:Q1-with-duality}%
Since $\Twar(\Delta^0) = \Delta^0$, we have $\Q_0(\C,\QF) = (\C,\QF)$. Similarly, we have $\Twar(\Delta^1) = \{(0 \leq 0) \leftarrow (0 \leq 1) \rightarrow (1 \leq 1)\}$, so $\Q_1(\C)$ is simply the $\infty$-category of spans in $\C$, with no condition imposed. Using Proposition~\refone{proposition:basic-properties-herm-diagrams}, the duality on $\Q_1(\C,\QF)$ takes
\[
\bigl(X \leftarrow Y \rightarrow Z\bigr) \quad \mapsto \quad \bigl(\Dual_\QF X \longleftarrow \Dual_\QF X \times_{\Dual_\QF Y} \Dual_\QF Z \longrightarrow \Dual_\QF Z\bigr).
\]
Following our explanation in the introduction to this section, we interpret $\Q_1(\C,\QF)$ as the \emph{$\infty$-category of cobordisms} in $(\C,\QF)$. 
\item 
\label{item:met-vs-Q-1}%
By \reftwoitem{item:Q1-with-duality}, the functor 
\[
d_1\colon \Q_1(\C, \QF)\lrar \Q_0(\C, \QF) = (\C,\QF), \quad (X\leftarrow Y\rightarrow Z)\mapsto X
\]
is duality-preserving so that its kernel is closed under the duality of $\Q_1(\C,\QF)$, and therefore a Poincaré \(\infty\)-category with the restricted Poincaré structure. And indeed, we have $\ker(d_1) \simeq \Met(\C,\QF)$ via
$(0\leftarrow W \to X) \mapsto (W\to X)$ and similarly for $\ker(d_0)$.
\item 
$\Q_2(\C,\QF)$ consists of those diagrams
\[
\begin{tikzcd}
[row sep=2.5ex, column sep=1ex]
 & & F(0 \leq 2) \ar[ld]\ar[rd]& & \\
 & F(0 \leq 1) \ar[ld]\ar[rd]& & F(1 \leq 2) \ar[ld]\ar[rd]& \\
F(0 \leq 0) && F(1 \leq 1) && F(2 \leq 2)
\end{tikzcd}
\]
in which the top square is bicartesian. It is therefore reasonable to think of $\Q_2(\C,\QF)$ as the $\infty$-category of two composable cobordisms equipped with a chosen composite.

\item 
\label{item:QnIn}%
We note that for the category $\I_n \subseteq \Twar(\Delta^n)$ spanned by the pairs $(i,j)$ with $j \leq i+1$ (the zig-zag along the bottom) the restriction functor 
\[
\Q_n(\C,\QF) \lrar (\Fun(\I_n,\C),\QF^{\I_n}) = (\C,\QF)^{\I_n}
\]
is an equivalence of hermitian \(\infty\)-categories: On underlying $\infty$-categories, it follows from \cite{HTT}*{Proposition 4.3.2.15} that the right Kan extension functor $\Fun(\I_n,\C) \rightarrow \Fun(\Twar(\Delta^n),\C)$ is both fully faithful and a left inverse to restriction. For $X \in \Fun(\Twar(\Delta^n),\C)$, it is then readily checked from the pointwise formulae  \cite{HTT}*{Lemma 4.3.2.13} that being in $\Q_n(\C)$ is equivalent to being right Kan extended from $\I_n$. For the quadratic functor, it follows since the inclusion $\I_n\op \subseteq \Twar(\Delta^n)\op$ is final. By Remark~\refone{remark:one-sided-inverse}, the arising hermitian structure on the right Kan extension functor $\Fun(\I_n,\C) \rightarrow \Fun(\Twar(\Delta^n),\C)$ is an instance of the exceptional functoriality of Construction~\refone{construction:exceptional}.

This description justifies us in thinking of $\Q_n(\C,\QF)$ as the $\infty$-category of $n$ composable cobordisms in $(\C,\QF)$ also for larger $n$.
\end{enumerate}
\end{examples}

Denoting the $\infty$-category of $\infty$-categories by $\Cat$, we obtain a functor
\[
\Cat\op \times \Cath \lrar \Cath, \quad (K,\C,\QF) \mapsto \Q_K(\C,\QF),
\]
from Proposition~\refone{corollary:functorial-cotensor}, since clearly induced maps preserve the cartesianness condition of Definition~\reftwo{definition:hermitianQ}. 

\begin{definition}
Restricting along the inclusion $\Delta \subseteq \Cat$ and adjoining the construction above, we obtain a functor $\Q \colon \Cath \rightarrow \sCath$, and in particular simplicial hermitian $\infty$-categories $\Q(\C,\QF) \in \sCath$ for all $(\C,\QF) \in \Cath$. We shall call both the initial construction as well as this restriction the \emph{hermitian $\Q$-construction}.
\end{definition}

We note that the underlying $\infty$-category of $\Q_n(\C,\QF)$ only depends on $\C$, and agrees with Barwick-Rognes' implementation $\Q_n(\C)$ of the $\Q$-construction, see \cite{BarwickRognesQexact}*{\S 3} upon restricting their set-up to stable \(\infty\)-categories. 
The following is at the heart of the present section:

\begin{lemma}
\label{lemma:segal}%
For every hermitian \(\infty\)-category $(\C,\QF)$ the simplicial hermitian \(\infty\)-category $\Qh(\C,\QF)$ is a Segal object of $\Cath$. Furthermore, it is complete in the sense that the diagram
\[
\begin{tikzcd}
\Q_0(\C,\QF) \ar[r,"s"] \ar[d,"\Delta"] & \Q_3(\C,\QF) \ar[d,"{(d_{02},d_{13})}"] \\ 
\Q_0(\C,\QF)^2 \ar[r,"{(s,s)}"] & \Q_1(\C,\QF)^2
\end{tikzcd}
\]
is cartesian in $\Cath$, with horizontal maps given by total degeneracies.
\end{lemma}
\begin{proof}
We will show that for every $0 \leq i \leq n$ the square
\[
\begin{tikzcd}
\Q_{n}(\C,\QF) \ar[r]\ar[d] & \Q_{[0,i]}(\C,\QF) \ar[d] \\
\Q_{[i,n]}(\C,\QF) \ar[r] & \Q_{[i,i]}(\C,\QF)
\end{tikzcd}
\]
is a pullback square of hermitian $\infty$-categories; the Segal condition then follows by iteration. The statement will follow readily from Example~\reftwo{examples:DescriptionsofQ} \reftwoitem{item:QnIn}. To this end, note that the inclusions $\Twar(\Delta^i) \rightarrow \Twar(\Delta^n)$ and $\Twar(\Delta^{\{i,\dots,n\}}) \rightarrow \Twar(\Delta^n)$ take the subcategories $\I_i$ and $\I_{[i,\dots,n]}$ to $\I_n$, and in fact the induced diagram 
\[
\begin{tikzcd}
\I_{[i,i]} \ar[r]\ar[d] & \I_{[i,n]} \ar[d] \\
\I_{i}\ar[r] & \I_{n}
\end{tikzcd}
\]
is readily checked to be cocartesian in $\Cat$, thus cartesian in $\Cat\op$. Since the functor 
\[
\Cat\op \lrar \Catp, \quad I \mapsto (\Fun(I,\C),\QF^I)
\]
is a right adjoint we obtain the first claim. 

To see that $\Q(\C,\QF)$ is complete, recall that limits in $\Catx$ may be computed in $\Cat$, 
so the map $P \rightarrow \Q_3(\C)$ from the pullback $P$ of the diagram 
\[
\Q_0(\C)^2 \lrar \Q_1(\C)^2 \longleftarrow \Q_3(\C)
\]
is fully faithful, since the degeneracy $\Q_0(\C) \rightarrow \Q_1(\C)$ is, and fully faithful functors are stable under pullback. Its essential image is given by the diagrams consisting entirely of equivalences, as one can check directly using the defining property of the $\Q$-construction, and these are precisely the constant diagrams, i.e.\ the totally degenerate ones, since $|\Twar(\Delta^3)| \simeq \ast$. 

The claim for the hermitian structure is immediate from Remark~\refone{remark:explicit-limits-cath},
since the diagram
\[
\QF_0(X)^2 \lrar \QF_1(sX)^2 \longleftarrow \QF_3(sX),
\]
whose pullback defines the hermitian structure on $P$, evaluates to
\[
\QF(X)^2 \xrightarrow{\id} \QF(X)^2 \xleftarrow{\Delta} \QF(X),
\]
so has pullback $\QF(X)$.
\end{proof}

\begin{lemma}
\label{lemma:Qlimitperserving}%
For fixed $(\C,\QF) \in \Cath$ the functor $\Q_{-}(\C,\QF) \colon \Cat\op \rightarrow \Cath$ preserves limits.
\end{lemma}

\begin{proof}
On underlying $\infty$-categories, this is \cite{HHLN-two-var-fil}*{Proposition~2.20}. We repeat the argument with hermitian structures tagging along. The preservation of limits can be recast as $\Q(\C,\QF) \colon \Delta\op \rightarrow \Cath$ being a complete Segal object (which we showed above), and the full functor $\Q_-(\C,\QF) \colon \Cat\op \rightarrow \Cath$ being right Kan extended from its restriction to $\Delta\op$; see \cite{HHLN-two-var-fil}*{Lemma 2.21}. By the pointwise formula this means that the natural map $\Q_J(\C,\QF) \rightarrow \lim_{n \in (\Delta/J)\op} \Q_n(\C,\QF)$ has to be an equivalence for every $J \in \Cat$. 

The diagram $\Delta/J \rightarrow \Delta \rightarrow \Cat$ is a typical example of a functor $I \colon K \rightarrow \Cat$ whose colimit is preserved by the nerve functor.
Such colimits are also preserved by the functor $\Twar \colon \Cat \rightarrow \Cat$, by a direct calculation at the level of Rezk nerves, and so are the subcategories making up the $\Q$-construction. Therefore, 
\[
\lim_{k \in K}\Q_{I_k}(\C,\QF) \quad \text{and} \quad \Q_{\colim_{k \in K} I_k}(\C,\QF)
\]
are the same hermitian subcategory of 
\[
\lim_{k \in K}(\C,\QF)^{\Twar(I_k)} \simeq (\C,\QF)^{\colim_{k \in K} \Twar(I_k)} \simeq (\C,\QF)^{\Twar(\colim_{k \in K} I_k)},
\]
and thus $\Q_-(\C,\QF)$ commutes with limits over diagrams that are compatible with the Rezk nerve, which is more than we need. 
\end{proof}

\begin{lemma}
\label{lemma:poincare}%
\label{lemma:boundarysplit}%
The functor $\Q \colon \Cat\op \times \Cath \rightarrow \Cath$ restricts to a functor $\Cat\op \times \Catp \rightarrow \Catp$. Moreover, for $(\C,\QF)$ Poincar\'e, $\Q(\C,\QF) \in \sCatp$ is a complete Segal object of $\Catp$, all whose structure maps induced by injections in $\Delta$ are split Poincar\'e-Verdier projections.
\end{lemma}

\begin{proof}
There are two good approaches to the statements about $\Q(\C,\QF)$. Either, one directly attacks them using the machinery developed in \S\refone{subsection:finite-complexes}, or one reduces the statement to explicit checks for small values of $n$ using the Segal condition. At the cost of being less elementary, we will here use the former route as it leads to shorter proofs. 

That the $\infty$-categories $\Q_n(\C,\QF)$ are Poincaré follows immediately from Proposition~\refone{proposition:posets-of-faces} and Examples~\reftwo{examples:DescriptionsofQ}, since $\I_n$ is the poset of faces for the triangulation of the interval using $n+1$ vertices. It also follows from \refone{proposition:posets-of-faces} that Poincaré functors $(\C,\QF) \rightarrow (\C',\QF')$ induce Poincaré functors $\Q_n(\C,\QF) \rightarrow \Q_n(\C',\QF')$.

To see that the induced hermitian functors $\alpha^* \colon \Q_n(\C,\QF) \rightarrow \Q_m(\C,\QF)$ for $\alpha \colon \Delta^m \rightarrow \Delta^n$ preserve the dualities, we distinguish two cases, namely the inner face maps on the one hand, and the outer face maps and degeneracies on the other. Since every morphism in $\Delta$ can be written as a composition of such, this will suffice for the claim. 

The latter maps all take the subset $\I_m \subseteq \Twar(\Delta^m)$ into $\I_n$, and the restriction is induced by a map of the simplicial complexes giving rise to $\I_m$ and $\I_n$. Thus, Proposition~\refone{proposition:poset-of-faces-functoriality} gives the claim. The interior faces do not preserves the subsets $\I_m$, however. Instead, we claim that they are instances of the exceptional functoriality of Construction~\refone{construction:exceptional} associated to a refinement among triangulations. Namely, one readily checks that $d_i \colon \Twar(\Delta^n) \rightarrow \Twar(\Delta^{n+1})$ admits a right adjoint $r_i \colon \Twar(\Delta^{n+1}) \rightarrow \Twar(\Delta^n)$ explicitly given by 
\[
(k \leq l) \longmapsto 
\begin{cases} 
        (k \leq l) & l < i \text{ or } k<l=i \\
	(k-1 \leq l) & k=l=i \\
	(k \leq l-1) & k<i<l \\								(k-1 \leq l-1) & i \leq k < l \text{ or } i <k=l.
\end{cases}
\]
As a right adjoint, $r_i$ is cofinal, so by Example~\refone{example:exceptional-is-direct}, the pullback functor
\[
(d_i)^* \colon (\C,\QF)^{\Twar(\Delta^{n+1})} \longrightarrow (\C,\QF)^{\Twar(\Delta^{n})}
\]
agrees with the exceptional functoriality along $r_i$. From the explicit formula, it is clear that $r_i$ takes $\I_{n+1}$ into $\I_n$, so we find a diagram
\[
\begin{tikzcd}
(\C,\QF)^{\I_{n+1}} \ar[r,"{(r_i)_*}"] \ar[d] & (\C,\QF)^{\I_n} \ar[d] \\
(\C,\QF)^{\Twar(\Delta^{n+1})} \ar[r,"{(r_i)_*}"] & (\C,\QF)^{\Twar(\Delta^n)}
\end{tikzcd}
\]
where the vertical maps are the exceptional functorialities associated to the inclusions $\I_n \subseteq \Twar(\Delta^n)$ which are also cofinal (the diagram commutes since exceptional functorialities compose by Remark~\refone{remark:compatible-with-composition}). But the vertical maps are equivalences onto $\Q_n(\C,\QF)$ by Example~\reftwo{examples:DescriptionsofQ} \reftwoitem{item:QnIn}. The claim now follows from Proposition~\refone{proposition:poset-of-faces-functoriality}, since the restriction of $r_i$ to $\I_{n+1} \rightarrow \I_n$ comes from the refinement of triangulation of the interval that adds a new $i$-th vertex. 

This shows that $\Q$ restricts to a functor $\Delta\op \times \Catp \rightarrow \Catp$, and in particular, it follows from Lemma~\reftwo{lemma:segal} above, that $\Q(\C,\QF)$ is a complete Segal object of $\Catp$ for every $(\C,\QF) \in \Catp$, because limits in $\Catp$ are computed in $\Cath$. It follows from Lemma~\reftwo{lemma:Qlimitperserving} that $\Q_{-}(\C,\QF)$ is right Kan extended from its restriction to $\Delta\op \subseteq \Cat\op$, so the hermitian $\Q$-construction indeed restricts to a functor $\Cat\op \times \Catp \to \Catp$ since $\Catp \to \Cath$ preserves and detects limits.
For the final claim, we only need to consider face maps, since split Poincaré-Verdier projections are stable under composition by the characterisation in Corollary~\reftwo{corollary:poincare-projection}. For the inner faces, this is immediate from Proposition~\reftwo{proposition:kan-extension-cofinal}, since $r_i \colon \I_{n+1} \rightarrow \I_n$ is evidently a localisation at the edges $(i-1 \leq i) \rightarrow (i \leq i)$ and $(i \leq i+1) \rightarrow (i \leq i)$. For the outer faces, it is an instance of Proposition~\reftwo{proposition:example-cosieve}.
\end{proof}

\subsection{The cobordism $\infty$-category of a Poincaré \(\infty\)-category}
\label{subsection:cobcat}%

We now proceed to extract the cobordism $\infty$-category from the hermitian $\Q$-construction. As mentioned in the introduction, it is useful to do this in the generality of an arbitrary additive $\F \colon \Catp \rightarrow \Sps$, see Definition \reftwo{definition:additive}, but the reader is encouraged to envision $\F = \Poinc$ throughout.

\begin{proposition}
\label{proposition:complete}%
Let $(\C,\QF)$ be a Poincaré \(\infty\)-category and $\F\colon \Catp\to \Sps$ an additive functor. Then $\F\Q(\C,\QF)$ is a Segal space which is complete if $\F$ preserves arbitrary pullbacks.
\end{proposition}

When $\F$ is the functor $\Core\colon\Catp \rightarrow \Sps$, completeness was established in \cite{BarwickRognesQexact}*{Proposition 3.4} by different means. For a general additive functor $\F$, the Segal space $\F\Q(\C,\QF)$ is not complete. For example, if $\F$ is group-like, then $\F\Q(\C,\QF)$ is complete if and only if $\F\Hyp(\C) \simeq 0$, see Remark~\reftwo{remarks:FQcompletebord}.

\begin{proof}
For the first part we show that 
\[
\begin{tikzcd}
\F\Q_{n}(\C,\QF) \ar[r]\ar[d] & \F\Q_{[0,i]}(\C,\QF) \ar[d] \\
\F\Q_{[i,n]}(\C,\QF) \ar[r] & \F\Q_{[i,i]}(\C,\QF)
\end{tikzcd}
\]
is cartesian for every $0 \leq i \leq n$. Indeed, before applying $\F$, the square is a Poincaré-Verdier square by Lemmas~\reftwo{lemma:boundarysplit} and \reftwo{lemma:segal}, and by assumption $\F$ preserves the cartesianness of such squares.

The assertion on completeness is immediate from the final part of Lemma~\reftwo{lemma:segal}.
\end{proof}

\begin{definition}
\label{definition:cobcat}%
Let $\Cob^\F(\C,\QF)$ denote the $\infty$-category associated to the Segal space $\F\Qh(\C,\QF\qshift{1})$. We shall write $\Cob(\C,\QF)$ for $\Cob^{\Poinc}(\C,\QF)$ and call it the \emph{cobordism $\infty$-category} of $(\C,\QF)$. Furthermore, we set $\Cob^\partial(\C,\QF) = \Cob(\Met(\C,\QF\qshift{1}))$, the \emph{cobordism $\infty$-category with boundaries}.
\end{definition}

We shall refer to $\Cob^\F(\C,\QF)$ as the $\F$-based cobordism $\infty$-category and hope the two possible superscripts ($\F$ and $\partial$) will not lead to confusion. By the functoriality of the $\Q$-construction and the previous discussion, the construction of these $\infty$-categories assemble into a functor $\Funadd(\Catp,\Sps) \times \Catp \lrar \Cat.$ An entirely analogous definition can be made for additive functors $\F \colon \Catx \rightarrow \Sps$ (i.e.\ reduced and sending split Verdier squares to cartesian squares), resulting in an $\infty$-category $\asscat (\F\Q(\C)) = \Span^\F(\C)$, with $\F = \Core$ giving rise to the usual span $\infty$-category $\Span(\C)$ considered for example in \cite{BarwickRognesQexact, HHLN-two-var-fil}.

\begin{examples}
\label{example:CobHyp-equal-Span}%
\ 
\begin{enumerate}
\item Straight from the definition we have $\Cob^{\Core}(\C,\QF) \simeq \Span(\C)$ for every Poincaré $\infty$-category $(\C,\QF)$.
\item Similarly, one obtains an equivalence
\[
\Cob^\F(\Hyp(\C)) \simeq \Span^{\F \circ \Hyp}(\C)
\]
by commuting the hyperbolic and $\Q$-constructions: From the natural equivalences of Remarks~\refone{remark:comparing-1} and \refone{remark:functor-category-to-hyp}, we find
\begin{align*}
\Funx((\E, \QF), \Q_n\Hyp(\C)) 
& \simeq \Funx((\E, \QF), \Hyp(\C)^{\I_n})  \\
& \simeq \Funx((\E, \QF)_{\I_n}, \Hyp(\C)) \\ 
& \simeq \Hyp(\Funx(\E_{\I_n}, \C)) \\
& \simeq \Hyp(\Funx(\E, \C^{\I_n})) \\
& \simeq \Funx((\E, \QF), \Hyp \Q_n(\C))
\end{align*}
so the natural map $\Q\Hyp(\C) \Rightarrow \Hyp \Q(\C)$ in $\sCatp$ is an equivalence.
\item In particular, $\Poinc\circ \Hyp \simeq \Core$ gives $\Cob(\Hyp(\C)) \simeq \Span(\C)$
for every stable $\infty$-category, see Proposition~\refone{proposition:poinc-of-hyp}.
\item There are canonical equivalences
\[
\Cob(\C,\QF^\sym)  \simeq \Span(\C)^\hC:
\]
By Remark~\refone{remark:action-on-core}, a Poincaré structure on an $\infty$-category $\D$ induces a natural $\Ct$-action on $\grpcr \D$. %
In particular, $\QF^\sym$ induces a $\Ct$-action on the simplicial space $\Core\Q(\C)$ and therefore a $\Ct$-action on the associated $\infty$-category $\Span(\C)$. By Proposition~\refone{proposition:basic-properties-hermitian-functor-cats}, the Poincaré structure $(\QF^\sym)^{\Twar[n]}$ is symmetric so that by Proposition~\refone{proposition:compare-symmetric} %
$\Poinc\Q_n(\C, \QF^\sym)\simeq \Core\Q_n(\C)^{\hC}$. As $\Core\Q(\C)$ is a complete Segal space, this implies the claim.
\item There is a natural equivalence
\[
\Cob^\F(\C,\QF) \simeq \Cob^\F(\C,\QF)\op
\]
since $\Q(\C,\QF)$ is naturally identified with the reversal simplicial object $\Q(\C,\QF)\op$ via the canonical (in fact unique) identification $\Twar(\Delta^n) \cong \Twar((\Delta^n)\op)$ of cosimplicial categories.
\end{enumerate}
\end{examples}

We now collect a few basic properties of the cobordism $\infty$-categories. Note that the inclusion of $0$-simplices of $\F\Q(\C,\QF\qshift{1})$ gives a natural map
\[
\F(\C,\QF\qshift{1}) \longrightarrow \grpcr\Cob^\F(\C,\QF)
\]
that is surjective on $\pi_0$. Informally, for $\F = \Poinc$, this map takes any Poincaré object to itself and an equivalence $f\colon X \rightarrow X'$ to the cobordism $X \xleftarrow{\id} X \xrightarrow{f} X'$. 
Proposition~\reftwo{proposition:complete} implies:

\begin{corollary}
\label{corollary:CobCatCore}%
The natural map 
$
\F(\C,\QF\qshift{1}) \rightarrow \grpcr\Cob^\F(\C,\QF)
$
is an equivalence whenever $\F$ preserves pullbacks. In particular, a Poincaré cobordism 
\[
(X,q) \longleftarrow (W,\eta) \lrar (X',q')
\]
considered as a morphism in $\Cob(\C,\QF)$ is invertible if and only if both underlying maps $W \rightarrow X$ and $W \rightarrow X'$ are equivalences in $\C$. 
\end{corollary}

\begin{remark}
There is an analogous result for the geometric cobordism category $\Cob_d$:
If a morphism $W$ in $\Cob_d$ is invertible, then it is an $h$-cobordism and the converse is true if $d\neq 4$, the inverse of $W$ being given by the $h$-cobordism with Whitehead torsion $-\tau(W) \in \mathrm{Wh}(\pi_1(\partial_0 W))$. 

Furthermore, the homotopy type of $\grpcr\Cob_d$ is closely related to the classifying space for $h$-cobordisms \cite{RaptisSteimlehcob}. 
\end{remark}

There is a simple way to produce diagrams in $\Cob^\F(\C,\QF)$, which will be heavily exploited in \cite{comparison}. Namely, for $K \in \Cat$ consider the map
\[
\F\Q_{K}(\C,\QF) \longrightarrow \Hom_{\sSps}(\snerv K,\F\Q(\C,\QF)) \longrightarrow \Hom_{\sSps}(\snerv K,\comp(\F\Q(\C,\QF)))\simeq \Hom_{\Cat}(K,\Cob^\F(\C,\QF))
\]
in which the first map arises by observing that the target is the right Kan extension of the restriction of the source, viewed as a functor in $K$, to $\Delta \subseteq \Cat$, 
the second map simply applies completion, and the final equivalence is an instance of the full faithfulness of $\snerv$.

\begin{proposition}
\label{proposition:functorstocob}%
If $\F \colon \Catp \rightarrow \Sps$ preserves all limits, the composite $\F\Q_{K}(\C,\QF) \rightarrow \Hom_{\Cat}(K,\Cob^\F(\C,\QF))$ above is an equivalence for all $K \in \Cat$ and $(\C,\QF) \in \Catp$.
\end{proposition}
\begin{proof}
All maps in the construction above are equivalences in this case, the first by Lemma~\reftwo{lemma:Qlimitperserving} and the second by Lemma~\reftwo{lemma:poincare}.
\end{proof}

Since the association $(\C,\QF) \mapsto \Q_n(\C,\QF)$ preserves finite products, as does completion of Segal spaces, it follows that the functor $\Cob^\F\colon \Catp \to \Cat$ preserves finite products. Since $\Catp$ is semi-additive (see Proposition~\refone{proposition:catp-pre-add}),
the $\infty$-categories $\Cob^\F(\C,\QF)$ acquire natural symmetric monoidal structures induced by the direct sum operation in $\C$ by \cite{GepGroNik}*{Corollary 2.5}. In particular, $\pi_0|\Cob^\F(\C,\QF)|$ is naturally a commutative monoid; explicitly when $\F = \Poinc$, it is the monoid of cobordism classes of Poincaré objects in $(\C,\QF)$ under orthogonal sum. The following is a description of it in general.

\begin{proposition}
\label{proposition:components-of-cob-I}%
\label{corollary:components-of-cob-II}%
For any additive functor $\F\colon \Catp\to \Sps$, the natural map $\pi_0\F(\C,\QF\qshift{1}) \rightarrow \pi_0|\Cob^\F(\C,\QF)|$ fits into a cocartesian square
\[
\begin{tikzcd}
\pi_0\F(\Met(\C,\QF\qshift{1})) \ar[r,"\met"] \ar[d] & \pi_0\F(\C,\QF\qshift{1}) \ar[d]\\
0 \ar[r] & \pi_0|\Cob^\F(\C,\QF)|
\end{tikzcd}
\]
of commutative monoids. Moreover, $\pi_0|\Cob^\F(\C,\QF)|$ is a group whose inversion is induced by the Poincaré functor $(\id_\C,-\id_\QF)\colon (\C,\QF) \to (\C,\QF)$.
\end{proposition}
In particular, we have $\pi_0|\Cob(\C, \QF)| \cong \L_{-1}(\C,\QF)$, see \S\refone{subsection:algebraic-thom} or \S\reftwo{subsection:L+tate} below for discussions of $\L$-groups in the present context. For the proof of Proposition~\reftwo{proposition:components-of-cob-I}, we will make use of the following.
\begin{construction}
\label{construction:bentcyl}%
Consider the hermitian functor $\bcyl \colon (\C,\QF) \to \Met(\C,\QF) \subseteq \Q_1(\C,\QF)$, representing a \emph{bent cylinder}, which consists of the functor
\[
X \mapsto [X \oplus X \xleftarrow{\Delta_X} X \rightarrow 0]
\]
and the map of quadratic functors induced by the diagram
\[
\begin{tikzcd}[cramped]
\QF(X) \ar[rd,dashed] \ar[rdd,bend right] \ar[rrrd,bend left=10,"{(\id,-\id,0)}"] & & & \\
  & \QF_1(\bcyl X) \ar[r] \ar[d] & \QF(X \oplus X) \ar[r] \ar[d,"{\Delta^*}"]& \QF(X) \oplus \QF(X) \oplus \Bil_\QF(X,X) \ar[dl,"\pr_1 + \pr_2"] \\
 & \ast \ar[r] & \QF(X) &
\end{tikzcd}
\]
whose left hand square is cartesian by definition of $\QF_1$, and whose right most horizontal map is an equivalence by definition of $\Bil_\QF$. One readily checks that $\bcyl$ is a Poincaré functor.
\end{construction}

\begin{proof}[Proof of Proposition~\reftwo{proposition:components-of-cob-I}]
We claim that
\[
\pi_0 \F(\Met(\C, \QF\qshift 1)) \xrightarrow\met \pi_0 \F(\C, \QF\qshift 1) \to \pi_0|\Cob^\F(\C,\QF)| \to 0
\]
is an exact sequence of commutative monoids. Granting this for a moment, the induced map 
\[
\pi_0\F(\C,\QF\qshift{1})/\mathrm{Im}(\met) \lrar \pi_0|\Cob^\F(\C,\QF)|
\]
is surjective and has trivial kernel. Furthermore, its source is a group since for each $x \in \pi_0\F(\C,\QF\qshift{1})$, the element $\bcyl_*(x)$ from Construction~\reftwo{construction:bentcyl} witnesses the relation 
\[
0 = x+ (\id_\C,-\id_\QF)_*(x)
\]
in $\pi_0\F(\C,\QF\qshift{1})/\mathrm{Im}(\met)$. Putting these together gives the proposition. 

To see the claim recall that $\pi_0|\Cob^\F(\C,\QF)|$ is the coequaliser of the two boundary maps $d_0, d_1\colon \pi_0 \F\Q_1(\C, \QF\qshift 1)\to \pi_0 \F(\C, \QF\qshift 1)$. Thus the second map in the sequence above is indeed surjective. For exactness at $\pi_0 \F(\C, \QF\qshift 1)$, we note that the image of $(d_0, d_1)$ is already a congruence relation on $\pi_0 \F(\C, \QF\qshift 1)$: It is clearly reflexive, transitive and compatible with addition, and symmetry follows from the evident automorphism of $\Q_1(\C,\QF)$ swapping source and target. Thus $x \in \pi_0\F(\C,\QF\qshift{1})$ vanishes in $\pi_0|\Cob^\F(\C,\QF)|$ if and only if there exists a $w \in \pi_0\F\Q_1(\C,\QF\qshift{1})$ with $d_1w = x$ and $d_0w = 0$. But since there is a fibre sequence
\[
\F(\Met(\C,\QF\qshift{1})) \lrar \F\Q_1(\C,\QF\qshift{1}) \xrightarrow{d_0} \F(\C,\QF\qshift{1})
\]
by Lemma~\reftwo{lemma:boundarysplit} and Example \reftwo{examples:DescriptionsofQ} \reftwoitem{item:met-vs-Q-1}, this is equivalent to $w$ lifting to $\pi_0 \F(\Met(\C, \QF\qshift 1))$. The claim thus follows from the fact that $d_1 \colon \Q_1(\C,\QF\qshift{1}) \to (\C,\QF\qshift{1})$ restricts to $\met \colon \Met(\C,\QF\qshift{1}) \rightarrow (\C,\QF\qshift{1})$.
\end{proof}

\begin{remark}
In particular, $|\Cob^\F(\C,\QF)|$ is canonically an $\Einf$-group. Its inversion map, however, is generally not induced by the Poincaré functor $(\id_\C,-\id_{\QF})$, see Corollary~\reftwo{corollary:inversion} for the general formula.
 \end{remark}

As the maps 
\[
\met \colon \Met(\Met(\C,\QF)) \lrar \Met(\C,\QF) \quad \text{and} \quad  \met \colon \Met(\Hyp(\C)) \lrar \Hyp(\C)
\]
are split by Remark~\refone{remark:comonad} and Corollary~\refone{corollary:met-hyp-splits}, we obtain:

\begin{corollary}
\label{corollary:pi0CorCobMet}%
For any Poincaré $\infty$-category $(\C,\QF)$, any small stable $\infty$-category $\D$ and any additive functor $\F \colon \Catp \rightarrow \Sps$, the $\infty$-categories
$\Cob^\F(\Met(\C,\QF))$ and $\Cob^\F(\Hyp(\D))$ are connected. 
\end{corollary}

Let us have a closer look at these two cobordism $\infty$-categories. We recorded in Example~\reftwo{example:CobHyp-equal-Span} that the forgetful functor $\Cob(\Hyp(\C)) \rightarrow \Span(\C)$ is an equivalence, so in particular we find:

\begin{observation}
\label{observation:CobHyp=SpanII}%
For every small stable $\infty$-category $\C$ there is a canonical equivalence
\[
|\Cob(\Hyp(\C))| \simeq \Omega^{\infty-1}\K(\C). 
\]
\end{observation}
Here, $\K(\C)$ denotes the connective algebraic $K$-theory spectrum of $\C$, defined for instance through the iterated $\Q$-construction for stable $\infty$-categories. 
\begin{proposition}
\label{proposition:CobMet=CobHyp}%
Ranicki's algebraic Thom construction gives a natural equivalence of $\infty$-categories
$\Cob(\Met(\C, \QF\qshift 1)) \to \Span(\catforms(\C,\QF))$ and the forgetful functor $\Span(\catforms(\C, \QF))\to \Span(\C)$ induces an equivalence 
\[
|\Cob(\Met(\C,\QF))| \simeq \Omega^{\infty-1}\K(\C).
\]
\end{proposition}
\begin{proof}
Commuting diagram categories, we find $\Q(\Met(\C,\QF)) \simeq \Met\Q(\C,\QF)$
so that Proposition~\refone{corollary:algebraic-thom-iso}, our incarnation of the algebraic Thom construction,
implies $\Poinc\Q(\Met(\C,\QF\qshift{1})) \simeq \spsforms\Q(\C,\QF)$.
But, without a non-degeneracy condition, hermitian objects in a diagram category are just diagrams of hermitian objects; see Corollary~\refone{corollary:forms-in-diagram-cats}. So the right hand side is equivalent to $\Core\Q(\catforms(\C,\QF))$. Passing to associated $\infty$-categories gives the first claim.

For the second claim, we will show that the projection $\pi\colon \Span(\catforms(\C,\QF)) \to \Span(\C)$
is cofinal and appeal to \cite{HTT}*{Corollary 4.1.1.12}. 
By \cite{HTT}*{Theorem 4.1.3.1}, it suffices to show that for every $X \in \Span(\C)$, the comma $\infty$-category $\Span(\catforms(\C,\QF))_{X/}$ is contractible.  To see this, we produce an adjunction 
\[
(\catforms(\C,\QF)_{/X})\op \adj \Span(\catforms(\C,\QF))_{X/},
\]
showing that the two sides have the same realisation, and note that $\catforms(\C,\QF)_{/X}$ has an initial object (the zero object of $\C$ with the trivial hermitian structure), so is contractible.
The adjunction itself is a special case of \cite{HLS}*{Lemma 4.6}, since the forgetful functor $\catforms(\C,\QF) \rightarrow \C$ is a right fibration by construction; the left adjoint is induced by the functor $\catforms(\C,\QF)\op \rightarrow \Span(\catforms(\C,\QF))$, which is the identity on objects and takes a morphism $X' \rightarrow X''$ to the span $X'' \leftarrow X' = X'$ %
and the right adjoint takes an element $(X \leftarrow W \rightarrow X', q \in \Omega^\infty \QF(X'))$ to $W \rightarrow X$ with the induced form on $W$.
\end{proof}

The equivalence $|\Cob(\Met(\C, \QF))| \simeq |\Cob(\Hyp(\C))|$ obtained from combining Observation~\reftwo{observation:CobHyp=SpanII} and Proposition~\reftwo{proposition:CobMet=CobHyp}
in fact holds more generally for the $\F$-based cobordism $\infty$-categories as a formal consequence of $|\Cob^\F-|$ being additive and group-like, see Corollary~\reftwo{corollary:FCobhyp=FCobMet}.

\subsection{Algebraic surgery}
\label{subsection:algebraic-surgery}%

In this subsection, we translate Ranicki's algebraic surgery to our set-up. This provides a simple way of producing cobordisms and gives an illuminating description of the comma categories of $\Cob(\C,\QF)$. We will approach these statements by translating them into assertions about certain Segal spaces derived from the $\Q$-construction, and for the present paper it is, in fact, the analysis thereof that plays the largest role. We follow the basic description of algebraic surgery given by Lurie in \cite{Lurie-L-theory}*{Lecture 11}.

Let $(\C, \QF)$ be a Poincaré \(\infty\)-category, and $(X,q)$ be a Poincaré object therein. A \emph{surgery datum} on $(X,q)$ consists of a map $f\colon T\to X$ and a nullhomotopy of $f^*q \in \Omega^\infty\QF(T)$. In other words, it is the extension of $(X,q)$ to a hermitian (but not necessarily Poincaré) nullbordism, i.e.\ to an object of $\catforms(\Met(\C,\QF))$. Surgery data organise into an $\infty$-category:

\begin{definition}
The \defi{$\infty$-category of surgery data in $(\C,\QF)$} is given by
\[
\Surg(\C,\QF) =\catforms(\Met(\C,\QF)) \times_{\catforms(\C,\QF)}  \Poinc(\C,\QF),
\]
where the left hand map in the pullback is induced by $\met \colon \Met(\C,\QF) \rightarrow (\C,\QF)$. The fibre of $\Surg(\C, \QF)$ over some $(X,q) \in \Poinc(\C,\QF)$ is called the \defi{category of surgery data on $(X,q)$} and denoted by $\Surg_{(X,q)}(\C,\QF)$. 

We shall refer to the groupoid cores of these $\infty$-categories as the \emph{spaces of surgery data}.
\end{definition}

\begin{remark}
In geometric topology, a surgery datum on a closed oriented $d$-dimensional manifold $M$ is an embedding $\amalg_{i \in I} S^k \rightarrow M$ with trivialised normal bundle (and $I$ finite). The induced map on singular chains inherits the structure of an algebraic surgery datum in $(\Dperf(\ZZ), {\QF^\sym}\qshift{-d})$ (the Poincaré form on $\mathrm{C}_*(M)$ arises via its identification with $\mathrm{C}^{d-*}(M)$ through Poincaré duality), for example by feeding the trace of the geometric surgery datum into the surgery equivalence of Proposition~\reftwo{proposition:surgequiv} below.

Let us warn the reader that therefore our presentation of algebraic surgery does not follow the overall convention of creating Poincaré chain complexes from manifolds via their cochains; that convention would require us to describe an algebraic surgery datum in a more cumbersome, though equivalent, fashion via a map $X \rightarrow S = \Dual_\QF T$, together with a null-homotopy of the form after pull-back along $\Dual_\QF S \to \Dual_\QF X \simeq X$. 
\end{remark}

Like in the geometric setting, surgery data can be used to produce cobordisms: Given a surgery datum $(f\colon  T\to X, h\colon f^* q\simeq 0)$, the composition
\[
T \xrightarrow f X \xrightarrow{q_\sharp} \Dual_\QF X \xrightarrow{\Dual_\QF f} \Dual_\QF T
\]
is identified with $(f^*q)_\sharp$ and therefore null via $h$. Therefore one can form the following diagram  
\[
\begin{tikzcd}
T \ar[d] \ar[r,equal] & T \ar[d,"f"] \ar[r] & 0 \ar[d] \\
\chi(f) \ar[r] \ar[d] & X \ar[r] \ar[d] & \Dual_\QF T \ar[d,equal] \\
X_f \ar[r] & X/T \ar[r] & \Dual_\QF T
\end{tikzcd}
\]
with exact rows and columns: Here, $\chi(f)$ is the fibre of the composition $X\simeq \Dual_\QF X \xrightarrow{\Dual_\QF f} \Dual_\QF T$ and $X_f$ is defined to be the cofibre of $T\to \chi(f)$. 

The resulting span $[X\leftarrow \chi(f)\to X_f] \in \Q_1(\C)$ is then the underlying object of the desired cobordism and we claim it upgrades to the simplest (but argueably most important) instance of the surgery equivalence, namely:
\begin{proposition}
\label{proposition:surgequiv}%
The association $\chi$ promotes to an equivalence
\[
\chi\colon \grpcr\Surg(\C,\QF) \longrightarrow  \Poinc\Q_1(\C,\QF),
\]
naturally in the Poincar\'e category $(\C,\QF)$.
\end{proposition}
To upgrade $\chi$ in this fashion (and also to extend the discussion to higher simplices in the $\Q$-construction) it is in fact easier to construct the inverse morphism which extracts a surgery datum from a cobordism, a process that seems to have no geometric analogue. Namely, there is a cartesian square
\[
\begin{tikzcd} 
\Q_1(\C,\QF) \ar[r] \ar[d,"d_1"] & \Met(\Met(\C,\QF\qshift{1})) \ar[d, "{\met}"] \\
(\C,\QF) \ar[r] & \Met(\C,\QF\qshift{1}) %
\end{tikzcd}
\]
where the upper horizontal map is  
given by
\[
(X \leftarrow Y \rightarrow Z) \quad \longmapsto \quad \begin{tikzcd} Y \ar[r]\ar[d] & X \ar[d] \\ Z \ar[r] & 0 \end{tikzcd} 
\]
with its tautological Poincar\'e refinement,
and the lower horizontal one is the inclusion as the objects without boundary. Applying the limit preserving functor $\Poinc$ now indeed gives an equivalence $s\colon \Poinc\Q_1(\C,\QF) \rightarrow \grpcr\Surg(\C,\QF)$ via the algebraic Thom isomorphism (see Corollary~\refone{corollary:algebraic-thom-iso}). Unwinding definitions $s$ takes $X \leftarrow W \rightarrow Y$ to $\mathrm{fib}(W \rightarrow Y) \rightarrow X$ at the level of underlying objects and this does indeed give an inverse construction to $\chi$. In this form, surgery generalises to higher simplices of the $\Q$-construction. We need a bit of notation:

\begin{definition}
\label{definition:higher-met}%
Let $(\C,\QF)$ be a Poincaré \(\infty\)-category. We define the simplicial object $\hMet(\C,\QF)$ in $\Catp$ as the fibre of the simplicial map $\dec(\Q(\C,\QF)) \rightarrow \Q_{0}(\C,\QF)=(\C,\QF)$. We will denote by $\hMet_n(\C)$ the underlying stable $\infty$-category of $\hMet_n(\C,\QF)$ (which only depends on $\C$).
\end{definition} 

In particular, $\hMet_0(\C,\QF) \simeq \Met(\C,\QF)$. By Lemma~\reftwo{lemma:boundarysplit}, the sequence
\[
\hMet_n(\C,\QF) \xrightarrow{i_n} \Q_{1+n}(\C,\QF) \xrightarrow{\mathrm{ev}_0} (\C,\QF) ,
\]
defining $\hMet_n(\C,\QF)$ is a split Poincaré-Verdier sequence. 
Here, the fully-faithful left and right adjoints 
$g_n,h_n$ of $d_0$ are given by left and right Kan extension along $\{(0\leq 0)\} \subseteq \Twar\Delta^n$, respectively. 
By Remark~\reftwo{remark:quasi-split}, the left and right adjoints $l_n$ and $r_n$ 
of \(i_n\) are then given by \(l_n(X) = \cof[g_n(X_{0 \leq 0}) \to X]\) and \(r_n(X) = \fib[X \to h_n(X_{0 \leq 0})]\). 
By Corollary~\reftwo{corollary:structure-square} and Remark~\reftwo{remark:explicit-met} we obtain a natural cartesian square
\[
\begin{tikzcd}
\Q_{1+n}(\C, \QF) \ar[r,"\varphi_n"] \ar[d, "\pi"] 
  & \Met\hMet_n(\C, \QF\qshift 1) \ar[d, "\met"] \\
 (\C, \QF) \ar[r, "\psi_n"] 
  & \hMet_n(\C, \QF\qshift 1) 
\end{tikzcd}
\]
classifying this split Poincar\'e-Verdier sequence, generalising the previous one from $n=0$ to arbitrary $n$; here the underlying exact functor of $\psi_n$ sends $Z \in \C$ to the diagram 
\[
\psi_n(Z)_{i \leq j} = l_n h_n(Z)_{i\leq j} = \left\{\begin{matrix}Z & i = 0, j \geq 1 \\ 0 & \text{otherwise.}\end{matrix}\right.
\]
and the underlying exact functor of $\varphi_n$ sends $X$ to the object $[l_n(X) \to \psi_n(X_{0 \leq 0})] \in \Met(\hMet_n(\C))$.

Our next goal is to exhibit the Poincar\'e $\infty$-categories $\Null_n(\C,\QF)$ as a Poincaré $\infty$-category of pairings. To this end we apply Proposition \refone{proposition:recognize-poincare} which requires us to identify a Lagrangian inside $\Null_n(\C,\QF)$; %
before diving in we recall some details regarding this recognition criterion: 

A Lagrangian $\Lag \subseteq \D$ in a Poincar\'e $\infty$-category $(\D,\QFD)$ is a full subcategory on which $\QFD$ vanishes, whose inclusion admits a right adjoint $p$ and such that the resulting inclusion
\[
\Lag \subseteq \Lag^{\perp} \coloneq \{Y \in \D \mid \Bil_\QFD(X,Y) \simeq 0 \text{ for all } Y \in \Lag\}
\]
is an equivalence.
In this case the Verdier quotient $\D/\Lag$ identifies with $\Lag\op$ by means of the functor $q = (p \circ  \Dual_\QFD)\op \colon \D \rightarrow \Lag\op$ resulting in a right split Verdier sequence $\Lag \rightarrow \D \rightarrow \Lag\op$ and consequently an exact counit-unit sequence $pX \rightarrow X \rightarrow \Dual_\QFD p \Dual_\QFD X$ for all $X \in \D$ by Remark \reftwo{remark:quasi-split}, compare Definition \refone{definition:lagrangian} or see \S\reftwo{subsection:isotrop} below for a generalisation.

With this in mind Proposition \refone{proposition:recognize-poincare} combined with Proposition \refone{proposition:first-adj-hermitian} and Corollary \refone{corollary:left-kan} implies:

\begin{observation}
\label{observation:recognizeWolfgang}%
For $(\D,\QFD) \in \Catp, (\E,\QFE) \in \Cath$ and an hermitian functor $(h,\eta) \colon (\D,\QFD) \rightarrow (\E,\QFE)$ the induced Poincar\'e functor $(\D,\QFD) \rightarrow \Pairings(\E,\QFE)$ is an equivalence if and only if 
\begin{enumerate}
\item $h$ is a Verdier projection,
\item $\mathrm{ker}(h) \subseteq \D$ is a Lagrangian, and
\item the induced map $\eta \colon h_! \QFD \rightarrow \QFE$ is an equivalence.
\end{enumerate}
\end{observation}

\begin{remark}
\label{remark:Wolfgangsformula}%
In the situation of Observation~\reftwo{observation:recognizeWolfgang}, the hermitian structure $\QFE \simeq h_!\QFD$ can also be identified naturally as the restriction of $\QFD\qshift{-1}$ along $\Omega\op \circ r\op \colon \E\op \rightarrow \D\op$, where $r$ is the right adjoint to $h$: 
For $X \in \D$ the exact counit-unit sequence $pX \to X \to rhX$ from Lemma~\reftwo{lemma:quasi-split} and Remark \reftwo{remark:quasi-split} provides a comparison map
\[
\QFD(X) \lrar \fib[\QFD(pX) \lrar \QFD(\Om rhX)] = \QFD\qshift{-1}(\Om rhX)
\]
and that it exhibits $\QFD\qshift{-1} \circ \Omega\op \circ r\op$ as the left Kan extension of $\QFD$ along $h$ can be seen as follows. By Lemma~\refone{lemma:goodwillie} applied to the rotated fibre sequence $\Omega rhX \to pX \to X$ there is a natural exact sequence
\[
 \Bil_{\QFD}(X,rhX) \longrightarrow \QFD(X) \longrightarrow \QFD\qshift{-1}(\Om rhX).
\]
Since $h$ is a localisation the left Kan extension of the middle term is the desired composite and the right hand side identifies with $\map_\D(X,\Dual_{\QFD}(rhX))$. The left Kan extension of this last term along $h$ is trivial: Since the slice categories occuring in the pointwise formula for this Kan extension are sifted (even filtered) it can be computed separately in both variables. But the left Kan extension in the first variable is given by $\map_\E(-, h\Dual_\QFD(rhX))$ and $h\Dual_{\QFD}(rhX)$ vanishes on account of the identification $\Dual_\QFD p\op \Dual_\QFD\op \simeq rh$ obtained by comparing the counit-unit sequence with the one coming from $\mathrm{ker}(h)$ being a Lagrangian.

Using this description of $h_!\QFD$ the functor underlying the equivalence $(\D,\QFD) \rightarrow \Pairings(\E,h_!\QFD)$ explicitly takes $X \in \D$ to $(hX,h\Dual_\QFD X, \partial)$ where $\partial$ is the boundary map of the counit-unit sequence using
\[
\Omega^\infty \Bil_{h_!\QFD}(hX,h\Dual_\QFD X) \simeq \Omega^{\infty+1} \Bil_\QFD(\Omega rhX, \Omega rh\Dual_\QFD X) \simeq \Omega^{\infty} \Bil_\QFD(rhX, \Omega \Dual_\QFD p X) \simeq  \Hom_{\D}(rhX, \Sigma pX).
\]
\end{remark}

Now, back to the matter at hand: The requisite Lagrangian in $\Null_n(\C,\QF)$ is given by the full subcategory spanned by those diagrams $X \colon \Twar(\Delta^{1+n}) \rightarrow \C$ which take all down-right maps \((i \leq j) \to (i' \leq j)\) to equivalences. Evidently, $\QF_{1+n}$ vanishes on this subcategory and to carry out the remaining analysis, let us, for a category \(\I\) and an object \(i \in \I\), denote by \(\C^{(\I,i)} \subseteq \C^\I\) the full subcategory spanned by those \(\I \to \C\) which send \(i\) to \(0\). Then the source and target projections \((s_n,t_n) \colon \Twar\Del^n \rightarrow \Del^n \times (\Del^n)\op\) send \((0 \leq 0)\) to \((0,0)\) and hence induce diagrams 
\[
\begin{tikzcd}
(\C,\QF)^{((\Del^{1+n})\op,0)} \ar[r, "\tau_n"]\ar[d, hook] & \hMet_n(\C,\QF) \ar[d, hook] & \ar[l,"\sig_n"'] (\C,\QF)^{(\Del^{1+n},0)} \ar[d, hook] \\
(\C,\QF)^{(\Del^{1+n})\op} \ar[r,"t^*_n"] & \dec_n(\Q(\C,\QF)) & \ar[l,"s^*_n"'] (\C,\QF)^{\Del^{1+n}}
\end{tikzcd}
\]
in $\Cath$, natural in $[n] \in \Delta\op$, where we note that any diagram \(\Twar\Del^{1+n} \to \C\) obtained via restriction along either the target or source projection automatically lies in \(\Q_{1+n}(\C)\). Then $\Lag_n(\C,\QF) \subseteq \Null_n(\C,\QF)$ is the essential image of $\tau_n$ and:

\begin{lemma}
\label{lemma:null-metabolic}%
The functors \(\tau_n\) and \(\sig_n\) admit adjoints fitting into a right split Verdier sequence
\[
\begin{tikzcd}
\C^{((\Del^{1+n})\op,0)} 
\ar[r,"\tau_n","\myperp"'] 
& \Null_n(\C)
\ar[l,"p_n",bend left=30,shift left=1.5ex,start anchor=west,end anchor=east]
\ar[r,"q_n","\myperp"'] 
& \C^{(\Del^{1+n},0)}
\ar[l,"\sig_n",bend left=30,shift left=1.5ex,start anchor=west,end anchor=east]
\end{tikzcd}
\]
compatibly with the simplicial transition functors. The adjoints %
are given by the formulas \(p_n(X)_j = X_{0 \leq j}\) and \(q_n(X)_i = \cof[X_{0\leq n+1} \to X_{i \leq n+1}]\).%
\end{lemma}

Note that the formula for $q_n$ is only a pointwise one (in $n \in \Delta\op$) and not obviously natural. Nevertheless, the lemma provides a natural refinement as $\tau_n$ and $p_n$ are evidently natural (compare Lemma \reftwo{lemma:adj-verdier}).

\begin{proof}
The source and target projections $s_n$ and $t_n$ both admit fully faithful adjoints and are hence localisations, e.g.\ \(t_n\) admits a left adjoint given by \(b_n(j) = (0 \leq j)\). It follows that the restriction functors \(t_n^*\) and \(s_n^*\) are both fully faithful, and their essential images consist of those diagrams \(X\colon \Twar\Del^n \to \C\) which send the down-right arrows \((i \leq j) \to (i' \leq j)\) to equivalences in the case of \(t_n^*\) and the down-left arrows \((i \leq j) \to (i \leq j')\) to equivalences in the case of \(s_n^*\). 

The functor \(b_n^*\) further provides a right adjoint to \(t_n^*\). 
Since \(b_n^*\) sends \(\hMet_n(\C)\) to
\(\C^{((\Del^{[1+n]})\op,0)}\) it restricts to give the right adjoint 
$p_n=b_n^*\colon \hMet_n(\C) \to \C^{((\Del^{1+n})\op,0)}$ 
to \(\tau_n\), visibly compatible with the simplicial transition functors.
Now the exactness conditions defining \(\Q_{1+n}(\C)\) and a 2-out-of-3 argument show that an \(X \in \Q_{1+n}(\C)\) sends all down-left arrows to equivalences if and only if the maps \(X_{0 \leq j} \to X_{0 \leq 0}\) are equivalences for every \(j\). For \(X \in \hMet_n(\C)\), this is equivalent to saying that \(X_{0 \leq j} = 0\) for every \(j\), that is, that \(p_n(X) = 0\). We conclude that 
\(\im(\sig_n)\) is the right orthogonal complement of \(\im(\tau_n)\), and so the rest of the statement is now a formal consequence of 
Lemma~\reftwo{lemma:quasi-split} and Lemma~\reftwo{lemma:adj-verdier}.
\end{proof}

Combined with Remark \reftwo{remark:quasi-split}, Lemma~\reftwo{lemma:null-metabolic} furnishes natural exact sequences $\tau_np_n(X) \rightarrow X \rightarrow \sigma_nq_n(X)$ and thus maps
\[
\eta_n\colon \QF_{1+n}(X) \lrar \fib\big[\QF_{1+n}(\tau_np_n(X)) \to \QF_{1+n}(\Om\sig_nq_n(X))\big] = \QF_{1+n}\qshift{-1}(\Om\sig_nq_n(X))
\]
refining \(\Om q_n\), naturally in \([n] \in \Del\op\), to a hermitian functor 
\[
(\Om q_n,\eta_n)\colon \hMet_n(\C,\QF) \lrar \left(\C^{(\Del^{1+n},0)},\sig_n^*\QF_{1+n}\qshift{-1}\right) \simeq (\C,\QF\qshift{-1})^{\Delta^n}
\]
where the second equivalence is by means of $d_0 \colon \Delta^n \rightarrow \Delta^{1+n}$. We thus find:

\begin{proposition}
\label{proposition:multiple-algebraic-Thom}%
The subcategory \(\Lag_n(\C,\QF) = \im(\tau_n)\) 
is a Lagrangian in \(\hMet_n(\C,\QF)\) and the above hermitian functor adjoins to an equivalence
\[
\hMet_n(\C,\QF) \longrightarrow %
\Pairings\left((\C,\QF\qshift{-1})^{\Del^n}\right)
\]
natural in $[n] \in \Del\op$.
\end{proposition}
\begin{proof}
Given Lemma \reftwo{lemma:null-metabolic} and the discussion preceeding it, to see that $\Lag_n(\C,\QF)$ is a Lagrangian it only remains to check that the right orthogonal of $\Lag_n(\C,\QF)$ agrees with $\Dual_{\QF_{1+n}}(\Lag_n(\C,\QF))$. But the latter consists exactly of those diagrams all of whose down-left maps are equivalences, so both agree with the essential image of $\sigma_n$. 

The remaining claim now follows from the recognition principle in Observation \reftwo{observation:recognizeWolfgang} together with Remark \reftwo{remark:Wolfgangsformula} applied to the Verdier projection $\Omega q_n$.
\end{proof}

In total the composite functor 
\[
(\C,\QF) \xrightarrow{\psi_n} \hMet_n(\C,\QF\qshift{1}) \xrightarrow{(\Om q_n,\eta_n)} (\C,\QF)^{\Del^n}
\]
unwinds to the functor that simply takes $Z \in \C$ to the constant diagram $\Delta^n \rightarrow \C$ with value $Z$, with its tautological hermitian refinement, so we find the following categorical version of the surgery equivalence:

\begin{corollary}
\label{proposition:catsurgequiv}%
We have the following Poincar\'e-Verdier square, natural in \((\C,\QF)\) and \([n] \in \Del\op\). 
\[
\begin{tikzcd}
\Q_{1+n}(\C, \QF) \ar[r, "s_n"] \ar[d, "\pi"] & \Met\Pairings((\C, \QF)^{\Delta^n})  \ar[d, "\met"] \ar[r,"\simeq"]&  \Pairings(\Met(\C, \QF)^{\Delta^n})\ar[d, "\Pairings(\met)"] \\
(\C, \QF) \ar[r, "\const"] & \Pairings((\C, \QF)^{\Delta^n}) \ar[r,equal] & \Pairings((\C, \QF)^{\Delta^n})
\end{tikzcd}
\]
\end{corollary}

Here, we abuse notation slightly in the upper right corner, by viewing the metabolic construction as an endofunctor on hermitian (as opposed to just Poincar\'e) $\infty$-categories. The top right horizontal functor is then obtained from the hermitian functor $\Met(c) \colon \Met\Pairings((\C,\QF)^{\Delta^n}) \rightarrow \Met((\C,\QF)^{\Delta^n})$ by the universal property of $\Pairings \colon \Cath \rightarrow \Catp$ as a right adjoint, see Proposition \refone{proposition:first-adj-hermitian}, with $c$ being the counit. 

The composite functor $\Met(c) \circ s_n$, which adjoins to the composite horizontal functor is per construction given by
\[
\Q_{1+n}(\C,\QF) \xrightarrow{\varphi_n} \Met\hMet_n(\C,\QF\qshift{1}) \xrightarrow{\Met(\Om q_n,\eta_n)} \Met((\C, \QF)^{\Delta^n})
\]
or explicitly by
\[
X \longmapsto \left(
  \begin{tikzcd} 
  \fib[X_{0 \leq n+1} \to X_{1 \leq n+1}] \ar[r]\ar[d] & \dots \ar[r]\ar[d] & \fib[X_{0 \leq n+1} \to X_{n+1 \leq n+1}] \ar[d]
  \\ X_{0 \leq 0} \ar[r] & \dots \ar[r]   & X_{0 \leq 0}
  \end{tikzcd}
 \right)
\]
on underlying categories, upgraded to a hermitian functor via the composite natural transformation
\[
\lim_{\Twar(\Delta^{1+n})\op} \QF\circ X \longrightarrow  \QF(X_{0 \leq 0}) \times_{\QF(X_{0 \leq n+1})}\QF(X_{n+1 \leq n+1}) \longrightarrow \fib\big[\QF(X_{0 \leq 0}) \to \QF(\fib[X_{0 \leq n+1} \to X_{n+1 \leq n+1}])\big],
\]
though again, this description is not visibly natural in $n \in \Delta\op$.

\begin{proof}[Proof of Corollary~\reftwo{proposition:catsurgequiv}]
It remains to argue that the top right horizontal map is an equivalence. This also follows from Observation~\reftwo{observation:recognizeWolfgang}: %
With $c \colon \Pairings((\C,\QF)^{\Delta^n}) \rightarrow (\C,\QF)^{\Delta^n}$ also $\Met(c)$ is a Verdier projection (e.g.\ by Proposition \reftwo{proposition:(co)tensor-Verdier}) whose kernel is a Lagrangian (by a simple direct check), and since $\QF^{\Delta^n}$ is left Kan extended from $(\QF^{\Delta^n})_\pair$ along $c\op$ by \reftwo{proposition:multiple-algebraic-Thom}, also $\QF^{\Delta^n}_\met$ is left Kan extended from $(\QF^{\Delta^n}_\pair)_\met$ along $\Met(c)\op$ (by the description of the induced Poincar\'e structures in \reftwo{remark:Wolfgangsformula}).	
\end{proof}

\medskip

Applying an additive functor $\F\colon \Catp\to \Sps$ to the diagram in Corollary \reftwo{proposition:catsurgequiv}, we obtain a cartesian diagram of (not necessarily complete) Segal spaces: Since $[n] \mapsto \Pairings((\C,\QF)^{\Delta^n})$ is a (complete) Segal object in $\Catp$, this follows from Lemma \reftwo{proposition:multiple-algebraic-Thom} together with Proposition \reftwo{lemma:boundarysplit}.

\begin{definition}
For additive $\F \colon \Catp \rightarrow \Sps$ and a Poincar\'e $\infty$-category $(\C,\QF)$ denote by $\catforms^\F(\C,\QF)$ the $\infty$-category associated to $[n] \mapsto \F(\Pairings((\C,\QF)^{\Delta^n}))$
and define the $\infty$-category of $\F$-based surgery data in $(\C,\QF)$ as
\[
\Surg^\F(\C,\QF) = \catforms^\F(\Met(\C,\QF)) \times_{\catforms^\F(\C,\QF)} \F(\C,\QF).
\]
For $X \in \F(\C,\QF)$ we denote by $\Surg^\F_X(\C,\QF)$ the fibre of $\Surg^\F(\C,\QF)$ over $X$.
\end{definition}

This notation is justified by the equivalence
$\snerv(\catforms(\C,\QF)) \simeq 
\spsforms((\C,\QF)^{\Delta^-}) \simeq \Poinc\Pairings((\C,\QF)^{\Delta^-})$,
via the algebraic Thom construction, which implies $\catforms^\Poinc(\C,\QF) \simeq \catforms(\C,\QF)$ and thus  
$\Surg^\Poinc(\C,\QF) \simeq \Surg(\C,\QF)$.

\begin{corollary}[Surgery equivalence]
\label{corollary:F-based-surgery-equivalence}%
For every additive $\mathcal F \colon \Catp \rightarrow \Sps$, Poincar\'e $\infty$-category $(\C,\QF)$, and $X \in \F(\C,\QF\qshift 1)$, the transformation $s\colon \dec(\Q(\C,\QF)) \to \Pairings(\Met(\C,\QF)^{\Delta^-})$ induces an equivalence
\[
s^\F_X\colon \Cob^\F(\C,\QF)_{X/} \longrightarrow \Surg^\F_X(\C,\QF).
\]
\end{corollary}
\begin{proof}
By Corollary~\reftwo{proposition:catsurgequiv} and Lemma~\reftwo{lemma:Segal_slice_asscat}, the category $\Cob^\F(\C, \QF)_{/X}$ is described as the associated category of the pullback of 
\[
 \F(\Pairings(\Met(\C, \QF\qshift 1)^{\Delta^-})) \lrar
 \F(\Pairings((\C, \QF\qshift 1)^{\Delta^-})) \longleftarrow
 \F(\C, \QF\qshift 1) \longleftarrow \{X\}.
\]
It remains to be shown that the pullback can be commuted with forming associated categories. This trivially holds true whenever the Segal spaces involved are complete, such as in the main case of interest $\F = \Poinc$. The general case of additive $\F$ requires the additional input that the left hand map is an isofibration of Segal spaces, see Lemma~\reftwo{lemma:asscat_pullbacks}. This can be checked directly, but for sake of brevity we deduce it from the additivity theorem in the next section (whose proof does not make use of surgery arguments): The map $\F(\Pairings(\Met(\C,\QF)^{\Delta^-})) \rightarrow \F(\Pairings((\C,\QF)^{\Delta^-}))$ is identified with $\F(\hMet\Met(\C,\QF)) \rightarrow \F(\hMet(\C,\QF))$ by \reftwo{proposition:catsurgequiv} and this is in turn obtained from the map $\F\Q\Met(\C,\QF) \rightarrow \F\Q(\C,\QF)$ by taking Segal slices under $0$. Theorem~\reftwo{theorem:fulladditivity} implies that this map is a cartesian fibrations of Segal spaces and cartesian fibrations are generally preserved under passage to Segal slices. Since cartesian fibrations are always isofibrations the claim follows.%
\end{proof}

\begin{remark}
\label{remark:underlyingpairingcotensordelta}%
The underlying stable $\infty$-category of $\Pairings((\C,\QF)^{\Delta^n})$ can be identfied naturally in $n \in \Delta\op$ and $(\C,\QF) \in \Catp$ as $\Fun(\Del^n \ast (\Del^n)\op,\C)$: Using that $\Hyp \colon \Catx \rightarrow \Catp$ is left adjoint to the forgetful functor and $\Pairings$ is right adjoint to the forgetful functor $\Catp \to \Cath$ we compute 
\begin{align*}
\Hom_{\Catx}(\E, \Pairings((\C,\QF)^{\Delta^n}) =& \Hom_{\Catp}(\Hyp(\E),\Pairings((\C,\QF)^{\Delta^n}))\\
=& \Hom_{\Cath}(\Hyp(\E), (\C,\QF)^{\Delta^n}) \\
=& \Hom_{\Cat}(\Del^n,\catforms(\Funx(\Hyp\E, (\C,\QF)))) \\
=& \Hom_{\Cat}(\Del^n,\catforms(\Hyp\Funx(\E, \C))) \\
=& \Hom_{\Cat}(\Del^n,\Twar(\Funx(\E,\C))) \\
=& \Hom_{\Cat}(\Del^n \ast (\Del^n)\op,\Funx(\E,\C)) \\
=& \Hom_{\Catx}(\E,\Fun(\Del^n \ast (\Del^n)\op,\C)) 
\end{align*}
for every $\E \in \Catx$, naturally in all three input variables; here we have used the natural equivalence $\Funx(\Hyp(\E),\C) = \Hyp\Funx(\E,\C)$ resulting from the projection formula of Corollary~\refone{corollary:hyp-projection} together with the fact that $\Hyp$ is also right adjoint to the forgetful functor, as well as the natural equivalence $\catforms\Hyp(-) \simeq \Twar(-)$, see \S\refone{subsection:hyp-and-sym-poincare-objects}. 

Under this equivalence the Poincar\'e structure on $\Pairings((\C,\QF)^{\Delta^n})$ translates (pointwise in $n \in \Delta$) to the functor assigning to $X_0 \rightarrow \dots \rightarrow X_n \xrightarrow{\beta} Y_n \rightarrow \dots \rightarrow Y_0$ the pullback of
\[
\hom_{\C^{\Delta^n}}(X,\Dual_\QF \circ Y) \xrightarrow{\Dual_\QF(\beta)_\ast} \hom_\C(X_n,\Dual_\QF(X_n)) \longleftarrow \QF(X_n).
\]
In fact, testing against a Poincar\'e $\infty$-category $(\E,\Psi)$ one can naturally identify $\Pairings((\C,\QF)^{\Delta^n})$ as the cotensor construction from \cite{Spitzweckreal} of $(\C,\QF)$ with $(\Delta^n \ast (\Delta^n)\op,\mathrm{flip})$, a (non-stable) category with duality, see \S\reftwo{subsection:HSV-compare} for a discussion thereof.

\end{remark}
\begin{examples}\
\begin{enumerate}
\item[i)] From Remark \reftwo{remark:underlyingpairingcotensordelta} we find that $\catforms^\core(\C,\QF)$ is associated to the complete Segal space $[n] \mapsto \Hom_{\Cat}(\Delta^n \ast (\Delta^n)\op,\C)$ resulting in an equivalence $\catforms^\core(\C,\QF) \simeq \Twar(\C)$. 
Unwinding the definitions, the surgery equivalence for $\mathcal F = \core$ then reduces to the equivalence
\[
\Span(\C)_{X/} \longrightarrow \Twar(\C_{/X}), \quad 
   (X \leftarrow Y \rightarrow Z) \longmapsto (\fib(Y \rightarrow Z) \rightarrow Y \rightarrow X) ,
\]
which is in fact valid for any stable $\infty$-category, and should be regarded as its non-hermitian analogue. %
\item[ii)] If $\F$ is group-like and additive we have $\catforms^\F(\C,\QF) \simeq \ast$ for all Poincar\'e $\infty$-categories as a consequence of the isotropic decomposition principle \reftwo{theorem:lagrangian} below, see Remark \reftwo{remarks:FQcompletebord}, and consequently $\Surg^\F(\C,\QF) \simeq \ast$. This shows that $\Cob^\F(\C,\QF)_{X/}$ is trivial in these cases, and thus in turn that $\Cob^\F(\C,\QF)$ is an $\infty$-groupoid; this is also a direct consequence of the work going into Remark~\reftwo{remarks:FQcompletebord}, see Corollary \reftwo{corollary:CobgrpdifFgrp}.
\item[iii)] The precursor to the surgery equivalence in Proposition \reftwo{proposition:surgequiv} has no analogue for general additive $\F$. More precisely, the induced map $\F\Q_1(\C,\QF) \to \iota\Surg^\F(\C,\QF)$ is an equivalence if and only if the Segal space $\F\dec(\Q(\C,\QF))$ is complete,
and this is generally not true, e.g., for $\F = \mathcal K$; compare also Remark \reftwo{remarks:FQcompletebord}.
\end{enumerate}
\end{examples}

We conclude this section by also giving an explicit description of the surgery functor $\chi_X^\F \colon \Surg^\F_X \rightarrow \Cob^\F(\C,\QF)_{X/}$ inverse to $s_X^\F$ at the level of underlying objects. Again, this functor is obtained as the functor associated to a map of Segal spaces, which, at level $n$ is induced by the Poincar\'e functor 
\[
\chi_n\colon \Pairings \Met((\C, \QF)^{\Delta^n})\times_{\Pairings((\C, \QF)^{\Delta^n})}(\C, \QF)\to \Q_{1+n}(\C, \QF)
\]
sending $(X\to \const_Z, Y\to \const_{\Dual_\QF(Z)}, b)\in \Pairings \Met((\C, \QF)^{\Delta^n})$ to the diagram
\[
(i\leq j) \longmapsto 
 \cof\bigl(X_{i-1} \to X_n \xrightarrow{s(b)} \fib(Z\to \Dual_\QF(Y_n)) \to \fib(Z\to \Dual_\QF(Y_{j-1}))\bigr)
\]
 where we set $X_{-1} = 0 = Y_{-1}$

\begin{lemma}
\label{lemma:inverse_surgery_equivalence}%
The functor $\chi$ induces the inverse of the equivalence $s_X^\F$ from Corollary~\reftwo{corollary:F-based-surgery-equivalence}. In particular, the induced equivalence of Poincar\'e objects
 
\[
\spsforms\Met(\C,\QF)^{\Delta^n}\times_{\spsforms(\C, \QF)^{\Delta^n}} \Poinc(\C, \QF) \longrightarrow \Poinc(\Q_{1+n}(\C, \QF))
\]
 sends a hermitian object $(X\to \const_Z, q)$ to the diagram 
 
\[
(i\leq j)\mapsto \cof\bigl(X_{i-1} \to X_n \xrightarrow{s(q_\sharp)} \fib(Z\to \Dual_\QF(X_n)) \to \fib(Z\to \Dual_\QF(X_{j-1}))\bigr)
\]
\end{lemma}

\begin{proof}
The second statement follows because the equivalence $\spsforms(\C, \QF)\to \Poinc(\Pairings(\C, \QF))$ sends $(X,q)$ to a Poincaré object with underlying object $(X, \Dual_\QF X, b)$ with $b\in \Bil_\QF(X, \Dual_\QF X)$ the canonical pairing. For the first statement, we use that the inverse equivalence (of stable $\infty$-categories) $\Pairings(\Dual_\QF \Lag\op, \QF\qshift{-1})\to \C$, for a metabolic Poincaré category $(\C,\QF)$ with Lagrangian $\Lag$, is explicitly given by the formula
\[
(X,Y,b)\mapsto \cof(b\colon X\to \Dual_\QF(Y));
\]
this follows straight from the proof of the equivalence in Proposition \refone{proposition:recognize-poincare}. Applied to the situation at hand, we find that the equivalence of stable $\infty$-categories \(\Pairings (\C, \QF)^{\Delta^n}\xrightarrow{\simeq} \hMet_n(\C)\)
is given by sending $(X, Y, b)$ to the diagram
\[
(i\leq j)\mapsto \cof(X_{i-1}\to X_{n} \xrightarrow b \Dual_\QF(Y_{n}) \to \Dual_\QF(Y_{j-1})),
\]
where we set again $X_{-1} = Y_{-1}=0$. 
The equivalence 
\[
\Pairings \Met((\C, \QF)^{\Delta^n})\xrightarrow{\simeq} \hMet_n(\Met(\C,\QF)) = \Met(\hMet_n(\C,\QF))
\]
is then given by applying the same formula to $\Met(\C,\QF)$ instead of $\C$. When applied to objects of the form $(X\to \const_Z, Y\to \const_{\Dual_\QF(Z)}, b)\in \Pairings \Met((\C, \QF)^{\Delta^n})$,
the outcome is a functor $\Twar(\Delta^{1+n})\to \Met(\C)$ described by the formula 
\[
(i\leq j) \mapsto \bigl((\cof(X_{i-1}\to \fib((\const_Z)_{i-1}\to \Dual_\QF(Y_{j-1})) \to \cof((\const_Z)_{i-1} \to (\const_Z)_{j-1})\bigr),
\]
where we set again $(\const_Z)_{-1}=0$. 
At the same time, the equivalence
\[
\Q_{1+n}(\C,\QF) \simeq (\C,\QF) \times_{\hMet_n(\C,\QF\qshift 1)} \Met\hMet_n(\C, \QF\qshift 1)
\]
induced by the classifying square of the split Poincar\'e-Verdier sequence defining $\hMet_n(\C,\QF\qshift{1})$, see the discussion after Definition \reftwo{definition:higher-met}, is given on underlying stable $\infty$-categories by \(X \mapsto (X_{0 \leq 0},[l_n(X) \to r_n(X_{0 \leq 0})])\) and its inverse equivalence is given by
\((Z,Y \to R_n(Z)) \mapsto h_n(Z) \times_{R_n(Z)} Y$. Composing this inverse with the formula above the desired result follows.
\end{proof}

\subsection{The additivity theorem}
\label{subsection:additvity}%

As we will see, the decisive step towards understanding the homotopy type of the cobordism $\infty$-categories $\Cob(\C,\QF)$ consists in analysing their behaviour under split Poincaré-Verdier sequences. To this end, we show:

\begin{theorem}[Additivity]
\label{theorem:additivity}%
Let $\F\colon \Catp\to\Sps$ be additive. Then the functor  $\vert \Cob^\F\vert$ is also additive. In particular, a split Poincaré-Verdier sequence 
$
(\C,\QF) \to (\D,\QFD) \to (\E,\QFE)
$
induces a fibre sequence
\[
|\Cob^\F(\C,\QF)| \lrar |\Cob^\F(\D,\QFD)| \lrar |\Cob^\F(\E,\QFE)|
\]
of $\Einf$-groups.
\end{theorem}

We heavily exploit this result in \S\reftwo{section:structure-theory} below. In particular, we use it to compute $\pi_1|\Cob^\F(\C,\QF)|$, produce deloopings of $|\Cob^\F(\C,\QF)|$ via the iterated $\Q$-construction, and it serves as the basis for Grothendieck-Witt theory in \S\reftwo{section:GW}. 
It also contains Waldhausen's additivity theorem for $\K$-theory as a special case, as we will detail in \S\reftwo{subsection:Kadd} below. 

On the other hand, Theorem~\reftwo{theorem:additivity} yields an algebraic analogue of Genauer's fibre sequence from geometric topology. To explain this analogy recall that there exists a fibre sequence
\[
|\Cob_{d+1}| \lrar |\Cob^\partial_{d+1}| \lrar |\Cob_{d}|
\]
relating cobordism categories of manifolds of different dimension (with the middle term allowing objects to have boundary). As mentioned in the introduction, this was originally proven by identifying the sequence term by term with the infinite loop spaces of certain Thom spectra, together with a direct verification that these Thom spectra form a fibre sequence; see \cite{Genauer}*{Proposition 6.2} and the main result of \cite{GMTW}. 

The connection to additivity arises by applying Theorem~\reftwo{theorem:additivity} for $\F = \Poinc$ to the metabolic Poincaré-Verdier sequence
\[
(\C,\QF\qshift{-1}) \lrar \Met(\C,\QF) \lrar (\C,\QF)
\]
from Example~\reftwo{example:metabolicfseq}, as we obtain the following algebraic analogue of the Genauer fibre sequence:

\begin{corollary}
\label{corollary:AlgGenauerSequence}%
For every Poincaré \(\infty\)-category $(\C,\QF)$, the metabolic Poincar\'e-Verdier sequence induces the fibre sequence
\[
|\Cob(\C,\QF)| \lrar |\Cob^\partial(\C,\QF)| \lrar |\Cob(\C,\QF\qshift{1})|
\]
of $\Einf$-groups.
\end{corollary}

Even more, our proof of the Additivity theorem will follow the strategy developed in \cite{Steimleadd} by the ninth author in his approach to Genauer's fibre sequence. It is based on a recognition criterion for realisation fibrations \cite{Steimleadd}*{Theorem 2.11}, which we shall give a brief account of in the proof of the Additivity theorem below, see also \cite{Steinebrunner} for a generalisation. The following result verifies its assumptions:

\begin{theorem}
\label{theorem:fulladditivity}%
Let $\F \colon \Catp \rightarrow \Sps$ be additive and $(p,\eta) \colon (\D,\QFD) \rightarrow (\E,\QFE)$ a split Poincaré-Verdier projection. Then the induced map
\[
(p,\eta)_* \colon \F\Q(\D,\QFD) \longrightarrow \F\Q(\E,\QFE)
\]
is a bicartesian fibration of Segal spaces and consequently
\[
(p, \eta)_*\colon \Cob^\F(\D, \QFD) \longrightarrow \Cob^\F(\E, \QFE)
\]
a bicartesian fibration of $\infty$-categories.
\end{theorem}

The proof of Theorem~\reftwo{theorem:fulladditivity} will show that an edge in $\F\Q(\D,\QFD)$ is $\F\Q(p)$-cocartesian if and only if it lies in the image of $\F(\C_p,\QF_1)$ where $\C_p \subseteq \Q_1(\D)$ is the subcategory spanned by those diagrams $X \leftarrow W \rightarrow Y$ with left hand map $p$-cartesian and right hand map $p$-cocartesian; the roles are reversed for $\Q(p)$-cartesian edges.

\begin{remarks}
\begin{enumerate}
\item A similar result in the context of $\infty$-categories of spans was given by Barwick as part of his unfurling construction in \cite{Barwick-MackeyI}*{Theorem 12.2}, but see \cite{HHLN-two-var-fil}*{Remark 3.3} for a small correction. While the main motivation for that construction is also $\K$-theoretic in nature, its use does not seem at all related to additivity in Barwick's work. Our proof, furthermore, proceeds rather differently than Barwick's combinatorial approach.
\item Neither Theorem~\reftwo{theorem:additivity} nor Theorem~\reftwo{theorem:fulladditivity} remain true upon assuming $\F$ Verdier-localising and the input Poincaré-Verdier, but not necessarily split. For example, with $\F = \Kspace \circ (-)^\natural$, which is Karoubi-localising and group-like, Corollary~\reftwo{corollary:components-of-cob-II} and Theorem~\reftwo{theorem:suspension} below in combination show that $|\Cob^{\Kspace \circ (-)^\natural}| \simeq \mathrm{B} \Kspace \circ (-)^\natural$ is the connected delooping, which is famously not (Poincaré-)Verdier-localising, since Verdier projections need not induce surjections on $\K_0\circ (-)^\natural$. 
\end{enumerate}
\end{remarks}

\begin{proof}[Proof of the Additivity theorem, assuming Theorem~\reftwo{theorem:fulladditivity}]
Suppose given a split Poincaré-Verdier square
\[
\begin{tikzcd}
	(\D,\QFD) \ar[r] \ar[d] & (\D',\QFD') \ar[d] \\
        (\E,\QFE) \ar[r] & (\E',\QFE').
\end{tikzcd}
\]
Since $\F\Q \colon \Catp \rightarrow \sCatp$ is additive, we find an associated cartesian square of Segal spaces, and it follows from Lemma~\reftwo{lemma:asscat_pullbacks}, that also
\[
\begin{tikzcd}
{\Cob^\F(\D,\QFD)} \ar[r] \ar[d] & {\Cob^\F(\D',\QFD')} \ar[d] \\
{\Cob^\F(\E,\QFE)} \ar[r] & {\Cob^\F(\E',\QFE')}
\end{tikzcd}
\]
is cartesian, since cartesian fibrations are isofibrations.

With the square above established as cartesian, we next argue that it remains cartesian after realisation. Writing $\mathrm{Un}(G)$ for Lurie's cocartesian unstraightening of a functor $G \colon \C \rightarrow \Cat$, it is generally true that if $|-| \circ G \colon \C \rightarrow \Cat \rightarrow \Sps$ factors over $|\C|$, the diagram
\[
\begin{tikzcd}
{|\mathrm{Un}(GH)|} \ar[r] \ar[d] & {|\mathrm{Un}(G)|} \ar[d] \\
            {|\C'|} \ar[r] & {|\C|}
\end{tikzcd}
\]
is cartesian for any $H \colon \C' \rightarrow \C$: Indeed, observe, for example via \cite{HTT}*{Corollary 3.3.4.3}, that the upper terms can also be regarded as the unstraightening of $|G| \colon |\C| \rightarrow \Sps$ and its precomposition with $|H|$, whence this follows from unstraightening translating compositions to pullbacks (by the adjoint of \cite{HTT}*{Proposition 3.2.1.4}).
Taking for $G$ the cocartesian straightening of $\Cob^\F(\D',\QFD') \rightarrow \Cob^\F(\E',\QFE')$ and for $H$ the functor $\Cob^\F(\E,\QFE) \rightarrow \Cob^\F(\E',\QFE')$ then gives the claim, since the cocartesian straightening of a bicartesian fibration takes values in $\Cat^\mathrm L$, so satisfies the assumption on $G$, as adjoint functors realise to equivalences.
\end{proof}

\subsection{Fibrations between cobordism categories}
\label{subsection:fibcob}%

The present section is devoted to the proof of Theorem~\reftwo{theorem:fulladditivity}. 
The strategy is as follows: After recording that a split Verdier projection (of stable $\infty$-categories) is a bicartesian fibration, we improve on this by showing that the maps
\[
p_*\colon (\D,\QFD)^{\Delta^n} \lrar (\E,\QFE)^{\Delta^n}
\]
behave like a bicartesian fibration between Segal objects in $\Cath$; we will not give a formal definition of this term, but instead formulate the relevant statements directly in Lemmas~\reftwo{lemma:cartesian-lifts-in-Cath} and \reftwo{lemma:cartesian-lifts-2-in-Cath}. We then use this to show that the map 
\[
\Q(p)\colon \Q(\D,\QFD)\lrar \Q(\E,\QFE)
\]
also behaves like such a bicartesian fibration; the cocartesian part is formulated in Lemmas~\reftwo{lemma:cocartesian-lifts-1-in-CobF} and \reftwo{lemma:cocartesian-lifts-2-in-CobF}, and the cartesian one follows by invariance of the $\Q$-construction under taking opposites. From there we will deduce the theorem by observing that any additive functor $\F$ can be used as a `cut-off' to obtain a bicartesian fibration $\F\Q(\D,\QFD) \rightarrow \F\Q(\E,\QFE)$ of Segal objects in $\Sps$, which implies the result.\\

To get started we need:

\begin{lemma}
\label{lemma:cartesian-in-bousfield}%
Let $p\colon \C\to \Ctwo$ be a functor with left adjoint $g$. Then:
\begin{enumerate}
 \item
\label{item:cartesian-morphism-in-bousfield}%
A morphism $\alpha\colon X\to Y$ in $\C$ is $p$-cocartesian if and only if the square 
\[
\begin{tikzcd}
[row sep=5ex,column sep=6ex]
gp(X) \ar[r,"gp(\alpha)"] \ar[d,"\mathrm{c}_X"'] & gp(Y) \ar[d,"\mathrm{c}_Y"] \\\
X \ar[r,"\alpha"] & Y,
\end{tikzcd} 
\]
obtained by applying the counit transformation to $\alpha$, is a pushout square.
\item
\label{item:cartesian-fibration-in-bousfield}%
If $\C$ admits pushouts which $p$ preserves and $g$ is fully faithful, then $p$ is a cocartesian fibration. 
\end{enumerate}
\end{lemma}

\begin{proof}
The first statement is immediate from the mapping space criterion for cocartesian morphisms \cite{HTT}*{Proposition 2.4.4.3}. For the second one readily checks that for $C \in \C$ and a map $p(C) \rightarrow D$ in $\Ctwo$ the edge $C \rightarrow C \cup_{gp(C)} g(D)$ is a $p$-cocartesian lift; here the pushout is formed using the counit $gp(C) \rightarrow C$ of the adjunction.
\end{proof}

Applying the previous lemma to the opposite $\infty$-category  as well we find:

\begin{corollary}
\label{corollary:splitverdbicart}%
Any split Verdier projection \(p\colon \D \to \E\) of stable $\infty$-categories is a bicartesian fibration. 
\end{corollary}

In fact, the converse also holds for exact functors $\C \rightarrow \C'$ by Remark \reftwo{remark:Verdiersplitcocartequiv} below.
Now, denote by 
\[
\Cart(p), \Cocart(p)\subseteq \Ar(\C)
\]
the full subcategories on $p$-cartesian, resp.~$p$-cocartesian morphisms. These are stable subcategories as a consequence of Lemma~\reftwo{lemma:cartesian-in-bousfield}, and the hermitian structure $\QF^{\Delta^1}$ endows $\Cart(\C)$ and $\Cocart(\C)$ with the structure of hermitian \(\infty\)-categories (we warn the reader that $\QF^{\Delta^1}(X \rightarrow Y) \simeq \QF(Y)$ is distinct from the Poincaré structures $\QF_\arr$ and $\QF_\met$). Finally, we denote by $s$ and $t\colon \Ar(\C) \to \C$ source and target functor, respectively.

\begin{lemma}
\label{lemma:cartesian-lifts-in-Cath}%
Let $p \colon (\D,\QFD) \rightarrow (\E,\QFE)$ be a split Poincaré-Verdier projection. Then the diagrams
\[
\begin{tikzcd}
[row sep=5ex]
(\Cart(p),\QFD^{\Delta^1})\ \ar[r,"t"] \ar[d,"p"]  & (\D,\QFD) \ar[d,"p"]  \\
(\E,\QFE)^{\Delta^1} \ar[r,"t"] & (\E,\QFE) 
\end{tikzcd}
\qquad
\begin{tikzcd}
[row sep=5ex]
(\Cocart(p),\QFD^{\Delta^1}) \ar[r,"s"] \ar[d,"p"] & (\D,\QFD) \ar[d,"p"] \\ 
(\E,\QFE)^{\Delta^1} \ar[r,"s"] & (\E,\QFE)
\end{tikzcd}
\]
in $\Cath$ are cartesian.
\end{lemma}

\begin{lemma}
\label{lemma:cartesian-lifts-2-in-Cath}%
Let $p \colon (\D,\QFD) \rightarrow (\E,\QFE)$ be a split Poincaré-Verdier projection. Then the square
\[
\begin{tikzcd}
[row sep=5ex,column sep=7ex]
(\Cart(p),\QFD^{\Delta^1})\times_{(\D,\QFD)^{\Delta^1}} (\D,\QFD)^{\Delta^2} \ar[r,"{(\id, d_1)}"] \ar[d,"p"] 
 & (\Cart(p),\QFD^{\Delta^1}) \times_{(\D,\QFD)}(\D,\QFD)^{\Delta^1} \ar[d,"p"]
 \\
(\E,\QFE)^{\Delta^2} \ar[r,"{(d_0, d_1)}"] 
 & (\E,\QFE)^{\Delta^1} \times_{(\E,\QFE)} (\E,\QFE)^{\Delta^1},
\end{tikzcd}
\]
where the pullback in the top left corner is formed using $d_0 \colon \Delta^1 \rightarrow \Delta^2$ and those on the right using the target functor is cartesian in $\Cath$. Similarly,  
\[
\begin{tikzcd}
[row sep=5ex,column sep=7ex]
(\Cocart(p),\QFD^{\Delta^1}) \times_{(\D,\QFD)^{\Delta^1}} (\D,\QFD)^{\Delta^2} \ar[r,"{(\id, d_1)}"] \ar[d,"p"] 
 & (\Cocart(p),\QFD^{\Delta^1})\times_{(\D,\QFD)} (\D,\QFD)^{\Delta^1} \ar[d,"p"]
 \\
(\E,\QFE)^{\Delta^2} \ar[r,"{(d_2, d_1)}"]
 & (\E,\QFE)^{\Delta^1}\times_{(\E,\QFE)} (\E,\QFE)^{\Delta^1}.
\end{tikzcd}
\]
with top left corner formed using $d_2 \colon \Delta^1 \rightarrow \Delta^2$ and right hand using the source functor, is cartesian in $\Cath$.
\end{lemma}

\begin{proof}[Proof of  Lemma~\reftwo{lemma:cartesian-lifts-in-Cath}]
One readily checks straight from the definitions and the mapping space criterion for cartesian edges \cite{HTT}*{Proposition 2.4.4.3} that the map $\Cart(p) \rightarrow \Ar(\Ctwo) \times_\Ctwo \C$ is essentially surjective and fully faithful for any cartesian fibration $p\colon \C \to \C'$. Now, apply Corollary~\reftwo{corollary:splitverdbicart} to obtain the statement at the level of underlying $\infty$-categories.

To see that this map is an equivalence $\Cath$, note first that by the discussion in \S\refone{subsection:limits}
it is enough to show that for a cartesian morphism $f\colon \Delta^1 \to \D$ the square 
\[
\begin{tikzcd}
[row sep=5ex]
\lim (\QFD\circ f\op) \ar[r] \ar[d] & \QFD(f(1)) \ar[d]\\
\lim (\QFE\circ (pf)\op) \ar[r] & \QFE(pf(1))
\end{tikzcd}
\]
is a pullback of spectra. But this is clear since the horizontal maps are equivalences, as $1$ is initial in $(\Delta^1)\op$.

Now we deal with the second square. That the underlying square of $\infty$-categories is cartesian is again easy (or indeed follows from the cartesian case applied to $p\op$). For the hermitian structure, we need to show that 
\[
\begin{tikzcd}
[row sep=5ex]
\QFD(f(1)) \ar[r] \ar[d] & \QFD(f(0)) \ar[d]\\
\QFE(pf(1)) \ar[r] & \QFE(pf(0)).
\end{tikzcd}
\]
is a pullback for every $p$-cocartesian morphism $f$. 
To see this, recall from Lemma~\reftwo{lemma:cartesian-in-bousfield} that $f(1)\simeq f(0) \cup_{lpf(0)} lpf(1)$, where $l$ is the left adjoint to $p$. Furthermore, the canonical map $\QFD\circ l\to \QFE$ is an equivalence, since $p$ is a split Poincaré-Verdier projection; see Corollary~\reftwo{corollary:split-poincare-projection-inclusion}. Thus, the square in question is equivalent to
\[
\begin{tikzcd}
\QFD(f(0)\cup_{lpf(0)} lpf(1)) \ar[r] \ar[d] & \QFD(f(0)) \ar[d]\\
\QFD(lpf(1)) \ar[r] & \QFD(lpf(0)).
\end{tikzcd}
\]
By Lemma~\refone{lemma:goodwillie} it is therefore enough to show that
\[
\Bil_\QFD(\cof(lpf), \cof(c)) \simeq 0,
\]
where $c \colon lpf(0) \rightarrow f(0)$ is the counit of the adjunction. 
For this, we compute 
\begin{align*}
\Bil_\QFD(\cof(lpf), \cof(c)) &\simeq \Hom_\D(l\cof(pf), \Dual_\QFD \cof(c)) \\
                             &\simeq \Hom_{\E}(\cof(pf), p \Dual_\QFD\cof(c)) \\
                             &\simeq \Hom_{\E}(\cof(pf), \Dual_{\QFE}\cof(pc)) = 0
\end{align*}
since $pc$ is an equivalence.
\end{proof}

For the proof of Lemma~\reftwo{lemma:cartesian-lifts-2-in-Cath}, we use the following observation; compare \cite{HTT}*{Corollary 2.4.2.5}:

\begin{observation}
\label{observation:restricted-cartesian-fibration}%
Let $p\colon \C\to \Ctwo$ be a cartesian fibration, and $\C_0\subseteq \C$ be a full subcategory that contains all $p$-cartesian morphisms whose target lies in $\C_0$. Then the restricted functor $p\colon \C_0\to \Ctwo$ is also a cartesian fibration. 
\end{observation}

\begin{proof}[Proof of Lemma~\reftwo{lemma:cartesian-lifts-2-in-Cath}]
We again start by showing that the upper square is a pullback of $\infty$-categories. We first claim that both vertical maps are cartesian fibrations. By \cite{HTT}*{3.1.2.1}, for a cartesian fibration $p\colon \C \to \C'$ and $K$ in $\Cat$, the functor $p_* \colon \Fun(K,\C) \rightarrow \Fun(K,\Ctwo)$ is again a cartesian fibration, with cartesian edges detected pointwise. Applying this with $K= \Delta^2$ and $\Lambda_2^2$, the claim easily follows from Observation~\reftwo{observation:restricted-cartesian-fibration} and the cancellability of cartesian edges \cite{HTT}*{Proposition 2.4.1.7}. The pointwise nature of cartesian edges also implies that the top horizontal map preserves cartesian edges, so to check that the underlying diagram is cartesian in $\Cat$ it suffices to check that the induced map on vertical fibres are equivalences by \cite{HTT}*{Corollary 2.4.4.4}.
But here again, one checks that the induced functors are fully faithful and essentially surjective from the mapping space criterion for cartesian edges \cite{HTT}*{Proposition 2.4.4.3} together with the description of spaces of natural transformation as iterated pullbacks arising from \cite{GHNfree}*{Proposition 5.1}.

This concludes the proof that the underlying diagram of the top square is a pullback in $\Cat$, and the argument for the bottom one is entirely analogous. To make the first square a pullback in $\Cath$, we need to show that for each $f\colon \Delta^2\to \D$ with $f(1) \to f(2)$ $p$-cartesian, the square of spectra
\[
\begin{tikzcd}
 \lim_{(\Delta^2)\op} \QFD\circ f\op \ar[r] \ar[d] & \lim_{(\Lambda_2^2)\op} \QFD\circ f\op \ar[d]\\
 \lim_{(\Delta^2)\op} \QFE\circ pf\op \ar[r] & \lim_{(\Lambda_2^2)\op} \QFE\circ pf\op
\end{tikzcd}
\]
is a pullback. But this is clear since $2$ is terminal in both $\Delta^2$ and $\Lambda_2^2$, so the inclusion $(\Lambda_2^2)\op \subset (\Delta^2)\op$ is final and the horizontal maps are equivalences. 

To see that the second square is a pullback in $\Cath$, we have to show that for each $f \colon \Delta^2 \to \D$ with $f(0) \to f(1)$ $p$-cocartesian, the following square is a pull-back:
\[
\begin{tikzcd}
 \lim_{(\Delta^2)\op} \QFD\circ f\op \ar[r] \ar[d] & \lim_{(\Lambda_0^2)\op} \QFD\circ f\op \ar[d]\\
 \lim_{(\Delta^2)\op} \QFE\circ pf\op \ar[r] & \lim_{(\Lambda_0^2)\op} \QFE\circ pf\op
\end{tikzcd}
\]
Since $2$ is terminal in $\Delta^2$, this reads
\[
\begin{tikzcd}
 \QFD(f(2)) \ar[r] \ar[d] & \QFD(f(1))\times_{\QFD(f(0))} \QFD(f(2)) \ar[d]\\
 \QFE(pf(2)) \ar[r]  & \QFE(pf(1))\times_{\QFE(pf(0))} \QFE(pf(2)).
\end{tikzcd}
\]
This square is indeed cartesian since $\QFD(f(1)) \simeq \QFE(pf(1)) \times_{\QFE(pf(0))}\QFD(f(0))$ as a consequence of Lemma~\reftwo{lemma:cartesian-lifts-in-Cath}; in fact this statement is also an explicit step in the proof of \reftwo{lemma:cartesian-lifts-in-Cath}.%
\end{proof}

Now let $\C_p\subseteq \Q_1(\D)$ denote the full subcategory on objects of the form $D \leftarrow W \to D'$
where the left arrow is $p$-cartesian, and the right arrow is $p$-cocartesian. This is a stable subcategory which inherits a hermitian structure from $\Q_1(\D)$.

\begin{lemma}
$\C_p\subseteq \Q_1(\D)$ is closed under the duality $\Dual_{\QFD_1}$. 
\end{lemma}

Therefore, $(\C_p, \QFD_1)$ is a Poincaré \(\infty\)-category and the inclusion functor $\C_p\to \Q_1(\D)$ canonically refines to a Poincaré functor.

\begin{proof}
Let $D\leftarrow W \rightarrow D'$ be an object of $\C_p$. The dual arrow is obtained by first completing the diagram to a pushout square; then applying $\Dual_\QFD$ termwise, and deleting the value at the terminal object of the square, see Proposition~\refone{proposition:basic-properties-herm-diagrams}. The claim now follows from the fact that $p$-(co-)cartesian morphisms are stable under (co-)base change, and that the dualities interchange $p$-cartesian with $p$-cocartesian morphisms since the diagram
\[
\begin{tikzcd}
\D\op \ar[r,"\Dual_\QFD"] \ar[d,"p\op"] & \D \ar[d,"p"] \\
\E\op \ar[r,"\Dual_\QFE"] & \E
\end{tikzcd}
\]
commutes as $p$ is Poincaré. 
\end{proof}

We are now ready to state the main technical results of this section, namely that $\Q(p) \colon \Q(\D,\QFD) \rightarrow \Q(\E,\QFE)$ 
behaves like a cocartesian fibration of Segal objects in $\Cath$, with cocartesian lifts given by $(\C_p,\QFD_1) \subseteq \Q_1(\D,\QFD)$. Since the $\Q$-construction is invariant under taking the opposite simplicial object, it follows that it also behaves like a cartesian fibration, see Example~\reftwo{example:CobHyp-equal-Span}.

\begin{lemma}
\label{lemma:cocartesian-lifts-1-in-CobF}%
\label{lemma:cocartesian-lifts-2-in-CobF}%
The squares
\[
\begin{tikzcd}
 (\C_p,\QFD_1) \ar[r,"d_1"] \ar[d,"p"] & (\D,\QFD) \ar[d,"p"] & (\C_p,\QFD_1) \times_{\Q_1(\D,\QFD)} \Q_2(\D,\QFD) \ar[r,"{(\id,d_1)}"] \ar[d,"p"]  & (\C_p,\QFD_1) \times_{(\D,\QFD)} \Q_1(\D,\QFD) \ar[d,"p"] \\
 \Q_1(\E,\QFE) \ar[r,"d_1"] & (\E,\QFE) & \Q_2(\E,\QFE) \ar[r,"{(d_2, d_1)}"] & \Q_1(\E,\QFE) \times_{(\E,\QFE)} \Q_1(\E,\QFE)
\end{tikzcd}
\]
are split Poincar\'e-Verdier squares. In the right hand square, the upper left pullback is formed using $d_2 \colon \Q_2(\D,\QFD) \rightarrow \Q_1(\D,\QFD)$ and the right hand ones using $d_1$. In particular, both diagrams are cartesian in $\Catp$.
\end{lemma}
\begin{proof}
We begin with the left hand square and factor it as follows.
\[
\begin{tikzcd}
 (\C_p,\QFD_1) \ar[r] \ar[d,"p"]
 & (\Cart(p),\QFD^{\Delta^1}) \ar[r,"t"] \ar[d,"p"]
 & (\D,\QFD) \ar[d,"p"]
 \\
 \Q_1(\E,\QFE) \ar[r] 
 & (\E,\QFE)^{\Delta^1} \ar[r,"t"]
 & (\E,\QFE)
\end{tikzcd}
\]
Here, the left horizontal maps are induced by including $\Delta^1$ into $\Twar\Delta^1$ as the  morphism $(0 \leq 1)\to (0 \leq 0)$. The right square is a pullback by Lemma~\reftwo{lemma:cartesian-lifts-in-Cath}. Now
\[
\Q_1(\D,\QFD) \simeq (\D,\QFD)^{\Lambda_0^2} \simeq (\D,\QFD)^{\Delta^1} \times_{(\D,\QFD)} (\D,\QFD)^{\Delta^1}
\]
using the source and target arrows for the pullback, and this equivalence restricts to an equivalence 
\[
\C_p\simeq \Cart(p)\times_\D \Cocart(p)
\]
by construction. So, the left square is obtained by pullback from the right hand square of Lemma~\reftwo{lemma:cartesian-lifts-in-Cath} %
and therefore cartesian as well (in $\Cath$, and hence in $\Catp$). 
Since $p$ is a split Poincaré-Verdier projection by assumption this implies the claim by Corollary~\reftwo{corollary:PV-proj-pullback}.

Next we treat the right hand square.
The \(\infty\)-category in its upper left corner is equivalent (as a hermitian \(\infty\)-category) to the full subcategory of $\Q_2(\D,\QFD)$ on those diagrams $F\colon \Twar\Delta^2\to \D$,
\[
\begin{tikzcd}[row sep=3ex,column sep=2ex]
 & & F(0 \leq 2) \ar[ld] \ar[rd,"{(\mathrm{III})}"']\\
  & F(0\leq 1) \ar[ld,"{(\mathrm{I})}"] \ar[rd,"{(\mathrm{II})}"'] & & F(1\leq 2) \ar[ld] \ar[rd]\\\
 F(0\leq 0) & & F(1\leq 1) & & F(2 \leq 2)
\end{tikzcd}
\]
such that (i) the map labelled by $(\mathrm{I})$ is $p$-cartesian, (ii) the map labelled by $(\mathrm{II})$ is $p$-cocartesian, and (iii) the middle square is cartesian. In view of Lemma~\reftwo{lemma:cartesian-in-bousfield}, one easily checks by pasting squares that condition (iii) is equivalent to the following two conditions: (iii') the map labelled by $(\mathrm{III})$ is $p$-cocartesian, and (iii'') the image of the middle square in $\E$ is cocartesian. In other words, if we denote by $(\D,\QFD)_p^{\Twar(\Delta^2)} \subseteq (\D,\QFD)^{\Twar(\Delta^2)}$ the full subcategory on diagrams satisfying (i), (ii), and (iii'), then the diagram
\[
\begin{tikzcd}
(\C_p,\QFD_1) \times_{\Q_1(\D,\QFD)} \Q_2(\D,\QFD) \ar[d,"p"] \ar[rr] & & (\D,\QFD)_p^{\Twar(\Delta^2)} \ar[d,"p"] \\
\Q_2(\E,\QFE) \ar[rr] && (\E,\QFE)^{\Twar(\Delta^2)}
\end{tikzcd}
\]
is a pullback in $\Cath$, since it is one in $\Cat$ and the hermitian structures on the left are the restrictions of those on the right. %
Now consider the following filtration 
\[
I_0 \to I_1 \to \dots I_4 = \Twar(\Delta^2)
\]
through (non-full) subposets, starting with 
\[
I_0 = d_2(\Twar\Delta^{1}) \cup d_1(\Twar\Delta^{1}).
\]
The remaining $I_i$ are obtained by adding relations in the order indicated in the following picture, in which circles indicate $2$-cells: 
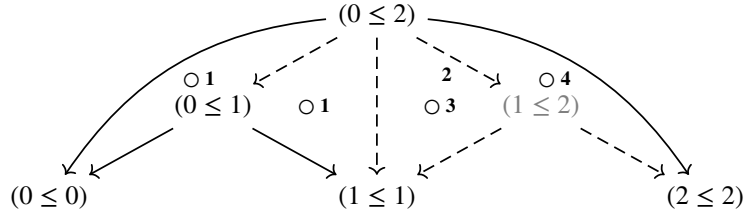
\[
\begin{tikzcd}
 & & (0\leq 2) \ar[ld,dashed] \ar[rd,dashed,"{\mathbf{2}}" description] \ar[lldd,bend right=25,"{\bigcircle\;\mathbf{1}}"] \ar[rrdd,bend left=25,"{\bigcircle\;\mathbf{4}}"'] \ar[dd,dashed] \\
 & (0 \leq 1) \ar[ld]  \ar[rd] \ar[r,phantom,"{\scriptstyle \bigcircle\;\mathbf{1}}"] & {} & {\color{gray} (1 \leq 2)} \ar[ld,dashed] \ar[rd,dashed] \ar[l,phantom,"{\scriptstyle \bigcircle\;\mathbf{3}}"] \\
 (0\leq 0) & & (1\leq 1) & & (2\leq 2) 
\end{tikzcd}
\]
Now, one readily checks that each $I_i\to I_{i+1}$ is obtained as a pushout of an outer horn inclusion (namely using $\Lambda_2^2, \Lambda_0^1$ and then  $\Lambda_0^2$ twice) in $\Cat$: This either follows from a simple direct argument by writing the posets involved as iterated pushouts of simplices, or from the corresponding statement at the level of simplicial sets using that homotopy pushouts in the Joyal model structure model pushouts in $\Cat$.

For $i\in \{0, \dots, 4\}$, let $(\D,\QFD)^{I_i}_p \subseteq (\D,\QFD)^{I_i}$ denote the full subcategory on functors that satisfy whichever of condition (i), (ii), and (iii') apply. Then for $i=0$ the map $(\D,\QFD)_p^{I_i} \rightarrow (\E,\QFE)^{I_i}$ induced by $p$ is equivalent to that in the right hand column of the statement of the Lemma, and for $i=4$ it is the right hand map in the square above.%

We then claim that the diagram 
\[
\begin{tikzcd}
(\D,\QFD)_p^{I_i} \ar[rr] \ar[d,"p"] & & (\D,\QFD)_p^{I_{i-1}} \ar[d,"p"] \\
(\E,\QFE)^{I_i} \ar[rr] & & (\E,\QFE)^{I_{i-1}},
\end{tikzcd}
\]
with horizontal maps given by restriction, is a pullback in $\Cath$ for $i = 1,2,3,4$. This establishes the lemma by pasting pullbacks.

Indeed, $I_2$ is obtained from $I_1$ by filling the $1$-horn $\Lambda^1_0\subset \Delta^1$ with a cocartesian edge, so that the restriction map $(\D,\QFD)^{I_2}\to (\D,\QFD)^{I_1}$ is pulled back from the restriction map $s\colon (\D,\QFD)^{\Delta^1} \to(\D,\QFD)$. It follows that the diagram in question is obtained from the second diagram of Lemma~\reftwo{lemma:cartesian-lifts-in-Cath} by base changes, and therefore is a pullback. 
Similarly, we see that the diagrams for $i=1,3,4$ are obtained by base changes from the diagrams of Lemma~\reftwo{lemma:cartesian-lifts-2-in-Cath} and therefore pullbacks. %

We are left to show that the map
\[
(\C_p,\QFD_1) \times_{(\D,\QFD)} \Q_1(\D,\QFD)  \longrightarrow \Q_1(\E,\QFE) \times_{(\E,\QFE)} \Q_1(\E,\QFE),
\]
is a split Poincaré-Verdier projection. But Lemma~\reftwo{lemma:cocartesian-lifts-1-in-CobF} identifies this map as a base change of the map $\Q_1(p) \colon \Q_1(\D,\QFD) \rightarrow \Q_1(\E,\QFE)$, which is a Poincaré-Verdier projection by Proposition~\reftwo{proposition:(co)tensor-Verdier}. The claim thus follows from Corollary~\reftwo{corollary:PV-proj-pullback} as well.
\end{proof}

\begin{proof}[Proof of Theorem~\reftwo{theorem:fulladditivity}]
Applying $\F$ to the squares of Lemma~\reftwo{lemma:cocartesian-lifts-1-in-CobF}, and using additivity, we deduce that the following squares are cartesian
\[
\begin{tikzcd}
\F(\C_p,\QFD_1) \ar[r,"d_1"] \ar[d,"p"] & \F(\D,\QFD) \ar[d,"p"] \\ 
\F(\Q_1(\E,\QFE)) \ar[r,"d_1"] & \F(\E,\QFE)
\end{tikzcd}
\begin{tikzcd}
\F(\C_p,\QFD_1)\times_{\F(\Q_1(\D,\QFD))} \F(\Q_2(\D,\QFD)) \ar[r,"{(\id, d_1)}"] \ar[d,"p"] & \F(\C_p,\QFD_1)\times_{\F(\D,\QFD)} \F(\Q_1(\D,\QFD)) \ar[d,"p"] \\
\F(\Q_2(\E,\QFE)) \ar[r,"{(d_2, d_1)}"] & \F(\Q_1(\E,\QFE)) \times_{\F(\E,\QFE)} \F(\Q_1(\E,\QFE));
\end{tikzcd}
\]
here, the pullbacks in the right hand square are formed using $d_2$ on the left and $d_1$ on the right. Now, the right hand square tells us that the image of $\pi_0\F(\C_p,\QFD_1) \rightarrow \pi_0\F(\Q_1(\D,\QFD))$ consists of $\F\Q(p)$-cocartesian arrows, whence the left hand square provides sufficiently many $\F\Q(p)$-cocartesian lifts to make $\F\Q(p) \colon \F\Q(\D,\QFD) \rightarrow \F\Q(\E,\QFE)$ into a cocartesian fibration of Segal spaces and therefore $\Cob^\F(p)$ into a cocartesian fibration of $\infty$-categories. 
Since $\Q(\D,\QFD)$ is naturally identified with $\Q(\D,\QFD)\op$ through the canonical equivalence $\Twar(\Delta^n) \cong \Twar((\Delta^n)\op)$, we conclude that both $\F\Q(p)$ and $\Cob^\F(p)$ are also cartesian fibrations.
\end{proof}

\subsection{Additivity in $\K$-Theory}
\label{subsection:Kadd}%

The arguments presented in the previous section work verbatim upon dropping hermitian structures and working with additive functors $\Catx \rightarrow \Sps$. In the present section, we briefly record the statements that are obtained this way.  
Let us first formally set terminology obviously analogous to that of Definition~\reftwo{definition:additive}. 
\begin{definition}
\label{definition:additive-on-catx}%
Let $\E$ be an $\infty$-category with finite limits and $\F\colon \Catx \to \E$ a reduced functor. We say that $\F$ is \defi{additive}, \defi{Verdier-localising} or \defi{Karoubi-localising} if it sends split Verdier squares, arbitrary Verdier squares or Karoubi squares to cartesian squares, respectively.
\end{definition}

Part of the following result also appears in \cite{BarwickRognesQexact}, though in incommensurable generality.

\begin{proposition}
For a stable \(\infty\)-category $\C$, the simplicial $\infty$-category $\Q(\C)$ is a complete Segal object in $\Cat$, whose boundary maps are split Verdier projections.
For an additive functor $\F \colon \Catx \rightarrow \Sps$, the simplicial space $\F\Q(\C)$ is a Segal space, which is complete if $\F$ preserves pullbacks.
\end{proposition}

\begin{proof}
The first two statements are obtained during the proofs of Lemmas~\reftwo{lemma:segal} (see also \cite{HHLN-two-var-fil}*{Lemma~2.17}) and \reftwo{lemma:boundarysplit}. The latter two statements are proven just as Proposition~\reftwo{proposition:complete}.
\end{proof}

In particular, we can set
$\Span^\F(\C) := \asscat(\F\Q(\C))$,
which inherits a symmetric monoidal structure since $\Catx$ is semi-additive and $\C \mapsto \Span^\F(\C)$ preserves finite products. The proof of Corollary~\reftwo{corollary:components-of-cob-II} gives the statement that the diagram of commutative monoids
\[
\begin{tikzcd}
\pi_0\F(\Ar(\C)) \ar[r,"t"] \ar[d] & \pi_0\F(\C) \ar[d] \\
0 \ar[r] & \pi_0 |\Span^\F(\C)|
\end{tikzcd}
\]
is cocartesian. In the present situation, the top horizontal map is, however, split surjective, hence we obtain: 

\begin{proposition}
\label{proposition:spanconn}%
The $\infty$-category $\Span^\F(\C)$ is connected for any stable $\C$ and additive $\F \colon \Catx \rightarrow \Sps$.
\end{proposition}

In particular, $|\Span^\F(\C)|$ is always an $\Einf$-group. Our notion of additive functor is geared to permit the following strong version of Waldhausen's additivity theorem:

\begin{theorem}[Additivity]
\label{theorem:addinkthm}%
If $\F\colon \Catx \rightarrow \Sps$ is additive, then so is $|\Span^\F-| \simeq |\F\Q-|$.
\end{theorem}

Just as the hermitian additivity theorem~\reftwo{theorem:additivity}, this theorem is deduced from the following statement and the fact that bicartesian fibrations are realisation fibrations.

\begin{theorem}
\label{theorem:classadd}%
Let $\F\colon \Catx \rightarrow \Sps$ be additive and $p \colon \D \rightarrow \E$ a split Verdier projection, then 
\[
\F\Q(p) \colon \F\Q(\D) \to \F\Q(\E)
\]
is a bicartesian fibration of Segal spaces and thus a realisation fibration.
\end{theorem}
\begin{proof}
The proof of Theorem~\reftwo{theorem:fulladditivity} in \S\reftwo{subsection:fibcob}, in particular, verifies Theorem~\reftwo{theorem:classadd} upon dropping all mention of Poincaré structures (which in fact made up the bulk of the work). 
\end{proof}

Waldhausen's additivity theorem now follows, by inserting the analogue of the metabolic sequence, i.e.\ the split Verdier sequence
$
\C \to \Ar(\C) \rightarrow \C
$
made up of the inclusion $X \mapsto (X \rightarrow 0)$ and the target projection $t$, into the additivity theorem (whence our terminology), and noting that either adjoint of $t$ give rise to splittings of the resulting fibre sequence
\[
\F(\C) \to \F(\Ar(\C)) \xrightarrow{t} \F(\C).
\]
This runs contrary to the situation of the metabolic sequence, where the adjoints of $\met \colon \Met(\C,\QF) \rightarrow (\C,\QF)$ are not compatible with the Poincaré structures. The splitting Lemma~\reftwo{lemma:retract} in $\Grp_\Einf(\Sps)$ then gives the equivalence
$
\F(\C)^2 \simeq \F(\Ar(\C)).
$
Applying this to $\Kspace = \Omega|\Span-|$, which is additive by Theorem~\reftwo{theorem:addinkthm} above, we find $\kk(\C)^2 \simeq \kk(\Ar(\C))$ as desired. 

In summary, the fibre sequence of cobordism $\infty$-categories induced by the metabolic fibre sequence is not just an algebraic analogue of Genauer's fibre sequence regarding geometric cobordism categories, but also of Waldhausen's additivity, the connection between which was first realised by the ninth author in \cite{Steimleadd}. 
The idea to systematically relate cobordism categories of manifolds to algebraic K-theory by viewing spans as a formal type of cobordism was originally put forward in joint work of his with Raptis, see \cite{RaptisSteimle-SmoothDWW, RaptisSteimle-CobcatK,RaptisSteimle-TopDWW}.

\begin{remarks}
\label{remark:leftsplitadd}%
\begin{enumerate}
\item We repeat the caveat that it is not generally true, that $|\F\Q-|$ is Verdier-localising whenever $\F \colon \Catx \rightarrow \Sps$ is, a counterexample being $\Kspace \circ (-)^\idem$. 
\item
\label{item:left-additivity}%
In the set-up of stable $\infty$-categories, any group-like additive functor $\F$ in fact takes all right split and all left split Verdier-sequences (i.e.\ those where the projection admits only one adjoint) to fibre sequences by unpacking an old argument of Waldhausen's, see e.g.\ \cite{HLS}*{Section 2}. An hermitian analogue of this fact is the isotropic decomposition principle in the next section, which in fact implies its non-hermitian counterpart,
see Example~\reftwo{example:left-split-isotropic-decomposition} below.
\item Let us also mention that the proof of Theorem~\reftwo{theorem:suspension} below also carries over without change to the setting of stable $\infty$-categories. As a consequence one obtains:
\begin{subenumerate}
\item
$\Kspace \simeq \Omega|\Span(-)| \colon \Catx \rightarrow \Sps$ is the initial group-like additive functor under $\Core$, compare Corollary~\reftwo{corollary:universal},
\item
the iterated $\Q$-construction defines a positive $\Omega$-spectrum $\mathbb{S}\mathrm{pan}(\C)$, with $\Omega^\infty \mathbb{S}\mathrm{pan}(\C) \simeq \Kspace(\C)$, compare Proposition~\reftwo{proposition:positive-om}, 
\item
whose spectrification $\K(\C)$ gives the initial additive functor $\Catx \rightarrow \Spa$ under $\SS[\Core]$, and in fact that $\K$ is the suspension spectrum of $\Core$ in the stabilisation of $\infty$-category of additive functors $\Catx \rightarrow \Sps$, compare Corollary~\reftwo{corollary:Cob=Susp} and Remark~\reftwo{remark:K-theorysuspension}.
\end{subenumerate}
This gives a simple and uniform approach to these fundamental results, which to the best of our knowledge have not been treated together in the literature.
\end{enumerate}
\end{remarks}

As we will have to make use of this result in the next section, let us also record the computation of $\pi_0\K(\C) = \pi_1|\Span(\C)|$ in the generality of an arbitrary additive $\F \colon \Catx \rightarrow \Sps$: The natural equivalence $\Hom_{\Span^\F(\C)}(0,0) \simeq \F(\C)$ provides maps
\[
\begin{tikzcd}
\pi_0\F(\C^2) \ar[d] & \ar[l,"{(\source,\cof)}"'] \pi_0\F(\Ar\C) \ar[r,"\target"] \ar[d] & \pi_0\F(\C) \ar[d]\\
\pi_1 |\Span^\F(\C^2)| & \ar[l] \pi_1|\Span^\F(\Ar\C)| \ar[r] & \pi_1|\Span^\F(\C)|,
\end{tikzcd}
\]
where $\source,\target$ and $\cof$ take the source, target, and cofibre of a morphism. The additivity theorem implies that the lower left horizontal map is an isomorphism. Inverting it produces a commutative square
\[
\begin{tikzcd}
\pi_0\F(\Ar\C) \ar[r,"\target"] \ar[d,"{(\source,\cof)}"'] & \pi_0\F(\C) \ar[d]\\
\pi_0\F(\C)^2 \ar[r] & \pi_1|\Span^\F(\C)|
\end{tikzcd}
\]
of commutative monoids natural in both $\C$ and $\F$.

\begin{proposition}
\label{proposition:pi0span}%
This square is cocartesian for every stable $\C$ and every additive $\F \colon \Catx \rightarrow \Sps$.
\end{proposition}

In particular, for $\F = \Core$ we recover the standard fact that $\K_0(\C)$ is given by $\pi_0\Core(\C)$ modulo extensions. While a proof internal to the $\Q$-construction is certainly possible, the quickest route is through the subdivision equivalence $|\F\Q(\C)| \simeq |\F\mathrm S(\C)|$ between Quillen's $\Q$- and Segal's $\rS$-constructions (the latter famously developed by Waldhausen). Let us briefly review this equivalence; it also features prominently in Appendix~\reftwo{appendix:AppIIB}. The $\rS$-construction can be viewed as functor $\rS\colon\Catx \to \sCatx$, with $\rS_n(\C) \subseteq \Fun(\Ar(\Delta^n),\C)$ given by the full subcategory on functors $F$ satisfying $F(i \leq i) = 0$ for all $i \in [n]$ and that for all $i \leq j \leq k\leq l$, the square 
\[
\begin{tikzcd}
	F(i\leq k) \ar[r] \ar[d] & F(i\leq l) \ar[d] \\
	F(j\leq k) \ar[r] & F(j \leq l)
\end{tikzcd}
\]
is cartesian. Now recall, that simplicial objects can be subdivided edgewise by precomposition with the functor $\Delta\op \to \Delta\op$ sending $[n]$ to $[n] \ast [n]\op$. In the case of simplicial spaces it takes the standard simplices to contractible simplicial spaces and preserves colimits, so preserves realisations. Writing $\rS^e(\C)$ for the edgewise subdivision of $\mathrm{S}(\C)$ we thus find $|\F\rS(\C)| \simeq |\F\rS^e(\C)|$ for any additive $\F \colon \Catx \rightarrow \Sps$. But it is easily checked that $\rS^e(\C) \simeq \Q(\C)$ via the map $\Twar(\Delta^n) \to \Ar(\Delta^n \ast (\Delta^n)\op)$, natural in $[n]$, taking $(i\leq j)$ to $(i_0 \leq j_1)$ where the subscripts indicate the join factor.%

\begin{proof}[Proof of Proposition~\reftwo{proposition:pi0span}]
In $\mathrm S(\C)$ the 0-,1- and 2-simplices are equivalently given by $*$, $\C$ and $\mathrm{Cof}(\C)$, respectively, where $\mathrm{Cof}(\C)$ denotes the $\infty$-category of cofibre sequences in $\C$. This is equivalent to $\Ar(\C)$ and under this identification the boundary maps of $\mathrm S(\C)$ are given by source, target and cofibre. Thus we find $\pi_1|\F\rS(\C)|$ given by $\pi_0\F(\C)$ modulo the relation $s(f) + \cof(f) = t(f)$ for every $f \in \pi_0\F(\Ar(\C))$, which is exactly the pushout above.
\end{proof}

Finally, let us mention that $\Kspace \colon \Catx \rightarrow \Sps$ is in fact Verdier localising and not just additive. Again this essentially goes back to work of Waldhausen, and (idempotent completions issues aside) is first explicit in the work of Blumberg, Gepner and Tabuada. A direct proof in the present setting can be found in \cite{HLS}*{Section 6}.

\section{Structure theory for additive functors}
\label{section:structure-theory}%

The objective of this section is to derive the fundamental theorems of Grothendieck-Witt theory from the additivity theorem. We will, however, do so in the generality of arbitrary additive functors $\Catp \rightarrow \Sps$. Even when only interested in Grothendieck-Witt spectra this additional layer of generality is useful, for example it enters our proof of the universal property of $\GW \colon \Catp \rightarrow \Spa$. The reader is encouraged to keep the two fundamental examples 
\[
\Poinc \quad \text{and} \quad  |\Cob(-)|\colon \Catp\to \Sps
\]
in mind throughout. In \S\reftwo{section:GW} below, we will specialise the results of this section to define Grothendieck-Witt theory and conclude the main theorems of this paper.

We begin by introducing the notion of a cobordism between Poincaré functors, and use this to establish some fundamental results for group-like additive functors. Chief among these is the agreeance of their values on hyperbolic and metabolic categories. In the case of $|\Cob(-)|$, we already proved this claim in Proposition~\reftwo{proposition:CobMet=CobHyp} by explicit identification of both sides. Using the general statement as a base case, we develop a general theory of isotropic decompositions of Poincaré \(\infty\)-categories, which allows for the computations of the values of a group-like additive functor $\F$ applied to many Poincaré \(\infty\)-categories $(\C,\QF)$ of interest, e.g.\ $\Q_n(\C,\QF)$ for all $n$, in terms of hyperbolic pieces and parts that are often simpler than the original $(\C,\QF)$.

We then use this machinery to establish precise relationships between additive functors taking values in the $\infty$-categories of $\Einf$-monoids, $\Einf$-groups, and spectra, in particular constructing left adjoints to the evident forgetful functors. The adjoint passing from $\Einf$-monoid- to $\Einf$-group-valued functors, the \defi{group-completion}, is given by $\F \mapsto \Omega|\Cob^\F(-)| = \Omega|\F\Q(-\qshift{1})|$, using the $\F$-based cobordism $\infty$-category from section \S\reftwo{section:cobcats}, and the adjoint from $\Einf$-group-valued to spectrum-valued functors, the \defi{spectrification}, is given by iterating the $\Q$-construction on $\F$. This generalises the work of Blumberg-Gepner-Tabuada on the universality of algebraic $\K$-theory \cite{BGT}. Many of our constructions also have geometric precursors in the work of Bökstedt-Madsen on the connection between iterated cobordism categories and algebraic K-theory \cite{BokstedtMadsen}. We will expand on these analogies in \S\reftwo{section:GW}.

We then turn to a more detailed analysis of spectrum-valued additive functors. To this end we introduce the notion of a bordism-invariant functor (i.e.\ one that vanishes on metabolic categories), the principal example being $\L \colon \Catp \rightarrow \Spa$, the $\L$-theory functor of Ranicki and Lurie. We show that the inclusion of bordism invariant  functors into all additive functors also admits a left adjoint $\bord$. It will then follow for rather formal reasons that there always is a natural cartesian square
\[
\begin{tikzcd}
\F(\C,\QF) \ar[r] \ar[d] & \F^\bord(\C,\QF) \ar[d] \\
\F(\Hyp(\C))^\hC \ar[r] & \F(\Hyp(\C))^\tC
\end{tikzcd}
\]
of spectra, which can in principle be used to compute $\F$ from its \emph{hyperbolisation} $\F^\hyp = \F \circ \Hyp$ and its \emph{bordification} $\F^\bord$, each of which may be easier to understand than $\F$. We also provide two direct formulas for $\F^\bord$, which again have precursors in manifold theory. We use these in \S\reftwo{section:GW} to identify the bordification of Grothendieck-Witt theory with L-theory, completing the proof of our main theorem.

\subsection{Cobordisms of Poincaré functors}
\label{subsection:isotropic}%

In the previous section, we introduced the concept of a cobordism in a Poincaré $\infty$-category. Using it, we make the following definition.

\begin{definition}
\label{definition:bordism-functors}%
Let $(\C,\QF)$ and $(\D,\QFD)$ be two Poincaré $\infty$-categories and let $f,g\colon (\C,\QF) \rightarrow (\D,\QFD)$ be two Poincaré functors. By a \defi{cobordism} from $f$ to $g$ we shall mean a cobordism in the Poincaré $\infty$-category $\Funx((\C,\QF),(\D,\QFD))$ between the Poincaré objects corresponding to $f$ and $g$. 
\end{definition}

We note that the data of such a cobordism can equivalently be encoded by a Poincaré functor $\phi\colon (\C,\QF) \rightarrow \Q_1(\D,\QFD)$ such that $d_0\phi=f$ and $d_1\phi=g$, in particular $d_0, d_1 \colon \Q_1(\C,\QF) \rightarrow (\C,\QF)$ are cobordant.

\begin{example}
\label{example:idmetbordnull}%
As an algebraic analogue of the fact that bordism groups of manifolds with boundaries vanish, we note that the identity of $\Met(\C,\QF)$ is cobordant to the null functor for any Poincar\'e $\infty$-category $(\C,\QF)$: In Remark~\refone{remark:comonad} we observed that the tautological Poincar\'e refinement of the functor $\mathrm{cm} \colon \Met(\C,\QF) \rightarrow \Met\Met(\C,\QF)$ given by sending $W \rightarrow X$ to the right square of
\[
\begin{tikzcd}0 \ar[d] &  \ar[l] W \ar[r] \ar[d] & W \ar[d] \\ 0 & \ar[l] W \ar[r] & X\end{tikzcd}
\]
is a common section of $\Met(\met_{(\C,\QF)})$ and $\met_{\Met(\C,\QF)}$. The former implies that the composition with the inclusion $\Met\Met(\C,\QF) \subseteq \Q_1(\Met(\C,\QF))$, in total taking $W \rightarrow X$ to the entire diagram above, gives a cobordism as desired; we will use that $\mathrm{cm}$ is also a section of $\met(\Met)$ in Proposition \reftwo{proposition:boundarymapGB} below. %
\end{example}
Our next goal is to describe the behaviour of group-like additive functors, see Definition~\reftwo{definition:grouplike}, under such cobordisms. 
We start by analysing the universal case of a Poincaré cobordism between functors with target $(\C,\QF)$. It is given by the two Poincaré functors $d_0,d_1\colon \Q_1(\C,\QF) \rightarrow \Q_0(\C,\QF)=(\C,\QF)$, which are equipped with a canonical cobordism between them. 

To this end, consider %
the functor 
\[
i\colon\C \lrar \Q_1(\C), \quad X \longmapsto \big[0 \leftarrow X \rightarrow X\big]
\]
and its right adjoint 
\[
p\colon\Q_1(\C) \lrar \C, \quad \big[X \leftarrow W \rightarrow Y\big] \longmapsto \fib(W \rightarrow X).
\]
Note that the unit transformation $\id \Rightarrow pi$ is an equivalence, so $i$ is fully faithful. By the universal property of the hyperbolic construction, Corollary~\refone{corollary:hyp-is-adjoint}, we obtain a pair of Poincaré functors
\[
\begin{tikzcd}
\Hyp(\C) \ar[r,"i_{\hyp}"] & \Q_1(\C,\QF) \ar[r,"p^{\hyp}"] & \Hyp(\C) 
\end{tikzcd}
\]
which is a retract diagram in $\Catp$. We also note that $i_{\hyp}\colon\Hyp(\C) \rightarrow \Q_1(\C,\QF)$ factors through $\Met(\C,\QF) \subseteq \Q_1(\C,\QF)$; the corresponding restriction of $i_{\hyp}$ agrees with $\ilag\colon \Hyp(\C) \rightarrow \Met(\C,\QF)$; see the recollection section for a review of notation. Similarly, the restriction of $p^{\hyp}$ to $\Met(\C,\QF) \subseteq \Q_1(\C,\QF)$ is $\lag\colon \Met(\C,\QF) \rightarrow \Hyp(\C)$. 
In particular, we obtain the diamond shaped diagram
\[
\begin{tikzcd}
& (\C,\QF) \ar[d,"\cyl"] \ar[dr,"\id"] & \\
\Met(\C,\QF) \ar[r] \ar[dr,"\lag"'] & \Q_1(\C,\QF) \ar[d,"p^{\hyp}"] \ar[r,"d_1"] & (\C,\QF) \\ 
& \Hyp(\C) 
\end{tikzcd}
\]
with $\cyl$ the inclusion of constant functors and $p^\hyp$ split, as a Poincaré functor, by $i_\hyp$.

\begin{lemma}
\label{lemma:met=hyp}%
For $(\C,\QF)$ a Poincaré \(\infty\)-category, both the horizontal and vertical sequences of the diamond above are split Poincaré-Verdier sequences.
\end{lemma}

\begin{proof}
For the horizontal sequence, this is immediate from Lemma~\reftwo{lemma:boundarysplit}. For the vertical sequence, we shall check that $p$ satisfies the assumptions of Lemma~\reftwo{lemma:PV-hyperbolic} to conclude that $p^\hyp$ is a split Poincaré-Verdier projection; the kernel of $p^\hyp$ is evidently given by the diagrams $\Twar(\Delta^1) \rightarrow \core\C$, and since $|\Twar(\Delta^1)|$ is contractible these are exactly the constant diagrams, which embed $\C$ fully faithfully into $\Q_1(\C)$ and evidently $\cyl^* \QF_1 \simeq \QF$. 

We already recorded above that $p$ admits a fully faithful left adjoint $i$ taking $X$ to $[0 \leftarrow X \rightarrow X]$, and
$
\Qoppa_1(0 \leftarrow X \rightarrow X) \simeq \QF(0) \simeq 0.
$
A right adjoint $r$ to $p$ is given by the formula 
$
X \mapsto [\Sigma X \leftarrow 0 \rightarrow 0]
$
and since
\[
\Dual_\QF\big([\Sigma X \leftarrow 0 \rightarrow 0]\big) \simeq [\Omega \Dual_\QF X \leftarrow \Omega \Dual_\QF X \rightarrow 0]
\]
we also find $\QF(\Dual_\QF(rX)) \simeq 0$ for all $X \in \C$ as desired.
\end{proof}
We thus obtain that both the vertical and horizontal fibre sequences in the diamond above are taken to fibre sequences by any additive functor. For group-like additive functors we even obtain the following strengthening from Proposition~\reftwo{proposition:retract}:

\begin{corollary}
\label{corollary:decomp-Q1}%
Let $\F\colon \Catp \rightarrow \E$ be a group-like additive functor. Then the following holds: 
\begin{enumerate}
\item
The Poincaré functor $\cyl\colon(\C,\QF) \rightarrow \Q_1(\C,\QF)$ and the inclusion $\Met(\C,\QF) \rightarrow \Q_1(\C,\QF)$ induce an equivalence
\[
\F(\C,\QF) \times\F(\Met(\C,\QF)) {\lrar} \F(\Q_1(\C,\QF)).
\]
\item
The functors $\cyl\colon(\C,\QF) \rightarrow \Q_1(\C,\QF)$ and $i_{\hyp}\colon\Hyp(\C) \rightarrow \Q_1(\C,\QF)$ induce an equivalence 
\[
\F(\C,\QF) \times \F(\Hyp(\C)) {\lrar} \F(\Q_1(\C,\QF)).
\]
\item
The functors $d_1\colon\Q_1(\C,\QF) \rightarrow (\C,\QF)$ and $p^{\hyp}\colon\Q_1(\C,\QF) \rightarrow \Hyp(\C)$ induce an equivalence 
\[
\F(\Q_1(\C,\QF)) {\lrar} \F(\C,\QF) \times \F(\Hyp(\C)).
\]
\end{enumerate}
\end{corollary}

As a consequence of the above, we obtain the following corollary, which plays a fundamental role throughout this paper. 

\begin{corollary}
\label{corollary:met-hyp}%
Let $\F\colon \Catp \rightarrow \E$ be a group-like additive functor. Then the functors $\lag\colon\Met(\C,\QF) \rightarrow \Hyp(\C)$ and $\ilag\colon\Hyp(\C) \rightarrow \Met(\C,\QF)$ induce inverse equivalences 
\[
\F(\Met(\C,\QF)) \simeq \F(\Hyp(\C)).
\]
\end{corollary}

\begin{proof}%
The composite 
\[
(\C,\QF) \times \Met(\C,\QF) \xrightarrow{(\cyl,\incl)} \Q_1(\C,\QF) \xrightarrow{(d_1,p_\hyp)} (\C,\QF) \times \Hyp(\C)
\]
is equivalent to the map $\id_{(\C,\QF)} \times \lag$. Since both constituents of this composite become equivalences after applying $\F$ by the previous corollary, $\F(\lag)$ is a retract of an equivalence and therefore an equivalence itself. 
Since the functor $\ilag$ is a one-sided inverse to $\lag$ at the level of Poincaré \(\infty\)-categories, it must induce the inverse equivalence after applying $\F$.
\end{proof}

Applying Corollary~\reftwo{corollary:met-hyp} to the group-like additive functor $(\C,\QF) \mapsto |\Cob^\F(\C,\QF)|$ for $\F$ a not necessarily group-like additive functor, we deduce immediately:

\begin{corollary}
\label{corollary:FCobhyp=FCobMet}%
The functors $\lag$ and $\ilag$ induce inverse equivalences 
\[
|\Cob^\F(\Met(\C,\QF))| \simeq |\Cob^\F(\Hyp(\C))|,
\]
for every additive functor $\F \colon \Catp \rightarrow \Sps$. 
\end{corollary}
This in particular gives an alternative proof of the second half of Proposition~\reftwo{proposition:CobMet=CobHyp} that does not use the algebraic Thom construction. 
To exploit Corollary~\reftwo{corollary:decomp-Q1} further, we need:

\begin{construction}
\label{construction:norm}%
Given two Poincaré \(\infty\)-categories $(\C,\QF),(\D,\QFD)$ and an exact functor $f\colon \C \rightarrow \D$ between the underlying $\infty$-categories, we obtain a Poincaré functor $\rN f\colon (\C,\QF) \rightarrow (\D,\QFD)$ by forming the composition
\[
\begin{tikzcd}
(\C,\QF) \ar[rr,"f^\hyp"] && \Hyp(\D) \ar[rr,"\id_\hyp"] && (\D,\QFD)
\end{tikzcd}
\]
using that the hyperbolic construction is both a left and a right adjoint to the forgetful functor $\Catp\to \Catx$. We refer to $\rN f$ as the \defi{norm} of $f$. Unwinding this construction, we find $(\rN f)(X) \simeq f(X) \oplus \Dual_\QFD f\op(\Dual_\QF X)$. 
\end{construction}

\begin{proposition}
\label{proposition:cobordism}%

Let $\F\colon \Catp \rightarrow \E$ be a group-like additive functor. Let $(\C,\QF)$ and $(\D,\QFD)$ be Poincaré $\infty$-categories and let $\phi = (f \leftarrow h \to g) \colon (\C,\QF) \to \Q_1(\D,\QFD)$ 
be a cobordism between two Poincaré functors $f,g\colon (\C,\QF) \rightarrow (\D,\QFD)$. Let $k= \fib(h \to f) \colon\C \rightarrow \D$ 
and let $\rN k\colon (\C,\QF) \rightarrow (\D,\QFD)$ be its norm. Then $\F(g)-\F(f) \sim \F(\rN k)$ as maps $\F(\C,\QF) \rightarrow \F(\D,\QFD)$.

\end{proposition}

\begin{proof}
By Corollary~\reftwo{corollary:decomp-Q1}, we have a pair of equivalences
\[
\begin{tikzcd}
\F(\D,\QFD) \oplus \F(\Hyp(\D)) \ar[rr,"\F(\cyl)\oplus \F(i_{\hyp})"] && \F(\Q_1(\D,\QFD)) \ar[rr,"(\F(d_1)\comma \F(p^{\hyp}))"] && \F(\D,\QFD) \oplus \F(\Hyp(\D)). 
\end{tikzcd}
\]
These equivalences are inverse to each other: indeed, the composite equivalence 
\[
\F(\D,\QFD) \oplus \F(\Hyp(\D)) \lrar \F(\D,\QFD) \oplus \F(\Hyp(\D))
\]
is equivalent to the identity since $p^{\hyp}i_{\hyp}$ and $d_1\cyl$ are equivalent to the respective identity functors while $d_1i_{\hyp}$ and $p^{\hyp}\cyl$ are equivalent to the respective zero functors. Consequently, we obtain a homotopy between the identity map $\id\colon\F(\Q_1(\D,\QFD)) \rightarrow \F(\Q_1(\D,\QFD))$ and the sum $\F(\cyl d_1)+\F(i_{\hyp} p^{\hyp})$, and hence a homotopy 
\[
\F(\phi) \sim \F(\cyl d_1\phi)+\F(i_{\hyp} p^{\hyp}\phi) = \F(\cyl f) + \F(i_{\hyp}k^{\hyp})
\]
of maps $\F(\C,\QF) \rightarrow \F(\Q_1(\D,\QFD))$. Post composing with the map $\F(d_0)\colon\F(\Q_1(\D,\QFD)) \rightarrow \F(\D,\QFD)$, we obtain a homotopy
\[
\F(g) = \F(d_0\phi) \sim \F(d_0\cyl f) + \F(d_0i_{\hyp}k^{\hyp}) = \F(f) + \F(\rN k)
\]
of maps $\F(\C,\QF) \rightarrow \F(\D,\QFD)$, as desired.
\end{proof}

\begin{corollary}
\label{corollary:inversion}%
For a group-like additive functor $\F \colon \Catp \rightarrow \E$, the inversion map on $\F(\C,\QF)$ is induced by the sum of the endofunctors $(\id_\C,-\id_\QF)$ and $\rN\Omega$ of $(\C,\QF)$.
\end{corollary}

\begin{proof} 
We resurrect the bent cylinder $\bcyl \colon (\C,\QF) \longrightarrow \Q_1(\C,\QF)$ with underlying functor
\[
X \longmapsto [X \oplus X \xleftarrow{\Delta} X \rightarrow 0]
\]
from Construction~\reftwo{construction:bentcyl}. By construction, it is a nullcobordism of $\id_{(\C,\QF)} + (\id_\C,-\id_\QF)$. We obtain the conclusion from Proposition~\reftwo{proposition:cobordism} by observing that the fibre of the diagonal $X \rightarrow X \oplus X$ is naturally equivalent to $\Omega X$.
\end{proof}

Next, we use Corollary~\reftwo{corollary:met-hyp} to determine the fundamental group of $|\Cob^\F(\C,\QF)|$. We base the calculation on the well-known analogue for the $\infty$-categories $\Span^\G(\C)$ for a small stable \(\infty\)-category $\C$ and an additive functor $\G \colon \Catx \rightarrow \Sps$
that we recalled in Proposition~\reftwo{proposition:pi0span}. 

Analogous to the construction in the non-hermitian case, we consider the diagram
\[
\begin{tikzcd}
\pi_0\F(\Hyp(\C)) \ar[d] & \ar[l,"\lag"'] \pi_0\F(\Met(\C,\QF)) \ar[r,"\met"] \ar[d] & \pi_0\F(\C,\QF) \ar[d] \\
\pi_1|\Cob^\F(\Hyp(\C))| & \ar[l] \pi_1|\Cob^\F(\Met(\C,\QF))| \ar[r] & \pi_1|\Cob^\F(\C,\QF)|,
\end{tikzcd}
\]
with the vertical maps induced by various instances of $\Hom_{\Cob^\F(\C,\QF)}(0,0) \simeq \F(\C,\QF)$. The lower left horizontal map is an isomorphism by Corollary~\reftwo{corollary:met-hyp}. Inverting it gives the square in the following:

\begin{theorem}
\label{theorem:pi1Cob}%
For a Poincaré \(\infty\)-category $(\C,\QF)$ and an additive functor $\F \colon \Catp \rightarrow \Sps$, the natural square
\[
\begin{tikzcd}
\pi_0\F\Met(\C,\QF) \ar[r,"\met"] \ar[d,"\lag"] & \pi_0\F(\C,\QF) \ar[d]\\
\pi_0\F\Hyp(\C) \ar[r] & \pi_1|\Cob^\F(\C,\QF)|
\end{tikzcd}
\]
of commutative monoids is cocartesian. 
\end{theorem}

Since the map $\lag$ is (split) surjective, this in particular describes $\pi_1|\Cob^\F(\C,\QF)|$ as the quotient monoid of $\pi_0\F(\C,\QF)$ identifying all metabolic objects with the hyperbolic objects on their Lagrangians. We thus, in particular, obtain an isomorphism $\pi_1|\Cob(\C,\QF)| \cong \GW_0(\C,\QF)$ with the Grothendieck-Witt group constructed in \S\refone{subsection:GW-group}. We discuss this further in \S\reftwo{section:GW} below.
For the proof we need:

\begin{proposition}
\label{proposition:boundarymapGB}%
The boundary map of the Bott-Genauer sequence 
\[
|\Cob^\F(\C,\QF\qshift{-1})| \longrightarrow |\Cob^\F(\Met(\C,\QF))| \longrightarrow |\Cob^\F(\C,\QF)|
\]
participates in a diagram
\[
\begin{tikzcd}
 & \F(\C,\QF) \ar[dl] \ar[dr] & \\
\Omega|\Cob^\F(\C,\QF)| \ar[rr,"\partial"] & & {|\Cob^\F(\C,\QF\qshift{-1})|}
\end{tikzcd}
\]
with the right hand map arising from the inclusion into the core, and the left hand map from the inclusion as the endomorphism of $0 \in \F(\C,\QF\qshift{1})$.%
\end{proposition}

Since the map $\pi_0\F(\C,\QF) \rightarrow \pi_1|\Cob^\F(\C,\QF)|$ is surjective by Theorem~\reftwo{theorem:pi1Cob}, this in particular determines the effect of the boundary map $\pi_1|\Cob^\F(\C,\QF)| \rightarrow \pi_0|\Cob^\F(\C,\QF\qshift{-1})|$.

\begin{proof}[Proof of Proposition~\reftwo{proposition:boundarymapGB}]
There is a commutative diagram of $\infty$-categories
\[
\begin{tikzcd}
	\F(\C,\QF) \ar[r] \ar[d] & \Cob^\F(\C,\QF\qshift{-1}) \ar[d] \ar[r] & \ast \ar[d] \\
	\Cob^\F(\Met(\C,\QF))_{0/} \ar[r] & \Cob^\F(\Met(\C,\QF)) \ar[r] & \Cob^\F(\C,\QF)
\end{tikzcd}
\]
The right square is simply obtained by applying $\Cob^\F(-)$ to the metabolic sequence. The left vertical map is given by the composite 
\[
\F(\C,\QF) \longrightarrow \F(\Met(\C,\QF\qshift{1})) \xrightarrow{\mathrm{cm}} \F(\Met\Met(\C,\QF\qshift{1})) \longrightarrow \Cob^\F(\Met(\C,\QF))_{0/},
\]
where the middle map is given by the one from Example \reftwo{example:idmetbordnull}. The commutativity of the left square comes from the fact that this map is a section of $\met(\Met)$. 

Applying realisation, using that $|\Cob^\F(-)|$ is additive, we obtain from the right square the boundary map $\Omega|\Cob^\F(\C,\QF))| \to |\Cob^\F(\C,\QF\qshift{-1})|$ of the Bott-Genauer sequence as the inclusion of the middle horizontal fibre. Since the lower left corner is contractible, the left hand square then factors the top horizontal map through the boundary map as desired. To finally see that the induced map $\F(\C,\QF) \to \Omega |\Cob^\F(\C,\QF)|$ is the one from the statement of the Proposition, one uses that $\mathrm{cm}$ is also a section of the map $\Met(\met)$: Indeed, this implies that in the commutative diagram
\[
\begin{tikzcd}
	\F(\C,\QF) \ar[r] \ar[dr] & \Cob^\F(\Met(\C,\QF))_{0/} \ar[r] \ar[d,"\met"] & \Cob^\F(\Met(\C,\QF)) \ar[d,"\met"] \\
	& \Cob^\F(\C,\QF)_{0/} \ar[r] & \Cob^\F(\C,\QF)
\end{tikzcd}
\]
the diagonal arrow is the fibre inclusion (over $0$) of the lower horizontal map. 
\end{proof}

\begin{proof}[Proof of Theorem~\reftwo{theorem:pi1Cob}]
Denote by $G(\C,\QF)$ the pushout of the diagram 
\[
\pi_0\F\Hyp(\C) \longleftarrow \pi_0\F\Met(\C,\QF) \longrightarrow \pi_0\F(\C,\QF),
\]
and similarly $W(\C,\QF)$ the pushout of 
\[
0 \longleftarrow \pi_0\F\Met(\C,\QF) \longrightarrow \pi_0\F(\C,\QF),
\]
giving a canonical map $G(\C,\QF) \rightarrow W(\C,\QF)$.
By construction, there is a natural map $G(\C,\QF) \rightarrow \pi_1|\Cob(\C,\QF)|$. Now, the discussion of the non-Poincaré case in Proposition~\reftwo{proposition:pi0span} implies that this map is an equivalence for hyperbolic categories: The square in Theorem~\reftwo{theorem:pi1Cob} for $\F \colon \Catp \rightarrow \Sps$ and input $\Hyp(\C)$ becomes that for $\F \circ \Hyp \colon \Catx \rightarrow \Sps$ and input $\infty$-category $\C$, under the equivalences $\Met(\Hyp(\C)) \simeq \Hyp(\Ar(\C))$ and $\Hyp(\C)^2 \simeq \Hyp(\Hyp(\C))$ from Corollary~\refone{corollary:met-hyp-splits} and Remark~\refone{remark:hyp-of-hyp}.

Let us now construct a diagram
\[
\begin{tikzcd}
\pi_1|\Cob^\F(\C,\QF\qshift{-1})| \ar[r] & \pi_1|\Cob^\F(\Met(\C,\QF))| \ar[r] & \pi_1|\Cob^\F(\C,\QF)| \ar[r,two heads] & \pi_0|\Cob^\F(\C,\QF\qshift{-1})| \\
G(\C,\QF\qshift{-1}) \ar[u] \ar[r] & G(\Hyp(\C)) \ar[r] \ar[u,"\cong"] & G(\C,\QF) \ar[u,two heads] \ar[r,two heads] & W(\C,\QF) \ar[u,"\cong"]
\end{tikzcd},
\]
whose upper sequence is induced by the metabolic fibre sequence via additivity and thus exact. Furthermore, the rightmost map of the top sequence is surjective as indicated, since the next term in the sequence is $\pi_0|\Cob^\F(\Met(\C,\QF))|$, which vanishes by Corollary~\reftwo{corollary:pi0CorCobMet}. The vertical maps are the evident ones (see Corollary~\reftwo{corollary:components-of-cob-II} for the right most one), except the second one, which is the composition
\[
G(\Hyp(\C)) \xrightarrow{\ilag} G(\Met(\C,\QF)) \longrightarrow \pi_1|\Cob^\F(\Met(\C,\QF))|.
\]
The left two horizontal maps in the lower sequence are
\[
G(\C,\QF\qshift{-1}) \longrightarrow G(\Met(\C,\QF)) \xrightarrow{\lag} G(\Hyp(\C))
\]
and
\[
\hyp \colon G(\Hyp(\C)) \xrightarrow{\ilag} G(\Met(\C,\QF)) \xrightarrow{\met} G(\C,\QF),
\]
respectively and the right one is the one constructed above. The second vertical map is an isomorphism by Corollary~\reftwo{corollary:met-hyp} and the claim for hyperbolic categories established above.

Now the middle square commutes by construction, the left one by Corollary~\reftwo{corollary:met-hyp} and the right by Proposition~\reftwo{proposition:boundarymapGB}. Furthermore, the lower sequence is exact at $G(\C,\QF)$ in the sense that two elements $x,y \in G(\C,\QF)$ have the same image in $W(\C,\QF)$ if and only if there are elements $w,z$ in the image of $G(\Hyp(\C))$ such that $x+w = z+y$: By the surjectivity of $\pi_0\F(\Hyp\C) \rightarrow G(\Hyp\C)$ this follows straight from the cocartesian diagram 
\[
\begin{tikzcd}
\pi_0\F\Met(\C,\QF) \ar[r] \ar[d] & \pi_0\F(\C,\QF) \ar[d]\\
\pi_0\F\Hyp(\C) \ar[r] & G(\C,\QF)
\end{tikzcd}
\]
by taking horizontal cokernels. It then follows formally that $G(\C,\QF)$ is in fact a group: Since $W(\C,\QF)$ is one, there is for every $a \in G(\C,\QF)$ an element $a' \in G(\C,\QF)$ such that $a+a'$ maps to $0$ in $W(\C,\QF)$. But then by exactness, there are $b,b' \in G(\Hyp(\C))$ with $a+a'+b'=b$, from which we can subtract $b$ to get an inverse to $a$, since $G(\Hyp(\C))$ is group.

Furthermore, the composition
\[
G(\C,\QF\qshift{-1}) \longrightarrow G(\Hyp(\C)) \longrightarrow G(\C,\QF)
\]
vanishes: By construction, the map $\met \colon G(\Met(\C,\QF)) \rightarrow G(\C,\QF)$ factors as 
\[
G(\Met(\C,\QF)) \xrightarrow{\lag} G(\Hyp(\C)) \xrightarrow{\hyp} G(\C,\QF)
\]
which identifies the composition above with
\[
G(\C,\QF\qshift{-1}) \longrightarrow G(\Met(\C,\QF)) \xrightarrow{\met} G(\C,\QF)
\]
which vanishes already at the level of $\infty$-categories. It is a bit tedious to check that the lower sequence is in fact exact at $G(\Hyp(\C))$. Luckily, we get away without doing so directly:

We deduce Theorem~\reftwo{theorem:pi1Cob} by two applications of the $4$-lemma. Applying one half of it to the right three columns (extended by $0$ to the right) gives surjectivity of the map $G(\C,\QF) \rightarrow \pi_1|\Cob^F(\C,\QF)|$ for every $(\C,\QF)$, in particular also for the left most column. This formally implies exactness at $G(\Hyp(\C))$ by a short diagram chase, whence the other half of the 4-lemma gives injectivity and thus the claim.
\end{proof}

\subsection{Isotropic decompositions of Poincaré \(\infty\)-categories}
\label{subsection:isotrop}%

We now describe a rather general situation which gives rise to cobordisms of Poincaré functors. We will use it to analyse the $\infty$-categories $\Q_n(\C,\QF)$; see Proposition~\reftwo{proposition:compQ} below. 

Let now $(\C,\QF)$ be a Poincaré \(\infty\)-category. %
Given a full subcategory $\Lag \subseteq \C$, we denote by $\Lag^{\perp} \subseteq \C$ the full subcategory spanned by the objects $\y \in \C$ such that $\Bil_\QF(x,y) \simeq 0$ for every $\x \in \Lag$.
Using $\Bil_\QF(x,y) \simeq \Hom_\C(x,\Dual_\QF(y))$, we immediately see that $\Dual_\QF(\Lag^{\perp}) \subseteq \C$ is the full subcategory of $\C$ consisting of the objects $\z \in \C$ that are right orthogonal to $\Lag$, i.e.\ for which $\Map_{\C}(\x,\z) \simeq 0$ for every $\x \in \Lag$.

\begin{definition}
\label{definition:isotropic}%
By an \defi{isotropic subcategory} of $(\C,\QF)$ we shall mean a full stable subcategory $\Lag \subseteq \C$ with the following properties:
\begin{enumerate}
\item
\label{item:1}%
$\QF$ vanishes on $\Lag$. 
\item
\label{item:2}%
The composite functor 
$
\Lag\op \longrightarrow \C\op \xrightarrow{\Dual_\QF} \C / \Lag^\perp
$
is an equivalence.
\end{enumerate}
\end{definition}

The first condition in particular implies $\Lag \subseteq \Lag^\perp$, and the second expresses a unimodularity condition on $\Lag$. In the ordinary theory of quadratic forms, the analogue of this condition is equivalent to the requirement that an isotropic subspace be a direct summand. It admits a convenient reformulation: 

\begin{lemma}
\label{lemma:unimod-vs-adjoint}%
For a stable subcategory $\Lag \subseteq \C$ the composite functor $\Lag\op \longrightarrow \C\op \xrightarrow{\Dual_\QF} \C / \Lag^\perp$ is an equivalence if and only if the inclusion of $\Lag$ into $\C$ admits a right adjoint. Furthermore, in this case $\Lag = (\Lag^\perp)^\perp$.
\end{lemma}

In particular, Lagrangians as appearing in the recognition criterion for Poinar\'e $\infty$-categories of pairings (see Definition~\refone{definition:lagrangian}) are examples of isotropic subcategories; we recall their definition in Definition~\reftwo{definition:lagrangianII} below.

\begin{proof}
The composite being an equivalence is equivalent to $\Lag \rightarrow \C \rightarrow \C / \Dual_\QF(\Lag^{\perp})$ being one. Since $\Dual_\QF(\Lag^{\perp})$ consists exactly of the right orthogonal of $\Lag$, it is closed under retracts in $\C$ and thus gives a Verdier inclusion into $\C$ by Proposition~\reftwo{proposition:inclusion-criterion}. Both the equivalence of the conditions in the statement and the last statement are then instances of Corollary~\reftwo{corollary:general-orthogonal}.
\end{proof}

\begin{remark}
\label{remark:kerp=DL}%
Applying the remainder of Corollary~\reftwo{corollary:general-orthogonal} in the situation at hand, we find that the kernel of the right adjoint $p \colon \C \rightarrow \Lag$ is given by $\Dual_\QF(\Lag^\perp)$, and thus $\Lag^\perp$ is the kernel of $p \circ \Dual_\QF$. 
\end{remark}

\begin{example}
\label{example:L=Lperpnotenough}%
The condition that $\Lag = (\Lag^\perp)^\perp$, or even that $\Lag=\Lag^\perp$, does not imply condition~\reftwoitem{item:2} of the definition of an isotropic category. For a concrete counterexample, fix a field $K$ of characteristic different from 2 and an element $a \in K \setminus \{0 \}$. Sending $T$ to $-T$ gives an involution on $K[T]$, which exchanges the ideals $(T-a)$ and $(T+a)$ in $K[T]$. Taking $R$ to be the localisation of $K[T]$ at the complement of $(T-a) \cup (T+a)$ we thus obtain a principal ideal domain with exactly two maximal ideals, and an involution which swaps them. Now consider the symmetric Poincar\'e structure on $\Dperf(R)$ provided by $R$ with this involution.

Denoting in this general situation the maximal ideals by $(m)$ and $(n)$, take $\Lag = \Dperf(R)_{m} \subseteq \Dperf(R)$ the subcategory spanned by the $m$-torsion complexes (i.e.\ those which become contractible upon inverting $m$, or equivalently whose homology is $m$-torsion). We claim that $\QF^\sym_{|\Lag} = 0$ and $\Lag = \Lag^\perp$, but $\Lag$ is not an isotropic subcategory: The inclusion $\Dperf(R)_{m} \subseteq \Dperf(R)$ does not admit a right adjoint (e.g.\ because such an adjoint would provide a retraction of the inclusion, but by devissage the induced map on $\K_0$ vanishes), so Lemma~\reftwo{lemma:unimod-vs-adjoint} applies. 

Since $\Dperf(R)_{m}$ it is generated as a stable $\infty$-category by $R/m$, see e.g.\ Example \reftwo{example:ore}, the claims boil down to the statement that for $X \in \Dperf(R)$ we have $\Bil_{\QF^\sym}(X,R/m) \simeq 0$ if and only if $X \in \Dperf(R)_{m}$. To see this in turn note that $\Dual_{\QF^\sym}(R/m) \simeq \Omega R/n$, which implies that $\Bil_{\QF^\sym}(X,R/m) \simeq 0$ if and only if $n$ acts invertibly on $X$, from which the statement follows by the classification of finitely generated modules over principal ideal domains.
\end{example}

\begin{definition}
\label{definition:homology}%
For an isotropic subcategory $\Lag$ of a Poincaré $\infty$-category $(\C,\QF)$, we define the \defi{homology $\infty$-category} $\Hlgy(\Lag)$ to be the cofibre of the inclusion $(\Lag,\QF) \rightarrow (\Lag^\perp,\QF)$ in $\Cath$. 
\end{definition}

Thus the underlying $\infty$-category is $\Lag^\perp/\Lag$ and per construction the hermitian structure is the left Kan extension of $\QF_{|(\Lag^\perp)\op}$ along the projection $(\Lag^\perp)\op \rightarrow (\Lag^\perp/\Lag)\op$. The next proposition, in particular, shows that this Kan extension is an optical illusion:%

\begin{proposition}
\label{proposition:homology}%
Let $\Lag$ be an isotropic subcategory of a Poincaré \(\infty\)-category $(\C,\QF)$. Then both $(\Bil_\QF)_{|(\Lag^\perp \times \Lag^\perp)\op}$ and $(\Lin_\QF)_{|(\Lag^\perp)\op}$ descend along the projection $(\Lag^\perp)\op \rightarrow (\Lag^\perp/\Lag)\op$ and give the bilinear and linear part of the hermitian structure on $\Hlgy(\Lag)$, which is Poincaré. The duality on $\Hlgy(\Lag)$ is induced by the functor $(\Lag^\perp)\op \rightarrow \Lag^\perp$ sending $X$ to $\fib(\Dual_\QF X \rightarrow \Dual_\QF pX)$, where $p$ denotes the right adjoint to $\Lag \subseteq \C$ and the arrow is induced by the counit.

In particular, the composite
\[
\Lag^{\perp} \cap \Dual_\QF(\Lag^{\perp}) \lrar \Lag^{\perp} \lrar \Lag^{\perp}/\Lag = \Hlgy(\Lag)
\]
canonically refines to an equivalence of Poincaré \(\infty\)-categories using the restriction of $\QF$ on the source.
\end{proposition}

In particular, $\Hlgy(\Lag)$ is equivalent to a full Poincaré subcategory of $(\C,\QF)$, which one may think of as the subcategory of harmonic objects for $\Lag$. We denote by  %
\[
\iota\colon \Hlgy(\Lag) \lrar (\C,\QF)
\]
the arising Poincaré-Verdier inclusion.

\begin{proof}
The first two statements follow from the general analysis of Kan-extended hermitian structures: By Lemma~\refone{proposition:left-kan-bilinear-linear} the linear and bilinear parts are given by the left Kan-extensions along $(\Lag^\perp)\op \rightarrow (\Lag^\perp/\Lag)\op$ of the restriction to $\Lag^\perp$. But they in fact descend along the projection: This is immediate from Condition~\reftwoitem{item:1} of Definition~\reftwo{definition:isotropic} for the linear part and from the definition of $\Lag^\perp$ in the case %
of the bilinear part. Note that this implies via the decomposition into linear and bilinear parts, that $\QF_{|(\Lag^\perp)\op}$ also descends along the projection $(\Lag^\perp)\op \rightarrow (\Lag^\perp/\Lag)\op$, as claimed above. It furthermore implies that the hermitian structure on $\Hlgy(\Lag)$ is also right Kan extended along this map, which we will use below.

For the equivalence of hermitian \(\infty\)-categories claimed in the statement, note first that by Lemma~\reftwo{lemma:quasi-split} and the comments thereafter the cofibre of the counit $pX \rightarrow X$ constitutes a right adjoint $q$ to the localisation $\Lag^\perp \rightarrow \Hlgy(\Lag)$. In fact, Lemma~\reftwo{lemma:quasi-split} implies that $q$ is an equivalence onto the kernel of $p \colon \Lag^\perp \rightarrow \Lag$, which is $\Dual_\QF(\Lag^\perp) \cap \Lag^\perp$ by Remark~\reftwo{remark:kerp=DL}. In particular, $q$ is also a right adjoint to the composite 
\[
c \colon \Lag^{\perp} \cap \Dual(\Lag^{\perp}) \rightarrow \Hlgy(\Lag)
\]
from the statement, which is thus also an equivalence. Now, right Kan extensions are computed by pullback along left adjoints, so the hermitian structure on $\Hlgy(\Lag)$ is given by $\QF \circ q\op$, which promotes $q$ and thus $c$ to an equivalence of hermitian \(\infty\)-categories. 

Finally, $\Lag^{\perp} \cap \Dual(\Lag^{\perp})$ is evidently closed under $\Dual_\QF$ so forms a Poincaré subcategory of $\C$, whence also $\Hlgy(\Lag)$ is Poincaré. The statement about the duality in $\Hlgy(\Lag)$ then follows from the formula for the inverse $q$ of $c$.
\end{proof}

\begin{definition}
\label{definition:lagrangianII}%
Let $\Lag \subseteq \C$ be an isotropic subcategory of a Poincaré $\infty$-category $(\C,\QF)$. We say that $\Lag$ is a \defi{Lagrangian} if $\Hlgy(\Lag) = 0$. We say that $(\C,\QF)$ is \defi{metabolic} if it contains a Lagrangian subcategory.
\end{definition}

As mentioned, Remark~\reftwo{lemma:unimod-vs-adjoint} shows that this definition of Lagrangian agrees with that discussed in Definition~\refone{definition:lagrangian}.

\begin{remark}
\label{remark:lagrangian}%
By Lemma~\reftwo{lemma:kernel-of-projection-to-verdier-quotient}, an isotropic subcategory $\Lag \subseteq \C$ is a Lagrangian if and only if the inclusion $\Lag \subseteq \Lag^{\perp}$ is an equivalence. Condition~\reftwoitem{item:2} of Definition~\reftwo{definition:isotropic} therefore yields a Verdier sequence
\[
\Lag \longrightarrow \C \longrightarrow \Lag\op
\]
exhibiting $\C$ as an extension of $\Lag$ by $\Lag\op$, where the right functor takes $X$ to $\fib(\Dual_\QF X \rightarrow \Dual_\QF pX)$ and $p$ is the right adjoint of the inclusion $\Lag\subseteq \C$. 
\end{remark}

\begin{examples}
\label{examples:hyp}%
\label{examples:Qn}%
\ %
\begin{enumerate}
\item
\label{item:metCQmeta}%
We showed in Proposition~\refone{proposition:recognize-poincare} that a Poincaré \(\infty\)-category $(\C,\QF)$ is metabolic if and only if it is of the form $\Pairings(\D,\QFD)$ for some hermitian \(\infty\)-category $(\D,\QFD)$. In fact, the Lagrangians in $(\C,\QF)$ are in one-to-one correspondence with representations of $(\C,\QF)$ as an $\infty$-category of pairings. Particular examples are the inclusion $\C \rightarrow \Met(\C,\QF)$ as the equivalences, and $\C \times 0 \subseteq \Hyp(\C)$.
\item
Extending the Lagrangian of the metabolic category, the full subcategory inclusion $i\colon\C \hrar \Q_1(\C)$
sending $X$ to $[0 \leftarrow X \rightarrow X]$ gives an isotropic subcategory $\Lag$, the adjoint $p$ witnessing Condition~\reftwoitem{item:2} of Definition~\reftwo{definition:isotropic} given by 
\[
[X \leftarrow Y \rightarrow Z] \longmapsto [0 \leftarrow \fib(Y \rightarrow X) \rightarrow \fib(Y \rightarrow X)].
\]
Thus $\Dual_{\QF_1}\Lag^\perp = \ker(p)$ is spanned by all diagrams with left pointing arrow an equivalence, whereas $\Lag^\perp$ itself consists of all diagram with right hand arrow an equivalence. Thus $\Hlgy(\Lag) \simeq (\C,\QF)$ embedded as the constant diagrams.
\item
\label{item:homlagQn-1}%
More generally, one can consider $\Lag$ given by the image of the inclusion $j_n \colon \C \rightarrow \Q_n(\C)$ as those diagrams which vanish away from $\{(i \leq n) \mid i \in \{0,...n\}\}$, and are constant on that subposet, i.e.
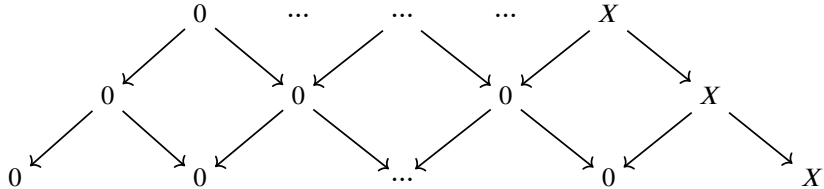
\[
\begin{tikzcd}[cramped]
 & & 0\ar[rd] \ar[ld] & \cdots & \ar[rd] \ar[ld] \cdots & \cdots & \ar[ld] \ar[rd,equal] X \\
 & 0 \ar[rd] \ar[ld] & & \ar[rd] \ar[ld] 0 & & 0 \ar[rd] \ar[ld] & & X \ar[rd,equal] \ar[ld] \\
0 & & 0 & & \cdots & & 0 & & X
\end{tikzcd}
\]
Formally this can be given by taking the embedding $\C \rightarrow \Q_1(\C)$ considered in the previous example and composing with the degeneracy $[n] \rightarrow [1]$ sending $n$ to $1$ and everything else to $0$. Using the Segal property of Lemma~\reftwo{lemma:segal} it is not difficult to see that the requisite adjoint $p$ is given by taking a diagram $\varphi \in \Q_n(\C)$ to the image of the fibre of the last left pointing arrow, namely $\varphi(n-1 \leq n) \rightarrow \varphi(n-1 \leq n-1)$. It follows that $\Dual_{\QF_n}(\Lag^\perp) = \ker(p)$ consists of all those diagrams $\varphi$ with the arrow $\varphi(n-1 \leq n) \rightarrow \varphi(n-1 \leq n-1)$ an equivalence (and thus all arrows $\varphi(i \leq n) \rightarrow \varphi(i \leq n-1)$ equivalences as well). From the explicit formula for the duality of $\Q_1(\C,\QF)$ from Example~\reftwo{examples:DescriptionsofQ} \reftwoitem{item:Q1-with-duality}, it then follows that $\Lag^\perp$ is spanned by the diagrams with the last right pointing arrow $\varphi(n-1 \leq n) \rightarrow \varphi(n \leq n)$ an equivalence, and so in total $\Hlgy(\Lag) \simeq \Q_{n-1}(\C,\QF)$ embedded in $\Q_n(\C,\QF)$ via the degeneracy $s_{n-1}$. 
\item
By Proposition~\reftwo{proposition:multiple-algebraic-Thom}, the full subcategory $\Lag_n \subseteq \Q_n(\C)$ spanned by those diagrams $\vphi\colon \Twar[n]\op \rightarrow \C$ for which $\vphi(0\leq 0)=0$ and $\vphi(0 \leq j) \rightarrow \vphi(i \leq j)$ is an equivalence for $i \leq j \in [n]$ is isotropic in $\Q_n(\C,\QF)$, and becomes a Lagrangian when considered as a subcategory of $\Null_n(\C,\QF)$, thus witnessing the fact that the latter is metabolic. This observation played a key role in the algebraic surgery equivalence of
section~\reftwo{subsection:algebraic-surgery}.
\end{enumerate}
\end{examples}

In generalisation of Corollary~\reftwo{corollary:decomp-Q1}, we now set out to prove: %

\begin{theorem}[Isotropic decomposition principle]
\label{theorem:lagrangian}%
\label{corollary:metabolic}%
Let $(\C,\QF)$ be a Poincaré $\infty$-category and $i\colon \Lag \rightarrow \C$ be the inclusion of an isotropic subcategory.
Let $\F\colon \Catp \rightarrow \E$ be a group-like additive functor. Then the Poincaré functors
$i_\hyp \colon \Hyp(\Lag) \to (\C,\QF)$ and $\iota \colon \Hlgy(\Lag) \to (\C,\QF)$
induce an equivalence 
\[
\F(\Hyp(\Lag)) \times \F(\Hlgy(\Lag)) \longrightarrow \F(\C,\QF).
\]
In particular, if $\Lag$ is a Lagrangian, $\F(\Hyp(\Lag)) \to \F(\C,\QF)$ is an equivalence.
\end{theorem}

We will explicitly construct an inverse. It relies on the following:

\begin{construction}
\label{construction:homology}%
Fix a Poincaré \(\infty\)-category $(\C, \QF)$
and the inclusion $i\colon \Lag \to \C$ of an isotropic subcategory with
right adjoint $p$. We note that the counit $ip \to \id_{\C}$
defines a surgery datum on the Poincaré object $\id_{(\C,\QF)}$ of $\Funx((\C, \QF), (\C, \QF))$. Performing surgery as in Proposition~\reftwo{proposition:surgequiv}, we obtain a
Poincaré object in $\Q_1(\Funx((\C, \QF), (\C, \QF)))$, in other
words, a Poincaré functor $\phi\colon (\C, \QF) \to \Q_1(\C, \QF)$. By
construction, $d_1\circ \phi=\id$, and we denote by $h$ the composite $
d_0 \circ \phi \colon (\C, \QF)\to (\C, \QF),$ that is, the result of surgery. In total, this gives a cobordism 
\[
\phi(X) = \bigl(X \leftarrow g(X) \to h(X) \bigr)
\]
where $g(X)$ is the fibre of the composite 
$X \simeq \Dual_\QF\Dual_\QF(X) \to \Dual_\QF(ip\Dual_\QF(X))$ whose second map is the dual of the counit, and $h(X)$ is the cofibre of the canonically induced map $ip(X) \to g(X)$.
\end{construction}

\begin{proof}[Proof of Theorem~\reftwo{theorem:lagrangian}]
We first claim that the Poincar\'e functor $h$ from Construction~\reftwo{construction:homology} takes values in $\Lag^{\perp} \cap \Dual_\QF(\Lag^{\perp}) \subseteq \C$, and so in particular factors as a composite 
\[
(\C,\QF) \xrightarrow{\wtl{h}} \Hlgy(\Lag) \xrightarrow{\iota} (\C,\QF)
\]
by Proposition \reftwo{proposition:homology}. To see this observe that both the cofibre of $ip(X) \rightarrow X$ and $\Dual_\QF (ip\Dual_\QF(X))$ belong to $\Dual_\QF(\Lag^\perp)$: The former because $\Dual_\QF(\Lag^\perp) = \ker(p)$ by Remark~\reftwo{remark:kerp=DL} and for the latter, we simply note $p \Dual_\QF(X) \in \Lag \subseteq \Lag^\perp$. Now, $h(X)$ participates in a cofibre sequence
\[
h(X) \longrightarrow \cof(ip(X) \rightarrow X) \longrightarrow \Dual_\QF(ip\Dual_\QF(X)),
\]
so also $h(X) \in \Dual_\QF(\Lag^\perp)$. Since $h$ commutes with the duality, its image is then also contained in $\Lag^\perp$, which gives the desired factorisation.

Now to start the proof proper, consider the composite
\[
\Hyp(\Lag) \oplus \Hlgy(\Lag) \xrightarrow{(i_\hyp,\iota)} (\C,\QF) \xrightarrow{(p^\hyp,\wtl{h})} \Hyp(\Lag) \oplus \Hlgy(\Lag)
\]
in $\Grp_{\Einf}(\E)$. We aim to show that it is the identity and to do so, we will analyse all four components in turn. Firstly, consider the composite $\wtl{h}\circ \iota$. Note that $ip(X) \simeq 0$ for $X \in  \Lag^\perp\cap \ker(p) = \Lag^\perp \cap \Dual_\QF(\Lag^\perp)$. Thus the cobordism $\phi(X)$ discussed above
consists of equivalences in this case, and so $\wtl{h}\iota \simeq \id_{\Hlgy(\Lag)}$ as desired.
Next, we consider $p^\hyp i_\hyp\colon \Hyp(\Lag) \to \Hyp(\Lag)$. Restricted to $\Lag \subseteq \Hyp(\Lag)$, we obtain the functor $X \mapsto (pi(X),p\Dual_\QF i(X))$, which is equal to $(X,0)$ since $pi = \id_\Lag$ and $\ker(p) = \Dual_\QF(\Lag^\perp) \supseteq \Dual_\QF(\Lag)$ as used before. From the adjointness properties of the hyperbolisation functor, we deduce $p^\hyp i_\hyp = \id_{\Hyp(\Lag)}$ as desired.

The map $p^\hyp \iota \colon \Hlgy(\Lag) \rightarrow \Hyp(\Lag)$ vanishes since $p^\hyp \iota = (p \iota)^\hyp$ and $\Lag^\perp \cap \Dual_\QF \Lag^\perp \subseteq \ker(p)$ by Remark~\reftwo{remark:kerp=DL} yet again.
Finally, we have $\wtl h i_\hyp \simeq (\wtl h i)_\hyp$ and $\wtl h  i \simeq 0$ as follows from the above cofibre sequence describing $h(X)$.

We have thus constructed a right inverse to the map of the theorem even at the level of categories. To finish the proof, it then remains to show that the composite map
\[
\begin{tikzcd}
\F(\C,\QF) \ar[rr,"(p^{\hyp}_*\comma \wtl{h}_*)"] && \F(\Hyp(\Lag)) \oplus \F(\Hlgy(\C)) \ar[rr,"((i_{\hyp})_*\comma \iota_*)"] && \F(\C,\QF)
\end{tikzcd}
\]
is homotopic to the identity. But this now readily follows by applying Proposition~\reftwo{proposition:cobordism} to the cobordism of Construction~\reftwo{construction:homology} and observing that $ip$ is identified by construction with the fibre of $g \rightarrow h$.
\end{proof}

\begin{example}
\label{example:left-split-isotropic-decomposition}%
Given a left split Verdier sequence of stable $\infty$-categories
\[
\begin{tikzcd}
\C \ar[r,"f"'] 
& \D \ar[r,"p"'] 
\ar[l,bend right=40,shift right=1.5ex,start anchor=west,end anchor=east,"g"']
\ar[l,phantom,shift right=1.3ex,start anchor=west,end anchor=east,"\myperp"]
& \E, 
\ar[l,bend right=40,shift right=1.5ex,start anchor=west,end anchor=east,"q"']
\ar[l,phantom,shift right=1.3ex,start anchor=west,end anchor=east,"\myperp"]
\end{tikzcd}
\]
the right split Verdier inclusion %
$
q \times f\op \colon \E \times \C\op \to \D \times \D\op
$
is %
a Lagrangian in $\Hyp(\D)$, as follows directly from the fact that $q(\E)$ and $f(\C)$ are respectively the left and right orthogonal complements of each other in $\D$. 
Theorem~\reftwo{theorem:lagrangian} then implies that 
\[
\F\Hyp(q) + \F\Hyp(f) \simeq \F[(q \times f\op)_\hyp] \colon \F\Hyp(\E) \oplus \F\Hyp(\C) \to \F\Hyp(\D)
\]
is an equivalence for any group-like additive functor $\F\colon \Catp \to \E$, so that $\F$ takes the sequence
\[
\Hyp(\C) \xrightarrow{\Hyp(f)} \Hyp(\D) \xrightarrow{\Hyp(p)} \Hyp(\E)
\]
to an exact sequence. Note that this is not a-priori evident since the above sequence is not a split Poincaré-Verdier sequence in general. Using Example~\reftwo{examples:additive}\;\reftwoitem{item:from-catx} and the fact that the original left split Verdier sequence is a retract of the last one on the level of underlying stable $\infty$-categories, this provides an alternative proof of the fact that any group-like additive functor on $\Catx$ sends left split (and similarly right split) Verdier sequences to exact sequences, see Remark~\reftwo{remark:leftsplitadd}\;\reftwoitem{item:left-additivity}.
\end{example}

Next, we use Theorem~\reftwo{theorem:lagrangian} to analyse the values of $\Q_n(\C,\QF)$ under group-like additive functors. To state the result consider the functor 
$f_n \colon \Q_n(\C,\QF) \rightarrow \C^n$ taking fibres of the left pointing maps along the bottom of a diagram $X$, i.e.
\[
X \longmapsto \big[\fib\big(X(0 \leq 1) \rightarrow X(0 \leq 0)\big), \dots, \fib\big(X(n-1 \leq n) \rightarrow X(n-1 \leq n-1)\big)\big].
\]
We then have:

\begin{proposition}
\label{proposition:compQ}%
The functors
\[
v_n \colon \Q_n(\C,\QF) \rightarrow (\C,\QF) \quad \text{and} \quad f_n^\hyp \colon \Q_n(\C,\QF) \rightarrow \Hyp(\C)^n
\]
the former induced by the inclusion $[0] \rightarrow [n]$, combine into an equivalence
\[
\F(\Q_n(\C,\QF)) \simeq \F(\Hyp(\C))^{n} \oplus \F(\C,\QF)
\]
for every group-like additive $\F \colon \Catp \rightarrow \E$. In fact, these equivalences give an identification of the simplicial $\Einf$-group  $\F\Q(\C,\QF)$ in $\E$ with the bar construction $\mathrm{B}(0,\F\Hyp(\C),\F(\C,\QF))$ of $\F(\Hyp(\C))$ acting on $\F(\C,\QF)$ via the hyperbolisation map $\hyp \colon \F(\Hyp(\C)) \longrightarrow \F(\C,\QF)$.
\end{proposition}

\begin{proof}
We proceed by induction. For $n=0$ there is nothing to show. Using the isotropic subcategory $j_{n+1} \colon \C \rightarrow \Q_{n+1}(\C)$ described in Example~\reftwo{examples:Qn} \reftwoitem{item:homlagQn-1}, we find an equivalence
\[
((j_{n+1})_\hyp, s_{n}) \colon \F\Hyp(\C) \oplus \F\Q_{n}(\C,\QF) \longrightarrow \F\Q_{n+1}(\C,\QF)
\]
as a consequence of Theorem~\reftwo{theorem:lagrangian}. It is readily checked that this equivalence translates the map $(f^\hyp_{n+1},v_{n+1})$ to the matrix
\[
\begin{pmatrix} 0 & f^\hyp_{n}\\ \id_{\Hyp(\C)} & 0 \\ 0 & v_{n}\\\end{pmatrix}\ \colon \F\Hyp(\C) \oplus \F\Q_{n}(\C,\QF) \longrightarrow \F(\Hyp(\C))^{n} \oplus \F\Hyp(\C) \oplus \F(\C,\QF).
\]
This matrix represents an equivalence by inductive assumption, which implies the first claim.

To obtain an identification with the bar construction, we first note that the bar construction $\mathrm{B}(M,R,N)$ of an action of $R$ on $N$ from the left and on $M$ from the right in a semi-additive $\infty$-category is the left Kan extension along the inclusion of the coequaliser diagram $(\Delta_{\leq 1})\op$ into $\Delta\op$ of the diagram
\[
\begin{tikzcd}
M \oplus R \oplus N \ar[r,shift left] \ar[r,shift right] & M \oplus N
\end{tikzcd}
\]
containing the two action maps; this follows directly by evaluation of the pointwise formulae for left Kan extensions. By the calculations above, $d_1 \colon \F\Q_1(\C,\QF) \rightarrow \F(\C,\QF)$ is identified with the projection $\F(\Hyp) \times \F(\C,\QF) \rightarrow \F(\C,\QF)$ and it is readily checked that $d_0 \colon \F\Q_1(\C,\QF) \rightarrow \F(\C,\QF)$ gets identified with the sum of the identity of $\F(\C,\QF)$ and the hyperbolisation map under the equivalence of Proposition~\reftwo{proposition:compQ}. We therefore obtain a map of simplicial objects 
\[
\mathrm{B}(0,\F(\Hyp(\C)),\F(\C,\QF)) \longrightarrow \F\Q(\C,\QF)
\]
and one readily unwinds the construction to find it given by the maps we just checked to be equivalences. 
\end{proof}

\begin{corollary}
\label{corollary:CobgrpdifFgrp}%
When $\F \colon \Catp \rightarrow \Sps$ is additive and group-like,
the $\infty$-category $\Cob^\F(\C,\QF)$ is an $\infty$-groupoid.
\end{corollary}
\begin{proof}
For a bar construction $\mathrm{B}(M,R,N)$ in $\Sps$ (which is automatically a Segal space), the equivalences in $\mathrm{B}(M,R,N)_1 = M \oplus R \oplus N$ consist precisely of $M \oplus R^\times \oplus N$ by direct inspection, where $R^\times \subseteq R$ consists of the units in the $\Eone$-monoid $R$. In particular, its associated $\infty$-category is an $\infty$-groupoid if and only if $R$ is an $\Eone$-group, as is the case in the situation of the corollary by assumption.
\end{proof}

\begin{remarks}
\label{remarks:FQcompletebord}%
\ %
\begin{enumerate}
\item The considerations in the proof of Corollary \reftwo{corollary:CobgrpdifFgrp} also show that for group-like $\F$, the Segal space $\F\Q(\C,\QF)$ is complete if and only if $\F(\Hyp\C)$ vanishes (see section~\reftwo{subsection:bordism} below for a detailed discussion of such functors). 

\item As a consequence of Proposition~\reftwo{proposition:compQ}, the functors $f^\hyp_{n+1} \colon \hMet_n(\C,\QF) \rightarrow \Hyp(\C)^{n+1}$ induce an equivalence
$\F(\hMet(\C,\QF)) \simeq \mathrm{B}(0,\F\Hyp(\C),\F\Hyp(\C))$
for every group-like additive $\F \colon \Catp \rightarrow \E$. This follows from the fact that the sequence defining $\hMet_n(\C,\QF)$ is in fact a split Poincaré-Verdier sequence by Lemma~\reftwo{lemma:boundarysplit}, so gives rise to a fibre sequence after applying $\F$.

\item
Proposition~\reftwo{proposition:compQ} can also be obtained directly using the Segal property of the simplicial space $\F\Q(\C,\QF)$ and the bar construction, together with the computation of $\F\Q_1(\C,\QF)$ from Corollary~\reftwo{corollary:decomp-Q1}.
\end{enumerate}
\end{remarks}

\subsection{The group-completion of an additive functor}
\label{subsection:group-complete}%

Our goal in this section is to study the behaviour of space-valued additive functors under the hermitian $\Q$-construction, or equivalently of the assignment $\F \mapsto |\Cob^\F(-)|$. %
This is based on the following observation: For any additive functor $\F \colon \Catp \rightarrow \Sps$, there is a natural cartesian square 
\[
\begin{tikzcd}
\F(\C,\QF) \ar[r] \ar[d] & \Cob^\F(\C,\QF)_{0/} \ar[d] \\
\{0\} \ar[r] & \Cob^\F(\C,\QF)
\end{tikzcd}
\]
in $\Cat$, since $\Hom_{\Cob^\F(\C,\QF)}(0,0) \simeq \F(\C,\QF)$, which is immediate from our discussion of Segal spaces in section~\reftwo{subsection:cobcat}. 
As an application of the isotropic decomposition principle, we saw in Corollary \reftwo{corollary:CobgrpdifFgrp} that $\Cob^\F(\C,\QF)$ is a groupoid if $\F$ is assumed group-like, and thus $\F \simeq \Omega|\Cob^\F|$ in this case. Together with the additivity theorem this suffices to recognise $\Omega|\Cob^\F-|$ as the group-completion of a not-necessarily group-like additive $\F$ by means of \cite{HTT}*{Proposition 5.2.7.4 (3)} and explicit inspection. In the present section, we shall go a bit further and prove: 

\begin{theorem}
\label{theorem:suspension}%
Let $\F\colon \Catp \rightarrow \Sps$ be an additive functor and consider the square of space valued functors
\[
\begin{tikzcd}
\F \ar[d] \ar[r]  &{ | \Cob^\F(-)_{0/}|} \ar[d] \\
\ast\ar[r] &  {|\Cob^\F(-)|}
\end{tikzcd}
\]
Then we have:
\begin{enumerate}
\item
\label{item:Funsaddsusp}%
The square is cocartesian in $\Funadd(\Catp,\Sps)$, and so exhibits $|\Cob^{\F}(-)|$ as the suspension of $\F$ in $\Funadd(\Catp,\Sps)$.
\item
\label{item:Funsaddcocart}%
If $\F$ is group-like then the square is also cartesian, so that the canonical map
$
\tau_\F\colon\F \to \Om|\Cob^\F(-)|
$
in $\Funadd(\Catp,\Sps)$ is an equivalence.
\end{enumerate}
\end{theorem}

Unwinding, one finds $\tau_\Poinc$ simply given by taking a Poincar\'e object $X \in \Poinc(\C,\QF)$ to the cobordism $[0 \leftarrow X \rightarrow 0] \in \Hom_{\Cob(\C,\QF)}(0,0)$ and then mapping it forward to $\Omega|\Cob(\C,\QF)|$.

\begin{remark}
\label{remark:alternativeproof}%
We have just argued that statement~\reftwoitem{item:Funsaddcocart} is a consequence of Corollary~\reftwo{corollary:CobgrpdifFgrp}. One can also give a more direct proof by making use of the Segal property of $\Q(\C,\QF)$ instead: By the translation into Segal spaces and \cite{HTT}*{Theorem 6.1.3.9}, it suffices to show that $\F\hMet(\C,\QF) \to \F\Q(\C,\QF)$ is an equifibred map of simplicial spaces, see also \cite{Rezkhatespistarkan}*{Proposition 2.4}. To do so, the Segal condition implies that it suffices to check that the squares
\[
\begin{tikzcd}
\F(\hMet_2(\C,\QF)\qshift{1}) \ar[d,"d_i"]\ar[r] & \F(\Q_2(\C,\QF\qshift{1})) \ar[d, "d_i"] \\
\F(\hMet_1(\C,\QF\qshift{1})) \ar[r] & \F(\Q_1(\C,\QF\qshift{1})) 
\end{tikzcd}
\]
are cartesian for $i = 0,1,2$. For $i=1,2$, these squares are split Poincaré-Verdier prior to applying $\F$ by Lemma~\reftwo{lemma:boundarysplit} and Corollary~\reftwo{corollary:PV-proj-pullback} and for $i=0$ the induced map on vertical fibres (over $0$, but this suffices since $d_0$ induces a surjection on $\pi_0$) identifies with 
\[
\ilag \colon \F(\Hyp(\C)) \longrightarrow \F(\Met(\C,\QF))
\]
which is an equivalence by Corollary \reftwo{corollary:met-hyp}.
\end{remark}

Before diving into the proof of \reftwo{theorem:suspension}~\reftwoitem{item:Funsaddsusp}, which will occupy most of this section, let us draw the desired consequences. We put:

\begin{definition}
Given an additive functor $\F \colon \Catp \rightarrow \Sps$ we denote $\Omega|\Cob^\F(-)|$ by $\F^\grp$ and refer to it as the \emph{group completion} of $\F$. 
\end{definition}

This is justified by:

\begin{corollary}
\label{corollary:universal}%
\label{corollary:Cobgrp}%
Let $\F\colon \Catp \to \Sps$ be an additive functor.
The natural map $\F \rightarrow \F^\grp$ exhibits $\F^\grp$ as universal among group-like additive functors receiving a map from $\F$; that is, the operation $\F \mapsto \F^\grp$ is left adjoint to the inclusion 
\[
\Funadd(\Catp,\Grp_{\Einf}(\Sps)) \subseteq \Funadd(\Catp,\Mon_{\Einf}(\Sps)) \simeq \Funadd(\Catp,\Sps).
\]
Moreover, the natural map
\[
|\Cob^\F(\C,\QF)| \longrightarrow |\Cob^{\F^\grp}(\C,\QF)| \simeq \Cob^{\F^\grp}(\C,\QF)
\]
is an equivalence.
\end{corollary}

\begin{proof}
Part~\reftwoitem{item:Funsaddcocart} of Theorem~\reftwo{theorem:suspension} implies that the unit $\id \Rightarrow \Omega|\Cob^{(-)}|$ is an equivalence on all group-like additive functors. Thus $|\Cob^{(-)}|$ restricts to a fully faithful functor on $\Funadd(\Catp,\Grp_{\Einf}(\Sps))$, so the corollary is a consequence of the following general Lemma~\reftwo{proposition:grpcomp} applied to the semi-additive $\infty$-category $\Fun^\add(\Catp,\Sps)$.
\end{proof}

\begin{lemma}
\label{proposition:grpcomp}%
Let $\E$ be a semi-additive $\infty$-category which admits suspensions and loops, and let $\E_{\grp} \subseteq \E$ be the full subcategory spanned by the group-like objects. Then the following holds:
\begin{enumerate}
\item
The full subcategory $\E_{\grp} \subseteq \E$ is closed under any limits and colimits that exist in $\E$, and the images of both the suspension and loop functors $\Sig,\Om\colon  \E \rightarrow \E$ are contained in $\E_{\grp}$. In particular, we may consider $\Om\Sig\colon  \E \rightarrow \E$ as a functor from $\E$ to $\E_{\grp}$.
\item
If the suspension functor $\Sig\colon  \E_{\grp} \rightarrow \E_{\grp}$ is fully faithful then the unit map $u\colon \id \Rightarrow \Om\Sig$ exhibits $\Om\Sig$ as left adjoint to the inclusion $\E_{\grp} \rightarrow \E$.
\item
\label{item:suspuniteq}%
If the suspension functor $\Sig\colon  \E_{\grp} \rightarrow \E_{\grp}$ is fully faithful then for every object $A \in \E$, the suspension of the unit $\Sigma u\colon \Sigma A \rightarrow \Sigma \Omega \Sigma A$ is an equivalence.
\end{enumerate}
\end{lemma}
\begin{proof}
The first claim follows from the fact that $X \in \E$ being group-like can be detected on the level of both the represented functor $\Map(-,X)$ and the corepresented functor $\Map(X,-)$ (which automatically take values in commutative monoid objects since $\E$ is semi-additive), and that loop spaces are always group-like. 

To prove the second claim, it suffices to check that under the given assumptions the natural transformations $u_{\Om\Sig X},\Om\Sig u_{X}\colon  \Om\Sig X \rightarrow \Om\Sig\Om\Sig X$ are both equivalences \cite{HTT}*{Proposition 5.2.7.4}. But since $\Om\Sig$ is a monad, these two natural transformations admit a common section (the multiplication of the monad) and $\Sig\colon \E_{\grp} \rightarrow \E_{\grp}$ being fully faithful implies that $u$ is a natural equivalence on all group-like objects of $\E$.

To see the final claim, it suffices to show that for all $B\in \E$, the induced map $(\Sigma u)^*\colon \Hom_\E(\Sigma\Omega\Sigma A,B) \to \Hom_\E(\Sigma A, B)$ is an equivalence. Under the adjunction equivalences, this map identifies with the map $u^*\colon \Hom_\E(\Omega\Sigma A,\Omega B) \to \Hom_\E(A,\Omega B)$ which is an equivalence by the first parts.
\end{proof}

\begin{remark}
\label{remark:unitcounitQ}%
While the unit map \(\F \rightarrow \Omega|\Cob^\F(-)| = \Omega|\F\Q(-\qshift{1})|\) of the adjunction arising from Theorem~\reftwo{theorem:suspension} is quite explicit, the counit \(\Omega| \F\Q(-\qshift{1})| \rightarrow \F\) is more elusive. However, the composite 
\[
|\Omega \F\Q(-\qshift{1})| \rightarrow \F \rightarrow \Omega|\F\Q(-\qshift{1})|
\]
of the counit and unit can be identified with the \emph{negative} of the canonical limit-colimit interchange map.
In case $\F$ is group-like, the unit map is an equivalence, so this determines the counit for such $\F$. In the cases $\F = \Poinc, \Core, \Core^\hC$, or more generally any additive $\F$ for which the component of $0$ in $\F(\C,\QF)$ is contractible, the source of the counit vanishes. These cases cover all additive functors of interest to us.

To see the claim, recall from Theorem~\reftwo{theorem:suspension} that the internal suspension functor $\Sig\colon \Funadd(\Catp,\Sps) \to \Funadd(\Catp,\Sps)$ is implemented by the formula $\F \mapsto |\F\Q(-\qshift{1})|$. Iterating this twice, we see that the double suspension functor $\Sig \circ \Sig$ is implemented by the formula $\F \mapsto |\F\Q(\Q(-\qshift{1})\qshift{1})|$.
Now the double suspension functor has an automorphism coming from the swap automorphism of the 2-sphere $S^2 = S^1 \wedge S^1$. On the level of the iterated $\Q$-construction this automorphism is induced by the natural equivalence $\Q_n(\Q_m(-\qshift{1})\qshift{1}) \simeq \Q_m(\Q_n(-\qshift{1})\qshift{1})$ swapping the two shifted $\Q$-constructions. %
Considering now the unit-counit composite and the limit exchange map as two maps from $|\Omega \F\Q(-\qshift{1})|$ to $\Omega|\F\Q(-\qshift{1})|$, we see that the former is homotopic to the mate transformation of the identify $\id\colon \Sig \circ \Sig \to \Sig \circ \Sig$, while the latter %
to the mate of the swap automorphism.
But the swap automorphism is homotopic to the negative of the identity in any semi-additive $\infty$-category with finite colimits by Yoneda's lemma and the analogous claim for double loop spaces.

\end{remark}

Let us now move towards the proof of Theorem~\reftwo{theorem:suspension}~\reftwoitem{item:Funsaddsusp}. We start by recalling a Segal space model for the square from its statement (prior to realisation). 
From 
Lemma~\reftwo{lemma:Segal_slice_asscat}, we have
\[
\Cob^\F(\C,\QF)_{0/} \simeq  \asscat(\F\hMet(\C,\QF\qshift{1}))
\]
where $\hMet_n(\C,\QF\qshift{1})$ is the fibre of the map
$\dec(\Q(\C,\QF\qshift{1})) \to \Q_0(\C,\QF\qshift{1}) = (\C,\QF\qshift{1})$.
We then recall that the maps $d_0$ in $\Delta$ induce a map $\dec(\Q(\C,\QF\qshift{1})) \to \Q(\C,\QF\qshift{1})$ which we may restrict to the map
$\pi \colon \hMet(\C,\QF\qshift{1}) \to \Q(\C,\QF\qshift{1})$,
modelling the right vertical map in the square upon applying $\F$ and passing to associated categories. Moreover, using \cite{HTT}*{Lemma 6.1.3.16}, which also defines the notion of split simplicial objects, we obtain:

\begin{observation}
\label{lemma:split}%
The simplicial objects $\hMet(\C,\QF)$ and $\dec(\Q(\C,\QF))$ extend to split simplicial objects over the zero Poincaré \(\infty\)-category and $(\C,\QF)$, respectively. 
\end{observation}

In particular, $|\F\Null(\C,\QF)| \simeq \ast$ and $|\F\dec(\Q(\C,\QF))| \simeq \F(\C,\QF)$; the former of which of course also follows from Lemma \reftwo{lemma:Segal_slice_asscat} and the latter can be established by similar means. 

At any rate, it follows from the discussion in section~\reftwo{subsection:rezk} that we have a commutative diagram of fibre sequences of simplicial Poincar\'e $\infty$-categories
\[
\begin{tikzcd}
	\const(\C,\QF\qshift{-1}) \ar[r,"\iota"] \ar[d,hook] & \hMet(\C,\QF) \ar[r,"\pi"] \ar[d,hook] & \Q(\C,\QF) \ar[d,equal] \\
	\const(\Met(\C,\QF)) \ar[r] & \dec(\Q(\C,\QF)) \ar[r] & \Q(\C,\QF)
\end{tikzcd}
\]
since $\ker(d_0) \colon \Q_1(\C,\QF) \to \Q_0(\C,\QF)$ canonically identifies with $\Met(\C,\QF)$, compare Example~\reftwo{examples:DescriptionsofQ}~\reftwoitem{item:met-vs-Q-1}. Unravelling definitions yields the following description of the functor $\iota_n \colon \C \to \hMet_n(\C,) \subseteq \Q_{1+n}(\C)$: It sends $X \in \C$ to the diagram $\vphi_{X}\colon \Twar[1+n] \rightarrow \C$ given by
\[
\vphi_{X}(i \leq j) = \left\{\begin{matrix} X & 0 = i < j \\ 0 & \text{otherwise}\end{matrix}\right.
\]
in which all the maps between the various $X$'s are identities.

\begin{observation}
\label{observation:for-annoying-sequence}%
The above sequences are split Poincar\'e-Verdier in each degree. Indeed, by Lemma~\reftwo{lemma:boundarysplit}, the maps $d_0 \colon \Q_{1+n}(\C,\QF) \rightarrow \Q_n(\C,\QF)$ are split Poincaré-Verdier projections, and the left adjoint of $d_0$ is given via extension by $0$ and thus factors through the underlying $\infty$-categories of $\Null_n(\C,\QF) \rightarrow \Q_{1+n}(\C,\QF)$, whence Corollary~\reftwo{corollary:split-poincare-projection-inclusion} gives the claim.
\end{observation}

Applying an additive functor $\F\colon \Catp \rightarrow \Sps$ levelwise to the upper horizontal sequence of the rectangle above then yields a sequence of Segal spaces which corresponds to the fibre sequence of $\infty$-categories
\[
\F(\C,\QF) \lrar \Cob^{\F}(\C,\QF)_{0/} \lrar \Cob^{\F}(\C,\QF).
\]
from the start of this section.

We will need to make use of the dual $\Q$-construction, denoted $\dualQ(\C,\QF)$, which we discuss next.
\begin{lemma}
\label{lemma:dqexists}%
The functors $\Q_n \colon \Catp \rightarrow \Catp$ and $\Q_n \colon \Catx \rightarrow \Catx$ admit left adjoints $\dualQ_n$, given by tensoring with the poset $\I_n$.
\end{lemma}
These adjoints participate in the two squares 
\[
\begin{tikzcd}
\Catp \ar[r,"{\dualQ_n}"] \ar[d,"\fgt"] & \Catp \ar[d,"\fgt"] \\
\Catx \ar[r,"{\dualQ_n}"] & \Catx 
\end{tikzcd}
\qquad\qquad
\begin{tikzcd}
\Catp \ar[r,"{\dualQ_n}"] & \Catp \\
\Catx \ar[r,"{\dualQ_n}"] \ar[u,"\Hyp "'] & \ar[u,"\Hyp"'] \Catx
\end{tikzcd}
\]
obtained by passing to left adjoints everywhere in the analogous diagrams involving the respective $\Q$-constructions.

\begin{proof}
Recall that for $[n] \in \Del$, we have denoted by $\I_n$ the full subposet of $\Twar(\Delta^n)$ spanned by the arrows of the form $(i \leq j)$ for $j \leq i+1$. From Examples~\reftwo{examples:DescriptionsofQ}, we find $\Q_n \simeq (-)^{\I_n}$, which by Proposition~\refone{proposition:tensor-mapping-property} has $(-)_{\I_n}$ as a left adjoint when regarded as a functor $\Cath \rightarrow \Cath$. As an application of Proposition~\refone{proposition:posets-of-faces}, we find that $(\C,\QF)_{\I_n}$ is Poincaré whenever $(\C,\QF)$ is and from Remark~\refone{remark:comparing-1} and Proposition~\refone{proposition:basic-properties-hermitian-functor-cats} we then find an equivalence of Poincaré \(\infty\)-categories
\[
\Funx((\C,\QF)_{\I_n},(\D,\QFD)) \simeq \Funx((\C,\QF),(\D,\QFD)^{\I_n})
\]
which according to Corollary~\refone{corollary:forms-in-functor-cats} gives the claim by passing to Poincaré objects.
\end{proof}

Now, recall that there is a canonical equivalence $\Fun^\mathrm{L}(\Catp,\Catp) \simeq \Fun^\mathrm{R}(\Catp,\Catp)\op$ for example as an immediate consequence of Lurie's straightening equivalences, which makes both $\infty$-categories equivalent to that of bicartesian fibrations over $\Delta^1$ with both fibres identified with $\Catp$; the superscripts $\mathrm{L}$ and $\mathrm{R}$ indicate left and right adjoint functors, respectively. In particular, as the $\Q$-construction is a simplicial object, the left adjoints above assemble into a cosimplicial object.

\begin{definition}
\label{definition:dualmet}%
Let $(\C,\QF)$ be an hermitian $\infty$-category. We will denote by 
$\dualQ(\C,\QF)$ the cosimplicial hermitian \(\infty\)-category obtained by applying the left adjoint of $\Q_n$ in each degree. %
Moreover, we define $\dMet_n(\C,\QF)$ to be the Poincaré-Verdier quotient of $\dualQ_{n+1}(\C,\QF)$ by the image of the functor $\C= \dualQ_{0}(\C) \rightarrow \dualQ_{1+n}(\C)$ induced by the inclusion $[0] \rightarrow [1+n]$.
\end{definition}

\begin{remark}
Note that Proposition~\reftwo{proposition:example-cosieve-2} shows that there is a Poincaré-Verdier sequence
\[
(\C,\QF) \longrightarrow \dualQ_{1+n}(\C,\QF) \longrightarrow \dMet_n(\C,\QF).
\]
In particular, the functor $\dMet_{n}$ is the cofibre of the natural transformation $\dualQ_{\{0\}} \Rightarrow \dualQ_{n+1}$, while the functor $\hMet_n(-)$ is the fibre of the natural transformation $\Q_{n+1} \rightarrow \Q_{0}$. We conclude that the association $(\C,\QF) \mapsto \dMet_n(\C,\QF)$ is left adjoint to $(\D,\QFD) \mapsto \hMet_n(\D,\QFD)$.
\end{remark}

\begin{proof}[Proof of Theorem~\reftwo{theorem:suspension}~\reftwoitem{item:Funsaddsusp}]

We have to show that the square 
\[
\begin{tikzcd}
{\Nat(|\F\Q(-\qshift{1})|,\G)} \ar[d] \ar[rr] & & \Nat({|\F\hMet(-\qshift{1})|,\G)}  \ar[d] \\
\ast\ar[rr] & & \Nat(\F,\G)
\end{tikzcd}
\]
is cartesian for every additive $\G \colon \Catp \rightarrow \Sps$. But we calculate 
\[
\Nat(|\F\Q(-\qshift{1})|,\G) \simeq \lim_{[n] \in \Delta} \Nat(\F\Q_n(-\qshift{1}),\G) \simeq \lim_{[n] \in \Delta} \Nat(\F,\G\dualQ_n(-\qshift{-1}))
\]
and similarly for the upper right hand term. Commuting limits, it thus suffices to show that
\[
\dualQ_n(\C,\QF\qshift{-1}) \longrightarrow \dMet_{n}(\C,\QF\qshift{-1}) \longrightarrow (\C,\QF)
\]
is split Poincaré-Verdier for every Poincaré $\infty$-category $(\C,\QF)$. 
Indeed, it is immediate from adjointness that it is a cofibre sequence in $\Catp$, so it remains to check that the composite 
\[
\dualQ_n(\C,\QF\qshift{-1}) \xrightarrow{d_0} \dualQ_{1+n}(\C,\QF\qshift{-1}) \lrar \dMet_n(\C,\QF\qshift{-1})
\]
is a Poincaré-Verdier inclusion. But from the equivalence $\dualQ_n(\C,\QF) \simeq (\C,\QF)_{\I_n}$, we find the first map to be such an inclusion by Proposition~\reftwo{proposition:example-cosieve-2}. Thus the Poincaré structure on $\dualQ_n(\C,\QF\qshift{-1})$ is obtained from that on $\dualQ_{1+n}(\C,\QF\qshift{-1})$ by pullback along $d_0$, or equivalently by left Kan extension along (the opposite of) the right adjoint $\dualQ_{1+n}(\C,\QF\qshift{-1}) \rightarrow \dualQ_n(\C,\QF\qshift{-1})$ to $d_0$. We are therefore done if we show that this right adjoint factors through $\dMet_{n}(\C,\QF\qshift{-1})$. But this follows from the corresponding statement for the left adjoint of $d_0 \colon \Q_{1+n}(\C,\QF\qshift{1}) \rightarrow \Q_n(\C,\QF\qshift{1})$ factoring through $\hMet_n(\C,\QF\qshift{1})$ in Observation~\reftwo{observation:for-annoying-sequence}, since the adjunction $\Q_n \vdash \dualQ_n$ is compatible with the passage to underlying $\infty$-categories by the discussion after Lemma~\reftwo{lemma:dqexists}; the same argument gives the claim for $\hMet_n \vdash \dMet_n$.
\end{proof}

\subsection{The spectrification of an additive functor}
\label{subsection:sigma-oo-Q}%

In section~\reftwo{subsection:group-complete} we showed that the association $\F \mapsto |\Cob^\F(-)| = |\F(\Q(-\qshift{1}))|$ identifies with the suspension functor in $\Funadd(\Catp,\Sps)$.
To iterate this observation, note that for each $n \geq 1$, we have an $n$-fold simplicial object in $\Catp$ given by
\[
\Q^{(n)}(\C,\QF)\colon(\Del\op)^n \lrar \Catp \quad\quad ([m_1],...,[m_n]) \longmapsto \Q_{m_1}\Q_{m_2}...\Q_{m_n}(\C,\QF).
\]
By Lemmas~\reftwo{lemma:poincare} and \reftwo{lemma:segal}, $\Q^{(n)}(\C,\QF)$ is an $n$-fold Segal object of $\Catp$, the $n$-fold iterated hermitian $\Q$-construction of $(\C,\QF)$. As such, it is an  $n$-fold $\infty$-category and presents an $(\infty,n)$-category, though we shall not attempt to make this precise. We simply set:

\begin{definition}
\label{definition:iteratedCob}%
For $\F\colon\Catp\rightarrow\Sps$ additive, we shall call the $n$-fold complete Segal space $\F \Qh^{(n)}(\C,\QF\qshift{n})$ the \defi{$\F$-based $n$-extended cobordism $\infty$-category} $\Cob^\F_n(\C,\QF)$ of $(\C,\QF)$. 
\end{definition}

The simplicial space $\Cob^\F_1(\C,\QF) = \F\Q(\C,\QF\qshift{1})$ really is the Segal space giving rise to the cobordism $\infty$-category $\Cob^\F(\C,\QF)$ and $\Cob^\F_0(\C,\QF) = \F(\C,\QF)$. Furthermore, there are canonical equivalences
\[
|\Cob_i^{|\Cob_j^\F|}(\C,\QF)| \simeq |\Cob_{j+i}^\F(\C,\QF)|.
\]
The $n$-fold Segal space $\Poinc\Q^{(n)}(\C,\QF\qshift{n})$ models the $(\infty,n)$-category informally described as having Poincaré objects of $(\C,\QF\qshift{n})$ as objects, their cobordisms as morphisms, cobordisms between cobordisms as squares and so on up to degree $n$.

\begin{remark}
The analogous $n$-fold topological category $\Cob_d^n$ (note the unfortunate index switch) for cobordism categories of $d$-manifolds first appeared in \cite{BokstedtMadsen}, see also \cite{SomPri-Invertible}, ironically inspired by the iterated version of Quillen's original $\Q$-construction, and served to produce cobordism theoretic deloopings of $|\Cob_d|$. In particular, Bökstedt and Madsen showed that $|\Cob_d^n| \simeq \Omega^{\infty-n} \MTSO(d)$, extending the theorem of Galatius, Madsen, Tillmann and Weiss from the case $n=1$. They used this description to give an entirely cobordism theoretic model for the spectrum $\MTSO(d)$, which endows it with an interesting map to $\A(\BSO(d))$, studied extensively by Raptis and the ninth author in \cite{RaptisSteimle-CobcattoK, RaptisSteimle-SmoothDWW, RaptisSteimle-TopDWW}, where it was used to give a short proof of the Dwyer-Weiss-Williams index theorem \cite{DwyerWW}. We will take up the study of the evident refinements of this map in a sequel to the present paper.

\end{remark}

Now, denote by $\PSp$ the $\infty$-category of pre-spectra, that is the lax limit of the diagram
\[
\cdots \xrightarrow{\Omega} \Sps_* \xrightarrow{\Omega} \Sps_* \xrightarrow{\Omega} \Sps_*,
\]
consisting of sequences $(X_n)_{n \in \NN}$ of pointed spaces together with structure maps $X_n \rightarrow \Omega X_{n+1}$. There is a fully faithful inclusion $\Spa \subseteq \PSp$, which admits a left adjoint that we will refer to as spectrification. It does not affect the homotopy groups. Furthermore, the evaluation functors $ev_n \colon \PSp \rightarrow \Sps_*$ commute with both limits and colimits, and restrict to the functors $\Omega^{\infty-n} \colon \Spa \rightarrow \Sps_*$ which still preserve limits, but only filtered colimits. Similarly, for any pointed $\infty$-category $\E$ which admits finite limits, one can define $\PSp(\E)$ as the analogous lax limit and obtains the spectrum objects $\Sp(\E)$ of $\E$ as a full subcategory of $\PSp(\E)$.

\begin{definition}
\label{definition:cob-pre-spectrum}%
Let $\F\colon \Catp \rightarrow \Sps$ be an additive functor. We denote by 
\[
\CCob^{\F} = [\Cob^\F_0(-),|\Cob^{\F}_1(-)|,|\Cob^{\F}_2(-)|,\ldots]
\]
the corresponding functor from $\Catp$ to pre-spectra with the structure maps determined by the square of Theorem~\reftwo{theorem:suspension} applied to the functors $|\Cob_i^\F|$. 
\end{definition}

\begin{remark}
\label{remark:suspension}%
By definition, the functor $|\Cob_{n}^\F(-)|$ is the colimit of an $n$-fold simplicial space.
As a result, there is a natural $\Sig_n$-action $|\Cob_{n}^\F(-)|$ induced by permuting the entries of this $n$-fold simplicial object, a structure which can be used to promote $\CCob^\F(\C,\QF)$ to a symmetric pre-spectrum object. Classically, the non-hermitian analogue of this observation was used in order to obtain the lax symmetric monoidal structure on algebraic $\K$-theory. We shall not use this approach, and argue instead by universal properties to construct multiplicative structures in \paperfour.
\end{remark}

By Theorem~\reftwo{theorem:suspension}, $|\Cob_n^\F|$ is a model for the $n$-fold suspension of $\F$ in $\Funadd(\Catp,\Sps)$. Considering $\CCob^\F$ as an object in $\Funadd(\Catp,\PSp) \simeq \PSp(\Funadd(\Catp,\Sps))$, i.e.\ as a pre-spectrum object in $\Funadd(\Catp,\Sps)$, it is thus the suspension pre-spectrum of $\F$.
Since the $0$-th object in the pre-spectrum $\CCob(\C,\QF)$ is $\F(\C,\QF)$ itself, we obtain a natural map 
\[
\F(\C,\QF) \lrar \Om^{\infty}{\CCob}^{\F}(\C,\QF),
\]
where the right hand side refers to the $0$-th space of the spectrification of $\CCob^\F(\C,\QF)$.

\begin{proposition}
\label{proposition:positive-om}%
Let $\F$ be an additive functor $\Catp \rightarrow \Sps$ and $(\C,\QF) \in \Catp$. Then:
\begin{enumerate}
\item
\label{item:CobFsplitadd}%
The functor $\CCob^\F \colon\Catp \rightarrow \PSp$ is again additive and takes values in positive $\Om$-spectra, i.e.\ the structure map $|\Cob^{\F}_n(\C,\QF)| \rightarrow \Om |\Cob^{\F}_{n+1}(\C,\QF)|$ is an equivalence for every $n \geq 1$. 
\item
\label{item:Fgroup-like}%
If $\F$ is group-like, then $\CCob^\F(\C,\QF)$ is in fact an ($\Om$-)spectrum, and $\CCob^\F$ is additive when considered as a functor 
$\CCob^\F \colon\Catp \rightarrow \Spa$.
\item
\label{item:Fspectrification}%
The natural map $\CCob^\F(\C,\QF) \rightarrow \CCob^{\F^\grp}(\C,\QF)$ exhibits the right hand side as the spectrification of the left. 
\end{enumerate}
In particular, for $n\geq 1$, we obtain equivalences 
\[
\F^\grp(\C,\QF) \simeq \Omega^\infty\CCob^\F(\C,\QF) \quad\text{and}\quad |\Cob_n^\F(\C,\QF)| \simeq \Omega^{\infty-n}\CCob^\F(\C,\QF).
\]
\end{proposition}
\begin{proof}
By Theorem~\reftwo{theorem:suspension} 
the $\infty$-category $\Funadd(\Catp,\Sps)$ satisfies that its suspension functor is fully faithful on group-like objects, and hence $\Omega\Sigma$ is a group completion by Lemma~\reftwo{proposition:grpcomp}. Said lemma then has the following immediate extension. Namely, for any $\E$ with these properties, we have
\begin{enumerate}
\item The functor $\E \to \PSp(\E)$ sending $E$ to its suspension pre-spectrum $\Sigma^\infty(E)$ takes values in positive $\Omega$-spectra.
\item If $E \in \E_\grp$ is group-like, then $\Sigma^\infty E$ is an $\Omega$-spectrum.
\item The natural map $\Sigma^\infty E \to \Sigma^\infty(\Omega \Sigma E)$ exhibits the target as the spectrification of the source.
\end{enumerate}
In particular, for $n\geq 1$, there are equivalences
\[
E^\grp\simeq \Omega\Sigma E \simeq \Omega^\infty \Sigma^\infty(E) \quad \text{ and } \quad \Sigma^n(E) \simeq \Omega^{\infty-n} \Sigma^\infty(E),
\]
proving the proposition.
\end{proof}

\begin{corollary}
\label{corollary:Cob=Susp}%
For a group-like additive functor $\F\colon \Catp \rightarrow \Sps$, the functor $\CCob^{\F} \colon \Catp \rightarrow \Spa$ is the initial additive functor under $\SS[\F]$, the pointwise suspension spectrum of $\F$. In other words, 
\[
\CCob \colon \Funadd(\Catp, \Grp_{\Einf}) \longrightarrow \Funadd(\Catp,\Spa)
\]
is left adjoint to the forgetful functor, i.e.\ composition with $\Omega^\infty$. Also,
\[
\CCob \circ (-)^\grp  \colon \Funadd(\Catp, \Sps) \longrightarrow \Funadd(\Catp,\Spa)
\]
is left adjoint to the forgetful functor.
\end{corollary}

Part~\reftwoitem{item:Fspectrification} of Proposition~\reftwo{proposition:positive-om}, applied to $\F\in \Funadd(\Catp,\Sps)$, identifies the non-negative homotopy groups of $\CCob^\F(\C,\QF)$ with those of $\F^\grp(\C,\QF)$. While these are generally very difficult to understand, we can determine the negative homotopy groups of the spectrum $\CCob^\F(\C,\QF)$ much more easily:

\begin{proposition}
\label{proposition:pinegCobF}%
For every additive $\F \colon \Catp \rightarrow \Sps$, Poincaré \(\infty\)-category $(\C,\QF)$, $n \geq 1$ and $0 \leq k < n$, the iterated bonding maps of the pre-spectrum $\CCob^\F(\C,\QF)$ induce isomorphisms
\[
\pi_k|\Cob^\F_n(\C,\QF)| \cong \pi_0|\Cob^\F(\C,\QF\qshift{n-k-1})| \quad \text{and} \quad \pi_{-n}\CCob^\F(\C,\QF) \cong \pi_0|\Cob^\F(\C,\QF\qshift{n-1})|.
\]
\end{proposition}

In other words, $\pi_k|\Cob_n^\F(\C,\QF)|$ for $k < n$ is just the $\F$-based cobordism group of $(\C,\QF\qshift{n-k})$ and similarly for the negative homotopy groups of $\CCob^\F$. 

\begin{proof}
Part~\reftwoitem{item:CobFsplitadd} of Proposition~\reftwo{proposition:positive-om} reduces the claim about the left hand side to the case $k=0$. By realising the $n$-fold simplicial object $\Cob_n^\F$ iteratively, this case follows from Corollaries~\reftwo{corollary:components-of-cob-II} and \reftwo{corollary:pi0CorCobMet} by induction on $n$. The statement for the right hand side is now immediate from Proposition~\reftwo{proposition:positive-om} \reftwoitem{item:Fspectrification} and Corollary~\reftwo{corollary:Cobgrp}.
\end{proof}

In particular, we find that $\CCob^\F(\C,\QF)$ is connective whenever $(\C,\QF)$ is metabolic by Proposition~\reftwo{proposition:pinegCobF}, Corollary~\reftwo{corollary:metabolic} and  Corollary~\reftwo{corollary:pi0CorCobMet}. In fact we have:

\begin{corollary}
\label{corollary:CCobcon}%
The functor
\[
\CCob \colon \Funadd(\Catp, \Grp_{\Einf}) \longrightarrow \Funadd(\Catp,\Spa)
\]
is fully faithful and its essential image consists precisely of the functors whose values on all metabolic Poincaré \(\infty\)-categories $(\C, \QF)$ are connective.
\end{corollary}

By Corollary~\reftwo{corollary:metabolic}, the essential image of $\CCob$ is equivalently described as those functors whose value on $\Hyp(\E)$ for all small stable $\infty$-categories $\E$ are connective.

\begin{proof}
That $\CCob$ is fully faithful follows from Corollary~\reftwo{corollary:Cob=Susp}, since the unit $\F \Rightarrow \Omega^\infty\CCob^\F$ is an equivalence by Proposition~\reftwo{proposition:positive-om} if $\F$ is group-like. Regarding the essential image, one inclusion was discussed before the statement of the corollary. For the converse consider an $\F\in \Fun^\add(\Catp,\Spa)$ whose values on metabolic categories are connective. By the triangle identities, the counit $\CCob^{\Omega^\infty \F}\to \F$ of the adjunction is an equivalence after applying $\Omega^\infty$, and therefore on non-negative homotopy groups. Applying this counit transformation to the metabolic fibre sequence 
\[
(\C, \QF) \to \Met(\C, \QF\qshift 1) \to (\C, \QF\qshift 1),
\]
and using that $\Met(\C,\QF\qshift{1})$ is indeed metabolic, see Example~\reftwo{examples:hyp}~\reftwoitem{item:metCQmeta}, we conclude inductively on $i$ that it is an equivalence on $\pi_{-i}$ for all $i\geq 0$.
\end{proof}

\begin{remark}
\label{remark:K-theorysuspension}%
Completely analogous definitions and arguments work in the non-hermitian set-up to give the $n$-fold Segal spaces $\Span_n^\F(\C)$ and (pre-)spectra $\mathbb{S}\mathrm{pan}^\F(\C)$, with the $\K$-theory functor $\Catx \rightarrow \Spa$ being the (pointwise) spectrification of $\mathbb{S}\mathrm{pan}^{\core}$ or equivalently $\mathbb{S}\mathrm{pan}^{\core^{\grp}}$. As a consequence of Proposition~\reftwo{proposition:spanconn}, one here finds that $\mathbb{S}\mathrm{pan}^\F(\C)$ is always a connective (pre-)spectrum. The analogue of the above corollary is the statement that
\[
\mathbb{S}\mathrm{pan} \colon \Funadd(\Catx, \Grp_{\Einf}) \longrightarrow \Funadd(\Catx,\Spa)
\]
is fully faithful with essential image the functors taking values in connective spectra. %
In particular, the non-connectivity of the iterated $\Q$-construction is an entirely hermitian phenomenon.
\end{remark}

Let us also record the relationship between the $\Q$-construction and suspension in $\Funadd(\Catp,\Spa)$. To this end, consider the diagram from Observation~\reftwo{observation:for-annoying-sequence} 
\[
\begin{tikzcd}
\const (\C,\QF\qshift{-1}) \ar[r]\ar[d] & \hMet(\C,\QF) \ar[r] \ar[d] & \Q(\C,\QF) \ar[d,equal] \\
\const(\Met(\C,\QF)) \ar[r] & \dec(\Q(\C,\QF)) \ar[r] & \Q(\C,\QF)
\end{tikzcd}
\]
consisting of split Poincaré-Verdier sequences in each simplicial degree. Applying an additive $\F\colon \Catp \rightarrow \Spa$, one obtains levelwise fibre sequences of simplicial spectra. As these are also cofibre sequences by stability, it follows that they remain fibre sequences after realisation. Moreover, as the simplicial objects in the top right corners are split by Lemma~\reftwo{lemma:split} over $(\C,\QF)$ and $0$, respectively, we obtain the following commutative diagram of horizontal fibre sequences.
\[
\begin{tikzcd}
	\F(\C,\QF\qshift{-1}) \ar[r] \ar[d] & 0 \ar[r] \ar[d] & {|\F\Q(\C,\QF)|} \ar[d,equal] \\
	\F(\Met(\C,\QF)) \ar[r] & \F(\C,\QF) \ar[r] & {|\F\Q(\C,\QF)|}
\end{tikzcd}
\]

\begin{corollary}
\label{corollary:Q=suspinsp}%
\label{remark:left-square-met}%
The endofunctor $\F \mapsto |\F\Q(-\qshift{1})|$ on $\Funadd(\Catp,\Spa)$ is the internal suspension functor, i.e.\ postcomposition with the suspension functor in $\Spa$, and the lower horizontal fibre sequence above is the rotation (to the right) of the metabolic fibre sequence.
\end{corollary}

Note that the geometric realisation $|\F\Q(-\qshift{1})|$ appearing in Corollary~\reftwo{corollary:Q=suspinsp} %
can be taken in either $\Fun(\Catp,\Spa)$ or $\Funadd(\Catp,\Spa)$ since
the latter is closed under colimits in the former by the stability of $\Spa$.

Finally, we study the effect of shifting the Poincaré structure on the $\Ct$-equivariant spectrum $\F(\Hyp(\C))$ acted on by the duality of $\QF$. This will ultimately lead to our generalisation of Karoubi's periodicity theorem in Corollary~\reftwo{corollary:Karoubiper} below. To this end, recall that the composite $\Catp \rightarrow \Catx \rightarrow \Catp$ of the forgetful and hyperbolic functor refines to a  
functor $\HypG \colon \Catp \rightarrow \Fun(\BC,\Catp)$ via the action of the duality, see Remark~\refone{remark:hyp-equivariant}.

\begin{definition}
\label{definition:Fhyp}%
Given a functor $\F \colon \Catp \rightarrow \E$ define the \defi{hyperbolisation} $\F^\hyp \colon \Catp \rightarrow \Fun(\BC,\E)$ of $\F$ as $\F \circ \HypG$.
\end{definition}

\begin{proposition}[Naive Karoubi periodicity]
\label{proposition:naive-karoubi}%
There is a canonical equivalence of $\Ct$-spectra 
\[
\F^{\hyp}(\C,\QF\qshift{-1}) \simeq \SS^{\sigma-1} \ssmash \F^{\hyp}(\C,\QF),
\]
natural in the Poincaré \(\infty\)-category $(\C,\QF)$ and the additive functor $\F\colon \Catp\rightarrow \Spa$. Furthermore, under this equivalence, the boundary map 
\[
\F^\hyp(\C,\QF) \rightarrow \SS^1 \ssmash \F^\hyp(\C,\QF\qshift{-1})
\]
of the metabolic fibre sequence is induced by the inclusion $\sph^0 \rightarrow \sph^\sigma$ as the fixed points.
\end{proposition}

Here, $\SS^\sigma$ denotes the $\Ct$-spectrum equivalently described as the suspension spectrum of $\sph^{\sigma}$, the $1$-sphere with complex conjugation action, or the functor 
\[
\BC = \BO(1) \rightarrow \BO \xrightarrow{\mathrm{J}} \Pic(\SS) \subseteq \Spa.
\]

\begin{proof}
We recall that under the equivalence $\F(\Met(\C,\QF)) \simeq \F(\Hyp(\C))$ induced by $\ilag\colon\Hyp(\C) \rightarrow \Met(\C,\QF)$, the map $\met \colon \Met(\C,\QF) \rightarrow (\C,\QF)$ identifies with $\hyp\colon \Hyp(\C) \rightarrow (\C,\QF)$. Using the metabolic Poincaré-Verdier sequence, we may therefore identify $\SS^1 \ssmash \F^\hyp(\C,\QF\qshift{-1})$ with the cofibre of the map 
\[
\F^{\hyp}(\hyp)\colon\F^\hyp(\Hyp(\C)) \lrar \F^\hyp(\C,\QF).
\]
Now, there is a natural equivalence
\[
\HypG(\Hyp(\C)) \simeq \Hyp(\C \times \C\op) \simeq \Hyp(\C) \otimes \Ct
\]
which translates the action of $\Dual_\QF$ on the left into the flip action on the right, see Remark~\refone{remark:hyp-of-hyp}. We may then identify the map $\F\hyp$ with the map
\[
\F(\Hyp(\C)) \ssmash \Ct \lrar \F(\Hyp(\C))
\]
obtained from the map $\Ct \rightarrow \ast$ of $\Ct$-spaces, whose cofibre is $\sph^{\sigma}$. We therefore obtain a natural equivalence
\[
\SS^1 \ssmash \F^\hyp(\C,\QF\qshift{-1}) \simeq \SS^\sigma \ssmash \F^\hyp(\C,\QF)
\]
which is the claim.
\end{proof}

\subsection{Bordism invariant functors}
\label{subsection:bordism}%

In the next two subsections, we introduce the notion of a bordism invariant functor out of $\Catp$, the main examples being various flavours of $\L$-theory. We then show that each additive functor $\F \colon \Catp \rightarrow \Sps$ admits an initial bordism invariant functor $\F^\bord$ equipped with a map $\F \rightarrow \F^\bord$, the bordification of $\F$, and that any group-like $\F$ can then be described in terms of this bordification and the hyperbolisation $\F^\hyp = \F \circ \Hyp$ from the previous section. This yields a version of the first part of our \reftwo{theorem:main-intro-two} for arbitrary additive{} functors in Corollary~\reftwo{corollary:tate-square}; it will be specialised to $\F=\Poinc$ in section~\reftwo{section:GW}.

To get started recall the notion of a cobordism between Poincaré functors from Definition~\reftwo{definition:bordism-functors}: It is a Poincaré functor $(\C,\QF) \rightarrow \Q_1(\Ctwo,\QFtwo)$ projecting correctly to the endpoints of $\Q_1$. 

\begin{definition}
\label{definition:bordism-invariant}%
A Poincaré functor $(F,\eta) \colon (\C,\QF) \rightarrow (\Ctwo,\QFtwo)$ is called a \defi{bordism equivalence} if there exists a Poincaré functor $(G,\vartheta) \colon (\Ctwo,\QFtwo) \rightarrow (\C,\QF)$ such that the composites $(F,\eta) \circ (G,\theta)$ and $(G,\theta) \circ (F,\eta)$ are cobordant to the respective identities. 
\end{definition}

\begin{example}
\label{example:isotropic}%
Let $(\C,\QF)$ be a Poincaré $\infty$-category and $\Lag \subseteq \C$ an isotropic subcategory (see Definition~\reftwo{definition:isotropic}). Then the inclusion $\Hlgy(\Lag) \subseteq (\C,\QF)$ of the homology $\infty$-category is a bordism equivalence. This follows directly from Construction~\reftwo{construction:homology}.
\end{example}

\begin{definition}
Given an $\infty$-category with finite products $\E$, we say that an additive functor $\F\colon \Catp \rightarrow \E$ is \defi{bordism invariant} if it sends bordism equivalences of Poincaré \(\infty\)-categories to equivalences in $\E$. We shall denote by $\Fun^\bord(\Catp,\E)$ the full subcategory of $\Funadd(\Catp,\E)$ spanned by the bordism invariant functors.
\end{definition}

\begin{example}
\label{example:bord-L}%
The formation of $\L$-spaces provides a bordism invariant functor $\Lspace \colon \Catp \rightarrow \Sps$. This was established by Lurie in \cite{Lurie-L-theory} based on fundamental work of Ranicki, see e.g.\ \cite{Ranickiblue}. We will discuss this functor and its spectral refinement in detail in section~\reftwo{subsection:L+tate} below. %
\end{example}

\begin{remark}
\label{remark:cob-bordism-invariance}%
If $\F \colon \Catp \rightarrow \Sps$ is additive and bordism invariant, then so is $|\Cob^\F|$. This follows straight from the definitions, as a cobordism of Poincaré functors $(\C,\QF) \rightarrow \Q_1(\Ctwo,\QFtwo)$, induces another such 
\[
\Q_n(\C,\QF) \rightarrow \Q_n(\Q_1(\Ctwo,\QFtwo)) \cong \Q_1(\Q_n(\Ctwo,\QFtwo)).
\]
Hence, a bordism equivalence $(\C,\QF) \rightarrow (\Ctwo,\QFtwo)$ gives an equivalence of simplicial objects $\F\Q(\C,\QF) \rightarrow \F\Q(\Ctwo,\QFtwo)$, and thus an equivalence on realisations.
\end{remark}

By Example~\reftwo{example:isotropic}, a bordism invariant functor vanishes on all metabolic Poincaré \(\infty\)-categories, i.e.\ those that admit a Lagrangian. For group-like additive functors the converse holds as well:

\begin{lemma}
\label{lemma:criterion}%
Let $\F\colon\Catp \rightarrow \E$ be a group-like additive functor. Then the following are equivalent:
\begin{enumerate}
\item 
\label{item:Fbordinv}%
$\F$ is bordism invariant.
\item 
\label{item:degmap}%
$\F$ takes $\cyl\colon (\C,\QF) \rightarrow \Q_1(\C,\QF)$ to an equivalence for every Poincaré \(\infty\)-category $(\C,\QF)$.
\item 
\label{item:Fmet}%
$\F$ vanishes on all metabolic Poincaré \(\infty\)-categories.
\item 
\label{item:FMetCQ}%
$\F(\Met(\C,\QF)) \simeq \ast$ for any Poincaré \(\infty\)-category $(\C,\QF)$.
\item 
\label{item:FHyp}%
$\F(\Hyp(\C)) \simeq \ast$ for any stable \(\infty\)-category $\C$.
\end{enumerate}
\end{lemma}

\begin{proof}
The functors in \reftwoitem{item:degmap} are bordism equivalences (essentially by definition) so \reftwoitem{item:Fbordinv} $\Rightarrow$ \reftwoitem{item:degmap} and it follows immediately from Corollary~\reftwo{corollary:decomp-Q1} that \reftwoitem{item:degmap} $\Rightarrow$ \reftwoitem{item:FMetCQ}. By Example~\reftwo{example:isotropic}, all metabolic categories are bordism equivalent to $0$, so \reftwoitem{item:Fbordinv} $\Rightarrow$ \reftwoitem{item:Fmet} and since $\Met(\C,\QF)$ really is metabolic, we have \reftwoitem{item:Fmet} $\Rightarrow$ \reftwoitem{item:FMetCQ}. To obtain \reftwoitem{item:FMetCQ} $\Rightarrow$ \reftwoitem{item:FHyp} observe that $\F\Hyp(\C)$ is a retract of $\F\Hyp(\C \times \C\op)$, which by Corollary~\reftwo{corollary:met-hyp} is equivalent to $\F\Met(\Hyp(\C)) \simeq \ast$. Finally, by Proposition~\reftwo{proposition:cobordism}, if $\F$ vanishes on hyperbolics, then cobordant Poincaré functors induce homotopic maps after applying $\F$, and so $\F$ is bordism invariant giving \reftwoitem{item:FHyp} $\Rightarrow$ \reftwoitem{item:Fbordinv}.
\end{proof}

To discuss another important example, recall from Definition~\reftwo{definition:Fhyp} the hyperbolisation $\F^\hyp(\C,\QF) = \F(\HypG(\C))$ taking values in the $\infty$-category $\Fun(\BC,\E)$ via the action of the duality $\Dual_\QF$ on $\Hyp(\C)$. 

\begin{example}
\label{example:tate}%
Given an additive functor $\F\colon\Catp \rightarrow \E$ with $\E$ stable and with sufficient (co)limits to form the Tate construction $(-)^\tC \colon \Fun(\BC,\E) \to \E$,  we find a functor $(\F^{\hyp})^{\tC}\colon\Catp \rightarrow \E$ which is bordism invariant; to see this, we invoke the natural equivalence 
\[
\HypG(\Hyp(\C)) \lrar \Hyp(\C) \otimes \Ct
\]
from~\refone{remark:hyp-of-hyp}
which shows that $\F^\hyp(\Hyp(\C))$ is an induced $\Ct$-object. It then follows that for every stable \(\infty\)-category $\C$, we have $(\F^{\hyp})^{\tC}\Hyp(\C) \simeq 0$, 
since the Tate construction generally vanishes on induced $\Ct$-objects.
\end{example}

\begin{proposition}
\label{proposition:bordism-invariant-sig}%
Suppose that $\F\colon \Catp \rightarrow \E$ is a bordism invariant additive functor. Then the natural map
$
\Omega \F(\C,\QF) \to \F(\C,\QF\qshift{-1})
$
arising from the metabolic Poincaré-Verdier sequence is an equivalence. In particular, $\F$ is automatically group-like.
If $\E$ is stable, an additive functor $\F\colon\Catp \rightarrow \E$ is bordism invariant if and only if this map is an equivalence for all Poincaré \(\infty\)-categories $(\C,\QF)$. 
\end{proposition} 
In particular, we find $\pi_i\F(\C,\QF) = \pi_0\F(\C,\QF\qshift{-i})$ for every space or spectrum valued bordism invariant functor. 
\begin{proof}
By Lemma~\reftwo{lemma:criterion}, $\F$ is bordism invariant if and only if $\F(\Met(\C,\QF)) \simeq \ast$ for all Poincaré \(\infty\)-categories $(\C, \QF)$, from which we obtain a fibre sequence
\[
\F(\C,\QF\qshift{-1}) \longrightarrow \ast \longrightarrow \F(\C,\QF),
\]
giving the first claim. If $\E$ is stable, then the map in question being an equivalence is equivalent to $\F\Met(\C,\QF) \simeq \ast$ giving the second claim.
\end{proof}

In particular, bordism invariant functors can be delooped simply by shifting the Poincaré structure, i.e.\ by considering
\[
\big[\F(\C,\QF), \F(\C,\QF\qshift{1}), \F(\C,\QF\qshift{2}), \ldots \big]
\]
with the structure maps provided by Proposition~\reftwo{proposition:bordism-invariant-sig}, and furthermore, by Corollary~\reftwo{corollary:inversion}, the inversion map on $\F(\C,\QF)$ is induced by the Poincaré functor $(\id_\C,-\id_\QF)$.
Next, we show that this delooping agrees with that from the previous section. In fact, we have as the main result of this subsection: %

\begin{theorem}
\label{theorem:lift}%
\label{theorem:unique-lift}%
The forgetful functor 
\[
\Fun^\bord(\Catp,\Spa) \longrightarrow \Fun^\bord(\Catp,\Sps),
\]
i.e.\ postcomposition with $\Omega^\infty$, is an equivalence with inverse $\F \mapsto \CCob^{\F}$. This equivalence restricts to Verdier- and Karoubi-localising functors.
In particular, any additive (or Verdier-localising) and bordism invariant functor $\F\colon \Catp \rightarrow \Sps$ admits an essentially unique lift to another such functor $\Catp \rightarrow \Spa$.
\end{theorem}

The same is not true for arbitrary group-like additive $\F \colon \Catp \to \Sps$ as the examples 
$
(\C,\QF) \to \K(\C^\natural), \KK(\C),
$
which have equivalent infinite loop spaces, show.

\begin{proof}[Proof of Theorem~\reftwo{theorem:lift}] 
By Proposition~\reftwo{proposition:bordism-invariant-sig} additive bordism invariant functors are group-like and so we may view $\Fun^\bord(\Catp,\Sps)$ as a full subcategory of $\Funadd(\Catp, \Grp_{\Einf})$. We now observe that the adjunction
\[
\CCob \colon \Funadd(\Catp, \Grp_{\Einf}) \adj \Fun^{\add}(\Catp,\Spa) \cocolon \Omega^\infty
\]
of Corollary \reftwo{corollary:Cob=Susp} restricts to bordism invariant functors: Indeed, this is clear for $\Omega^\infty$, and for the left adjoint it follows by induction from Remark~\reftwo{remark:cob-bordism-invariance}. The restricted adjunction is then an equivalence by Corollary~\reftwo{corollary:CCobcon} since any $\F \in \Fun^{\bord}(\Catp,\Spa)$ vanishes on metabolic Poincaré $\infty$-categories, and so in particular sends them to connective spectra. It remains to observe that $\Cob^\F$ is Verdier- or Karoubi-localising whenever $\F\colon \Catp \to \Sps$ is Verdier- or Karoubi-localising and bordism invariant, which follows readily from the delooping explained before the proof together with the fact that $(-)\qshift{1}$ preserves Poincar\'e-Verdier and Poincar\'e-Karoubi sequences.
\end{proof}

\subsection{The bordification of an additive functor}
\label{subsection:bordification}%

In this subsection, we establish the following theorem and deduce a general version of our \reftwo{theorem:main-intro-two} in Corollary~\reftwo{corollary:tate-square}.

\begin{theorem}
\label{theorem:localisation}%
The inclusions 
\[
\Fun^\bord(\Catp,\Spa) \subseteq \Funadd(\Catp,\Spa) \quad \text{and} \quad \Fun^{\bord}(\Catp,\Sps) \subseteq \Funadd(\Catp,\Sps)
\]
of the bordism invariant into all additive functors admit left and right adjoints.
\end{theorem}

\begin{definition}
\label{definition:bordify}%
We will refer to these left adjoint functors as \defi{bordification} and denote their values on an additive functor $\F \colon \Catp \rightarrow \Spa$ or $\F\colon \Catp \rightarrow \Sps$ by $\F^\bord$. The right adjoint we shall call \defi{cobordification} and denote by $(-)^\rbord$.
\end{definition}

The existence of the left adjoints, Theorem~\reftwo{theorem:lift} and Corollary~\reftwo{corollary:Cob=Susp} may be summarised by the following square of forgetful functors and their dotted left adjoints, whose left hand vertical arrows are inverse equivalences, as follows:
\[
\begin{tikzcd}
[column sep=7ex]
\Fun^\bord(\Catp, \Spa) \ar[d, shift left=1ex] \ar[r,shift right=.5ex] & \Funadd(\Catp, \Spa) \ar[d,shift left=2ex] \ar[l,dotted,shift left=.5ex,bend right=10,"{(-)^\bord}"'] \\
\Fun^\bord(\Catp, \Sps) \ar[r,shift right=.5ex] \ar[u,dotted,bend left=10,"\CCob"] & \Funadd(\Catp, \Sps) \ar[u, dotted,shift right=1ex,bend left=10,"{\CCob\circ(-)^\grp}"] \ar[l,dotted,shift left=.5ex,bend right=10,"{(-)^\bord}"']
\end{tikzcd}
\]
Thus the existence of the upper horizontal adjoint implies the existence of the lower one, or more precisely for $\F \colon \Catp \rightarrow \Sps$ additive, we have
\[
\F^\bord \simeq \Omega^\infty((\CCob^{\F^\grp})^\bord).
\]
Similarly, we find 
\[
\F^\rbord \simeq \Omega^\infty((\CCob^{\F^\times})^\rbord),
\]
where $\F^\times$ denotes the pointwise units of $\F$: For $\G$ bordism invariant one computes
\begin{align*}
\Nat(\G,\Omega^\infty((\CCob^{\F^\times})^\rbord)) &\simeq \Nat(\CCob^\G,(\CCob^{\F^\times})^\rbord) \\
                                                  &\simeq \Nat(\CCob^\G,\CCob^{\F^\times}) \\
                                                  &\simeq \Nat(\G,\F^\times) \\
                                                  &\simeq \Nat(\G,\F)
\end{align*}
where the first and third identities follow from Corollary~\reftwo{corollary:Cob=Susp}. In particular, the (co)bordifications at the spectrum level determine those at the space level, so we will restrict attention to the case of functors taking values in spectra in this section. 

\begin{remark}
\label{remark:bordvsomegainf}%
Using Theorem~\reftwo{theorem:unique-lift} it is also easy to see that $\F^\rbord \simeq \CCob^{(\Omega^\infty \F)^\rbord}$ for additive $\F\colon \Catp \rightarrow \Spa$, and so $\F^{\rbord}$ depends only on $\Om^{\infty}\F$. We explicitly warn the reader that $\F^\bord$ cannot in general be reconstructed from $\Omega^\infty\F$ in similar manner. A counterexample is given by the Karoubi-Grothendieck-Witt functor, discussed in detail in \paperfour, see Remark~\reftwo{remark:compare-bordification}.
\end{remark}

We give three distinct formulae for the spectral bordification functor in Proposition~\reftwo{proposition:phi-bord}, Corollary~\reftwo{corollary:rho-bord} and Corollary~\reftwo{corollary:Q-oo-bord}, and it really is the comparison between these that is most relevant for our work. While this comparison can be established by direct calculations, that route does not lead to shorter arguments and the present framework allows for a more conceptual interpretation; see Proposition~\reftwo{proposition:bordrecoll}. Along the way we also establish the cobordification in Proposition \reftwo{proposition:phi-bord}. We begin with the following lemma.

\begin{lemma}
\label{lemma:hyptobord}%
If $\F \colon \Catp \rightarrow \E$ is arbitrary and $\G \colon \Catp \rightarrow \E$ is bordism-invariant, then the spaces
$\Nat(\F^\hyp,\G), \Nat(\F^\hyp_\hC,\G), \Nat(\G,\F^\hyp)$ and $\Nat(\G,(\F^\hyp)^\hC)$ are contractible (assuming $\E$ admits sufficient (co)limits to form the homotopy orbits and fixed points appearing).
\end{lemma}
\begin{proof}
Since $\Hyp\colon \Catx \rightarrow \Catp$ is both left and right adjoint to the forgetful functor by Corollary~\refone{corollary:hyp-is-adjoint}, the composite $\Catp \to \Catx \to \Catp$ is both left and right adjoint to itself. Hence the association $\F \mapsto \F^{\hyp}$ is both left and right adjoint to itself as well. Since $\G^{\hyp} \simeq *$ for any bordism invariant functor, it follows that the mapping space from $\F^{\hyp}$ to and from any bordism invariant functor is trivial. 
The computations
\[
\Nat(\F^\hyp_\hC,\G) \simeq \Nat(\F^\hyp,\G)^\hC \simeq * \quad \text{and} \quad \Nat(\G,(\F^\hyp)^\hC) \simeq \Nat(\G,\F^\hyp)^\hC \simeq *
\]
give the remaining claims.
\end{proof}

Next, recall from Corollary~\refone{corollary:hyperbolic-mackey}
that the hyperbolic and forgetful maps
refine to $\Ct$-equivariant maps 
\[
\HypG(\C,\QF) \lrar (\C,\QF) \lrar \HypG(\C,\QF),
\]
where $(\C,\QF)$ is considered with the trivial $\Ct$-action. It then follows that the induced natural maps
$\F^{\hyp} \to \F \to \F^{\hyp} $
refine to maps of the form
\[
\F^{\hyp}_{\hC} \lrar \F \lrar (\F^{\hyp})^{\hC} .
\]
As part of \refone{corollary:hyperbolic-mackey}, we also showed that the composite of these two maps coincides with the norm map $\F^{\hyp}_{\hC} \rightarrow (\F^{\hyp})^{\hC}$ associated to the $\Ct$-action on $\F^{\hyp}$; we showed this at the level of Poincaré $\infty$-categories and it persists since $\F$ preserves finite products.
From Lemma~\reftwo{lemma:hyptobord}, we find $\Nat(\F^{\hyp}_{\hC},\G) \simeq \ast$
if $\G$ is bordism invariant. In particular, assuming the existence of a bordification, there is a sequence
\[
\F^{\hyp}_{\hC} \lrar \F \lrar \F^\bord
\]
whose composite admits an essentially unique null-homotopy. There is a universal way to produce such a sequence:

\begin{proposition}
\label{proposition:phi-bord}%
Consider the functor $\Phi\colon \Funadd(\Catp,\Spa) \rightarrow \Funadd(\Catp,\Spa)$ given by the formula
\[
\Phi\F = \cof(\F^{\hyp}_{\hC} \lrar \F).
\]
Then, the canonical transformations $\F \Rightarrow \Phi \F$ exhibit $\Phi$ as a bordification. Dually, $\fib(\F \rightarrow (\F^\hyp)^\hC)$ is a cobordification of $\F$.
\end{proposition}

\begin{proof}
Let $\F$ be an additive functor. We first verify that $\Phi\F$ is bordism invariant. By Lemma~\reftwo{lemma:criterion}, we need to check that $\Phi\F(\Hyp(\C)) \simeq 0$, or, equivalently, that the canonical transformation
\[
\F^{\hyp}_{\hC}(\Hyp(\C)) \lrar \F(\Hyp(\C))
\]
is an equivalence, which follows, similarly as in Example~\reftwo{example:tate}, from the observation that $\HypG(\Hyp(\C))$ is an induced $\Ct$-object.
Now, for any bordism invariant functor $\G$, we have a fibre sequence
\[
\Nat(\Phi\F,\G) \lrar \Nat(\F,\G) \lrar \Nat(\F^{\Hyp}_{\hC},\G)
\]
and $\Nat(\F^{\Hyp}_{\hC},\G) \simeq *$ by Lemma~\reftwo{lemma:hyptobord}. Consequently, the first map is an equivalence, showing that $\F \to \Phi\F$ is a bordification.
The argument for the cobordification is entirely dual.
\end{proof}

\begin{examples}
\label{lemma:vanishing-property}%
\label{lemma:bord-tate}%
\label{example:underlyingbord}%
Let $\F \colon \Catp\to \Spa$ and $\G \colon \Catx \to \Spa$ be additive. Using Example~\reftwo{example:tate} we find:
\begin{enumerate}
\item $(\F^{\hyp})^\bord \simeq 0 \simeq (\F^\hyp_\hC)^\bord$ and the natural map $(\F^\hyp)^\hC \Rightarrow (\F^\hyp)^\tC$ descends to an equivalence $((\F^\hyp)^\hC)^\bord \simeq (\F^\hyp)^\tC$.
\item Dually, $(\F^\hyp)^{\rbord} \simeq 0 \simeq ((\F^\hyp)^\hC)^{\rbord}$ and $(\F^\hyp_\hC)^{\rbord} \simeq \SS^{-1} \otimes (\F^\hyp)^\tC$.
\item On also generally has $(\G\circ \fgt)^\bord \simeq 0 \simeq (\G\circ\fgt)^{\rbord}$. In particular, the bordification of $\Core, \K,\KK, \mathrm{THH}, \mathrm{TC}$ and all similar invariants of Poincar\'e $\infty$-categories vanish: (Co)bordification is a genuinely hermitian concept hat has no classical counterpart.
\end{enumerate}
\end{examples}

We find an general form of our main result, the \emph{fundamental fibre sequence} and \emph{fundamental fibre square}: %
\begin{corollary}
\label{corollary:bord-seq-unique}%
\label{corollary:tate-square}%
For every additive functor $\F \colon \Catp \rightarrow \Spa$ and Poincaré \(\infty\)-category $(\C,\QF)$, there is a canonical fibre sequence 
\[
\F^{\hyp}_{\hC}(\C,\QF) \xrightarrow{\hyp} \F(\C,\QF) \xrightarrow{\bord} \F^\bord(\C,\QF).
\]
Moreover, applying bordification to the natural map $\F \to (\F^\hyp)^{\hC}$, it extends to a cartesian square
\[
\begin{tikzcd}
\F(\C,\QF) \ar[r] \ar[d] & \F^\bord(\C,\QF) \ar[d] \\
(\F^\hyp)^\hC(\C,\QF) \ar[r] & (\F^\hyp)^\tC(\C,\QF).
\end{tikzcd}
\]
\end{corollary}

\begin{proof}
The first part is a reformulation of Proposition~\reftwo{proposition:phi-bord}. For the second claim we use Example~\reftwo{lemma:bord-tate} to identify the bordification of $(\F^\hyp)^\hC$ with $(\F^\hyp)^\tC$ via the canonical map. It follows that both horizontal fibres identify with $\F^\hyp_\hC$ and furthermore, the induced map between them identifies with the identity: By definition a map into the lower fibre is determined by the map to $(\F^\hyp)^\hC$ together with the null-homotopy of the composite to $(\F^\hyp)^\tC$. And indeed, we have shown in Corollary~\refone{corollary:hyperbolic-mackey} that the composite of the natural maps 
\[
\F^\hyp_\hC \to \F \to (\F^\hyp)^\hC
\]
is the norm map of $\F^\hyp$ and there is only one null-homotopy of the further composite by another application of Lemma~\reftwo{lemma:bord-tate}.
\end{proof}

Corollary~\reftwo{corollary:tate-square} can also be recast in terms of stable recollements and their associated fracture squares, see section~\reftwo{section:appendix-split-verdier}, particularly Definition~\reftwo{definition:stable-recollement} and the discussion thereafter. 

\begin{proposition}
\label{proposition:bordrecoll}%
The diagram
\[
\begin{tikzcd}
[column sep=7ex]
\Fun^\bord(\Catp,\Spa) \ar[r] & 
\Fun^\add(\Catp,\Spa) 
\ar[r] 
\ar[l,bend left=25,shift left=1.5ex,start anchor=west,end anchor=east, "{\rbord}"] 
\ar[l,bend right=25,shift right=1.5ex,start anchor=west,end anchor=east,"{\bord}"'] 
\ar[l,phantom,shift left=1.2ex,start anchor=west,end anchor=east,"\myperp"] 
\ar[l,phantom,shift right=1.3ex,start anchor=west,end anchor=east,"\myperp"]
& \Fun^\add(\Catx,\Spa)^\hC,
\ar[l,bend left=25,shift left=1.5ex,start anchor=west,end anchor=east,"{(- \circ \fgt)^\hC}"] 
\ar[l,bend right=25,shift right=1.5ex,start anchor=west,end anchor=east,"{(- \circ \fgt)_\hC}"']
\ar[l,phantom,shift left=1.2ex,start anchor=west,end anchor=east,"\myperp"] 
\ar[l,phantom,shift right=1.3ex,start anchor=west,end anchor=east,"\myperp"]
\end{tikzcd}
\]
where the right unlabelled map is pullback along $\Hyp \colon (\Catx)_\hC \rightarrow \Catp$, constitutes a stable recollement, whose associated fracture square is the fundamental fibre square from Corollary~\reftwo{corollary:tate-square}. The images of the left and right adjoints on the right hand side consist of those additive functors $\F$ such that
$\F^\hyp_\hC \to \F$ or $\F \to (\F^\hyp)^\hC$
is an equivalence, respectively.
\end{proposition}
One might call functors satisfying the properties characterising the essential images above \emph{left} and \emph{right hyperbolic}, whence we find that $\F_\hC^\hyp \rightarrow \F$ is the terminal approximation to $\F$ by a left hyperbolic, and $\F \rightarrow (\F^\hyp)^\hC$ is the initial approximation by a right hyperbolic functor.

\begin{proof}
Recall from the discussion above (or see \refone{corollary:hyperbolic-mackey}) that the functor
\[
\Hyp \colon \Catx \rightarrow \Catp
\]
is a both sided adjoint to $\fgt \colon \Catp \rightarrow \Catx$ and that this adjunction is equivariant for the trivial action of $\Ct$ on $\Catp$ and the opponing action on $\Catx$ (i.e.\ it is a relative adjunction over $\mathrm B\Ct$). It follows formally that 
\[
\Nat((\F\circ\fgt)_\hC,\G) \simeq \Nat((\F\circ\fgt),\G)^\hC \simeq \Nat(\F,\G \circ \Hyp)^\hC
\]
and similarly $\Nat(\G,(\F\circ\fgt)^\hC) \simeq \Nat(\G \circ \Hyp,\F)^\hC$, so that we obtain the right hand adjunctions. Unwinding definitions, one finds that in the former case the counit evaluated at some stable $\infty$-category $\C$ is the canonical equivalence
\[
\F(\C) \longrightarrow (\Ct \otimes \F(\C))^\hC
\]
since the inclusion $\F(\C) \rightarrow \F(\C) \oplus \F(\C\op) \longrightarrow \F(\fgt(\Hyp(\C)))$ gives rise to an equivalence $\Ct \otimes \F(\C) \simeq \F(\fgt(\Hyp(\C)))$ by direct inspection. Thus $(-)^\hyp \colon \Fun^\add(\Catp,\Spa) \rightarrow \Fun^\add(\Catx,\Spa)^\hC$ has a fully faithful adjoint, and is therefore a split Verdier projection. Its kernel is precisely $\Fun^\bord(\Catp,\Spa)$ by Lemma~\reftwo{lemma:criterion}, whence Proposition~\reftwo{proposition:criterion-split} and the discussion thereafter give the recollement. 

The claim about the fracture square is true by construction, and the conditions given in the final statement unwind to the counit and unit being equivalences, respectively, which characterise the essential images by the triangle identities.
\end{proof}

The construction of (co)bordifications via the hyperbolisation map $\F^\hyp_\hC \rightarrow \F$ or forgetful map $\F \rightarrow (\F^\hyp)^\hC$ discussed so far are, however, not very suitable for computations of $\F^\bord$ or $\F^\rbord$. Therefore we present two more formulae for the bordification, both of which we put to use in the next section. These also have dual versions which yield cobordifications, though we will not pursue this direction here.
Our approach relies on the following general recognition principle for bordifications:

\begin{lemma}
\label{lemma:bordification}%
Suppose that $\B\colon  \Funadd(\Catp,\Spa) \rightarrow \Funadd(\Catp,\Spa)$ is a functor equipped with a natural transformation $\bet: \id \Rightarrow \B$. Suppose the following conditions hold:
\begin{enumerate}
\item
\label{item:colim}%
$\B$ commutes with colimits;
\item
\label{item:preserve}%
if $\F$ is bordism invariant then $\beta_\F:\F \Rightarrow \B\F$ is an equivalence;
\item
\label{item:metabolic}%
$\B(\F^{\hyp}) \simeq 0$ for every additive $\F \colon \Catp \rightarrow \Spa$.
\end{enumerate}
Then $\beta$ exhibits $\B$ as a bordification functor.
\end{lemma}
Let us explicitly point out that we do not assume a priori that $\B$ takes values in bordism invariant functors. The price is that we have to invest that we already know that there exists a bordification functor into the proof. Direct arguments are also certainly possible, but slightly more cumbersome.

\begin{proof}
Let $\F \colon \Catp \rightarrow \Spa$ be an additive functor. Applying $\B$ to the fibre sequence $\F^{\hyp}_{\hC} \rightarrow \F \rightarrow \F^\bord$ from Corollary~\reftwo{corollary:bord-seq-unique} yields a commutative diagram
\[
\begin{tikzcd}
\F^{\hyp}_{\hC} \ar[r]\ar[d] & \F \ar[r]\ar[d] & \F^\bord \ar[d] \\
\B(\F^{\hyp}_{\hC}) \ar[r] & \B\F \ar[r] & \B(\F^\bord)
\end{tikzcd}
\]
in which both rows are bifibre sequences and the vertical maps are all the respective components of $\beta$. By definition,  $\F^\bord$ is bordism invariant and hence, by property~\reftwoitem{item:preserve}, we get that the right most vertical map above is an equivalence. This implies that the left square is cartesian. On the other hand, by properties~\reftwoitem{item:colim} and~\reftwoitem{item:metabolic}, the lower left corner in the above diagram is $0$, hence the lower right map is an equivalence as well. The right hand square thus exhibits $\B$ as equivalent to $\bord$ under the identity of $\Funadd(\Catp,\Spa)$. %
\end{proof}

Our second formula for bordification is modelled on the classical definition of $\L$-theory spectra via ad-spaces. Its starting point is the \emph{$\ads$-construction}: For $[n]\in\Delta$ we denote by $\T_n = \mathcal{P}_0([n])\op$ the opposite of the poset of nonempty subsets of $[n]$. We observe that $\T_n$ depends functorially on $[n] \in \Del$, giving rise to a simplicial $\infty$-category $\rho(\C,\QF)$: Given a Poincaré \(\infty\)-category $(\C,\QF)$ denote
\[
\rho_n(\C,\QF) = (\Fun(\T_n,\C),\QF^{\T_n})
\]
the cotensor of $(\C,\QF)$ by $\T_n$. Since $\T_n$ is the reverse face poset of $\Delta^n$, we find from Proposition~\refone{proposition:posets-of-faces} that the hermitian \(\infty\)-categories $\rho_n(\C,\QF)$ are Poincaré for every $[n] \in \Del$ and from Proposition~\refone{proposition:poset-of-faces-functoriality} that the hermitian functor $\sig^*\colon \rho_n(\C,\QF) \rightarrow \rho_m(\C,\QF)$ is Poincaré for every $\sig\colon [m] \rightarrow [n]$ in $\Del$. We may hence consider $\rho{(\C,\QF)}$ as a simplicial object in $\Catp$.

\begin{definition}
\label{definition:rho-construction}%
Let $\E$ be an $\infty$-category with sifted colimits. Given a functor $\F\colon \Catp \rightarrow \E$, we denote by $\ads\F\colon \Catp \rightarrow \E$ the functor given by
\[
\ads \F(\C,\QF) = |\F\rho(\C,\QF)|.
\]
\end{definition}

Using the functoriality of the cotensor construction, we may promote the association $\F \mapsto \ads\F$ to a functor
\[
\ads\colon \Fun(\Catp,\E) \lrar \Fun(\Catp,\E) .
\]
The inclusion of vertices then equips $\ads$ with a natural transformation
$ b_{\F}\colon \F \rightarrow \ads\F$.

In this section, we consider the $\ads$-construction only in the case when $\E=\Spa$, as this entails great simplifications (though the case $\E=\Sps$ is fundamental for the discussion of $\L$-theory in section~\reftwo{subsection:L+tate}). The key is that for a stable $\E$, the collection of additive functors from $\Catp$ to $\E$ is closed under colimits inside the $\infty$-category of all functors. Since in addition, the functor $\rho_n:\Catp \rightarrow \Catp$ preserves split Poincaré-Verdier sequences by Proposition~\reftwo{proposition:(co)tensor-Verdier}, it follows that $\ads\F$ is additive whenever $\F:\Catp \rightarrow \Spa$ is. 
In particular, we may consider $\ads$ as an (evidently colimit preserving) functor
\[
\ads\colon \Funadd(\Catp,\Spa) \longrightarrow \Funadd(\Catp,\Spa) .
\]

\begin{remark}
\label{remark:adsadd}%
The analogous statement with target category $\Sps$ requires an additivity theorem for the $\ads$-construction. For the L-space functor, such a result is indeed available. In fact, Lurie showed in \cite{Lurie-L-theory}*{Lecture~8,~Corollary~9} (see Theorem~\reftwo{theorem:L-is-localising}) that $\Lspace = \ads(\Poinc)$ is even Verdier-localising, generalising results of  Ranicki in more classical language; see e.g.\ \cite{Ranickiblue}*{Proposition 13.11}. 
\end{remark}

Next, we set out to show:

\begin{proposition}
\label{proposition:rho-bordism}%
\label{corollary:rho-bord}%
Let $\F\colon \Catp \rightarrow \Spa$ be an additive functor. Then
\begin{enumerate}
\item
\label{item:FadsF}%
if $\F$ is bordism invariant, then the map $b_\F\colon  \F \Rightarrow \ads\F$ is an equivalence.
\item
\label{item:adsFhyp}%
$\ads(\F^{\hyp}) \simeq 0$.
\end{enumerate}
In particular, the natural transformation $b$ exhibits $\ads$ as a bordification functor.
\end{proposition}

For the proof of Proposition~\reftwo{proposition:rho-bordism}, we denote by $\mathcal P([n])$ the full power set of $[n]$ and endow $\Fun(\mathcal P([n])\op,\C)$ with the hermitian structure $\QF^\mathrm{tf}$ that sends a cubical diagram
$\varphi\colon \mathcal P([n])\op\to \C$ to the total fibre of
$\QF\qshift{1}\circ \varphi\op$; through the isomorphism
\[
\mathcal P([n])\op\cong \mathcal \prod_{i= 0}^n[1]
\]
the hermitian \(\infty\)-category
$(\Fun(\mathcal P([n])\op, \C),\QF^{\mathrm{tf}})$
is equivalent to $\Met^{(n+1)}(\C, \QF\qshift{1})$, but in the form given, it is clear that it assembles into a functor $\Cath \rightarrow \sCath$. Through the identification as an iterated metabolic category, it is, however, easy to check that it restricts to $\Catp \rightarrow \sCatp$.

\begin{lemma}
\label{lemma:rho_vs_power_set}%
The sequence
\[
\rho_n(\C, \QF) \longrightarrow (\Fun(\mathcal P([n])\op, \C), \QF^\mathrm{tf})\xrightarrow{\ev_\emptyset} (\C, \QF\qshift 1)
\]
is a Poincaré-Verdier sequence for all Poincaré $\infty$-categories $(\C, \QF)$ and $n \in \mathbb N$. Furthermore, there are equivalences
\[
(\Fun(\mathcal  P([-])\op, \Ar(\C)), \QF_\met^\mathrm{tf}) \simeq \dec(\Fun(\mathcal  P([-])\op, \C), \QF^\mathrm{tf}) \simeq \Met(\Fun(\mathcal  P([-])\op, \C), \QF^\mathrm{tf})
\]
of simplicial Poincaré $\infty$-categories.
\end{lemma}

Note that the left term in the second display is just $\Fun(\mathcal  P([-])\op, \D), \QFD^\mathrm{tf})$ for $(\D,\QFD) = \Met(\C,\QF)$.

\begin{proof}
The map $\ev_\emptyset$ is given by evaluation of the cubical diagram
$\varphi$ at $\emptyset$, together with the canonical projection of
hermitian functors. Under the equivalence of the middle term with
$\Met^{(n+1)}(\C, \QF\qshift 1)$, the second map is the $(n+1)$-fold
iteration of the map $\met\colon \Met(\C, \QF\qshift 1)\to (\C,\QF\qshift 1)$;
thus it is a split Poincaré-Verdier projection. 
Its kernel is equivalent to the first term by restriction along $\T_n = \mathcal
P_0([n])\op\subset \mathcal P([n])\op$ and the equivalence
\[
\lim_{\emptyset \neq A \subseteq [n]} \QF\circ \varphi\op(A) \simeq
\fib\bigl(0 \longrightarrow \lim_{\emptyset \neq
A \subseteq [n]} \QF\qshift 1\circ
\varphi\op(A)\bigr),
\]
since for $\varphi \in \ker(\ev_\emptyset)$ the second term is equivalent to the total fibre of $\QF\qshift 1\circ\varphi\op$.

For the second claim, note that commuting limits and functor categories gives equivalences
\[
\Fun(\mathcal  P([-])\op , \Ar(\C)), \QF_\met^\mathrm{tf}) \simeq \Fun(\mathcal  P([-])\op \times \Delta^1,\C),\QF^\mathrm{tf}) \simeq \Met(\Fun(\mathcal P([-])\op,\C),\QF^\mathrm{tf})
\]
and the middle term is the requisite décalage by inspection.
\end{proof}

\begin{proof}[Proof of Proposition~\reftwo{proposition:rho-bordism}]
If $\F$ is bordism invariant, then it vanishes on the middle term of the sequence of Lemma~\reftwo{lemma:rho_vs_power_set} (which is an iterated
metabolic construction); we conclude that the left term becomes constant in $n$, after application of $\F$. This shows \reftwoitem{item:FadsF}. 
To show \reftwoitem{item:adsFhyp}, we note that for a Poincaré $\infty$-category $(\C,\QF)$ 
\[
\F^{\hyp}\left(\Fun(\mathcal P[-], \C)\right) = \F\left(\Hyp(\Fun(\mathcal P[-], \C))\right)\simeq \F\left(\Met(\Fun(\mathcal P[-], \C),\QF^\mathrm{tf})\right)
\]
by Corollary~\reftwo{corollary:met-hyp}. But this is the décalage of $\F(\Fun(\mathcal P[-], \C),\QF^\mathrm{tf})$ by Lemma~\reftwo{lemma:rho_vs_power_set} with augmentation induced by $\ev_\emptyset$. Interpreting this map as a map of split simplicial objects, we conclude that its fibre $\F^\hyp(\rho_n(\C,\QF))$ is split over $0$, and therefore has contractible realisation. The final claim follows from Lemma~\reftwo{lemma:bordification}.
\end{proof}

\begin{remarks}
\begin{enumerate}
\item
In \paperfour we show that, more generally, the construction $\ads$ actually gives the universal bordism invariant replacement of an arbitrary (in particular, not necessarily additive) functor $\Catp \rightarrow \Sps$.
\item 
We expect an additivity statement similar to Theorem \reftwo{theorem:additivity} to hold for $\ads(\F)$ for functors $\F \colon \Catp \rightarrow \Sps$, but do not pursue this here, since we will be mainly interested in the case $\F=\Poinc$ in which case Lurie has provided as stronger result which we record in \reftwo{theorem:L-is-localising} below, and on which we will base our analysis of $\L$-spectra.
\item
If $\F\colon \Catp \to \Spa$ is an additive functor which factors through the forgetful functor $\Catp \to \Catx$ then $\F$ is a retract of $\F^{\hyp}$ (after forgetting the $\Ct$-action)
and hence Proposition~\reftwo{proposition:rho-bordism}\;\reftwoitem{item:adsFhyp} implies that $\ads(\F) \simeq 0$. More generally, this argument shows that if $\F$ is any additive functor on $\Catx$ then the non-hermitian analogue of the $\ads$-construction vanishes on $\F$.
\item
In \cite{WW-duality}*{Lemma 9.3}, Weiss and Williams give a direct verification that $\ads(\K) \simeq 0$. Their proof immediately generalises to give a different argument for the vanishing of $\ads$ %
on additive functors on $\Catx$.
\end{enumerate}
\end{remarks}

Finally, we present a third formula for bordification, obtained by iterating the boundary map $\F(\C,\QF) \rightarrow \SS^1 \ssmash \F(\C,\QF\qshift{-1})$ of the metabolic fibre sequence. 

\begin{definition}
\label{definition:stab-construction}%
Let $\F\colon \Catp \rightarrow \Spa$ be an additive functor. We define its \defi{stabilisation} $\metstab \F$ by the formula
\[
(\metstab \F)(\C,\QF) = \colim(\F(\C,\QF) \longrightarrow \SS^1 \ssmash \F(\C,\QF\qshift{-1}) \longrightarrow  \SS^2 \ssmash \F(\C,\QF\qshift{-2}) \longrightarrow \dots),
\]
with structure maps the shifts of the boundary map for $\F$, and we denote by
$
\sig^{\infty}_{\F}\colon\F \to \metstab \F 
$
the arising natural transformation.
\end{definition}

Recall from the discussion preceding Corollary~\reftwo{corollary:Q=suspinsp} that the boundary map $\F(\C,\QF) \to \SS^1 \ssmash \F(\C,\QF\qshift{-1})$ of the metabolic fibre sequence
\[
\F(\C,\QF\qshift{-1}) \longrightarrow \F(\Met(\C,\QF)) \xrightarrow\met \F(\C,\QF)
\]
is also modelled by the inclusion of vertices
$
\sigma_\F \colon \F(\C,\QF) \to |\F\Q(\C,\QF)|.
$
So we equally well find that
\[
\metstab\F(\C,\QF) \simeq \colim(\F(\C,\QF) \xrightarrow{\sig_{\F}}  |\F\Q(\C,\QF)| \xrightarrow{|\sig_{\F} \Q|}  |\F\Q^{(2)}(\C,\QF)| \longrightarrow \dots).
\]
arises from iteration of the (unshifted!) $\Q$-construction.

\begin{proposition}
\label{proposition:Q-bordism}%
\label{corollary:Q-oo-bord}%
Let $\F\colon \Catp \rightarrow \Spa$ be an additive functor. Then
\begin{enumerate}
\item
\label{item:FmetstabF}%
if $\F$ is bordism invariant, the map $\sig^\infty_{\F}\colon \F \rightarrow \metstab\F$ is an equivalence.
\item
\label{item:metstabFhyp}%
$\metstab(\F^{\hyp}) \simeq 0$. %
\end{enumerate}
In particular, the transformation $\sig^{\infty}$ exhibits $\metstab$ as a bordification.
\end{proposition}

Dually, the limit over the stabilisation maps give a cobordification, as is easy to check, but we shall not make use of that fact.

\begin{proof}
Property~\reftwoitem{item:FmetstabF} follows immediately from Corollary~\reftwo{proposition:bordism-invariant-sig}. To prove \reftwoitem{item:metstabFhyp}, it will suffice to show that for any additive $\F\colon \Catp \rightarrow \Spa$ and any stable $\infty$-category $\C$, the boundary map 
\[
\F(\Hyp(\C)) \lrar \SS^1 \ssmash \F(\Hyp(\C)\qshift{-1})
\]
is null-homotopic. But this follows immediately from the metabolic functor $\met \colon \Met(\Hyp(\C)) \rightarrow \Hyp(\C)$ being split by Corollary~\refone{corollary:met-hyp-splits}. Since $\metstab$ evidently commutes with colimits and preserves additivity, we can apply Lemma~\reftwo{lemma:bordification} to obtain the final claim.
\end{proof}

\begin{corollary}
\label{corollary:filtration}%
For every group-like additive $\F\colon \Catp \rightarrow \Sps$ the %
map 
\[
\pi_{i}\CCob^\F(\C,\QF) \longrightarrow \pi_i(\CCob^{\F})^\bord(\C,\QF)
\]
is an isomorphism for $i<0$ and for $i=0$ becomes the canonical projection $\pi_0\F(\C,\QF) \to \pi_0|\Cob^{\F}(\C,\QF\qshift{-1})|$. %
In particular, for every $i \in \ZZ$ we have canonical isomorphisms
\[
\pi_i(\CCob^{\F})^\bord(\C,\QF) \cong \pi_{-1}(\CCob^{\F})^\bord(\C,\QF\qshift{-(i+1)}) \cong \pi_0|\Cob^\F(\C,\QF\qshift{-(i+1)})|.
\]

\end{corollary}

The group $\pi_i(\CCob^{\F})^\bord(\C,\QF)$ is thus the $\F$-based cobordism group of $(\C,\QF\qshift{-i})$. In fact, the proof below will show that the colimit description for $\F^\bord(\C,\QF) \simeq \metstab\F(\C,\QF)$ stabilises on $\pi_i$ after step $i$. %

\begin{proof}
By Corollary~\reftwo{corollary:CCobcon} the spectra $\SS^k \otimes \CCob^\F(\Met(\C,\QF\qshift{-k}))$ are $k$-connective, and so the maps in the colimit sequence
\[
(\CCob^{\F})^\bord(\C,\QF) = \colim(\CCob^\F(\C,\QF) \longrightarrow \SS^1 \ssmash \CCob^\F(\C,\QF\qshift{-1}) \longrightarrow  \SS^2 \ssmash \CCob^\F(\C,\QF\qshift{-2}) \longrightarrow \dots),
\]
have increasingly connective fibres. In particular, all of these 
induce an isomorphism on negative homotopy groups, and all but the first induce an isomorphism on $\pi_0$, in which case we furthermore apply %
Corollary~\reftwo{remark:left-square-met}.
We conclude using Proposition~\reftwo{proposition:bordism-invariant-sig}.
\end{proof}

\begin{remark}
\label{remark:compare-bordification}%
Given an arbitrary additive functor $\F\colon \Catp \to \Spa$, the fundamental fibre sequence of Corollary~\reftwo{corollary:tate-square} combined with Corollary~\reftwo{corollary:CCobcon} shows that %
the fibre of the comparison map 
\[
(\Om^{\infty}\F)^{\bord} = \Om^{\infty}((\CCob^{\Omega^\infty\F})^{\bord}) \longrightarrow \Om^{\infty}(\F^{\bord})
\]
is naturally equivalent to $\Om^{\infty}(\G^{\hyp}_{\hC})$, where 
$\G = \cof(\CCob^{\Om^{\infty} \F} \to \F)$. Corollary \reftwo{corollary:filtration} always computes the homotopy groups of the source, but $\G$ can generally be hard to access: For example, for
 $\F = \KGW$ the Karoubi-Grothendieck-Witt spectrum we study in \paperfour, the defect term is $\Omega^\infty ((\tau_{< 0} \KK)_\hC)$. This observation ultimately encodes part of the classical Rothenberg sequences for decorated $\L$-groups, and we shall upgrade the statement to the spectrum level in \paperfour.%

\end{remark}

We finally mention that the defining filtration of $\sig^{\infty}_{\F}\colon \F\to \metstab(\F)\simeq \F^\bord$ is compatible with the canonical filtration of $(-)^\hC\to (-)^\tC$
in the sense that from naive Karoubi periodicity \reftwo{proposition:naive-karoubi}, we obtain a diagram
\[
\begin{tikzcd}[column sep=small]
\F(\C, \QF) \ar[r] \ar[d]
& \SS^1\ssmash \F(\C, \QF\qshift{-1}) \ar[r] \ar[d]
& \SS^2\ssmash \F(\C, \QF\qshift{-2}) \ar[r] \ar[d]
& \dots \ar[r] 
& \metstab(\F)(\C, \QF) \ar[d] \\
(\F^{\hyp})^\hC(\C, \QF) \ar[r] 
& (\SS^\sigma \ssmash \F^{\hyp})^\hC(\C, \QF) \ar[r]
& (\SS^{2\sigma} \ssmash \F^{\hyp})^\hC(\C, \QF) \ar[r]
& \dots \ar[r] 
& \colim_{n \in \mathbb N} (\SS^{n\sigma} \ssmash \F^{\hyp})^\hC(\C, \QF) 
\end{tikzcd}
\]
with right most vertical arrow being equivalent to $ \F^\bord(\C, \QF) \to (\F^{\Hyp})^\tC(\C, \QF)$. Considering horizontal fibres each of the occuring squares is cartesian as an application of Corollary \reftwo{corollary:met-hyp}, which provides an alternative proof of Corollary~\reftwo{corollary:tate-square}. 
This filtration of the fundamental fibre square plays a major role in the works of Weiss--Williams \cite{WWIII} in case of LA-theory (which we identify with Grothendieck-Witt theory of parametrised spectra in section~\reftwo{subsection:Weiss-Williams} below).

\subsection{The genuine hyperbolisation of an additive functor}
\label{subsection:Mackey}%

In this final subsection, we recast the fundamental fibre square from Corollary~\reftwo{corollary:tate-square} of an additive functor $\F \colon \Catp \rightarrow \Spa$ as the isotropy separation square of a genuine $\Ct$-spectrum, that is a spectral Mackey functor for the group $\Ct$, refining the hyperbolisation $\F^\hyp(\C,\QF) \in \Spa^\hC$. This allows for a convenient way of combining Karoubi periodicity with the shifting behaviour of bordism invariant functors; see Theorem~\reftwo{theorem:periodicity} below.

Let us briefly recall the notion of a spectral Mackey functor, see also \paperone section~7.4 for further details. For a finite discrete group $G$, we denote by $\Span(\mathrm{Fin}_G)$ the $\infty$-category of spans of finite $G$-sets, first used for the purposes of equivariant homotopy theory in \cite{Barwick-MackeyI}*{Df.~3.6} (under the name effective Burnside category). 
We set $\Spa^{\mathrm{g}G}=\Fun^\times(\Span(\mathrm{Fin}_G),\Spa)$, the $\infty$-category of product preserving functors, and refer to \cite{Nardin-stability}*{Appendix A}, \cite{GepnerMeier}*{Appendix C}, or \cite{CMNN2} for comparisons with other models of genuine $G$-spectra. 

Restricting to the relevant case $G=\Ct$, we have associated to $E \in \Spagc$ an underlying $\Ct$-spectrum $uE \in \Spa^\hC$, the evaluation of $E$ at the $\Ct$-set $\Ct$ with $\Ct$-action through
\[
\Ct \xleftarrow{\id} \Ct \xrightarrow{\operatorname{flip}} \Ct
\]
and the genuine fixed points $E^{\g\Ct} \in \Sp$, given by the evaluation at the $\Ct$-set $\ast$, connected by restriction $\mathrm{res} \colon E^{\gCt}\to E^\hC$ and transfer maps $\mathrm{tr}\colon E_\hC\to E^{\gCt}$ coming from the spans
\begin{equation*}
* \leftarrow \Ct \xrightarrow{\id} \Ct \qquad \text{and} \qquad \Ct \xleftarrow{\id} \Ct \rightarrow *
\end{equation*}
Finally, to $E$ is associated $E^{\geofix} \in \Sp$, the geometric fixed points, which for us are most easily defined as the cofibre of the transfer; they can also more conceputally be defined by left Kan extending $E$ along the fixed points functor $\Span(\mathrm{Fin}_{\Ct}) \rightarrow \Span(\mathrm{Fin})$ and then evaluating at $\ast$, see \cite{Barwick-MackeyI}*{B.7}). %
In fact, there is a stable recollement
\[
\begin{tikzcd}
[column sep=7ex]
\Spa
\ar[r,"R" description] 
& \Spa^\gCt \ar[r] 
\ar[l,bend left=25] 
\ar[l,bend right=25,pos=.4,"{(-)^\geofix}"']
\ar[l,phantom,shift left=1.2ex,start anchor=west,end anchor=east,"\myperp"] 
\ar[l,phantom,shift right=1.3ex,start anchor=west,end anchor=east,"\myperp"]
& \Spa^\hC
\ar[l,bend left=20,"{(-)^\sym}"]
\ar[l,bend right=25,"{(-)^\qdr}"']
\ar[l,phantom,shift left=1.2ex,start anchor=west,end anchor=east,"\myperp"] 
\ar[l,phantom,shift right=1.3ex,start anchor=west,end anchor=east,"\myperp"]
\end{tikzcd}
\]
where $R$ is given by restriction along the fixed point functor $\Span(\mathrm{Fin}_{\Ct}) \rightarrow \Span(\mathrm{Fin})$, under the identification $\Spa \simeq \Fun^\times(\Span(\mathrm{Fin}),\Spa)$, and the lower left functor takes $X$ to the fibre of $X^\gC \rightarrow X^\hC$, and $(-)^\sym$ and $(-)^\qdr$ are Borel completion and cocompletion, respectively, defined as the adjoints to the forgetful functor $u$ described above.
In particular, there results a cartesian square
\[
\begin{tikzcd}
E^\gC \ar[r] \ar[d] & E^\geofix \ar[d] \\
E^\hC \ar[r] & E^\tC,
\end{tikzcd}
\]
as a special case of the general theory in section~\reftwo{appendix:AppIIA}, specifically Remark \reftwo{remark:explicit-classifyingsplit}. In fact, \reftwo{corollary:classificationsplitverdierdumdidum} implies that the entire $\infty$-category $\Spa^\gCt$ can equivalently be described via such isotropy separation squares. %

In Corollary~\refone{corollary:hyperbolic-mackey},
we showed:
\begin{theorem}
\label{theorem:genuine}%
The construction of hyperbolic Poincar\'e $\infty$-categories canonically refines to a functor
$
\gHyp\colon \Catp\to \Fun^{\times}(\Span(\mathrm{Fin}_{\Ct}),\Catp)
$
together with natural equivalences :
\begin{enumerate}
\item $\gHyp(\C,\QF)^\gCt \simeq (\C,\QF)$ of Poincaré $\infty$-categories,
\item $\undl(\gHyp(\C,\QF)) \simeq \HypG(\C,\QF)$ of $\Ct$-Poincaré $\infty$-categories
\end{enumerate}
under which transfer and restriction%
\[
\gHyp(\C,\QF)_\hC \rightarrow \gHyp(\C,\QF)^\gC \quad \text{and} \quad \gHyp(\C,\QF)^\gC \rightarrow \gHyp(\C,\QF)^\hC,
\]
are naturally identified with
\[
\hyp \colon \HypG(\C,\QF)_\hC \rightarrow (\C,\QF) \quad \text{and} \quad \fgt \colon (\C,\QF) \rightarrow \HypG(\C,\QF)^\hC.
\]
\end{theorem}

\begin{definition}
    Let $\F\colon\Catp\to \Spa$ be an additive functor. Then we call the composite
\[
    \Catp\xrightarrow{\gHyp} \Fun^\times(\Span(\Ct),\Catp)
    \xrightarrow{\F}\Fun^\times(\Span(\Ct),\Spa) = \Spagc\,.
\]
   the \emph{genuine hyperbolisation} $\F^\ghyp$ of $\F$. 
\end{definition}

Now from the identification of the transfer in Theorem~\reftwo{theorem:genuine} and Proposition~\reftwo{proposition:phi-bord}, there results an identification $\F^{\ghyp}(\C,\QF)^\geofix \simeq \F^\bord(\C,\QF)$ and by the universal property of bordifications, this determines the entire isotropy separation square. We conclude:

\begin{corollary}
\label{corollary:geometric-fixed-points-are-bordism-invariant}%
The isotropy separation square of the genuine $\Ct$-spectrum $\F^\ghyp(\C,\QF)$ is naturally identified with the fundamental fibre square, in symbols
\[
\begin{tikzcd}
\F^\ghyp(\C,\QF)^\gC \ar[r] \ar[d] & \F^\ghyp(\C,\QF)^\geofix \ar[d] \\
\F^\ghyp(\C,\QF)^\hC \ar[r] & \F^\ghyp(\C,\QF)^\tC
\end{tikzcd} 
\quad \simeq \quad 
\begin{tikzcd}
\F(\C,\QF) \ar[r] \ar[d] & \F^\bord(\C,\QF) \ar[d] \\
(\F^\hyp)^\hC(\C,\QF) \ar[r] & (\F^\hyp)^\tC(\C,\QF),
\end{tikzcd}
\]
for any additive functor $\F \colon \Catp \rightarrow \Sps$ and Poincaré \(\infty\)-category $(\C,\QF)$.
\end{corollary}

Finally, we use the genuine $\Ct$-spectrum $\F^\ghyp(\C,\QF)$ to combine naive periodicity with the behaviour of bordism invariant functors under shifting.

\begin{lemma}
\label{lemma:genuineofmet}%
Let $\F\colon\Catp\to \Spa$ be an additive functor and $\C$ a stable $\infty$-category. Then the map of genuine $\Ct$-spectra
\[
\Ct \otimes \F(\Hyp\C) \longrightarrow \F^\ghyp(\Hyp\C)
\]
where $\Ct \otimes - \colon \Spa \rightarrow \Spa^\gCt$ is left adjoint to the forgetful functor, and the map is adjoint to the diagonal
\[
\F(\Hyp\C)\to \F(\Hyp\C)\oplus \F(\Hyp\C) \simeq \F(\Hyp\C\times\Hyp\C) \simeq \F(\HypG(\Hyp\C))
\]
is an equivalence. In particular, $\F^\ghyp(\Met(\C,\QF)) \simeq \Ct \ssmash \F(\Hyp\C)$.
\end{lemma}
\begin{proof}
The map is an equivalence both on underlying spectra and on geometric fixed points: On underlying spectra, this follows immediately from the corresponding statement
\[
\HypG(\Hyp(\C)) \simeq \Ct \otimes \Hyp(\C)
\]
on underlying $\Ct$-Poincaré $\infty$-categories from~\refone{remark:hyp-of-hyp}.
Furthermore, both spectra have vanishing geometric fixed points: The left hand side by the symmetric monoidality of geometric fixed points together with $\Ct^\gCt = \emptyset$, the right hand side by bordism invariance. 
\end{proof}

As a direct generalisation of Proposition~\reftwo{proposition:naive-karoubi}, we then have:
\begin{theorem}[Genuine Karoubi periodicity]
\label{theorem:periodicity}%
There is a canonical equivalence of genuine $\Ct$-spectra
\[
\F^\ghyp(\C,\QF\qshift{-1}))\simeq \SS^{\sigma-1} \ssmash \F^\ghyp(\C,\QF),
\]
natural in the Poincar\'e $\infty$-category $(\C,\QF)$ and the additive functor $\F \colon \Catp \to \Sp$. Furthermore, under this equivalence, the boundary map
\[
\F^\ghyp(\C,\QF)) \longrightarrow \SS^1 \ssmash \F^\ghyp(\C,\QF\qshift{-1}))
\]
of the metabolic fibre sequence is induced by the inclusion $S^0 \rightarrow S^\sigma$ as the fixed points.
\end{theorem}

In particular, passing to geometric fixed points we recover the equivalence $\F^\bord(\C,\QF\qshift{i}) \simeq \SS^i \otimes \F^\bord(\C,\QF)$ from Proposition~\reftwo{proposition:bordism-invariant-sig}.

\begin{proof}
Given the previous lemma, the proof of Proposition~\reftwo{proposition:naive-karoubi} applies essentially verbatim, when interpreted in the $\infty$-category of genuine $\Ct$-spectra: Lemma~\reftwo{lemma:genuineofmet} identifies the once-rotated metabolic fibre sequence
\[
\F^\ghyp(\Met(\C,\QF)) \xrightarrow{\met} \F^\ghyp(\C,\QF) \xrightarrow{\partial} \SS^1 \ssmash \F^\ghyp(\C,\QF\qshift{-1})
\]
with
\[
\Ct \otimes \F(\Hyp(\C)) \longrightarrow \F^\ghyp(\C,\QF) \longrightarrow \SS^\sigma \ssmash \F^\ghyp(\C,\QF)
\]
obtained by tensoring $\F^\ghyp(\C,\QF)$ with $\SS[\Ct] \rightarrow \SS \rightarrow \SS^\sigma$.
\end{proof}

\section{Grothendieck-Witt theory}
\label{section:GW}%

We start this section by defining the Grothendieck-Witt space $\gw(\C,\QF)$ and record its properties, as specialisations of the general results of the previous section to the case $\F = \Poinc \in \Funadd(\Catp,\Sps)$. We then proceed to analyse the Grothendieck-Witt spectrum in the same manner, in particular identifying its hyperbolisation as $\K$-theory and its bordification as $\L$-theory. 

This leads to the identification of the homotopy type of the algebraic cobordism categories in Corollary~\reftwo{corollary:algGMTW}, the fundamental fibre square in Corollary~\reftwo{corollary:tate-square-L}, localisation sequences for Grothendieck-Witt spectra of discrete rings in Corollary~\reftwo{corollary:localisation-sequences-papertwo} and our generalisation of Karoubi periodicity in Corollaries~\reftwo{corollary:karoubi-fundamental} and~\reftwo{corollary:Karoubiper}, constituting the main results of the present paper.

In the final subsection, we spell out the relation of our constructions to the $\LA$-theory of Weiss and Williams from \cite{WWIII}. In particular, our results provide a cycle model for the infinite loop spaces of their spectra. %

\subsection{The Grothendieck-Witt space}
\label{subsection:gw-space}%

In this section, we define the Grothendieck-Witt space of a Poincaré $\infty$-category $(\C,\QF)$, whose homotopy groups are by definition the higher Grothendieck-Witt groups of $(\C,\QF)$. 

Recall that a functor $\Catp\to \E$ into an $\infty$-category with finite limits is called additive if it carries split Poincaré-Verdier squares to cartesian squares, see section~\reftwo{subsection:additive-localizing}. Additive functors automatically take values in $\Einf$-monoids (with respect to the cartesian product in $\E$) but may well not be group-like; the functor $\Poinc\colon \Catp\to \Sps$ taking Poincaré objects being the first example. Denoting by $\Funadd(\Catp, \E)\subseteq \Fun(\Catp, \E)$ the full subcategory of additive functors, Corollary~\reftwo{corollary:universal} asserts that the forgetful functor  
\[
\Funadd(\Catp, \Grp_{\Einf}(\Sps)) \longrightarrow \Funadd(\Catp, \Sps)
\]
admits a left adjoint $(-)^\grp$, the group completion functor.

\begin{definition}
\label{definition:gw-space}%
We define the \defi{Grothendieck-Witt space} functor $\GWspace\colon  \Catp \lrar \Grp_{\Einf}(\Sps)$ to be the group-completion
\[
\GWspace(\C,\QF) = \Poinc^\grp(\C,\QF),
\]
of the functor $\Poinc \in \Funadd(\Catp,\Sps)$. Furthermore, for a Poincaré $\infty$-category $(\C,\QF)$ and $i\geq 0$, we set
\[
\GW_i(\C,\QF) = \pi_i\GWspace(\C,\QF),
\]
the \emph{Grothendieck-Witt-groups} of $(\C,\QF)$.
\end{definition}

Let us also mention that we frequently regard $\GWspace$ as a functor $\Catp \rightarrow \Sps$ without further comment.
As a direct reformulation of the definition of the Grothendieck-Witt functor, we record:

\begin{observation}
\label{observation:univ-GW-space}%
The functor $\GWspace \colon \Catp \rightarrow \Sps$ is additive and group-like, and it is the initial such functor under $\Poinc \colon \Catp \rightarrow \Sps$.
\end{observation}

We will show in Corollary~\reftwo{corollary:gwadd} that $\GWspace$ is in fact Verdier localising (and not just additive); that is, it takes all Poincaré-Verdier squares to cartesian squares (and not just the split ones).

Also note that we have already introduced a group $\GW_0(\C,\QF)$ explicitly in \S\refone{subsection:GW-group}, which by 
Theorem~\reftwo{theorem:pi1Cob} indeed agrees with $\pi_0\GWspace(\C,\QF)$:
\begin{corollary}
\label{corollary:pi-0gw}%
The natural map 
$\pi_0\Poinc(\C,\QF) \to \pi_0\GWspace(\C,\QF)$
exhibits the target as the quotient of the source by the congruence relation generated by
\[
[X,\qone] \sim [\hyp(W)],
\]
where $(X,\qone)$ runs through the Poincaré objects of $(\C,\QF)$ with Lagrangian $W\to X$.
\end{corollary}

In particular, $\GW_0(\C,\QF)$ is the quotient of $[\pi_0\Poinc(\C,\QF)]^\grp$ by the subgroup spanned by the differences $[X,\qone] - [\hyp(W)]$, but part of the statement is that one does not need to complete $\pi_0\Poinc(\C,\QF)$ to a group in order to obtain $\GW_0(\C,\QF)$ as a quotient. From Corollary~\reftwo{corollary:inversion} or indeed from Lemma~\refone{lemma:basic-GW-group},
we find that for $[X,q] \in \GW_0(\C,\QF)$ we have
\[
-[X,q] = [X,-q] + \hyp(\Omega X).
\]

\begin{remark}
Let us explicitly warn the reader that $\GWspace(\C,\QF)$ is %
usually \emph{not} the group-completion $\Poinc(\C,\QF)^\grp$ of $\Poinc(\C,\QF)$.%
Indeed, this termwise group completion of $\Poinc$ does not yield an additive functor. %
For the relation between Grothendieck-Witt spectra defined by group-completions of forms and those considered here, see \cite{comparison}.
\end{remark}

From the results of the previous section, we obtain several formulae for $\GWspace(\C,\QF)$: Recall from Definition~\reftwo{definition:cobcat} the cobordism $\infty$-category $\Cob(\C,\QF)$ associated to the Segal space $\Poinc\Q(\C,\QF\qshift{1})$ given by the hermitian $\Q$-construction. From Corollary~\reftwo{corollary:universal}, we find:

\begin{corollary}
There are canonical equivalences
\[
\GWspace(\C,\QF) \simeq \Omega|\Cob(\C,\QF)| \simeq \Omega|\Poinc\Q(\C,\QF\qshift{1})|
\]
natural in the Poincaré $\infty$-category $(\C,\QF)$.
\end{corollary}

These formulae are in accordance with the usual definition of the $\K$-space $\Omega|\Span(\C)| \simeq \Omega|\core\Q(\C)|$ of $\C$.

\begin{corollary}
\label{corollary:filtcolimpresgwspace}%
The functor $\GWspace\colon \Catp\to \Sps$ preserves filtered colimits.
\end{corollary}

\begin{proof}
It is a composite of functors 
\[
\Catp\xrightarrow{(-)\qshift 1} \Catp \xrightarrow{\Q} \sCatp \xrightarrow{\Poinc} \sSps \xrightarrow{\vert-\vert} \Sps \xrightarrow\Omega \Sps
\]
each of which preserves filtered colimits: The first is an equivalence, for $\Poinc$ this is \refone{proposition:forms-poinc-compact}, for $\Q_n$ this follows from \refone{corollary:limit-preservation}, and for the remaining two functors this is well-known.
\end{proof}

Classically, the Grothendieck-Witt space is often defined as the fibre of the forgetful functor from the hermitian to the usual $\Q$-construction. We obtain such a description from the fibre sequence, as we explain next. Applying the hermitian $\Q$-construction and $\Poinc$ to the (split!) metabolic Poincaré-Verdier sequence
\[
(\C,\QF) \longrightarrow \Met(\C,\QF\qshift{1}) \xrightarrow{\met} (\C,\QF\qshift{1})
\]
results in the fibre sequence
\[
|\Poinc\Q(\C,\QF)| \longrightarrow |\Poinc\Q\Met(\C,\QF\qshift{1})| \xrightarrow{\met} |\Poinc\Q(\C,\QF\qshift{1})|.
\]
Now from Proposition~\reftwo{proposition:CobMet=CobHyp} and Example~\reftwo{example:CobHyp-equal-Span}, we obtain  canonical equivalences 
\[
|\Poinc\Q\Met(\C,\QF\qshift{1})| \xrightarrow{\mathrm{can}} |\Poinc\Q\Hyp(\C)| \xrightarrow{\fgt} |\core\Q(\C)|.
\]
Inserting these and rotating the resulting fibre sequence once to the left we find:
\begin{corollary}
\label{corollary:classicaldefofgw}%
\label{corollary:GWspace-is-fiber-Pn-to-Cr}%
The metabolic fibre sequence identifies with
\[
|\Poinc\Q(\C,\QF)| \xrightarrow\fgt |\core\Q(\C)| \xrightarrow{\hyp} |\Poinc\Q(\C,\QF\qshift{1})|
\]
for every Poincar\'e $\infty$-category $(\C,\QF)$ and there is a natural equivalence
\[
\GWspace(\C,\QF) \simeq \fib(|\Poinc\Q(\C,\QF)| \xrightarrow\fgt |\core\Q(\C)|).
\]
\end{corollary}

This formula for the Grothendieck-Witt space is the transcription of the classical definition for example from \cite{schlichting-exact} into our framework. We will say more about the metabolic fibre sequence in section~\reftwo{subsection:Genauer-Bott-Karoubi} below.

\subsection{The Grothendieck-Witt spectrum}
\label{subsection:The-GW-spectrum-Paper2}%

Our next goal is to deloop the Grothendieck-Witt-space into a spectrum valued additive functor. To this end recall from Corollary~\reftwo{corollary:universal} that the forgetful functor
\[
\Omega^\infty \colon \Funadd(\Catp,\Spa) \longrightarrow \Funadd(\Catp,\Grp_{\Einf})
\]
admits a left adjoint $\CCob$.

\begin{definition}
\label{definition:gw-spectrum}%
We define the \defi{Grothendieck-Witt spectrum} functor $\GW \colon \Catp \rightarrow \Spa$ by 
\[
\GW(\C,\QF) =\CCob^{\GWspace}(\C,\QF),
\]
and denote by $\GW_i(\C,\QF)$ its homotopy groups for $i \in \mathbb Z$.
\end{definition}

We will see in Corollary~\reftwo{corollary:algGMTW} below, that for $i \geq 0$ this conforms with Definition~\reftwo{definition:gw-space}.
Again, we list the properties that are immediate from the results of the previous section. As a reformulation of the definition we find:

\begin{corollary}
\label{corollary:univ-GW-spectrum}%
The functor $\GW \colon \Catp \rightarrow \Spa$ is additive, and it is the initial such functor equipped a transformation $\Poinc \Rightarrow \Omega^\infty \GW$ of functors $\Catp \rightarrow \Sps$.%
\end{corollary}

By adjunction $\GW$ is also the initial additive functor $\GW \colon \Catp \rightarrow \Spa$ under $\SS[\Poinc]$. Moreover, we show in Corollary~\reftwo{corollary:gwadd} below, that $\GW$ is Verdier-localising and not just additive.

\begin{proof}
By Corollary~\reftwo{corollary:Cob=Susp}, the functor $\GW$ is the initial additive functor to spectra equipped with a transformation $\GWspace \Rightarrow \Omega^\infty\GW$ of $\Einf$-groups. By the universal property of $\GWspace$ established in Observation~\reftwo{observation:univ-GW-space}, the functor $\GW$ is therefore also the initial additive functor to spectra with a transformation $\Poinc \Rightarrow \Omega^\infty \GW$ of $\Einf$-monoids. 
\end{proof}

Next, we identify the spaces $\Omega^{\infty-i}\GW(\C,\QF)$. To this end, recall the $i$-fold simplicial space 
\[
\Cob_i(\C,\QF) = \Poinc\Q^{(i)}(\C,\QF\qshift{i}),
\]
given by the iterated hermitian $\Q$-construction from Definition~\reftwo{definition:hermitianQ}. These model the extended cobordism categories of $(\C,\QF)$ and by Proposition~\reftwo{proposition:positive-om} form a positive $\Omega$-spectrum $\CCob^{\Poinc}(\C,\QF)$. The natural map
\[
\CCob^{\Poinc}(\C,\QF) \longrightarrow \CCob^{\GWspace}(\C,\QF) = \GW(\C,\QF)
\]
exhibits the right hand side as the spectrification of the left. From Proposition~\reftwo{proposition:positive-om}, we also find:

\begin{corollary}
\label{corollary:algGMTW}%
For any Poincaré $\infty$-category $(\C,\QF)$, there are canonical equivalences
\[
\GWspace(\C,\QF) \simeq \Omega^\infty\GW(\C,\QF) \quad \text{and} \quad 
|\Cob_i(\C,\QF)| \simeq \Omega^{\infty-i}\GW(\C,\QF)
\]
for any $i \geq 1$, that are natural in $(\C,\QF)$. In particular, we obtain isomorphisms
\[
\pi_i\GW(\C,\QF) \cong \GW_i(\C,\QF)
\]
for all $i \geq 0$.
\end{corollary}

\begin{corollary}
The functor $\GW\colon \Catp\to \Spa$ preserves filtered colimits.
\end{corollary}

\begin{proof}
The functor $\CCob^{\Poinc}$ from $\Catp$ to the $\infty$-category of prespectra does, by our description through the iterated $\Q$-construction, and the spectrification functor preserves filtered colimits.
\end{proof}

\begin{remark}
\label{remark:geometric}%
We propose to view the equivalences 
\[
|\Cob_i(\C,\QF)| \simeq \Omega^{\infty-i}\GW(\C,\QF)
\]
for $i \geq 1$ as a close analogue of the equivalence
\[
|\Cob^i_d| \simeq \Omega^{\infty-i}\MTSO(d),
\]
established by Galatius, Madsen, Tillmann and Weiss for $i = 1$, and Bökstedt and Madsen in general \cite{GMTW, BokstedtMadsen}. In particular, the sequence of spectra $\GW(\C,\QF\qshift{-d})$ can be considered as an algebraic analogue of the Madsen-Tillmann-spectra $\MTSO(d)$. 

Of course, our arguments so far correspond only to the statement that the higher cobordism categories $\Cob_d^n$ deloop one another, i.e.\ that $|\Cob_d^n| \simeq \Omega|\Cob_d^{n+1}|$ for $n \geq 1$. The identification of the resulting spectrum via a Pontryagin-Thom construction has no direct analogue in our work. We describe the homotopy type of $\GW(\C,\QF)$ by different means in Corollary~\reftwo{corollary:tate-square-L} below. 

We shall make this connection more than an analogy in future work by promoting the association of its cochains or stable normal bundle to a manifold into functors
\[
\sigma \colon \Cob_d \longrightarrow \Cob\big((\Spa/\BSO(d))^\cp, \QF^\vis_{-\gamma_d}\big) \longrightarrow \Cob\big(\Dperf(\ZZ), (\QF^\sym)\qshift{-d}\big)
\]
from the geometric to our algebraic cobordism categories. The Grothendieck-Witt spectrum of the middle term has already appeared in manifold topology, see \S\reftwo{subsection:Weiss-Williams}, and we expect the comma $\infty$-category of the composite functor over $0$ to be closely related to the $\infty$-category $\Cob^\mathcal L_{2n+1}$ from \cite{Hebestreit-Perlmutter} for $d = 2n+1$.
\end{remark}

Just as the negative homotopy groups of the Madsen-Tillmann spectra are given by the cobordism groups, so are the negative homotopy groups of the Grothendieck-Witt spectrum. From Definition~\refone{definition:L-groups}, we recall:

\begin{definition}
\label{definition:L-group-the-fourth}%
The \defi{$\L$-group} $\L_0(\C,\QF)$ of a Poincaré $\infty$-category is defined as the quotient monoid of $\pi_0\Poinc(\C,\QF)$ by the submonoid of forms $(X,q)$ admitting a Lagrangian $W \rightarrow X$.
\end{definition}

Also, $\L_0(\C,\QF)$ really is a group: We showed in Proposition~\reftwo{corollary:components-of-cob-II} that there is a canonical isomorphism
\[
\pi_0|\Cob(\C,\QF\qshift{-1})| \cong \L_0(\C,\QF),
\]
and consequently, we get 
\[
[x,-q] + [x,q] = 0
\]
in $\L_0(\C,\QF)$ from Proposition~\reftwo{proposition:components-of-cob-I}. In other words, $\L_0(\C,\QF)$ is the cobordism group of Poincaré forms in $(\C,\QF)$ and inverses are given be reversing the orientation. From Proposition~\reftwo{proposition:pinegCobF}, we obtain:

\begin{corollary}
For $i > 0$ there are canonical isomorphisms
\[
\pi_{-i}\GW(\C,\QF) \cong \L_0(\C,\QF\qshift{i})
\]
natural in the Poincaré \(\infty\)-category $(\C,\QF)$.
\end{corollary}

\subsection{The Bott-Genauer sequence and Karoubi's fundamental theorem}
\label{subsection:Genauer-Bott-Karoubi}%

In the present section, we analyse the behaviour of the metabolic Poincaré-Verdier sequence
\[
(\C,\QF\qshift{-1}) \longrightarrow \Met(\C,\QF) \xrightarrow{\met} (\C,\QF)
\]
under the Grothendieck-Witt functor. From Example~\reftwo{example:CobHyp-equal-Span} and Corollary~\reftwo{corollary:met-hyp}, we obtain:

\begin{corollary}
\label{corollary:GWHyp}%
For a Poincar\'e $\infty$-category $(\C,\QF)$, the functors $\lag \colon \Met(\C,\QF) \leftrightarrow \Hyp(\C) \colon \ilag$ induce inverse equivalences
$
\GW(\Met(\C,\QF)) \simeq \GW(\Hyp(\C))
$
and for a stable $\infty$-category $\D$, we have a canonical equivalence
$
\GW(\Hyp(\D)) \simeq \K(\D).
$
In particular, 
we find $\GW^\hyp \simeq \K.$
\end{corollary}

Consequently, applying $\GW$ to the metabolic sequence gives a fibre sequence
\[
\GW(\C,\QF\qshift{-1}) \xrightarrow{\fgt} \K(\C) \xrightarrow\hyp \GW(\C,\QF),
\]
of spectra, which we call the \defi{Bott-Genauer sequence}. It is a general version of the Bott-sequence appearing for example in \cite{schlichting-derived}*{Section 6}. 

\begin{remark}
\label{remark:geometric-2}%
We chose the present terminology to highlight the analogy with the fibre sequence
\[
\MTSO(d+1) \lrar \SS[\BSO(d+1)] \lrar \MTSO(d)
\]
originally appearing in \cite{GMTW}*{Section 3}, which complemented Genauer's theorem in \cite{Genauer}*{Section 6} that
\[
|\Cob^\partial_d| \simeq \Omega^{\infty-1}\SS[\BSO(d)].
\]
In particular, in the Bott-Genauer sequence the algebraic $\K$-theory spectrum really arises via the metabolic category, encoding objects with boundary, rather than the hyperbolic category. From this perspective, the connectivity of the algebraic $\K$-theory spectrum corresponds to the fact that the bordism groups of manifolds with boundary vanish.
\end{remark}

Finally, we observe that the Bott-Genauer sequence gives a vast extension of Karoubi's fundamental theorem: Following Karoubi and Schlichting \cite{Karoubi-Le-theoreme-fondamental, schlichting-derived}, we define functors 
\[
\U(\C,\QF) = 
\fib(\K(\C) \stackrel{\hyp}{\longrightarrow} \GW(\C,\QF))
\quad\text{and}\quad 
\V(\C,\QF) = 
\fib(\GW(\C,\QF) \stackrel{\fgt}{\longrightarrow} \K(\C)).
\]
Karoubi's fundamental theorem \cite{Karoubi-Le-theoreme-fondamental}*{p. 260} compares these functors in the setting of discrete rings with involution. In the setting of Poincaré $\infty$-categories, this statement is a direct consequence of the Bott-Genauer sequence (we will specialise this abstract version to discrete rings in Corollary~\reftwo{corollary:karoubi-fundamental}
below).

\begin{corollary}[Karoubi's fundamental theorem]
\label{corollary:karoubi-fundamental-theorem}%
There is a canonical equivalence
\[
\U(\C,\QF\qshift{2}) \simeq \SS^1 \otimes \V(\C,\QF)
\]
natural in the Poincaré \(\infty\)-category $(\C,\QF)$.
\end{corollary}
\begin{proof}
Simply note that the Bott-Genauer sequence identifies both sides with $\GW(\C,\QF\qshift{1})$.
\end{proof}

Now, we showed in \refone{proposition:shiftofqingeneral} that the Poincaré $\infty$-category $(\C,\QF\qshift{2})$ can be described in another way: By \refone{corollary:quadratic-genuine} $\QF$ canonically lifts to a functor $\wtl{\QF} \colon \C\op \rightarrow \Spagc$ with $\wtl{\QF}^\gC \simeq \QF$, and $\Omega^k \colon \C \rightarrow \C$ upgrades to an equivalence
\[
(\C,\QF\qshift{2k}) \simeq (\C,(\SS^{k-k\sigma} \otimes \wtl{\QF})^\gC),
\]
where $\SS^{\sigma}$ is the suspension spectrum of the sign representation of $\Ct$ on $\mathbb R$. While more complicated at first glance, the right hand side is actually much easier to analyse directly, and corresponds to the idea of "inserting a sign" into a Poincaré structure.

Let us spell this out for the Grothendieck-Witt theory of rings. Recall that these are integrated into our set-up via their derived categories of modules. More generally, consider an $\Einf$-ring $k$, and
an $\Eone$-algebra $R$ over $k$ and an invertible $k$-module with genuine involution $(M, \alpha \colon N \rightarrow M^\tC)$ over $R$. Then by Proposition~\refone{proposition:general-equivalence-of-poincare-infty-categories}, the general equivalence above becomes 
\[
(\Modp{R},(\QF^\alpha_M)\qshift{2n}) \simeq (\Modp{R},\QF^{\SS^n \otimes \alpha}_{\SS^{n-n\sigma} \otimes M}).
\]
Now, if $k$ is equipped with a $2\sigma$-orientation, i.e.\ a factorisation of its unit map as
\[
\SS \longrightarrow (\SS^{2-2\sigma} \otimes k)^\hC \xrightarrow{\fgt} k,
\]
then we learn from \refone{remarks:basic} that the module with involution $\SS^{n - n\sigma} \otimes M$ only depends on the parity of $n$ modulo $2$ up to a canonical equivalence. We denote the common value for odd $n$ by $-M$, so that we get
\[
(\Modp{R},(\QF^\alpha_M)\qshift{2}) \simeq (\Modp{R},\QF^{\SS^1 \otimes \alpha}_{-M})
\]
as a special case of the above. This situation occurs most importantly, whenever $k$ is discrete, e.g.\ $k = \mathrm H\ZZ$, and if $R$ and $M$ are discrete as well $-M$ really is given by changing its involution by a sign. The situation also occurs more generally, however, when $k$ is complex oriented (e.g.\ even periodic) or $2 \in \pi_0(k)^\times$, see \refone{example:signperiodicmodules} for details.

A $2\sigma$-orientation on $k$ in particular produces an equivalence $(-M)^\tC \simeq \SS^1 \otimes M^\tC$, so we obtain equivalences 
\[
(\Modp{R},(\QF^\sym_M)\qshift{2}) \simeq (\Modp{R},\QF^\sym_{-M}), \qquad (\Modp{R},(\QF^\qdr_M)\qshift{2}) \simeq (\Modp{R},\QF^\qdr_{-M})
\]
and, whenever $R$ is furthermore connective, 
\[
(\Modp{R}, (\Qgen m M)\qshift{2}) \simeq (\Modp{R},\Qgen{m+1}{-M}).
\]
Also note that if $k$ admits a $\sigma$-orientation, e.g.\ it is discrete of characteristic $2$ or more generally real oriented, then we even find $M \simeq -M$. 

Recall furthermore, that for $\mathrm c \in \K_0(R) = \K_0(\Modp{R})$ a subgroup, we denote by $\Mod^\mathrm{c}(R) \subseteq \Modp{R}$ the full subcategory spanned by those $R$-modules $X$ with $[X] \in \mathrm c$, the most interesting special cases being 
\[
\Mod^{\K_0(S)}(\mathrm HS) = \Dperf(S) \quad \text{and} \quad \Mod^{\langle \mathrm HS \rangle}({\mathrm HS}) = \Dfree(S)
\]
for $S$ an ordinary ring. We shall need to assume that $\mathrm c$ is closed under the involution induced by $M$. This is clearly always true in the former case, and in the latter amounts to $[M] = \pm [S] \in \K_0(S)$.

\begin{corollary}
\label{corollary:karoubi-fundamental}%
For $R$ an $\Eone$-algebra over an $\Einf$-ring $k$, which carries a $2\sigma$-orientation, $M$ an invertible $k$-module with involution over $R$, and $\mathrm c \subseteq \K_0(R)$ a subgroup closed under the involution induced by $M$, there are canonical equivalences
\[
\U(\Modc{\mathrm c}{R},\QF^\qdr_{-M}) \simeq \SS^1 \otimes \V(\Modc{\mathrm c}{R},\QF^\qdr_{M}) \quad \text{and} \quad \U(\Modc{\mathrm c}{R},\QF^\sym_{-M}) \simeq \SS^1 \otimes \V(\Modc{\mathrm c}{R},\QF^\sym_{M})
\]
and if $R$ is furthermore connective, then also
\[
\U(\Modc{\mathrm c}{R},\Qgen {m+1}{-M}) \simeq \SS^1 \otimes \V(\Modc{\mathrm c}{R},\Qgen{m}{M})
\]
for arbitrary $m \in \ZZ$. %
\end{corollary}

Now, these equivalence can be specialised further to a discrete ring $R$, $c = \K_0(R)$ and a discrete invertible module with involution over $R$, in which case one can take $k = \mathbb Z$ together with its canonical $2\sigma$-orientation. Choosing either $m=1$ or $m=2$, we then obtain the following extension of Karoubi's fundamental theorem:

\begin{corollary}
\label{corollary:karoubi-even-more-fundamental}%
For a discrete ring $R$ and a discrete invertible module with involution $M$ over $R$, there are canonical equivalences
\[
 \U(\Dperf(R),\QF^\gq_{-M}) \simeq \SS^1 \otimes \V(\Dperf(R),\QF^\gev_{M}),
 \quad 
 \U(\Dperf(R),\QF^\gev_{-M}) \simeq \SS^1 \otimes \V(\Dperf(R),\QF^\gs_{M}),
\]
and
\[
 \U(\Dperf(R),\QF^\qdr_{-M}) \simeq \SS^1 \otimes \V(\Dperf(R),\QF^\qdr_{M}),
 \quad 
 \U(\Dperf(R),\QF^\sym_{-M}) \simeq \SS^1 \otimes \V(\Dperf(R),\QF^\sym_{M}).
\]
\end{corollary}

Given the comparisons in Appendix~\reftwo{appendix:AppIIB}, all of these equivalences collapse into the classical formulation of Karoubi's fundamental theorem upon restricting to discrete rings in which $2$ is invertible; if $2$ is not assumed invertible they are, however, distinct. We will explore their uses for discrete rings in the third paper of this series.
\medskip

Let us also apply the fundamental theorem \reftwo{corollary:karoubi-fundamental-theorem} to the case of a form parameter, as originally defined by Bak \cite{bak-form-parameter} and generalised by Schlichting \cite{SchlichtinghigherI}. Recall from \refone{subsection:discrete-rings} that given a discrete ring $R$ and a discrete invertible module with involution $(M,\sigma)$ over $R$, a form parameter $\fpm$ on $M$ is the data of an abelian group $\fpmg$, equipped with an action by the multiplicative monoid of $R$ through group homomorphisms, and two $R$-equivariant homomorphisms
\[
M_{\Ct} \xrightarrow{\tau} \fpmg \xrightarrow{\rho} M^{\Ct}
\]
whose composition is the norm map (i.e.\ the map induced by $\id_M + \sigma$) and such that
\[
(a+b)x= ax + \tau((a\otimes b)\rho(x)) + bx
\]
 for all $a,b \in R$ and $x \in Q$. As explained in \S\refone{subsection:discrete-rings}, the $2$-polynomial functor $\QF^{\fpm}_{\proj}\colon\Proj(R)\op \to \Ab$ sending $P$ to its group of $\fpm$-hermitian forms, followed by the canonical Eilenberg-Mac Lane inclusion $\Ab \rightarrow \Spa$ extends (essentially uniquely) by non-abelian derivation to a hermitian structure $\QF^\gfpm_M\colon \Dperf(R)\op \to \Spa$. 

Any form parameter $(\fpmg,\tau,\rho)$ on $M$ admits a dual form parameter $\dfpm$ %
\[
(-M)_{\Ct} \xrightarrow{\mathrm{pr}} M/\fpmg \xrightarrow{\id - \sigma} (-M)^{\Ct}
\]
on $(-M)$ in which the first map is surjective. For example we have $\check{\pm\s} = \mp \ev$ and $\check{\pm \ev} = \mp\qdr$, in what is hopefully evident notation. We showed in Proposition~\refone{proposition:form-parameter-shift} that the loop functor refines to an equivalence of Poincaré \(\infty\)-categories
\[
\big(\Dperf(R), (\QF_{M}^{\gfpm})\qshift{2}\big) \longrightarrow \big(\Dperf(R), \QF^{\gdfpm}_{-M}\big),
\]
whenever $\rho$ is injective in the original form parameter. Applying Corollary~\reftwo{corollary:karoubi-fundamental-theorem} to this equivalence immediately implies: 
\begin{corollary}
\label{corollary:karoubi-periodicity-form-parameters}%
For a discrete ring $R$, a discrete invertible module with involution $M$, $c \subseteq \K_0(R)$ a subgroup closed under the involution induced by $M$ and a form parameter $\fpm=(\fpmg,\tau,\rho)$ on $M$ with $\rho$ injective and dual $\dfpm$, there is a canonical equivalence
\[
\U(\D^c(R),\QF^{\g\dfpm}_{-M}) \simeq \SS^1 \otimes \V(\D^c(R),\QF^\gfpm_M).
\]
\end{corollary}
Together with \cite{comparison}*{Theorem A}, which identifies the Grothendieck-Witt space of $(\Dperf(R),\QF^\gfpm_M)$ considered in the present paper with the Grothendieck-Witt space considered in \cite{karoubi-periodicity} (which is defined as the group completion of the $\Einf$-monoid of the corresponding Poincaré forms on projective modules), this proves Conjectures 1\ and 2\ in \S3.4 and \S4.3 of loc.\ cit, respectively. Note for translation purposes that Karoubi's hermitian and quadratic modules for the form parameter $\fpm$ correspond to Poincaré forms for the form parameters $\fpm$ and $\dfpm$, respectively, as explained in \cite{SchlichtinghigherI}*{Example~3.9}.

\subsection{L-theory and the fundamental fibre square}
\label{subsection:L+tate}%

In the present section, we prove our main result on the homotopy type of the Grothendieck-Witt spectrum. In \S\reftwo{subsection:bordification}, we studied the bordification of an additive functor $\F \colon \Catp \rightarrow \Spa$ and in Corollary~\reftwo{corollary:tate-square-L}, we produced a cartesian square reconstructing $\F$ from its hyperbolisation $\F^\hyp$ and its bordification $\F^\bord$. In the previous subsection, we obtained an equivalence $\GW^\hyp \simeq \K$, and in the present section, we show $\GW^\bord \simeq \L$. To set the stage, recall the $\rho$-construction from Definition~\reftwo{definition:rho-construction}.

\begin{definition}
The \defi{$\L$-theory space} is the functor $\Catp \rightarrow \Sps$ given by
\[
\Lspace(\C,\QF) = |\Poinc\rho(\C,\QF)|
\]
obtained by applying $\Poinc$ to the $\rho$-construction.
\end{definition}

Since $\rho_0(\C,\QF) = (\C,\QF)$, there is a canonical map 
$
\Poinc(\C,\QF) \rightarrow \Lspace(\C,\QF).
$
and by construction, the $1$-skeleta of the $\rho$ and $\Q$ construction agree, so, from Corollary~\reftwo{corollary:components-of-cob-II}, we find that the natural map  
$
\pi_0\Poinc(\C,\QF) \to \pi_0 \Lspace(\C,\QF)
$
descends to an isomorphism 
$
\L_0(\C,\QF) \to \pi_0 \Lspace(\C,\QF)
$
for all Poincaré $\infty$-categories $(\C,\QF)$. 

But much more is true: Generalising a classical result of Ranicki, Lurie showed in \cite{Lurie-L-theory}*{Lecture 7, Theorem 9} that there are canonical isomorphisms
\[
\pi_i\Lspace(\C,\QF) = \L_0(\C,\QF\qshift{-i})
\]
for all $i \geq 0$. While analogous to our results on bordifications, this is more difficult and fundamentally rests on the fact that $\Poinc\rho(\C,\QF)$ is a Kan simplicial space. In fact:

\begin{theorem}
\label{theorem:L-is-localising}%
Given a Poincaré-Verdier sequence $(\C,\QF) \rightarrow (\D,\QFD) \rightarrow (\E,\QFE)$, the map $\Poinc\rho(\D,\QFD) \rightarrow \Poinc\rho(\E,\QFE)$ is a Kan fibration of simplicial spaces with fibre $\Poinc\rho(\C,\QF)$. 
In particular, the functor $\Lspace \colon \Catp \rightarrow \Sps$ is Verdier-localising, bordism invariant, and it also preserves filtered colimits.
\end{theorem}

The above identification of homotopy groups is then a consequence of Proposition~\reftwo{proposition:bordism-invariant-sig}. The result itself is the main content of \cite{Lurie-L-theory}*{Lectures 8 \& Lecture 9}; we give the proof here for completeness' sake. It rests on the following lemma:

\begin{lemma}
\label{lemma:hermitian-lifts-verdier}%
Given a Poincaré-Verdier projection $(\D,\QFD) \xrightarrow{(p,\eta)} (\E,\QFE)$, an object $X \in \D$, a map $f \colon Y \rightarrow p(X)$ in $\E$ and a diagram
\[
\begin{tikzcd}
K \ar[rr]\ar[d] && * \ar[d,"q"] \\
\Omega^\infty \QFD(X) \ar[r,"\eta"] & \Omega^\infty\QFE(p(X)) \ar[r,"f^*"] & \Omega^\infty\QFE(Y)
\end{tikzcd}
\]
with $K \in \Sps^\cp$, there exists an arrow $g \colon Z \rightarrow X$ in $\D$ lifting $f$ together with a lift
\[
\begin{tikzcd}
K \ar[r] \ar[d] & * \ar[d,"r"]\\ 
\Omega^\infty \QFD(X) \ar[r,"g^*"] & \Omega^\infty\QFD(Z)
\end{tikzcd}
\]
of the original rectangle.
\end{lemma}
\begin{proof}
Since $p$ is essentially surjective by Corollary~\reftwo{corollary:verdier-surjective}, there exists a $V$ in $\D$ with $p(V) \simeq Y$, and applying \cite{NS}*{Theorem I.3.3 ii)}, we can then modify $V$ to find $h \colon W \rightarrow Y$ lifting $f$. From Remark~\reftwo{remark:compute}, setting $\C = \ker(p)$, we furthermore find
\[
\QFE(Y) \simeq \colim_{C \in \C_{W/}} \QFD(\fib(W \rightarrow C))
\]
so, putting $U = \fib(W \rightarrow C)$ for an appropriate $C$, we find a lift $s \in \Omega^\infty\QFD(U)$ lifting $q$, and the composite $U \rightarrow W \rightarrow Y$ still lifts $f$. To find a lift of the homotopy of maps $K \rightarrow \Omega^\infty\QFE(Y)$, note that the colimit above is filtered so, since $K$ is assumed compact, we also have
\[
\Hom_\Sps(K,\Omega^{\infty+1}\QFE(Y)) \simeq \colim_{C' \in \C_{U/}} \Hom_\Sps(K,\Omega^{\infty+1}\QFD(\fib(U \rightarrow C')),
\]
which for appropriate $C'$ yields all the desired data for $Z = \fib(U \rightarrow C')$ and $g$ the composite $Z \rightarrow U \rightarrow Y$.
\end{proof}

\begin{proof}[Proof of Theorem~\reftwo{theorem:L-is-localising}]
We need to show that each solid diagram
\[
\begin{tikzcd}
\Lambda_i^n \ar[r] \ar[d] & \rho(\D,\QFD) \ar[d]\\ 
\Delta^n \ar[r] & \rho(\E,\QFE)
\end{tikzcd}
\]
admits a dotted filler up to homotopy; here we regard $\Delta^n$ and $\Lambda_i^n$ as simplicial spaces via the inclusion $\mathrm{Set} \subset \Sps$.

To unwind this, recall that $\rho_n(\D,\QFD) = (\D,\QFD)^{\T_n}$, where $\T_n = \mathcal P_0([n])\op$ is the opposite of the barycentric subdivision $\sd(\Delta^n)$ of $\Delta^n$. Denote then by $H_n^i \subseteq \T_n$ the opposite of the subdivision of the $i$-horn, i.e.\ the collection of subsets missing an element besides $i$. Then the lifting problem above translates to showing that the canonical map
\[
\Poinc\left((\D,\QFD)^{\T_n}\right) \longrightarrow \Poinc\left((\E,\QFE)^{\T_n}\right) \times_{\Poinc\left((\E,\QFE)^{H^i_n}\right)} \Poinc\left((\D,\QFD)^{H^i_n}\right)
\]
is surjective on $\pi_0$.
To this end, we first show the corresponding statement on spaces of hermitian objects, and then explain how to adapt a lift in $\spsforms\left((\D,\QFD)^{\T_n}\right)$ to a Poincaré one, provided its images in $\spsforms\left((\E,\QFE)^{\T_n}\right)$ and $\spsforms\left((\D,\QFD)^{H^i_n}\right)$ are Poincaré. The first claim even holds for boundary inclusions instead of horn inclusions, so denote by $B_n$ the opposite of the subdivision of $\partial \Delta^n$ and consider hermitian objects $(F \colon \T_n \rightarrow \E,q)$ and $(G \colon B_n \rightarrow \D,r)$ and an equivalence between their images in $\spsforms\left((\E,\QFE)^{B_n}\right)$. Let $X \in \D$ be the limit of $G$. By construction, there is then a canonical map $f \colon Y \rightarrow p(X)$, where $Y = F(b)$ with $b$ the barycentric vertex $[n]$ in $\T_n$. Furthermore, regarding $r \in \QFD^{B_n}(G)$ as a map $r \colon \ast \rightarrow \lim_{B_n\op} \Omega^\infty\QFD \circ G\op$ it is adjoint to a transformation $\const_\ast \Rightarrow \Omega^\infty\QFD \circ G\op$, which gives rise to a map
\[
|B_n\op| \simeq \colim_{B_n\op} \ast \longrightarrow \colim_{B_n\op} \Omega^\infty\QFD \circ G\op \longrightarrow \Omega^\infty \QFD(\lim_{B_n} G) = \Omega^\infty\QFD(X)
\]
whose composition down to $\Omega^\infty\QFE(Y)$ is canonically identified with the constant map with value $q \in \Omega^\infty\QFE^{\T_n}(F) \simeq \Omega^\infty\QFE(F(b))$, since $b = [n]$ is initial in $\T_n\op$, so $|\T_n\op| \simeq \ast$.

We can therefore apply the previous lemma to obtain a lift $g \colon Z \rightarrow X$ of $f$, together with a lift $s \in \Omega^\infty \QFD(Z)$ of $q$ and an identification of the composite 
\[
|B_n\op| \longrightarrow \Omega^\infty\QFD(X) \xrightarrow{g^*} \Omega^\infty\QFD(Z)
\]
with the constant map on $s$, that lifts the identification above. Since $\T_n$ is the cone on $B_n$, the map $g$ precisely defines an extension of $G$ to a map $\widetilde G \colon \T_n\op \rightarrow \D$, on which $s$ defines a hermitian form, and the remainder of the data produced bears witness to $(\widetilde G,s)$ being a lift as desired.

For the second step, we need to modify a hermitian lift $(\widetilde G,r) \in \spsforms\left((\D,\QFD)^{\T_n}\right)$ of $(F,q) \in \Poinc\left((\E,\QFE)^{\T_n}\right) \in $ and $(G,s) \in \Poinc\big((\D,\QFD)^{H_i^n}\big)$ into a Poincaré lift. This is achieved by performing surgery as follows: The algebraic Thom construction from Corollary~\refone{corollary:algebraic-thom-iso} gives an equivalence 
\[
\spsforms\left((\D,\QFD)^{\T_n}\right) \simeq \Poinc(\Met\left((\D,\QFD\qshift{1})^{\T_n}\right)
\]
refining the map taking $(\widetilde G,s)$ to 
\begin{equation}\tag{$\ast$}
\label{equation:arrow-in-proof-of-L}%
\Dual_{\QFD^{\T_n}}\big(\widetilde G\big) \longrightarrow \cof\big(\widetilde G \xrightarrow{s_\sharp}  \Dual_{\QFD^{\T_n}}\big(\widetilde G\big)\big).
\end{equation}
In particular, the Poincaré objects in $(\D,\QFD)^{\T_n}$ correspond precisely to those arrows with vanishing target (the target is the boundary of $(\widetilde G,s)$ in the sense of Definition~\reftwo{definition:Ranicki-boundary} below). Since $(F,q)$ and $(G,r)$ and the boundary maps in the $\rho$-construction are Poincaré (see the discussion before Definition~\reftwo{definition:rho-construction}), it follows that the target in our case already lies in the kernels of both
\[
\D^{\T_n} \longrightarrow \D^{H_i^n} \quad \text{and} \quad \D^{\T_n} \longrightarrow \E^{\T_n}.
\]
We claim that the intersection of these kernels is equivalent to $\Met(\C,\QF\qshift{1-n})$ as a Poincaré $\infty$-category. This is clear on underlying $\infty$-categories, and follows for the hermitian structures from the iterative formulae for limits of cubical diagrams, i.e.
\[
\lim_{\T_n\op} X \simeq X([n]) \times_{\lim\limits_{\T_{n-1}\op} X \circ (- \cup n)} \lim_{\T_{n-1}\op} X
\]
which is easily verified using \cite{HTT}*{Corollary 4.2.3.10} by decomposing $\T_n$ as the pushout of $\T_{n-1}$ and $\T_n\setminus\{0,...,n-1\}$ over their intersection. We thus find that the cofibre of $s_\sharp$ admits a Lagrangian $L$, since objects in metabolic Poincaré \(\infty\)-categories are canonically metabolic by Remark~\refone{remark:comonad}. We can thus perform surgery on \eqreftwo{equation:arrow-in-proof-of-L} with the surgery datum $0 \rightarrow L$, see Proposition~\reftwo{proposition:surgequiv}. The resulting arrow has vanishing target, and by design the surgery changes neither the image in $\E^{\T_n}$ nor the restriction to $\D^{H_i^n}$. Translating back along the algebraic Thom construction thus provides the desired Poincaré lift of $(F,q)$ and $(G,q)$.

To deduce the remaining claims, note that the statement about the fibre is immediate from both cotensors and $\Poinc$ preserving limits. That $\Lspace$ is Verdier-localising now follows, since colimits of simplicial fibre sequences with second map a Kan fibration are again fibre sequences, see e.g.\ \cite{SAG}*{Theorem A.5.4.1}. 

To obtain bordism invariance, one can 
proceed by observing that on account of the Kan property the $i$-th homotopy groups of $\L(\C,\QF) = |\Poinc\rho(\C,\QF)|$ can be described as the quotient of 
\[
\pi_0\fib\left(\Hom_{\sSps}(\Delta^i, \Poinc\rho(\C,\QF)) \longrightarrow \Hom_{\sSps}(\partial \Delta^i, \Poinc\rho(\C,\QF))\right)
\]
by the equivalence relation generated by a pair of such elements admitting an extension to 
\[
\pi_0\fib\left(\Hom_{\sSps}(\Delta^1 \times \Delta^i, \Poinc\rho(\C,\QF)) \rightarrow \Hom_{\sSps}(\Delta^1 \times \partial \Delta^i, \Poinc\rho(\C,\QF))\right).
\]
This quotient is readily checked to be exactly $\L_0(\C,\QF\qshift{-i})$. Since these $\L$-groups evidently take bordism equivalences to isomorphisms, bordism invariance of $\L$ follows. Finally, that $\Lspace$ preserves filtered colimits follows by the same argument as for $\GWspace$ in Corollary \reftwo{corollary:filtcolimpresgwspace} above.
\end{proof}

It now follows from Theorem~\reftwo{theorem:lift} that $\Lspace$ admits an essentially unique lift to a functor with values in spectra.

\begin{definition}
We define the \defi{$\L$-theory spectrum} $\L \colon \Catp \rightarrow \Spa$ by
$$ \L(\C,\QF) = \CCob^{\Lspace}(\C,\QF)$$
with $(\C,\QF)$ a Poincaré \(\infty\)-category, and denote by $\L_i(\C,\QF)$ its homotopy groups for $i \in \mathbb Z$.
\end{definition}

The following is then immediate:

\begin{corollary}
\label{corollary:Ladd}%
The functor $\L \colon \Catp \rightarrow \Spa$ is Verdier-localising, bordism invariant, and preserves filtered colimits. 
\end{corollary}

The above definition of $\L$-groups agrees with Definition~\reftwo{definition:L-group-the-fourth}, since from Proposition~\reftwo{proposition:positive-om} and Proposition~\reftwo{proposition:bordism-invariant-sig}, we obtain:

\begin{corollary}
There are canonical equivalences
$
\Omega^{\infty-i}\L(\C,\QF) \simeq \Lspace(\C,\QF\qshift{i})
$
for all $i \in \ZZ$. In particular, there are isomorphisms
$
\pi_i\L(\C,\QF) \cong \L_0(\C,\QF\qshift{-i})
$
also for negative $i$.
\end{corollary} 

In fact, the equivalence
\[
\L(\C,\QF) \simeq 
\left[\Lspace(\C,\QF), \Lspace(\C,\QF\qshift{1}), \Lspace(\C,\QF\qshift{2}), \dots\right]
\]
with structure maps arising from Proposition~\reftwo{proposition:bordism-invariant-sig} is a direct generalisation of the classical definition of $\L$-theory spectra due to Ranicki; see for example \cite{Ranickiblue}*{Section 13}, and it is rather more elegant than our definition which iterates the $\Q$-construction on top of the $\rho$-construction. %

One can also directly describe the boundary operator of the long exact sequence on the $\L$-groups of a Poincaré-Verdier sequence. To that end, we make the following definition.
\begin{definition}
\label{definition:Ranicki-boundary}%
Given a Poincaré \(\infty\)-category $(\C, \QF)$ and a hermitian object $(X,q)\in \spsforms(\C, \QF)$, the \defi{boundary} of $(X,q)$ is the Poincaré object $\partial(X,q)\in \Poinc(\C, \QF\qshift{1})$ obtained as the result of surgery on $(X \rightarrow 0,q) \in \Surg_0(\C,\QF\qshift{1})$.
\end{definition}

Note that by the discussion preceding Proposition~\reftwo{proposition:surgequiv}, the object underlying $\partial(X,q)$ is given by the cofibre of $q_\sharp \colon X \rightarrow \Dual_{\QF}X$. 
\begin{proposition}
\label{proposition:boundary-L}%
Given a Poincaré-Verdier sequence 
$
(\C,\QF) \to (\D,\QFD) \to (\E,\QFE),
$
the boundary operator $\L_i(\E,\QFE) \rightarrow \L_{i-1}(\C,\QF)$ of the resulting long exact sequence takes a Poincaré object $(X,q) \in \Poinc(\E,\QFE\qshift{-i})$ to $\partial(Y, q') \in \Poinc(\C,\QF\qshift{1-i})$, where $(Y,q')\in \spsforms(\D,\QFD\qshift{-i})$ is any lift of $(X,q)$. 
\end{proposition}

In particular, the proposition asserts that such a hermitian lift of $X$ can always be found, and its image in $\L_{i-1}(\C,\QF)$ is the obstruction against finding a Poincaré lift of $X$.

\begin{proof}
From Proposition~\reftwo{proposition:boundarymapGB}, we find that the inverse to the boundary isomorphism $\pi_1\L(\E,\QFE) \rightarrow \pi_0\L(\E,\QFE\qshift{-1})$ takes a Poincaré object $X$ in the target to the loop $w$ represented by $0 \leftarrow X \rightarrow 0 \in \Poinc\rho_1(\E,\QFE)$.
We now compute the map $\L_1(\E,\QFE) \rightarrow \L_0(\C,\QF)$, the case of general $i \in \ZZ$ follows by shifting the quadratic functor. That any Poincaré object $(X,q) \in \Poinc(\E,\QFE\qshift{-1})$ can be lifted to some $(Y,q') \in \spsforms(\D,\QFD\qshift{-1})$ is an application of Lemma~\reftwo{lemma:hermitian-lifts-verdier} (with $K=\emptyset$). %

Now, regarding the map $(Y \rightarrow 0,q')$ as a surgery datum in $\Surg_0(\D,\QFD)$, we can apply Proposition~\reftwo{proposition:surgequiv} to obtain a cobordism
from $0$ to the result of surgery, which is $\partial(Y, q')$. We regard this cobordism 
as an element of $\Poinc(\rho_1(\D,\QFD))$ and thus as a path in $\L(\D,\QFD)$. By construction, this path lifts the loop in $\L(\E,\QFE)$ defined by $X$ via the consideration in the first paragraph. Therefore, its endpoint $\cof(Y \rightarrow \Dual_{\QF\qshift{-1}} Y)$ represents the image of $(X,q)$ under the boundary map as claimed.
\end{proof}

Now, by construction there is a natural transformation $\Poinc \Rightarrow \Lspace \simeq \Omega^\infty \L$, which uniquely extends to a transformation 
$\bord \colon \GW \Rightarrow \L$ 
of functors $\Catp \rightarrow \Spa$ by Corollary~\reftwo{corollary:univ-GW-spectrum}. We record, see Corollary~\reftwo{corollary:filtration}:

\begin{corollary}
\label{corollary:inducedgwtol}%
Under the identifications of Theorem~\reftwo{theorem:pi1Cob} the map 
\[
\bord \colon \pi_0\GW(\C,\QF) \rightarrow \pi_0\L(\C,\QF)
\]
becomes the canonical projection $\GWgroup(\C,\QF) \rightarrow \L_0(\C,\QF)$. Similarly, for $i>0$, the induced map $\pi_{-i}\GW(\C,\QF) \rightarrow \pi_{-i}\L(\C,\QF)$ is identified with the identity of $\L_0(\C,\QF\qshift{i})$ by Proposition~\reftwo{proposition:pinegCobF}.
\end{corollary}

\begin{remark}
While the map $\bord \colon \GW(\C,\QF) \rightarrow \L(\C,\QF)$ is most easily constructed via the universal property of $\GW$, it is also easy to obtain a direct map between these spectra when defining them via the $\Q$- and $\rho$-constructions: Consider the map of cosimplicial objects $\eta \colon (\sd\Delta^n)\op \rightarrow \Twar(\Delta^n)$, that sends a non-empty subset $T \subseteq [n]$ to the pair $(\mathrm{min} T \leq \mathrm{max} T)$. It is an isomorphism in degrees $0$ and $1$, 
and $\Q_2(\C,\QF\qshift{1}) \to \rho_2(\C,\QF\qshift{1})$ is depicted as follows: 
\[
\begin{tikzcd}[row sep={5ex,between origins},column sep={5ex,between origins},
nodes={rectangle,inner sep=0.2ex}]
& & \ar[ld] \ar[rd] U & & \\
& V \ar[rd] \ar[ld] & & W \ar[rd]\ar[ld] & \\
X & & Y && Z
\end{tikzcd} 
\mapsto
\begin{tikzpicture}
[nodes={rectangle,inner sep=0.1ex,minimum height=3ex},baseline=(current bounding box.west)]
\node (y) at (90:14ex) {$Y$};
\node (x) at (195:14ex) {$X$};
\node (z) at (-15:14ex) {$Z$};
\node (center) at (barycentric cs:x=1,y=0.9,z=1) {$U$};
\node (v) at (barycentric cs:x=1,y=1,z=0) {$V$};
\node (w) at (barycentric cs:x=0,y=1,z=1) {$W$};
\node (d) at (barycentric cs:x=1,y=0,z=1) {$U$};
\draw [->,thick] (v) -- (y);
\draw [->,thick] (v) -- (x);
\draw [->,thick] (w) -- (y);
\draw [->,thick] (w) -- (z);
\draw [->,thick] (d) -- (x);
\draw [->,thick] (d) -- (z);
\draw [->,thick] (center) -- (v); 
\draw [->,thick] (center) -- (w); 
\draw [->,thick] (center) -- (d);
\end{tikzpicture}
\]

The analogous operation on manifold cobordisms takes two composable cobordisms to the $2$-ad given by the cartesian product of their composition with an interval; the ad-structure is given (after smoothing corners) by decomposing the boundary into the original two cobordisms, represented along the diagonal edges, and their composite given by the lower horizontal edge. In general then, the transformation $\eta \colon \Q \Rightarrow \rho$ regards $n$ composable $1$-ads as a special case of an $n$-ad.

Now, $\eta$ induces a map 
\[
\GWspace(\C,\QF) = \Omega |\Poinc\Q(\C,\QF\qshift{1})|  \xrightarrow{\Omega|\eta|} \Omega |\Poinc\rho(\C,\QF\qshift{1})| \xrightarrow{\partial} |\Poinc\rho(\C,\QF)| = \Lspace(\C,\QF)
\]
and thus a map $\eta \colon \GW = \CCob^{\GWspace} \Rightarrow \CCob^{\Lspace} = \L$. Using Proposition~\reftwo{proposition:boundarymapGB}, it is not difficult to check, that this map satisfies the universal property defining $\bord$. Since we shall not have to make use of that statement, we leave the details to the reader.
\end{remark}

We now turn to the main result of this section:

\begin{theorem}
\label{theorem:lisbordgw}%
The transformation $\bord$ exhibits $\L$ as the bordification of $\GW$. In particular, $\L \colon \Catp \rightarrow \Spa$ is the initial bordism invariant, additive functor equipped with a transformation $\Poinc \Rightarrow \Omega^\infty \L$ of functors $\Catp \rightarrow \Sps$.
\end{theorem}

From Theorem~\reftwo{theorem:lift}, we also find that $\Lspace \colon \Catp \rightarrow \Sps$ is the initial bordism invariant, additive functor under either $\Poinc$ or $\GWspace$.

\begin{proof}
The second statement follows from the first and 
Corollary~\reftwo{corollary:univ-GW-spectrum}. 
For the first statement, note that the canonical map $\GW \to \L$ induces a map $\GW^\bord \to \L$. We have recorded in \reftwo{corollary:filtration} that $\GW \to \GW^\bord$ induces an isomorphism on negative homotopy groups, hence by Corollary~\reftwo{corollary:inducedgwtol} the map $\GW^\bord \to \L$ also induces an isomorphism on negative homotopy groups, and thus by Proposition~\reftwo{proposition:bordism-invariant-sig} an isomorphism on all homotopy groups.
\end{proof}

\begin{remark} 
We can also employ the $\metstab$-construction for a proof of Theorem~\reftwo{theorem:lisbordgw}: By Proposition~\reftwo{proposition:bordism-invariant-sig}, the map $\bord$ factors over a map 
\[
\colim_d\SS^d \otimes \GW(\C,\QF\qshift{-d}) \lrar \L(\C,\QF).
\]
and it follows from Corollary~\reftwo{corollary:inducedgwtol} and Corollary~\reftwo{corollary:filtration} that this map is an isomorphism on homotopy groups. 
Under the analogy between $\GW(\C,\QF\qshift{-d})$ and $\MTSO(d)$ (see Remarks~\reftwo{remark:geometric} and~\reftwo{remark:geometric-2}) this equivalence
corresponds to the left of the equivalences
\[
\colim_d \SS^d\otimes \CCob_d \simeq \Omega^\SO \quad \text{and} \quad \colim_d\SS^d\otimes\MTSO(d) \simeq \MSO,
\]
which is equated with the right hand one by the theorems of Galatius, Madsen, Tillman and Weiss on the one and Thom on the other hand; . %
even the definition of $\L(\C,\QF)$ in terms of the $\rho$-construction is modelled on Quinn's construction of the ad-spectrum of manifolds $\Omega^\SO$, whose homotopy groups by construction are the cobordism groups. 
\end{remark}

Now since the functor $(\C,\QF) \mapsto \K(\C,\QF)^\tC$ is bordism invariant by Example~\reftwo{example:tate}, the composite $\GW \rightarrow \K^\hC \rightarrow \K^\tC$ factors uniquely over a map $\Xi \colon \L \rightarrow \K^\tC$ and we obtain the main result of this paper:
\begin{corollary}[The fundamental fibre square]
\label{corollary:tate-square-L}%
The natural square 
\[
\begin{tikzcd}
\GW(\C,\QF) \ar[r,"\bord"] \ar[d,"\fgt"] & \L(\C,\QF) \ar[d,"\Xi"] \\
\K(\C,\QF)^\hC \ar[r] & \K(\C,\QF)^\tC 
\end{tikzcd}
\]
is cartesian for every Poincaré \(\infty\)-category $(\C,\QF)$ and in particular, there is a natural fibre sequence
\[
\K(\C,\QF)_\hC \xrightarrow{\hyp} \GW(\C,\QF) \xrightarrow{\bord} \L(\C,\QF).
\]
\end{corollary}

\begin{proof}
Apply Corollary~\reftwo{corollary:tate-square} in combination with Corollary~\reftwo{corollary:GWHyp} and Theorem~\reftwo{theorem:lisbordgw}.
\end{proof}

\begin{corollary}
\label{corollary:GW-split-one-half}%
There is a canonical equivalence
\[
\GW(\C,\QF)[\tfrac 1 2] \simeq \K(\C,\QF)[\tfrac 1 2]_\hC \oplus \L(\C,\QF)[\tfrac 1 2]
\]
natural in the Poincaré $\infty$-category $(\C,\QF)$ and in particular
\[
\GW_i(\C,\QF)[\tfrac 1 2] \cong \K_i(\C)[\tfrac 1 2]_\Ct \oplus \L_i(\C,\QF)[\tfrac 1 2].
\]
\end{corollary}

\begin{proof}
Note that the homotopy orbit spectral sequence collapses for spectra in which $2$ is invertible. Thus the first claim implies the second. For the first note that the composite
\[
\K(\C,\QF)[\tfrac{1}{2}]_\hC \simeq \K(\C,\QF)_\hC[\tfrac{1}{2}] \xrightarrow{\hyp} \GW(\C,\QF)[\tfrac{1}{2}] \xrightarrow{\fgt} \K(\C,\QF)^\hC[\tfrac{1}{2}] \to \K(\C,\QF)[\tfrac{1}{2}]^\hC
\]
is the norm of $\K(\C,\QF)[\tfrac{1}{2}]$ which is an equivalence, and thus splits $\K(\C,\QF)[\tfrac{1}{2}]_\hC$ off $\GW(\C,\QF)[\tfrac{1}{2}]$ as claimed.
\end{proof}

 \subsection{Localisation and the Mayer-Vietoris principle}

We next use the fundamental fibre sequence to establish localisation sequences. Fundamentally we have:

\begin{corollary}
\label{corollary:gwadd}%
The functor $\GW \colon \Catp \rightarrow \Spa$ is Verdier-localising.
\end{corollary}

\begin{proof}
Given Corollary~\reftwo{corollary:Ladd}, we need only recall that $\K$-theory is a Verdier-localising functor $\K \colon \Catx \rightarrow \Spa$ (as by Proposition~\reftwo{proposition:criterion-poincare} the underlying sequence of a Poincaré-Verdier sequence is indeed a Verdier sequence). The statement about $\K$-theory in large parts goes back to Waldhausen's fibration theorem and a direct proof in the present context was recently recorded in \cite{HLS}*{Section 6}. 
\end{proof}
This result is an extension of Schlichting's localisation theorem \cite{schlichting-derived}*{Theorem 6.6} in the context of differential graded categories on which $2$ acts invertibly.
As a consequence 
we now deduce localisation properties of Grothendieck-Witt spectra, which will form the basis of our analysis of the Grothendieck-Witt groups of Dedekind rings in the third paper of this series; see e.g.\ Corollary~\refthree{corollary:decomposition-s}. We recall that in order to obtain Poincar\'e-Verdier sequences from localisations, in general we need to introduce suitable subcategories of $\Dperf(-)$, controlled by a subgroup $c \subseteq \K_0(-)$ closed under the involution. 
\begin{remark}
\label{remark:cofinality-GW}%
For algebraic K-groups the map $\K_i(\D^\mathrm c(R)) \to \K_i(\Dperf(R)) = \K_i(R)$ induces an isomorphism $i>0$ and is the inclusion $\mathrm c \to \K_0(R)$ for $i = 0$ by the cofinality theorem. We extend Schlichting's hermitian analogue thereof in \paperfour by showing that for any pair of involution-closed subgroups $\mathrm c \subseteq \mathrm d \subseteq \K_0(R)$ the squares
\[
\begin{tikzcd}
\GW(\D^\mathrm c(R),\QF) \ar[r] \ar[d] & \GW(\D^\mathrm d(R),\QF) \ar[d] & \L(\D^\mathrm c(R),\QF) \ar[r] \ar[d] & \L(\D^\mathrm d(R),\QF) \ar[d] \\
\K(\D^\mathrm c(R),\Dual_\QF)^{\hC} \ar[r] & \K(\D^\mathrm d(R),\Dual_\QF)^{\hC} &\K(\D^\mathrm c(R),\Dual_\QF)^{\tC} \ar[r] & \K(\D^\mathrm d(R),\Dual_\QF)^{\tC}
\end{tikzcd}
\]
are cartesian; see Theorem~\reffour{theorem:strong-cofinality}. It follows that there are fibre sequences
\[
\GW(\D^\mathrm c(R),\QF) \lto \GW(\D^\mathrm d(R),\QF) \lto \mathrm{H}(\mathrm d/\mathrm c)^{\hC}
\]
\[
\L(\D^\mathrm c(R),\QF) \lto \L(\D^\mathrm d(R),\QF) \lto \mathrm{H}(\mathrm d/\mathrm c)^{\tC}.
\]
In particular, the map $\GW_i(\D^\mathrm c(R),\QF) \lto \GW_i(\D^\mathrm d(R),\QF)$ is an isomorphism for positive $i$ and injective for $i=0$. On the $\L$-theoretic side, we recover Ranicki's \emph{Rothenberg-sequences} 
\[
\dots \lto \L_i(\D^\mathrm c(R),\QF) \lto \L_i(\D^\mathrm d(R),\QF) \lto \widehat{\mathrm H}^{-i}(\Ct;d/c) \lto \L_{i-1}(\D^\mathrm c(R),\QF) \lto \dots
\]
\cite{RanickiATS1}*{Proposition 9.1}. 
\end{remark}

Next, recall the compatibility condition between a multiplicative subset and an invertible module, Definition~\reftwo{definition:module-compatible-with-induction} and Example~\reftwo{example:involution-induction-compatible}; in the setting of discrete rings with $S \subseteq A$ satisfying an Ore condition, it means that inverting the action of $S$ on an invertible module $M$ with involution over $A$ using one of the $A$-module structures on $M$ also inverts it for the other one.  
From Corollary~\reftwo{corollary:ore-localisation},
we immediately obtain:

\begin{corollary}
\label{corollary:localisation-sequences-papertwo}%
Let \(A\) be a discrete ring, $M$ a discrete invertible module with involution over $A$, $\mathrm c \subseteq \K_0(A)$ a subgroup closed under the involution induced by $M$ and \(S \subseteq A\) a multiplicative subset compatible with $M$, such that \((A,S)\) satisfies the left Ore condition. Let $\D^\mathrm c(A)_S$ denote the full subcategory of $\D^\mathrm c(A)$ spanned by the $S$-torsion complexes. Then the inclusion and localisation functors fit into fibre sequences
\[
\GW(\D^\mathrm{c}(A)_S,\Qgen m M) \longrightarrow \GW(\D^\mathrm{c}(A),\Qgen m M) \longrightarrow \GW(\D^{\im(\mathrm c)}(A[S^{-1}]),\Qgen m{M[S^{-1}]})
\]
\[
\L(\D^\mathrm{c}(A)_S,\Qgen m M) \longrightarrow \L(\D^\mathrm{c}(A),\Qgen m M) \longrightarrow \L(\D^{\im(\mathrm c)}(A[S^{-1}]),\Qgen m{M[S^{-1}]}),
\]
for all $m \in \ZZ \cup \{\pm \infty\}$.
\end{corollary}

Granting the hermitian cofinality theorem, the above fibre sequence of Grothendieck--Witt spectra corresponds to the affine version of Schlichting's localisation sequence \cite{schlichting-derived}*{Theorem 9.5}. The fibre sequence of $\L$-spectra recovers localisation sequences of Ranicki's from \cite{Ranickiyellowbook}*{Section 3.2}, upon investing the identification of the genuine $\L$-spectra in Theorem~\refthree{theorem:main-theorem-L-theory}.

\begin{remark}
By Corollary~\reftwo{corollary:quadratic-Verdier-sequence}, the quadratic version, i.e.\ the case $m= \infty$, of Corollary~\reftwo{corollary:localisation-sequences-papertwo} holds true for an arbitrary $\Eone$-ring spectrum $A$ and an invertible module $M$ with involution over $A$, but for the symmetric and genuine variants, one has to require further conditions, see e.g.\ Example~\reftwo{example:2-inverted-sphere}. 
\end{remark}

In a similar vein, one can compare localisations along a ring homomorphism. 
In order to do this, we have to introduce a bit of notation. Given a map $(f,g) \colon A \to A' \times B$ of rings, and an involution closed (with respect to suitable modules with involution) subgroup $\mathrm{c} \subseteq \K_0(A)$ as before, we denote by $\bar{\mathrm{c}} \subseteq \K_0(A)$ the subgroup $f^{-1}(f(\mathrm{c})) \cap g^{-1}(g(\mathrm{c}))$. 
Note that $\mathrm{c} \subseteq \bar{\mathrm{c}} \subseteq \K_0(A)$ where both inclusions can be strict.
\begin{proposition}
\label{proposition:analytic-isomorphism}%
Let \(p: A \rightarrow B\) be a homomorphism of discrete rings, $M$ and $N$ discrete invertible modules with involution over $A$ and $B$, respectively, $\eta \colon M \rightarrow N$ a group homomorphism that is $p \otimes p$-linear, $S \subseteq A$ a subset and $m \in \ZZ \cup \{\pm\infty\}$. Then if
\begin{enumerate}
\item the map $B \otimes_A M \rightarrow N$ induced by $\eta$ is an isomorphism,
\label{item:eta-is-an-iso}%
\item the subset $S$ is compatible with $M$,
\label{item:S-compatible-M}%
\item 
\label{item:right-mult-equiv}%
for every $s \in S$ the induced map $p \colon A\sslash s \rightarrow B \sslash p(s)$ on cofibres of right multiplication by $s$ and $p(s)$, respectively, is an equivalence in $\D(A)$,
\item 
\label{item:both-left-ore}%
the pairs $(S,A)$ and $(p(S),A)$ both satisfy the left Ore condition, and
\item 
\label{item:boundary-map}%
the boundary map $\widehat \mathrm H^{-m}(\Ct,N[p(S)^{-1}]) \rightarrow \widehat \mathrm H^{-m+1}(\Ct,M)$ in Tate cohomology of the short exact sequence
\[
M \xrightarrow{(-\eta,\mathrm{can})} N \oplus M[S^{-1}] \xrightarrow{(\mathrm{can},\eta)} N[p(S)^{-1}]
\]
vanishes,
\end{enumerate}
the square
\[
\begin{tikzcd}
 (\D^\mathrm{\bar{c}}(A),\Qgen{m}{M}) \ar[r] \ar[d] & (\D^{\im(\mathrm{c})}(A[S^{-1}]),\Qgen{m}{M[S^{-1}]}) \ar[d] \\
 (\D^{p(\mathrm{c})}(B),\Qgen{m}{N}) \ar[r] & (\D^{\im(p(\mathrm{c}))}(B[p(S)^{-1}]),\Qgen{m}{N[p(S)^{-1}]})
\end{tikzcd}
\]
is a Poincaré-Verdier square for every subgroup $\mathrm{c} \subseteq \K_0(A)$ stable under the involution induced by $M$, and so in particular becomes cartesian after taking $\GW$-, $\K$- or $\L$-spectra. 
\end{proposition}

Here, condition \reftwoitem{item:boundary-map} is to be interpreted as vacuous if $m = \pm \infty$. In \cite{motives}*{Proposition~4.4.1}, it is furthermore shown to hold for $B$ a flat, commutative (discrete) $A$-algebra, in case $M$ is an invertible $B$-module and $N$ is an invertible $A$-module. Note also that condition \reftwoitem{item:both-left-ore} is equivalent to requiring that $p$ induces an isomorphism on kernels and cokernels of right multiplication by any $s \in S$.

The $\K$-theoretic part is a classical result of Karoubi, Quillen and Vorst, see \cite{Vorst-localization}*{Proposition 1.5}, and investing the identification of the $\L$-spectra from the third instalment in this series, see Theorem~\refthree{theorem:main-theorem-L-theory}, the $\L$-theoretic part recovers an analogous result of Ranicki \cite{Ranickiyellowbook}*{Section~3.6}.

\begin{proof}
Let us start out by observing that the diagram 
\[
\begin{tikzcd}
	A \ar[r] \ar[d] & A[S^{-1}] \ar[d] \\
	B \ar[r] & B[p(S)^{-1}]
\end{tikzcd}
\]
is cartesian in $\D(A)$: Denoting the top horizontal fibre by $F$, this is equivalent to the assertion that $F \rightarrow B \otimes^\mathbb L_A F$ is an equivalence in $\D(A)$, but combining Example~\reftwo{example:ore} with assumptions \reftwoitem{item:right-mult-equiv} and \reftwoitem{item:both-left-ore}, this holds for any object of $\D(A)_S$. Tensoring the square with $M$ (over $A$) then produces the short exact sequence appearing in \reftwoitem{item:boundary-map}. Furthermore, from the Ore conditions, we also find that the natural map $B \otimes^\mathbb L_A A[S^{-1}] \rightarrow B[p(S)^{-1}]$ is an equivalence. It is then readily checked that $p(S)$ is compatible with $N$.

Now, the rows of the diagram of Poincaré $\infty$-categories
\[
\begin{tikzcd}
(\D^\mathrm{\bar{c}}(A)_S,\Qgen{m}{M}) \ar[r] \ar[d] & (\D^\mathrm{\bar{c}}(A),\Qgen{m}{M}) \ar[r] \ar[d] & (\D^{\im(\mathrm{c})}(A[S^{-1}]),\Qgen{m}{M[S^{-1}]}) \ar[d] \\
(\D^{p(\mathrm{c})}(B)_{p(S)},\Qgen{m}{N}) \ar[r] & (\D^{p(\mathrm{c})}(B),\Qgen{m}{N}) \ar[r] & (\D^{\im(p(\mathrm{c}))}(B[p(S)]^{-1}),\Qgen{m}{N[p(S)^{-1}]})
\end{tikzcd}
\]
are Poincaré-Verdier sequences by Proposition~\reftwo{proposition:left-Kan-extension-specific-structures}, and the vertical maps are Poincaré functors on account of assumption \reftwoitem{item:eta-is-an-iso}, see Lemma~\refone{lemma:Poincresscalars}.
We next show that the right square is cartesian by verifying the criteria from Proposition~\reftwo{proposition:poincare-verdier-square-char}. Adjointability of the square is a special case of  Example~\reftwo{example:adjointable-square-E1-rings}.
It remains to show that the left hand vertical map is an equivalence of Poincaré $\infty$-categories. The fact that the underlying functor of stable $\infty$-categories is an equivalence follows from assumption \reftwoitem{item:right-mult-equiv}: Since the $\infty$-categories $\Dperf(A)_S$ and $\Dperf(B)_{p(s)}$ are generated by the objects $A\sslash s$ and $B \sslash p(s)$ under shifts, retracts and finite colimits, see Example~\reftwo{example:ore}, essential surjectivity is clear and it suffices to compute%
\[
\Hom_A(A \sslash s,A \sslash t) \simeq \Hom_A(A \sslash s,B \sslash p(t)) \simeq  \Hom_\B(B \sslash p(s),B \sslash p(t)).
\]
for full faithfulness. See also \cite{LT}*{Proposition 1.17} for an alternative argument.
It remains to check that the natural map $\Qgen{m}{M}(X) \rightarrow \Qgen{m}{N}(p_!X)$ induced by $\eta$ is an equivalence for all $X \in \Dperf(A)_S$. For $m = \pm \infty$ this follows from the fact that $p_!$ is a Poincaré functor and an equivalence on underlying $\infty$-categories, as this evidently implies that $p_!$ induces an equivalence on bilinear parts. We are thus reduced to considering the linear parts for finite $m$. 
Using the adjunction $p_! \vdash p^*$, we have to show that for every $S$-torsion perfect complex of $A$-modules $X$, the map 
\[
\map_A(X,\tau_{\geq m}(M^{\tC})) \lto \map_A(X,p^\ast\tau_{\geq m}(N^{\tC}))
\]
induced by $\eta$ is an equivalence. Since the $\infty$-category $\Dperf_S(A)$ in generated under finite colimits and desuspensions by objects of the form $A\sslash s = \cof(A \stackrel{\cdot s}{\to} A)$ one can equivalently show that every element $s \in S$ acts invertibly on $F_m = \cof\left(\tau_{\geq m}(M^{\tC}) \to f^\ast\tau_{\geq m}(N^{\tC})\right)$, i.e.\ that the canonical map 
\[
F_m \lto F_m[S^{-1}]
\]
is an equivalence. We note that $F_m \to F_{-\infty}$ induces an isomorphism on homology groups in degrees larger than $m$, and that there is an exact sequence
\[
0 \lto \mathrm H_m(F_m) \lto \mathrm H_m(F_{-\infty}) \lto K \lto 0,
\]
where 
\[
K = \ker\left( \widehat{\mathrm H}^{-m+1}(\Ct;M )\to \widehat{\mathrm H}^{-m+1}(\Ct;N) \right).
\]
From the fact that the bilinear parts of the two functors agree, we find that $S$ acts invertibly on $F_{-\infty}$. Hence it remains to show that $S$ acts invertibly on $\mathrm H_m(F_m)$. 
The above short exact sequence maps into its localisation at $S$. Since this localisation is an exact functor, the snake lemma implies that it suffices to check that the map $K \to K[S^{-1}]$ is injective. 
Writing $M[S^{-1}]$ as $(R[S^{-1}]\otimes R[S^{-1}])\otimes_{R\otimes R} M$ and likewise for $N$, using assumption~\reftwoitem{item:S-compatible-M}, we find that 
\[
K[S^{-1}] = \ker\left( \widehat{\mathrm H}^{-m+1}(\Ct;M[S^{-1}] )\to \widehat{\mathrm H}^{-m+1}(\Ct;N[p(S^{-1})]) \right),
\]
since Tate cohomology commutes with filtered colimits in the coefficients (see the discussion in the proof of Proposition~\reftwo{proposition:left-Kan-extension-specific-structures}). The kernel of $K \to K[S^{-1}]$ therefore canonically identifies with the kernel of
\[
\widehat{\mathrm{H}}^{-m+1}(\Ct;M) \lto \widehat{\mathrm{H}}^{-m+1}(\Ct;N\oplus M[S^{-1}])
\]
which vanishes by assumption~\reftwoitem{item:boundary-map}.
\end{proof}

For example, with same notations as in Proposition~\reftwo{proposition:analytic-isomorphism}, we obtain:

\begin{corollary}
\label{corollary:open-cover-GW}%
Let $R$ be a discrete commutative ring, $M$ an invertible $R$-module with an $R$-linear involution, $f,g \in R$ elements spanning the unit ideal and $c \subseteq \K_0(R)$ closed under the involution associated to $M$. Then the square
\[
\begin{tikzcd}
\GW(\D^{\overline{\mathrm{c}}	}(R),\Qgen{m}{M}) \ar[r] \ar[d] & \GW\Huge(\D^{\im(\mathrm{c})}(R[\frac 1 f]),\Qgen{m}{M[1/f]}\Huge) \ar[d] \\
\GW\Huge(\D^{\im(\mathrm{c})}(R[\frac 1 g]),\Qgen{m}{M[1/g]}\Huge) \ar[r] & \GW\Huge(\D^{\im(\mathrm{c})}(R[\frac 1 {fg}]),\Qgen{m}{M[1/{fg}]}\Huge)
\end{tikzcd}
\]
and the analogous squares in $\K$ and $\L$-theory are cartesian.
\end{corollary}

\begin{proof}
We verify conditions \reftwoitem{item:eta-is-an-iso} through \reftwoitem{item:boundary-map} of the previous proposition. The first and fourth are obvious and the second is implied by the two $R$-module structures on $M$ agreeing. For the third one simply notes that $g$ acts invertibly on $R\sslash f$, since with $f$ and $g$ also any powers thereof span the unit ideal.
To verify the final condition, recall that Tate cohomology groups over $\Ct$ with coefficients in $M$ are $2$-periodic with values alternating between the kernels of the norm maps $\id_M \pm \sigma \colon M_\Ct \rightarrow M^\Ct$. Thus, we may check that $M \rightarrow M[1/f] \oplus M[1/g]$ induces injections on both these kernels. But taking coinvariants commutes with localisation at both $f$ and $g$, so the map in question is injective on the entire coinvariants.
\end{proof}

\begin{example}
If $A$ is a regular coherent ring, i.e.\ any finitely presented $A$-module is perfect, $\K_0(A)\rightarrow \K_0(A[1/h])$ is surjective for all $h \in A$ satifsying the Ore condition, so Corollary~\reftwo{corollary:open-cover-GW} holds true for the full $\infty$-categories of perfect complexes in all four corners; indeed, the finitely generated projective $A[1/h]$-modules are a generating system of $\K_0(A[1/h])$ and writing some such $P$ as a cokernel of a map between finitely generated free modules and clearing denominators, one sees that $P$ can be lifted to a finitely presented $A$-module $P'$. By assumption $P'$ is then perfect and thus a preimage of $P$. Examples of regular coherent rings include regular Noetherian rings (in the case of finite Krull dimension see \cite{stacks}*{Lemma 10.110.8}, for the general case one can run the same argument using the upper-semicontinuity of the projective dimension of a finitely generated $A$-module, thought of as a function on the quasi-compact space $\mathrm{Spec}(A)$, as was explained to us by Tamme), and Bézout rings, in particular valuation rings by \cite{stacks}*{Lemma 10.50.15} (as then every finitely generated submodule of a finite free module is again free by a simple induction).

\end{example}

\subsection{The real algebraic $\K$-theory spectrum and Karoubi-Ranicki periodicity}
\label{subsection:real-alg-K-Karoubi}%

Just as in \S\reftwo{subsection:Mackey}, the fundamental fibre square can be cleanly encapsulated as the isotropy separation square of a genuine $\Ct$-spectrum:

\begin{definition}
We define the \emph{real algebraic $\K$-spectrum} $\KR(\C,\QF)$ of a Poincaré \(\infty\)-category $(\C,\QF)$ to be the genuine $\Ct$-spectrum $\GW^\ghyp(\C,\QF)$. 
\end{definition}

In particular, from Corollary~\reftwo{corollary:geometric-fixed-points-are-bordism-invariant},
we obtain:

\begin{corollary}
The real algebraic $\K$-spectra define a Verdier-localising functor
\[
\KR \colon \Catp \longrightarrow \Spagc,
\]
such that 
\[
u\KR \simeq \K, \quad \KR^{\gCt} \simeq \GW \quad \text{and} \quad \KR^{\geofix} \simeq \L,
\]
where $u \colon \Spa^\gCt \rightarrow \Spa^\hC$ denotes the functor extracting the underlying $\Ct$-spectrum, and $(-)^\gCt$ and $(-)^\geofix \colon \Spa^\gCt \rightarrow \Spa$ denote the genuine and geometric fixed points, respectively. Furthermore, the isotropy separation square associated to $\KR(\C,\QF)$ is naturally equivalent to the fundamental fibre square of $(\C,\QF)$.
\end{corollary}

From Theorem~\reftwo{theorem:periodicity}, we succinctly find:

\begin{corollary}
\label{corollary:real-k-theory-shift}%
There are canonical equivalences
\[
\KR(\C,\QF\qshift{1}) \simeq \SS^{1-\sigma} \otimes \KR(\C,\QF)
\]
natural in the Poincaré \(\infty\)-category $(\C,\QF)$. In particular, any equivalence $(\C,\QF) \rightarrow (\C,\QF\qshift{n})$ induces a periodicity equivalence
\[
\KR(\C,\QF) \simeq \SS^{n-n\sigma} \otimes \KR(\C,\QF).
\]
\end{corollary}

Recall from \refone{proposition:shiftofqingeneral} and the discussion preceding Corollary~\reftwo{corollary:karoubi-fundamental} that, generally,
\[
(\C,\QF\qshift{2n}) \simeq (\C,(\SS^{n-n\sigma} \otimes \wtl{\QF})^\gC)
\]
via the $n$-fold shift in $\C$, and that there are many examples where $\QF$ and $(\SS^{n-n\sigma} \otimes \wtl{\QF})^\gC)$ are closely related, if not equal for some $n$. The simplest such situation is:

\begin{corollary}
\label{corollary:karoubi-skew-periodicity}%
Let $R$ be an $\Eone$-algebra over an $\Einf$-ring $k$ equipped with an $n\sigma$-orientation, $M$ an invertible $k$-module with involution over $R$ and $\mathrm c \subseteq \K_0(R)$ a subgroup closed under the involution induced by $M$. Then there are canonical equivalences
\[
\KR(\Modc{\mathrm c}{R},\QF^\sym_{M}) \simeq \SS^{2n-2n\sigma} \otimes \KR(\Modc{\mathrm c}{R},\QF^\sym_{M}) \quad \text{and} \quad \quad \KR(\Modc{\mathrm c}{R},\QF^\qdr_{M}) \simeq \SS^{2n-2n\sigma} \otimes \KR(\Modc{\mathrm c}{R},\QF^\qdr_{M}).
\]
If $n=2$, e.g. if $k$ is complex oriented or $2 \in \pi_0(k)^\times$, this refines to
\[
\KR(\Modc{\mathrm c}{R},\QF^\sym_{-M}) \simeq \SS^{2-2\sigma} \otimes \KR(\Modc{\mathrm c}{R},\QF^\sym_{M}) \quad \text{and} \quad \quad \KR(\Modc{\mathrm c}{R},\QF^\qdr_{-M}) \simeq \SS^{2-2\sigma} \otimes \KR(\Modc{\mathrm c}{R},\QF^\qdr_{M}),
\]
and if $R$ is furthermore connective we also have
\[
\KR(\Modc{\mathrm c}{R},\Qgen {m+1}{-M}) \simeq \SS^{2-2\sigma} \otimes \KR(\Modc{\mathrm c}{R},\Qgen m M).
\]
\end{corollary}

In particular, we obtain the following periodicity result:

\begin{corollary}[Karoubi-Ranicki periodicity]
\label{corollary:Karoubiper}%
Let $R$ be a discrete ring, $M$ a discrete invertible module with involution over $R$ and $\mathrm c \subseteq \K_0(R)$ a subgroup closed under the involution induced by $M$. Then the genuine $\Ct$-spectra
\[
\KR(\D^{\mathrm c}(R),\QF^\sym_M) \quad \text{and} \quad \KR(\D^{\mathrm c}(R),\QF^\qdr_M)
\]
are $(4-4\sigma)$-periodic, and even $(2-2\sigma)$-periodic if $R$ has characteristic $2$. Furthermore, we have
\[
\SS^{2\sigma-2} \otimes \KR(\D^{\mathrm c}(R),\QF^{\gq}_M) \simeq \KR(\D^{\mathrm c}(R),\QF^{\gev}_{-M}) \simeq \SS^{2-2\sigma} \otimes \KR(\D^{\mathrm c}(R),\QF^{\gs}_M).
\]
\end{corollary}

Passing to geometric fixed points in Corollary~\reftwo{corollary:karoubi-skew-periodicity} extends Ranicki's classical periodicity results for $\L$-groups from the case of discrete rings:

\begin{corollary}
\label{corollary:Ranickipergeneral}%
Let $R$ be an $\Eone$-algebra over an $\Einf$-ring $k$ equipped with an $n\sigma$-orientation, $M$ an invertible $k$-module with involution over $R$ and $\mathrm c \subseteq \K_0(R)$ a subgroup closed under the involution induced by $M$. Then there are canonical equivalences 
\[
\L(\Modc{\mathrm c}{R},\QF^\sym_{M}) \simeq \SS^{2n} \otimes \L(\Modc{\mathrm c}{R},\QF^\sym_{M}) \quad \text{and} \quad \quad \L(\Modc{\mathrm c}{R},\QF^\qdr_{M}) \simeq \SS^{2n} \otimes \L(\Modc{\mathrm c}{R},\QF^\qdr_{M}).
\]
If $n=2$, e.g. if $k$ is complex oriented or $2 \in \pi_0(k)^\times$, this refines to
\[
\L(\Modc{\mathrm c}{R},\QF^\sym_{-M}) \simeq \SS^{n} \otimes \L(\Modc{\mathrm c}{R},\QF^\sym_{M}) \quad \text{and} \quad \quad \L(\Modc{\mathrm c}{R},\QF^\qdr_{-M}) \simeq \SS^{n} \otimes \L(\Modc{\mathrm c}{R},\QF^\qdr_{M}).
\]
\end{corollary}

Of course, this corollary can also easily be obtained straight from the shifting behaviour of bordism invariant functors. Let us also recall briefly from \refone{example:signperiodicmodules} that there are also many examples with $n\sigma$-orientations for $n>2$:  For example, any Spin-orientation induces a $4\sigma$-orientation and a String-orientation an $8\sigma$-orientation (to see the latter, one can use that $\frac{p_1}{2}$ reduces to $w_4$ in $H^*(\mathrm{BSpin};\mathbb{F}_2)$ or that $\frac{p_1}{2} \in H^*(\mathrm{BSpin};\mathbb{Z})$ is primitive). We thus find that for example that
\[
\KR(\Mod^\omega({\mathrm{ko}}),\QF^\sym) \quad \text{and} \quad \KR(\Mod^\omega({\mathrm{tmf}}),\QF^\sym)
\]
are $(8-8\sigma)$- and $(16-16\sigma)$-periodic, respectively.\medskip

Let us finally discuss the case of form parameters as discussed at the end of \S \reftwo{subsection:Genauer-Bott-Karoubi}. 
\begin{corollary}
For a discrete ring $R$, a discrete invertible module with involution $M$ over $R$, a subgroup $c \subseteq \K_0(R)$ closed under the involution induced by $M$ and a form parameter $\fpm =(M_\Ct \rightarrow{\tau} \fpmg \xrightarrow{\rho} M^\Ct)$ on $M$ with $\rho$ injective and dual form parameter $\dfpm$, there is a canonical equivalence
\[
\KR(\D^c(R),\QF^{\g\dfpm}_{-M}) \simeq \SS^{2-2\sigma} \otimes \KR(\D^c(R),\QF^\gfpm_M).
\]
\end{corollary}

In particular, we find 
\[
\L(\D^c(R),\QF^{\g\dfpm}_{-M}) \simeq \SS^{2} \otimes \L(\D^c(R),\QF^\gfpm_M)
\]
by passing to geometric fixed points (or directly from bordism invariance of $\L$-theory). Note also, that this corollary can be applied twice whenever the dual form parameter $(-M)_\Ct \rightarrow M/\fpmg \rightarrow (-M)^\Ct$ again has its second map injective. This is the case if and only if $\rho \colon \fpmg \rightarrow M^\Ct$ is an isomorphism, in which case $\QF^\gfpm_M = \QF_M^\gs$ per construction. We thus find equivalences
\[
\SS^4 \otimes \L(\D^c(R),\QF^{\gq}_{M}) \simeq \SS^2 \otimes \L(\D^c(R),\QF^{\gev}_{-M}) \simeq \L(\D^c(R),\QF^{\gs}_{M}),
\]
as claimed in the introduction.

\begin{remarks}
\begin{enumerate}
\item The genuine $\L$-spectra really are not periodic in general, as we will show in \paperthree of this series by explicit computation of $\L(\Dperf(\ZZ),\QF^{\gs})$. 
\item It follows from \cite{WWIII}*{Theorem 4.5}, that $\L(\Modp{\SS},\QF^\sym)$ is not periodic as we will explain in Remark~\reftwo{remark:Lsymsphere} below. 
\item The higher periodicities for the $\KR$- and $\L$-spectra of ring spectra such as $\mathrm{ko}$ and $\mathrm{tmf}$, or their periodic versions, have not been studied prior to our work. We cannot rule out that they obey lower periodicities than what the respective symmetries provide. To illustrate that this can actually occur, let us mention that the algebraic $\pi$-$\pi$-Theorem of Weiss and Williams from \cite{WWII}, see also \refthree{corollary:pi-pi}, gives an equivalence 
\[
\L(\Mod^\omega(R),\QF^\qdr) \longrightarrow \L(\Dperf(\pi_0 R),\QF^\qdr)
\]
 for all connective $\Eone$-rings $R$, so that these spectra are indeed always $4$-periodic.
\end{enumerate}
\end{remarks}

\subsection{$\LA$-theory after Weiss and Williams}
\label{subsection:Weiss-Williams}%

In this final subsection, we will relate Grothendieck-Witt spectra to the $\LA$-spectra arising in the work of Weiss and Williams \cite{WWIII} via the fundamental fibre square. We start by comparing the map $\Xi\colon \L \rightarrow \K^{\tC}$ appearing in Corollary~\reftwo{corollary:tate-square-L} with the map $\L \rightarrow \K^{\tC}$ constructed by Weiss and Williams in \cite{WW-duality}*{Section 9}. Translated to our set-up, they consider the map %
\[
\Lspace(\C,\QF^\sym) = |\core\rho(\C,\QF^\sym)^\hC| \longrightarrow \Omega|\core\Q\rho(\C,(\QF^\sym)\qshift{1})^\hC|  \longrightarrow \Omega^\infty \ads(\K^\hC)(\C,\QF^\sym),
\]
where the second map is a colimit-limit interchange and the first is the realisation (in the $\rho$-direction) of the structure maps for the group completions of the additive functor $\core^\hC$; here $\QF^\sym$ denotes the symmetrisation of a hermitian structure $\QF$ on $\C$, given by
$\QF^\sym(X)=\Bil_\QF(X,X)^\hC$ as in Example~\refone{example:quadratic-symmetric}. For the reader filling in the details of the translation, we also note that Weiss and Williams use $\Omega^\sigma|\core\rS^e\rho(\C,\QF)|$ in place of $\Omega|\core\Q\rho(\C,(\QF^\sym)\qshift{1})^\hC|$, where $\rS^e$ is the edgewise subdivision of Segal's $\rS$-construction. We shall discuss the comparison between the $\Q$- and $\rS$-constructions in some detail in Appendix \reftwo{subsection:Spitzweck-compare}, and the equivalence between these terms is then an instance of Karoubi-periodicity, for example it is obtained by applying $\Omega^\infty$ to the equivalence of \reftwo{corollary:real-k-theory-shift} above.

Precomposing with the composite
\[
\Poinc(\C,\QF) \longrightarrow \Lspace(\C,\QF) \xrightarrow{\mathrm{pol}} \Lspace(\C,\QF^\sym)
\]
and unwinding definitions this is the same as 
\[
\Poinc(\C,\QF) \longrightarrow \GWspace(\C,\QF) \longrightarrow |\GWspace\rho(\C,\QF)| \xrightarrow{\fgt} |\Kspace\rho(\C,\QF^\sym)^\hC| \longrightarrow \Omega^\infty\ads(\K^\hC)(\C,\QF^\sym).
\]
The latter part of this composite can in turn be rewritten as
\[
\GWspace(\C,\QF)\simeq \Omega^\infty\GW(\C,\QF) \longrightarrow \Omega^\infty\ads\GW(\C,\QF) \xrightarrow{\fgt} \Omega^\infty\ads(\K^\hC)(\C,\QF^\sym).
\]
Now, the canonical map $\ads(\K^\hC)(\C,\QF) \rightarrow \ads(\K^\hC)(\C,\QF^\sym)$ is an equivalence, so the forgetful map is nothing but $\Omega^\infty \Xi \colon \Omega^\infty \Lspace(\C,\QF) \rightarrow \Omega^\infty \K(\C,\QF)^\tC$ under the identifications of Corollary~\reftwo{corollary:rho-bord} and Theorem~\reftwo{theorem:lisbordgw}. By the universal property of $\L$-theory in Theorem~\reftwo{theorem:lisbordgw}, we conclude that the Weiss-Williams map $\L \Rightarrow \K^\tC$ agrees with ours.

\begin{corollary}
\label{corollary:LA}%
For a space $B$ and a stable spherical fibration $\xi$ over $B$ the spectrum $\GW(\Spa_B^\cp,\QF^r_\xi)$, identifies with Weiss' and Williams' $\LA^r(B,\xi)$, where $r \in \{\sym,\vis,\qdr\}$, i.e.\ either of symmetric, visible or quadratic. In particular, we find equivalences
\[
\Omega^{\infty-1} \LA^r(B,\xi) \simeq |\Cob(\Spa_B^\cp,\QF^r_\xi)|.
\]
\end{corollary}

The displayed equivalence 
provides a cycle model for the left hand object, which seems to be new. In particular, specialising to $B=\ast$ we find that the $(-1)$-st infinite loop spaces of
\[
\GW(\Spa^\cp,\QF^\sym) \simeq \LA^\sym(*) \quad \text{and} \quad \GW(\Spa^\cp,\QF^\uni) \simeq \LA^\vis(*),
\]
where $\QF^\uni \colon (\Spa^\cp)\op \rightarrow \Spa$ is the universal hermitian structure of \S\refone{subsection:universal},
are the homotopy types of the cobordism $\infty$-categories of Spanier-Whitehead selfdual spectra, and selfdual spectra equipped with a lift along
\[
\Dual_\SS X \rightarrow (\Dual_\SS X)_2^\wedge \simeq \hom_\Spa(X,\Dual_\SS X)^\tC
\]
of the image of the selfduality map, respectively.

\begin{remark}
Here we applied a naming scheme similar to Lurie's suggestion of writing $\L^\qdr(R)$ and $\L^\sym(R)$ instead of Ranicki's $\L_\bullet(R)$ and $\L^\bullet(R)$ for what we would systematically call $\L(\Dperf(R),\QF^\qdr_R)$ and $\L(\Dperf(R),\QF^\sym_R)$.

In \cite{WWIII} the spectra $\LA^\qdr(B,\xi), \LA^\vis(B,\xi)$ and $\LA^\sym(B,\xi)$ are called $\LA_\bullet(B,\xi \otimes \SS^{d},d)$, $\mathrm{VLA}^\bullet(B,\xi \otimes \SS^{d},d)$ and $\LA^\bullet(B,\xi \otimes \SS^{d},d)$, where $d$ is the dimension of $\xi$.
\end{remark}

\begin{proof}[Proof of Corollary~\reftwo{corollary:LA}]
The spectra $\LA^r(B,\xi)$ of Corollary~\reftwo{corollary:LA} are defined by certain pullbacks \cite{WWIII}*{Definition 9.5}
\[
\begin{tikzcd}
\LA^r(B,\xi) \ar[r] \ar[d] & \L^r(B,\xi) \ar[d,"\Xi"] \\
\mathrm{A}(B,\xi)^\hC \ar[r] & \mathrm{A}(B,\xi)^\tC,
\end{tikzcd}
\]
which we claim correspond precisely to our fundamental fibre square from Corollary~\reftwo{corollary:tate-square-L} for $(\Spa_B^\cp,\QF^r_\xi)$. 

As we detailed in Section \refone{subsection:parametrised-spectra}, the sets of quadratic, symmetric and visible Poincaré objects that Weiss and Williams consider canonically map to $\Poinc(\Spa_B^\omega,\QF_\xi^r)$ for the appropriate value of $r$. Since they define their $\L$-spaces by a point-set implementation of the $\rho$-construction, and then deloop them by shifting the duality, see \cite{WW-duality}*{Sections 10 \& 11}, there result comparison maps between the $\L$-spectra, that are equivalences by \refone{corollary:compareWWquadsymL} and \refone{corollary:compare-visible-WW}. As we identified the map $\Xi$ occurring in the definition of the $\LA$-spectra with ours above, we obtain the claim from the well-known equivalence $\mathrm{A}(B) \simeq \K(\Spa_B^\cp)$. 
\end{proof}

\begin{remark}
In subsequent work, we will construct for $\xi$ a stable $(-d)$-dimensional vector bundle over $B$ a functor 
\[
\Cob_d^{\xi}\to \Cob(\Spa_B^\cp,\QF_{\xi}^\vis)
\]
from the geometric, normally-$\xi$ oriented cobordism category into the algebraic cobordism $\infty$-category of parametrised spectra over $B$. Through the equivalence
\[
\Omega^{\infty} \LA^\vis(B,\xi) \simeq \Omega|\Cob(\Spa_B^\cp,\QF^\vis_\xi)|
\]
this provides a factorisation of the Weiss-Williams index map
\[
\mathrm{BTop}^\xi(M) \longrightarrow \Omega^\infty \LA^\vis(B,\xi)
\]
whenever $M$ is a closed $\xi$-oriented manifold with stable normal bundle $\nu_M$, through the (topological) cobordism category $\Cob_d^{\xi}$; here $\mathrm{Top}^\xi(M)$ denotes the $\Eone$-group of $\xi$-oriented homeomorphisms of $M$. Now the homotopy type of the cobordism category is excisive in the bundle data by the main result of \cite{GomezLopez-Kupers}, which was exploited by Raptis and the ninth author in the $\K$-theoretic context for a new proof of the Dwyer-Weiss-Williams index theorem \cite{RaptisSteimle-SmoothDWW}. One can now follow their strategy so as to provide a canonical lift of the map $\Omega|\Cob_d^{\xi}| \rightarrow \Omega^\infty\LA^\vis(B,\xi)$, into the source of the assembly map of $\LA^\vis$. Inserting $\xi = \nu_M$, one obtains a lift of the Weiss-Williams index, resulting in a new perspective on substantial parts of \cite{WWIII} and by compatibility with the classical comparison between block homeomorphism and $\L$-spaces also on Weiss--Williams' map 
\[
\widetilde{\mathrm{Top}}(M)/\mathrm{Top}(M) \longrightarrow \Omega^{\infty+1}( \mathrm{Wh}(M)_\hC),
\]
into the topological Whitehead spectrum of $M$, by investing the fundamental fibre sequence.
\end{remark}

We offer one application of the identification $\GW(\Spa^\cp, \QF^\uni) \simeq \LA^\vis(*)$. To this end, recall from Proposition~\refone{proposition:GW-L-lax} that the functors $\GW_0$ and $\L_0$ and $\K_0 \colon \Catp \rightarrow \Ab$ are compatibly lax symmetric monoidal for the tensor product of $\Catp$ and that $(\Spa^\cp,\QF^\uni)$ is the unit of the tensor product on $\Catp$. Hence, there are rings maps
\[
\K_0(\Spa^\cp) \xleftarrow{\fgt} \GW_0(\Spa^\cp,\QF^\uni) \xrightarrow{\bord} \L_0(\Spa^\cp,\QF^\uni).
\]
Abbreviating the underlying spectra to $\K(\SS)$, $\GW^\uni(\SS)$ and $\L^\uni(\SS)$, and similarly their homotopy groups, we have:

\begin{proposition}
\label{proposition:pi0gwuniv}%
There is a commutative diagram with vertical maps isomorphisms 
\[
\begin{tikzcd}
[column sep=9ex]
\ZZ \ar[d] & \ZZ[e,h]/I \ar[l,"{8 \mapsfrom e,\ 2 \mapsfrom h}"'] \ar[d] \ar[r,"{e \mapsto e,\ h \mapsto 0}"] & \ZZ[e]/(e^2-8e) \ar[d] \\
\K_0(\SS) & \GW^\uni_0(\SS) \ar[l] \ar[r] & \L^\uni_0(\SS)
\end{tikzcd}
\]
where $e$ and $h$ denote the classes of the spherical $E_8$-lattice and $\hyp(\SS)$, respectively, and $I$ is the ideal generated by $e^2-8e, h^2-2h$ and $eh-8h$.
\end{proposition}

The calculation of $\K_0(\SS) = \pi_0\mathrm A(*)$ is of course due to Waldhausen and the calculation of $\L^\uni_0(\SS)$ is due to Weiss-Williams (due to the identification $\L^\uni(\SS) = \L(\Spa^\cp,\QF^\vis_\SS)$). The spherical lift of the $\mathrm E_8$-lattice comes from the canonical map $\pi_0\QF^\qdr(\SS^{\oplus i}) \rightarrow \pi_0\QF^\qdr(\ZZ^{\oplus i})$ being an isomorphism, which in fact implies that every quadratic form over $\ZZ$ lifts to $\SS$ uniquely up to homotopy.

Without multiplicative structures the result says that $\GW^\uni_0(\SS)$ is free of rank $3$ generated by the Poincaré spectra $\hyp(\SS)$, $(\SS,\id_\SS)$ and the spherical lift of the $E_8$-lattice. In particular, as already observed by Weiss and Williams, the equality $[E_8] = 8[\ZZ,\id_\ZZ]$ in the symmetric (Grothendieck-)Witt-group of the integers (a consequence of the classification of indefinite forms over $\ZZ$ through rank, parity and signature) does not lift to the sphere spectrum.
\begin{proof}
We first identify the underlying abelian groups.
From Corollary~\reftwo{corollary:tate-square-L}, 
we deduce a short exact sequence
\[
0 \longrightarrow \K_0(\SS) \xrightarrow{\hyp} \GW^\uni_0(\SS) \longrightarrow \L^\uni_0(\SS) \longrightarrow 0;
\]
here the injectivity of the map $\hyp$ follows e.g.\ from the observation that the composite with the forgetful map $\GW^\uni_0(\SS) \to \K_0(\SS)$ is multiplication by $2$ map on $\ZZ$.
Next, we argue that $\L_0^\uni(\SS)$ is free of rank two. To do so, we recall that
Weiss and Williams constructed a fibre sequence
\[
\L^\qdr(\SS) \longrightarrow \L^\uni(\SS) \longrightarrow \SS \oplus \MTO(1),
\]
by identifying the latter term with visible, normal (or hyperquadratic) $\L$-theory of the sphere in \cite{WWIII}*{Theorem 4.3}; a proof of this fibre sequence in the language of the present paper will appear as part of \cite{Harpaz-Nikolaus-Shah}, which derives a general formula for relative $\L$-spectra. By the algebraic $\pi$-$\pi$-theorem the base change map
$\L^\qdr(\SS) \to \L^\qdr(\ZZ)$
is an equivalence; this appears for example as \cite{WWII}*{Proposition 6.2}, a proof in the present language is given in \cite{Lurie-L-theory}*{Lecture 14} and we will also derive it in the third installment of this series, see Corollary~\refthree{corollary:pi-pi}. Using that $\L^\qdr_0(\ZZ) \to \L^\sym_0(\ZZ)$ identifies with the inclusion $8\ZZ \subseteq \ZZ$ and that $\L^\qdr_{-1}(\ZZ) = 0$, the Weiss--Williams sequence induces a short exact sequence on $\pi_0$. To deduce that
$\L_0^\uni(\SS) \cong \ZZ^2$, generated by the spherical $E_8$-lattice and $(\SS,\id_\SS)$ (compare the discussion following \cite{WWIII}*{Theorem 4.3}) we will show $\pi_0(\MTO(1)) = 0$. Then, since by Example~\refone{example:quadratic-module} the map
$\L_0^\qdr(\SS) \to \L_0^\uni(\SS)$
is an $\L_0^\uni(\SS)$-module map, we find that the map $\L_0^\uni(\SS) \to \pi_0(\SS \oplus \MTO(1)) \cong \mathbb Z$ admits the unit as a section as needed. Finally, $\pi_0(\MTO(1))=0$ follows e.g.\ from the Genauer fibre sequence
\[
\MTO(1) \longrightarrow \SS[\BO(1)] \longrightarrow \SS,
\]
whose right hand map is the transfer map for $\mathrm{O}(1)$ \cite{GMTW}*{Remark 3.2}, which is evidently injective on $\pi_0$ and surjective on $\pi_1$ by the Kahn--Priddy theorem \cite{KahnPriddy} or direct calucation. %

We are left to calculate the ring structures on $\GW^\uni_0(\SS)$ and $\L^\uni_0(\SS)$. We start with the latter. Using again that the map
$
\L_0^\qdr(\SS) \to \L_0^\uni(\SS)
$
is an $\L_0^\uni(\SS)$-module map, we obtain $e^2 = n e$ for some $n \in \ZZ$. Mapping to the integers shows that $n=8$, giving the claim. For the ring structure of $\GW^\uni_0(\SS)$, we similarly observe that the exact sequence at the start of the proof
consists of $\GW^\uni_0(\SS)$-modules by Corollary~\refone{corollary:hyp-fgt-GW}. This immediately gives $eh = 8h$ and $h^2 = 2h$, and also that $e^2 = 8e + kh$ for some $k \in \ZZ$. But then we find
\[
64h = 8he = he^2 = h(8e+kh) = 64h + 2kh
\]
which forces $k=0$. 
\end{proof}

\begin{remark}
\label{remark:lower-universal-GW-groups}%
As a consequence of the fundamental fibre sequence, we find that the map $\GW^\uni(\SS) \to \L^\uni(\SS)$ induces an isomorphism on negative homotopy groups. In Example~\refthree{example:pi-pi} we show that the canonical map $\L^\uni(\SS) \to \L^\gs(\ZZ)$ also induces an isomorphism on negative homotopy groups. Combined with the calculation of the latter, this results in
\[
\GW_{n}^\uni(\SS) \cong \begin{cases} 0 & \text{ for } n= -1,-2 \\ \L^\qdr_n(\ZZ) & \text{ for } n< -2. \end{cases}
\]
The case $n<-2$ is also a direct consequence of the Weiss--Williams fibre sequence 
\[
\L^\qdr(\SS) \longrightarrow \L^\uni(\SS) \longrightarrow \SS \oplus \MTO(1)
\]
and the $\pi$-$\pi$-theorem. Regarding the remaining degrees one finds an exact sequence
\[
0 \longrightarrow \mathrm{L}_{-1}^\uni(\mathbb S) \longrightarrow \pi_{-1}\MTO(1) \longrightarrow \L_{-2}^\qdr(\mathbb S) \longrightarrow \L_{-2}^\uni(\mathbb S) \longrightarrow 0
\]
with both middle groups isomorphic to $\mathbb Z/2$, and one can indeed check by hand that a generator of $\L_{-2}^\qdr(\mathbb S) = \mathbb Z$ admits a Lagrangian for $\QF^\uni$, just as the corresponding element of $\L_{-2}^\qdr(\mathbb Z)$ admits one for $\QF^\sym$.

\end{remark}

\begin{remark}
\label{remark:Lsymsphere}%
Similar to the sequence used in the previous remark, Weiss and Williams produce a fibre sequence
\[
\L^\qdr(\ZZ) \longrightarrow \L^\sym(\SS) \longrightarrow (\SS_2^\wedge \otimes \SS_2^\wedge) \oplus \MTO(1),
\]
in \cite{WWIII}*{Theorem 4.5}, which rules out any sort of periodicity for $\L^\sym(\SS)$.
\end{remark}

\begin{appendices}
\renewcommand{\appendixtocname}{{\bf Appendices}}
\addappheadtotoc

\section{Verdier sequences, Karoubi sequences and stable recollements}
\label{appendix:AppIIA}%
In this appendix we investigate in detail the \(\infty\)-categorical variants of the notion of Verdier sequences, i.e.\ sequences in \(\Catx\) that are simultaneously fibre and cofibre sequences, and the same notion up to idempotent completion, called Karoubi sequences. The results are mostly well-known and various parts can be found in the literature, but we do not know of a coherent account at the level of detail we need. In the hope that it can serve as a general reference for this material, we have attempted to kept this appendix mostly self-contained.

For the reader familiar with \cite{BGT}, here is a comparison of the terminology we will develop in this appendix: A Karoubi sequence is called an \emph{exact sequence} in~\cite{BGT}, while our notion of a Verdier sequence corresponds to that of a \emph{strict-exact} sequence in~\cite{BGT}; this follows from Corollary~\reftwo{corollary:criterion-projection}, Proposition~\reftwo{proposition:inclusion-criterion} and Proposition~\reftwo{proposition:criterion-karoubi}. 
Our notion of a split Verdier sequence is however stricter than the corresponding notion of \emph{split-exact} sequence in~\cite{BGT}, since we require the projection to have both adjoints (in which case these adjoints are automatically fully faithful, and the injection has both adjoints as well, see Proposition~\reftwo{proposition:criterion-split}), while in the corresponding notion in~\cite{BGT} only the right adjoints are assumed to exist. In particular, the notion of a split-exact sequence in~\cite{BGT} corresponds to what we call a right-split Verdier sequence.

We write \(\Cat\) for the \(\infty\)-category of (small) \(\infty\)-categories, and \(\Catx\) for its (non-full) subcategory spanned by stable \(\infty\)-categories and exact functors. %
\newline

\subsubsection*{Brief recollection of inductive and projective completions}
\label{subsection:recall-ind-pro}%
For an \(\infty\)-category \(\C\) its inductive completion \(\Ind(\C)\) is given by the smallest full subcategory of \(\Fun(\C\op,\Sps)\) containing the representable functors and closed under filtered colimits. If \(\C\) admits finite colimits then we may also identify \(\Ind(\C)\) with the \(\infty\)-category of left exact functors \(\C\op \to \Sps\), i.e.\ those preserving finite limits, see~\cite{HTT}*{Corollary 5.3.5.4}. If \(\C\) is furthermore stable then this is equivalent, by the universal property of \(\Sp\), to the \(\infty\)-category of exact functors \(\C\op \to \Sp\) (and in particular $\Ind(\C)$ is stable as soon as $\C$ is so). Under this equivalence the covariant dependence of $\Ind(\C)$ in $\C$ is via left Kan extension (an operation which preserves exactness, see Lemma~\refone{lemma:kan-extension-exact-quadratic}), so that each induced functor $\Ind(\C) \to \Ind(\C')$ admits a right adjoint given by restriction.
We note that for $\C$ small $\Ind(\C)$ is generally large, but still locally small. The dual notion is that of projective completion, which, for a stable \(\C\), is given by \(\Pro(\C) = \Ind(\C\op)\op = \Funx(\C,\Sp)\op\).

For $\C$ small and $\E$ possibly large and possessing filtered colimits, restriction and left Kan extension along the Yoneda embedding give inverse equivalences between functors \(\C \to \E\) and filtered colimits preserving functors \(\Ind(\C) \to \E\), 
see~\cite{HTT}*{Proposition 5.3.5.10}. For $\C$ and $\E$ in addition stable, they restrict to equivalences between exact functors $\C \to \E$ and colimit preserving functors $\Ind(\C) \rightarrow \E$.
Finally, for \(g\colon \C \to \C'\) an exact functor the left Kan extension \(g_!\vphi\) of any (not necessarily exact) functor \(\vphi\colon \C \to \E\) 
can be computed as the composite
\[
\C' \lrar \Ind(\C') \xrightarrow{\,r\,} \Ind(\C) \xrightarrow{\wtl{\vphi}} \E, 
\]
where \(\wtl{\vphi}\) is the essentially unique filtered colimits preserving extension of \(\vphi\) and \(r\) the right adjoint of \(\Ind(g)\). In particular, this left Kan extension sends exact functors to exact functors for any target $\E$. Similarly, right Kan extensions can be computed via Pro-completions, and preserve exact functors.

For an \(\infty\)-category \(\C\) we denote by \(\C^{\natural}\) its idempotent completion. It can be constructed explicitly as the full subcategory of \(\Ind(\C)\) spanned by compact objects (\cite{HTT}*{Lemma 5.4.2.4}), or the full subcategory of \(\Pro(\C)\) spanned by the cocompact objects. In particular, if \(\C\) is stable then \(\C^{\natural}\) is again stable.

\subsection{Verdier sequences}
\label{section:appendix-verdier}%
We begin by analysing in detail the notion of a Verdier sequence. We recall the definition:
\begin{definition}
\label{definition:verdier}%
Let 
\[
\C \xrightarrow{f} \D \xrightarrow{p} \E
\]
be a sequence in \(\Catx\) with vanishing composite. We will say that it is a \defi{Verdier sequence} if it is both a fibre and a cofibre sequence in \(\Catx\). In this case we will refer to \(f\) as a \defi{Verdier inclusion} and to \(p\) as a \defi{Verdier projection}.
\end{definition}

\begin{remark}
Note that the extension of a sequence $\C \rightarrow \D \rightarrow \E$ with vanishing composte to a square with additional corner $0$ is essentially unique, since $0 \in \Fun(\C,\E)$ is initial (and terminal). The condition of being a fibre or cofibre sequence in $\Catx$ refers to such an extension being a cartesian or cocartesian square, respectively.
\end{remark}

\begin{example}
If \(R\) is a commutative ring, \(a \in R\) is an element, and $R[a^{-1}]$ is the associated localisation of $R$ obtained by inverting $a$, then
\[
\Dfree_a(R) \lrar \Dfree(R) \lrar \Dfree(R[a^{-1}])
\]
is a Verdier sequence, where $\Dfree(R)$ and $\Dfree(R[a^{-1}])$ denote the corresponding finitely presented derived categories (that is, the full subcategories of the respective derived categories spanned by the bounded complexes of finitely generated free modules), %
and \(\Dfree_a(R)\) 
is the stable subcategory of \(\Dfree(R)\) generated by $R\sslash a$. %
Being of particular interest in applications, we have dedicated \S\reftwo{section:appendix-module-examples} to a systematic study of Verdier sequences among module \(\infty\)-categories (the above statement, for example, is a special case of Corollary~\reftwo{corollary:ore-2}).
\end{example}

To make the notion of a Vedier sequence more concrete, let us recall how to compute fibres and cofibres in \(\Catx\): The fibre of an exact functor \(f\colon \C\to \D\) is computed in \(\Cat\) and given by the kernel \(\ker(f)\), which is the full subcategory of \(\C\) on all objects mapping to a zero object in \(\D\). Cofibres, in turn, are described by Verdier quotients: 

\begin{definition}
\label{definition:quotient}%
Let \(f\colon \C \to \D\) be an exact functor between stable \(\infty\)-categories. We say that a map in \(\D\) is an \defi{equivalence modulo \(\C\)} if its fibre (equivalently, its cofibre) lies in the smallest stable subcategory spanned by the essential image of \(f\). We write \(\D/\C\) for the localisation of \(\D\) with respect to the collection \(W\) of equivalences modulo \(\C\), and refer to \(\D/\C\) as the \defi{Verdier quotient} of \(\D\) by \(\C\).
\end{definition}

By definition, the Verdier quotient \(\D/\C\) only depends on the essential image of \(f\), and is characterised inside \(\Cat\) by the universal mapping property of localisations, that is, for every $\infty$-category $\E$, functors \(\D/\C \to \E\) identify via restriction with functors \(\D \to \E\) sending equivalences mod \(\C\) to equivalences.
While \(\D/\C\) is a priori an arbitrary \(\infty\)-category, it turns out to always be stable, and is in fact also characterised by a universal property as such:

\begin{proposition}
\label{proposition:verdier-localisation}%
Let \(f\colon \C \to \D\) be an exact functor between stable \(\infty\)-categories. The \(\infty\)-category \(\D/\C\) is stable, the localisation functor \(\D \to \D/\C\) is exact, and the sequence 
\[
\C \lrar \D \lrar \D/\C
\]
is a cofibre sequence in \(\Catx\). In addition, for every \(x,y \in \D\), the map
\[
\displaystyle\mathop{\colim}_{[\bet\colon z \to y] \in \C_{/y}}\map_{\D}(x,\cof(\bet)) \lrar \map_{\D/\C}([x],[y])
\]
is an equivalence.
\end{proposition}
\begin{proof}
This is \cite{NS}*{Theorem I.3.3}.
\end{proof}

\begin{remark}
\label{remark:kan-along-verdier-projection}%
In the situation of Proposition~\reftwo{proposition:verdier-localisation}, the mapping spectra formula also holds on the level of mapping spaces, as can be seen by applying to both sides the filtered colimits preserving functor $\Om^{\infty}$. Since colimits in spaces commutes with base change we can pass to fibres over a given $\alp\colon [x] \to [y]$ in $\D/\C$ to deduce that the comma category of the functor
\[
\C_{/y} \lrar \D_{/[y]} \quad\quad [\bet\colon z \to y] \mapsto (\cof(\beta),[\cof(\beta)] \simeq [y])
\]
over $(x,[\alp]) \in \D_{/[y]}$ is contractible, where on the right hand side we have the equivalence $[\cof(\beta)] \simeq [y]$ inverse to the one induced by the mod $\C$ equivalence $y \to \cof(\bet)$. %
We hence conclude that the above functor is cofinal, and so for every functor $\varphi\colon \D \to \A$ where $\A$ is some (not necessarily small) stable $\infty$-category with filtered colimits, the left Kan extension $p_!\varphi$ of $\varphi$ along $p$ is given by the formula
\[
(p_!\varphi)([y]) = \displaystyle\mathop{\colim}_{[\bet\colon z \to y] \in \C_{/y}}\varphi(\cof(\bet)) .
\]
Passing to opposites we also deduce the dual formula
\[
(p_*\varphi)([y]) = \displaystyle\mathop{\lim}_{[\bet\colon y \to z] \in \C_{y/}}\varphi(\fib(\bet))
\]
for the right Kan extension of $\varphi$ along $p$.
\end{remark}

\begin{corollary}
\label{corollary:criterion-projection}%
\label{corollary:verdier-surjective}%
Let \(p\colon \D \to \E\) be an exact functor between stable \(\infty\)-categories.
Then the following are equivalent:
\begin{enumerate}
\item
\label{item:verdier-projection}%
The functor \(p\) is a Verdier projection.
\item
\label{item:verdier-quotient}%
The induced map $\D/\mathrm{ker}(p) \rightarrow \E$ is an equivalence.%
\item
\label{item:localisation}%
The functor \(p\) is a localisation (at the maps it takes to equivalences).
\item
\label{item:surjective-strong}%
The functor \(p\) is essentially surjective and for every exact functor \(\varphi\colon \D \to \A\) where $\A$ is a cocomplete stable $\infty$-category the fibre of the map \(\varphi \Rightarrow p^*p_!\varphi\) is left Kan extended from its restriction to \(\ker(p)\).
\item
\label{item:surjective}%
The functor \(p\) is essentially surjective and for every exact functor \(\varphi\colon \D \to \Sp\) the fibre of the map \(\varphi \Rightarrow p^*p_!\varphi\) is left Kan extended from its restriction to \(\ker(p)\).
\end{enumerate}
\end{corollary}

In \reftwoitem{item:surjective} and \reftwoitem{item:surjective-strong}, 
note that $p_!\varphi$ is exact (see discussion in the recollection on Ind- and Pro-completions)
and that
the restriction of $p^*p_!\varphi$ to the kernel of $p$ vanishes, so that for $f \colon \ker(p) \rightarrow \D$ the inclusion we obtain a fibre sequence
\[
f_!f^*\varphi \Longrightarrow \varphi \Longrightarrow p^*p_*\varphi
\]
for every Verdier projection $p$. 
\begin{proof}
By Proposition~\reftwo{proposition:verdier-localisation} both \reftwoitem{item:verdier-projection} and \reftwoitem{item:verdier-quotient} are equivalent to \(p\) exhibiting its target as the cofiber of $\ker(p) \to \D$,
and so 
\reftwoitem{item:verdier-projection}\(\Leftrightarrow\)\reftwoitem{item:verdier-quotient}.
Next \reftwoitem{item:verdier-quotient}\(\Leftrightarrow\)\reftwoitem{item:localisation} holds 
because a morphism in \(\D\) maps to an equivalence in \(\E\) if and only if its cofibre lies in $\ker(p)$ (by stability of \(\E\) and exactness of \(p\)). 
Now if $p$ is a Verdier projection and $\varphi\colon \D \to \A$ is an exact functor then by the formula of Remark~\reftwo{remark:kan-along-verdier-projection} we have that 
\[
\fib[\varphi \Rightarrow p^*p_!\varphi](y) = \colim_{[z \to y] \in \ker(p)_{/y}} \varphi(z) = f_!f^*\varphi,
\]
so \reftwoitem{item:verdier-projection}$\Rightarrow$\reftwoitem{item:surjective-strong},
while clearly \reftwoitem{item:surjective-strong}$\Rightarrow$\reftwoitem{item:surjective}.
Finally, let us now show that~\reftwoitem{item:surjective}$\Rightarrow$\reftwoitem{item:verdier-quotient}.
For this, note that taking $\varphi = \map_{\D}(x,-)$ in~\reftwoitem{item:surjective} gives that the functor \(\fib[\map_{\D}(x,-) \Rightarrow \map_{\E}(p(x),p(-))]\) is left Kan extended from \(\ker(p)\), that is, that the sequence
\[
\displaystyle\mathop{\colim}_{[\bet\colon z \to y] \in \ker(p)_{/y}}\map_{\D}(x,z) \lrar \map_{\D}(x,y) \lrar \map_{\E}(p(x),p(y))
\]
is a fibre sequence.
By the last part of Proposition~\reftwo{proposition:verdier-localisation} we hence see that~\reftwoitem{item:surjective} for all representable \(\varphi\) is equivalent to saying that \(\D \to \E\) is essentially surjective and that the induced map \(\D/\ker(p) \to \E\) is fully faithful. Since \(\D \to \D/\ker(p)\) is essentially surjective (since localisations are so) we conclude that~\reftwoitem{item:surjective}$\Rightarrow$\reftwoitem{item:verdier-quotient}.
\end{proof}

We now examine the notion of a Verdier inclusion. For this, we need the following result:

\begin{lemma}
\label{lemma:kernel-of-projection-to-verdier-quotient}%
The kernel of the canonical map \(p\colon \D\to \D/\C\) consists of all objects of \(\D\) which are retracts of objects in the essential image of an exact functor \(f\colon \C \to \D\) between stable $\infty$-categories. 
\end{lemma}

\begin{proof}
Clearly, any retract of an object in \(\C\) lies in the kernel of \(p\). For the converse inclusion, let \(\x\) be an object of \(\ker(p)\). We note that by Proposition~\reftwo{proposition:verdier-localisation}, every exact functor \(\D \to \Spa\) that vanishes on \(\C\) also vanishes on \(\ker(p)\).
In particular, we may consider the exact functor \(\vphi_x\colon \D \to \Spa\) given by the formula 
\[
\vphi_{\x}(y) = \map_{\D/\C}(x,y) = \displaystyle\mathop{\colim}_{[\bet\colon \z \to \y] \in \C_{/\y}}\map_{\D}(\x,\cof(\bet)) \ .
\]
Then \(\vphi_{\x}\) vanishes on \(\C\), and so by the above vanishes on \(\ker(p)\).
In particular, \(\vphi_{\x}\) vanishes on \(\x \in \ker(p)\) itself, which implies the existence of a map \(\bet\colon \z \to \x\) for some \(\z \in \C\) such that \(\id\colon \x \to \x\) is in the kernel of the composite map \(\pi_0\map_{\D}(\x,\cof(0 \to \x)) \to \pi_0\map_{\D}(\x,\cof(\bet))\). We may then conclude that \(\id\colon \x \to \x\) factors through \(\z\) and hence \(\x\) is retract of \(\z\), as desired.
\end{proof}

\begin{proposition}
\label{proposition:inclusion-criterion}%
Let \(f\colon \C \to \D\) be an exact functor between stable \(\infty\)-categories.
Then the following are equivalent:
\begin{enumerate}
\item
\label{item:incl-verdier-inclusion}%
\(f\) is a Verdier inclusion.
\item 
\label{item:blubverdierin}%
the induced map $\C \rightarrow \mathrm{ker}(\D \rightarrow \D/\C)$ is an equivalence.
\item
\label{item:incl-ff-retracts}%
\(f\) is fully faithful and its essential image is closed under retracts in \(\D\).
\end{enumerate}
\end{proposition}
\begin{proof}
If \(f\) is a Verdier inclusion, then it is a kernel so that \reftwoitem{item:incl-ff-retracts} holds. On the other hand, if \reftwoitem{item:incl-ff-retracts} holds, then \(f\) extends to a cofibre sequence \(\C\to \D\to \D/\C\), and by Lemma~\reftwo{lemma:kernel-of-projection-to-verdier-quotient} this is also a fibre sequence. The same lemma shows \reftwoitem{item:blubverdierin}$\Leftrightarrow$\reftwoitem{item:incl-ff-retracts}.
\end{proof}

Summarising our discussion, we obtain: 
\begin{corollary}
\label{corollary:verdier-criterion}%
For a sequence \(\C\xrightarrow f \D \xrightarrow p \E\) in \(\Catx\) with vanishing composite, the following are equivalent:
\begin{enumerate}
\item The sequence is a Verdier sequence.
\item \(f\) is fully faithful with essential image closed under retracts in \(\D\), and \(p\) exhibits \(\E\) as the Verdier quotient of \(\D\) by \(\C\). 
\item \(p\) is a localisation, and \(f\) exhibits \(\C\) as the kernel of \(p\). 
\end{enumerate}
\end{corollary}

The collection of Verdier sequences is naturally organised into an \(\infty\)-category, see Definition~\reftwo{definition:Ver-category} below. The notion of a morphism we use in that definition is constructed out of that of adjointable squares. 
For this, let us first say that a square of stable \(\infty\)-categories
\[
\begin{tikzcd}
\D \ar[d, "f"']\ar[r, "p"] & \E  \ar[d, "f'"] \\
\D'\ar[r, "p'"] & \E'
\end{tikzcd}
\]
is vertically/horizontally \defi{inducitvely adjointable} (or $\Ind$-adjointable for short) if it becomes vertically/horizontally right adjointable after inductive completion, and that it is vertically/horizontally \defi{projectively adjointable} (or $\Pro$-adjointable) if it becomes vertically/horizontally left adjointable after projective completion; here, horizontally left adjointable means that the Beck-Chevalley transformation imposed after replacing the horizontal arrows with their left adjoints is an equivalence, 
see~\cite{HTT}*{Definition 7.3.1.1} (and similarly for vertical and and/or right adjointability).

\begin{remarks}
\label{remark:adj-concrete}%
\label{remark:adj-beck-chevalley}%
\
\begin{enumerate}
\item
\label{item:concrete}%
Concretely, a square as above is horizontally Ind-adjointable if and only if for every \(e \in \E\), the induced map 
\[
\colim_{p(d) \to e \in \D_{/e}}f(d) \longrightarrow \colim_{p'(d') \to f'(e) \in \D'_{/f'(e)}} d'
\]
is an equivalence in \(\Ind(\D')\). Since equivalences in \(\Ind(\D')\) are detected by mapping out of objects in \(\D'\), this is the same as saying that for every \(e \in \E\) and \(d'\in \D'\), the induced map of spectra
\[
\colim_{p(d) \to e \in \D'_{/e}} \map_\D(d',f(d)) \longrightarrow \map_{\E'}(p'(d'),f'(e))
\]
is an equivalence. In particular, if \(f\) is an equivalence then the square is horizontally Ind-adjointable if and only if the map \(\map_\E(p(d),e) \to \map_{\E'}(f'p(d),f'(e))\) is an equivalence for every \(d \in \D,e \in \E\).
\item
\label{item:equivalent-def}%
Another way of phrasing horizontal $\Ind$-adjointability is that the Beck-Chevalley transformation \(p_!f^*\vphi \Rightarrow (f')^*p'_!\vphi\) is an equivalence for every (not necessarily exact) functor $\varphi \colon \D' \rightarrow \A$ with $\A$ (not necessarily stable) admitting filtered colimits to form the requisite Kan extensions; indeed, the definition above is precisely this statement in the ``universal case'' of $\varphi$ being the canonical embedding 
$\D' \rightarrow \Ind(\D')$,
while the case of a general $\varphi$ follows by applying its unique filtered colimit extension $\wtl{\varphi}\colon \Ind(\D') \to \A$ to the map of inducitve systems in
\reftwoitem{item:concrete} and using the pointwise formula for left Kan extensions.

The same relation holds between horizontal Pro-adjointability and the Beck-Chevalley transformation involving the relevant right Kan extensions. 
\end{enumerate}	
\end{remarks}

\begin{lemma}
\label{lemma:adjointable-equivalent}%
For a square of stable $\infty$-categories the following conditions are equivalent:
\begin{enumerate}
\item
\label{item:hor-ind-adj}%
The square is horizontally Ind-adjointable.
\item
\label{item:ver-ind-adj}%
The square is vertically Pro-adjointable.
\item
\label{item:op-hor-pro-adj}%
Applying \((-)\op\) to the square yields a horizontally Pro-adjointable square.
\item
\label{item:op-ver-ind-adj}%
Applying \((-)\op\) to the square yields a vertically Ind-adjointable square.
\end{enumerate}
\end{lemma}

\begin{definition}
\label{definition:adjointable}%
We say that a square of stable $\infty$-categories is \defi{adjointable} if the square and its square of opposites both satisfy the equivalent conditions of Lemma~\reftwo{lemma:adjointable-equivalent}.
\end{definition}

\begin{proof}[Proof of Lemma~\reftwo{lemma:adjointable-equivalent}]
Firstly, the equivalence \(\Ind(-)\op \simeq \Pro(-\op)\) implies that \reftwoitem{item:hor-ind-adj}\(\Leftrightarrow\)\reftwoitem{item:op-hor-pro-adj} and \reftwoitem{item:ver-ind-adj}\(\Leftrightarrow\)\reftwoitem{item:op-ver-ind-adj}.
Secondly, using Remark~\reftwo{remark:adj-beck-chevalley}\;\reftwoitem{item:equivalent-def} and the identification $\Pro(-) =\Funx(-,\Spa)\op$ with functoriality given by left Kan extension (see \S\reftwo{subsection:recall-ind-pro}) we have that horizontal $\Ind$-adjointability implies vertical $\Pro$-adjointability.
Similarly, the same argument after taking opposites and flipping the square shows that \reftwoitem{item:op-ver-ind-adj}\(\Rightarrow\)\reftwoitem{item:op-hor-pro-adj}.
\end{proof}

\begin{example}
\label{example:adjointable-square-E1-rings}%
For a square 
\[
\begin{tikzcd}
A \ar[r] \ar[d] & A' \ar[d] \\
B \ar[r] & B'
\end{tikzcd}
\]
of \(\Eone\)-rings, the associated square of $\infty$-categories of compact module spectra
is vertically Ind-adjointable if and only if the induced map \(\alp\colon A' \otimes_A B \rightarrow B'\) is an equivalence, and is horizontally Ind-adjointable if and only if the induced map \(\bet\colon B \otimes_A A' \to B'\) is an equivalence. Indeed, this is a simple unwinding of the definitions coupled with the fact that natural transformations of exact functors out of \(\Dperf(B)\) or \(\Dperf(A')\) are equivalences if and only if they are so on \(B \in \Dperf(B)\) and \(A'\in \Dperf(A')\), respectively. In particular, the square of compact module categories is adjointable if and only if both \(\alp\) and \(\bet\) are equivalences. 
\end{example}

\begin{lemma}
\label{lemma:adj-verdier}%
Given a commutative diagram 
\[
\begin{tikzcd}
\C \ar[r,"f"]\ar[d, "h"] & \D \ar[r, "p"]\ar[d, "k"] & \E \ar[d, "g"] \\
\C'\ar[r, "f'"] & \D'\ar[r,"p'"] & \E'
\end{tikzcd}
\]
whose rows are Verdier sequences, the left square is horizontally projectively or inductively adjointable if and only if the right one is. In particular, the left square is adjointable if and only if the right square is adjointable.
\end{lemma}

We will say that a diagram as in Lemma~\reftwo{lemma:adj-verdier} is \defi{adjointable} if it satisfies the equivalent conditions of Lemma~\reftwo{lemma:adj-verdier}.

\begin{definition}
\label{definition:Ver-category}%
We define the \(\infty\)-category of Verdier sequences to be the (non-full) subcategory  \(\Ver \subseteq \Fun(\Del^2,\Catx)\) spanned by the Verdier sequences and the natural transformations between them corresponding to adjointable diagrams as above.
\end{definition} 

\begin{remark}
By Lemma~\reftwo{lemma:adj-verdier}, the \(\infty\)-category \(\Ver\) is equivalent by the obvious forgetful functor to the subcategory of \(\Ar(\Catx)\) 
spanned by Verdier projections 
and adjointable squares between them, and also to the subcategory of \(\Ar(\Catx)\) 
spanned by the Verdier inclusions and adjointable squares between them.
\end{remark}

\begin{proof}[Proof of Lemma~\reftwo{lemma:adj-verdier}]
Since taking opposites sends Verdier sequences to Verdier sequences it will suffice by Lemma~\reftwo{lemma:adjointable-equivalent} to prove that the left square is horizontally Ind-adjointable if and only if the right square is so. We use the characterisation of Remark~\reftwo{remark:adj-concrete}~\reftwoitem{item:equivalent-def}, and hence fix a cocomplete stable $\infty$-category $\A$. Since any functor \(\C'\to \A\) restricts from some functor \(\D'\to \Sp\) (e.g., from its own left Kan extension), the condition that the left square is horizontally Ind-adjointable can be equivalently stated as saying that for every exact \(\vphi\colon \D'\to \A\), the natural map 
\[
f_!f^*k^*\vphi = f_!h^*(f')^*\vphi \Longrightarrow k^*(f')_!(f')^*\vphi
\]
is an equivalence. On the other hand, since restriction along \(p\) is fully faithful we have that the right square is horizontally Ind-adjointable if and only if for every \(\vphi\colon \D'\to \A\), the natural map 
\[
p^*p_!k^*\vphi \Longrightarrow p^*g^*p'_!\vphi = k^*(p')^*p'_!\vphi
\]
is an equivalence. Now by \reftwoitem{item:surjective-strong} of Corollary~\reftwo{corollary:criterion-projection} 
these two transformations fit in a commutative diagram
\[
\begin{tikzcd}
f_!f^*k^*\vphi \ar[r,Rightarrow]\ar[d,Rightarrow] & k^*\vphi \ar[d,equal]\ar[r,Rightarrow] & p^*p_!k^*\vphi \ar[d,Rightarrow] \\
k^*f'_!(f')^*\vphi \ar[r,Rightarrow] & k^*\vphi \ar[r,Rightarrow] & k^*(p')^*p'_!\vphi
\end{tikzcd}
\]
whose rows are fibre sequences, and so one is an equivalence if and only if the other one is.
\end{proof}

Our next goal is to prove that Verdier projections are closed under base change, and that such base changes are automatically adjointable. The Poincaré analogue of such squares plays an important role in the present paper, and the following statement is used in the characterisation of such squares. 

\begin{proposition}
\label{proposition:verdier-square-char}%
For a commutative square %
\[
\begin{tikzcd}
\D \ar[r, "f"]\ar[d, "p"'] & {\D'} \ar[d, "p'"] \\
\E \ar[r, "f'"] & {\E'}
\end{tikzcd}
\]
of stable $\infty$-category and exact functors the following are equivalent:
\begin{enumerate}
\item
\label{item:cartesian}%
The square is cartesian and \(p'\) is a Verdier projection.
\item
\label{item:Achimscriterion}%
Both $p$ and $p'$ are Verdier projections, the induced functor $\ker(p) \rightarrow \ker(p')$ is essentially surjective and the induced map
\[
\hom_\D(x,y) \longrightarrow \hom_{\D'}(f(x),f(y))
\]
is an equivalence if at least one of $x,y$ lies in $\ker(p)$. 
\item
\label{item:adj-ker-equiv}%
The square is adjointable, both \(p\) and \(p'\) are Verdier projections and the induced functor \(\ker(p) \to \ker(p')\) is an equivalence.
\end{enumerate}

Similarly, the following are equivalent:
\begin{enumerate}
\item
\label{item:cocartesian1}%
The square is cocartesian and $f$ is a Verdier inclusion.
\item
\label{item:cocartesian2}%
The square is adjointable, both $f$ and $f'$ are Verdier inclusions and the induced functor $\mathrm{coker}(f) \to \mathrm{coker}(f')$ is an equivalence.
\end{enumerate}

\end{proposition}

The second criterion of the first three above is due to Krause, see~\cite{Krause-Picard}*{Lemma~3.9}. 

\begin{corollary}
\label{corollary:pullback-non-split-verd-proj}%
The collection of Verdier projections is closed under pullbacks and the collection of Verdier inclusions is closed under pushouts. 
\end{corollary}

\begin{proof}[Proof of Proposition~\reftwo{proposition:verdier-square-char}]
Throughout the proof, let us set $h := f' \circ p = f \circ p'$. Suppose first that~\reftwoitem{item:cartesian} holds. 
The condition that the square is cartesian implies, on the one hand, that 
the induced functor \(\ker(p) \to \ker(p')\) is an equivalence, and on the other hand that the map 
\[
\map_{\D}(x,y) = \map_{\E}(p(x),p(y)) \times_{\map_{\E'}(h(x),h(y))} \map_{\D'}(f(x),f(y)) \to \map_{\D'}(f(x),f(y))
\]
is an equivalence if at least one of \(x,y\) is in \(\ker(p)\). We now complete the implication \reftwoitem{item:cartesian}\(\Rightarrow\)\reftwoitem{item:Achimscriterion} by showing that \(p\) is indeed a Verdier projection. We first note that from Remark~\reftwo{remark:adj-concrete}~\reftwoitem{item:concrete} we get that the square
\[
\begin{tikzcd}
\ker(p) \ar[r,"g"]\ar[d,"i"] & \ker(p') \ar[d,"i'"] \\
\D \ar[r,"f"] & \D'
\end{tikzcd}
\]
is adjointable. Now we use the characterisation of Verdier projections given by \reftwoitem{item:surjective} of Corollary~\reftwo{corollary:criterion-projection}. In particular, by this characterisation \(p'\) is essentially surjective, and hence \(p\) is essentially surjective since the square is assumed cartesian. Next suppose given an exact functor \(\varphi\colon \D \to \Sp\). We need to show that \(\fib[\varphi \Rightarrow p^*p_!\varphi]\) is left Kan extended from \(\ker(p)\). 
Writing \(\varphi\) as a colimit of representable functors, the assumption that the square from the statement is cartesian implies that the square of exact functors
\[
\begin{tikzcd}
\varphi \ar[r, Rightarrow]\ar[d, Rightarrow] & f^*f_!\varphi \ar[d, Rightarrow] \\
p^*p_!\varphi \ar[r,Rightarrow] & h^*h_!\varphi
\end{tikzcd}
\]
is cartesian. It will hence suffice to show that \(\fib[f^*f_!\varphi \to h^*h_!\varphi]\) is left Kan extended from \(\ker(p)\).
Now since \(p'\) is a Verdier projection, the characterisation in~\reftwoitem{item:surjective} of Corollary \reftwo{corollary:criterion-projection} implies that the \(\fib[f_!\varphi \to (p')^*h_!\varphi]\) is left Kan extended from \(\ker(p')\). Since the upper square is adjointable 
we now conclude that \(f^*\fib[f_!\varphi \to (p')^*h_!\varphi] =\fib[f^*f_!\varphi \to h^*h_!\varphi]\) is left Kan extended from \(\ker(p)\), as desired.

To see that \reftwoitem{item:Achimscriterion}$\Rightarrow$\reftwoitem{item:adj-ker-equiv}, note that \reftwoitem{item:Achimscriterion} still implies that the induced functor $\ker(p) \rightarrow \ker(p')$
is an equivalence and hence that the upper square is adjointable by Remark~\reftwo{remark:adj-concrete}~\reftwoitem{item:concrete}, and by Lemma~\reftwo{lemma:adj-verdier}, this implies that the square from the statement is adjointable as well.

We now prove \reftwoitem{item:adj-ker-equiv}\(\Rightarrow\)\reftwoitem{item:cartesian}, and so we assume that \(p\) and \(p'\) are Verdier projections, that the square in the statement is adjointable, and that the induced functor \(g\colon \ker(p) \to \ker(p')\) is an equivalence. By Lemma~\reftwo{lemma:adj-verdier}, the upper square above is then adjointable as well.
We now run the argument of the previous part of the proof backwards. Given an exact functor \(\varphi\colon \D \to \Sp\), the characterisation \reftwoitem{item:surjective} of Corollary~\reftwo{corollary:criterion-projection} for both \(p\) and \(p'\) together with the adjointability of the upper square implies that the lower square is cartesian: Indeed, the induced map on vertical fibres can be identified with the composite
\[
i_!i^*\varphi \Rightarrow i_!g^*g_!i^*\varphi \Rightarrow f^*i'_!g_!i^*\varphi \Rightarrow f^*i'_!(i')^*f_!\varphi
\]
of the unit of $g_! \dashv g^*$ (which is an equivalence since $g$ is) and the two Beck-Chevalley transformations (which are equivalences by adjointability). 
Taking $\varphi$ to be representable, we thus conclude that the induced functor
\[
\D \lrar \E \times_{\E'} \D'
\]
is fully faithful. We now claim that the essential image of this functor is closed under retracts, that is, this functor is a Verdier inclusion. For this, 
suppose given an object \((e,d',\eta\colon f'(e)\simeq p'(d')) \in \E \times_{\E'} \D'\) 
together with a retract diagram
\[
(e,d',\eta) \longrightarrow (p(d),f(d),\id) \longrightarrow (e,d',\eta)
\]
with \(d \in \D\). Then \(e\) is a retract of \(p(d)\) and by the last part of Proposition~\reftwo{proposition:verdier-localisation}
we can find an \(x \in \D\) together with a map \(d \rightarrow x\) covering the projection \(p(d) \rightarrow e\). But then the fibre of the composite \((e,d',\eta) \rightarrow (p(d),f(d),\id) \rightarrow (p(x),f(x),\id)\) lies in the kernel of \(\E \times_{\E'} \D'\to \E\), and is hence in the image of \(\D\). We conclude that \((e,d',\eta)\) is in the image of \(\D\) as well, and hence that \(\D \to \E \times_{\E'} \D'\) is a Verdier inclusion.

Now, by the implication \reftwoitem{item:cartesian}\(\Rightarrow\)\reftwoitem{item:adj-ker-equiv} proven above, the sequence \(\ker(p') \to \E \times_{\E'} \D'\to \E\) is a Verdier sequence. 
Taking the levelwise Verdier quotient of the two sequences we hence obtain a cofibre sequence
\[
0 = \ker(p')/\ker(p) \lrar (\E \times_{\E'} \D') / \D \lrar \E/\E = 0 ,
\]
and so \((\E \times_{\E'} \D') / \D = 0\). By Lemma~\reftwo{lemma:kernel-of-projection-to-verdier-quotient} we conclude that the Verdier inclusion 
\(\D \to \E \times_{\E'} \D'\)
is an equivalence, as desired.

To see the claims about cocartesian squares, we reduce to cartesian ones by taking inductive completions and passing to right adjoints everywhere (this will require some result from later in this appendix, but is not used in their deduction, or in fact anywhere in this appendix). This entire process%
preserves and detects adjointable squares, sends Verdier inclusions to Verdier projections (by Example \reftwo{example:ind-pro-split}) and takes cocartesian squares to cartesian ones. Furthermore, it detects Verdier inclusions and cocartesian squares up to idempotent completion. %
We thus deduce from the first part, that the two statements are equivalent up to idempotent completions.
For the implication \reftwoitem{item:cocartesian1}$\Rightarrow$\reftwoitem{item:cocartesian2}, it then remains to show that $f'$ is a Verdier inclusion, i.e.\ its image is closed under retracts in $\E'$ or equivalently by Lemma \reftwo{lemma:kernel-of-projection-to-verdier-quotient} that it agrees with $\overline \E \coloneq \mathrm{ker}(\E' \rightarrow \E'/\E)$. Using Thomason's classification of dense subcategories \reftwo{theorem:thomason-classification} it suffices to show that $\K_0(\E) \rightarrow \K_0(\overline \E)$ is surjective, and this in turn follows via a diagram chase from comparing the long exact sequences of K-groups associated to the Verdier sequences $\D \rightarrow \D' \rightarrow \E'/\E$ and $\overline \E \rightarrow \E' \rightarrow \E'/\E$ together with the fact that %
$\K_0(\E) \oplus \K_0(\D') \to \K_0(\E')$ is surjective. For the implication \reftwoitem{item:cocartesian2}$\Rightarrow$\reftwoitem{item:cocartesian1} we need to argue that $\E'$ is the smallest subcategory of its idempotent completion generated from the essential images of $\E$ and $\D'$ in it, or equivalently, that $\K_0(\E) \oplus \K_0(\D') \to \K_0(\E')$ is surjective, which follows again from a small diagram chase.
\end{proof}

\begin{remark}
We have recorded several constructions preserving Poincar\'e-Verdier sequences in \S\reftwo{subsection:preservation-poincare-verdier}; upon dropping hermitian structures, one obtains a number of constructions that preserve Verdier-sequences (and their relatives discussed in the next two sections) using the tensor product of stable $\infty$-categories and the (co)tensor constructions of $\Catx$ over $\Cat$; we shall not repeat that material in this appendix.
\end{remark}

\subsection{Split Verdier sequences and stable recollements}
\label{section:appendix-split-verdier}%

In this section we discuss the notion of \emph{split Verdier sequence} (Definition~\reftwo{definition:split-verdier}), and its relationship with stable recollements (Definition~\reftwo{definition:stable-recollement} and Proposition~\reftwo{proposition:criterion-split}). The hermitian analogue of this notion (see \S\reftwo{subsection:poincare-recollements}) plays a central role in the present paper.

\begin{definition}
\label{definition:split-verdier}%
A Verdier sequence
\[
\C \xrightarrow{f} \D \xrightarrow{p} \E
\]
is \defi{right split} if \(p\) admits a right adjoint, \defi{left split} if $p$ admits a left adjoint und \defi{split} if $p$ admits both adjoints.
In this case we will refer to \(f\) as a (right or left) \defi{split Verdier inclusion} and to \(p\) as a (right or left) \defi{split Verdier projection}. We write \(\Ver^{\perp} \subseteq \Ver\) for the full subcategory spanned by the split Verdier sequences. 
\end{definition}

To obtain criteria similar to Corollary~\reftwo{corollary:criterion-projection} and \reftwo{proposition:inclusion-criterion} for split Verdier projections and inclusions 
we first investigate further the existence of adjoints for localisation functors.

\begin{lemma}
\label{lemma:bousfieldgenI}%
Let \(\D\) be a small \(\infty\)-category and \(W\) a collection of morphisms in \(\D\). Write \(\D_W \subseteq \D\) for the full subcategory spanned by the left \(W\)-local objects, that is, the objects \(d \in \D\) such that \(\Hom_{\D}(d,-)\) sends the morphisms in \(W\) to equivalences. Then the projection \(p\colon \D \to \D[W^{-1}]\) admits a left adjoint if and only if the composite functor \(\D_W \hrar \D \to \D[W^{-1}]\) is essentially surjective. In addition, in that case the arising left adjoint \(g\colon \D[W^{-1}] \to \D\) is fully faithful and its essential image is exactly \(\D_W\). 

The analogous statement holds if we ask instead for a right adjoint to \(p\) and we replace left \(W\)-local objects by right \(W\)-local objects, that is, the objects such that \(\Hom_{\D}(-,d)\) inverts \(W\).
\end{lemma}

\begin{corollary}
\label{corollary:bousfieldgenII}%
For a functor \(p\colon \D \to \E\) the following are equivalent:
\begin{enumerate}
\item
The functor \(p\) is a localisation functor and admits a left (resp. right) adjoint \(\E \to \D\).
\item
The functor \(p\) admits a fully faithful left (resp. right) adjoint \(\E \to \D\).
\end{enumerate}
\end{corollary}
\begin{proof}
The first statement implies the second by Lemma~\reftwo{lemma:bousfieldgenI}. The second implies the first by~\cite{HTT}*{Proposition 5.2.7.12}.  
\end{proof}

We refer to functors satisfying the equivalent conditions of Corollary~\reftwo{corollary:bousfieldgenII} as right (resp. left) Bousfield localisations. 

\begin{proof}[Proof of Lemma~\reftwo{lemma:bousfieldgenI}]
The dual final statement is obtained from the main statement by passing to opposites, and so we prove only the main statement. Since \(p \colon \D \rightarrow \D[W^{-1}]\) is essentially surjective, Yoneda's lemma implies that \(p\) admits a left adjoint if and only if, for each \(x \in \D\) the functor 
\[
\Hom_{\D[W^{-1}]}(p(x),p(-)) \colon \D \rightarrow \Sps
\]
is representable in \(\D\), in which case the desired left adjoint necessarily sends \(p(x)\) to the representing object. We also note that once such a left adjoint \(g\colon \D[W^{-1}] \to \D\) exists then it is automatically fully faithful by the full faithfulness of Yoneda 
and the fact that the restriction functor \(\Fun(\D[W^{-1}]\op,\Sps) \to \Fun(\D\op,\Sps)\) is fully faithful by the universal property of localisation.

To finish the proof it will hence suffice to show that an object \(y \in \D\) represents \(\Hom_{\D[W^{-1}]}(p(x),p(-))\) if and only if \(y\) is \(W\)-colocal and \(p(y) \simeq p(x)\).
For this, note that since \(\Hom_{\D[W^{-1}]}(p(x),p(-))\) inverts \(W\) such a representing object is necessarily \(W\)-colocal, and, on the other hand, any \(W\)-colocal object \(y\) represents \(\Hom_{\D[W^{-1}]}(p(y),p(-))\); indeed, since \(p\colon \D \to \D[W^{-1}]\) is a localisation we have that for any functor \(\vphi\colon \D \to \Sps\) that inverts \(W\) the unit map \(\vphi \Rightarrow p^*p_!\vphi\) is an equivalence, and the left Kan extension functor \(p_!\) sends \(\Hom_{\D}(y,-)\) to \(\Hom_{\D[W^{-1}]}(p(y),-)\).
\end{proof}

Now the statement of Lemma~\reftwo{lemma:bousfieldgenI} simplifies a bit 
in the stable case.
For this, given a stable \(\infty\)-category \(\D\) and a full subcategory \(\C \subseteq \D\), let us say that an object \(\y \in \D\) is \emph{right orthogonal} to \(\C\) if \(\map_{\D}(\x,\y) \simeq 0\) for every \(\x \in \C\) and that \(\y\) is \emph{left orthogonal} to \(\C\) if \(\map_{\D}(\y,\x) \simeq 0\) for every \(\x \in \C\). 
Let us write \(\C^r\) and \(\C^l\) for the full subcategories spanned by these objects, respectively. 
Applying Lemma~\reftwo{lemma:bousfieldgenI} in the case where \(\D\) is stable and \(W\) consist of equivalences mod some full subcategory \(\C \subseteq \D\) we immediately obtain:

\begin{corollary}
\label{corollary:bousfield}%
Let \(\D\) be a stable \(\infty\)-category and \(\C \subseteq \D\) a full subcategory. Then the Verdier projection \(\D \to \D/\C\) admits a left adjoint if and only if the composite functor \(\C^l \to \D \to \D/\C\) is essentially surjective. In addition, in this case the left adjoint \(g\colon \D/\C \to \D\) is fully faithful and 
its essential image is exactly \(\C^l\). Similarly, the analogous statement for \(\C^r\) and a right adjoint to \(p\) holds. 
\end{corollary}

\begin{corollary}
\label{corollary:criterion-split}%
An exact functor \(p\colon \D \to \E\) is a (left or right) split Verdier projection if and only if it admits fully faithful (left or right) adjoints. 
\end{corollary}

Put differently, a left (resp. right) split Verdier projection is precisely a right (resp. left) Bousfield localisation (hence automatically exact) functor among stable $\infty$-categories,
and similarly a split Verdier projection is precisely a functor among stable $\infty$-categories which is both a left and a right Bousfield localisation.

\begin{corollary}
\label{corollary:split-verdier-proj-pullback}%
The collection of (right of left) split Verdier projections is closed under pullbacks.
\end{corollary}
\begin{proof}
This is a direct consequence of Corollary~\reftwo{corollary:criterion-split}: Using the universal property of the pullback, one readily constructs the requisite adjoints from the original one (using the fact that these are fully faithful, and therefore sections of the original Verdier projection). That these are again fully faithful adjoints follows immediately from the description of mapping spaces in pullbacks of \(\infty\)-categories as pullbacks of mapping spaces. 
\end{proof}

\begin{remark}
\label{remark:Verdiersplitcocartequiv}%
There is another characterisation of left/right split Verdier projections
worth mentioning: 
an exact functor \(\C \rightarrow \D\) is a right split Verdier projection if and only if it is a cartesian fibration, and similarly for left split Verdier projections and cocartesian fibrations. 
One implication we proved in Lemma \reftwo{lemma:cartesian-in-bousfield} above, and the other is recorded in \cite{HHLN-two-var-fil}*{Example 5.6}.
\end{remark}

Passing now to split Verdier inclusions, we observe that the existence of adjoints on the side of the projection implies the same for the inclusion, and vice-versa:

\begin{lemma}
\label{lemma:quasi-split}%
Let 
\[
\C \xrightarrow{f} \D \xrightarrow{p} \E
\]
be a sequence in \(\Catx\) with vanishing composite. Then the following are equivalent:
\begin{enumerate}
\item
\label{item:split-verdier}%
The sequence is a left (right) split Verdier sequence.
\item
\label{item:fiber-left-to-right}%
The sequence is a fibre sequence, and \(p\) admits a fully faithful left (right) adjoint \(q\).
\item
\label{item:cofiber-left-to-right}%
The sequence is a cofibre sequence, and \(f\) is fully faithful and admits a left (right) adjoint \(g\).
\end{enumerate}
Furthermore, when these equivalent conditions hold the sequence  
\[
\C \xleftarrow g \D\xleftarrow q \E,
\]
formed by passing to left (right) adjoints is a right (left) split Verdier sequence.
\end{lemma}

\begin{remark}
\label{remark:quasi-split}%
The proof of Lemma \reftwo{lemma:quasi-split} also shows the following important addendum: In the case of a left split Verdier sequence the adjoints participate in a canonical exact sequence %
\[
qp \Rightarrow \id \Rightarrow fg
\]
with the first map the adjunction counit and the second the adjunction unit. Dually, in the case of right adjoints the unit and counit fit into a sequence in the opposite direction.
\end{remark}

\begin{corollary}
\label{corollary:equivalence-split}%
In a Verdier sequence \(\C \to \D \to \E\), the projection admits a left (resp. right) adjoint if and only if the inclusion admits a left (resp. right) adjoint.
\end{corollary}
\begin{proof}
By Proposition~\reftwo{proposition:inclusion-criterion} any Verdier inclusion is fully faithful and by 
Corollary~\reftwo{corollary:criterion-projection} and Corollary~\reftwo{corollary:bousfieldgenII}
any adjoint to a Verdier projection is fully faithful. The desired statement hence follows from Lemma~\reftwo{lemma:quasi-split}.
\end{proof}

\begin{corollary}
\label{corollary:criterion-split-inclusion}%
An exact functor \(f\colon \C \to \D\) is a (left or right) split Verdier inclusion if and only if it is fully faithful and admits (left or right) adjoints.
\end{corollary}

\begin{proof}[Proof of Lemma~\reftwo{lemma:quasi-split}]
We prove the claim left split Verdier sequences. The claim for right split ones follows by the dual argument (or by replacing all \(\infty\)-categories by their opposites). Now clearly~\reftwoitem{item:split-verdier}\(\Rightarrow\)\reftwoitem{item:fiber-left-to-right}, and the combination of~\reftwoitem{item:fiber-left-to-right} and \reftwoitem{item:cofiber-left-to-right} implies~\reftwoitem{item:split-verdier}. We now show that \reftwoitem{item:fiber-left-to-right} and \reftwoitem{item:cofiber-left-to-right} are equivalent to each other.

Suppose first that \reftwoitem{item:fiber-left-to-right} holds. Then we obtain a left adjoint \(g\) of \(f\) by considering the exact functor 
\[
\tilde g = \cof[qp \to \id]\colon \D \lrar \D
\]
given by the cofibre of the counit. Since \(q\) is fully faithful, the unit map \(\id\to qp\) is an equivalence, from which we can conclude that \(p\circ \tilde g\) vanishes. Thus, \(\tilde g\) factors uniquely through \(f\), giving rise to a functor \(g\colon \D\to \C\). We now claim that the canonical transformation \(\id \to \tilde g = f\circ g\) acts as a unit exhibiting \(g\) as left adjoint to \(f\). Given objects \(X\in \D\) and \(Y \in \C\), it will suffice to check that the composite map
\[
\map_{\C}(g(X),Y) \lrar \map_{\D}(fg(X),f(Y)) \lrar \map_{\D}(X,f(Y))
\]
is an equivalence of spectra. Indeed, the first map is an equivalence since \(f\) is fully faithful and the second map is an equivalence because its cofibre is \(\map_{\D}(qp(X),f(Y)) \simeq \map_{\E}(p(X),pf(Y)) \simeq 0\). 

In this situation, \(p\) is a localisation by Corollary~\reftwo{corollary:bousfieldgenII}, so the sequence formed by \(f\) and \(p\) is a Verdier sequence by Corollary~\reftwo{corollary:verdier-criterion}, in particular a cofibre sequence. Also, the kernel of \(g\) consists, by the adjunction rule, of those objects that are left orthogonal to \(\C\), and by Corollary~\reftwo{corollary:bousfield}, this agrees with the essential image of \(q\). So, the sequence formed by the adjoints satisfies \reftwoitem{item:fiber-left-to-right} (in the version with right adjoints), and is therefore also a Verdier sequence by what we have just shown.

On the other hand, suppose that \reftwoitem{item:cofiber-left-to-right} holds. Then \(g\) is a localisation by Corollary~\reftwo{corollary:bousfieldgenII} and the essential image of \(f\) is given by the right orthogonal of \(\ker(g)\) by Lemma~\reftwo{lemma:bousfieldgenI}. It is therefore, in particular, closed under retracts in \(\D\). But according to Corollary~\reftwo{corollary:criterion-projection}, \(p\) exhibits \(\E\) as the Verdier quotient of \(\D\) by this image so it equals \(\ker(p)\) by Lemma~\reftwo{lemma:kernel-of-projection-to-verdier-quotient}. This shows that the sequence in the statement is a fibre sequence. To see that \(p\) admits a left adjoint, we can appeal to Lemma~\reftwo{lemma:bousfieldgenI}: For \(X \in \D\) the fibre of the unit map \(X \rightarrow fg(X)\) clearly projects to \(p(X)\) under \(p\), and for \(C \in \C\), we have
\[
\Hom_\D\left(\fib(X \rightarrow fg(X)),f(C)\right) \simeq \cof\left[\Hom_\D(X,f(C)) \rightarrow \Hom_\D(fg(X),f(C))\right]
\]
and since \(f\) is fully faithful the latter term is also given by \(\Hom_\C(g(X),C)\), which identifies the map on the right as the adjunction equivalence. It follows that $\fib(X \rightarrow fg(X))$ is left local for the equivalences mod $\C$ and thus left adjoint to $p(X)$.
\end{proof}

\begin{corollary}
\label{corollary:split-verdier-incl-pushout}%
The collection of (right or left) split Verdier inclusion is closed under pushouts.
\end{corollary}

\begin{proof}
Given a diagram $\C' \leftarrow  \C \rightarrow \D$ with pushout $\D'$ whose right leg is a Verdier inclusion, the collection of objects of the pushout admitting left adjoint objects for $\C' \rightarrow \D'$ is evidently stable. It thus suffices to see that it contains the essential images of $\C' \rightarrow \D'$ and $\D \rightarrow \D'$ as these generate $\D'$. But the former of these is fully faithful by Corollary \reftwo{corollary:pullback-non-split-verd-proj}, so in this case there is nothing to do, and for the latter the projective adjoint $\Pro(\D') \rightarrow \Pro(\C')$ takes values $\C'$ once restricted to $\D$ by the adjointability statement from \reftwo{corollary:pullback-non-split-verd-proj}. 
\end{proof}

\begin{example}
\label{example:universal-split-verdier}%
The prototypical example of a split Verdier sequence is the sequence
\[
\C \xrightarrow{c \mapsto \id_c} \Ar(\C) \xrightarrow{\cof} \C
\]
for a stable \(\infty\)-category \(\C\), where \(\Ar(\C) = \Fun(\Del^1,\C)\) is the associated arrow category and \(\cof\) sends an arrow to its cofibre, so that \(\ker(\cof)\) consists of the equivalences. The fully faithful right and left adjoints of \(\cof\) are given by \(c \mapsto [0 \to c]\) and \([\Om c \to 0]\), respectively, and the right and left adjoints of \(\id_{(-)}\) send an arrow \(x \to y\) to \(x\) and \(y\), respectively. This split Verdier sequence is in fact the universal split Verdier sequence with kernel \(\C\), see Corollary~\reftwo{corollary:universal-split-verdier} below.
\end{example}

\begin{example}
\label{example:fun-split}%
If \(\C \to \D \to \E\) is a Verdier sequence (or in fact a Karoubi sequence as in the next section) then the sequence of restriction functors
\[
\Funx(\C,\Sp) \longleftarrow \Funx(\D,\Sp) \longleftarrow \Funx(\E,\Sp)
\]
is a split Verdier sequence (among large \(\infty\)-categories). Indeed, since \(\C \to \D \to \E\) is a cofibre sequence the sequence of restriction functors is a fibre sequence, but since \(\C\op \to \D\op\) is fully faithful the restriction functor on the left 
admits fully faithful left and right adjoints in the form of left and right Kan extensions (operations which preserve exact functors, see Lemma~\refone{lemma:kan-extension-exact-quadratic}),
and is hence a split Verdier projection by Corollary~\reftwo{corollary:criterion-split}.
\end{example}

\begin{example}
\label{example:ind-pro-split}%
If \(\C \to \D \to \E\) is a Verdier sequence then the sequences
\[
\Ind(\C) \lrar \Ind(\D) \lrar \Ind(\E) \quad\text{and}\quad  \Pro(\C) \lrar \Pro(\D) \lrar \Pro(\E)
\]
are right split and left split Verdier sequences (among large \(\infty\)-categories), respectively. 
Indeed, under the equivalence \(\Ind(\C\op) = \Funx(-,\Sp) = \Pro(\C)\op\) these examples are obtained from the split Verdier sequence of Example~\reftwo{example:fun-split} by passing to the corresponding sequence of left adjoints (before or after taking opposites) using Lemma~\reftwo{lemma:quasi-split}.
\end{example}

As a final method for producing split Verdier sequences, let \(\C \rightarrow \D \to \E\) be a (not necessarily split) Verdier sequence. Write 
\(\D^{\spl} \subseteq \D\) for the full subcategory spanned by those objects whose images under the inductive right adjoint $\D \rightarrow \Ind(\C)$ and projective left adjoints $\D \rightarrow \Pro(\C)$ lie in \(\C\). Then \(\D^{\spl}\) contains \(\C\) and we set \(\E^{\qspl} \coloneq \D^{\spl}/\C\), and we call the arising sequence the \emph{split core} of a Verdier sequence. To state its universal property fully, we finally let $\Fun^{\Ver}(\D\rightarrow \E,\D' \rightarrow \E')$ denote the full subcategory of the pullback of 
\[
\Funx(\D,\D') \longrightarrow \Funx(\C,\D') \longleftarrow \Funx(\C,\C')
\]
spanned by those functors giving rise to adjointable squares for any two Verdier sequences $\C \rightarrow \D \rightarrow \E$ and $\C' \rightarrow \D' \rightarrow \E'$, so that in particular
\[
\grpcr(\Fun^{\Ver}(\D\rightarrow \E,\D'\rightarrow \E')) \simeq \Hom_{\Ver}(\D\rightarrow \E,\D' \rightarrow \E').
\]

\begin{lemma}
\label{lemma:split-core}%
For every Verdier sequence $\C \rightarrow \D \rightarrow \E$ the sequence
\[
\C \lrar \D^{\spl} \longrightarrow \E^{\mathrm{q}\spl}
\]
is a split Verdier sequence. Furthermore, for every split Verdier sequence 
\(\C' \to \D' \to \E'\)
the map
\[
\Fun^{\Ver}(\D'\to\E',\D^{\spl}\to\E^{\mathrm{q}\spl})  \longrightarrow \Fun^{\Ver}(\D'\to\E',\D\to\E)
\]
is an equivalence, so that the formation of split cores constitues a right adjoint to the inclusion $\Ver^\perp \subseteq \Ver$.
\end{lemma}

In complete analogy one can also form right split and left split cores of a Verdier sequence, but these will not play a large role in the sequel.

\begin{proof}
Since \(i\colon \C \to \D\) is fully faithful the Pro-left and Ind-right adjoints of \(i\) send \(i(c)\) to \(c\) for every \(c \in \C\), and so \(\D^{\spl}\) contains \(\C\). In addition, the inclusion \(\C \to \D^{\spl}\) admits both a left and a right adjoint by construction and is hence a split Verdier inclusion by Corollary~\reftwo{corollary:criterion-split}, so that the split core is indeed a split Verdier sequence. To obtain the universal mapping property it suffices by Lemma~\reftwo{lemma:adj-verdier} to note that in any adjointable square
\[
\begin{tikzcd}
\C' \ar[r]\ar[d] & \C \ar[d] \\
\D'\ar[r, "f"] & \D \ ,
\end{tikzcd}
\]
with \(\C' \to\D'\) a split Verdier inclusion the map \(f\) must factor through \(\D^{\spl}\) by adjointability.
\end{proof}

Finally, let us record the following recognition principle for left and right split Verdier sequences, which we will use in Section~\reftwo{subsection:isotrop}. 
\begin{corollary}
\label{corollary:general-orthogonal}%
Let \(\D\) be a stable \(\infty\)-category and \(\C,\E \subseteq \D\) two full stable subcategories such that \(\map_{\D}(\x,\y) \simeq 0\) for every \(\x \in \C,\y \in \E\). Then the following are equivalent:
\begin{enumerate}
\item
\label{item:right-adj-incl-eq}%
\(\C \subseteq \D\) admits a right adjoint \(p\colon \D \to \C\) and the inclusion \(\E \subseteq \C^r\) is an equivalence.
\item
\label{item:verdier-incl-proj-eq}%
\(\E \subseteq \D\) is a Verdier inclusion and the projection \(\C \to \D/\E\) is an equivalence.
\item
\label{item:left-adj-incl-eq}%
\(\E \subseteq \D\) admits a left adjoint \(q\colon \D \to \E\) and the inclusion \(\C \subseteq \E^l\) is an equivalence.
\item
\label{item:verdier-incl-other-proj-eq}%
\(\C \subseteq \D\) is a Verdier inclusion and the projection \(\E \to \D/\C\) is an equivalence.
\end{enumerate}
Furthermore, when either of these equivalent conditions holds, the resulting sequences 
\[
\C\lrar \D\lrar \E\quad \text{and} \quad \E\lrar \D\lrar \C
\]
formed by the inclusions and their adjoints are right-split and left-split Verdier sequences, respectively.
\end{corollary}
\begin{proof}
The implications \reftwoitem{item:right-adj-incl-eq} \(\Rightarrow\) \reftwoitem{item:verdier-incl-proj-eq} and \reftwoitem{item:left-adj-incl-eq} \(\Rightarrow\) \reftwoitem{item:verdier-incl-other-proj-eq}  are dual to each other, and the same for the implications \reftwoitem{item:verdier-incl-proj-eq} \(\Rightarrow\) \reftwoitem{item:left-adj-incl-eq} and \reftwoitem{item:verdier-incl-other-proj-eq} \(\Rightarrow\) \reftwoitem{item:right-adj-incl-eq}. It will hence suffice to show \reftwoitem{item:right-adj-incl-eq} \(\Rightarrow\) \reftwoitem{item:verdier-incl-proj-eq} \(\Rightarrow\) \reftwoitem{item:left-adj-incl-eq}, along with the last claim.

To prove the first of these implications, suppose that \(i\colon\C \subseteq \D\) admits a right adjoint \(p\colon \D \to \C\) and that \(\E \subseteq \C^r\) is an equivalence. By the adjunction rule, \(\C^r\) agrees with the kernel of \(p\) so we have a right-split Verdier sequence
\[
\E\lrar \D\xrightarrow{\,p\,} \C
\]
from which we conclude that the map \(\D/\E\to \C\) induced by \(p\) is an equivalence. The projection \(\C\to \D/\E\) is a one-sided inverse and therefore also an equivalence. 

On the other hand, if \reftwoitem{item:verdier-incl-proj-eq} holds then by Corollary~\reftwo{corollary:bousfield}, the projection \(\D\to \D/\E\) has a left adjoint, and the inclusion of \(\C\) into \(\E^l\) is an equivalence (since both project to \(\D/\E\) by an equivalence); the existence of the left adjoint \(q\) follows from Lemma~\reftwo{lemma:quasi-split}.
\end{proof}

We now come back to the notion of a split Verdier sequence and show that it is essentially equivalent to that of a recollement in the sense of~\cite{HA}*{Section A.8} in the setting of stable \(\infty\)-categories. Specialising the definition to this case, we have:
\begin{definition}
\label{definition:stable-recollement}%
A stable \(\infty\)-category \(\D\) is a \emph{stable recollement} of a pair of stable subcategories \(\C\) and \(\E\) if 
\begin{enumerate}
\item the inclusions of both \(\C\) and \(\E\) admit left adjoints \(L_\C\) and \(L_\E\), 
\item the composite \(\C \rightarrow \D \xrightarrow{L_\E} \E\) vanishes, and
\item \(L_\E\) and \(L_\C\) are jointly conservative.
\end{enumerate}
\end{definition}

\begin{proposition}
\label{proposition:criterion-split}%
If \(\D\) is a stable recollement of \(\C\) and \(\E\), then the sequence \(\C\to \D\xrightarrow{L_\E} \E\) is a split Verdier sequence.
Conversely, if \(\C\xrightarrow{f} \D\xrightarrow{p} \E\) is a split Verdier sequence, then \(\D\) is a stable recollement of the essential images \(f(\C)\) and \(q(\E)\), where \(q\) denotes the right adjoint of \(p\). 
\end{proposition}

\begin{proof}
Consider the first statement. We claim that the sequence under consideration is a fibre sequence, so that it is split Verdier by Lemma~\reftwo{lemma:quasi-split}. Since the composite is zero by assumption, we are left to show that every object \(x\) in \(\ker(L_\E)\) already belongs to the essential image of \(\C\). Denoting by \(L_\C\) the left adjoint of the inclusion of \(\C\), then the unit \(x\to L_\C(x)\) is mapped to an equivalence under both \(L_\E\) and \(L_\C\). By assumption \(L_\E\) and \(L_\C\) are jointly conservative, so the unit \(x\to L_\C(x)\) is an equivalence and therefore \(x\) indeed lies in the essential image of \(\C\).

For the second statement, \(f\) admits a left adjoint \(g\) by Lemma~\reftwo{lemma:quasi-split} , since \(p\) does and it remains to see that \(p\) and \(g\) are jointly conservative. Since we are in the stable setting it will suffice to show that the functors \(p\) and \(g\) together detect zero objects. Indeed, if \(\x \in \D\) is such that \(p(\x)\simeq 0\) then \(\x\) belongs to the essential image of \(f\). In this case, if \(g(\x)\) is zero as well, then \(\x \simeq 0\) because the counit of \(g \dashv f\) is an equivalence. 
\end{proof}

Another definition of stable recollements is given by Barwick and Glasman in \cite{barwick2016note}. In pictures, the definitions of \cite{HA}*{Section A.8}, \cite{barwick2016note}, and that of split Verdier sequences above correspond to 
\[
\begin{tikzcd}
\C \ar[r] 
& \D \ar[r,"L_\E"] 
\ar[l,bend right=40,shift right=1.5ex,start anchor=west,end anchor=east,"L_\C"']
\ar[l,phantom,shift right=1.3ex,start anchor=west,end anchor=east,"\myperp"]
& \E \ar[l,bend left=40,shift left=1.5ex,start anchor=west,end anchor=east]
\ar[l,phantom,shift left=1.2ex,start anchor=west,end anchor=east,"\myperp"]
\end{tikzcd}\quad,\quad 
\begin{tikzcd}
\C \ar[r] 
& \D \ar[r,"p"] 
\ar[l,bend left=40,shift left=1.5ex,start anchor=west,end anchor=east]
\ar[l,phantom,shift left=1.2ex,start anchor=west,end anchor=east,"\myperp"]
\ar[l,bend right=40,shift right=1.5ex,start anchor=west,end anchor=east]
\ar[l,phantom,shift right=1.3ex,start anchor=west,end anchor=east,"\myperp"]
& \E
\end{tikzcd}\quad\text{and}\quad
\begin{tikzcd}
\C \ar[r,"f"] 
& \D \ar[r] 
& \E \ar[l,bend left=40,shift left=1.5ex,start anchor=west,end anchor=east]
\ar[l,phantom,shift left=1.2ex,start anchor=west,end anchor=east,"\myperp"]
\ar[l,bend right=40,shift right=1.5ex,start anchor=west,end anchor=east]
\ar[l,phantom,shift right=1.3ex,start anchor=west,end anchor=east,"\myperp"],
\end{tikzcd}
\]
respectively.
Our results above show that all of these types of diagrams can be completed to the full
\[
\begin{tikzcd}
[column sep=7ex]
\C \ar[r, "f" description] & 
\D 
\ar[r,"p" description] 
\ar[l,bend left=25,shift left=1.5ex,start anchor=west,end anchor=east,"r"] 
\ar[l,bend right=25,shift right=1.5ex,start anchor=west,end anchor=east,"q"'] 
\ar[l,phantom,shift left=1.2ex,start anchor=west,end anchor=east,"\myperp"] 
\ar[l,phantom,shift right=1.3ex,start anchor=west,end anchor=east,"\myperp"]
& \E,
\ar[l,bend left=25,shift left=1.5ex,start anchor=west,end anchor=east,"h"] 
\ar[l,bend right=25,shift right=1.5ex,start anchor=west,end anchor=east,"g"']
\ar[l,phantom,shift left=1.2ex,start anchor=west,end anchor=east,"\myperp"] 
\ar[l,phantom,shift right=1.3ex,start anchor=west,end anchor=east,"\myperp"]
\end{tikzcd}
\]
in which both the top and the bottom left pointing maps also form Verdier sequences. In particular, all three definitions are equivalent.

\begin{remark}
There are many examples of split Verdier sequences, or stable recollements, arising ``in nature'', often involving large \(\infty\)-categories. One example we established in \paperone is given by
\[
\begin{tikzcd}
[column sep=7ex]
\Funx(\C\op,\Spa) \ar[r] & 
\Funq(\C\op,\Spa) 
\ar[r] 
\ar[l,bend left=25,shift left=1.5ex,start anchor=west,end anchor=east] 
\ar[l,bend right=25,shift right=1.5ex,start anchor=west,end anchor=east,"{\Lin}"'] 
\ar[l,phantom,shift left=1.2ex,start anchor=west,end anchor=east,"\myperp"] 
\ar[l,phantom,shift right=1.3ex,start anchor=west,end anchor=east,"\myperp"]
& \Funs(\C\op,\Spa),
\ar[l,bend left=25,shift left=1.5ex,start anchor=west,end anchor=east,"{(-)^\sym}"] 
\ar[l,bend right=25,shift right=1.5ex,start anchor=west,end anchor=east,"{(-)^\qdr}"']
\ar[l,phantom,shift left=1.2ex,start anchor=west,end anchor=east,"\myperp"] 
\ar[l,phantom,shift right=1.3ex,start anchor=west,end anchor=east,"\myperp"]
\end{tikzcd}
\]
and a closely related example is that of genuine $\Ct$-spectra discussed in \S\reftwo{subsection:Mackey}. Many more examples of split Verdier sequences among large \(\infty\)-categories can be obtained as special cases of Example~\reftwo{example:fun-split}. Another standard example is the case where $\D = \Spa$ and $f$ is the inclusion of those spectra on which a prime $l$ acts invertibly:
\[
\begin{tikzcd}
[column sep=7ex]
\Spa [\oneoverl] \ar[r]  & 
\Spa \ar[r] 
\ar[l,bend left=25,shift left=1.5ex,start anchor=west,end anchor=east,"{\mathrm{div}_l}"] 
\ar[l,bend right=25,shift right=1.5ex,start anchor=west,end anchor=east,"{(-)[1/l]}"'] 
\ar[l,phantom,shift left=1.2ex,start anchor=west,end anchor=east,"\myperp"] 
\ar[l,phantom,shift right=1.3ex,start anchor=west,end anchor=east,"\myperp"]
& \Spa[\text{l-adic equiv's}^{-1}],
\ar[l,bend left=25,shift left=1.5ex,start anchor=west,end anchor=east,"{(-)^\wedge_l}"] 
\ar[l,bend right=25,shift right=1.5ex,start anchor=west,end anchor=east,"{(-)[l^\infty]}"']
\ar[l,phantom,shift left=1.2ex,start anchor=west,end anchor=east,"\myperp"] 
\ar[l,phantom,shift right=1.3ex,start anchor=west,end anchor=east,"\myperp"]
\end{tikzcd}
\]
with $X[l^\infty] \simeq \fib(X \rightarrow X[1/l])$ and $\mathrm{div}_l(X) = \fib(X \rightarrow X_l^\wedge)$.
\end{remark}

\subsection{Karoubi sequences}
\label{subsection:Karoubi}%

We now proceed to the more general notion of Karoubi sequences, which is a version of Verdier sequences invariant under idempotent completion.

\begin{definition}
\label{definition:karoubi-equivalence}%
We call an exact functor \(\C \to \D\) between stable \(\infty\)-categories a \defi{Karoubi equivalence} if it is fully faithful and has dense image, in the sense that every object of \(\D\) is a retract of an object in the essential image. 
\end{definition}

The most important example of Karoubi equivalences are of course idempotent completions \(\C \rightarrow \C^\natural\). When fixing the target, Karoubi equivalences can be entirely classified, see \cite{thomason-classification}*{Theorem 2.1}:

\begin{theorem}[Thomason]
\label{theorem:thomason-classification}%
Karoubi equivalences induce injections on \(\K_0 \colon \Catx \rightarrow \Ab\), and Karoubi equivalences to a fixed stable \(\infty\)-category \(\C\), that is, the poset of dense subcategories of \(\C\), is isomorphic to the poset of subgroups of \(\K_0(\C)\), by 
sending a dense subcategory \(\Ctwo \subseteq \C\) to the image of \(\K_0(\Ctwo)\) in \(\K_0(\C)\), and conversely sending a subgroup \(c \subseteq \K_0(\C)\) to the full subcategory on objects \(x\) with \([x]\in c\).
\end{theorem}

Note that the statement in \cite{thomason-classification} is for triangulated categories, but the proof works verbatim in the setting of stable \(\infty\)-categories.

\begin{definition}
Stable \(\infty\)-categories \(\C\) with the property that \(\K_0(\C)\) vanishes we will call \defi{minimal}, and refer to the assignment sending \(\C\) to the full subcategory \(\C^\mn \subseteq \C\) spanned by the objects \(x \in \C\) with \(0 = [x] \in \K_0(\C)\), as \defi{minimalisation}.
\end{definition}

\begin{proposition}
\label{proposition:minimal-cats}%
The localisation of \(\Catx\) at the collection of Karoubi equivalences is both a left and a right Bousfield localisation. The left and right local objects are the minimal and idempotent complete stable $\infty$-categories, respectively, and the left and right adjoints of the localisation at the Karoubi equivalences are given by 
$[\C] \mapsto \C^{\mn}$ and \([\C] \mapsto \C^\natural\), respectively.
Furthermore, an exact functor is a Karoubi equivalence if and only if it induces an equivalence on minimalisations, and if and only if it induces an equivalence on idempotent completions.
\end{proposition} 

Let \(\Catxi \subseteq \Catx\) denote the full subcategory spanned by the idempotent complete stable \(\infty\)-categories. We %
then conclude that \((-)^\natural \colon \Catx \rightarrow \Catxi\) exhibits $\Catxi$ as the localisation of $\Catx$ by the collection of Karoubi equivalences. In particular, by Proposition~\reftwo{proposition:minimal-cats} it is both a left and right Bousfield localisation, and hence preserves both limits and colimits.

\begin{proof}[Proof of Proposition~\reftwo{proposition:minimal-cats}]
It is an exercise in pasting retract diagrams to check that Karoubi equivalences are closed under 2-out-of-3. The characterisation in the last statement then follows immediately from the fact that both inclusions \(\C^\mn \subseteq \C \subseteq \C^\natural\) are Karoubi equivalences, the former since every \(X \in \C\) is a retract of \(X \oplus \Sigma X \in \C^\mn\). 

Together with \cite{HTT}*{Lemma 5.1.4.7} and the preservation of exactness by Kan extensions, this also means that given a Karoubi equivalence \(i \colon \C \rightarrow \D\) and a functor \(f \colon \D \rightarrow \E\), the exactness of \(f\) is equivalent to that of \(fi\). Combining this with \cite{HTT}*{Proposition 5.1.4.9} we hence conclude that in this case the restriction map \(i^*\colon \Funx(\D,\E) \to \Funx(\C,\E)\) is an equivalence if in addition \(\E\) is idempotent complete.

The statement about the adjoints now follows from Lemma~\reftwo{lemma:bousfieldgenI}: That idempotent completion satisfies the requisite conditions was just established 
and that minimalisations do is immediate from the functoriality of \(\K_0\).
\end{proof}

\begin{corollary}
\label{corollary:karoubi-equivalences-pullbacks}%
The collection of Karoubi equivalences is closed under pullbacks and pushouts.
\end{corollary}

Let us now define our main object of study in this section.
\begin{definition}
\label{definition:kar-seq}%
A sequence 
$
\C \xrightarrow{f} \D \xrightarrow{p} \E
$
of exact functors with vanishing composite 
is a \defi{Karoubi sequence} if the sequence
$
\C^\natural \rightarrow \D^\natural \rightarrow \E^\natural
$
is both a fibre and cofibre sequence in \(\Catxi\). In this case, we refer to \(f\) as a \defi{Karoubi inclusion} and to \(p\) as a \defi{Karoubi projection}. 
\end{definition}

\begin{remark}Equivalently, by Proposition~\reftwo{proposition:minimal-cats}, we might ask the sequence
\[
\C^\mn \rightarrow \D^\mn \rightarrow \E^\mn
\]
to be both a fibre and a cofibre sequence in the full subcategory of \(\Catx\) spanned by the minimal stable \(\infty\)-categories, or more symmetrically, that the original sequence give a fibre and cofibre sequence in the localisation of \(\Catx\) at the Karoubi equivalences. 

We have chosen the present formulation as the idempotent completion plays a disproportionally more important role, both in the detection of Karoubi sequences and in applications.
\end{remark}

We also have a concrete characterisation of Karoubi sequences, analogous to the one for Verdier sequences of Corollary~\reftwo{corollary:verdier-criterion}.
\begin{proposition}
\label{proposition:criterion-karoubi}%
Let \(\C \xrightarrow{f} \D \xrightarrow{p} \E\) be a sequence of exact functors between stable \(\infty\)-categories with vanishing composite. Then
\begin{enumerate}
\item 
\label{item:fiber-catex-idem}%
the sequence \(\C^\natural \to \D^\natural \to \E^\natural\)is a fibre sequence in \(\Catxi\) if and only if \(f\) becomes a Karoubi equivalence when regarded as a functor \(\C \rightarrow \ker(p)\).
\item 
\label{item:cofiber-catex-idem}%
the sequence \(\C^\natural \to \D^\natural \to \E^\natural\) is a cofibre sequence in \(\Catxi\) if and only if the induced functor from the Verdier quotient of \(\D\) by the stable subcategory generated by the image of \(f\) is a Karoubi equivalence to \(\E\).
\item 
\label{item:karoubi-catex-idem}%
the sequence \(\C \to \D \to \E\) is a Karoubi sequence if and only if \(f\) is fully faithful and the induced map \(\D/\C \to \E\) is a Karoubi equivalence. 
\end{enumerate}
In particular, every Verdier sequence is a Karoubi sequence.
\end{proposition}

Let us explicitly warn the reader, however, that the Verdier quotient of two idempotent complete, stable \(\infty\)-categories need not be idempotent complete.

\begin{proof}[Proof of Proposition~\reftwo{proposition:criterion-karoubi}]
By Proposition~\reftwo{proposition:minimal-cats}, the functor \((-)^\natural\colon \Catx\to \Catxi\) preserves both limits and colimits, and \(\Catxi\) is closed under limits in \(\Catx\). 
This yields an equivalence
\[
\ker(p^\natural) \simeq \ker(p)^\natural,
\]
which proves \reftwoitem{item:fiber-catex-idem}. 

Similarly, 
\reftwoitem{item:cofiber-catex-idem} follows from the description of cofibres in \(\Catx\) as Verdier quotients together with the preservation of cofibres under idempotent completion.

Finally, the forward direction of \reftwoitem{item:karoubi-catex-idem} follows directly from the previous two statements. On the other hand, if \(f\) is fully faithful, and \(\D/\C\to \E\) is a Karoubi equivalence, then the kernel of \(p\) agrees with the kernel of the projection \(q\colon \D \to \D/\C\). Thus, by Lemma~\reftwo{lemma:kernel-of-projection-to-verdier-quotient}, the map \(f\colon \C\to \ker(q)\) has dense essential image and therefore is a Karoubi equivalence. The reverse claim thus also follows from the first two statements.
\end{proof}

\begin{corollary}
\label{corollary:characterisation-of-Karoubi-inclusions-and-projections}%
An exact functor \(f\colon \C\to \D\) is a Karoubi inclusion if and only if it is fully faithful. It is a Karoubi projection if and only if it has dense essential image \(f(\C)\subseteq \D\), and the induced functor \(f\colon \C\to f(\C)\) is a Verdier projection.
\end{corollary}

Combining this statement with Thomason's result above, we find:

\begin{corollary}
\label{proposition:karoubi-surjective}%
Let \(p\colon \D \to \E\) be a Karoubi projection. Then the following are equivalent:
\begin{enumerate}
\item
\label{item:karoubi-is-verdier}%
\(p\) is a Verdier projection.
\item
\label{item:p-ess-surj}%
\(p\) is essentially surjective.
\item
\label{item:kzero-surj}%
The induced group homomorphism \(\K_0(\D) \to \K_0(\E)\) is surjective.
\end{enumerate}
\end{corollary}

\begin{corollary}
\label{corollary:Karoubi-proj-pullback}%
The collection of Karoubi projections is closed under pullbacks.
\end{corollary}

\begin{proof}
By Corollary~\reftwo{corollary:characterisation-of-Karoubi-inclusions-and-projections} the Karoubi projections are exactly the functors that can be written as a composite of a Verdier projection followed by a Karoubi equivalence. Both these classes of maps are closed under pullback by 
Corollaries \reftwo{corollary:pullback-non-split-verd-proj} and \reftwo{corollary:karoubi-equivalences-pullbacks}, respectively. 
\end{proof}

Next, we record the following detection criterion for Karoubi-sequences, often called the Thomason-Neeman localisation theorem in the context of triangulated categories, see \cite{neeman-localisation}*{Theorem 2.1}. 

\begin{theorem}
\label{theorem:indkaroubi}%
A sequence \(\C \rightarrow \D \rightarrow \E\) of stable \(\infty\)-categories and exact functors with vanishing composite is a Karoubi sequence if and only if the induced sequence
\[
\Ind(\C) \lrar \Ind(\D) \lrar \Ind(\E)
\]
is a Verdier sequence (of not necessarily small \(\infty\)-categories).
In addition, in this case the above sequence is automatically right split and the corresponding sequence of right adjoints
\[
\Ind(\C) \longleftarrow \Ind(\D) \longleftarrow \Ind(\E)
\]
is a split Verdier sequence, i.e.\ both functors admit further right adjoints.
\end{theorem} 

\begin{proof}
 
By Proposition~\cite{HTT}*{5.3.5.11} the functor \(\Ind\) preserves full faithfulness and sends Karoubi equivalences to equivalences. At the same time, by Proposition~\reftwo{proposition:criterion-karoubi} every Karoubi sequence relates via a zig-zag of levelwise Karoubi equivalences to a Verdier sequence, and hence \(\Ind(-)\) takes Karoubi sequences to split Verdier sequences by combining Examples~\reftwo{example:fun-split} and \reftwo{example:ind-pro-split}. 

In the other direction, since the Yoneda functor \(\C \to \Ind(\C)\) is fully faithful we have that \(\Ind(-)\) also detects full faithfulness, and by \cite{HTT}*{Lemma 5.4.2.4} 
\(\Ind(-)\) also detects Karoubi equivalences, that is, a functor is a Karoubi equivalence if and only if it is sent to an equivalence by \(\Ind\). In particular, if \(\C \to \D \to \E\) is now a sequence which becomes Verdier after inductive completion then the \(\C \to \D\) is fully faithful and the induced functor \(\D/\C \to \E\) is an equivalence on inductive completions and hence a Karoubi equivalence. We conclude that \(\C \to \D \to \E\) 
is a Karoubi sequence by Proposition~\reftwo{proposition:criterion-karoubi}.
\end{proof}

\begin{remark}
To avoid confusion, note that while the characterisation of Theorem~\reftwo{theorem:indkaroubi} works on the level of inclusions alone (\(f\) is a Karoubi inclusion if and only if \(\Ind(f)\) is a Verdier inclusion), it does not work on the level of projections: there exist exact functors \(\D \to \E\) which are sent to Verdier projections under inductive completion, but that are not themselves Karoubi projections: The kernel of $\Ind(\D) \rightarrow \Ind(\E)$ may fail to be of the form $\Ind(\C)$, see, e.g.\ Example~\reftwo{example:ore} below. The kernel is, however, automatically a dualisable stable $\infty$-category, and can be taken into account by using these as a foundational set-up for localising Karoubi-localising invariants in place of small stable $\infty$-categories. This framework has recently been gaining traction as a preferred setup for (non-connective) \(\K\)-theory, topological cyclic homology and the like, see \cite{Efimov}.
\end{remark}

\subsection{Verdier sequences among $\infty$-categories of modules}
\label{section:appendix-module-examples}%

Let \(\phi\colon A\to B\) be a map of \(\Eone\)-ring spectra. Extension of scalars induces an exact functor
\[
\phi_!\colon \Mod(A) \lrar \Mod(B), \quad M \mapsto B\ssmash_A M
\]
on the \(\infty\)-categories of (left) modules, which is left adjoint to the restriction of scalars functor \(\phi^*\colon \Mod(B)\to \Mod(A)\). Extension of scalars restricts to functors
\[
{\phi}_!\colon\Modp{A} \lrar \Modp{B} \quad \text{and} \quad {\phi}_!\colon\Modc{\mathrm{c}}{A} \lrar \Modc{\phi(\mathrm{c})}{B},
\]
where \(\mathrm c \subseteq \K_0(A)\) is a subgroup and \(\Modc{\mathrm{c}}{A}\) the full subcategory of \(\Modp{A}\) spanned by those \(A\)-modules \(X\) with \([X] \in \mathrm c \subseteq \K_0(A)\). The most important special case of the latter construction is the case where \(\mathrm c\) is the image of the canonical map \(\ZZ \rightarrow \K_0(A), 1 \mapsto A\), in which case \(\Modc{\mathrm{c}}{A} = \Modf{A}\) is the stable subcategory of \(\Modp{A}\) generated by \(A\). In this section we analyse when these functors are Verdier or Karoubi projections.

Now let \(\Mod(A)_B \subseteq \Mod(A)\) denote the kernel of the functor \(\phi_!\colon\Mod(A) \to \Mod(B)\) extending scalars. For the next lemma, we note that the restriction of scalars functor \(\phi^*\colon \Mod(B)\to \Mod(A)\) admits $\hom_A(B,-)$ as a right adjoint, which we denote by \(\phi_*\colon \Mod(A) \to \Mod(B)\).
\begin{lemma}
\label{lemma:unital}%
Let \(\phi \colon A \to B\) be a map of \(\Eone\)-ring spectra and denote by \(I\) the fibre of \(\phi\), considered as an \(A\)-bimodule. Then the following are equivalent:
\begin{enumerate}
\item
\label{item:unital-equivalence}%
The multiplication \(B \ssmash_A B \rightarrow B\) is an equivalence.
\item
\label{item:unital-fibre}%
We have \(B \otimes_A I \simeq 0\).
\item
\label{item:idempotent-kernel}%
The map \(I \otimes_A I \to A \otimes_A I = I\) is an equivalence.
\item
\label{item:unital-bousfield}%
\label{item:unital-verdier}%
The diagram
\[
\begin{tikzcd}
[column sep=10ex]
\Mod(B)
\ar[r, "\phi^*" description]
& \Mod(A) 
\ar[r, "I \otimes_A (-)" description]
\ar[l,bend left=35,shift left=1.5ex,start anchor=west,end anchor=east,"{\phi_*}"] 
\ar[l,bend right=35,shift right=1.5ex,start anchor=west,end anchor=east,"{\phi_!}"'] 
\ar[l,phantom,shift left=2ex,start anchor=west,end anchor=east,"\myperp"] 
\ar[l,phantom,shift right=2ex,start anchor=west,end anchor=east,"\myperp"]
& \Mod(A)_B 
\ar[l,bend left=35,shift left=1.5ex,start anchor=west,end anchor=east,"{\hom_A(I,-)}"]
\ar[l,bend right=35,shift right=1.5ex,start anchor=west,end anchor=east,"{\incl}"']
\ar[l,phantom,shift left=2ex,start anchor=west,end anchor=east,"\myperp"] 
\ar[l,phantom,shift right=2ex,start anchor=west,end anchor=east,"\myperp"]
\end{tikzcd}
\]
\end{enumerate}
is a stable recollement.
\end{lemma}

Note that \reftwoitem{item:unital-bousfield} in particular contains the statement that \(I \otimes_A - \colon \Mod(A) \rightarrow \Mod(B)\) has image in \(\Mod(A)_B\) as indicated. Of course, one may as well replace \(\Mod(A)_B\) in the statement by the kernel of \(\phi_*\) and extrapolating from the example \(B = A[s^{-1}]\) one might call such modules \(\phi\)-complete: In this case the lower adjoint becomes the inclusion, the right pointing map \(\hom_A(I,-)\) and the top adjoint \(I \otimes_A -\).

\begin{proof}
For the equivalence between the first two items, simply note that
\[
B \simeq B \otimes_A A \xrightarrow{\id \otimes \phi} B \otimes_A B
\]
is always a right inverse to the multiplication map of \(B\). So the latter is an equivalence if and only if the fibre \(B \otimes_A I\) of the former vanishes. Similarly, \(B \otimes_A I\) is also the cofibre of the map in~\reftwoitem{item:idempotent-kernel}, and so the first three items are equivalent. 
The statement of \reftwoitem{item:unital-verdier} contains \reftwoitem{item:unital-fibre}, since \( I = I \otimes_A A \in \Mod(A)_B\). Finally, assuming the first two items, we first find that 
\[
B \otimes_A I \otimes_A X \simeq 0
\]
so that \(I \otimes_A X \in \Mod(A)_B\) for all \(X \in \Mod(A)\) and the diagram in \reftwoitem{item:unital-bousfield} is well-defined. 
Furthermore, it follows that \(\phi^*\) is fully faithful: For this one needs to check that the counit transformation \(B \otimes_A Y \rightarrow Y\) is an equivalence for every \(B\)-module \(Y\). 
But as both sides preserve colimits and \(B\) generates \(\Mod(B)\) under colimits and shifts, it suffices to check this for \(Y=B\), where we have assumed it. 

We have thus established that under the first two equivalent items the functor \(\phi^*\) is fully faithful. Since it admits both a left and a right adjoint it is a split Verdier inclusion by Corollary~\reftwo{corollary:criterion-split}.
It then follows from 
Proposition~\reftwo{proposition:criterion-split} that the diagram 
\[
\begin{tikzcd}
[column sep=7ex]
\Mod(B) \ar[r] 
& \Mod(A) 
\ar[l,bend left=30,shift left=1.5ex,start anchor=west,end anchor=east,"{\phi_*}"] 
\ar[l,bend right=30,shift right=1.5ex,start anchor=west,end anchor=east,"{\phi_!}"'] 
\ar[l,phantom,shift left=1.2ex,start anchor=west,end anchor=east,"\myperp"] 
\ar[l,phantom,shift right=1.2ex,start anchor=west,end anchor=east,"\myperp"]
\end{tikzcd}
\]
can be completed to a stable recollement and the fibre sequences connecting the various adjoints are easily checked to give the formulae from the statement.
\end{proof}

\begin{definition}
\label{definition:solid}%
We will call a map \(\phi\colon A \to B\) of \(\Eone\)-ring spectra satisfying the equivalent conditions of the previous lemma a \defi{\solid}. 

A map \(R \rightarrow S\) between discrete rings is called a \defi{derived localisation} if the associated map \(\GEM R \rightarrow \GEM S\) is a localisation in the sense above.
\end{definition}

Before we give examples, we introduce one more bit of terminology that we will need below. Given a full subcategory \(\C \subseteq \Mod(A)\), let us write \(\C_B = \C \cap \Mod(A)_B\). For the following definition, note that by Lemma~\reftwo{lemma:unital}, a map \(\phi \colon A \rightarrow B\) is a \solid if and only if its fibre \(\fib(A \to B)\) lies in \(\Mod(A)_B\).

\begin{definition}
\label{definition:perfectly-generated-fibre}%
We say that a \solid \(\phi\colon A \rightarrow B\) has perfectly generated fibre if 
\(I := \fib(A \to B) \in \Mod(A)_B\) lies in the smallest subcategory of \(\Mod(A)_B\) containing \((\Modp{A})_B\) and closed under colimits. 
\end{definition}

\begin{examples}
\label{remark:examples-discrete}%
\label{example:ore}%
\label{example:ore-2}%
\begin{enumerate}
\item
\label{item:locisderloc}%
If \(A\) is an \(\Eone\)-ring and \(S \in \pi_*(A)\) is a subset of homogeneous elements satisfying the left or right Ore condition then \(A \to A[S^{-1}]\) is a \solid with perfectly generated fibre. Indeed, it is a \solid by Lemma~\reftwo{lemma:unital} since the forgetful functor \(\Mod({A[S^{-1}]}) \to \Mod(A)\) is fully faithful,
and under the Ore condition, the fibre \(\Mod(A)_S \coloneq \Mod(A)_{A[S^{-1}]}\) is generated under colimits and shifts by the perfect modules \(A/s = \cof[A \xrightarrow{s} A]\) for \(s \in S\), see~\cite{HA}*{Lemma 7.2.3.13}. 
\item 
The discrete counterpart of Definition~\reftwo{definition:solid} for ordinary rings and ordinary tensor products was studied by Bousfield and Kan in~\cite{bk-solid}, where they have classified all commutative rings \(R\) whose multiplication \(R \otimes_{\ZZ} R \rightarrow R\) is an isomorphism, a property which is called \defi{solidity} in \cite{bk-solid}. We note that for a map of connective \(\Eone\)-rings \(A \rightarrow B\), being a \solid implies the solidity of \(\pi_0A \rightarrow \pi_0B\) but even for discrete \(A\) and \(B\), the converse is not automatic, as it additionally entails that \(\Tor_i^A(B,B) = 0\) for all \(i > 0\). Note that for non-commutative ordinary rings this condition is not even automatic for localisations.
\item 
If $R$ is a discrete commutative ring and $U \subseteq \spec(R)$ is a quasi-compact open subset then it follows from \cite{stacks}*{Lemma 0FLQ} that the map $R \to \RR\Gamma(U,\cO_R)$ is a derived \solid. More generally, if $U$ is an intersection of quasi-compact open subsets subsets of $\spec(R)$ then $R \to \RR\Gamma(U,\cO_R) = \colim_{V \in \mathcal U} \RR\Gamma(V,\cO_R)$ is again a \solid, where $\mathcal U$ is the (filtered) family of all quasi-compact open subsets containing $U$. By~\cite{thomason-classification}*{Theorem 3.15} these are exactly all the localisations of $R$ with compactly generated fibre, a result previously proven in \cite{hopkins-global, neeman-chromatic} for Noetherian rings. In fact, in the latter case all localisations of $R$ have perfectly generated fibre, and are hence all accounted for by this construction, see \cite{neeman-chromatic}*{Theorem 3.4}. 
\item
\label{item:almost-math}%
Besides these examples Falting's almost ring theory provides another important case of derived localisations, namely the projection $R \rightarrow R/J$ for a flat and idempotent ideal $J$ in a commutative ring $R$. The most prominent example of this situation is the ideal $J$ of topologically nilpotent elements inside the ring $R$ of power bounded elements in a perfectoid field, see e.g.\ \cite{Bhatt}. Note that these examples have highly non-finitely generated $J$, and indeed Nakayama's lemma implies that a finitely generated idempotent ideal $J$ is principal on an idempotent element $e$. In this case $R/J \cong R[(1-e)^{-1}]$. Examples with non-principal $J$ usually do not have perfectly generated fibre.
This is the case, for example, if \(R = k[t,t^{1/2},t^{1/4},...]\) is %
obtained from the polynomial algebra \(k[t]\) over a field $k$ by adding all iterated square roots of \(t\), and \(J \subseteq R\) is the ideal generated by the \(t^{1/2^i}\) for \(i\geq 0\), see~\cite{keller1994remark}; a quick argument using Thomason's classification is as follows: The $\Einf$-rings $\mathbb R\Gamma(U,\mathcal O_R)$ for $U \subset \spec(R)$ occuring in it all inherit the property from open affine subsets that tensoring with them preserves coconnectivity. In particular, if discrete they must be flat as $R$-modules. But in the example $k$ is not flat as an $R$-module as $\Tor_1^R(k,(t)) = k$. 

A generalisation of this set-up was recently developed in \cite{almostpaper}, showing that to every idempotent ideal $J \subseteq \pi_0(A)$ in a connective $\Eone$-ring $A$ corresponds a unique localisation $A/J^\infty$ of $A$ which is also connective, and for which $\pi_0(A) \rightarrow \pi_0(A/J^\infty)$ is surjective with kernel $J$, in particular removing the flatness assumption on $J$ entirely.
\item 

The maps $A \rightarrow \mathrm{L}_{n,p}A$ and $A \rightarrow \mathrm L^{\mathrm{fin}}_{n,p} A$ in chromatic homotopy theory are localisations for every $\Eone$-ring spectrum; here, $\mathrm{L}_{n,p}$ denotes the localisation at $\bigoplus_{i=0}^n \mathrm K(i,p)$, where $\mathrm K(i,p)$ denotes the Morava K-spectrum at height $i$ and prime $p$ (equivalently it is the localisation at a Lubin-Tate spectrum at height $n$ and prime $p$) and $\mathrm L^{\mathrm{fin}}_{n,p}$ is similarly the localisation at $\bigoplus_{i=0}^n \mathrm T(i,p)$ for some choice of $v_i$-telescope $\mathrm T(i,p)$ at the prime $p$ (the localisation is independent of the exact choices), see e.g.\ \cite{ravenel-book, miller-finite}. The fibre of $A \to \mathrm L^{\mathrm{fin}}_{n,p} A$ is always perfectly generated, and for $A = \mathbb S$ it is a theorem of Hopkins and Smith \cite{nilpotenceII} that all localisations of $\mathbb S$ with perfectly generated fibre are tensor products of the $\mathrm L^{\mathrm{fin}}_{n,p} \mathbb S$ %
over different primes. Since the compact objects in the categories of $\mathrm L_{n,p}$- and $\mathrm L_{n,p}^\mathrm{fin}$-acyclic spectra %
agree, the disproof of the telescope conjecture in \cite{telescopeisfalse} shows that $\mathbb S \rightarrow \mathrm L_{n,p}\mathbb S$ does not have perfectly generated fibre unless $n=0,1$.

\end{enumerate}
\end{examples}

\begin{proposition}
\label{proposition:app-verdier-free}%

Let \(\phi \colon A \rightarrow B\) be a \solid of \(\Eone\)-rings with perfectly generated fibre. Then 
\[
\Modp{A} \xrightarrow{\phi_!} \Modp{B} \quad \text{ and } \quad \Modc{\mathrm{c}}{A} \xrightarrow{\phi_!} \Modc{\phi(\mathrm{c})}{B}
\]
are Karoubi and Verdier projections for every $\mathrm{c} \subseteq \K_0(A)$, respectively.
\end{proposition}
\begin{proof}
Combining Theorem~\reftwo{theorem:indkaroubi} and Lemma~\reftwo{lemma:unital}, it only remains to show that \(\Ind((\Modp{A})_B) \simeq \Mod(A)_B\) to obtain the first claim. But by \cite{HTT}*{Proposition 5.3.5.11}, the former term is equivalent to the smallest subcategory of \(\Mod(A)_B\) containing \((\Modp{A})_B\) and closed under colimits, so by assumption it contains \(I\). But the smallest stable subcategory of \(\Mod(A)\) containing \(I\) and closed under colimits is \(\Mod(A)_B\), as follows immediately from the stable recollement of \reftwoitem{item:unital-bousfield} (since \(A\) generates \(\Mod(A)\) under colimits and shifts). 

Since the inclusion \(\Modc{\mathrm c}{A} \rightarrow \Modp{A}\) is a Karoubi equivalence and similarly for \(B\), it follows that \(\phi_!\colon\Modc{\mathrm c}{A} \to \Modc{\phi(\mathrm c)}{B}\) is a Karoubi projection as well. But the essential image of this functor is then the Verdier quotient by its kernel, see Corollary~\reftwo{corollary:characterisation-of-Karoubi-inclusions-and-projections}, and therefore a dense stable subcategory of \(\Modp{B}\). The second claim follows from the classification of dense subcategories \reftwo{theorem:thomason-classification} and Proposition~\reftwo{proposition:karoubi-surjective}.
\end{proof}

Combining Proposition~\reftwo{proposition:app-verdier-free} with Example~\reftwo{example:ore} we get:

\begin{corollary}
\label{corollary:ore-2}%

Given an \(\Eone\)-ring \(A\) and a subset \(S \in \pi_*(A)\) of homogeneous elements satisfying the left Ore condition, for example \(\pi_*(A)\) could be (skew-)commutative, then 
\[
\Modp{A} \xrightarrow{-[S^{-1}]} \Modp{A[S^{-1}]} \quad \text{ and } \quad \Modc{\mathrm c}{A} \xrightarrow{-[S^{-1}]} \Modc{\im(\mathrm c)}{A[S^{-1}]}
\]
are Karoubi and Verdier projections for every $\mathrm{c} \subseteq \K_0(A)$, respectively.
\end{corollary}

\begin{remark}
\label{example:app-Verdier-not-karoubi}%
In the situation of Proposition~\reftwo{proposition:app-verdier-free}, the map \(\Modp{A} \to \Modp{B}\) can fail to be a Verdier projection. For example, suppose that \(k\) is a field and \(R \subseteq k[s]\) is the nodal curve we have also considered in Remark~\reftwo{remark:no-verdier}. There we noted that the fundamental theorem in K-theory implies that the cokernel of the map $\K_0(R[t]) \to \K_0(R[t^{\pm 1}])$ contains $\K_{-1}(R) \cong \ZZ$ as a direct summand, showing that $\Dperf(R[t]) \to \Dperf(R[t^{\pm1}])$ is not a Verdier projection.
\end{remark}

\subsection{The classification of Verdier sequences}
\label{section:classverdier}%

Our goal in this subsection is to develop the structure theory of Verdier, split Verdier and Karoubi sequences.
In the split case the theory is rather well-documented, often from the point of view of classifying recollements, see e.g.\ \cite{HA}*{\S A.8}, \cite{barwick2016note} or~\cite{ShahQuigley}.
For Verdier and Karoubi sequences this is new, originally conceived of by the eighth author.

To begin, let us extend the $\infty$-category $\mathrm{Ver}$ to an $\infty$-category $\textsc{Ver}$ consisting of Verdier sequences (that is, fibre-cofibre sequences) $\C \rightarrow \D \rightarrow \E$ of locally small stable $\infty$-categories in which $\C$ is assumed small. 
We will call such sequences \defi{locally small Verdier sequences}, and similarly call locally small Vedier projections/inclusions for the projections and inclusions appearing in such sequences (so that by definition a locally small Verdier inclusion has a small domain and a locally small Verdier projection has a small kernel). To apply the theory in this generality, we first note that by a direct adaptation of Proposition~\reftwo{proposition:verdier-localisation} 
a Verdier quotient of a locally small stable $\infty$-category by a small stable subcategory is again locally small, and similarly by Remark~\reftwo{remark:kan-along-verdier-projection} we have that if $p\colon \D \to \E$ is a locally small Verdier projection then left Kan extensions along $p$ exist for all functors $p\colon \D \to \A$ with $\A$ locally small and possessing filtered colimits, and similarly for right Kan extensions and $\A$ possessing cofiltered limits. 
Corollary~\reftwo{corollary:criterion-projection} then applies as stated when $\D$ and $\E$ are locally small and $\ker(p)$ is assumed small.
Taking the characterisation of Remark~\reftwo{remark:adj-concrete}~\reftwoitem{item:equivalent-def} as the definition of horizontal Ind- and Pro-adjointability we can then extend these notions to squares of locally small stable $\infty$-categories whose horizontal arrows are locally small Verdier inclusions or projections.
Defining adjointability as the conjunction of vertical Ind- and Pro-adjointability, Lemma~\reftwo{lemma:adj-verdier} then still holds, and one can define morphisms in $\textsc{Ver}$ as adjointable diagrams analogously to $\mathrm{Ver}$. Finally, let us also note that with these definitions Proposition~\reftwo{proposition:verdier-square-char}
applies verbatim when the $\infty$-categories in the first part of its statement are only assumed locally small, as long as the vertical fibres are assumed small, and in the second part the lower row is only assumed locally small.

Now given a small stable \(\infty\)-category \(\cB\), there are associated two locally small stable $\infty$-categories $\Tate(\cB)$ and $\Latt(\cB)$, abstracting the classical notions of Tate vector spaces and lattices therein, originally constructed in \cite{Hennion-Tate}. We define $\Tate(\cB)$ to be the smallest stable subcategory of $\Ind\Pro(\cB)$ spanned by $\Pro(\cB)$ and $\Ind(\cB)$, and $\Latt(\cB)$ to be the full subcategory of $\Ar(\Ind\Pro(\cB))$ spanned by the arrows starting in $\Ind(\cB)$ and ending in $\Pro(\cB)$. 
Let us warn the reader that Hennion works in the realm of idempotent complete stable $\infty$-categories and therefore considers the idempotent completion of $\Tate(\cB)$ throughout; what we have called $\Tate(\cB)$ is said to consist of the elementary Tate objects in \cite{Hennion-Tate}. At any rate, assigning to an object $b \in \cB$ the identity arrow $\id_b \in \Ar(\cB) \subseteq \Latt(\cB)$ gives a fully faithful functor $\cB \rightarrow \Latt(\cB)$. Sending an arrow to its cofibre then gives a functor $\Latt(\cB) \rightarrow \Tate(\cB)$ which vanishes on the image of $\cB$. The goal for the remainder of this section is to prove and exploit the following result: 

\begin{theorem}
\label{theorem:universal-verdier}%
\label{corollary:structure-square-verdier}%
For every (small) stable $\infty$-category \(\cB\) the sequence
\[
\cB^\natural \longrightarrow \Latt(\cB) \longrightarrow \Tate(\cB)
\]
is a Verdier sequence, and for every locally small Verdier sequence \(\C \to \D\to \E\) the resulting map
\[
\Fun^{\textsc{Ver}}(\D \to \E,\Latt(\cB) \to \Tate(\cB)) \longrightarrow \Funx(\C,\cB^\natural)
\]
extracting kernels is an equivalence of (small) stable $\infty$-categories. Consequently, there 
exists an essentially unique adjointable square
\[
\begin{tikzcd}
\D \ar[r,"\vphi"]\ar[d, "p"'] & \Latt(\C) \ar[d] \\
\E \ar[r,"\psi"] & \Tate(\C) \ ,
\end{tikzcd}
\]
inducing the inclusion $\C \rightarrow \C^\natural$ on vertical fibres. This square is cartesian if and only if $\C$ is idempotent complete, so that in this case the Verdier projection \(p\) is pulled back from the universal Verdier projection with fibre $\C$ on the right. 
\end{theorem}

The exact functor \(\psi\colon \E \to \Tate(\C)\) is called the \defi{classifying functor} of the given Verdier sequence, and the square its \defi{classifying square}. Theorem~\reftwo{theorem:universal-verdier} then implies the following classification result for Verdier projections (among small stable $\infty$-categories), for the statement of which we note that $\K(\Tate(\C)) \simeq \mathbb S^1 \otimes \K(\C^\natural)$ since $\Latt(\C)$ admits an Eilenberg swindle, compare \cite{Hennion-Tate}*{Corollary 4.3}:
	
\begin{corollary}
\label{corollary:classverdierfull}%
Given stable $\infty$-categories $\C$ and $\E$, extracting classifying functors and boundary maps on $\K$-spectra induces an equivalence between the space of Verdier sequences $\C \rightarrow \D \rightarrow \E$ and the space of pairs $(\psi,\eta)$ where $\psi \colon \E \rightarrow \Tate(\C)$ is an exact functor and $\eta \colon \K(\E) \rightarrow \mathbb S^1 \otimes \K(\C)$ is a lift of
\[
\K(\E) \xrightarrow{\psi_\ast} \K(\Tate(\C)) \xrightarrow{\partial} \mathbb S^1 \otimes \K(\C^\natural).
\]
Equivalently, it is the space of pairs $(\psi,s)$ where $\psi$ is as above, such that in the exact sequence 
\[
\K_0(\C^\natural)/\K_0(\C) \longrightarrow \K_0\big(\Latt(\C) \times_{\Tate(\C)} \E\big)/\K_0(\C) \longrightarrow \K_0(\E)\longrightarrow 0
\]
the left hand map is injective and $s$ is a splitting of it as abelian groups.
\end{corollary}

A version for Verdier sequences among stably dualisable $\infty$-categories will appear in \cite{Nikolaus-Verdier}, and a similar result was also obtained independently by Efimov, see \cite{Efimov}*{Theorem 3.4}.

\begin{proof}
Theorem \reftwo{theorem:universal-verdier} implies that the inverse construction to taking classifying functors is given by sending $(\psi, \eta)$ with $\psi \colon \E \rightarrow \Tate(\C)$ and $\eta \colon \K(\E) \rightarrow \mathbb S^1 \otimes \K(\C)$ to
\[
\D = \{X \in \Latt(\C) \times_{\Tate(\C)} \E \mid [X] \in \pi_0(\fib(\eta))\},
\]
equipped with the evident inclusion and projection from $\C$ and to $\E$, respectively, where $\pi_0(\fib(\eta))$ is regarded as a subgroup of $\K_0(\Latt(\C) \times_{\Tate(\C)} \E)$ via
\[
\pi_0\fib(\eta \colon \K(\E) \rightarrow \mathbb S^1 \otimes \K(\C)) \subseteq \pi_0\fib(\psi_\ast \colon \K(\E) \rightarrow \mathbb S^1 \otimes \K(\C^\natural)) = \K_0(\Latt(\C) \times_{\Tate(\C)} \E).
\]
But through this description it is also equivalent to the space of pairs $(\psi, \widetilde\eta)$ where $\psi$ is as before and $\widetilde \eta$ is a null-homotopy of the composite
\[
\widetilde \psi \colon \K(\E) \xrightarrow{\psi} \K(\Tate(\C)) = \mathbb S^1 \otimes \K(\C^\natural) \longrightarrow \mathbb S^1 \otimes \K(\C^\natural)/\K(\C) = \mathbb S^1 \otimes \GEM(\K_0(\C^\natural)/\K_0(\C)),
\]
where the final equality follows from the cofinality theorem.
On general grounds the space of such nullhomotopies is discrete with components corresponding to splittings of the induced sequence
\[
\K_0(\C^\natural)/\K_0(\C) \longrightarrow \pi_0\fib(\widetilde\psi) \longrightarrow \pi_0\K(\E)
\]
which identifies with the sequence from the statement.
\end{proof}

\begin{remarks}
\label{remark:concrete}%
\
\begin{enumerate}
\item The functor \(\map_\cB \colon \cB\op \times \cB \to \Sp\) extends in an essentially unique way to a functor 
\[
\wtl{\map}_\cB \colon \Ind(\cB)\op \times \Pro(\cB) \lrar \Sp
\]
which preserves limits in each variable. We claim that $\Latt(\cB)$ coincides with the associated pairing \(\infty\)-category $\Pairings(\Ind(\cB),\Pro(\cB),\wtl{\map}_\cB)$, that is, with the bivariant unstraightening of the space valued functor $\Om^{\infty}\wtl{\map}_{\cB}(-,-)$, see \S 7.1-7.2 in \paperone. Indeed, we generally have $\Pairings(\C,\C,\hom_\C) \simeq \Ar(\C)$ for stable $\C$, see Example 7.1.1 in \paperone,
and applying this twice we find cartesian squares
\[
\begin{tikzcd}
\Ar(\cB^\natural) \ar[r] \ar[d,"{(s,t)}"] & \Pairings(\Ind(\cB),\Pro(\cB),\wtl{\map}_\cB) \ar[r]\ar[d] & \Ar(\Tate(\cB)) \ar[d,"{(s,t)}"] \\
\cB^\natural \times \cB^\natural \ar[r] & \Ind(\cB) \times \Pro(\cB) \ar[r] & \Tate(\cB) \times \Tate(\cB)
\end{tikzcd}
\]
the right one of which implies the desired identification.
\item 
The Ind-right and Pro-left adjoints of \(i \colon \cB \rightarrow \Latt(\cB)\) send an Ind-to-Pro arrow \(f\colon x \to y\) to \(x\) and \(y\), respectively: This follows from the fact that the right and left adjoints of \(\id_{(-)}\colon \Ind\Pro(\cB) \to \Ar(\Ind\Pro(\cB))\) are given by the source and target projections, respectively, and that this projection maps \(\Latt(\cB) \subseteq \Ar(\Ind\Pro(\cB))\) to \(\Ind(\cB)\) and \(\Pro(\cB)\), respectively.
\item Similarly, the projection $\Latt(\cB) \rightarrow \Tate(\cB)$ admits a partial left adjoint on $\Ind(\cB) \subseteq \Tate(\cB)$, taking $z$ to $\Omega z \rightarrow 0$ and a partial right adjoint on $\Pro(\cB) \subseteq \Tate(\cB)$ taking $z$ to $0 \rightarrow z$. 
\end{enumerate}
\end{remarks}

\begin{proof}[Proof of Theorem \reftwo{theorem:universal-verdier}]
The proof naturally splits into two independent parts:
We first verify that $\cB \rightarrow \Latt(\cB) \rightarrow \Tate(\cB)$ is indeed a Verdier sequence. This is a special case of Clausen's discussion of $\infty$-categories of cones, see \cite{Clausen-Artin}*{Section 3.1} and particularly \cite{Clausen-Artin}*{Remark 3.23}. For the convenience of the reader we briefly summarise the proof. 
We first note that the fibre of $\Latt(\cB) \to \Tate(\cB)$ is $\Ind(\cB) \cap \Pro(\cB) = \cB^\natural$, and that the image of the induced map $\Latt(\cB)/\cB \to \Tate(\cB)$ generates $\Tate(\cB)$ under finite limits and colimits by the definition of $\Tate(\cB)$. It will hence suffice to show that this map is fully-faithful. 
Writing an $[i \to p] \in \Latt(\cB)$ as the middle object of the exact sequence 
\[
\begin{tikzcd}
i\ar[d,equal]\ar[r,equal] & i \ar[d]\ar[r] & 0 \ar[d] \\
i \ar[r] & p \ar[r] & \cof[i\to p]
\end{tikzcd}
\]
in $\Ar(\Ind\Pro(\cB))$ 
shows that the fibre of
\[
\hom_{\Latt(\cB)}(i' \rightarrow p',i \rightarrow p) \lrar \hom_{\Tate(\cB)}(\cof(i' \rightarrow p'), \cof(i \rightarrow p))
\]
is $\map_{\Ar(\Ind\Pro(\cB))}(i' \to p',i =i) = \map_{\Ind\Pro(\cB)}(p',i)$. 
At the same time, the analogous fibre with $\Latt(\cB)/\cB$ instead of $\Tate(\cB)$ can be computed using the colimit formula for mapping spectra in Verdier projections (Proposition~\reftwo{proposition:verdier-localisation}). Comparing the two expressions 
it will now suffice to show that the induced map 
\[
\colim_{[a \rightarrow i] \in \cB_{/i}}\hom_{\Ind\Pro(\cB)}(p',a) \lrar \hom_{\Ind\Pro(\cB)}(p',i)
\]
is an equivalence. But this is clear, since $p'$ is compact in $\Ind(\Pro(\cB))$. 
 
We now address the second statement, and first show that it holds after applying cores. Write \(j,\hat{j},i,\hat{i}, h\) for the inclusions appearing in the commutative square
\[
\begin{tikzcd}
\cB^\natural\ar[r, "j"]\ar[d, "\hat{j}"']\ar[dr,"h"] & \Ind(\cB) \ar[d,"\hat{i}"] \\
\Pro(\cB) \ar[r,"i"] & \Ind\Pro(\cB).
\end{tikzcd}
\]
Given an exact functor \(g\colon \C \to \cB^\natural\), we show that the homotopy fibre of the map in question over \(g \in \Map_{\Catx}(\C,\cB^\natural)\) is contractible. By Lemma~\reftwo{lemma:adj-verdier}, this homotopy fibre can equivalently be described as the space of refinements of \(g\) to an adjointable commutative square of the form appearing on the top left of the diagram
\[
\begin{tikzcd}
\C\ar[r, "g"]\ar[d, "f"'] & \cB^\natural \ar[d]\ar[r,"h"] & \Ind\Pro(\cB) \ar[d,"{\id_{(-)}}"]\\
\D\ar[r, dashed]\ar[dr,dashed,"{(\phi,\rho)}"'] & \Latt(\cB) \ar[r]\ar[d] & \Ar(\Ind\Pro(\cB)) \ar[d] \\
&\Ind(\cB) \times \Pro(\cB)  \ar[r] & \Ind\Pro(\cB) \times \Ind\Pro(\cB).
\end{tikzcd}
\]
The bottom right square is cartesian by the definition of $\Latt(\cB)$, and hence, temporarily ignoring 
the adjointability constraint, the data of such an extension 
corresponds to the data of a tuple \((\phi,\rho,\tau)\), 
where \(\phi\colon \D \to \Ind(\cB)\) is an extension of \(jg\) to \(\D\), \(\rho\colon \D \to \Pro(\cB)\) is an extension of \(\hat{j}g\) to \(\D\), and \(\eta \colon \hat{i}\phi\Rightarrow i\rho\) is a natural transformation between the resulting composites \(\D \to \Ind\Pro(\B)\), extending the identity transformation \(\hat{i}\phi|_{\C} = hg = i\rho|_{\C}\) on \(\C\). 
As for the adjointability constraint, we note that since the Ind-left and Pro-right adjoints of \(\cB \to \Latt(\cB)\) are given by the source and target projections respectively (see the last part of Remark~\reftwo{remark:concrete}), we have that 
the top left square in the larger diagram above determined by a tuple \((\phi,\rho,\tau)\) 
is vertically Ind-adjointable if and only if \(\phi\) is left Kan extended from \(\C\), and is vertically Pro-adjointable if and only if \(\rho\) is right Kan extended from \(\C\) (see Remark~\reftwo{remark:adj-beck-chevalley}).
In addition, since \(\hat{i}\) preserves  
filtered colimits and \(i\) preserves 
cofiltered limits these conditions imply that \(\hat{i}\phi\) is left Kan extended from \(\C\) and \(i\rho\) is right Kan extended from $\C$. 
The space of natural transformations \(\hat{i}\phi \Rightarrow i\rho\) extending the identity is 
hence contractible for such \(\phi\) and \(\rho\), and  
so by uniqueness of Kan extensions we conclude that the space of triples \((\phi,\rho,\eta)\) satisfying the adjointability constraint  
is contractible, as desired.

To finally obtain the statement at the level of functor categories, insert the Verdier sequences $\C_{\Delta^n} \rightarrow \D_{\Delta^n} \rightarrow \E_{\Delta^n}$ formed by the tensor with $\Delta^n \in \Cat$, confer (the proof of) Proposition \reftwo{proposition:(co)tensor-Verdier}, into the already proven statement and use conservativity of Rezk's nerve.
\end{proof}

\begin{remark}
\label{remark:explicit-classifying}%
The functors appearing in the classifying square of $\C$ can be explicitly described: 
Let \(r\colon \D \to \Ind(\C)\) and \(q\colon \D \to \Pro(\C)\) be the Ind-right and Pro-left adjoints, respectively, of \(f\colon \C \to \D\).
Unwinding the proof of Theorem~\reftwo{theorem:universal-verdier} in this case, 
we see that \(\vphi\) is determined by the unique natural transformation \(r \Rightarrow q\) extending the identity transformation on \(\C\), where both \(r\) and \(q\) are considered as functors \(\D \to \Ind\Pro(\C)\). 
But this transformation is easily described: it is necessarily the preimage under \(f\) of the composite  
of the unit and counit \(fr \to \id \to fq\) in \(\Ind\Pro(\D)\). 
Calling this canonical transformation \(\eta\),  
we then conclude that 
\[
\varphi(d) = [r(d) \xrightarrow{\eta_d} q(d)] \in \Latt(\C) .
\]
Consequently, \(\psi\) must send \(p(d) \in \E\) to \(\cof[\eta_d] \in \Tate(\C)\). 

Writing \(g\colon \E \to \Pro(\D)\) and \(h\colon \E \to \Ind(\D)\) for the Pro-left and Ind-right adjoints of \(p\), respectively, we have in \(\Ind\Pro(\D)\) unit and counit maps \(gp(d) \to d \to hp(d)\) such that \(q(d) = \cof[gp(d) \to d]\) and \(r(d) = \fib[d \to hp(d)]\). For \(e \in \E\) we hence conclude that 
\[
\psi(e) = \cof[g(e) \to h(e)],
\]
an expression which a-priori describes an object of \(\Tate(\D)\), but that by the above actually lies in \(\Tate(\C) \subseteq \Tate(\D)\).
\end{remark}

\begin{example}
For the Verdier projection $\Modf{A}\rightarrow \Modf{A[s^{-1}]}$ for some $s \in \pi_k(A)$, $A$ an $\Einf$-ring spectum, the classifying functor $\varphi \colon \Modf{A} \rightarrow \Latt(\Modf{A}_s)$ takes $M \in \Modf{A}$ to the composite
\[
\colim_{i \in \mathbb N} \mathbb S^{-1} \otimes M/s^i \xrightarrow{\beta} M \xrightarrow{\mathrm{pr}} \lim_{j \in \mathbb N\op} M/s^j
\]
where $\beta$ denotes the Bockstein operation associated with multiplication by the various $s^i$. Identifying $\Tate(\Modf{A}_s)$ as a full subcategory of $\Pro(\Mod(A)_s)$
we thus find that $\psi \colon \Modf{A[s^{-1}]} \rightarrow \Pro(\Mod(A)_s)$ is determined by taking $M[s^{-1}]$ for $M \in\Modf{A}$ to
\[
\cof(\colim_{i \in \mathbb N} \mathbb S^{-1} \otimes M/s^i \xrightarrow{\beta} \lim_{j \in \mathbb N\op} M/s^j) \simeq \lim_{j \in \mathbb N\op} \colim_{i \in \mathbb N} M/s^{i+j} \simeq \lim_{j \in \mathbb N\op} M/s^\infty
\]
where $M/s^\infty = \cof(M \rightarrow M[s^{-1}])$ and the transition maps in the rightmost term are given by multiplication with $s$. Consequently, the composite with $\lim \colon \Pro(\Mod(A)_s) \rightarrow \Mod(A)$ is given by $N \longmapsto  A_s^\wedge \otimes_A N$ on account of the fibre sequence $M_s^\wedge\rightarrow \lim_{j \in \mathbb N\op} M/s^\infty \rightarrow M/s^{\infty}$ in $\Mod(A)$.
\end{example}

Let us now explain how the above structure theory specialises to give a classification of split Verdier sequences.
Namely, suppose that \(\cB\) is an idempotent complete stable \(\infty\)-category and consider the universal Verdier sequence with kernel \(\cB\). Since by Remark~\reftwo{remark:concrete} the Ind-right and Pro-left adjoints of \(\cB \to \Latt(\cB)\) are given exactly by \((i \rightarrow p) \mapsto i\) and \((i \rightarrow p) \mapsto p\), respectively, we see that the split core of the universal Verdier sequence is given by the split Verdier sequence
\[
\cB \xrightarrow{b \mapsto \id_b} \Ar(\cB) \xrightarrow{\cof} \cB,
\]
of Example~\reftwo{example:universal-split-verdier}.
This operation gives us a universal split Verdier sequence with kernel \(\cB\). In fact, assuming \(\cB\) is idempotent complete is not necessary:

\begin{corollary}
\label{corollary:universal-split-verdier}%
Let \(\cB\) be a stable \(\infty\)-category. Then, 
for every split Verdier sequence
\(\C \to \D \to \E\)
the map
\[
\Fun^{\Ver}(\D \to \E,\Ar(\cB) \to \cB) \longrightarrow \Funx(\C,\cB)
\]
extracting kernels is an equivalence. Consequently, there exists an essentially unique adjointable square
\[
\begin{tikzcd}
\D \ar[r,"\vphi"]\ar[d, "p"'] & \Ar(\C) \ar[d, "\cof"] \\
\E \ar[r,"\psi"] & \C
\end{tikzcd}
\]
inducing the identity on vertical fibres. It is cartesian, so that the split Verdier projection \(\D\to \E\) is pulled back from the universal split Verdier projection \(\cof \colon \Ar(\C) \to \C\).
\end{corollary}

\begin{proof}
In the case where \(\cB\) is idempotent complete this follows from Theorem~\reftwo{theorem:universal-verdier} via Lemma~\reftwo{lemma:split-core}. For arbitrary \(\cB\), it suffices by Lemma~\reftwo{lemma:adj-verdier} to note that in any adjointable square
\[
\begin{tikzcd}
\C \ar[r, "\rho"]\ar[d] & \cB^{\natural} \ar[d,"\id_{(-)}"] \\
\D \ar[r, "\vphi"] & \Ar(\cB^{\natural}) \ ,
\end{tikzcd}
\]
the condition that \(\rho\) factors through \(\cB \subseteq \cB^{\natural}\) is equivalent to the condition that \(\vphi\) factors through \(\Ar(\cB)\). 
\end{proof}

We refer to the split Verdier sequence \(\cB \xrightarrow{b \mapsto \id_b} \Ar(\cB) \xrightarrow{\cof} \cB\) of Corollary~\reftwo{corollary:universal-split-verdier} as the universal split Verdier sequence, and to its inclusion and projection as the universal split Verdier inclusion and projection, respectively. We obtain the following well-known consequence:

\begin{corollary}
\label{corollary:classificationsplitverdierdumdidum}%
Given two stable $\infty$-categories $\C$ and $\E$ extracting classifying functors induces an equivalence between the space of Verdier sequences $\C \rightarrow \D \rightarrow \E$ and $\Hom_{\Catx}(\E,\C)$.
\end{corollary}

\begin{remark}
\label{remark:explicit-classifyingsplit}%
\label{corollary:classification-split-verdier}%
Again the classifying functor can be described explicitly: Given a split Verdier sequence
\[
\begin{tikzcd}
[column sep=7ex]
\C \ar[r, "f" description] & 
\D 
\ar[r,"p" description] 
\ar[l,bend left=25,shift left=1.5ex,start anchor=west,end anchor=east,"r"] 
\ar[l,bend right=25,shift right=1.5ex,start anchor=west,end anchor=east,"q"'] 
\ar[l,phantom,shift left=1.2ex,start anchor=west,end anchor=east,"\myperp"] 
\ar[l,phantom,shift right=1.3ex,start anchor=west,end anchor=east,"\myperp"]
& \E,
\ar[l,bend left=25,shift left=1.5ex,start anchor=west,end anchor=east,"h"] 
\ar[l,bend right=25,shift right=1.5ex,start anchor=west,end anchor=east,"g"']
\ar[l,phantom,shift left=1.2ex,start anchor=west,end anchor=east,"\myperp"] 
\ar[l,phantom,shift right=1.3ex,start anchor=west,end anchor=east,"\myperp"]
\end{tikzcd}
\]
the functor \(\vphi\colon \D \to \Ar(\C)\) sends \(d \in \D\) to the canonical unit-counit composite map \(\eta_d\colon r(d) \to q(d)\)
as in 
Remark~\reftwo{remark:explicit-classifying}. Since \(rh = qg = 0\) we consequently have that 
the functor \(\psi\colon \E \to \C\) is given by 
\[
\psi(e) = \cof(\eta_{h(e)}) = qh(e),
\]
and a similar calculation using the other adjoints gives $\psi(e) = \cof(\eta_{g(e)}) = \Sig rg(e)$. These formulas then yield two cartesian squares
\[
\begin{tikzcd}
d \ar[r] \ar[d] & fq(d) \ar[d] &&   \Omega f \psi p (d) \ar[r] \ar[d] & fr(d) \ar[d] \\ 
hp(d) \ar[r] & f\psi p(d) &&  gp(d) \ar[r] & d
\end{tikzcd}
\]
for every $d \in \D$.
\end{remark}

\begin{remark}
\label{remark:references}%
Applying the left/right split cores to the universal Verdier sequence, one obtains similar classifications of left/right split Verdier sequences, respectively. The universal right and left split Verdier sequences with fibre $\C$ (not necessarily idempotent complete) are given by
\[
\C \longrightarrow \Latt(\C) \times_{\Ind(\C)} \C \longrightarrow \Pro(\C) \quad \text{and} \quad \C \longrightarrow \Latt(\C) \times_{\Pro(\C)} \C \longrightarrow \Ind(\C)
\]
and the spaces of right and left split Verdier sequences with fibre $\C$ and base $\E$ are equivalent to the spaces of functors $\E \rightarrow \Pro(\C)$ and $\E \rightarrow \Ind(\C)$, respectively. 

Interpreting a functor $\E \rightarrow \Ind(\C)$ instead as a functor $\C\op \times \E \rightarrow \Sps$ the inverse process to the formation of classifying functors for left split Verdier sequences is furthermore simply given by forming the category of pairings as in \refone{subsection:bifibrations}, and similarly for the right split case by taking appropriate opposites.
\end{remark}

Finally, let us explain how to deduce a classification result for Karoubi sequences, essentially a direct consequence of Theorem \reftwo{theorem:universal-verdier} and the universal property of idempotent completions. In what follows, $\textsc{Kar}$ denotes the $\infty$-category analogous to $\textsc{Ver}$: It has objects the Karoubi sequences among locally small stable $\infty$-categories with small kernel and morphisms those transformations inducing adjointable squares on the level of Karoubi inclusions/projections.

\begin{corollary}
\label{corollary:universal-karoubi}%
For $\cB$ a stable $\infty$-category, the restriction map
\[
\Fun^{\textsc{Kar}}(\C \to \D \to \E, \cB \to \Latt(\cB) \to \Tate(\cB)^\natural) \longrightarrow \Funx(\C,\cB)
\]
is an equivalence for every locally small Karoubi sequence \(\C\to \D\to \E\), and there is consequently an essentially unique adjointable square
\[
\begin{tikzcd}
\D \ar[r,"\vphi"]\ar[d, "p"'] & \Latt(\C) \ar[d] \\
\E \ar[r,"\psi"] & \Tate(\C)^\natural \ ,
\end{tikzcd}
\]
inducing the inclusion $\C \rightarrow \C^\natural$ on vertical fibres. This square is cartesian if and only if $\C$ is idempotent complete, so that in this case the Karoubi projection \(p\) is pulled back from the universal Karoubi projection with fibre $\C^\natural$ on the right. 
\end{corollary}

\begin{corollary}
Given stable $\infty$-categories $\C$ and $\E$ extracting classifying functors and boundary maps on $\K$-spectra induces an equivalence between the space of Karoubi sequences $\C \rightarrow \D \rightarrow \E$ and the space of quadruples $(\psi, \E',\C',\eta)$, where $\psi \colon \E \rightarrow \Tate(\C)^\natural$ is an exact functor, $\E^\mathrm{min} \subseteq \E' \subseteq \E$ and $\C \subseteq \C' \subseteq \C^\natural$ are stable subcategories and $\eta \colon \K(\E') \rightarrow \mathbb S^1 \otimes \K(\C')$ is a lift of the induced map $\K(\E') \rightarrow \K(\E) \rightarrow \K(\Tate(\C)^\natural)$ along
\[
\mathbb S^1 \otimes \K(\C') \longrightarrow\mathbb S^1 \otimes \K(\C^\natural) \xrightarrow{\partial} \K(\Tate(\C)) \longrightarrow \K(\Tate(\C)^\natural).
\]
Equivalently, it is the space of quadruples $(\psi,e,c,s)$, where $\psi$ is as before, 
\[
e \subseteq \mathrm{ker}(\psi_\ast \colon \K_0(\E) \longrightarrow \K_{0}(\Tate(\C)^\natural)) \quad \text{and} \quad \K_0(\C) \subseteq c \subseteq \K_0(\C^\natural)
\]
are subgroups, such that
in the exact sequence
\[
\K_0(\C^\natural)/c \longrightarrow \K_0\big(\Latt(\C) \times_{\Tate(\C)} \E^e\big)/c \longrightarrow e \longrightarrow 0
\]
the first map is injective, and $s$ is a splitting of it as abelian groups, where $\E^e = \{X \in \E \mid [X] \in e\}$.
\end{corollary}

The subgroups $e$ and $c$ simply record the essential image of $\D \rightarrow \E$ and the kernel of $\D \rightarrow \E$. Armed with this information the statement follows easily from Corollary \reftwo{corollary:classverdierfull}.

\begin{proof}[Proof of Corollary~\reftwo{corollary:universal-karoubi}]
By Proposition~\reftwo{proposition:criterion-karoubi}
and Corollary~\reftwo{corollary:characterisation-of-Karoubi-inclusions-and-projections}, both $\C \xrightarrow{f} \D \xrightarrow{p} \E$ and $\ker(p) \rightarrow \D \rightarrow \mathrm{im}(p)$ admit canonical pointwise dense inclusions into $\ker(p) \rightarrow \D \rightarrow \E$,
which then induce equivalences on mapping spaces into $\cB^\natural \rightarrow \Latt(\cB) \rightarrow \Tate(\cB)^\natural$, since the latter consists of idempotent complete $\infty$-categories. Similarly, the inclusion $\C \subseteq \ker(p)$ induces an equivalence on mapping spaces into $\cB^{\natural}$.
But $\ker(p) \rightarrow \D \rightarrow \mathrm{im}(p)$ is a Verdier sequence by Corollary~\reftwo{corollary:characterisation-of-Karoubi-inclusions-and-projections},
and so combining the above  
with Theorem \reftwo{theorem:universal-verdier} we conclude that the map
\[
\Fun^{\textsc{Kar}}([\C \to \D \to \E], [\cB^{\natural} \to \Latt(\cB) \to \Tate(\cB)^\natural]) \lrar \Funx(\C,\cB^{\natural}) ,
\]
induced by extracting kernels, is indeed an equivalence.
\end{proof}

\section{Comparisons to previous work}
\label{appendix:AppIIB}%

For a ring $R$ equipped with an invertible module with involution $M$ and a form parameter $\lambda$ on $(R,M)$, the main result of~\cite{comparison} 
asserts that the canonical map
\[
\Unimod(R,\lambda)^\grp \rightarrow \GWspace(\Dperf(R),\QF^{\g\lambda}_M)
\]
is an equivalence; here, $\Unimod(R,\lambda)$ denotes the groupoid of unimodular, 
$\lambda$-hermitian forms on finitely generated projective $R$-modules, symmetric monoidal under orthogonal sum, and $(-)^{\grp}$ denotes group completion in the $\infty$-category of $\Einf$-monoids in spaces. 
For the symmetric and quadratic form parameters, the left hand side is essentially the classical definition of the Grothendieck-Witt spaces of Karoubi and Villamayor.
When $2$ is invertible in $R$, Schlichting~\cite{schlichting-derived} refined these Grothendieck-Witt-spaces to spectra, whose negative homotopy groups are Balmer's Witt groups. Keeping the assumption that $2$ is invertible, 
these coincide with the classical symmetric $\L$-groups of Ranicki, 
and hence with the $\L$-groups of $\Dperf(R)$ with respect to the symmetric Poincaré structures associated to $M$. 
Unwinding the definition of the relevant maps, 
the above comparison of Grothendieck-Witt groups thus extends to negative groups for the quadratic/symmetric form parameter when $2$ is invertible. 

In this appendix we extend such comparisons beyond the case of rings and to the spectrum level, 
and equate our construction with previous work in two new cases: 
First, in \S\reftwo{subsection:HSV-compare}
we compare our Grothendieck-Witt space with the 
Grothendieck-Witt space of Heine, Spitzweck and Verdugo in \cite{Spitzweckreal}, which is defined 
in the setting of Waldhausen \(\infty\)-categories with genuine duality. This in particular also implies an equivalence with Spitzweck's definition of Grothendieck-Witt spaces for stable $\infty$-categories with duality \cite{Spitzweck-GW} in the case where the Poincaré structure is the symmetric one associated to a given duality. 
Then, in \S\reftwo{subsection:schlichting} we show that Hornbostel and Schlichting's Grothendieck-Witt spectrum of a exact category with duality on which $2$ act invertibly \cite{hornbostel,schlichting-exact} can be recovered as the symmetric Grothendieck-Witt spectrum of its derived $\infty$-category (also known as its stable envelope, or stable hull, see~\cite{klemenc}).

\begin{Rmk}
\begin{enumerate}
\item The comparison results of \S\reftwo{subsection:schlichting} can sometimes be extended to situations in which $2$ is not invertible using the associated genuine symmetric Poincaré structure on the derived \(\infty\)-category. The case of the exact category of vector bundles on a sufficiently nice scheme is treated in \cite{motives}*{Section 4.6}, for example. \footnote{During the revision process of this paper, these results were improved quite drastically in \cite{Schlichting-genuine}.}
\item An equivalence of the real $\K$-spaces (and thus in particular the Grothendieck-Witt spaces) defined in \cite{HeineLopez-AvilaSpitzweck} with our genuine symmetric ones is given in \cite{comparison}*{Section 8.4}, and Grothendieck-Witt spaces of split-exact categories with form parameters from \cite{SchlichtinghigherI} are treated in \cite{comparison}*{Section 8.2}.
\item We do not attempt here a comparison of our work to either the Grothendieck-Witt space of \cite{SchlichtinghigherI} for non-split exact categories or the entire real algebraic $K$-space defined in \cite{Spitzweckreal}. \footnote{During the revision process of this paper, a comparison in the former case was provided in \cite{Schlichting-Marlowe}, and in the latter in \cite{realSvshermQ}.} %
\end{enumerate}
\end{Rmk}

\subsection{The Heine-Spitzweck-Verdugo Grothendieck-Witt space} 
\label{subsection:HSV-compare}%
\label{subsection:Spitzweck-compare}%

In \cite{HeineLopez-AvilaSpitzweck, Spitzweckreal} the authors study the Grothendieck-Witt theory of what they call Waldhausen \(\infty\)-categories with genuine duality. As they show in~\cite{Spitzweckreal}*{\S 6}, this is a direct generalisation of the notion of a Poincaré \(\infty\)-category. In particular, the \(\infty\)-category \(\Waldgd\) of Waldhausen \(\infty\)-categories with genuine duality contains \(\Catp\) as a full subcategory spanned by those objects whose underlying \(\infty\)-category is stable and whose Waldhausen structure is maximal. 
In loc.\ cit., the authors construct in particular a Grothendieck-Witt functor, which they denote by $\mathrm{KH}$, and which we will denote by $\GWspace^{\gd}$ for the sake of notational consistency. Our goal in this subsection is to show that both definitions of the Grothendieck-Witt space are compatible for Poincaré $\infty$-categories.

As the Waldhausen structure shall play no role in this appendix, we suppress it from further discussion, and recall that the \(\infty\)-category \(\Catgd\) of \(\infty\)-categories with genuine duality is defined via the cartesian square
\[
\begin{tikzcd}
\Catgd \ar[r]\ar[d] & \Cat^{\hC} \ar[d, "\catforms^{\sym}"] \\
\RFib \ar[r,"t"]& \Cat,
\end{tikzcd}
\]
where \(\RFib \subseteq \Ar(\Cat)\) is the full subcategory spanned by the right fibrations and \(\catforms^{\sym}\) 
is the functor sending an \(\infty\)-category with duality \((\C,\Dual)\) to the \(\infty\)-category \(\catforms^{\sym}(\C,\Dual)\coloneq  \Twar(\C)^{\hC}\) of symmetric hermitian objects in \((\C,\Dual)\). 

Poincar\'e $\infty$-categories are then regarded as such objects by means of the functor 
\[
\Catp \longrightarrow \Catgd, \quad (\C,\QF) \longmapsto (\C,\Dual_{\QF},\catforms(\C,\QF) \to \catforms(\C,\QF^\sym) = \catforms^{\sym}(\C,\Dual_{\QF})).
\]
In \cite{Spitzweckreal}*{\S 6} the authors show that the essential image consists of exactly those $\infty$-categories with genuine duality whose underlying $\infty$-category is stable and for two Poincaré \(\infty\)-categories \((\C,\QF),(\D,\QFD)\) the map
\[
\Hom_{\Catp}((\C,\QF),(\D,\QFD)) \to \Hom_{\Catgd}((\C,\Dual_{\QF},\catforms(\C,\QF),(\D,\Dual_{\QFD},\catforms(\D,\QFD)))
\]
is a full inclusion of spaces whose essential image contains all duality preserving functors whose underlying functor is exact. 

The construction of \(\GWspace^{\gd}\) uses in an essential that \(\Catgd\) is cartesian closed \cite{Spitzweckreal}*{Corollary 2.43}. An explicit model for the internal mapping object %
from $(\I,\Dual,\phi)$ to $(\J,\Dual',\phi'))$ is given by the functor \(\infty\)-category \(\Fun(\I,\J)\), equipped with the conjugation duality $f \mapsto \Dual' f\op \Dual\op\), and with its genuine refinement
classified by 
\[
\catforms^{\sym}(\Fun(\I,\J), \Dual' \circ - \circ  \Dual\op)\op \longrightarrow \Sps, \quad f \longmapsto \Hom_{/\catforms^{\sym}(\I,\Dual)}(\phi,\catforms^{\sym}(f)^*\phi').
\]

Since the functor \(\Fun(\I,-)\) preserves stable \(\infty\)-categories and exact functors this operation restricts to a cotensoring operation \((\Catgd)\op \times \Catp \to \Catp\), naturally in \(\I\). Let us denote its output by \((\C,\QF)^{(\I,\Dual,\phi)}\), which explicitly is the Poincaré \(\infty\)-category whose underlying stable \(\infty\)-category is \(\C^\I = \Fun(\I,\C)\), whose duality is the conjugation duality \(\mathrm{Cf}_\Dual^{\Dual_\QF}\), and whose Poincaré structure \(\QF^{(\I,\Dual,\phi)}\colon \Fun(\I,\C)\op \to \Sp\) fits into the cartesian square
\begin{equation}
\label{equation:poincare-cotensor}%
\begin{tikzcd}
\QF^{(\I,\Dual,\phi)}(f) \ar[r]\ar[d] & \lim_{\mathcal H}\pi^*f^*\QF \ar[d]\\
\mathrm{nat}(f,\Dual_{\QF} f\op\Dual\op)^{\hC} \ar[r] & \lim_{\mathcal H}\pi^*f^*\Bil_{\QF}^{\hC},
\end{tikzcd}
\end{equation}
where \(\pi\) is the composite of \(\phi \colon \mathcal H \to \catforms^{\sym}(\I,\Dual)$ with the projection to $\I$ and the bottom horizontal map is the composite
\[
\mathrm{nat}(f,\Dual_{\QF}\op f\op\Dual)^{\hC} = \big[\lim_{i \to j}\Bil_{\QF}(f(i),f(\Dual j))\big]^{\hC} 
\to \big[\lim_{i \in \mathcal H}\Bil_{\QF}(f(\pi(i)),f(\pi(i)))\big]^{\hC} = \lim_{\mathcal H}\pi^*f^*\Bil_{\QF}^{\hC},
\]
with the left limit ranging over $\Twar(\I)$.

This cotensor operation is used to define $\GWspace^\gd(\C,\QF)$ for a Poincaré \(\infty\)-category $(\C,\QF)$ as follows. 
One considers the poset \([n] \ast [n]\op\) equipped with (unique) duality \(\Dual_{\fl}\) flipping the two components, and consider \(\Ar([n] \ast [n]\op) = \Fun([0) \ast [0]\op, [n] \ast [n]\op)\) as equipped with the associated conjugation duality. Viewing \(\Ar([n]\ast [n]\op)\) as a category with genuine duality via genuine refinement the identity of \(\catforms^{\sym}(\Ar([n]\ast [n]\op))\) we can take the cotensor Poincaré \(\infty\)-category \((\C,\QF)^{\Ar([n] \ast [n]\op)}\) in the above sense, and define the (edgewise subdivided) $\rS$-construction 
\(\rS^e_n(\C,\QF)\) be the full subcategory of \((\C,\QF)^{\Ar([n] \ast [n]\op)}\) spanned by those diagrams \(\vphi\colon \Ar(\Del^n \ast (\Del^n)) \to \C\) such that \(\varphi(i \leq i) \simeq 0\) for every \(i \in [n] \ast [n]\op\) and such that the square  
\[
\begin{tikzcd}
\varphi(i \leq k) \ar[r] \ar[d] & \varphi(i \leq l) \ar[d] \\
\varphi(j \leq k) \ar[r] & \varphi(j \leq l)
\end{tikzcd}
\]
is exact for every \(i \leq j \leq k \leq l\) in \([n] \ast [n]\op\), equipped with the restricted hermitian structure, %
which is Poincaré since this full subcategory is duality invariant. We note that the simplicial object in stable \(\infty\)-categories underlying \(\rS^e(\C,\QF)\) is the edgewise subdivision of the usual \(\rS\)-construction \(\rS(\C)\) used to define the algebraic \(\K\)-theory of \(\C\), see the discussion after Propositon \reftwo{proposition:pi0span} for a reminder. Finally, by \cite{Spitzweckreal}*{Proposition 10.2} the $\GWspace^\gd(\C,\QF)$ agrees with the fibre of $|\core\rS^{e}(\C,\QF)^\Ct| \to |\core\rS^{e}(\C)|$, where $(-)^\Ct \colon \Cat^\gd \rightarrow \Sps$ is the genuine fixed point functor. Restricted to Poincar\'e $\infty$-categories it agrees (by construction) with $\Poinc$, so we finally arrive at
\[
\GWspace^{\gd}(\C,\QF) \simeq  \fib[|\Poinc\rS^{e}(\C,\QF)| \longrightarrow |\core\rS^{e}(\C,\QF)|].
\]

\begin{remark}
\label{remark:recovers-spitzweck}%
When $(\C,\QF)$ is the symmetric Poincaré category associated to a given stable $\infty$-category with duality $(\C,\Dual)$ the vertical maps in the square defining $\QF^{(\I,\Dual,\phi)}$ are equivalences, and so the cotensor construction of \cite{Spitzweckreal} reduces to the usual cotensor construction $(\C,\Dual)^{(\I,\Dual')} = (\Fun(\I,\C),\Dual' \circ - \circ \Dual\op)$
of $\infty$-categories with duality. Consequently, the edgewise subdivided $\rS$-construction $\rS^{e}(\C,\QF^\sym_\Dual)$ coincides with Spitzweck's edgewise subdivided $\rS$-construction for stable $\infty$-categories with duality. In particular, for these types of Poincaré $\infty$-categories the Grothendieck-Witt spectrum of~\cite{Spitzweckreal} recovers that of~\cite{Spitzweck-GW}, see~\cite{Spitzweckreal}*{Corollary 10.4}.
\end{remark}

\begin{proposition}
\label{proposition:S-is-Q}%
For a Poincaré $\infty$-category $(\C,\QF)$ there is a natural equivalence of simplicial Poincar\'e $\infty$-categories \(\rS^e(\C,\QF) \simeq \Q(\C,\QF)\).
\end{proposition}

\begin{corollary}
For a Poincaré $\infty$-category $(\C,\QF)$ there is a natural equivalence
$\GWspace^{\gd}(\C,\QF) \simeq \GWspace(\C,\QF)$.
\end{corollary}

\begin{proof}
Combine Proposition~\reftwo{proposition:S-is-Q} with Corollary~\reftwo{corollary:classicaldefofgw}.
\end{proof}

When $(\C,\QF)$ is the symmetric Poincaré $\infty$-category associated to a given stable $\infty$-category with duality $(\C,\Dual)$, the above comparison reduces to a comparison between $\GWspace(\C,\QF)$ and Spitzweck's Grothendieck-Witt space of an $\infty$-category with duality~\cite{Spitzweck-GW}, see Remark~\reftwo{remark:recovers-spitzweck}.

\begin{proof}[Proof of Proposition~\reftwo{proposition:S-is-Q}]
We construct an equivalence between the \(\Q\)- and \(\rS\)-construction in two steps. Let \(\I_n \subseteq \Ar([n] \ast [n]\op)\) be the duality invariant subcategory spanned by those arrows of \([n] \ast [n]\op\) whose source is \([n]\) and whose target is in \([n]\op\). Then \(\I_n = [n] \times [n]\op\) and we write objects in \(\I_n\) as pairs \((i,j)\) with \(i \in [n],j \in [n]\op\).
We then see that \(\I_n\) is closed under duality and the induced duality on \(\I_n\) is given by \((i,j) \mapsto (j,i)\). 
We can view \(\I_n\) as an \(\infty\)-category with 
genuine duality using the identity of $\catforms^{\sym}(\I_n)$, and view the inclusion 
\(\I_n \hrar \Ar([n] \ast [n]\op)\) as a morphism of \(\infty\)-categories with genuine duality. Let \(\rT_n(\C,\QF) \subseteq (\C,\QF)^{\I_n}\) denote the full subcategory spanned by those 
diagrams \(\vphi\colon \I_n\to \C\) for which the square 
\[
\begin{tikzcd}
\varphi(i, k) \ar[r] \ar[d] & \varphi(i, l) \ar[d] \\
\varphi(j, k) \ar[r] & \varphi(i, l)
\end{tikzcd}
\]
is exact for every \(i \leq j \in [n]\) and \(l \geq k \in [n]\op\), with the restricted Poincar\'e structure and, adhering to our conventions, denote by \(\rT_n(\C)\) the underlying stable \(\infty\)-category of \(\rT_n(\C,\QF)\). The induced functor on cotensors (as $\infty$-categories with genuine duality) \((\C,\QF)^{\Ar([n] \ast [n]\op)} \to (\C,\QF)^{\I_n}\) maps \(\rS^e_n(\C)\) to \(\rT_n(\C)\), and since
the full subcategory \(\I_n \subseteq \Ar([n]\ast[n]\op)\) is clearly stable under all simplicial transition functors we obtain a map of simplicial  Poincaré \(\infty\)-categories $\rS^e_n(\C,\QF) \to \rT_n(\C,\QF)$.

Next, we observe that 
\[
\catforms^{\sym}(\I_n) = \Twar([n] \times [n]\op)^{\hC} = (\Twar([n]) \times \Twar([n]\op))^{\hC} = \Twar([n])
\]
with the equivalence realised by sending \((i \leq j) \in \Twar([n])\) to \((i \leq j,j \geq i) \in \Twar([n]) \times \Twar([n]\op)\).
Now the cotensor construction comes with a canonical evaluation map
\[
\I_n \to \Funx((\C,\QF)^{\I_n},(\C,\QF)) ,
\]
in \(\Catgd\). Taking symmetric hermitian objects we then obtain a functor
\[
\Twar[n] = \catforms^{\sym}(\I_n) \to \catforms^{\sym}(\Funx((\C,\QF)^{\I_n},(\C,\QF))) = \Funh((\C,\QF)^{\I_n},(\C,\QF)) ,
\]
which induces an hermitian functor
$q_n\colon (\C,\QF)^{\I_n} \to (\C,\QF)^{\Twar[n]}$.
To avoid confusion, let us emphasise that the source and target of this map are cotensors of two different types: the former is the cotensor of Poincaré $\infty$-categories over $\infty$-categories with genuine duality constructed in~\cite{Spitzweckreal}, while the target is the cotensor of hermitian $\infty$-categories over $\infty$-categories constructed in~\S 6 of \paperone; the map $q_n$ is constructed using the universal property of its target, and is hence a-priori only an hermitian functor.

Comparing the definitions of \(\rT^e_n\) and \(\Q_n\), we see that $q_n$ 
sends the full subcategory \(\rT_{n}(\C)\) to \(\Q_n(\C)\), and hence restricts to an hermitian functor 
\[
q_n' \colon \rT_n(\C,\QF) \longrightarrow \Q_n(\C,\QF)
\]
We now claim that both functors in the sequence
\begin{equation}
\label{equation:S-T-Q}%
\rS^e_n(\C,\QF) \to \rT_n(\C,\QF) \to \Q_n(\C,\QF) 
\end{equation}
are equivalences in $\Catp´h$ (and thus $\Catp$). Let us first consider the situation on the level of underlying stable \(\infty\)-categories. 
To see that the second map is also an equivalence on the level of stable \(\infty\)-categories note that by construction it is given by restriction along the composite \(\Twar([n]) = \catforms^{\sym}([n] \times [n]\op) \to [n] \times [n]\op\), which is just the full embedding of posets \((i \leq j) \mapsto (i,j)\). 
It is then straightforward to verify that left Kan extension along this full subposet inclusion sends \(\Q_n(\C)\) to \(\rT_n(\C)\) and constitutes an inverse to the restriction functor \(\rT_n(\C) \to \Q_n(\C)\). 
For the map $\rS^e_n(\C) \to \rT_n(\C)$ one can argue in a similar manner, this time breaking the argument into several steps which alternate between left and right Kan extension. We also point out that the composite functor agrees with that used in Proposition \reftwo{proposition:pi0span}, in turn the stable analogue of the classical identification of the $\Q$-construction and the edgewise subdivided $\rS$-construction established by Waldhausen \cite{waldhausen}*{\S 1.9}. 

Now for the first functor in~\eqreftwo{equation:S-T-Q}, we note that since it is Poincaré by construction and an equivalence on underlying \(\infty\)-categories, it will suffice to show that it induces an equivalence on colinear parts. By the formula in~\eqreftwo{equation:poincare-cotensor}, this amounts to showing that for every \(\vphi \in \rS^e_n(\C)\) the map
\[
\lim_{\stackrel{(x \leq \Dual x) \in}{\catforms^{\sym}(\Ar([n] \ast [n]\op))\op}} \Lam^{\vee}_{\QF}(\vphi(x)) \to 
\lim_{\stackrel{(x \leq \Dual x) \in}{\catforms^{\sym}(\I_n)\op}} \Lam^{\vee}_{\QF}(\vphi(x))
\]
is an equivalence, where \(\Lam^{\vee}_{\QF}(-) = \fib[\QF(-) \to \Bil^{\hC}(-,-)]\) is the colinear part of \(\QF\). We prove this by showing that the induced map 
\begin{equation}
\label{equation:cofinal-map-twar}%
\Twar([n]) = \catforms^{\sym}(\I_n) \to \catforms^{\sym}(\Ar([n]\ast[n]\op)) 
\end{equation}
is cofinal. 
For this, note that for $\infty$-categories with duality whose underlying categories are posets, the category of symmetric hermitian objects embeds as the full subposet spanned by the elements which are smaller than their dual. Applying this principle to the pair of posets $\I_n \subseteq \Ar([n] \ast [n]\op)$, we may view both sides of~\eqreftwo{equation:cofinal-map-twar} as subposets of \(\Ar([n] \ast [n]\op)\), the source corresponding to the full subposet  
spanned by the arrows of the form \(i_0 \leq j_1\) for \(i\leq j \in [n]\) (where \(i_0\) is considered in \([n]\) and \(j_1\) in \([n]\op\)), 
and the target containing all these objects 
and in addition the objects of the form \(i_0 \leq j_0\) with \(i,j \in [n]\). We then see that the comma categories of the inclusion~\eqreftwo{equation:cofinal-map-twar} 
all have initial objects and so this map is cofinal, as desired.

We now show that the second map in~\eqreftwo{equation:S-T-Q} is an equivalence on hermitian structures. Observing the formula in~\eqreftwo{equation:poincare-cotensor} and using the isomorphism $\Twar(\I_n) = \Twar([n] \times [n]\op) = \Twar([n]) \times \Twar([n]\op)$, this amounts to showing that for every \(\vphi \in \rT_n(\C)\) the map
\[
\Nat(\vphi,\Dual_{\conj}(\vphi))^{\hC} = \left[\lim_{\stackrel{[(i \leq l,j \geq k)] \in}{(\Twar([n]) \times \Twar([n]\op))\op}}\Bil_{\QF}(\vphi(i,j),\vphi(k,l))\right]^{\hC} \to \lim_{(i \leq j)\in \Twar([n])\op}\Bil_{\QF}^{\hC}(\vphi(i,j),\vphi(i,j))
\]
induced by restriction along the inclusion 
\[
\Twar([n]) \to \Twar([n]) \times \Twar([n]\op) \quad\quad (i \leq j) \mapsto (i \leq j,j \geq i) ,
\]
is an equivalence.
Consider the intermediate poset 
\[
\Twar([n])\subseteq \J_n \subseteq \Twar([n]) \times \Twar([n]\op)
\]
spanned by those objects \((i \leq l,j \geq k)\) 
such that \(i \leq j\) and \(k \leq l\). 
The exactness conditions on \(\vphi \in \rT_n(\C)\) imply that \((i \leq l,j \geq k) \mapsto \Bil_{\QF}(\vphi(i, j),\vphi(k , l))\) is right Kan extended from its restriction to \(\J_n\).
To finish the proof we now claim that the inclusion \(\Twar([n]) \to \J_n\) is cofinal.
Indeed, note that for a given \((i \leq l, j \geq k) \in \J_n\), we have 
\(\Twar([n])_{(i \leq l,j \geq k)/} = \Twar([i,l] \cap [k,j])\), 
and \([i,l] \cap [k,j]\) is non-empty by the definition of \(\J_n\).
\end{proof}

\subsection{Schlichting's Grothendieck-Witt spectrum of an exact category with duality}
\label{subsection:schlichting}%

We now turn to the 
comparison with the classical set-up of exact categories with weak equivalences and duality from ~\cite{schlichting-exact,schlichting-mv,hornbostel}. It consists of an additive (ordinary) category \(\E\), equipped with two special types of arrows, namely, inflations and deflations, satisfying suitable properties, as well as a duality \(\Dual\colon \E \to \E\op\), which switches between the inflations and deflations. A symmetric object in \(\E\) is then an object \(X \in \E\) equipped with a self-dual map \(\phi\colon X \to \Dual X\). Such a symmetric object is said to be Poincaré \(\phi\) is an isomorphism. More generally, Schlichting also considers the case where $\E$ is endowed in addition with a collection $W$ of morphisms called weak equivalences, and satisfying suitable axioms. In the presence of such a $W$ the notion of a Poincaré object $(X,\phi)$ is adapted so that the self-dual map \(\phi\colon X \to \Dual X\) is only assumed to be a weak equivalence. One then defines $\cor(\E,W)$ to be the category whose objects are the objects of $\E$ and whose morphisms are the weak equivalences, and $\Poi(\E,\Dual,W)$ to be the category whose objects are the Poincaré objects in $\E$ and whose morphisms are the form preserving weak equivalences.
The Grothendieck-Witt space of an exact category with duality and weak equivalences is then given by
\[
\GWspace(\E,\D, W) := \fib[|\Poi(\rS^e(\E,\Dual,W))| \to |\cor(\rS^e(\E,W))|],
\]
where $\rS^e(\E,\Dual,W)$ is the simplicial object in exact categories with duality and weak equivalences obtained by taking the edgewise subdivided $\rS$-construction of $(\E,\Dual,W)$, and $|-|$ means that we take nerves levelwise and then geometrically realise the resulting bisimplicial set (equivalently, we geometrically realise the categories levelwise and then take homotopy colimit over $\Del\op$). 
In the case when $W$ consists just of isomorphisms one speaks of the Grothendieck-Witt space of an exact category with duality.

In~\cite{schlichting-mv}*{Proposition 6}, Schlichting proves a Quillen type resolution theorem, showing that the Grothendieck-Witt space of an exact category with duality $(\E,\Dual)$ 
is naturally equivalent to the Grothendieck-Witt space of the category $\Ch(\E)$ of bounded complexes in $\E$, considered as an exact category with duality and weak equivalences. Here, the inflations and deflations on $\Ch(\E)$ are defined pointwise and the weak equivalences are the quasi-isomorphisms, that is, the maps which induce isomorphisms on homology groups, where these are computed in the abelian hull $\A(\E)$ of $\E$. The duality $\Dualb$ on $\Ch(\E)$ is then induces from the duality on $\E$.

The $\infty$-categorical localisation $\Der^{\mathrm{b}}(\E) := \Ch(\E)[\qiso^{-1}]$ is stable; it is the bounded derived $\infty$-category of $\E$. To avoid confusion, let us point out that the notion of a quasi-isomorphism in $\Ch(\E)$, and hence also the formation of $\Der^{\mathrm{b}}(\E)$, depend on the exact structure of $\E$ (as they depend on the abelian hull). 
Since the duality preserves the weak equivalences it induces a duality $\ovl{\Dual}$ after localisation, and we may consider the associated Poincaré $\infty$-category $(\Der^{\mathrm{b}}(\E),\QF^{\sym}_{\ovl{\Dual}})$. This leads to the question of comparing Schlichting's Grothendieck-Witt space of $(\E,\Dual)$, or, equivalently, that of $(\Ch(\E),\Dualb,\qiso)$, with the Grothendieck-Witt space of $(\Der^{\mathrm{b}}(\E),\QF^{\sym}_{\ovl{\Dual}})$ defined in the present paper.

Let us say that an exact category with duality and weak equivalences is \defi{homotopically sound} if the collection of deflations and weak equivalences exhibits it as a category of fibrant objects in the sense of~\cite{cisinski}*{Definition 7.5.7}. Since the duality preserves weak equivalences and switches between inflations and deflations, this is also equivalent to saying that the collection of inflations and weak equivalences exhibit it as a category of cofibrant objects. 

\begin{example}
\label{example:homotopically-sound}%
If $(\E,\Dual)$ is any exact category with duality, then the exact category with duality and weak equivalences $(\Ch(\E),\Dualb,\qiso)$ is homotopically sound. Indeed, since quasi-isomorphisms are detected in the abelian hull $\A(\E)$ of $\E$, the collection of trivial fibrations (that is, deflations which are also weak equivalences) is closed under base change by the long exact sequence in ($\A(\E)$-valued) homology groups, and the factorisation axiom follows from the mapping cylinder construction.
\end{example}

Now suppose given an exact category with duality and weak equivalences $(\E,\Dual,W)$ such that $\E[W^{-1}]$ is stable. Write $\ovl{\Dual}$ for the induced duality on $\E[W^{-1}]$ and $\QF^{\sym}_{\ovl{\Dual}}$ for the associated symmetric Poincaré structure. Then we have a natural map from (Poincaré) objects in $(\E,\Dual,W)$ to (Poincaré) objects in $(\E[W^{-1}],\QF^{\sym}_{\ovl{\Dual}})$, and similarly from (Poincaré) objects in $\rS^e_n(\E,\Dual,W)$ to (Poincaré) objects in $\rS^e_n(\E[W^{-1}],\QF^{\sym}_{\ovl{\Dual}})$, 
assembling to give a commutative square
\begin{equation}
\label{equation:comparison-square}%
\begin{tikzcd}
{|\Poi(\rS^e(\E,\Dual,W))|} \ar[r]\ar[d] & {|\Poinc(\rS^e(\E[W^{-1}],\QF^{\sym}_{\ovl{\Dual}}))|} \ar[d]\ar[r,"\simeq"] & {|\Poinc(\Q(\E[W^{-1}],\QF^{\sym}_{\ovl{\Dual}}))|} \ar[d] \\
{|\cor(\rS^e(\E,W))|}  \ar[r] & {|\core(\rS^e(\E[W^{-1}]))|}  \ar[r,"\simeq"] & {|\core(\Q(\E[W^{-1}]))|} \ ,
\end{tikzcd}
\end{equation}
where we have used Proposition~\reftwo{proposition:S-is-Q} and Remark~\reftwo{remark:recovers-spitzweck} to identify the $\Q$- and edgewise subdivided $\rS$-construction.
The induced map on vertical fibres then gives, in light of Corollary~\reftwo{corollary:GWspace-is-fiber-Pn-to-Cr}, a map
\[
\GWspace(\E,\Dual,W) \to \GWspace(\E[W^{-1}], \QF_{\ovl{\Dual}}^\sym).
\]

\begin{proposition}
\label{proposition:comp-schlichting-space}%
Let $\E$ be a $\ZZ[1/2]$-linear exact category with duality $\Dual$ and weak equivalences $W$. Assume that $(\E,\Dual,W)$ is homotopically sound and that $\E[W^{-1}]$ is stable. Then the above comparison map of $\GW$-spaces is an equivalence.
\end{proposition}

By Example~\reftwo{example:homotopically-sound} we consequently conclude:

\begin{corollary}
\label{corollary:comp-schlichting}%
Let $\E$ be an exact category with duality $\Dual$. 
Then we have a canonical equivalence of Grothendieck-Witt spaces
\[
\GWspace(\E,\Dual) \simeq \GWspace(\Ch(\E),\qiso,\Dualb) \simeq \GWspace(\Der^{\mathrm{b}}(\E),\QF^{\sym}_{\ovl{\Dual}}) .
\]
\end{corollary}

In \cite{schlichting-derived}*{Section 5} Schlichting also defines a Grothendieck-Witt spectrum in the context of $\ZZ[1/2]$-linear dg-categories with duality. In particular, this can be applied to $\Ch(\E)$, which admits a natural structure of a dg-category. 
One then considers the shifted duality \((\Dualb)\qshift{n}\colon \Ch(\E) \to \Ch(\E)\op\) obtained by post-composing \(\Dualb\) with the \(n\)-th suspension functor. 
Schlichting's Grothendieck-Witt (pre-)spectrum \(\GW(\Ch(\E),\Dualb)\) is then defined as the sequence of spaces
\[
\GW(\Ch(\E),\Dualb) = \big(|\Poi(\Ch(\E),\qiso,\Dualb)|, |\Poi(\rS^{e}(\Ch(\E),\qiso,\Dualb\qshift{1}))|, |\Poi((\rS^{e})^{(2)}(\Ch(\E),\qiso,\Dualb\qshift{2}), \ldots \big),
\]
with bonding maps induced by the functor \(\Poi(\Ch(\E),\qiso,\Dualb\qshift{n})\to \Poi(\rS^{e}_1(\Ch(\E),\qiso,\Dualb\qshift{n+1}))\)
as in the construction of the Grothendieck-Witt spectrum of a stable \(\infty\)-category using the iterated \(\Q\)-construction. Using Proposition~\reftwo{proposition:S-is-Q} and Remark~\reftwo{remark:recovers-spitzweck} to identify the $\Q$- and edgewise subdivided $\rS$-construction the natural maps
\[
|\Poi((\rS^{e})^{(n)}(\Ch(\E),\qiso,(\Dualb)\qshift{n}))| \to |\Poinc((\rS^{e})^{(n)}(\Der^{\mathrm{b}}(\E),(\Dualb)\qshift{n}))| = |\Poinc(\Q^{(n)}(\Der^{\mathrm{b}}(\E),(\QF^{\sym}_{\ovl{\Dual}})\qshift{n}))|
\]
then fit to give a natural map
\[
\GW(\Ch(\E),\qiso,\Dualb) \to \GW(\Ch(\E)[\qiso^{-1}],\QF^{\sym}_{\ovl{\Dual}}) .
\]

\begin{corollary}
\label{corollary:compare-schlichting-spectrum}%
The above comparison map among \(\GW\)-spectra is an equivalence.
\end{corollary}

The remainder of this section is devoted to the proof of Proposition~\reftwo{proposition:comp-schlichting-space} and Corollary~\reftwo{corollary:compare-schlichting-spectrum}.

\begin{lemma}
\label{lemma:passagetoSn}%
Let \(\E\) be an exact category with weak equivalences which is homotopically sound. Then for every \(n \geq 0\), the exact category with weak equivalences \(\rS_n\E\) is homotopically sound and the natural functor \(\rS_n\E \to \rS_n(\E[W^{-1}])\) exhibits the \(\infty\)-category \(\rS_n(\E[W^{-1}])\) as the localisation of \(\rS_n\E\) with respect to the pointwise weak equivalences.
\end{lemma}

\begin{proof}
We note that \(\rS_n\E\) is equivalent to the category of sequences of inflations
\[
X_1 \hrar X_2 \hrar ... \hrar X_n
\]
with the inflations in \(\rS_n\E\) being the Reedy inflations. It is then standard that if \(\E\) is category of cofibrant objects, then the collection of Reedy inflations exhibit \(\rS_n\E\) as a category of cofibrant objects, see e.g.~\cite{cisinski}*{Theorem 7.4.20 \& Example 7.5.8}. On the other hand, \(\rS_n(\E[W^{-1}])\) is equivalent to the \(\infty\)-category \(\Fun(\Delta^n,\E[W^{-1}])\) of sequences of \(n-1\) composable maps in \(\E[W^{-1}]\). The fact that \(\Fun(\Delta^n,\E[W^{-1}])\) is the \(\infty\)-categorical localisation of the category of Reedy sequences of inflations then follows from \cite{cisinski}*{Theorems 7.5.18 \& 7.6.17}.
\end{proof}

\begin{proposition}
\label{proposition:homotopicallysoundpoinc}%
Suppose that \(\E\) is a $\ZZ[1/2]$-linear exact category with duality $\Dual$ and weak equivalences $W$. Suppose that $(\E,W)$ is a homotopically sound and that $\E[W^{-1}]$ is stable. Then the natural map
\[
|\Poi(\E,W,\Dual)| \to \Poinc(\E[W^{-1}],\QF_{\ovl{\Dual}}^\sym)
\]
is an equivalence, where $\ovl{\Dual}$ denotes the induced duality on $\E[W^{-1}]$. 
\end{proposition}
\begin{proof}
Put $\E[W^{-1}] = \E_\infty$ for readability. Consider then the diagram
\begin{equation}
\label{equation:square-comp-sch}%
\begin{tikzcd}
\Poi(\E,W,\Dual) \ar[r] \ar[d] & \Poinc(\E_\infty,\QF_{\ovl{\Dual}}^\sym) \ar[d] \\
W \ar[r] & \core\E_\infty
\end{tikzcd}
\end{equation}
where the vertical functors are forgetful. By~\cite{cisinski}*{Corollary 7.6.9}, the bottom horizontal map becomes an equivalence upon realisation. By Quillen's Theorem B, it hence suffices to show that for every \(X \in \E\), the map
\[
\begin{tikzcd}
\Poi(\E,W,\Dual) \times_{W} W_{/X} \ar[r] & \Poinc(\E_{\infty},\QF_{\ovl{\Dual}}^{\sym}) \times_{\core\E_{\infty}} (\core\E_{\infty})_{/X} 
\end{tikzcd}
\]
is an equivalence after realisation. Let \(\I_X \subseteq W_{/X}\) be the full subcategory spanned by the deflations \(Y \twoheadrightarrow X\) that are also weak equivalences. We claim that the map
\begin{equation}
\label{equation:vertical-commas-2}%
i \colon \Poi(\E,W,\Dual) \times_{W} \I_X \longrightarrow  \Poi(\E,W,\Dual) \times_{W} W_{/X} 
\end{equation}
induces an equivalence on realisations. To see this, let \(X \to X^I \to X \times X\) be a path object for \(X\), whose existence is guaranteed by our assumption that \(\E\) is a category of fibrant objects with respect to deflations. Construct a functor
\[
q\colon \Poi(\E,W,\Dual) \times_{W} W_{/X} \to \Poi(\E,W,\Dual) \times_{W} \I_X
\]
by sending \((q \colon Y \to \Dual Y, Y {\to} X)\) to \((\Dual(\pr) \circ q \circ \pr, Y \times_{X} X^{I} {\twoheadrightarrow} X)\), where \(\pr\colon Y \times_{X} X^{I} \to Y\) is the projection to the first component. The natural map \((q,Y {\to} X) \to (\Dual(\pr) \circ q \circ \pr, Y \times_{X} X^{I})\) induced by the structure map \(X \to X^I\) then determines natural transformations \(\id \Rightarrow q \circ i\) and \(\id \Rightarrow i \circ q\), showing that~\eqreftwo{equation:vertical-commas-2} is an equivalence after realisation. It will hence suffice to show that the map
\begin{equation}
\label{equation:vertical-commas-3}%
\begin{tikzcd}
{|\Poi(\E,W,\Dual) \times_{W} \I_X|} \ar[r,"i"] &  \Poinc(\E_{\infty},\QF^{\sym}_{\ovl{\Dual}}) \times_{\core\E_{\infty}} (\core\E_{\infty})_{/X}
\end{tikzcd}
\end{equation}
is an equivalence. We now observe that the left vertical map in~\eqreftwo{equation:square-comp-sch} is a right fibration classified by the functor \(X \mapsto \Hom_{W}(X,\Dual X)^{C_2} \simeq \Hom_{W}(X,\Dual X)^{\hC}\); recall that $\Set \subset \Sps$ is closed under limits. Similarly, the right vertical map is classified by \(X \mapsto \Map_{\core\E_\infty}(X,\ovl{\Dual} X)^{\hC}\); since $\core(\E_\infty)_{/X}$ is contractible, we do not need the full statement here, but rather only that the fibre of $\Poinc(\E_\infty,\QF_{\ovl{\Dual}}^{\sym}) \rightarrow \core(\E_\infty)$ over a point $X$ is given by $\Map_{\core\E_\infty}(X,\ovl{\Dual} X)^{\hC}$. This follows from the general fact that for a $\Ct$-space $X \in \Sps^\hC$ the fibre of $X^\hC \rightarrow X$ over some $x \in X$ si either empty or may be computed as $\Map_x(\mathrm S^\sigma, X)^\hC$ from the fibre sequence 
\[
\Map_x(\mathrm S^\sigma, X) \longrightarrow \Map(*, X) \longrightarrow \Map(\Ct,X).
\]
Since total spaces of right fibrations are given as the opposites of the colimits in $\Cat$ of their classified functors by \cite{HTT}*{Corollary 3.3.4.6}, and thus their realisation as the colimits in $\Sps$, we may identify~\eqreftwo{equation:vertical-commas-3} with the natural map 
\[
\colim_{[Y {\twoheadrightarrow} X] \in \I_X\op }\Hom_{W}(Y,\Dual Y)^{\hC} \longrightarrow \Hom_{\core\E_{\infty}}(X,\ovl{\Dual} X)^{\hC} 
\]
in $\Sps$. Now, since $2$ is assumed invertible in \(\E\), multiplication by $2$ acts invertibly on the $\Einf$-groups \(\Hom_{\E}(Y,\Dual Y)\), which is of course an ordinary abelian group, and \(\Hom_{\E_{\infty}}(Y,\ovl{\Dual} Y)\). It follows that the norm map identifies their homotopy fixed points with their homotopy orbits (in $\Einf$-groups). In particular, the homotopy fixed point functor commutes with colimits of $\E_\infty$-groups in which $2$ is invertible. Now note that the category \(\I_X\) admits products (given by fibre products in \(\E\) over \(X\)), and so \(\I_X\op\) is sifted in the $\infty$-categorical sense. Since the forgetful functor from $\E_\infty$-groups to spaces preserves sifted colimits by \cite{HA}*{Proposition 1.4.3.9}, we conclude that
\[
\colim_{[Y {\twoheadrightarrow} X] \in \I_X\op }\Hom_{W}(Y,\Dual Y)^{\hC} \simeq \left[\colim_{[Y {\twoheadrightarrow} X] \in \I_X\op }\Hom_{W}(Y,\Dual Y)\right]^\hC,
\]
so it suffices to establish that 
\[
\begin{tikzcd}[cramped, sep=small]
\displaystyle\colim_{[Y \to X] \in \I_X\op }\Hom_{W}(Y,\Dual Y) \ar[rr] && \Hom_{\core\E_{\infty}}(X,\ovl{\Dual} X)
\end{tikzcd}
\]
is an equivalence. Since $\I_X\op$ is sifted, the map induced by the diagonal
\[
\begin{tikzcd}
\displaystyle\colim_{[Y {\twoheadrightarrow} X] \in \I_X\op }\Hom_{W}(Y,\Dual Y) \ar[r] & \displaystyle\colim_{[Y {\twoheadrightarrow} X,Z {\twoheadrightarrow} X] \in \I_X\op \times \I_X\op } \hspace{-3ex} \Hom_{W}(Y,\Dual Z)
\end{tikzcd}
\]
is an equivalence. We have thus reduced to showing that the natural map
\[
\begin{tikzcd}[cramped, sep=small]
\displaystyle\colim_{[Y {\twoheadrightarrow} X,Z {\twoheadrightarrow} X] \in \I_X\op \times \I_X\op }\hspace{-4ex} \Hom_{W}(Y,\Dual Z) \ar[r] 
\ar[r, phantom, "\simeq" description] 
& \Hom_{\Core\E_{\infty}}(X,\ovl{\Dual} X) 
\end{tikzcd}
\]
is an equivalence. Since the duality switches inflations and deflations, we may rewrite this as
\[
\begin{tikzcd}[cramped, sep=small]
\displaystyle\colim_{[Y {\twoheadrightarrow} X,\Dual X {\hrar} Z'] \in \I_X\op \times \J_{\Dual X}} \hspace{-3ex} \Hom_{W}(Y,Z') \ar[rr] && \Hom_{\core\E_{\infty}}(X,\ovl{\Dual} X) ,
\end{tikzcd}
\]
where $\J_{\Dual X}$ denotes the subcategory of $W_{\Dual X/}$ spanned by the inflations. Now, this last map is an equivalence on general grounds; it is one formula for derived mapping spaces in categories of fibrant/cofibrant objects~\cite{cisinski-k-theorie}*{Proposition 3.23}.
\end{proof}

\begin{proof}[Proof of Proposition~\reftwo{proposition:comp-schlichting-space}]
Proposition~\reftwo{proposition:homotopicallysoundpoinc} and Lemma~\reftwo{lemma:passagetoSn} together imply that the horizontal maps in the left square in~\eqreftwo{equation:comparison-square} are equivalences, and so induce an equivalence on vertical fibres.
\end{proof}

\begin{proof}[Proof of Corollary~\reftwo{corollary:compare-schlichting-spectrum}]
Combining Proposition~\reftwo{proposition:homotopicallysoundpoinc} and Lemma~\reftwo{lemma:passagetoSn} we get that each of the maps
\[
|\Poi((\rS^{e})^{(n)}(\Ch(\E),\qiso,(\Dualb)\qshift{n}))| \to |\Poi((\rS^{e})^{(n)}(\Der^{\mathrm{b}}(\E),(\QF^{\sym}_{\ovl{\Dual}})\qshift{n}))|
\]
is an equivalence of spaces, and so induce an equivalence on the level of pre-spectra, and eventually of spectra.
\end{proof}

\end{appendices}

\end{document}